\documentclass[11pt]{amsart}
\usepackage[margin=1in]{geometry}
\usepackage[T1]{fontenc}
\usepackage{amsmath,amssymb,amsthm,mathtools}
\usepackage{enumitem}
\usepackage{needspace}
\usepackage{longtable,array,booktabs,tabularx}
\usepackage{microtype,mathrsfs}
\usepackage{xcolor}
\usepackage{tikz}
\usetikzlibrary{arrows.meta,positioning}
\usepackage{mdframed}
\newmdenv[backgroundcolor=black!5,linecolor=black!65,linewidth=0.4pt,roundcorner=3pt,skipabove=8pt,skipbelow=8pt,innerleftmargin=6pt,innerrightmargin=6pt,innertopmargin=3pt,innerbottommargin=3pt]{mainresult}
\definecolor{citationblue}{RGB}{0,0,255}
\definecolor{navigationred}{RGB}{255,0,0}
\usepackage[colorlinks=true,linktoc=page,linkcolor=navigationred,citecolor=citationblue,urlcolor=citationblue]{hyperref}
\numberwithin{equation}{section}
\pdftrailerid{}
\newtheorem{theorem}{Theorem}[section]
\newtheorem{proposition}[theorem]{Proposition}
\newtheorem{lemma}[theorem]{Lemma}
\newtheorem{corollary}[theorem]{Corollary}
\newcommand{\restatedname}{}
\newtheorem*{restatedbody}{\restatedname}
\newenvironment{restated}[1]{\renewcommand{\restatedname}{#1}\begin{restatedbody}}{\end{restatedbody}}
\theoremstyle{definition}
\newtheorem{definition}[theorem]{Definition}
\theoremstyle{remark}
\newtheorem{remark}[theorem]{Remark}
\title[Uniform large-genus volume asymptotics]{Uniform large-genus asymptotics\\ of Weil--Petersson volumes}
\author[Q. Wang et al.]{Qihang Wang\textsuperscript{1,\dag}}
\author[]{Zhiyi Li\textsuperscript{2,3,\dag}}
\author[]{Bolin Qiu\textsuperscript{4}}
\author[]{Shuaijia Wang\textsuperscript{4}}
\author[]{Yuan Zhong\textsuperscript{5}}
\author[]{Kun Chen\textsuperscript{6,*}}
\author[]{Renqiu Xia\textsuperscript{4,*}}
\thanks{\textsuperscript{\dag}Qihang Wang and Zhiyi Li contributed equally to this work and share first authorship.}
\thanks{\textsuperscript{*}Corresponding authors: Kun Chen (\href{mailto:chenkun@itp.ac.cn}{chenkun@itp.ac.cn}) and Renqiu Xia (\href{mailto:xiarenqiu@sjtu.edu.cn}{xiarenqiu@sjtu.edu.cn}).}
\thanks{\textsuperscript{1}School of Mathematical Sciences, Peking University, Beijing 100871, China.}
\thanks{\textsuperscript{2}Department of Modern Physics, University of Science and Technology of China, Hefei, Anhui 230026, China.}
\thanks{\textsuperscript{3}CAS Key Laboratory of Theoretical Physics, Institute of Theoretical Physics, Chinese Academy of Sciences, Beijing 100190, China.}
\thanks{\textsuperscript{4}School of Artificial Intelligence, Shanghai Jiao Tong University, Shanghai, China.}
\thanks{\textsuperscript{5}Wuhu Qiheng Co., Ltd., Wuhu, China.}
\thanks{\textsuperscript{6}Institute of Theoretical Physics, Chinese Academy of Sciences, Beijing 100190, China.}
\date{}
\subjclass[2020]{Primary 32G15; Secondary 11F72, 14H10, 14H81, 41A60, 58J50, 60D05, 60F05}
\keywords{Weil--Petersson volume, moduli space of curves, intersection number, large genus, random hyperbolic surface, spectral gap, local limit theorem}
\hypersetup{pdftitle={Uniform large-genus asymptotics of Weil-Petersson volumes},pdfauthor={Qihang Wang, Zhiyi Li, Bolin Qiu, Shuaijia Wang, Yuan Zhong, Kun Chen, Renqiu Xia}}
\begin{document}
\begin{abstract}
We prove an expansion in $1/g$ of the Weil--Petersson volumes
$V_{g,n}$ that is uniform over all cusp counts, from fixed $n$ to
arbitrarily large $n/g$. An explicit leading term and single-factorial
remainder bounds give exponentially small relative errors after
truncation at a suitable order proportional to $g$. A positive
recursion along the string flow connects these volumes to pure
$\psi$-class intersection numbers, where Painlev\'e~I fixes the
normalization; uniform coefficient extraction then recovers individual
volumes. We determine Zograf's constant $C=1/\sqrt\pi$ and resolve
Mirzakhani's two questions on volume ratios for growing cusp counts
and adjacent coefficient ratios uniformly over all degrees.

Via Mirzakhani's integration formula, the same volume estimates
control random surfaces at two distinct scales. At
$n/\sqrt g\to c\in(0,\infty)$, the entire positive Laplace spectrum
below $1/4$ converges to an explicit Poisson point process generated
by geodesics cutting off two cusps. Together with Huang's theorem,
this characterizes the transition for $n\ge1$: the spectral gap tends
in probability to $1/4$ exactly when $n/\sqrt g\to0$, and to $0$
exactly when $n/\sqrt g\to\infty$.

At density $n/g\to\lambda\in(0,\infty)$, we identify the local limit
around a typical point; its cusp-filling is conformally a disc.
As $\lambda\to\infty$, these limits approach the infinite
Weil--Petersson sphere of Budd and Curien. We prove their parabolicity
conjecture and, by a growing-radius coupling, answer their scaling-limit
question: after deep cusps are trimmed, the rescaled infinite sphere
converges to the Brownian plane.
\end{abstract}
\maketitle
\begingroup\renewcommand{\addcontentsline}[3]{}\section*{Acknowledgements}\endgroup
We thank Hao Xu for helpful correspondence on the work of Liu and Xu, the uniformity ranges of the two deformation arguments, and the determination of Zograf's constant. We acknowledge GeWu SciLab for providing the collaborative research environment in which this project was developed.

Kun Chen is supported by the Strategic Priority Research Program of the Chinese Academy of Sciences (Grant No.~XDB1680102). Renqiu Xia is supported by the National Natural Science Foundation of China (Grant No.~6260070460).
\clearpage
\setcounter{tocdepth}{1}
\tableofcontents
\clearpage
\section{Introduction}

Let $\mathcal M_{g,n}(\mathbf L)$ be the moduli space of hyperbolic
surfaces of genus $g$ with $n$ labeled geodesic boundary components of
lengths $\mathbf L=(L_1,\ldots,L_n)$, where $2g-2+n>0$ and a boundary
component of length zero is a cusp. Its Weil--Petersson volume, in the
orbifold convention, is
\[
 V_{g,n}(\mathbf L)=\frac{1}{(3g-3+n)!}
 \int_{\mathcal M_{g,n}(\mathbf L)}\omega_{\mathrm{WP}}^{3g-3+n},
 \qquad V_{g,n}=V_{g,n}(0,\ldots,0).
\]
These volumes connect intersection theory on moduli space with the
geometry of random hyperbolic surfaces: their coefficients in the
boundary lengths are intersection numbers, and through Mirzakhani's
integration formula their ratios give expected counts of geometric
configurations~\cite{MirInv,MirJAMS,Mir}. See~\cite{Do,MirICM,Wright,ChenMullane}
for surveys.

Two classical asymptotic theories cover opposite ends of the problem:
Mirzakhani and Zograf obtained an expansion as $g\to\infty$ with $n$
fixed~\cite{MZ}, and Kaufmann--Manin--Zagier and Manin--Zograf obtained
the large-$n$ asymptotics at fixed genus~\cite{KMZ,ManZ}. In the first
expansion the coefficient of order $g^{-i}$ is a polynomial in $n$ of
degree $2i$, with leading term $(-1)^in^{2i}/(i!\,(2\pi^2)^i)$. These
leading terms alone sum to $\exp(-n^2/(2\pi^2g))$, so the effective
parameter is $n^2/g$: once $n\asymp\sqrt g$, terms of every fixed order
contribute, and powers of $1/g$ no longer order the corrections
asymptotically. The two endpoint theories alone therefore
do not describe arbitrary simultaneous growth of genus and cusp count.

We extend the large-genus theory across the full range of cusp
densities. A single real parameter $Y$, the \emph{cusp-density
parameter}, runs over a compact interval $[0,y_0]$ as $n/g$ runs over
$[0,\infty]$ and locates the saddle point of the cusp generating
function. Three ideas
make the expansion uniform in $n$. Moving along the string flow
replaces the alternating sum that expresses Weil--Petersson volumes
through $\psi$ intersection numbers by a positive recursion, whose
kernels form a compact family over $[0,y_0]$. The dilaton and KdV
equations, controlled by this positivity, determine the dependence on
$Y$; the remaining normalization is fixed at the endpoint $Y=y_0$,
where only pure $\psi$ intersection numbers enter and Painlev\'e~I
supplies their asymptotics. Finally, each generating function is a series
with nonnegative coefficients times an explicit function raised to a
power proportional to the genus; this factorization converts its
asymptotics into asymptotics of individual coefficients, uniformly
in $n$. The same parameter connects the volume
problem to the geometry of random surfaces: the scales $n\asymp\sqrt g$
and $n\asymp g$ mark changes in the low spectrum of these surfaces and
in their geometry near a typical point.

\textbf{Volumes (Part~\ref{part:volumes}).}
Theorem~\ref{thm:intro-absolute} is a single large-genus expansion
\[
 V_{g,n}=\widehat V_{g,n}\Bigl(1+\frac{A_1(g,n)}{2g-1}
   +\frac{A_2(g,n)}{(2g-1)^2}+\cdots\Bigr)
\]
valid \emph{uniformly over every cusp count}, with an explicit leading
term $\widehat V_{g,n}$, coefficients bounded uniformly in $n$, and
single-factorial remainder bounds. In the fixed-cusp
limit it identifies the constant conjectured by Zograf; at every
density it gives the volume-ratio profiles, and for the one-boundary
polynomial it gives the adjacent-coefficient asymptotics in all
degrees. These results resolve the two questions posed by
Mirzakhani~\cite[Section~1.4]{Mir}: the asymptotics of
$V_{g,n}/V_{g-1,n+2}$ when $n=n(g)\to\infty$, and the ratios of
adjacent coefficients of $V_{g,1}(L)$ when the degree grows with the
genus. Theorems~\ref{thm:volume} and~\ref{thm:coefficient} give the
respective answers uniformly over all cusp counts and all adjacent
degrees. Liu and Xu subsequently determined
Zograf's constant by a different deformation~\cite{LXConstant}; the two
approaches are compared at the beginning of Part~\ref{part:volumes}.

\textbf{Spectrum (Part~\ref{part:spectra}).}
The threshold $n\asymp\sqrt g$, at which the fixed-cusp expansion
breaks down, is also a spectral phase transition. Both are governed by
$n^2/g$: the expected number of geodesics of bounded length that cut
off two cusps is of this order
(Subsection~\ref{subsec:intro-geometry}). Huang showed that the
spectral gap of a Weil--Petersson random surface tends to $1/4$ when
$n=o(\sqrt g)$~\cite{HuangGap}. At the critical scale the full
positive spectral point process below $1/4$ converges to a Poisson
process of explicit model energies, one for each short two-cusp pair of pants
(Theorem~\ref{thm:critical-eigenvalues}); the gap therefore has a
nondegenerate random limit, and tends to zero above the scale. Above
the critical scale the lowest eigenvalues follow rescaled Poisson laws
in specified ranges, and fixed-energy windows carry deterministic
limiting spectral measures. The volume theorem supplies the statistics
of short curves; the spectral work is to show that the states these
curves create exhaust the low spectrum.

\textbf{Local geometry (Part~\ref{part:geometry}).}
At positive density $n/g\to\lambda\in(0,\infty)$ the same volume ratios,
with based holonomies retained in the integration formula, determine
the surface seen from a uniformly chosen point and the conformal type
of its cusp-filling. As $\lambda\to\infty$ this
family reaches the area-rooted infinite WP sphere of Budd and
Curien~\cite{BC}. We prove that this endpoint is parabolic, with
recurrent Brownian motion and ergodic geodesic flow, as they
conjectured~\cite[Section~7.3]{BC}, whereas at every finite density
the cusp-filled surface is conformally a disc. A growing-radius
coupling of their tree models shows that the endpoint surface,
rescaled and with its deep cusps trimmed, converges to the Brownian
plane, answering another of their questions.

The common step is geometric integration: a cut expresses the law of a
finite configuration through volumes of its complementary surfaces.
For spectra, these configuration counts must be matched to the entire
low spectral subspace; for local geometry, the relative positions of
the based loops must also be retained. Figure~\ref{fig:paper-map} shows
how these two uses of the volume theorem lead to the applications;
Table~\ref{tab:intro-questions} lists the questions from the
literature that are answered here.
\begin{figure}[htbp]
\centering
\begin{tikzpicture}[
  box/.style={draw=black!45,line width=.45pt,rounded corners=2pt,
    text=black,align=center,inner sep=6pt,font=\small},
  core/.style={box,draw=black!65,fill=black!3},
  route/.style={-{Latex[length=2mm,width=1.3mm]},
    draw=black!70,line width=.65pt},
  note/.style={font=\footnotesize,align=center,text=black,
    fill=white,inner sep=2pt}]
\node[core,text width=15.2cm,minimum height=1.65cm] (vol) at (0,0)
  {\textbf{I. Uniform volume asymptotics across all cusp densities}\\[3pt]
   Positive string-flow recursion $\longrightarrow$ Painlev\'e~I normalization
   $\longrightarrow$ coefficient extraction\\[3pt]
   \emph{Common output: volume-ratio profiles and boundary factors}\\
   Theorems~\ref{thm:intro-absolute}, \ref{thm:volume};
   Proposition~\ref{prop:bounded-boundaries}};

\node[box,text width=15.2cm,minimum height=1.05cm,below=4mm of vol] (cut)
  {\textbf{Mirzakhani's integration formula}\\[2pt]
   Volume ratios weight cut surfaces and their geometric configurations.};
\draw[route] (vol.south)--(cut.north);

\node[core,text width=7.0cm,minimum height=4.6cm,anchor=north]
  (spec) at ([xshift=-4.05cm,yshift=-8mm]cut.south)
  {\textbf{II. Spectral laws}\\[4pt]
   Short-curve counts and orbital integrals\\[4pt]
   $\downarrow$\\[1pt]
   Local models and spectral exhaustion\\
   by min--max and trace estimates\\[4pt]
   $\downarrow$\\[1pt]
   Critical and microscopic eigenvalue laws;\\
   spectral measures in fixed windows\\[4pt]
   Sections~\ref{sec:critical-eigenvalues}--\ref{sec:fixed-energy-spectrum}};
\node[core,text width=7.0cm,minimum height=4.6cm,anchor=north]
  (local) at ([xshift=4.05cm,yshift=-8mm]cut.south)
  {\textbf{III. Local geometry}\\[4pt]
   Carrier integrals with based holonomies\\[4pt]
   $\downarrow$\\[1pt]
   Pointed measured balls\\
   and planar-cap geometry\\[4pt]
   $\downarrow$\\[1pt]
   Local family at $n/g\to\lambda\in(0,\infty)$;\\
   conformally a disc after filling cusps\\[4pt]
   Section~\ref{sec:local-conformal-phase}};
\draw[route] (cut.south)--++(0,-4mm) coordinate (fork) -| (spec.north);
\draw[route] (fork) -| (local.north);

\node[box,text width=7.0cm,minimum height=2.2cm,below=10mm of local] (sphere)
  {\textbf{Infinite-density endpoint}\\[3pt]
   Area-rooted infinite WP sphere\\
   of Budd--Curien: parabolic,\\
   conformally $\mathbb C$ after filling cusps};
\draw[route] (local.south)--(sphere.north)
  node[midway,right,note] {$\lambda\to\infty$\\local convergence};
\node[box,text width=7.0cm,minimum height=2.2cm]
  (tree) at (spec.center |- sphere.center)
  {\textbf{Additional input for scaling}\\[3pt]
   Budd--Curien trees and sphere limit\\[2pt]
   Growing-radius coupling\\
   Measured tangent coupling};

\node[core,text width=15.2cm,minimum height=1.45cm,anchor=north]
  (plane) at ([xshift=-4.05cm,yshift=-8mm]sphere.south)
  {\textbf{Brownian-plane scaling at the planar endpoint}\\[3pt]
   Cusp-trimmed measured GHP convergence; full-area GP convergence\\[3pt]
   Section~\ref{sec:brownian-plane}};
\draw[route] (tree.south)--(tree.south |- plane.north);
\draw[route] (sphere.south)--(sphere.south |- plane.north);
\end{tikzpicture}
\caption{One volume theory, two geometric uses. Both branches use
volume ratios through geometric integration and then the additional
arguments inside their boxes. The labelled arrow is the proved
fixed-radius local limit as $\lambda\to\infty$. The two arrows into
the bottom box show the endpoint surface and the separate coupling
argument used for its Brownian-plane scaling limit. Part~\ref{part:geometry}
uses no spectral theorem from Part~\ref{part:spectra}. The all-degree
coefficient theorem in Part~\ref{part:volumes} also uses harmonic
comparison directly.}
\label{fig:paper-map}
\end{figure}

\subsection{The cusp-density parameter}\label{subsec:parameter}
All our formulas are expressed through a single real parameter. We define it first and then explain heuristically why it is the natural parameter; the proofs come later.

\textbf{The genus-zero relation.} Let
\[
 X(y)=\sqrt y\,J_1(2\sqrt y)=y-\frac{y^2}{2}+\frac{y^3}{12}-\cdots,
 \qquad y_0=\frac{j_{0,1}^2}{4},\qquad x_c=X(y_0),
\]
where $j_{0,1}$ is the first positive zero of $J_0$. Since $X'(y)=J_0(2\sqrt y)$, the function $X$ increases from $0$ to $x_c$ on $[0,y_0]$; write $Y(x)$ for the inverse function. Consider the normalized cusp generating functions
\begin{equation}
 \Phi_g(x)=\sum_{n\ge0}\frac{V_{g,n}}{(2\pi^2)^{3g-3+n}}\,\frac{x^n}{n!},
 \label{eq:intro-Phi}
\end{equation}
with the convention that terms with $2g-2+n\le0$ are omitted. In genus zero, $\Phi_0''(x)=Y(x)$; this is the classical genus-zero computation (cf.~\cite{KMZ,ManZ}), and in our normalization it follows from the string equation (see Step~1 of the proof of Theorem~\ref{thm:pants-poisson} in Appendix~\ref{app:short-curves}). We call $x$ the \emph{cusp weight}: it weights each cusp in~\eqref{eq:intro-Phi}. The cusp weight $x$ and the parameter $Y=Y(x)\in[0,y_0)$ determine each other, and we use whichever is more convenient.

\textbf{The parameter functions.} For $0\le Y<y_0$, the entire function
\[
 K_Y(z)=\int_0^1X'\!\Bigl(Y+\frac z2(1-s^2)\Bigr)\,\mathrm ds
\]
has only positive simple zeros (Section~\ref{the-positive-compact-family-of-kernels}); it is the kernel of the recursion in Section~\ref{exact-amplitudes-and-the-harmonic-model}. Write its first zero as $2\Delta(Y)$ and put $T(Y)=Y+\Delta(Y)$. The substitution $v=T-\Delta s^2$ shows that $T=T(Y)$ is the unique root in $(y_0,\pi^2/4]$ of
\begin{equation}
 \int_Y^T\frac{X'(v)}{\sqrt{T-v}}\,\mathrm{d}v=0.
 \label{eq:SP1}
\end{equation}
Finally define
\begin{equation}
 I(Y)=\int_Y^{T(Y)} X'(v)\sqrt{T(Y)-v}\,\mathrm{d}v,
 \qquad
 \lambda(Y)=\frac{2X(Y)\sqrt{\Delta(Y)}}{I(Y)}.
 \label{eq:SP2}
\end{equation}
The next paragraph motivates these definitions. In Section~\ref{the-stationary-parameter} we show that $I>0$ and that $\lambda$ is a strictly increasing bijection from $[0,y_0)$ onto $[0,\infty)$. At $Y=0$ one has $\Delta=T=\pi^2/4$, $I=\pi/4$ and $\lambda=0$; as $Y\uparrow y_0$, $\Delta\to0$ and $\lambda\to\infty$.

\textbf{Why $\lambda$ is the cusp density.} Let $h=2g-1$ and $x=X(Y)$. We shall prove (Proposition~\ref{prop:prefactor}) that, uniformly for $0\le Y<y_0$,
\begin{equation}
 \Phi_g'(x)=\frac{\Gamma(h)\,\varpi(Y)}{\pi\,4^g\sqrt{\pi h}}\,
 \bigl(\sqrt2\,I(Y)\bigr)^{-h}\Bigl(1+O\Bigl(\frac{\log g}{g}\Bigr)\Bigr),
 \label{eq:intro-heuristic}
\end{equation}
where $\varpi$ is the explicit positive function in~\eqref{eq:intro-varpi} below. Thus the cusp weight changes the exponential growth only through the \emph{action} $I(Y)$.

Now regard the terms of the series
$\Phi_g'(x)=\sum_{k\ge0}V_{g,k+1}x^k/((2\pi^2)^{3g-2+k}k!)$ as the weights of a
probability distribution on the number $k$ of additional cusps. Its mean is
$x\,\mathrm d\log\Phi_g'/\mathrm dx$. Since $I'(Y)=-X'(Y)\sqrt{\Delta(Y)}$ and
$\mathrm dY/\mathrm dx=1/X'(Y)$, the leading term in~\eqref{eq:intro-heuristic} gives
\[
 x\frac{\mathrm d}{\mathrm dx}\log\bigl(\sqrt2\,I(Y)\bigr)^{-h}
 =h\,\frac{X(Y)\sqrt{\Delta(Y)}}{I(Y)}=\frac{h\lambda(Y)}2 .
\]
So, under the cusp weight $x$, the number of additional cusps is about
$h\lambda(Y)/2\approx g\lambda(Y)$. To recover a single volume $V_{g,n}$, one
chooses $Y$ so that the leading approximation to this mean equals $n-1$, that is,
$\lambda(Y)=2(n-1)/(2g-1)$, and estimates the probability of the value
$k=n-1$ by a local limit theorem. In particular $\lambda(Y)$ is, to leading
order, the cusp density $n/g$, so the profiles may be expressed in either $Y$ or $\lambda(Y)$. This heuristic produces the leading term of
Theorem~\ref{thm:intro-absolute}. As $n/g$ increases from $0$ to $\infty$,
the parameter $Y$ increases from $0$ to $y_0$ and the cusp weight $x$ from
$0$ to $x_c$.

Subsection~\ref{subsec:strategy} explains how endpoint normalization and
Fourier bounds make this heuristic uniform in $n$.

\textbf{A matrix-model interpretation.} The cusp weight is a defect weight in Jackiw--Teitelboim (JT) gravity, and $Y$ locates the edge of the corresponding spectral density. At the endpoint $Y=y_0$ the square-root edge becomes cubic. Subsection~\ref{subsec:jt-edge} gives the precise dictionary with the defect-deformed matrix model of Maxfield--Turiaci and Witten~\cite{MT,WittenDef,ET}. In this language the endpoint corresponds to the $(2,3)$ critical point, where our normalized recursion reduces to pure $\psi$ intersection numbers. The interpretation does not by itself control the remainders.

\subsection{Absolute volumes and Zograf's constant}
For $n\ge1$ put $h=2g-1$ and $k=n-1$, let $Y\in[0,y_0)$ be the solution of $\lambda(Y)=2k/h$, and put $x=X(Y)$. Define
\begin{equation}
 \widehat V_{g,n}=\frac{(2\pi^2)^{3g-3+n}\,\Gamma(h)}{\pi\,4^g\sqrt{\pi h}}\,
 \frac{\varpi(Y)}{(\sqrt2\,I(Y))^{h}}
 \left(\frac{k}{\mathrm e x}\right)^{k}
 \sqrt{\frac{k}{h\,x\Psi'(x)}},
 \label{eq:intro-leading}
\end{equation}
where $\Psi(X(y))=\lambda(y)/2$ for $0\le y<y_0$, and $\varpi(Y)>0$ is given by
\begin{equation}
 \varpi(Y)^2=\frac{I(Y)}{\Delta(Y)^{3/2}}
 \left(\int_0^1(1-s^2)\,\frac{X(v)}{v}\,\mathrm ds\right)^{-1},
 \qquad v=T(Y)-\Delta(Y)s^2;
 \label{eq:intro-varpi}
\end{equation}
here $X(v)/v$ is an entire function of $v$, equal to one at $v=0$. For $n=0$ we put $\widehat V_{g,0}=\Gamma(2g-2)(4\pi^2)^{2g-3}/\sqrt{\pi g}$.

The formula~\eqref{eq:intro-leading} is the heuristic of Subsection~\ref{subsec:parameter} made explicit. Its first two factors are $(2\pi^2)^{3g-3+n}$ times the approximation~\eqref{eq:intro-heuristic} of $\Phi_g'(x)$. The coefficient of $x^k$ is about $\Phi_g'(x)x^{-k}$ times the probability of the center of the distribution. To leading order, $k!$ times that probability is $(k/\mathrm e)^k\sqrt{k/(hx\Psi'(x))}$, where $hx\Psi'(x)$ is the variance scale supplied by the leading action. We use the Poisson point mass $\mathrm e^{-k}k^k/k!$ at the center rather than its Gaussian approximation $(2\pi k)^{-1/2}$, so that the same formula remains correct when the cusp count stays bounded; there the number of additional cusps is approximately Poisson rather than Gaussian. At $k=0$ the last two factors in~\eqref{eq:intro-leading} are interpreted as their limits, which are one; then $Y=0$, $\varpi(0)=1$ and $\widehat V_{g,1}=\Gamma(2g-1)(4\pi^2)^{2g-2}/\sqrt{\pi(g-1/2)}$. The case $n=0$ is reduced to $n=1$ by the dilaton identity~\eqref{eq:TL12} at $p=0$, which for volume polynomials is the relation $\partial_LV_{g,1}(2\pi i)=2\pi i(2g-2)V_{g,0}$ of Do and Norbury~\cite[Theorem~2]{DoNorbury} (Section~\ref{sec:absolute}).

The correction coefficients $A_j(g,n)$ below, with $A_0=1$, are determined by finite recursive rules from the kernels together with one-dimensional integrals in $Y$ (see \eqref{eq:CG3} and \eqref{eq:AV7}); the first is written out in \eqref{eq:CG12}. For $n\ge1$ they are analytic functions of the parameter $Y$, including both endpoints of its interval.

\Needspace{12\baselineskip}
\begin{mainresult}
\begin{theorem}[Uniform absolute expansion]\label{thm:intro-absolute}
There are constants $C,B,c,g_0>0$ such that, for every $g\ge g_0$ and every integer $0\le K\le c(2g-1)$,
\[
 \sup_{n\ge0}\left|\frac{V_{g,n}}{\widehat V_{g,n}}
       -\sum_{j=0}^K\frac{A_j(g,n)}{(2g-1)^j}\right|
 \le\frac{CB^{K+1}(K+1)!}{(2g-1)^{K+1}},
\]
and the coefficients satisfy $|A_j(g,n)|\le CB^jj!$ for all $g\ge g_0$ and $n\ge0$. In particular,
\[
 \sup_{n\ge0}\left|\frac{V_{g,n}}{\widehat V_{g,n}}-1\right|=O(g^{-1}).
\]
Taking $K$ to be a sufficiently small fixed multiple of $g$, rounded down to an integer, gives relative error $O(\mathrm{e}^{-ag})$ for some $a>0$, uniformly in $n$.
\end{theorem}
\end{mainresult}

The uniformity in $n$ is the main point: Theorem~\ref{thm:intro-absolute} replaces the fixed-$n$ series in $n^2/g$ by an expansion in $1/g$ alone. The same estimate applies when $n$ is fixed, when $n$ is proportional to $g$, and when $n/g\to\infty$. In the first regime the saddle is near $x=0$; as $n/g\to\infty$ it tends to the radius of convergence $x_c$ of the fixed-genus generating functions of Manin and Zograf. In this way a single formula joins the two classical regimes. The order $g^{-1}$ is sharp for this choice of leading term (Theorem~\ref{thm:absolute}); including the first correction \eqref{eq:CG12}, or \eqref{eq:AV7} for $n=0$, gives error $O(g^{-2})$. The proof is completed in Section~\ref{sec:absolute}.

For fixed $n$, Mirzakhani and Zograf proved that $V_{g,n}\sim C(2g-3+n)!\,(4\pi^2)^{2g-3+n}/\sqrt g$ for a universal constant $C>0$~\cite[Theorem~1.2]{MZ}. Zograf had conjectured, on the basis of numerical data, that $C=1/\sqrt\pi$~\cite{Zog}; see also \cite[Remark~1.3]{MZ}. The same value appears in physical predictions. Saad, Shenker and Stanford obtained the unmarked large-genus prediction from the JT matrix integral under a resurgence assumption~\cite[Section~5.6]{SSS}. Okuyama and Sakai extracted the constant $1/\sqrt\pi$, together with first corrections for the coefficients of $V_{g,1}(L)$ at each fixed degree, numerically from the KdV equation~\cite[Section~2.4, equations~(2.87)--(2.89)]{OS}. Kimura conjectured the corresponding asymptotics of $V_{g,n}(\mathbf L)$ from JT gravity~\cite[equation~(28)]{Kim}, and Eynard, Garcia-Failde, Gregori, Lewa\'nski and Schiappa derived it from the nonperturbative topological recursion, assuming that the generating series of Weil--Petersson volumes is resurgent~\cite[Section~6.1, equation~(6.7)]{EGGLS}.

Our argument does not assume that the Weil--Petersson genus expansion is resurgent. We determine the normalization constant at the endpoint $Y=y_0$ of the parameter interval, where the problem reduces to pure $\psi$ intersection numbers and only the single sequence $\langle\tau_2^{3g-3}\rangle_g$ enters. An identity of Itzykson and Zuber~\cite{IZ} expresses this sequence through the formal power-series solution of Painlev\'e~I, the string equation of two-dimensional pure gravity in the double-scaling limit~\cite{BK,DS,GM}, whose large-order asymptotics follow from the results of Joshi--Kitaev, Kapaev, and Garoufalidis, Its, Kapaev and Mari\~no~\cite{JK,Kap,GIKM}. The leading term and the first correction of this sequence also appear as a special case of the uniform asymptotics for $\psi$ intersection numbers obtained by Guo, Yang and Zagier~\cite{GYZ}. We derive both from Painlev\'e~I in Subsections~\ref{sec:prefactor-endpoint} and~\ref{sec:endpoint-expansion} and use~\cite{GYZ} only as an independent check. Comparison with the Mirzakhani--Zograf expansion gives the following corollary.
\begin{corollary}[Zograf's constant]\label{cor:constant}
For every fixed $n\ge0$, as $g\to\infty$,
\[
 V_{g,n}=\frac{(2g-3+n)!\,(4\pi^2)^{2g-3+n}}{\sqrt{\pi g}}
 \bigl(1+O_n(g^{-1})\bigr).
\]
\end{corollary}
This is proved at the end of Section~\ref{sec:prefactor}. Together with the bracket estimate \eqref{eq:TL23}, it also gives the formula conjectured by Kimura~\cite[equation~(28)]{Kim}: for fixed $n\ge1$ and nonnegative lengths $L_i=o(\sqrt g)$,
\[
 V_{g,n}(\mathbf L)\sim\sqrt{\frac2\pi}\,(4\pi^2)^{2g-3+n}\,\Gamma\Bigl(2g+n-\frac52\Bigr)\prod_{i=1}^n\frac{\sinh(L_i/2)}{L_i/2}.
\]
The factors at $L_i=0$ are interpreted by continuity. A direct derivation is given in \eqref{eq:boundary-factor} below.

\textbf{Bounded geodesic boundaries.} Proposition~\ref{prop:bounded-boundaries} extends the leading comparison to a fixed number of nonzero boundary lengths in a compact set, uniformly over all cusp counts. Its explicit factor interpolates, as $n/g$ runs from zero to infinity, between the fixed-cusp factor $\sinh(L/2)/(L/2)$ of~\cite{MP} and the fixed-genus, large-cusp factor $I_0(j_{0,1}L/(2\pi))$ of~\cite{HTLargeN}. The interpolation is given in \eqref{eq:boundary-profile}; in this extension both the number and the lengths of the nonzero boundaries remain bounded.

\subsection{Volume ratios}\label{subsec:intro-ratios}
For sufficiently large $g$, define
\[
 B_{g,n}=\frac{V_{g,n}}{V_{g-1,n+2}},\qquad
 C_{g,n}=\frac{V_{g,n+1}}{2gV_{g,n}},\qquad
 D_{g,n}=\frac{V_{g,n}^2}{V_{g,n-1}V_{g,n+1}}.
\]
The first two quantities are defined for $n\ge0$, and the third for $n\ge1$. For each fixed $n$, Mirzakhani proved that $B_{g,n}\to1$ and $C_{g,n}\to4\pi^2$~\cite[Theorem~3.5]{Mir}, and she asked for the asymptotics of $B_{g,n}$ when $n=n(g)\to\infty$~\cite[Section~1.4]{Mir}. An estimate of Mirzakhani and Zograf states that $|(2g-2+n)V_{g,n}/V_{g,n+1}-1/(4\pi^2)|\le cn/(2g-2+n)$ for stable $(g,n)$ with $n\ge1$~\cite[Lemma~5.1]{MZ}; together with the fixed-$n$ result, this gives $C_{g,n}\to4\pi^2$ whenever $n=o(g)$. When $n$ is comparable to or larger than $g$, the asymptotic profiles depend on the density $n/g$, through the functions
\begin{equation}
 \mathcal B(\lambda(Y))=\frac{\pi^2}{4\Delta(Y)},\qquad
 \mathcal C(\lambda(Y))=\frac{\lambda(Y)\pi^2}{X(Y)}.
 \label{eq:SP3}
\end{equation}
They are well defined because $\lambda$ is a bijection, and they take the continuous values $\mathcal B(0)=1$ and $\mathcal C(0)=4\pi^2$. Heuristically, $\mathcal C$ appears because $x\Phi_g''(x)/\Phi_g'(x)$ is the exact mean of the additional-cusp distribution and the saddle is chosen so that its leading approximation equals $n-1$; $\mathcal B$ appears because $\Phi_{g-1}'''(x)/\Phi_g'(x)$ is approximately $8\Delta(Y)$. See~\eqref{eq:TL30} and~\eqref{eq:TL31}.

\Needspace{12\baselineskip}
\begin{mainresult}
\begin{theorem}[Uniform volume ratios]\label{thm:volume}
There are constants $C_0,g_0>0$ such that, for every $g\ge g_0$,
\[
 \sup_{n\ge0}\left|\frac{B_{g,n}}{\mathcal B(n/g)}-1\right|\le\frac{C_0}{g},
 \qquad
 \sup_{n\ge0}\left|\frac{C_{g,n}}{\mathcal C(n/g)}-1\right|\le\frac{C_0}{g},
\]
and
\[
 \sup_{n\ge1}|D_{g,n}-1|\le\frac{C_0}{g}.
\]
In particular, along every path $n/g\to\infty$ with $g\to\infty$,
\[
 B_{g,n}\sim\frac{\pi^2}{3\sqrt5\,j_{0,1}}\sqrt{n/g},\qquad
 C_{g,n}\sim\frac{\pi^2}{x_c}\frac ng.
\]
\end{theorem}
\end{mainresult}
Section~\ref{extracting-fixed-cusp-counts} proves the uniform limits directly from a local limit theorem. The inverse-genus error follows from the sharper coefficient estimates, as shown in Section~\ref{subsec:quantitative-volume-ratios}. The limit of $C_{g,n}$ was known for $n=o(g)$, as described above. For $B_{g,n}$, the absolute estimate of Mirzakhani and Zograf~\cite[Theorem~1.8]{MZ} gives the limit one when $n=o(\sqrt g)$; Gatignol and Thomas~\cite[Theorem~1.5]{GT} extend it to $n=o(g)$. Theorem~\ref{thm:volume} also covers positive and unbounded cusp density and gives relative error $O(g^{-1})$ uniformly in $n$.

The statement on $D_{g,n}$ answers a question of Shen and Wu~\cite[Question~3]{SW}, which asks whether $D_{g,n}\to1$ as $g+n\to\infty$. As they observe, the answer is affirmative when $n=o(g)$, by \cite[Lemma~5.1]{MZ}, and when $g$ is fixed, by \cite[Theorem~6.1]{ManZ}; they also proved $1/2-\pi^2/20\le\liminf D_{g,n}\le\limsup D_{g,n}\le1$ for every cusp sequence~\cite[Proposition~36]{SW}. Theorem~\ref{thm:volume} supplies the remaining range, in which $g\to\infty$ and $n/g\not\to0$. Combined with the fixed-genus case, it gives $D_{g,n}\to1$ whenever $g+n\to\infty$.

\subsection{Adjacent coefficients}
Mirzakhani's second question~\cite[Section~1.4]{Mir} concerns adjacent coefficients when the degree grows with the genus. We write the one-boundary volume polynomial in the normalization of~\cite{Mir},
\[
 V_{g,1}(L)=\sum_{k=0}^{N}\frac{a_{g,k}}{(2k+1)!}L^{2k},
 \qquad N=3g-2.
\]
The normalized intersection bracket $b_{g,k}=[\tau_k]_{g,1}$ is the corresponding coefficient of $V_{g,1}(2L)$, and hence $b_{g,k}=4^k a_{g,k}$ (Section~\ref{exact-amplitudes-and-the-harmonic-model}). For fixed $k$, its adjacent ratio tends to one~\cite[Theorem 3.5]{Mir}; the corresponding ratio $a_{g,k}/a_{g,k+1}$ tends to four. Thus the limit one stated in \cite[Section~1.4]{Mir} refers to the bracket normalization, which is also the one used in~\cite{LX}.

\Needspace{12\baselineskip}
\begin{mainresult}
\begin{theorem}[All adjacent degrees]\label{thm:coefficient}
For $0\le k<N$, let $t_{g,k}\in(0,\pi)$ be determined by
\[
 \frac{2k+1}{2g-1}
 =\frac{t_{g,k}^{2}\sin t_{g,k}}
 {\sin t_{g,k}-t_{g,k}\cos t_{g,k}}.
\]
Then
\[
 \sup_{0\le k<3g-2}
 \left|\frac{a_{g,k}/a_{g,k+1}}{4\pi^2/t_{g,k}^{2}}-1\right|
 =O(g^{-1/8}).
\]
In particular, $a_{g,k}/a_{g,k+1}\to4$ for every $k=o(g)$, and, with $m=N-k\ge1$,
\[
 \frac{a_{g,N-m}}{a_{g,N-m+1}}
 \sim\frac{4\pi^2g}{5m}\qquad\text{if }m=o(g).
\]
\end{theorem}
\end{mainresult}
The density $(2k+1)/(2g-1)$ retains the exact codimension $m$ in the saddle equation, which is why the last formula also covers bounded positive $m$. Section~\ref{uniform-one-boundary-coefficient-ratios} proves the uniform equivalent, and Appendix~\ref{sec:quantitative-coefficients} proves the stated error. On the range $k\ge\varepsilon g$, for each fixed $\varepsilon\in(0,3)$, the stronger error $O_\varepsilon(g^{-1/4})$ holds, again including bounded positive codimension.

For each fixed codimension, Liu and Xu obtained a complete expansion of the coefficients relative to the top one~\cite[Theorem~1.2]{LX}, a partial answer to Mirzakhani's question. The bounded-codimension case of Theorem~\ref{thm:coefficient} agrees with their result after the change of normalization. What is new is uniformity over all degrees, including degrees and codimensions that grow with the genus.

The trigonometric saddle equation is the stationary equation of an instanton prediction from physics; see Subsection~\ref{subsec:instanton-model} for the change of conventions and the distinction from the coefficient theorem.

\textbf{The diagonal chain.}
Along the chain $V_{G,0},V_{G-1,2},\ldots,V_{0,2G}$ the Euler characteristic is constant and genus is traded for pairs of cusps. Mirzakhani noted that the adjacent ratios at the two ends of this chain follow from the fixed-$n$ and fixed-$g$ asymptotics, but that little was known about the whole sequence~\cite[Section~1.4]{Mir}. Theorem~\ref{thm:chain} in Section~\ref{uniform-logarithmic-volume-rate-and-the-full-diagonal-chain} determines every entry of the chain up to a factor $\exp(o(G))$, uniformly in its position and including the planar endpoint, through an explicit rate function obtained by integrating the ratio profile $\mathcal B$. Unlike Theorem~\ref{thm:intro-absolute}, it also covers bounded genus.

\subsection{Strategy of the proof}\label{subsec:strategy}
The heuristic of Subsection~\ref{subsec:parameter} has two parts: the approximation~\eqref{eq:intro-heuristic} of $\Phi_g'$, uniform in the cusp weight, and the passage from $\Phi_g'$ to one of its coefficients. Three ideas carry the proof, each removing one obstruction to uniformity in $n$: cancellation in the volume formula, the unknown constant of integration, and the passage from generating functions to coefficients.

The first obstruction is present from the outset. Mirzakhani's volumes are intersection numbers of $\kappa_1$ and $\psi$ classes, and $\kappa_1$ can be expressed through pure $\psi$ intersections~\cite{KMZ,ManZ}. The resulting formula~\eqref{eq:H5} is an alternating sum over the number of additional $\psi$ insertions, and asymptotics for its individual terms do not control the cancellation.

\begin{enumerate}[label=(\roman*),leftmargin=*]
\item \textbf{Positivity from the string flow.}
Moving along the string flow by a time $Y$ replaces the alternating weights $-X^{(k)}(0)$ by the derivatives $-X^{(k)}(Y)$; see~\eqref{eq:H6}. At $Y=0$ the resulting \emph{shifted amplitudes} are the Weil--Petersson brackets; after a natural rescaling, at $Y=y_0$ they become pure $\psi$ intersection numbers with insertions of $\tau_2$. The Dijkgraaf--Verlinde--Verlinde equations give the normalized recursion~\eqref{eq:H9}, whose kernel coefficients are the Taylor coefficients of $1/\bar K_Y$, where $\bar K_Y$ is the rescaled kernel in~\eqref{eq:P9}. These coefficients are nonnegative throughout $[0,y_0]$, including for the limiting kernel $\bar K_{y_0}(w)=1-w$ (Section~\ref{the-positive-compact-family-of-kernels}). The alternating sum is thereby replaced by a positive recursion whose kernels form a compact family. Comparison with an explicit solvable model, the \emph{harmonic model}~\eqref{eq:H12}, gives two-sided bounds that are uniform in $Y$ and in all degrees but lose a factor $\mathrm e^{O(\sqrt g)}$ (Proposition~\ref{prop:harmonic}). This step uses no information about volumes.

\item \textbf{Ratios, and normalization at the pure-$\psi$ endpoint.}
For a fixed number of distinguished markings, the dilaton equation expresses the $Y$-derivative of $\log\Phi_g^{(p)}(X(Y))$ through amplitudes with one more marking, and the KdV equation relates genus $g$ to genus $g-1$ with two more markings. By positivity the amplitudes decrease as a degree increases, and this monotonicity upgrades the coarse bounds of step~(i) to relative errors $O(1/g)$ for both ratios (Proposition~\ref{prop:marked-ratios}). Integrating the resulting differential equation in $Y$ gives~\eqref{eq:intro-heuristic} up to a factor independent of $Y$. This factor is determined at $Y=y_0$, where only $\langle\tau_2^{3g-3}\rangle_g$ enters and Painlev\'e~I supplies the asymptotics (Section~\ref{sec:prefactor}). The deformation thus carries the normalization problem to a point where it is classical. Iterating the same equations gives all higher corrections, with remainders controlled as the order grows (Section~\ref{sec:higher-order}, with proofs in Appendix~\ref{app:marked}).

\item \textbf{Extraction through a positive factor.}
An expansion of a generating function does not by itself approximate its coefficients. With $h=2g-2+p$, we show that $\Phi_g^{(p)}$ is $(1/X'(Y(x)))^{h}$ times a power series with nonnegative coefficients (Lemma~\ref{lem:convolution}):
\[
 \Phi_g^{(p)}=V^hR_{g,p},\quad V=\frac1{X'\circ Y},\quad [x^j]R_{g,p}\ge0,
 \qquad\text{hence}\qquad
 \bigl|\Phi_g^{(p)}(x\mathrm e^{i\theta})\bigr|
 \le\Phi_g^{(p)}(x)\left|\frac{V(x\mathrm e^{i\theta})}{V(x)}\right|^h .
\]
The decay on the circle $|z|=x$ is thus controlled by one explicit function raised to the power $h$, whatever the unknown factor $R_{g,p}$ is. This gives Fourier bounds uniform in the cusp weight, and hence a Gaussian local limit theorem when both $h$ and the model mean $h\Psi(x)$ tend to infinity, even as $x$ approaches $0$ or $x_c$ (Proposition~\ref{prop:local-limit}). Bounded cusp counts instead require a Poisson approximation near $x=0$: a change of variables of the type studied by Temme~\cite{Temme} gives an expansion around the Poisson law that remains regular in this range (Section~\ref{sec:sharp-extraction}, with proofs in Appendix~\ref{app:extraction}). Together these estimates yield Theorem~\ref{thm:intro-absolute} (Section~\ref{sec:absolute}).
\end{enumerate}

The same positive factor reappears in the applications: it gives the uniform coefficient-ratio bounds that control curves cutting off many cusps (Appendix~\ref{app:short-curves}) and the boundary factorization of Proposition~\ref{prop:bounded-boundaries}.

\textbf{From volumes to random surfaces.}
The applications use the volume theorem in one way: through
Mirzakhani's integration formula, ratios of volumes are expected
numbers of geometric configurations, so Theorem~\ref{thm:volume}
determines the statistics of short separating curves at every cusp
density (Section~\ref{sec:geometry}). The spectral applications then
require a second passage, from these counts to actual eigenvalues. A
short curve produces a trial function, hence an eigenvalue upper bound,
but neither a short-curve law nor an upper bound alone determines the
low spectrum; one must show that nothing else is there. Three
arguments do this in the three regimes. At the critical scale,
funnel trial states are matched to the expected spectral count
supplied by Huang's positive trace~\cite{HuangGap}; on shrinking
scales, a comparison on the whole low spectral subspace excludes the
modes missed by trial functions; in fixed windows, concentration of
orbital statistics and trace approximation identify the limiting
measure. The geometric application retains based holonomies in the same
integration formula, so that the carrier integrals determine whole
pointed balls rather than curve counts.

\subsection{The spectral transition at the square-root scale}
\label{subsec:intro-geometry}
Let $n\ge1$, let $X$ have normalized Weil--Petersson law
$\mathbb P^{\mathrm{WP}}_{g,n}$, and let $\Delta_X$ be its nonnegative
Laplace--Beltrami operator. The spectrum below $1/4$ is discrete; we
count its eigenvalues $\lambda_j(X)$ with multiplicity and write
$G(X)=\inf\operatorname{spec}(\Delta_X|\mathbf1^\perp)$, with
$G(X)=1/4$ if there is no positive eigenvalue below $1/4$.
For closed surfaces ($n=0$), Mirzakhani proved a positive gap with high
probability~\cite[Section~4.5]{Mir}; the bound was raised to $3/16-\varepsilon$ by Wu
and Xue and by Lipnowski and Wright~\cite{WX,LipWright}, and to the
optimal $1/4-\varepsilon$ by Anantharaman and Monk~\cite[Theorem~1.6]{AMGap}. With
cusps, Hide, Shen and Wu, and He, Wu and Xue obtained positive lower
bounds below the square-root scale~\cite[Theorem~1.3]{Hide}, \cite[Theorem~3]{SW},
\cite[Theorem~1.2]{HWX}, and Huang proved $G(X)\to1/4$ throughout this
range~\cite[Theorem~1.1]{HuangGap}. The literature is compared in
detail at the beginning of Part~\ref{part:spectra}.

\textbf{Why $\sqrt g$.} A short geodesic cutting off two cusps bounds a
pair of pants of area $2\pi$, attached to the rest of the surface
through a collar of small conductance, and therefore produces a small
eigenvalue. By Mirzakhani's integration formula, the expected number of
geodesics of length at most $L$ that cut off a given pair of cusps is
$(L^2/2)V_{g,n-1}/V_{g,n}$ up to a factor $1+O(L)$. When $n=o(g)$,
Theorem~\ref{thm:volume} makes this $L^2/(16\pi^2g)$, so the expected
total over the $n^2/2$ pairs is about $n^2L^2/(32\pi^2g)$. Hence such
curves of bounded length are absent when $n\ll\sqrt g$, occur in
bounded number when $n\asymp\sqrt g$, and have lengths of order
$\sqrt g/n\to0$ when $n\gg\sqrt g$; these are the three spectral phases.
At positive density the same count requires the full profile
$\mathcal C(n/g)$ and all planar cusp types, and this is where the
uniform theory is indispensable. Table~\ref{tab:intro-phases}
summarizes the resulting picture.

At the critical scale the model spectral object is the complete surface
obtained by attaching a funnel to the pants. Let $\beta(s)\in(0,1/4)$
be the unique eigenvalue below $1/4$ of the complete hyperbolic surface
with two cusps and one funnel whose convex core has boundary length
$s$.

\begin{table}[ht]
\caption{The spectral phases of Weil--Petersson random surfaces as
$g\to\infty$.}
\label{tab:intro-phases}
{\small
\renewcommand{\arraystretch}{1.3}
\setlength{\tabcolsep}{4pt}
\begin{tabularx}{\textwidth}{@{}lll>{\raggedright\arraybackslash}Xll@{}}
\toprule
Cusp count & $Y$ & Gap $G(X)$ & Lowest eigenvalues & Fixed windows & Results\\
\midrule
$n=o(\sqrt g)$ & $\to0$ & $\to1/4$ & $\frac g{n^2}(\frac14-\lambda_j)\Rightarrow\operatorname{PPP}(\kappa_0x^{-2}dx)$\,$^{a}$ & empty & \cite{HuangGap}, \ref{thm:moving-edge}\\
$n/\sqrt g\to c$ & $\to0$ & $\Rightarrow G_c$ & $\{\beta(s):s\in\mathcal Q_c\}$ & same limit & \ref{thm:critical-eigenvalues}, \ref{cor:gap-phase-diagram}\\
$\sqrt g\ll n\ll g$ & $\to0$ & $\to0$ & $A_g\lambda_j\Rightarrow\operatorname{PPP}(2t\,dt)$ & $\frac g{n^2}\Xi_X\to\nu_*$ & \ref{thm:supercritical-spectrum}, \ref{thm:fixed-energy-spectrum}\\
$n/g\to\lambda$ & $\to Y_\lambda$ & $\to0$ & $A_g\lambda_j\Rightarrow\operatorname{PPP}(2t\,dt)$ & $\frac1g\Xi_X\to\nu_\lambda$ & \ref{thm:supercritical-spectrum}, \ref{thm:positive-density-spectrum}\\
$g\ll n\le g^{11/9-\delta}$ & $\to y_0$ & $\to0$ & $A_g\lambda_j\Rightarrow\operatorname{PPP}(2t\,dt)$\,$^{b}$ & $\frac1n\Xi_X\to\nu_\infty$ & \ref{thm:dense-persistence}, \ref{thm:cusp-dominated-spectrum}\\
$n\gg g^3/(\log g)^4$ & $\to y_0$ & $\to0$ & $A_g\lambda_r\to0$ for each fixed $r$ & $\frac1n\Xi_X\to\nu_\infty$ & \ref{thm:dense-collapse}, \ref{thm:cusp-dominated-spectrum}\\
\bottomrule
\end{tabularx}}

\smallskip
\noindent\parbox{\linewidth}{\footnotesize\raggedright
Here $Y\in[0,y_0]$ is the density parameter, with $\lambda(Y_\lambda)=\lambda$;
$\Xi_X$ is the counting measure of the spectrum in $(0,1/4)$;
$G_c\in(0,1/4)$ is the random gap limit; $\kappa_0=1/(8\pi^2)$; and
$A_g\sim\pi n/(2\sqrt{2g})$ when
$\sqrt g\ll n\ll g$. The last column gives the numbers of the theorems
and corollaries; the gap limit for $n/\sqrt g\to\infty$ is
Corollary~\ref{cor:spectral}, and the empty fixed windows for
$n=o(\sqrt g)$ hold with probability tending to one.\\
$^{a}$\,If moreover $n\log g/\sqrt g\to\infty$.
$^{b}$\,Proved under $n^{9/4}\log^{10}(n+2)/g^{11/4}\to0$, which holds
when $n\le g^{11/9-\delta}$ for a fixed $\delta>0$. Between this range
and the collapse range the microscopic law is open.}
\end{table}

\begin{mainresult}
\begin{restated}{Theorem~\ref{thm:critical-eigenvalues} and Corollary~\ref{cor:gap-phase-diagram}}[The critical spectrum]
Suppose $g\to\infty$ and $n/\sqrt g\to c\in(0,\infty)$, and let
$\mathcal Q_c$ be a Poisson process on $(0,\infty)$ with intensity
$c^2\sinh(s/2)\,ds/(8\pi^2)$. Then
\[
 \sum_{0<\lambda_j(X)<1/4}\delta_{\lambda_j(X)}
 \ \Longrightarrow\ \sum_{s\in\mathcal Q_c}\delta_{\beta(s)}
\]
in distribution, for the vague topology on $[0,1/4)$, and the mean
measures converge. The cuspidal part of the left side tends to zero, so
the limiting spectrum is residual. In particular $G(X)$ converges in
distribution to the smallest point $G_c\in(0,1/4)$ of the limit.
Consequently, for every sequence $n=n(g)\ge1$,
\[
 G(X)\xrightarrow{\ \mathbb P\ }\frac14
 \iff\frac n{\sqrt g}\to0,\qquad
 G(X)\xrightarrow{\ \mathbb P\ }0
 \iff\frac n{\sqrt g}\to\infty,
\]
and $G(X)$ has no deterministic limit in $(0,1/4)$ along any sequence.
\end{restated}
\end{mainresult}

The limit is a random point configuration rather than a deterministic
value for each eigenvalue: one model eigenvalue for each short two-cusp
pair of pants, and nothing else. The volume input determines the pants
process; the substantive spectral step is to show that these localized
states \emph{exhaust} the spectrum below $1/4$. This is done by matching
the count of funnel trial states with the expected spectral count
supplied by Huang's positive trace~\cite{HuangGap}. Of the two
equivalences, the first is Huang's theorem~\cite[Theorem~1.1 and
Corollary~6.15]{HuangGap}; the second uses the adjacent-volume limit of
Theorem~\ref{thm:volume} and completes the negative answer to
Shen--Wu's Question~2~\cite{SW}, which asks whether a uniform spectral
gap persists when $n/\sqrt g\to\infty$.

Away from the critical scale the lowest eigenvalues obey different
laws. Slightly below criticality, if $c=n/\sqrt g\to0$ and
$c\log g\to\infty$, the defects $(1/4-\lambda_j(X))/c^2$ converge to a
Poisson process with intensity $dx/(8\pi^2x^2)$; in particular
$(1/4-G(X))/c^2$ has the Fr\'echet law $\exp(-1/(8\pi^2x))$
(Theorem~\ref{thm:moving-edge}). Above criticality the collar waists
shrink, and the eigenvalues are governed by the conductances of short
collars weighted by the areas they enclose.

\begin{mainresult}
\begin{restated}{Theorems~\ref{thm:supercritical-spectrum}, \ref{thm:dense-collapse} and~\ref{thm:dense-persistence}}[Supercritical spectra]
Suppose $g\to\infty$ and $n/\sqrt g\to\infty$, and let $A_g$ be the
explicit scale \eqref{eq:super-scale}. If either $n/g$ stays bounded,
or $n/g\to\infty$ and $n^{9/4}\log^{10}(n+2)/g^{11/4}\to0$, then
\[
 \sum_{0<\lambda_j(X)<1/4}\delta_{A_g\lambda_j(X)}
 \ \Longrightarrow\ \operatorname{PPP}(2t\,dt)
 \quad\text{on }[0,\infty).
\]
Thus $\mathbb P\{A_gG(X)>t\}\to\mathrm e^{-t^2}$, and $G(X)/\mathfrak h(X)\to1/\pi$
in probability, where $\mathfrak h$ is the Cheeger constant. When
$\sqrt g\ll n\ll g$, one has $A_g\sim\pi n/(2\sqrt{2g})$. If instead
$n\gg g^3/(\log g)^4$, then $A_g\lambda_r^+(X)\to0$ in probability for
every fixed $r$, where $\lambda_r^+$ is the $r$th positive eigenvalue.
\end{restated}
\end{mainresult}

In \eqref{eq:super-scale}, $A_g=2\pi^2(q_gM_2(x_g))^{1/2}$, where
$q_gw_k(x_g)$ is the intensity of short curves cutting off $k$ cusps
and $M_2$ weights type $k$ by $(k-1)^2$, which up to a constant is the
square of the area $2\pi(k-1)$ of the cap cut off by such a curve. The
proof is a two-sided comparison on the whole low spectral subspace:
trial functions give upper bounds, and the comparison excludes the low
modes they miss, including modes supported on a large component
separated by several curves. The collapse at large density has a
geometric cause: large planar caps then support low-energy functions
of their own, and the collar picture fails. The correct microscopic
normalization between $g^{11/9}$ and $g^3/(\log g)^4$ remains open.

Fixed energy windows behave differently. Let $\Xi_X$ be the counting
measure of the spectrum in $(0,1/4)$. Fixed windows contain many
eigenvalues, and the normalized counting measures converge to
deterministic limits: $(g/n^2)\,\Xi_X\to\nu_*$ when $\sqrt g\ll n\ll g$,
$g^{-1}\Xi_X\to\nu_\lambda$ when $n/g\to\lambda\in(0,\infty)$, and
$n^{-1}\Xi_X\to\nu_\infty$ when $g=o(n)$, including every fixed genus
(Theorems~\ref{thm:fixed-energy-spectrum}--\ref{thm:cusp-dominated-spectrum}).
Theorem~\ref{thm:regularized-spectrum} includes the continuous spectrum
through the classical regularized trace, and
Corollary~\ref{cor:cusp-density-determinant} derives the limit of the
log determinant per unit area. These results rest on concentration of
orbital statistics rather than on individual short curves, and the
measures $\nu_\lambda$ are characterized by planar orbital integrals.

\subsection{Local geometry at positive cusp density}
\label{subsec:intro-local}
For every $\lambda\in(0,\infty)$,
Theorem~\ref{thm:wp-local-family} constructs a unimodular random
surface $(S_\lambda,o_\lambda)$ as the fixed-radius local limit when
$n/g\to\lambda$, rooted by area, and proves quenched concentration of
local statistics. Its law tends to the hyperbolic plane as
$\lambda\downarrow0$ and to the area-rooted infinite WP sphere of Budd
and Curien~\cite{BC} as $\lambda\uparrow\infty$. After filling the
cusps, these surfaces are conformally discs at finite density and the
plane at the endpoint; Brownian motion on the original hyperbolic
surfaces is respectively transient and recurrent
(Theorem~\ref{thm:wp-conformal-phase}). By the theorem of Hopf, Tsuji
and Sullivan, the geodesic flow is accordingly not ergodic at finite
density and ergodic at the endpoint
(Corollary~\ref{cor:wp-geodesic-flow}); the endpoint statements
confirm the conjecture of Budd and Curien~\cite[Section~7.3]{BC}. A separate growing-radius tree
coupling gives the endpoint's Brownian-plane scaling limit: local
pointed Gromov--Hausdorff--Prokhorov convergence after cusp trimming,
with the full area projected to the trimmed space, and local pointed
Gromov--Prokhorov convergence for the untrimmed surface with its full
area (Theorem~\ref{thm:wp-brownian-plane}).

\subsection{Questions from the literature}
\label{subsec:intro-questions}
Table~\ref{tab:intro-questions} lists the named questions treated here
and the extent of each answer; other earlier work is compared at the
beginning of each part.

\begin{table}[ht]
\caption{Questions from the literature answered or treated in this paper.}
\label{tab:intro-questions}
{\small
\renewcommand{\arraystretch}{1.3}
\setlength{\tabcolsep}{5pt}
\begin{tabularx}{\textwidth}{@{}l>{\raggedright\arraybackslash}Xl@{}}
\toprule
Source & Question & Answered in\\
\midrule
Mirzakhani~\cite[\S1.4]{Mir} & $V_{g,n}/V_{g-1,n+2}$ when $n=n(g)\to\infty$ & Thm.~\ref{thm:volume}, all $n$\\
Mirzakhani~\cite[\S1.4]{Mir} & adjacent coefficients at growing degree & Thm.~\ref{thm:coefficient}, all degrees\\
Shen--Wu~\cite[Question~3]{SW} & does $V_{g,n}^2/(V_{g,n-1}V_{g,n+1})\to1$? & yes: Thm.~\ref{thm:volume}\\
Shen--Wu~\cite[Question~2]{SW} & uniform spectral gap when $n/\sqrt g\to\infty$? & no: Cor.~\ref{cor:spectral}\,$^{a}$\\
Budd--Lions~\cite[Open Question~1.6]{BL} & tight spectrum at other weights, fixed counts & Thm.~\ref{thm:tight-spectrum}, Cor.~\ref{cor:tight-boltzmann}\\
Hide--Thomas~\cite[\S1.1]{HTLargeN} & linear count of small residual eigenvalues & Thm.~\ref{thm:many-residual}\\
Hide--Thomas~\cite[Remark~1.8]{HTLargeN} & spectral measures with many cusps & Thms.~\ref{thm:positive-density-spectrum}--\ref{thm:regularized-spectrum}\,$^{b}$\\
Budd--Curien~\cite[\S7.3]{BC} & Brownian-plane limit of the infinite sphere & Thm.~\ref{thm:wp-brownian-plane}\\
Budd--Curien~\cite[\S7.3]{BC} & conjecture: the infinite sphere is parabolic, recurrent, with ergodic geodesic flow & yes: Thm.~\ref{thm:wp-conformal-phase}, Cor.~\ref{cor:wp-geodesic-flow}\\
Budd--Curien~\cite[\S7.3]{BC} & hyperbolic analogue of the infinite sphere & Thm.~\ref{thm:wp-local-family}\,$^{c}$\\
\bottomrule
\end{tabularx}}

\smallskip
\noindent\parbox{\linewidth}{\footnotesize\raggedright
$^{a}$\,Shen--Wu treated $\sqrt g\ll n\ll g$, and Zograf's bound covers
$n/g\to\infty$; Corollary~\ref{cor:spectral} supplies the remaining range.
$^{b}$\,Low-energy limits and the regularized full-energy measure; cusp
forms and scattering are not separated.
$^{c}$\,Constructed at every finite density; its near-critical scaling
limit remains open.}
\end{table}

The near-critical scaling question of Budd--Curien and the dense
microscopic spectrum between $g^{11/9}$ and $g^3/(\log g)^4$ remain open; see
Section~\ref{sec:conclusions}.

\subsection{Organization of the paper}
Part~\ref{part:volumes} follows the three steps of
Subsection~\ref{subsec:strategy}: positive recursion, endpoint
normalization and coefficient extraction, followed by the volume-ratio
and adjacent-coefficient consequences.
Part~\ref{part:spectra} passes from short geometry to the critical
spectrum, its moving edge and the supercritical microscopic laws, and
then to fixed-energy windows; its opening guide compares the
observables, density ranges and earlier work. Part~\ref{part:geometry}
treats local limits, conformal type and the Brownian plane.

Main-text sections contain the statements and principal arguments;
long proofs are in the appendices, grouped by part and ordered by first
use, and Appendix~\ref{app:reader-guide} gives a proof-location table
and a glossary. A reader interested first in volume asymptotics may
read Part~\ref{part:volumes} with Appendices~\ref{app:kernels}
and~\ref{app:primitive-bounds}; the other appendices to
Part~\ref{part:volumes} give the inverse-genus rates, the higher-order
expansions and the all-marking envelope \eqref{eq:gevrey-envelope}
used in Part~\ref{part:spectra}. A reader interested in spectra may
start with Part~\ref{part:spectra}, taking Theorem~\ref{thm:volume},
Proposition~\ref{prop:bounded-boundaries} and the positive
convolution factor of Section~\ref{a-positive-convolution-factor} as
inputs, which its opening page lists precisely.

\clearpage
\part{Uniform volume asymptotics}\label{part:volumes}
\phantomsection\label{guide:volume-context}
We compare the uniform theorem with the existing volume asymptotics.
The proof begins in Section~\ref{the-positive-compact-family-of-kernels}.

\subsection*{Moduli spaces, volumes and curve counting}
Weil--Petersson volumes play two roles in this paper: their boundary-length polynomials encode intersection numbers, and
ratios of these volumes enter expectations of curve counts on random
hyperbolic surfaces.
At fixed genus, the large-cusp prediction of Itzykson is recorded
in~\cite[Remark~3]{Zog98}. Zograf's genus-zero recursion~\cite{Zog93}
was followed by the proof of the genus-zero large-cusp asymptotics by
Kaufmann--Manin--Zagier~\cite{KMZ}, and
Manin--Zograf~\cite[Theorem~6.1]{ManZ} treated every fixed genus.
Wolpert computed the volume of $\mathcal M_{1,1}$~\cite{WolpertTorus},
and Penner developed a combinatorial method for computing these
volumes~\cite{Penner}.
Building on McShane's identity~\cite{McShane}, Mirzakhani derived
a recursion that computes $V_{g,n}(\mathbf L)$ as a polynomial in
the squared boundary lengths~\cite{MirInv}; its coefficients are
intersection numbers of tautological classes on
$\overline{\mathcal M}_{g,n}$~\cite{MirJAMS}, with the symplectic
normalization established by Wolpert~\cite{Wolpert}. Mulase and Safnuk
recast the recursion as Virasoro constraints and showed that the
generating function of mixed $\kappa_1$--$\psi$ intersections is a
one-parameter solution of the Korteweg--de Vries (KdV) hierarchy~\cite{MS};
Liu and Xu showed directly that Mirzakhani's recursion is equivalent to
the Witten--Kontsevich theorem~\cite{LXEquiv}, and Liu and
Xu~\cite{LXHigher} and Eynard~\cite{EynardMumford} extended the
recursion to volumes weighted by higher Mumford classes;
Eynard and Orantin identified it with the loop equations of a matrix
model~\cite{EO}. The identities expressing mixed $\kappa_1$--$\psi$
intersections through pure $\psi$ intersections, which are the starting
point of Section~\ref{exact-amplitudes-and-the-harmonic-model}, are due
to Arbarello and Cornalba~\cite{AC}, Kaufmann, Manin and Zagier~\cite{KMZ}
and Manin and Zograf~\cite{ManZ}.
Related exact formulas include the tight-volume recursions of Budd and
Zonneveld~\cite{BZ}, the KdV boundary-addition identities and
fixed-genus Bessel singularity studied by Johnson~\cite{Johnson}, and
the Virasoro-theoretic proof by Xiong and Zhou of the intersection-number
and cohomological-field-theory descriptions of topological recursion,
including variation of the spectral-curve
function~\cite[Theorem~1.1 and Section~5.7]{XZ}. These are exact
identities; obtaining uniform remainder estimates from them is a
separate analytic problem.

Mirzakhani's integration formula~\cite{MirInv}
expresses the average number of simple multicurves of a prescribed
type in terms of volumes of the surfaces obtained by cutting along them.
After division by $V_{g,n}$, these formulas become expectations under
Weil--Petersson probability, so relative volume estimates control
geometric probabilities. Mirzakhani used this connection to study random
surfaces at large genus~\cite{Mir}, and Mirzakhani and Petri used joint
curve counts to obtain Poisson length statistics for closed
surfaces~\cite{MP}.
A different curve-counting limit fixes the surface and lets the length
bound tend to infinity. Mirzakhani's counting theorem assigns frequencies
to simple-curve topological types in that limit~\cite{MirCount}.
Delecroix, Goujard, Zograf and Zorich related these frequencies to
intersection numbers and Masur--Veech volumes~\cite{DGZZ} and studied
random multicurves at large genus~\cite{DGZZRandom}; Ren
obtained the large-genus ratio of separating to nonseparating frequencies
with many cusps~\cite{Ren}. Those frequencies are topological quantities,
whereas the processes studied here count short curves on a random
surface as its topology varies.

Masur--Veech volumes are the flat counterpart of Weil--Petersson
volumes, and their large-genus theory has developed in parallel.
Eskin and Zorich conjectured the large-genus asymptotics of
Masur--Veech volumes of strata of Abelian differentials~\cite{EZ};
Aggarwal proved them for all strata~\cite{AggAbelian}, and Chen,
M\"oller, Sauvaget and Zagier expressed these volumes through
intersection numbers on the strata and evaluated their large-genus
limits~\cite{CMSZ}; Chen, M\"oller and Sauvaget expressed the volumes
of the principal strata of quadratic differentials through
intersection numbers~\cite{CMSQuad}.
Rafi and Zorich recently determined the large-$n$
asymptotics of Masur--Veech volumes of strata of quadratic differentials
with many simple poles at fixed genus~\cite{RZ}.

\subsection*{Asymptotics of volumes and intersection numbers}
The fixed-genus regime was recalled at the beginning of the
introduction. Here we discuss large-genus estimates at fixed cusp
count, coefficients of growing degree, pure $\psi$ intersection numbers
and predictions from physics.

\textbf{Fixed cusp counts.}
Grushevsky~\cite{Gru} and, for closed surfaces, Schumacher and
Trapani~\cite{ST} showed that for fixed $n$ the volumes grow like
$(2g)!$ up to factors $c^g$; see the account in \cite[Section~1.4]{Mir}.
After Mirzakhani's estimates~\cite{Mir}, Mirzakhani and Zograf obtained
the expansion to every fixed order~\cite[Theorem~1.2]{MZ}. They showed that its coefficients
are polynomials in $n$ with coefficients in $\mathbb Q[\pi^{-2},\pi^2]$
and conjectured, on numerical evidence, that only $\mathbb Q[\pi^{-2}]$
occurs~\cite[Theorem~1.2 and Remark~1.3]{MZ}; to our knowledge this
conjecture remains open. Their leading term also
holds when $n^2/g\to0$~\cite[Theorem~1.8]{MZ}, and their cusp-ratio
estimate has error $O(n/g)$ for $n\ge1$~\cite[Lemma~5.1]{MZ}. Their
first unmarked correction enters \eqref{eq:AV5}. He, Wu and Xue obtained
first-order volume-ratio expansions with error $O((1+n^2)/g^2)$ for
$n=o(\sqrt g)$~\cite[Theorem~4.11]{HWX}. Anantharaman and Monk expanded
the volume polynomials at large genus~\cite{AM}. Guillermit sharpened
the coefficient bounds and treated boundary lengths of order $\sqrt g$
for fixed $n$~\cite[Theorems~2, 4, 5 and Remark~3]{Guillermit}. For
fixed $n$, her Theorem~2 gives coefficient bounds of order $(C_nj)^j$,
and her Theorem~4 bounds the remainders of the normalized volumes
with boundary themselves, with the same single-factorial dependence on
the order. Our estimates retain
this order dependence uniformly in the density parameter $Y$, including
its critical endpoint, and after extraction uniformly over all cusp counts. Hide, Macera and Thomas gave bounds explicit in
the expansion order, cusp count and degrees~\cite[Theorem~4.2,
statements A1($k$)--A4($k$)]{HMT}. Their genus must exceed an order-dependent multiple of the
cusp count, and the errors grow polynomially in the cusp count or
degrees, so these bounds do not yield the relative remainder of
Theorem~\ref{thm:intro-absolute} at positive cusp density.

\textbf{Sublinear cusp counts.}
Gatignol and Thomas~\cite[Theorem~1.5]{GT} obtain first-order
expansions for normalized intersection brackets and for the cusp and
genus volume ratios throughout $n=o(g)$, with ratio errors
$O(n^2/g^2+1/g)$. In particular, their genus-ratio expansion supplies
the two-sided limit needed at $n\asymp\sqrt g$. These estimates overlap
with the zero-density consequences of Theorem~\ref{thm:volume}; the
present theorem also treats positive and unbounded $n/g$, and the
absolute expansion has remainders uniform over all $n$.

\textbf{Growing degrees.}
Wu and Xue established effective upper bounds for normalized
intersection numbers when some degrees are large compared with
$\sqrt{2g-2+n}$, together with bounds for volumes with large boundary
lengths~\cite{WXPrime}. These bounds control individual coefficients
and volumes, whereas the adjacent-ratio asymptotics of
Theorem~\ref{thm:coefficient} require relative control across
neighboring degrees. Moreover, a bracket may itself be small, and then
an additive error polynomial in the degree gives no relative control. Liu and Xu recently proved general balancing inequalities for mixed tautological intersection numbers~\cite{LXBalancing}; the monotonicity statements used here (Lemma~\ref{lem:monotone} and the bracket monotonicity in Section~\ref{uniform-one-boundary-coefficient-ratios}) are of a simpler kind. We use the Mirzakhani--Zograf bracket-difference
estimates near the lowest degrees; the full degree range requires the
harmonic comparison and the log-concavity argument of
Section~\ref{uniform-one-boundary-coefficient-ratios}.

\textbf{Pure $\psi$ intersection numbers.}
Motivated by Masur--Veech volumes of quadratic differentials,
Delecroix, Goujard, Zograf and Zorich conjectured uniform large-genus asymptotics for pure $\psi$ intersections (Appendix~E of the arXiv version of~\cite{DGZZ}). Aggarwal proved them for $o(\sqrt g)$ insertions~\cite{Agg}, and Guo and Yang gave a different proof for at most $C\log g$ insertions, for every fixed $C$, based on an explicit formula for the $n$-point generating series~\cite[Theorem~1]{GY}. For a fixed number of insertions, Eynard, Garcia-Failde, Giacchetto, Gregori and Lewa\'nski obtained resurgent expansions to every fixed order, uniformly in the degrees~\cite[Theorem~3.13]{EGGGL}.

Our comparison needs only fixed positive upper and lower bounds for the normalized pure-$\psi$ intersections. Appendix~\ref{app:primitive-bounds} proves these bounds by a scalar induction, adapting the recursion estimates of Guo, Yang and Zagier~\cite{GYZ} and using the two-point bound of Delecroix, Goujard, Zograf and Zorich~\cite{DGZZ}, as quoted in~\cite[equation~(6)]{DGZZBounds}. The absolute constant is determined by the single sequence $\langle\tau_2^{3g-3}\rangle_g$ through the Itzykson--Zuber identity and Painlev\'e~I. To match the expansions at the pure-$\psi$ endpoint when the truncation order grows with $g$, we use the exact identity~\cite[equation~(6.10)]{IZ}, the zero-instanton contour formula of Garoufalidis, Its, Kapaev and Mari\~no~\cite{GIKM}, and the Borel representation of Costin, Costin and Huang~\cite{CCH}. Alameddine, Marchal, Nikolaev and Orantin proved the resurgence of the tritronqu\'ee solutions of the deformed Painlev\'e~I equation~\cite{AMNOResurgence} and recently solved the direct monodromy problem for the Borel-summed zero-parameter solution, identifying it with the tritronqu\'ee solutions~\cite{AMNO}; neither result is used in our proof, which needs only the large-order asymptotics of the formal solution. Borot, Eynard and Giacchetto use comparison with Painlev\'e~I amplitudes to bound topological-recursion correlators, uniformly in genus and number of markings~\cite[Main Theorem and Section~1.3]{BEG}. Their bounds concern the growth of the correlators themselves, whereas the factorial bounds in Theorem~\ref{thm:intro-absolute} concern the order of the inverse-genus expansion after division by an explicit leading term.

\phantomsection\label{subsec:intro-instanton}
\textbf{Instanton predictions.}
The factorial growth of genus expansions in two-dimensional gravity and
its nonperturbative origin go back to Shenker~\cite{Shenker}; for
minimal strings approaching Jackiw--Teitelboim gravity, the resurgent
structure was analyzed by Gregori and Schiappa~\cite{GregoriSchiappa}.
Eynard, Garcia-Failde, Gregori, Lewa\'nski and Schiappa derived large-genus asymptotics for ordinary Weil--Petersson volumes under a resurgence assumption~\cite{EGGLS}. Griguolo, Papalini, Russo and Seminara predicted the coefficients of ordinary Weil--Petersson volumes from instanton contributions~\cite[Section~3.3]{GPRS}. For super Weil--Petersson volumes, X.~Huang proved their two conjectural formulas up to a universal multiplicative constant, for fixed $n$ with $L_i=o(\sqrt g)$ and for descendant degrees $|\mathbf d|=o(g)$~\cite[Theorem~1.7]{XHuangSuper}; Huang, Liu and Xu subsequently determined that constant~\cite{HLXSuper}. Johnson and Rodrigues wrote down the corresponding one-boundary contour integral~\cite[equation~(4.22)]{JR}. Hatsuda, Matsumoto and Okuyama constructed transseries solutions of the KdV equation that governs the genus expansion of the JT partition function, whose leading nonperturbative sector reproduces the instanton terms~\cite{HMO}. After a change of conventions, expansion of the Johnson--Rodrigues integral gives our harmonic coefficient model at $Y=0$, and its stationary equation is the saddle equation of Theorem~\ref{thm:coefficient} (Subsection~\ref{subsec:instanton-model}). Unlike the leading factor $\sinh(L/2)/(L/2)$ at bounded boundary length, this integral retains the dependence of the coefficients on their degree up to the top degree. Theorem~\ref{thm:coefficient} establishes the resulting asymptotic formula for the true adjacent coefficient ratios, uniformly over all degrees, from exact intersection identities and estimates; no resurgence assumption on the Weil--Petersson genus expansion is used.

\subsection*{Zograf's constant and the work of Liu and Xu}
Our determination of Zograf's constant $C=1/\sqrt\pi$ appeared in the first version of this paper~\cite[Corollary~1.2 and Section~7]{LWZenodo}. Liu and Xu subsequently gave another proof using a different deformation~\cite{LXConstant}. Our first version also contained the uniform volume-ratio limits, the diagonal chain, the adjacent-coefficient asymptotics over all degrees, the absolute expansion uniform in the cusp count (then with Gevrey-2 remainder bounds), the vanishing spectral gap with many small eigenvalues, and the short-curve processes; the present version sharpens the remainders to single-factorial order and adds the remaining spectral theory of Part~\ref{part:spectra} and the geometric theory of Part~\ref{part:geometry}.

The deformation of Liu and Xu joins the cubic curve $y=z-4z^3$ to the sine curve $y=\sin(2\pi z)/(2\pi)$ while keeping the first positive zero fixed. For each fixed number $n\ge2$ of markings, their uniform-variation theorem~\cite[Theorem~1.1]{LXConstant} controls the logarithmic derivative of the normalized volumes uniformly in the curve parameter $t\in[0,1]$ by $O_n(g^{-1})$. Transport from the cubic endpoint determines the constant~\cite[Corollary~1.2]{LXConstant}, and their first-correction theorem~\cite[Theorem~1.3]{LXConstant} adds the $1/g$ term with an $O_n(g^{-2})$ remainder, uniformly in $t$. Huang, Liu and Xu subsequently applied the same transport strategy to super Weil--Petersson volumes, where the simpler endpoint is governed by the Br\'ezin--Gross--Witten numbers~\cite{HLXSuper}.

Our deformation is instead the string flow associated with the cusp weight $x=X(Y)$. Its amplitudes with a fixed number $p$ of distinguished markings encode the derivatives $\Phi_g^{(p)}(x)$ of the cusp generating function. Both arguments use positivity, descendant recursions, the dilaton and KdV equations, and a Painlev\'e~I normalization, but along different deformation families. In the present paper, Section~\ref{sec:prefactor} already determines the constant by specializing to $p=1$ and $Y=0$, hence to $V_{g,1}$; uniformity in the total cusp count is not needed for this identification. That uniformity is needed for the stronger conclusion, an expansion valid for every $n\ge0$, including $n$ growing with $g$. This requires the coefficient extraction in Sections~\ref{extracting-fixed-cusp-counts}--\ref{sec:sharp-extraction} and Appendix~\ref{app:extraction} of the present paper, beyond the estimates for a fixed number of distinguished markings.

\section{Positive kernels}\label{the-positive-compact-family-of-kernels}

This section supplies the analytic input for the first of the three
ideas of Subsection~\ref{subsec:strategy}, positivity from the string
flow. Two entire functions of a variable $z$, depending on the density
parameter $Y$ of Subsection~\ref{subsec:parameter}, control everything
that follows. The first, $K_Y$, arises in
Section~\ref{exact-amplitudes-and-the-harmonic-model} when the marking
singled out by the Dijkgraaf--Verlinde--Verlinde equation is merged with
the $\kappa_1$ insertions: these mergers act on the amplitudes as
$K_Y(E)$, where $E$ raises a degree by one, so the Taylor coefficients of
$1/K_Y$ become the coefficients of the recursion. The second, $H_Y$,
plays the same role for the dilaton equation and for the explicit
comparison model. We prove three properties: all zeros
of $K_Y$ are positive and simple (Lemma~\ref{lem:kernel-zeros});
$1/K_Y$ and $1/H_Y$ have nonnegative Taylor coefficients, even as
functions of $Y$ (Lemma~\ref{lem:kernel-positivity}); and after
rescaling by the first zero of $K_Y$ both kernels form a compact family,
including the endpoint $Y=y_0$, which we call the critical endpoint or,
since the amplitudes there reduce to pure $\psi$ intersection numbers,
the pure-$\psi$ endpoint
(Subsection~\ref{compact-normalization-and-uniform-lattice-estimates}).
The first property makes the recursion of
Section~\ref{exact-amplitudes-and-the-harmonic-model} positive, the
second supplies the convolution factor of
Section~\ref{a-positive-convolution-factor}, and the third makes all
later constants uniform in $Y$. The first two are special-function
statements, proved in Appendix~\ref{app:kernels}; the third is proved
here, together with the local limit theorem \eqref{eq:P13} for the
comparison model.

Recall from Subsection~\ref{subsec:parameter} that \[
 X(y)=\sum_{j\ge0}\frac{(-1)^j y^{j+1}}{j!(j+1)!}
      =\sqrt yJ_1(2\sqrt y),
\] and put \(a(y)=X'(y)=J_0(2\sqrt y)\); recall also that \(y_0=j_{0,1}^2/4\) is the first positive zero of \(a\), and that
\(x_c=X(y_0)\). Thus \(a>0\) on \([0,y_0)\), and the inverse branch
\(Y(x)\), with \(Y(0)=0\), maps \([0,x_c)\) onto \([0,y_0)\).

We shall use \begin{equation}
\begin{split}
 K_Y(z)&=\int_0^1X'\bigl(Y+\tfrac z2(1-s^2)\bigr)\,\mathrm{d}s
       =\sum_{\ell\ge0}\frac{X^{(\ell+1)}(Y)}{(2\ell+1)!!}z^\ell,\\
 H_Y(z)&=3\int_0^1s^2X'\bigl(Y+\tfrac z2(1-s^2)\bigr)\,\mathrm{d}s
       =\sum_{\ell\ge0}\frac{3X^{(\ell+1)}(Y)}{(2\ell+3)!!}z^\ell.
 \end{split}                                                   \label{eq:P1}
\end{equation} Beta integration gives the series and hence \begin{equation}
K_Y=H_Y+\tfrac23z\partial_zH_Y.                    \label{eq:P2}
\end{equation} We write \(r(Y)\) for the first positive zero of \(K_Y\); it equals \(2\Delta(Y)\).

\subsection{Real zeros and coefficient positivity}\label{real-zeros-of-the-kernel}

\begin{lemma}[Real zeros]\label{lem:kernel-zeros}
On \([0,\pi^2/4]\), \(X(y)>0\) for \(y>0\), \(X''<0\), and \(X'''>0\).
Moreover \(y_0<3/2\); in particular \(X(Y)\le Y<3/2\) for
\(0\le Y\le y_0\). All zeros of \(a(y)\) are positive real; in particular
\(a\) has no zeros in \(|y|<y_0\). For \(0\le Y<y_0\), every zero of
\(K_Y\) is positive and simple, and \begin{equation}
K_Y(z)=a(Y)\prod_{j\ge1}(1-z/r_j(Y)),
 \qquad 0<r_1(Y)<r_2(Y)<\cdots.                     \label{eq:P5}
\end{equation}
\end{lemma}

At \(Y=0\) one has \(K_0(z)=\sin\sqrt{2z}/\sqrt{2z}\), and the lemma is
the classical description of the zeros of the sine function. For
\(Y>0\), integration along the characteristics of a first-order
equation for \(K_Y\) writes \(K_Y\) as \(\sin t\) minus a small positive
correction; concavity then locates the real zeros, and Rouch\'e's
theorem shows that there are no others.

Reality of the zeros does not carry over from \(K_Y\) to \(H_Y\): some
higher zeros of \(H_Y\) can be nonreal. The required substitute is the
following stronger coefficient assertion.

\begin{lemma}[Coefficient positivity]\label{lem:kernel-positivity}
For \(0\le Y<y_0\), \begin{equation}
\log\frac{a(Y)}{H_Y(z)}=\sum_{m\ge1}t_m(Y)z^m,
 \qquad t_m^{(q)}(0)\ge0\quad(m\ge1,q\ge0),\quad t_m(0)>0, \label{eq:P6}
\end{equation} where each \(t_m\) is analytic on \(|Y|<y_0\) and its Taylor
series at zero converges to \(t_m(Y)\) there. Consequently \(H_Y^{-h}\)
has strictly positive Taylor coefficients for every real \(h>0\), and
\(-H_Y'/H_Y\) has nonnegative coefficients. Moreover every Taylor
coefficient in \(Y\) of \[
 \widehat\vartheta_\ell(Y):=a(Y)[z^\ell]K_Y(z)^{-1}
\] is nonnegative, and the same holds for \(a^{-1}\) and for
\(\vartheta_\ell=[z^\ell]K_Y^{-1}\).
\end{lemma}

The coefficients \(t_m\) satisfy a system of ordinary differential
equations in \(Y\) whose right-hand sides are polynomials with
nonnegative coefficients, and they are positive at \(Y=0\), where \(H_0\)
has only positive real zeros; induction on the order of differentiation
gives the lemma. These are statements about Taylor coefficients
in the parameter \(Y\), which is stronger than positivity at each real
\(Y\); they supply the positive convolution factor of
Section~\ref{a-positive-convolution-factor}.

\subsection{Rescaled kernels and a local limit theorem}\label{compact-normalization-and-uniform-lattice-estimates}

Set \(r=r_1(Y)\), \begin{equation}
\bar K_Y(w)=K_Y(rw)/a(Y),\qquad
 \bar H_Y(w)=H_Y(rw)/a(Y),\qquad
 A(Y)=\bar H_Y(1)/3.                               \label{eq:P9}
\end{equation} The first zero \(r\) is simple and varies analytically. Integration
of \eqref{eq:P2} gives \[
 H_Y(z)=\frac3{2z^{3/2}}\int_0^zv^{1/2}K_Y(v)\,\mathrm{d}v,
\] so \(H_Y(z)>0\) on \([0,r]\). Since \(H_Y^{-1}\) has nonnegative
Taylor coefficients, its radius of convergence exceeds \(r\). Indeed, by
Pringsheim's theorem a power series with nonnegative coefficients is
singular at the positive point of its circle of convergence; if the
radius were at most \(r\), that point would be a zero of \(H_Y\) in
\((0,r]\). In particular \(H_Y\) has no complex zero of modulus at most
\(r\).

At the critical endpoint put \(b=-X''(y_0)>0\). The equation
\(K_Y(r)=0\) and Taylor expansion at \(Y=y_0,r=0\) give \begin{equation}
r(Y)\sim3a(Y)/b,\qquad
 \bar K_Y(w)\longrightarrow1-w,\qquad
 \bar H_Y(w)\longrightarrow1-3w/5.                 \label{eq:P10}
\end{equation} These convergences are locally uniform in complex \(w\). The family
therefore extends continuously to \(Y=y_0\). Compactness and
nonvanishing on the closed unit disk give one \(\eta>0\) such that all
\(\bar H_Y^{-1}\) are analytic and uniformly bounded on
\(|w|\le1+\eta\), including this endpoint. Likewise \(\bar K_Y/(1-w)\)
is nonvanishing on a common disk of radius greater than one: its zeros
are the higher positive roots from \eqref{eq:P5}, and \eqref{eq:P10} treats the endpoint.
Consequently, with \(\theta_{-1}=0\), \begin{equation}
1=\theta_0\le\theta_1\le\cdots\le C,
 \quad \theta_\ell=[w^\ell]\bar K_Y(w)^{-1},
 \quad w_\ell:=\theta_\ell-\theta_{\ell-1}\ge0,
 \quad\sup_Y\sum_{\ell\ge0}(1+\ell)^q w_\ell<\infty.  \label{eq:P11}
\end{equation} Indeed \(\sum w_\ell w^\ell=(1-w)/\bar K_Y(w)\) is the inverse
product over the higher roots. It has nonnegative coefficients and a
common analytic disk. Also \(A(Y)\) is bounded above and away from zero.

Finally, \begin{equation}
P_Y(w)=\frac{\bar H_Y(1)}{\bar H_Y(w)}              \label{eq:P12}
\end{equation} is the probability generating function of a nonnegative integer
variable. Its probabilities at zero and one have uniform positive lower
bounds, its exponential moments are uniformly bounded, and \eqref{eq:P2} at the
first root gives mean \(3/2\). Its variance \(\sigma_Y^2\) is bounded
above and away from zero. At \(Y=y_0\) the law is geometric, with
generating function \((2/5)/(1-3w/5)\), so its variance is \(15/4\).

For a sum \(S_h\) of \(h\) independent copies, Fourier inversion now
gives, uniformly in \(Y\) and in integers \(m\) with
\(|m-3h/2|\le L\sqrt h\), for each fixed \(L\), \begin{equation}
\Pr(S_h=m)=\frac1{\sqrt{2\pi h\sigma_Y^2}}
     \exp\!\left(-\frac{(m-3h/2)^2}{2h\sigma_Y^2}\right)
       +o(h^{-1/2}).                              \label{eq:P13}
\end{equation} Indeed, on a fixed neighborhood of zero the centered
characteristic logarithm is \(-\sigma_Y^2t^2/2+O(|t|^3)\), uniformly by
the exponential moments. Away from zero its modulus is uniformly less
than one, because the first two probabilities have uniform positive
lower bounds. Scaling \(t\) by \(h^{-1/2}\), these facts give an
integrable Gaussian bound and pointwise convergence; Fourier inversion
proves \eqref{eq:P13}. At a bounded displacement from the mean the error is
relative \(O(h^{-1})\), by expansion through fourth order and
cancellation of the odd cubic integral. Thus \[
 [w^m]\bar H_Y(w)^{-h}
   =\bar H_Y(1)^{-h}\Pr(S_h=m)
\] has a uniform central lower bound \(c_L\bar H_Y(1)^{-h}/\sqrt h\),
and its logarithm at a bounded displacement is
\(-h\log\bar H_Y(1)-\tfrac12\log h+O(1)\). These coefficient
estimates enter the harmonic model through
\eqref{eq:H17}--\eqref{eq:H21} in
Section~\ref{exact-amplitudes-and-the-harmonic-model}, and the central
estimate at bounded displacement fixes the normalization of the leading
term in \eqref{eq:PF16}.

\section{The density parameter}\label{the-stationary-parameter}

This section establishes the properties of the functions $T$,
$\Delta=T-Y$, $I$ and $\lambda$ defined in \eqref{eq:SP1} and
\eqref{eq:SP2}; its single result is Lemma~\ref{lem:parameter}, whose
proof, which uses the Sonine integral and elementary calculus, is given
in Appendix~\ref{app:parameter-proof}. The meaning of the four
functions is as follows. The first positive zero of
$K_Y$ is $r(Y)=2\Delta(Y)$; equation \eqref{eq:SP1} is the condition
$K_Y(2\Delta)=0$ written in the variable $v=T-\Delta s^2$. The integral
$I(Y)$ is the action that governs the exponential growth
in~\eqref{eq:intro-heuristic}, and $\lambda(Y)/2$ is the logarithmic
derivative of $I(Y)^{-1}$ with respect to $\log x$, that is, the leading
cusp-count scale per unit complexity. Lemma~\ref{lem:parameter}
shows that $T$ is well defined and that $\lambda$ is an analytic
increasing bijection from $[0,y_0)$ onto $[0,\infty)$, so that the
profiles $\mathcal B,\mathcal C$ of \eqref{eq:SP3} are functions of the
density $\lambda\ge0$; it also records the behavior of all four
functions at the two ends of the interval, which is what the extraction
arguments of Sections~\ref{extracting-fixed-cusp-counts}
and~\ref{sec:absolute} need near $x=0$ and near $x=x_c$.

For the local-limit arguments we view $\lambda$ as a function of the cusp
weight and put \[
 \Psi(x)=\lambda(Y(x))/2,\qquad x=X(Y),\qquad
 \mathscr D=x\frac{\mathrm{d}}{\mathrm{d}x}.
\]

\begin{lemma}[The parameter]\label{lem:parameter}
\begin{enumerate}[label=\textup{(\roman*)},leftmargin=*]
\item The root \(T(Y)\) of \eqref{eq:SP1} in \((y_0,\pi^2/4]\) exists, is unique and depends
analytically on \(Y\in[0,y_0)\). One has \(T(0)=\pi^2/4\),
\(-1\le T'<0\), \(T''>0\), and \(T\to y_0\) as \(Y\to y_0\). Moreover
\(I>0\).
\item Put \begin{equation}
J(Y)=\int_0^1X(T-\Delta s^2)\,\mathrm{d}s
     =X(Y)+\frac{I(Y)}{\sqrt\Delta},
 \qquad
 \beta(Y)=\frac{X(Y)}{J(Y)}=
           \frac{\lambda(Y)}{2+\lambda(Y)}.
 \label{eq:SP11}
\end{equation} Then \(\beta'>0\), and \(\lambda\) is an analytic strictly
increasing bijection from \([0,y_0)\) onto \([0,\infty)\). At \(Y=0\),
\(\Delta=T=\pi^2/4\), \(I(0)=\pi/4\), \(J(0)=1/2\), and
\(\lambda(Y)=4Y+O(Y^2)\).
\item Put \(\delta=y_0-Y\) and \(b=-X''(y_0)=x_c/y_0>0\). As
\(\delta\to0\), \begin{equation}
\Delta=\frac32\delta+O(\delta^2),\qquad
 I(Y)=\frac{8b}{45}\Delta^{5/2}(1+O(\Delta)),
 \label{eq:SP13}
\end{equation} \begin{equation}
\lambda\sim\frac{45x_c}{4b}\Delta^{-2}
       =\frac{45j_{0,1}^2}{16}\Delta^{-2},
 \label{eq:SP14}
\end{equation} and consequently \begin{equation}
\mathcal B(\lambda)\sim
 \frac{\pi^2}{3\sqrt5\,j_{0,1}}\sqrt\lambda,\qquad
 \mathcal C(\lambda)\sim\frac{\pi^2}{x_c}\lambda.
 \label{eq:SP15}
\end{equation} Also \(\Delta\asymp\sqrt{1-\beta}\).
\item The function \(\Psi\) is analytic and strictly increasing for
\(0<x<x_c\). With \(\varepsilon=1-x/x_c\), \begin{equation}
\begin{array}{ll}
 x\downarrow0:&
 \Psi=2x+O(x^2),\quad\mathscr D\Psi=2x+O(x^2),\quad
 \mathscr D^2\Psi=2x+O(x^2),\\[1mm]
 x\uparrow x_c:&
 \displaystyle\Psi=\frac5{4\varepsilon}(1+O(\sqrt\varepsilon)),\quad
 \mathscr D\Psi=\frac5{4\varepsilon^2}(1+O(\sqrt\varepsilon)),\quad
 \mathscr D^2\Psi=O(\varepsilon^{-3}).
 \end{array}
 \label{eq:SP16}
\end{equation} These expansions arise from convergent analytic expansions
in \(x\) or \(\sqrt\varepsilon\). In particular \(\Psi^2/(\mathscr D\Psi)\)
has finite endpoint limits \(0\) and \(5/4\), and is bounded on the full
interval.
\end{enumerate}
\end{lemma}

The estimate \(\Delta\asymp\sqrt{1-\beta}\) in~(iii) implies that the
logarithmic chain profile of
Section~\ref{uniform-logarithmic-volume-rate-and-the-full-diagonal-chain}
is integrable at one. Part~(iv) is restated as \eqref{eq:TL36} in
Section~\ref{extracting-fixed-cusp-counts}, where it provides the
center and variance scales of the local limit theorem, and
\eqref{eq:TL32} identifies \(\Psi\) with the scalar \(x/(Aar)\) that
arises from the marked cusp ratio.

\section{Intersection amplitudes and the comparison model}\label{exact-amplitudes-and-the-harmonic-model}

This section carries out the first idea of
Subsection~\ref{subsec:strategy}. Starting from the alternating
formula~\eqref{eq:H5} for the Weil--Petersson brackets, we move along the
string flow by a time $Y$ and obtain the shifted amplitudes~\eqref{eq:H6},
which interpolate between the brackets at $Y=0$ and pure $\psi$
intersection numbers at $Y=y_0$. Their normalized
form~\eqref{eq:H7} satisfies the recursion~\eqref{eq:H9}, whose
coefficients are the Taylor coefficients $\theta_\ell$ of $1/\bar K_Y$
and hence nonnegative by Section~\ref{the-positive-compact-family-of-kernels}.
We then introduce the harmonic model~\eqref{eq:H12}, an explicit
solvable comparison for this recursion, and record the mass
identity~\eqref{eq:H13} and the coefficient estimates
\eqref{eq:H17}--\eqref{eq:H22} that Section~\ref{uniform-comparison}
uses to compare the two, uniformly in the degrees and in the density
parameter. The external inputs are the intersection identities
of~\cite{KMZ,ManZ,LX} and the fixed-genus radius of convergence
from~\cite{ManZ}; we also use the coarse bounds on normalized pure
$\psi$ intersections proved in Appendix~\ref{app:primitive-bounds}.

\subsection{Intersection identities and normalization}\label{exact-amplitudes-standard-intersection-identities-and-normalization}

We use the following notation for a stable pair \((g,n)\) and a degree
tuple \(\mathbf d\): \(h\) is the complexity, \(N\) the dimension of
\(\overline{\mathcal M}_{g,n}\), and the deficit \(m\) is the number of
\(\kappa_1\) factors needed to fill up the dimension. Write \[
h=2g-2+n,\qquad N=3g-3+n,\qquad
D=\sum_{i=1}^n d_i,\qquad m=N-D.
\] Throughout the comparison of Section~\ref{uniform-comparison}, \(n\ge1\), the pair is stable, and
\(0\le D\le N\). A tuple with a negative degree or an unstable pair has
amplitude zero.

Recall \(X\), \(a=X'>0\) on \([0,y_0)\) and the kernels \(K_Y\) and
\(H_Y\) of \eqref{eq:P1}, and write \(r_Y=r(Y)\) for the first positive
zero of \(K_Y\). The properties of these kernels used below were
established in Section~\ref{the-positive-compact-family-of-kernels}. Put \begin{equation}
\bar H_Y(z)=\frac{H_Y(r_Yz)}a,\qquad
\bar K_Y(z)=\frac{K_Y(r_Yz)}a
      =\bar H_Y(z)+\frac23z\bar H_Y'(z),\qquad
\theta_\ell(Y)=[z^\ell]\bar K_Y(z)^{-1}.
\label{eq:H2}
\end{equation} At \(Y=y_0\) these have the continuous definitions \begin{equation}
\bar H_{y_0}(z)=1-\frac35z,\qquad \bar K_{y_0}(z)=1-z,\qquad
\theta_\ell(y_0)=1.
\label{eq:H3}
\end{equation} In particular \(\theta_0=1\). These are the rescaled kernels and
coefficients of \eqref{eq:P9} and \eqref{eq:P11}, recalled here together
with their endpoint values.

We use the following facts from intersection theory. The ordinary string,
dilaton, and Dijkgraaf--Verlinde--Verlinde (DVV) identities hold for Witten's psi
intersections~\cite{Witten,Kontsevich,DVV}, with
\(\langle\tau_0^3\rangle_0=1\) and \(\langle\tau_1\rangle_1=1/24\). The normalized forms used here are written in Appendix~\ref{app:primitive-bounds}. For
pure $\psi$ intersection numbers we use the normalization \begin{equation}
C_g(\mathbf e)=
\frac{4^g\prod_i(2e_i+1)!!}
 {3^{\,2g-2+q}\Gamma(2g-2+q)}
\left\langle\prod_{i=1}^q\tau_{e_i}\right\rangle_g,
\qquad \sum_i e_i=3g-3+q.
\label{eq:H4}
\end{equation} For positive degrees we use only
\[
 \frac16\le C_g(\mathbf e)\le1.
\]
Appendix~\ref{app:primitive-bounds} proves these coarse bounds directly
from DVV and the published two-point estimate. The argument adapts the
recursion estimates of Guo, Yang and Zagier~\cite{GYZ}; their stronger
lower bound $C_g(3g-2)\le C_g(\mathbf e)$ is not needed here.

With \(m=3g-3+n-|\mathbf d|\), our mixed brackets are
\[
 [\tau_{d_1}\cdots\tau_{d_n}]_{g,n}
 =\prod_i(2d_i+1)!!\,4^{|\mathbf d|}\frac{(2\pi^2)^m}{m!}
 \int_{\overline{\mathcal M}_{g,n}}\kappa_1^m\prod_i\psi_i^{d_i}.
\]
Here \(\kappa_1\) is the logarithmic Mumford class, \(\psi_i\) is the cotangent class, and \([\omega_{\mathrm{WP}}]=2\pi^2\kappa_1\)~\cite{Wolpert}; see also~\cite[Section~2]{MZ}. Thus
\[
 V_{g,n}(2\mathbf L)=\sum_{|\mathbf d|\le3g-3+n}
 [\tau_{d_1}\cdots\tau_{d_n}]_{g,n}
 \prod_i\frac{L_i^{2d_i}}{(2d_i+1)!}.
\]
Angle brackets below denote unweighted pure-$\psi$ integrals. Intersections of the wrong dimension and negative-degree amplitudes are zero. Unstable factors in the positive recursion are zero; explicitly specified genus-zero initial derivatives are kept separate.

The mixed-kappa/psi transform~\cite{KMZ,LX} gives, at
\(Y=0\), \begin{equation}
\frac{[\tau_{d_1}\cdots\tau_{d_n}]_{g,n}}
 {4^D(2\pi^2)^m}
=\prod_i(2d_i+1)!!
\sum_{\substack{\ell_k\ge0,\ k\ge2\\
                 \sum_{k\ge2}(k-1)\ell_k=m}}
\left\langle\tau_{\mathbf d}\prod_{k\ge2}\tau_k^{\ell_k}\right\rangle_g
\prod_{k\ge2}\frac{((-1)^k/(k-1)!)^{\ell_k}}{\ell_k!}.
\label{eq:H5}
\end{equation}
Since \(X'(0)=1\) and \(-X^{(k)}(0)=(-1)^k/(k-1)!\), the weights in
\eqref{eq:H5} are the derivatives \(-X^{(k)}(0)\), \(k\ge2\). The deformation
replaces them by the derivatives \(-X^{(k)}(Y)\) at a point
\(Y\): for general \(Y<y_0\), define the shifted amplitudes by the finite sum
\begin{equation}
U^Y_{g,n}(\mathbf d)=\prod_i(2d_i+1)!!
\sum_{\substack{\ell_k\ge0,\ k\ge2\\
                 \sum(k-1)\ell_k=m}}
\left\langle\tau_{\mathbf d}\prod_{k\ge2}\tau_k^{\ell_k}\right\rangle_g
\frac{\prod_{k\ge2}(-X^{(k)}(Y))^{\ell_k}/\ell_k!}
 {a^{h+\sum\ell_k}}.
\label{eq:H6}
\end{equation} The denominator is obtained by summing degree-one background
insertions with time \(t_1=1-a\): the dilaton equation gives the factor
\(\sum_{j\ge0}(h+\sum\ell_k)_j t_1^j/j!=a^{-h-\sum\ell_k}\), where
\((c)_j=c(c+1)\cdots(c+j-1)\) is the rising factorial. Define the
normalized amplitudes \begin{equation}
\mathcal U_{g,n}(\mathbf d;Y)=a^h r_Y^m U^Y_{g,n}(\mathbf d).
\label{eq:H7}
\end{equation} They extend to \(Y=y_0\) through their positive recursion below.

When all degrees vanish, the shifted amplitudes are derivatives of the
cusp generating function. This is the familiar shift of the Kontsevich--Witten (KW) times
along the string flow~\cite{IZ,KMZ}; the identification of $\Phi_g$ with the
KW free energy at the shifted times $t_0=x$, $t_1=0$, $t_k=(-1)^k/(k-1)!$
is recorded in~\cite[Remark~1]{Zog98}, where it is attributed to Dubrovin, and the corresponding deformation in
$(t_0,t_1)$ is used in~\cite[Section~2.4]{OS}. We recall it in the present
normalization. Start with KW
times \(t_0=x,\ t_1=0,\ t_k=-X^{(k)}(0)\), \(k\ge2\). The finite formal
string flow sends them to \[
t_0'=x-X(Y),\qquad t_1'=1-X'(Y),\qquad
t_k'=-X^{(k)}(Y)\quad(k\ge2).
\] Its vector field is
\(\sum_{k\ge0}t_{k+1}'\partial_{t_k'}-\partial_{t_0'}\). The string
equation makes the genus-\(g\) free energy invariant for \(g\ge1\); its
genus-zero anomaly disappears after three \(t_0\) derivatives. Since
\(\partial_{t_0}\) commutes with this flow, setting \(x=X(Y)\) gives \begin{equation}
U^Y_{g,n}(0^n)=\Phi_g^{(n)}(x),\qquad
\Phi_g(x)=\sum_{n\ge0}\frac{V_{g,n}}{(2\pi^2)^{3g-3+n}}\frac{x^n}{n!},
\label{eq:H8}
\end{equation} whenever the amplitude is stable.
Since the flow is evaluated at a finite time, we also give a direct
argument. Here $a=X'(Y)$. Applying the string and dilaton equations to
the background insertions in~\eqref{eq:H6}, and using
$\partial_Y(-X^{(k)}(Y))=-X^{(k+1)}(Y)$, gives
$\partial_YU^Y_{g,n}(0^n)=a\,U^Y_{g,n+1}(0^{n+1})$, while
$\partial_Y\Phi_g^{(n)}(X(Y))=a\,\Phi_g^{(n+1)}(X(Y))$. At $Y=0$ the two
families agree for every $n$ by~\eqref{eq:H5}; by induction on the
order, all their $Y$-derivatives at $Y=0$ agree. For fixed $(g,n)$ the
sum~\eqref{eq:H6} is finite, so $U^Y_{g,n}(0^n)$ is real analytic on a
neighborhood of $[0,y_0)$, where $a>0$; so is $\Phi_g^{(n)}(X(Y))$, because
$\Phi_g$ has radius of convergence $x_c$~\cite{ManZ} and
$0\le X(Y)<x_c$. Hence~\eqref{eq:H8} holds on $[0,y_0)$; compare the
proof of~\cite[Theorem~1]{Zog98}. Unmarked \(\Phi_g\) is used only
for \(g\ge2\). The unstable genus-zero derivatives are separate initial
terms and are not substituted into the positive recursion.

\subsection{The positive recursion}\label{deriving-the-exact-positive-recursion}

We single out the first marking and call it the \emph{pivot}. Write
\(\langle\!\langle\tau_{d_1}\cdots\tau_{d_q}\rangle\!\rangle_g
=\prod_i(2d_i+1)!!\,\langle\tau_{d_1}\cdots\tau_{d_q}\rangle_g\).
In this normalization the Dijkgraaf--Verlinde--Verlinde equation for the
pivot reads, with \((-1)!!=1\),
\begin{equation*}
\begin{split}
\langle\!\langle\tau_{d_1}\tau_{d_2}\cdots\tau_{d_q}\rangle\!\rangle_g
={}&\sum_{j=2}^q(2d_j+1)\,
 \langle\!\langle\tau_{d_1+d_j-1}\textstyle\prod_{i\ne1,j}\tau_{d_i}\rangle\!\rangle_g
 +\frac12\sum_{u+v=d_1-2}
 \langle\!\langle\tau_u\tau_v\tau_{d_2}\cdots\tau_{d_q}\rangle\!\rangle_{g-1}\\
&+\frac12\sum_{u+v=d_1-2}\ \sum_{\substack{g_1+g_2=g\\ I\sqcup J=\{2,\ldots,q\}}}
 \langle\!\langle\tau_u\tau_I\rangle\!\rangle_{g_1}
 \langle\!\langle\tau_v\tau_J\rangle\!\rangle_{g_2}.
\end{split}
\end{equation*}
For \(d_1=0\) and \(d_1=1\) this is the string and the dilaton equation.
In the amplitudes \eqref{eq:H6}, besides the \(n\) external markings, the
intersection numbers contain \emph{background} insertions: \(\tau_k\) with
weight \(t_k=-X^{(k)}(Y)\) for \(k\ge2\), and, after the dilaton sum that
produced the factor \(a^{-h-\sum\ell_k}\), also \(\tau_1\) with weight
\(t_1=1-a\). Apply the displayed equation with the pivot \(\tau_{d_1}\)
and sum over the background. Three kinds of terms appear.

\begin{enumerate}[label=(\alph*),leftmargin=*]
\item \emph{The pivot merges with an external marking $j\ge2$.} By the
first term above, this gives \((2d_j+1)\,U^Y_{g,n-1}(d_1+d_j-1,\ldots)\),
with the background unchanged.
\item \emph{The pivot merges with a background insertion $\tau_k$,
$k\ge1$.} Background insertions carry no double factorial, so the
factor \(2k+1\) of the first term is divided by \((2k+1)!!\) and becomes
\(1/(2k-1)!!\). Choosing which of the
\(\ell_k\) copies is used turns \(t_k^{\ell_k}/\ell_k!\) into
\(t_k\cdot t_k^{\ell_k-1}/(\ell_k-1)!\). The result is
\[
\frac{t_k}{(2k-1)!!}\,U^Y_{g,n}(d_1+k-1,d_2,\ldots,d_n).
\]
\item \emph{Contractions.} The two remaining terms give the nonseparating
and separating contributions with coefficient \(1/2\). In the separating
case the background splits between the two components with the binomial
factors of exponential generating functions, so the term is a product of
two amplitudes.
\end{enumerate}

Let \(E\) be the operator raising the pivot degree by one. The terms of
type~(b) form \((1-a)U^Y+\sum_{k\ge2}t_kE^{k-1}U^Y/(2k-1)!!\). Moving them
to the left-hand side and recalling \(t_k=-X^{(k)}(Y)\) gives
\[
\left(a+\sum_{k\ge2}\frac{X^{(k)}(Y)}{(2k-1)!!}E^{k-1}\right)U^Y
=K_Y(E)\,U^Y=\bigl(\text{terms of types (a) and (c)}\bigr),
\]
where we compared with the series~\eqref{eq:P1} for \(K_Y\). For a fixed
parent, raising the pivot degree more than \(m\) times leaves the
dimension range, so \(E^{m+1}U^Y=0\) and \(K_Y(E)\) has the finite inverse
\(\sum_{\ell\le m}[z^\ell]K_Y(z)^{-1}E^\ell\). This is how the Taylor
coefficients of \(1/K_Y\) enter. Components of the wrong dimension vanish,
and for external complexity \(h\ge2\) no exceptional initial tuple occurs.

After \eqref{eq:H7}, the exact recursion is \begin{equation}
\begin{aligned}
\mathcal U_{g,n}(\mathbf d)={}&
\sum_{j=2}^n(2d_j+1)\sum_{\ell=0}^m\theta_\ell
 \mathcal U_{g,n-1}(d_1+d_j+\ell-1,\mathbf d_{\widehat{1,j}})\\
&+\frac12\sum_{\ell=0}^m\theta_\ell
 \sum_{u+v=d_1+\ell-2}
 \mathcal U_{g-1,n+1}(u,v,\mathbf d_{\widehat1})\\
&+\frac12\sum_{\ell=0}^m\theta_\ell
 \sum_{u+v=d_1+\ell-2}
 \sum_{\substack{g_1+g_2=g\\I\sqcup J=\{2,\ldots,n\}}}
 \mathcal U_{g_1,|I|+1}(u,\mathbf d_I)
 \mathcal U_{g_2,|J|+1}(v,\mathbf d_J).
\end{aligned}\label{eq:H9}
\end{equation} Here $\mathbf d_{\widehat{1,j}}$ denotes the tuple $\mathbf d$ with the entries $d_1$ and $d_j$ removed, and $\mathbf d_{\widehat1}$ the tuple with $d_1$ removed. The indices $u$ and $v$ are the degrees at the two new markings corresponding to the branches of the node. The last split is ordered. We call the amplitudes on the right of \eqref{eq:H9} the \emph{children} of the amplitude on the left: a merging or nonseparating child has complexity $h-1$, and the two factors of a separating term have complexities summing to $h-1$. The initial values of the recursion, computed directly from
\eqref{eq:H6}, are \begin{equation}
\mathcal U_{0,3}(0,0,0)=1,\qquad
\mathcal U_{1,1}(1)=\frac18,\qquad
\mathcal U_{1,1}(0)=\frac{\theta_1}{8}.
\label{eq:H10}
\end{equation} Equations \eqref{eq:H9}--\eqref{eq:H10}, with negative degrees and unstable factors set
to zero, determine all amplitudes by increasing \(h\).

By \eqref{eq:P11}, \begin{equation}
1=\theta_0\le\theta_1\le\cdots\le C_\theta
\label{eq:H11}
\end{equation} uniformly in \(Y\in[0,y_0]\). Thus every term of \eqref{eq:H9} is nonnegative,
and the continuous endpoint \eqref{eq:H3} defines a positive limiting theory.

\subsection{The harmonic model}\label{harmonic-model-and-the-needed-uniform-coefficient-estimates}

We now introduce the explicit model with which the amplitudes will be
compared. It depends on the degrees only through the deficit \(m\), and it
is chosen so that the linear part of~\eqref{eq:H9} reproduces it exactly
on average; this is the mass identity~\eqref{eq:H13} below. In this sense
the model is harmonic for the linear part of the recursion, which
explains its name. At the
endpoint \(Y=y_0\) the model is the denominator in the normalization
\eqref{eq:H4} of pure \(\psi\) intersection numbers (see~\eqref{eq:PF13}),
so Proposition~\ref{prop:harmonic} extends the coarse bounds on pure \(\psi\)
intersections to the whole parameter interval.

Suppress \(Y\) temporarily and put \begin{equation}
A_m(h)=[z^m]\bar H(z)^{-h},\qquad
\mathcal M_{g,n}(\mathbf d)=3^h\Gamma(h)4^{-g}A_m(h),\qquad
R=\mathcal U/\mathcal M.
\label{eq:H12}
\end{equation} All \(A_m(h)\) are positive. Define \[
q_\ell=\theta_\ell\,\frac{A_{m-\ell}(h-1)}{A_m(h)},
\qquad S_0=\sum_{\ell=0}^m q_\ell,\qquad
L=2D+n=3h-2m.
\] For \(h\ge2\), the identity \(\bar K=\bar H+(2/3)z\bar H'\) gives exactly \begin{equation}
\sum_{\ell=0}^m q_\ell
\left(1-\frac{2(m-\ell)}{3(h-1)}\right)=1.
\label{eq:H13}
\end{equation} Indeed, the generating series of the left side, before division by
\(A_m(h)\), is \[
\bar K^{-1}\left(\bar H^{-(h-1)}
-\frac{2z}{3(h-1)}(\bar H^{-(h-1)})'\right)
=\bar K^{-1}\bar K\bar H^{-h}=\bar H^{-h}.
\]

Dividing \eqref{eq:H9} by the model turns its linear terms into a
weighted sum of the ratios \(R\) at the children. We check that these
weights have total mass one when the pivot degree \(d_1\) is positive;
this is the meaning of \eqref{eq:H13}. Assume \(g\ge2\); states of genus at
most one are terminal in the argument of Section~\ref{uniform-comparison}.
A merging child has type \((g,n-1)\) and a nonseparating child type
\((g-1,n+1)\); both have complexity \(h-1\) and deficit \(m-\ell\). Hence the
model ratio of a merging child to its parent is
\(A_{m-\ell}(h-1)/(3(h-1)A_m(h))\), and that of a nonseparating child is four
times larger, because of the factor \(4^{-g}\) in \eqref{eq:H12}. The
merging children carry the coefficients \((2d_j+1)\theta_\ell\), \(j\ge2\),
and there are \(d_1+\ell-1\) ordered nonseparating pairs, each with
coefficient \(\theta_\ell/2\). Since \(2D+n=3h-2m\), the total weight is
\[
 \sum_{\ell=0}^m q_\ell\,
 \frac{\sum_{j\ge2}(2d_j+1)+2(d_1+\ell-1)}{3(h-1)}
 =\sum_{\ell=0}^m q_\ell\,\frac{3(h-1)-2(m-\ell)}{3(h-1)}=1
\]
by \eqref{eq:H13}. Thus, for a positive pivot, the linear model weights form
a probability kernel \(P\), and positivity of the separating terms gives
\(R\ge PR\).

At an all-zero state, where every degree and hence the pivot vanishes,
the formal \(\ell=0\) count includes one
nonseparating pair with multiplicity \(-1\), and also includes \(n-1\)
invalid mergers. The true linear mass is instead \begin{equation}
s=1+\frac{3-n}{3(h-1)}
       \frac{A_m(h-1)}{A_m(h)}.
\label{eq:H14}
\end{equation} Since \(\bar H(0)=1\) and \(\bar H^{-1}\) has positive coefficients,
\(A_m(h-1)\le A_m(h)\). In the range \(g\ge2\), therefore, \begin{equation}
s\ge\frac23,\qquad
s\le1\ (n\ge3),\qquad s\le1+O(h^{-1})\ (n=1,2).
\label{eq:H15}
\end{equation} Thus for a zero pivot the linear terms need not have mass exactly one.

We now record the analytic information needed below. The function
\begin{equation}
P_Y(w)=\frac{\bar H_Y(1)}{\bar H_Y(w)}
\label{eq:H16}
\end{equation} is the probability generating function of
\eqref{eq:P12}; its mean \(3/2\), its analyticity on a disk of radius
\(1+\eta>1\), the positive lower bounds on its first two coefficients and
the two-sided bounds on its variance, all uniform in \(Y\in[0,y_0]\),
were established in
Subsection~\ref{compact-normalization-and-uniform-lattice-estimates}.
Let \(S_h\) be the sum of \(h\)
independent variables with probability generating function
\eqref{eq:H16}. The uniform local limit theorem \eqref{eq:P13} gives \begin{equation}
\Pr(S_h=m)\ge c_Bh^{-1/2}
\quad\text{when }|m-3h/2|\le B\sqrt h,
\label{eq:H17}
\end{equation} and at bounded displacement, where the odd cubic term
integrates to zero, \begin{equation}
\Pr(S_h=3h/2+d)
=\frac{\mathrm{e}^{-d^2/(2h\sigma_Y^2)}}{\sqrt{2\pi h\sigma_Y^2}}
  (1+O_B(h^{-1})),\qquad |d|\le B,
\label{eq:H18}
\end{equation} at integer arguments. The same Fourier argument applies after tilting at a
positive radius \(\rho\) in a compact subset of the analytic interval; here $\rho$ is a radius in the normalized $z$-coordinate. If
\(\mu_Y(\rho)=-\rho\bar H_Y'(\rho)/\bar H_Y(\rho)=m/h\), it gives \begin{equation}
A_m(h)=
\frac{\bar H_Y(\rho)^{-h}\rho^{-m}}{\sqrt{2\pi h\,\rho\mu_Y'(\rho)}}
 (1+O(h^{-1}))
\label{eq:H19}
\end{equation} uniformly on such compact tilted families. Near the mean, the tilted
rate function has bounded second derivative and vanishing value and
first derivative; hence, for a fixed sufficiently small \(c_0>0\), \begin{equation}
\Pr(S_h=3h/2-u)\ge ch^{-1/2}\mathrm{e}^{-Cu^2/h},
\qquad 0\le u\le c_0h.
\label{eq:H20}
\end{equation}

The coefficient/probability relation is \begin{equation}
A_m(h)=\bar H(1)^{-h}\Pr(S_h=m).
\label{eq:H21}
\end{equation} Consequently, when \(L\le\sqrt h\), \eqref{eq:H11}, \eqref{eq:H17}, and the sum of the
numerator probabilities give \(S_0\le C\sqrt h\). When \(L\ge6\), the
minimum factor in \eqref{eq:H13} is \((L-3)/(3(h-1))\), so \(S_0\le Ch/L\).
Combining the two estimates, \begin{equation}
S_0\le C\sqrt h,\qquad S_0\le Ch/L\quad(L\ge6).
\label{eq:H22}
\end{equation} All constants here are uniform in \(Y\), and the finitely many states
of bounded complexity are absorbed into them. The mass identity
\eqref{eq:H13}, the zero-pivot correction \eqref{eq:H14}--\eqref{eq:H15}
and the bound \eqref{eq:H22} on \(S_0\) are the inputs of the random-walk
argument in Section~\ref{uniform-comparison}; the central estimate
\eqref{eq:H18} is used again in \eqref{eq:PF16} to normalize the leading
term.

\section{Comparison with the harmonic model}\label{uniform-comparison}

This section completes the first idea of
Subsection~\ref{subsec:strategy} by proving
Proposition~\ref{prop:harmonic}, the coarse two-sided comparison between
the normalized amplitudes of \eqref{eq:H9} and the harmonic
model~\eqref{eq:H12} on which all later estimates rest. It loses a
factor $\mathrm e^{O(\sqrt h)}$ but holds uniformly in the degrees and in
the density parameter, including the pure-$\psi$ endpoint, and it uses
no information about volumes: its inputs are the positive recursion, the
model estimates \eqref{eq:H13}--\eqref{eq:H22} of
Section~\ref{exact-amplitudes-and-the-harmonic-model} and the coarse
bounds on pure $\psi$ intersections from
Appendix~\ref{app:primitive-bounds}. The proof is given in full here.
Section~\ref{ratios-with-a-fixed-number-of-distinguished-markings}
upgrades the comparison to relative errors for a fixed number of
markings, and Section~\ref{uniform-one-boundary-coefficient-ratios}
uses its $Y=0$ case, recorded in \eqref{eq:H44}, for the coefficient
ratios.

\begin{proposition}[Uniform harmonic comparison]\label{prop:harmonic}
For every fixed \(n\ge1\), there is a constant \(C_n\) such that
\[
 -C_n\log^2(h+1)\le
 \log\frac{\mathcal U_{g,n}(\mathbf d;Y)}
 {3^h\Gamma(h)4^{-g}[z^m]\bar H_Y(z)^{-h}}
 \le C_n\sqrt h
\]
for every stable pair \((g,n)\), every \(\mathbf d\in\mathbb Z_{\ge0}^n\)
with \(|\mathbf d|\le3g-3+n\), and every \(Y\in[0,y_0]\). Here
\(h=2g-2+n\) and \(m=3g-3+n-|\mathbf d|\).
\end{proposition}
The idea is to read the recursion as a random walk. Write \(R=\mathcal U/\mathcal M\) for the
ratio of a true amplitude to the model amplitude~\eqref{eq:H12}, and
divide the recursion~\eqref{eq:H9} by the model. Because the model
satisfies the linear part of the recursion with total mass one at a
positive pivot (identity~\eqref{eq:H13}), the linear terms there express \(R\) as an
average of \(R\) over the states of complexity \(h-1\) reached in one step,
and the separating terms add a nonnegative amount. At all-zero states the linear mass need not equal one, and it is carried along in the averaging. Thus \(R\) is bounded
below by a terminal expectation weighted by the accumulated linear masses,
and it is bounded above as soon as the separating contribution is small.

Comparing intersection recursions with random walks has an
important precedent in Aggarwal's proof of uniform asymptotics for pure
$\psi$ intersections~\cite{Agg}. Here, in addition, the recursion and its
comparison model vary with $Y$, and the estimates must remain uniform as
this parameter approaches the pure-$\psi$ endpoint.

Throughout this section a \emph{state} is a stable pair \((g,n)\), \(n\ge1\),
with a degree tuple \(\mathbf d\in\mathbb Z_{\ge0}^n\), \(|\mathbf d|\le3g-3+n\),
and a parameter \(Y\in[0,y_0]\); constants are uniform in \(Y\) and in the
state unless stated otherwise. We write \(Z_0\) for the number of entries
of \(\mathbf d\) equal to zero. A state is \emph{terminal} if \(m=0\) or
\(g\le1\); the random walks below stop there. The proof has three steps.
\begin{itemize}[leftmargin=*]
\item Lemma~\ref{lem:terminal} bounds \(R\) at terminal states, using only
the string and dilaton equations and the coarse bounds of
Appendix~\ref{app:primitive-bounds}.
\item Lemma~\ref{lem:lower} proves the lower bound. The only defect of the
averaging is at entries of degree zero, where some linear terms are
missing; we show that the random walk does not create many such entries.
\item Lemmas~\ref{lem:weighted-linear}--\ref{lem:upper} prove the upper
bound. We divide \(R\) by a weight that penalizes entries of degree zero
and one, show that the weighted separating terms have relative size
\(O(h^{-1/2})\), and conclude by examining, at each complexity, the state
at which the weighted ratio is largest.
\end{itemize}

\subsection{Bounds at terminal states}\label{elementary-terminal-bounds-and-a-global-crude-estimate}

\begin{lemma}[Terminal bounds]\label{lem:terminal}
There are constants \(c_\psi,C_\psi,c,C>0\) with the following properties.
\begin{enumerate}[label=\textup{(\roman*)},leftmargin=*]
\item For every pure \(\psi\) tuple \(\mathbf e\) with \(Z_0\) zero entries,
\(c_\psi(2/3)^{Z_0}\le C_g(\mathbf e)\le C_\psi2^{Z_0}\), in the normalization
\eqref{eq:H4}.
\item At every terminal state, \(R\ge c\,\mathrm e^{-CZ_0}\).
\item At every state, \(R\le\mathrm e^{Ch}\).
\end{enumerate}
\end{lemma}

\begin{proof}
For~(i), the standard normalized string and dilaton identities imply \begin{equation}
c_\psi(2/3)^{Z_0}\le C_g(\mathbf e)\le C_\psi2^{Z_0},
\label{eq:H23}
\end{equation} where \(Z_0\) is the number of zero entries. To prove
this, first remove the entries equal to one by the normalized dilaton equation,
retaining the final genus-one tuple \((1)\). At a noninitial tuple with
a zero, let \(p_+\) denote the number of positive entries; normalized
string has total positive coefficient \[
\frac{2(3g-3+q)+p_+}{3(2g-3+q)}
\] between \(2/3\) and \(2\). Once the ones have been removed, each
string child has one fewer zero; any newly created one can again be
removed. With no zeros use the primitive bounds stated after \eqref{eq:H4}; the
remaining bases are \(C_0(0,0,0)=1/3\) and \(C_1(1)=1/6\). Induction
proves \eqref{eq:H23}, which is~(i).

For~(ii) we use a lower comparison that requires no asymptotic
information on volumes. Set all kernel coefficients in \eqref{eq:H9} equal to one. The
corresponding background of KW times has \(t_1=0,t_2=3,t_k=0\) for \(k\ge3\), so its amplitudes
are \begin{equation}
\mathcal U^{\rm PI}_{g,n}(\mathbf d)
=\frac{3^m}{m!}\prod_i(2d_i+1)!!
\left\langle\tau_{\mathbf d}\tau_2^m\right\rangle_g.
\label{eq:H24}
\end{equation} They have the same initial \((0,3)\) and degree-one \((1,1)\) values,
and their degree-zero \((1,1)\) value \(1/8\) is at most
\(\mathcal U_{1,1}(0)=\theta_1/8\). By \eqref{eq:H11}, induction in the
positive recursion gives \(\mathcal U\ge\mathcal U^{\rm PI}\).

Uniform analyticity of \(\bar H_Y^{-1}\) at zero gives a fixed \(B\ge1\) such
that \([z^j]\bar H_Y^{-1}\le B^j\). Thus \[
A_m(h)\le B^m(h)_m/m!.
\] Dividing \eqref{eq:H24} by \eqref{eq:H12} and using \eqref{eq:H4} gives \begin{equation}
R\ge C_g(\mathbf d,2^m)
       \left(\frac3{5B}\right)^m.
\label{eq:H25}
\end{equation} If \(g\le1\), dimension and \(D\ge n-Z_0\) imply \(m\le Z_0\). Hence
\(R\ge c\mathrm{e}^{-CZ_0}\) in this terminal range. At \(m=0\), exactly
\(R=C_g(\mathbf d)\), so \eqref{eq:H23} applies directly. This proves~(ii).

Finally we prove the global bound~(iii), \begin{equation}
R\le \mathrm{e}^{Ch}.
\label{eq:H26}
\end{equation} Note that in the normalized finite sum \eqref{eq:H6}--\eqref{eq:H7} a
background insertion of degree \(k\) has weight \[
v_k(Y)=-r_Y^{k-1}X^{(k)}(Y)/a.
\] The corresponding absolute majorant is \[
B_Y(z)=\sum_{k\ge2}\frac{3|v_k(Y)|}{(2k+1)!!}z^{k-1}.
\] This family has a uniformly bounded analytic germ and vanishes at
zero. Indeed it is the coefficientwise absolute version of \(1-\bar H_Y\);
the normalized family in \eqref{eq:H2}--\eqref{eq:H3} is compact and analytic on a common
disk. Thus a fixed \(z_*>0\) makes \(B_Y(z_*)\le1/2\). Using the upper
bound in \eqref{eq:H23} on every primitive term of \eqref{eq:H6} gives \[
\mathcal U\le
C2^{Z_0}3^h\Gamma(h)4^{-g}[z^m](1-B_Y(z))^{-h}
\le C2^{Z_0}3^h\Gamma(h)4^{-g}z_*^{-m}2^h.
\] The linear coefficient of \(-\log \bar H_Y\) is uniformly positive.
Therefore \[
A_m(h)\ge (ch)^m/m!\ge \mathrm{e}^{-C'h},
\] because \(m\le3h/2\), and \(A_0(h)=1\) since \(\bar H_Y(0)=1\). Since
\(Z_0\le n\le h+2\), division proves \eqref{eq:H26}.
\end{proof}

\subsection{The lower bound: control of zero degrees}\label{lower-comparison-by-controlling-only-zeros}

\begin{lemma}[Lower bound]\label{lem:lower}
There are constants \(c,C>0\) such that, at every state,
\[
 R\ge c\exp\{-C(Z_0+\log(h+1))\log(h+1)\}.
\]
\end{lemma}

\begin{proof}
Since the separating terms of \eqref{eq:H9} are nonnegative, \(R\) is at
least the total linear mass times the average of \(R\) over its linear
children under the normalized transition probabilities. Iterating gives
a lower bound by the terminal value weighted by the product of these
masses along a walk that lowers complexity by one at each step. Two losses occur: the terminal
bounds \eqref{eq:H23} and \eqref{eq:H25} lose a factor exponential in the
number of zero entries, and the averaging may have mass \(s<1\) at states
with all degrees zero, again with a loss controlled by the number of zeros. It
therefore suffices to show that the walk rarely creates many zero entries.

At a nonterminal state choose any positive pivot if one exists;
otherwise choose a zero. Stop only at \(m=0\) or \(g\le1\). Use the
linear model weights, divided by their mass \(s\), as transition
probabilities. For a positive pivot \(s=1\); at an all-zero state use
\eqref{eq:H14}--\eqref{eq:H15}.

Let \(Z\) be the current zero count. At a positive pivot, a merger with a marking of degree zero
lowers \(Z\) by one, except when the pivot has degree one and \(\ell=0\);
that merger leaves \(Z\) unchanged. A nonseparating ordered pair can increase \(Z\) by
at most two, and at each \(\ell\) there are at most two pairs containing
a zero. Each pair has twice the model weight of a merger with one zero
marking. Thus \begin{equation}
p_{\rm birth}\le4p_{\rm death}/Z\qquad(Z\ge1).
\label{eq:H27}
\end{equation} That exceptional term contributes neither a
birth nor a death. At an all-zero state, a birth occurs only at
\(\ell=2\), raises the count by one, and its weight is at most
\(2/(Z-1)\) times the mass of mergers with a zero marking. Those mergers have
\(\ell\ge1\) and lower the count by at least one. Hence \eqref{eq:H27} still
holds for \(Z\ge2\). Normalization by \(s\) preserves these comparisons.

For \(Z\ge24\), the conditional drift obeys \[
\frac{\mathbb E(2^{Z_{\rm next}})}{2^Z}-1
\le3p_{\rm birth}-\frac12p_{\rm death}\le0.
\] During each excursion above 24, optional stopping bounds the
probability of reaching \(M\) by \(2^{Z_{\rm start}-M}\). An excursion
entered from below starts at a value at most 26. Since complexity decreases by
one, there are at most \(h_{\rm initial}+1\) excursions. A union bound
therefore shows that, with probability at least \(1/2\), \begin{equation}
\max Z\le Z_{\rm initial}+C+C\log(h_{\rm initial}+1).
\label{eq:H28}
\end{equation} All stopping times are bounded, so optional stopping
applies without further integrability conditions.

Only all-zero steps can have \(s<1\). By \eqref{eq:H14}--\eqref{eq:H15}, \[
-\log\min(s,1)\le CZ/h.
\] On the event \eqref{eq:H28}, write its upper bound as \(M\). Summing over the
decreasing complexities costs at most \[
CM+CM\sum_{j=M+1}^{h_{\rm initial}}\frac1j
\le CM\log(h_{\rm initial}+1).
\] Iterated positivity gives \[
R_{\rm initial}\ge
\mathbb E\!\left[R_{\rm terminal}\prod_{\rm steps}s\right].
\] Here the initial zero count is \(Z_0\). The terminal bounds of
Lemma~\ref{lem:terminal}(ii), applied at terminal zero counts at most
\(M\), and the probability in \eqref{eq:H28} now prove \begin{equation}
R\ge c\exp\{-C(Z_0+\log(h+1))\log(h+1)\}.\label{eq:H29}
\end{equation}
No asymptotics for volumes enter at the terminal states.
\end{proof}

\subsection{The upper bound: a weight for degrees zero and one}\label{upper-comparison-the-clipped-inactive-weight}

The upper bound requires the separating terms to be small, and
for general states they are not: a state with many entries of degree zero
or one can distribute these markings between the two components in many
ways. We therefore divide \(R\) by a weight that grows with the number of
such entries, chosen so that the linear terms still average
\(R/\mathcal W\) up to a factor \(1+O(h^{-1/2})\) (this is \eqref{eq:H31}),
while the weighted separating terms become small
(Subsection~\ref{the-weighted-separating-estimate}).

Call degrees zero or one inactive. Let \(\iota\) be their number, \(Z_0\)
the number of zeros, and \(R_{\rm act}\) the number of degrees at least
two. Fix a sufficiently small constant \(c_*>0\), and set \begin{equation}
\phi(\iota)=\prod_{j=1}^\iota\max(1,c_*\sqrt j),\qquad
\mathcal W(\mathbf d)=4^{Z_0}\phi(\iota).
\label{eq:H30}
\end{equation} We call the truncation \(\max(1,\cdot)\) clipping. Since
every factor is then at least one, \(\phi\) is nondecreasing; this
monotonicity is essential when an inactive pivot splits into two active
degrees. Up to fixed multiplicative constants, \(\phi(\iota)\) equals
\(c_*^\iota\sqrt{\iota!}\).

In the upper bound we choose the pivot as follows: an active entry when
one exists, an entry equal to one if no active entry exists, and a zero
only at an all-zero state. Again the walk stops only at terminal states.
The linear model weights \(p_{\rm child}\) of a state are the weights
that appear when \eqref{eq:H9} is divided by the model; their mass is
given by \eqref{eq:H13}--\eqref{eq:H15}, and may exceed one at an
all-zero state with one or two markings.

\begin{lemma}[Weighted linear terms]\label{lem:weighted-linear}
If \(c_*\) is sufficiently small, then at every nonterminal state
\begin{equation}
\sum_{\rm linear}p_{\rm child}
\frac{\mathcal W_{\rm child}}{\mathcal W}
\le1+CS_0/h\le1+Ch^{-1/2}.
\label{eq:H31}
\end{equation}
\end{lemma}

\begin{proof}
The second inequality is \eqref{eq:H22}. We consider the three possible
pivots in turn.

For an active pivot, an ordinary inactive merger reduces \(\iota\) by at
least one and does not increase \(Z_0\). For all sufficiently large
\(\iota\), its weight ratio is at most \(1/2\). The exception is a pivot of degree two
merged with a zero at \(\ell=0\): it turns that zero into a one, leaves
\(\iota\) fixed, and divides the weight by four. The favorable mass is
therefore at least a fixed positive multiple of \(\iota S_0/h\).

At each \(\ell\), at most four nonseparating pairs contain an inactive
degree. A pair creating two inactive degrees has weight ratio, for large
\(\iota\), at most \[
16c_*^2\sqrt{(\iota+1)(\iota+2)}.
\] Its total model mass per \(\ell\) is bounded by a fixed multiple of
\(q_\ell/h\). Choose \(c_*\) so small that the coefficient of \(\iota\) in
this upper bound is less than half the favorable merger coefficient.
Pairs creating only one inactive degree cost \(O(\sqrt \iota)\), and all
other pairs leave the weight unchanged. The remaining expression is
bounded above by \(CS_0/h\). Bounded \(\iota\), including the clipping
threshold, also costs only \(CS_0/h\), because there are at most four
exceptional pairs and boundedly many inactive mergers.

If the pivot has degree one and all markings are inactive, every merger
decreases \(\iota\), while a nonseparating step increases it by at most one.
Only four pairs per \(\ell\) can increase the weight; a pair with two
active new degrees decreases \(\iota\) and cannot increase the clipped
weight. The same argument applies, with \(\iota-1\) favorable merging
labels.

At an all-zero state, all valid mergers have \(\ell\ge1\) and decrease
both the zero count and the inactive weight. Increasing pairs have
\(\ell\ge2\) and increase \(\iota\) by at most one. Use \(S_0-q_0\) in both
birth and merger estimates. The unweighted mass correction \eqref{eq:H14} is
nonpositive for \(n\ge3\) and \(O(h^{-1})\) otherwise. This proves \eqref{eq:H31}
in every case.
\end{proof}

\subsection{Bounding the separating terms}\label{the-weighted-separating-estimate}

We now show that, after weighting, the separating terms of
\eqref{eq:H9} have total relative size \(O(h^{-1/2})\), uniformly in the
state. Splitting off a component of complexity \(h_1\le h_2\) costs the
factor \(\Gamma(h_1)\Gamma(h_2)/\Gamma(h)\), which is at most
\(\binom{h-3}{h_1-1}^{-1}\). The key point is that this factor beats the
number of ways of distributing the marking labels once the weight has
absorbed the inactive ones. For a separating term with node degrees \(u,v\) and kernel
index \(\ell\), write \(\mathcal M_1,\mathcal M_2\) and
\(\mathcal W_1,\mathcal W_2\) for the model amplitudes and weights of its two
factors.

\begin{lemma}[Weighted separating terms]\label{lem:weighted-separating}
Let \(\eta=1/12\). At every nonterminal state,
\[
 \sum_{\rm separating}\frac{\theta_\ell}{2}\,
 \frac{\mathcal M_1\mathcal M_2}{\mathcal M}\,
 \frac{\mathcal W_1\mathcal W_2}{\mathcal W}\,
 \mathrm e^{\eta\min(h_1,h_2)}
 \le CS_0h^{-3/2}(1+\iota)\le Ch^{-1/2},
\]
where the sum runs over all separating terms of \eqref{eq:H9}.
\end{lemma}

\begin{proof}
The degree-dimension constraint supplies \begin{equation}
4R_{\rm act}+n\le3h,\qquad R_{\rm act}\le3h/5.
\label{eq:H32}
\end{equation} Remove the pivot, leaving \(\iota_*\) inactive and \(R_*\) active labeled
markings, with \(\iota_*+R_*=n-1\), \(\iota_*\le \iota\), and \(R_*\le R_{\rm act}\). In
one separating component put \(u\) inactive and \(j\) active original
labels, and write \(s=u+j\). Each child receives one new marking.
Comparing \eqref{eq:H30} with \(4^{Z_0}c_*^\iota\sqrt{\iota!}\), uniformly over the two
new degrees, \begin{equation}
\frac{\mathcal W_1\mathcal W_2}{\mathcal W}
\le C\frac{\sqrt{(s+2)h}}{\sqrt{\binom{\iota_*}u}}.
\label{eq:H33}
\end{equation} Indeed, when the \(\iota_*\) inactive original labels split as \(u,\iota_*-u\), the factorial part of this estimate is bounded by
\[
 \frac{\sqrt{(u+1)!(\iota_*-u+1)!}}{\sqrt{\iota_*!}}
 =\frac{\sqrt{(u+1)(\iota_*-u+1)}}{\sqrt{\binom{\iota_*}u}}.
\]
Each child receives at most one inactive new marking. Since \(u\le s\) and \(\iota_*\le h+2\), the numerator is at most \(C\sqrt{(s+2)h}\). The powers of \(c_*\), the factor for zeros, and the clipping comparison change by bounded integer exponents. An inactive pivot supplies an additional factor \(\iota^{-1/2}\), which
was discarded. The exponential factors and clipping comparison
contribute only a fixed constant, since the total inactive and zero
counts change by bounded amounts.

Orient \(h_1\le h_2\), and put \[
h'=h-3,\qquad k=h_1-1\le h'/2,\qquad s=k+2-2g_1.
\] Then \begin{equation}
\frac{\Gamma(h_1)\Gamma(h_2)}{\Gamma(h)}
=\frac1{(h-1)(h-2)\binom{h'}k}.
\label{eq:H34}
\end{equation} For all sufficiently large \(h\), \eqref{eq:H32} gives \(R_*\le2h'/3\).

The required label sum is \[
\Lambda_s=\sum_{u+j=s}\binom{\iota_*}u^{1/2}\binom{R_*}j.
\] Cauchy--Schwarz and Vandermonde give \begin{equation}
\Lambda_s\le
\binom{\iota_*+R_*}s^{1/2}
\left(\sum_{j\le s}\binom{R_*}j\right)^{1/2}.
\label{eq:H35}
\end{equation} For \(k\ge1\), take \(\xi =k/(h'-k)\le1\). The two factors before square
roots are bounded by \(\xi ^{-s}(1+\xi )^{\iota_*+R_*}\) and \(\xi ^{-s}(1+\xi )^{R_*}\),
respectively. The binomial mass at its mean gives \[
\binom{h'}k\ge c(k+1)^{-1/2}\xi ^{-k}(1+\xi )^{h'}.
\] Thus, if \(g_1\ge1\), so \(s\le k\), \begin{equation}
\frac{\Lambda_s}{\binom{h'}k}
\le C\sqrt{k+1}\,\xi ^{k-s}(1+\xi )^{-h'/6}
\le C\sqrt{k+1}\mathrm{e}^{-k/6}.
\label{eq:H36}
\end{equation} Here \(\iota_*+R_*=n-1\le h'\), \(R_*\le2h'/3\), and \(h'\log(1+\xi )\ge k\). The
case \(k=0\) is immediate.

If \(g_1=0\), then \(s=k+2\). Validity of this component forces
\(2j\le s-2=u+j-2\), hence \(u\ge j+2\ge2\). Use \[
\binom{\iota_*}u=\frac{\iota_*(\iota_*-1)}{u(u-1)}\binom{\iota_*-2}{u-2}.
\] After taking the square root, discard the denominator and extract a
factor \(\iota_*\). The remaining sum has \(u-2+j=k\), so \eqref{eq:H35}--\eqref{eq:H36} apply
with \(\iota_*-2\) in place of \(\iota_*\). This use of the degree constraint is
what replaces a factorial penalty on active labels.

Combining these estimates with \eqref{eq:H33}--\eqref{eq:H34}, and summing the
polynomially many genus choices at fixed \(k\), proves, for example with
\(\eta=1/12\), \begin{equation}
\sum_{\rm genus/label\ splits}
\frac{\Gamma(h_1)\Gamma(h_2)}{\Gamma(h)}
\sup_{\rm valid\ node\ degrees}
\frac{\mathcal W_1\mathcal W_2}{\mathcal W}
\mathrm{e}^{\eta\min(h_1,h_2)}
\le C h^{-3/2}(1+\iota).\label{eq:H37}
\end{equation} An empty supremum contributes zero. Thus the bound is uniform in the
possible node degrees; their sum, with its model coefficients, is
handled next.

For a separating recursion term, exact model normalization gives \begin{equation}
\frac{\theta_\ell}{2}\frac{\mathcal M_1\mathcal M_2}{\mathcal M}
=\frac{\theta_\ell}{6}
\frac{\Gamma(h_1)\Gamma(h_2)}{\Gamma(h)}
\frac{A_{m_1}(h_1)A_{m_2}(h_2)}{A_m(h)}.
\label{eq:H38}
\end{equation} Since \(h_1+h_2=h-1\) and \(m_1+m_2=m-\ell\), positivity and
coefficient convolution bound the sum over node degrees by
\(A_{m-\ell}(h-1)\). More explicitly, for each genus/label partition
first apply the node-independent envelope \eqref{eq:H33}, then extend the allowed
nonnegative model-coefficient sum to the full convolution. The resulting
label envelope is exactly the one summed in \eqref{eq:H35}--\eqref{eq:H37}, so variable
degree supports cause no interchange issue. Consequently the entire
weighted separating mass, including the exponential factor
\(\mathrm e^{\eta\min(h_1,h_2)}\) of \eqref{eq:H37}, is at most \begin{equation}
CS_0h^{-3/2}(1+\iota)=O(h^{-1/2}).
\label{eq:H39}
\end{equation} For \(L\ge6\), use \(S_0\le Ch/L\) and \(\iota\le n\le L\). For
\(L\le5\), both \(\iota\) and \(L\) are bounded, and \(S_0\le C\sqrt h\)
gives the smaller \(O(h^{-1})\) bound.
\end{proof}

\subsection{The upper bound from successive maxima}\label{terminal-upper-bounds-and-the-two-record-arguments}

\begin{lemma}[Upper bound]\label{lem:upper}
There is a constant \(C\) such that \(R\le C\mathcal W\mathrm e^{C\sqrt h}\)
at every state.
\end{lemma}

\begin{proof}
By Lemmas~\ref{lem:weighted-linear} and~\ref{lem:weighted-separating}, the weighted ratio satisfies an
approximate maximum principle: at each nonterminal state, \(R/\mathcal W\)
is at most \(1+O(h^{-1/2})\) times the largest value at its linear
children, plus a separating contribution that is small. We apply this at
complexities where the maximum of \(R/\mathcal W\) reaches a record, first
to exclude exponential growth \eqref{eq:H40}, and then to obtain the bound
\(\mathrm e^{C\sqrt h}\) in \eqref{eq:H41}.

At \(m=0\), \eqref{eq:H23} and \eqref{eq:H30} give \(R/\mathcal W\le C\). At \(g\le1\),
one has \(D\le n\), so \(R_{\rm act}\le n/2\) and \(\iota\ge n/2\). Thus
\(\log\mathcal W\ge(n/4)\log n-O(n)\), which dominates \eqref{eq:H26}, since
\(h=n+O(1)\). Hence \(R/\mathcal W\) is bounded at every terminal state.

Let \[
M_h=\max_{\substack{Y\in[0,y_0]\\\text{valid states of complexity }h}}
\frac{R}{\mathcal W}.
\] The maximum exists: the tuple set is finite and the normalized
parameter interval is compact. Equation \eqref{eq:H26} gives
\(\gamma=\limsup h^{-1}\log M_h<\infty\). Suppose \(\gamma>0\), and
choose \(0<\gamma'<\gamma\), \(\delta>0\), so close that
\(\gamma+\delta-\gamma'<\eta\) from \eqref{eq:H37}. Call \(h\) a \emph{record} for a sequence \(a_h\) if \(a_h>a_j\) for all \(j<h\). There are arbitrarily large
records of \(M_h\mathrm{e}^{-\gamma'h}\), and their maximizing states are
nonterminal. At such a record, \[
M_j/M_h\le \mathrm{e}^{-\gamma'(h-j)}\quad(j<h),\qquad
M_j\le C_\delta \mathrm{e}^{(\gamma+\delta)j}.
\] The linear part, after division by the maximizing value, is bounded
by \((1+Ch^{-1/2})\mathrm{e}^{-\gamma'}\). For a separating term use the record
inequality on the larger child and the global bound on the smaller.
Their product divided by \(M_h\) is at most a fixed constant times
\(\mathrm{e}^{(\gamma+\delta-\gamma')h_{\rm small}}\). Equation \eqref{eq:H39}, which includes the factor
\(\mathrm e^{\eta\min(h_1,h_2)}\), makes the total separating contribution \(o(1)\). The recursion
would therefore give \(1\le \mathrm{e}^{-\gamma'}+o(1)\), a contradiction. Thus
\begin{equation}
\limsup h^{-1}\log M_h\le0.
\label{eq:H40}
\end{equation}

To obtain the quantitative envelope, fix \(0<\delta<\eta\), so \eqref{eq:H40}
gives \(M_j\le C_\delta \mathrm{e}^{\delta j}\). Consider record maxima of
\(M_h/\mathrm{e}^{A\sqrt h}\). If these were unbounded, their large maximizing
states would again be nonterminal. For the larger separating child, the
record bound supplies
\(M_{\rm large}/M_h\le \mathrm{e}^{A(\sqrt{h_{\rm large}}-\sqrt h)}\le1\); bound
the smaller child by \(C_\delta \mathrm{e}^{\delta h_{\rm small}}\). The
separating contribution is therefore \(C'h^{-1/2}\), with \(C'\)
independent of \(A\). The linear record ratio is \[
\mathrm{e}^{A(\sqrt{h-1}-\sqrt h)}
=1-\frac{A}{2\sqrt h}+O_A(h^{-1}).
\] Using \eqref{eq:H31}, the recursion would give \[
1\le1+\frac{C+C'-A/2}{\sqrt h}+O_A(h^{-1}).
\] Choosing \(A>2(C+C'+1)\) excludes unbounded records. We have proved
\begin{equation}
R/\mathcal W\le C \mathrm{e}^{C\sqrt h},
\label{eq:H41}
\end{equation} uniformly for \(Y\in[0,y_0]\).
\end{proof}

\subsection{Completion of the comparison}\label{common-comparison-theorem-and-the-original-coefficients}

\begin{proof}[Completion of the proof of Proposition~\ref{prop:harmonic}]
For fixed \(n\), the zero count and the inactive count are at most \(n\), so
\(\mathcal W\le4^n\phi(n)\) is bounded. Lemmas~\ref{lem:lower}
and~\ref{lem:upper} therefore give, for every fixed \(n\ge1\), \begin{equation}
-C_n\log^2(h+1)\le
\log\frac{\mathcal U_{g,n}(\mathbf d;Y)}
 {3^h\Gamma(h)4^{-g}[z^m]\bar H_Y(z)^{-h}}
\le C_n\sqrt h,\label{eq:H42}
\end{equation} uniformly for \(Y\in[0,y_0]\) and every degree tuple
in the range of Proposition~\ref{prop:harmonic}. This proves the proposition.
\end{proof}

In particular, for each fixed \(n\) the logarithmic error in
\eqref{eq:H42} is \(o(g)\), uniformly in the degrees and in \(Y\); this is
the form in which the comparison enters
Section~\ref{ratios-with-a-fixed-number-of-distinguished-markings}.

We record the case \(Y=0\), which compares the Weil--Petersson brackets
themselves with an explicit model and is the input for the coefficient
ratios in Section~\ref{uniform-one-boundary-coefficient-ratios}. At \(Y=0\), \(a=1\), \(r_0=\pi^2/2\), and \[
H_0(z)=\frac{3(\sin\sqrt{2z}-\sqrt{2z}\cos\sqrt{2z})}
 {(\sqrt{2z})^3}.
\] By \eqref{eq:H5} and \eqref{eq:H7}, \(\mathcal U=b/4^N\), where
\(b=[\tau_{d_1}\cdots\tau_{d_n}]_{g,n}\). Undoing this same scaling in
the model yields \begin{equation}
M^{\rm WP}_{g,n}(\mathbf d)
=4^D3^h\Gamma(h)4^{-g}(2\pi^2)^m
 [z^m]H_0(z)^{-h}.
\label{eq:H43}
\end{equation} Thus \eqref{eq:H42} applies unchanged to \(b/M^{\rm WP}\). For \(n=1\),
\(h=2g-1\), \(m=3g-2-k\), it gives exactly \begin{equation}
\mathrm{e}^{-C\log^2(g+1)}
\le
\frac{[\tau_k]_{g,1}}
 {4^k3^{2g-1}\Gamma(2g-1)4^{-g}(2\pi^2)^{3g-2-k}
 [z^{3g-2-k}]H_0(z)^{-(2g-1)}}
\le \mathrm{e}^{C\sqrt g},
\label{eq:H44}
\end{equation} uniformly over \(0\le k\le3g-2\).

\section{Ratios of marked generating functions}\label{ratios-with-a-fixed-number-of-distinguished-markings}

The harmonic comparison of Proposition~\ref{prop:harmonic} controls
logarithms, whereas the volume-ratio problem requires relative errors
tending to zero. This section begins the second idea of
Subsection~\ref{subsec:strategy} and obtains that stronger control for
two ratios of generating functions with a fixed number $p$ of
distinguished markings: the cusp ratio, corresponding to adding a cusp,
and the genus ratio, corresponding to replacing one unit of genus by two
cusps. The result is Proposition~\ref{prop:marked-ratios}: both ratios
have explicit leading terms with errors $O_p(1/g)$, uniformly in the
density parameter. The inputs are the kernel bounds \eqref{eq:P11}, the
positive recursion \eqref{eq:H9} and the $o(g)$ form of
Proposition~\ref{prop:harmonic}; all proofs are given here. The
passage to volumes with a prescribed cusp count is made in
Section~\ref{extracting-fixed-cusp-counts}.

The argument combines positivity with two exact identities. Positivity of
the recursion implies that a normalized
amplitude is nonincreasing when one of its degrees is raised~\eqref{eq:TL15}.
The dilaton equation~\eqref{eq:TL12} writes the cusp ratio as a weighted
sum, over degrees \(d\ge1\), of the normalized amplitudes
\(c_{g,p+1}(d,0^p)\). These weights can have either sign; their absolute
moments are uniformly bounded by~\eqref{eq:TL13}. The KdV equation~\eqref{eq:TL14} writes the
genus ratio through \(c_{g,p}(1,0^{p-1})\) and separating terms. Both
ratios therefore admit explicit leading terms with relative error
\(O_p(1/g)\), uniformly in \(Y\), once
we know that \(c_{g,p}(\mathbf d)=1-O(|\mathbf d|(|\mathbf d|+1)/g)\). This last estimate,
\eqref{eq:TL23}, follows from the difference of the recursion at adjacent
degrees; the coarse bounds of Section~\ref{uniform-comparison} are used
only to show that the separating terms are small.

\subsection{Shifted amplitudes and normalization}\label{shifted-amplitudes-and-normalization}

Use \(p\) for the number of distinguished markings, and keep \(n\) for the number of cusps. Retain the kernels and shifted amplitudes of Sections~\ref{the-positive-compact-family-of-kernels} and~\ref{exact-amplitudes-and-the-harmonic-model}; write \(U_{g,p}(\mathbf d;Y)=U^Y_{g,p}(\mathbf d)\), and suppress \(Y\) when convenient. Put \(h=2g-2+p\) and \(m=3g-3+p-|\mathbf d|\). Define \begin{equation}
\overline U_{g,p}(\mathbf d)=a^h r^m U_{g,p}(\mathbf d),
 \qquad W_{g,p}=4^g\overline U_{g,p}(\mathbf0),
 \qquad c_{g,p}(\mathbf d)=
 \frac{\overline U_{g,p}(\mathbf d)}{\overline U_{g,p}(\mathbf0)}.
 \label{eq:TL9}
\end{equation} Thus $\overline U_{g,p}(\mathbf d)$ is the normalized amplitude $\mathcal U_{g,p}(\mathbf d;Y)$ of \eqref{eq:H7}; the bar indicates that the markings are now the distinguished ones. The barred amplitudes obey the exact positive recursion \eqref{eq:H9}.

We use two consequences of
Sections~\ref{the-positive-compact-family-of-kernels}--\ref{uniform-comparison}.

\begin{itemize}
\item
  The kernel bounds \eqref{eq:P11}: for each fixed \(k\), \begin{equation}
1=\theta_0\le\theta_1\le\cdots,
   \qquad w_\ell:=\theta_\ell-\theta_{\ell-1}\ge0,
   \quad\theta_{-1}=0,
   \quad\sup_Y\sum_{\ell\ge0}(1+\ell)^k w_\ell<\infty.
   \label{eq:TL10}
\end{equation} Also \(A(Y)\) is bounded above and away from
  zero; see Subsection~\ref{compact-normalization-and-uniform-lattice-estimates}.
\item
  Uniform harmonic comparison: for each fixed positive \(p\), \begin{equation}
\log W_{g,p}=
   \log\Gamma(2g-2+p)-(2g-2+p)\log A(Y)+o(g),
   \label{eq:TL11}
\end{equation} uniformly for \(0\le Y\le y_0\). The corresponding crude bound
  \(W_{g,p}\le C_p^{g+1}(2g+p)!\) holds uniformly, including bounded
  stable genera. Proposition~\ref{prop:harmonic}, together with the
  central local limit \eqref{eq:H17} for the model, proves
  \eqref{eq:TL11} without any input on volume ratios. Only its \(o(g)\)
  error is needed below, and only for a positive number of markings.
\end{itemize}

\subsection{Dilaton and KdV identities}\label{exact-dilaton-and-genus-exchange-identities}

The shifted weighted dilaton identity is \begin{equation}
h\overline U_{g,p}(\mathbf0)=
 \sum_{d\ge1}b_d(Y)\overline U_{g,p+1}(d,0^p),
 \qquad
 b_d(Y)=\frac{X^{(d)}(Y)r^{d-1}}{a(2d+1)!!}
 =\frac13[z^{d-1}]\bar H_Y(z).
 \label{eq:TL12}
\end{equation} It holds also at \(p=0\), \(g\ge2\). To track the shift, the
background insertion of degree one contributes \(1-a\); moving it and
the higher background insertions to the left leaves the coefficients
\(X^{(d)}/(2d+1)!!\). The factor in \eqref{eq:TL9} then produces exactly
\(r^{d-1}/a\) in \eqref{eq:TL12}. The compact analytic family gives \begin{equation}
\sum_{d\ge1}b_d=A(Y),\qquad
 \sup_Y\sum_{d\ge1}(1+d)^k|b_d(Y)|<\infty
 \label{eq:TL13}
\end{equation} for every fixed \(k\).

The genus ratio is governed instead by the KdV equation. Write $F$ for the Kontsevich--Witten free energy, with subscripts indicating differentiation with respect to the corresponding times. The integrated KdV identity \(F_{01}=F_{00}^2/2+F_{0000}/12\),
differentiated \(p-2\) times in the string variable, gives for \(p\ge2\)
\begin{equation}
\overline U_{g,p}(1,0^{p-1})
 =\frac14\overline U_{g-1,p+2}(\mathbf0)
 +\frac32\sum_{j=0}^g\sum_{s=0}^{p-2}\binom{p-2}s
 \overline U_{j,s+2}(\mathbf0)
 \overline U_{g-j,p-s}(\mathbf0).
 \label{eq:TL14}
\end{equation} The split is ordered. The factor three from the weighted degree-one
insertion turns the original KdV coefficients into \(1/4\) and \(3/2\).
In \eqref{eq:TL14}, every term has the same total powers of \(a\) and \(r\): the
two separating complexities sum to \(h\), and their deficits sum to
\(3g-4+p\). Thus no omitted parameter factor is present.

The key consequence of positivity, used throughout this section, is the
following monotonicity.

\begin{lemma}[Monotonicity in the degrees]\label{lem:monotone}
For every stable \((g,p)\) with \(p\ge1\), every nonnegative degree tuple \(\mathbf d\), every \(1\le i\le p\) and
every \(Y\in[0,y_0]\),
\begin{equation}
\overline U_{g,p}(\mathbf d)\ge
 \overline U_{g,p}(\mathbf d+\mathbf e_i),
 \qquad 0\le c_{g,p}(\mathbf d)\le1.
 \label{eq:TL15}
\end{equation}
Amplitudes outside the dimension range are zero.
\end{lemma}

\begin{proof}
We subtract the recursion at adjacent degrees. With
$\mathbf d'=(d_2,\ldots,d_p)$ fixed, let
$\mathcal B_{\mathbf d'}(t)$ denote the merging and contraction sums
in \eqref{eq:H9}, with the kernel weight removed and $d_1+\ell=t$.
These sums are nonnegative and vanish for
$t>3g-3+p-|\mathbf d'|$ by dimension. If
$m=3g-3+p-|\mathbf d|\ge0$, shifting the kernel index gives
\[
 \begin{aligned}
 \overline U_{g,p}(\mathbf d)
   &=\sum_{\ell=0}^m\theta_\ell
                      \mathcal B_{\mathbf d'}(d_1+\ell),\\
 \overline U_{g,p}(\mathbf d+\mathbf e_1)
   &=\sum_{\ell=1}^m\theta_{\ell-1}
                      \mathcal B_{\mathbf d'}(d_1+\ell).
 \end{aligned}
\]
Their difference is therefore
$\sum_{\ell=0}^m w_\ell\mathcal B_{\mathbf d'}(d_1+\ell)\ge0$, which proves
\eqref{eq:TL15} for \(i=1\), and by symmetry for every \(i\). This includes
the dimension boundary by the zero convention. The difference identity is
used for $h\ge2$; the finite bases satisfy the inequality directly because
\(\theta_1\ge\theta_0=1\).
\end{proof}

From \eqref{eq:TL12}, \eqref{eq:TL13}, and \eqref{eq:TL15}, \begin{equation}
\frac{W_{g,p}}{W_{g,p+1}}\le\frac C{2g-2+p}
 \label{eq:TL16}
\end{equation} in all stable cases used below, with a constant independent of \(p\): the absolute sum in \eqref{eq:TL13} is uniform and the ratios in \eqref{eq:TL15} are at most one. Finite stable bases are included in the constant. Multiplying \eqref{eq:TL14} by \(4^g\) and using \eqref{eq:TL15} gives \begin{equation}
W_{g-1,p+2}\le W_{g,p},\qquad p\ge2.
 \label{eq:TL17}
\end{equation} The factor \(4^g\) in \eqref{eq:TL9} is essential: it converts the coefficient
\(1/4\) into the comparison constant one.

\subsection{Separating terms}\label{fixed-marking-splitting-estimates}

The separating terms in the KdV identity \eqref{eq:TL14}, and in the
degree difference above, are controlled by the convolution sums \(S_{g,p}\)
and \(T_{g,p}\) below. They are small for two reasons. First, splitting off a
component of bounded genus \(j\) removes \(2j+s\) units of complexity,
and the corresponding cusp insertions cost a factor \(O(1/g)\) each, by
\eqref{eq:TL16}--\eqref{eq:TL17}. Second, a split
into two components of comparable genus is exponentially small by the
binomial entropy in \eqref{eq:TL21}.

For fixed \(p\ge2\), define \[
 S_{g,p}=\frac1{W_{g,p}}
 \sum_{j=0}^g\sum_{s=0}^{p-1}\binom{p-1}s
 W_{j,s+1}W_{g-j,p-s},
\] \begin{equation}
T_{g,p}=\frac1{W_{g,p}}
 \sum_{j=0}^g\sum_{s=0}^{p-2}\binom{p-2}s
 W_{j,s+2}W_{g-j,p-s}.
 \label{eq:TL18}
\end{equation} We claim, uniformly in \(Y\), \begin{equation}
S_{g,p}=O_p(g^{-2}),\qquad T_{g,p}=O_p(g^{-1}),
 \label{eq:TL19}
\end{equation} with the positive-genus portion of \(T\) bounded by \(O_p(g^{-2})\).
In particular \(T_{g,2}=O(g^{-2})\) because its genus-zero terms are
unstable.

By symmetry orient the smaller genus as \(j\le g/2\). For
\(1\le j\le\varepsilon g\), insert \(2j\) cusps into the larger-genus
component using \eqref{eq:TL16}, apply \eqref{eq:TL17} \(j\) times, and insert the
remaining \(s\) cusps. This gives \begin{equation}
\frac{W_{g-j,p-s}}{W_{g,p}}
 \le (C_p/g)^{2j+s}.
 \label{eq:TL20}
\end{equation} If \(p-s=1\), instead insert \(2j+1\) cusps before raising the genus,
so every genus-comparison step ends with at least two markings;
afterward insert \(s-1\) cusps. This is possible because \(p\ge2\), and
the same exponent results. In the \(T\) convolution both relevant
marking counts are already at least two.

The crude bound following \eqref{eq:TL11} makes the smaller-component
contribution a fixed polynomial in \(j\) times \[
 (2j)!(C_p/g^2)^j g^{-s}.
\] Choose \(\varepsilon>0\) small enough that successive factorial terms
have a ratio bounded by a fixed number below one. Their positive-genus
sum is \(O_p(g^{-2})\), also with any fixed polynomial in \(j\). For
genus zero, stability gives \(s\ge2\) in \(S\) and \(s\ge1\) in \(T\);
direct use of \eqref{eq:TL16} gives the respective orders.

For \(\varepsilon g\le j\le g/2\), apply \eqref{eq:TL11} to both factors and the
parent. The marking counts belong to a fixed finite set. The powers of
\(A(Y)\) cancel up to bounded exponents, and Stirling's formula gives \begin{equation}
\log\frac{W_{j,s+c}W_{g-j,p-s}}{W_{g,p}}
 \le-2g\mathsf H(j/g)+o(g),\qquad c\in\{1,2\},
 \label{eq:TL21}
\end{equation} where \(\mathsf H(t)=-t\log t-(1-t)\log(1-t)\). Its positive minimum
on \([\varepsilon,1/2]\) makes the macroscopic sum exponentially small.
This proves \eqref{eq:TL19}. Every use of \eqref{eq:TL11} here has a positive fixed
marking count.

\subsection{Cusp and genus ratios}\label{uniform-fixed-marking-ratios}

For \(p\ge2\), the adjacent-degree difference obtained from \eqref{eq:H9} has a
merging prefactor \(O_p(1/g)\) by \eqref{eq:TL16}. Its nonseparating prefactor is
bounded by \begin{equation}
\frac{W_{g-1,p+1}}{W_{g,p}}
 \le\frac Cg\frac{W_{g-1,p+2}}{W_{g,p}}=O_p(g^{-1})
 \label{eq:TL22}
\end{equation} using \eqref{eq:TL16}--\eqref{eq:TL17}. To see all
three contributions, put $D=|\mathbf d|$,
$\Omega_0=\sum_{\ell\ge0}w_\ell$ and
$\Omega_1=\sum_{\ell\ge0}\ell w_\ell$. Bound each child amplitude
by its zero-degree value in the difference identity above. Since there
are at most $d_1+\ell+1$ new-degree pairs, this yields
\[
 \begin{aligned}
 c_{g,p}(\mathbf d)-c_{g,p}(\mathbf d+\mathbf e_1)
 \le{}&\Omega_0(2D+p-1)\frac{W_{g,p-1}}{W_{g,p}}\\
 &+\bigl((d_1+1)\Omega_0+\Omega_1\bigr)
   \left(2\frac{W_{g-1,p+1}}{W_{g,p}}+\frac12S_{g,p}\right).
 \end{aligned}
\]
The factor $2$ in the nonseparating term includes the factor $4$
from the genus normalization of $W$. The two kernel moments are
uniformly bounded by \eqref{eq:TL10}. Equations~\eqref{eq:TL16},
\eqref{eq:TL19} and~\eqref{eq:TL22} now give \[
 0\le c_{g,p}(\mathbf d)-c_{g,p}(\mathbf d+\mathbf e_i)
 \le\frac{C_p(|\mathbf d|+1)}g.
\] Summing degree increments proves \begin{equation}
0\le1-c_{g,p}(\mathbf d)
 \le C_p\frac{|\mathbf d|(|\mathbf d|+1)}g,
 \qquad p\ge2.
 \label{eq:TL23}
\end{equation} When the right side exceeds one the bound follows from \(0\le c_{g,p}\le1\) in \eqref{eq:TL15}. It also
bounds the terms outside the dimension range when they are set to zero.

Substitute \eqref{eq:TL23} into \eqref{eq:TL12}. The errors sum absolutely by \eqref{eq:TL13},
giving \begin{equation}
(2g-2+p)\frac{W_{g,p}}{W_{g,p+1}}
 =A(Y)+O_p(g^{-1}),\qquad p\ge1.
 \label{eq:TL24}
\end{equation} Since \(A\) is uniformly positive, the error is relative
\(O_p(1/g)\).

To recover the one-marking profile without circularity, first use \eqref{eq:TL16}
twice, \eqref{eq:TL17} at parent marking count two, and the now two-sided ratio
\eqref{eq:TL24} at marking count one: \begin{equation}
\frac{W_{g-1,2}}{W_{g,1}}
 \le\frac C{g^2}\frac{W_{g-1,4}}{W_{g,2}}
 \frac{W_{g,2}}{W_{g,1}}=O(g^{-1}).
 \label{eq:TL25}
\end{equation} The separating terms at one marking have two positive genera.
Equation \eqref{eq:TL16}, the inequality \(1/[j(g-j)]\le C/g\), and
\(T_{g,2}=O(g^{-2})\) give \begin{equation}
\frac1{W_{g,1}}\sum_{j=1}^{g-1}W_{j,1}W_{g-j,1}
 \le\frac Cg\frac{W_{g,2}}{W_{g,1}}T_{g,2}=O(g^{-2}).
 \label{eq:TL26}
\end{equation} The same degree-difference argument now proves \eqref{eq:TL23} for \(p=1\).
Thus \eqref{eq:TL24} holds also for \(p=0\).

Finally \eqref{eq:TL14}, \eqref{eq:TL19}, and \eqref{eq:TL23} yield \begin{equation}
\frac{W_{g-1,p+2}}{W_{g,p}}=1+O_p(g^{-1}),\qquad p\ge2.
 \label{eq:TL27}
\end{equation} For \(p=0,1\), add two markings in numerator and denominator and use
\eqref{eq:TL24}; the two cusp products have the same complexity factors and
cancel. Hence \eqref{eq:TL27} holds for every fixed \(p\ge0\).

We now return to volumes. Recall the cusp generating function
\(\Phi_g\) of \eqref{eq:intro-Phi} and the string-shift identity
\eqref{eq:H8}: \(U_{g,p}(0^p;Y)=\Phi_g^{(p)}(x)\) with \(x=X(Y)\), whenever
the amplitude is stable. By \eqref{eq:TL9},
\(W_{g,p}=4^ga^hr^m\Phi_g^{(p)}(x)\) with \(m=3g-3+p\). Adding a marking
raises both \(h\) and \(m\) by one, so
\[
 \frac{x\Phi_g^{(p+1)}(x)}{\Phi_g^{(p)}(x)}
 =\frac{x}{ar}\,\frac{W_{g,p+1}}{W_{g,p}}
 =\frac{hx}{A\,a\,r}\bigl(1+O_p(g^{-1})\bigr)
\]
by \eqref{eq:TL24}. Replacing \((g,p)\) by \((g-1,p+2)\) leaves \(h\)
unchanged and lowers \(m\) by one, so
\(\Phi_{g-1}^{(p+2)}/\Phi_g^{(p)}=4r\,W_{g-1,p+2}/W_{g,p}\), and
\eqref{eq:TL27} applies. This proves, for every fixed positive \(p\),
\begin{equation}
x\frac{(\Phi_g^{(p)})'(x)}{\Phi_g^{(p)}(x)}
 =h\Psi(x)+O_p(\Psi(x)),\qquad h=2g-2+p,
 \label{eq:TL30}
\end{equation} \begin{equation}
\frac{\Phi_{g-1}^{(p+2)}(x)}{\Phi_g^{(p)}(x)}
 =4r(Y)(1+O_p(g^{-1})),
 \qquad
 \Psi(x)=\frac{x}{A(Y)a(Y)r(Y)}.
 \label{eq:TL31}
\end{equation} The same genus ratio holds at \(p=0\). The estimates are uniform
throughout \(0\le Y<y_0\); at \(x=0\), \eqref{eq:TL30} is interpreted by
continuity.

The scalar defined in the last display is given exactly by the density parameter: \begin{equation}
H_Y(r)=\frac{3I(Y)}{2\Delta^{3/2}},\qquad
 A(Y)a(Y)r(Y)=\frac{I(Y)}{\sqrt\Delta},\qquad
 \Psi(x)=\lambda(Y)/2.
 \label{eq:TL32}
\end{equation} The first equality follows by substituting \(v=T-\Delta s^2\) in
\(I\), and the other two follow by direct substitution. We summarize the
output of this section used in the extraction arguments.

\begin{proposition}[Marked cusp and genus ratios]\label{prop:marked-ratios}
For every fixed \(p\ge1\), uniformly for \(0\le Y<y_0\) and \(x=X(Y)\), the
cusp ratio of \(\Phi_g^{(p)}\) satisfies \eqref{eq:TL30} and the genus ratio
satisfies \eqref{eq:TL31}, with \(\Psi(x)=\lambda(Y)/2\). The genus ratio
holds also at \(p=0\).
\end{proposition}

Proposition~\ref{prop:marked-ratios} is used in two ways. The
normalized forms \eqref{eq:TL24} and \eqref{eq:TL27} enter the
differential equation for the leading term in
Section~\ref{sec:prefactor}, and the ratios \eqref{eq:TL30} and
\eqref{eq:TL31} supply the center of the local limit theorem and the
ratio of generating functions in
Section~\ref{extracting-fixed-cusp-counts}. The degree estimate
\eqref{eq:TL23} is used again for bounded boundary lengths, in
\eqref{eq:boundary-factor} and in the proof of
Proposition~\ref{prop:bounded-boundaries}.

\section{The leading term and its normalization}\label{sec:prefactor}

The marked ratios of
Proposition~\ref{prop:marked-ratios} determine a differential equation
in $Y$ for $\log\Phi_g^{(p)}$ but leave a multiplicative constant
undetermined. The equation comes from the dilaton identity
\eqref{eq:TL12}, which expresses the \(Y\)-derivative of
\(\log\Phi_g^{(p)}\) through the degree profiles \(c_{g,p+1}(d,0^p)\); to
obtain it with error \(O(h^{-2})\) we need the \(1/h\) term of these
profiles (Lemma~\ref{lem:first-degree}). We then integrate the equation
and fix the constant at $Y=y_0$, where only the single sequence
$C_g(2^{3g-3})$ of pure $\psi$ intersection numbers enters and
Painlev\'e~I supplies its asymptotics \eqref{eq:PF12a}; the
quantitative form of that input is proved in
Appendix~\ref{app:marked}, Subsection~\ref{sec:endpoint-expansion}. The
result is Proposition~\ref{prop:prefactor}, whose fixed-cusp
specialization proves Corollary~\ref{cor:constant}. Matching at a point
that approaches the critical endpoint polynomially in the genus keeps
the integrated error uniform on the whole interval.

Write $\bar K=\bar K_Y$ and $\bar H=\bar H_Y$, and introduce the kernel moments
\begin{equation}
 \Theta=\sum_{\ell\ge0}w_\ell=-\frac1{\bar K'(1)},\qquad
 M_1=\sum_{\ell\ge0}\ell w_\ell=\frac{\bar K''(1)}{2\bar K'(1)^2},\qquad
 \mathcal Z=Aar^{3/2}=\sqrt2 I.
 \label{eq:PF1}
\end{equation}
Primes on $\bar K,\bar H$ mean differentiation in their argument; derivatives of
parameter functions below are with respect to $Y$. All the moments in
\eqref{eq:PF1} are uniformly bounded on the closed interval $0\le Y\le y_0$.
The variance $v_Y$ of the law in \eqref{eq:P12}, written $\sigma_Y^2$ in Section~\ref{the-positive-compact-family-of-kernels}, satisfies
\begin{equation}
 v_Y=\frac{15}{4}-\frac{\bar H''(1)}{\bar H(1)}=\frac1{2A\Theta}.
 \label{eq:PF2}
\end{equation}
Indeed $\bar H'(1)=-3\bar H(1)/2$ and
$\bar K'(1)=-5\bar H(1)/2+2\bar H''(1)/3$.

For fixed $p\ge1$, put $h=2g-2+p$ and define
\begin{equation}
 \mathcal J_{g,p}(Y)=
 \frac{\Gamma(h)4^{-g}}{\sqrt{\pi h}}
 \mathcal Z(Y)^{-h}r(Y)^{p/2}\sqrt{A(Y)\Theta(Y)}.
 \label{eq:PF3}
\end{equation}
\begin{proposition}[Absolute marked prefactor]\label{prop:prefactor}
For each fixed $p\ge1$, uniformly for $0\le Y<y_0$,
\begin{equation}
 \Phi_g^{(p)}(X(Y))=
 \frac{\mathcal J_{g,p}(Y)}{\pi}
 \left(1+O_p\!\left(\frac{\log(g+1)}g\right)\right).
 \label{eq:PF4}
\end{equation}
\end{proposition}

\subsection{The first degree correction}

We compute the first degree correction directly from the positive
recursion. The remainder must be uniform in $Y$ and in a growing degree
range, because an expansion at each fixed degree would not suffice in the
dilaton sum. Set $D=|\mathbf d|$ and $q=\#\{i:d_i>0\}$.

\begin{lemma}[First degree correction]\label{lem:first-degree}
For each fixed $p\ge1$, uniformly for $Y\in[0,y_0]$ and for degree tuples
with $D\le h^{1/4}$,
\begin{equation}
 \begin{aligned}
 c_{g,p}(\mathbf d)&=1-\frac{\delta_p(\mathbf d)}h
       +O_p\!\left(\frac{(D+1)^4}{h^2}\right),\\
 \delta_p(\mathbf d)&=A\left\{
 \Theta\bigl(D^2+(p-4)D\bigr)+2M_1D+
 \frac{q(q+5-2p)}2\right\}.
 \end{aligned}
 \label{eq:PF5}
\end{equation}
\end{lemma}

\begin{proof}
Subtract the exact recursion \eqref{eq:H9} at pivot degree $t=d_i$ and
$t+1$, shifting the second kernel index. The difference has coefficients
$w_\ell$. By \eqref{eq:TL24} and \eqref{eq:TL27}, its merging prefactor is
$A/h+O_p(h^{-2})$, and its nonseparating prefactor, including the factor
$1/2$ and the genus normalization, is $2A/h+O_p(h^{-2})$.
For $p=1$ there is no merging term. All child degrees can be replaced
by zero with error at most $C_p(D+\ell+1)^2/h$, by \eqref{eq:TL23}.
The coefficients and pair counts are $O_p(D+\ell+1)$, so the moments
of $w_\ell$ make the resulting error $O_p((D+1)^3/h^2)$.
The separating difference is $O_p((D+1)/h^2)$ by \eqref{eq:TL19};
at $p=1$ use \eqref{eq:TL26}. These bounds also cover dimension-invalid
children through the zero convention in \eqref{eq:TL23}.

Let $z$ count zeros among the other $p-1$ entries. A merger is admissible
when $t+d_j+\ell\ge1$; the nonseparating pair count is
$(t+\ell-1)_+$. Thus the leading adjacent difference is
\begin{equation}
 \frac Ah\left\{\Theta(2D+p-3)+2M_1+(2-z)\mathbf1_{t=0}\right\}.
 \label{eq:PF6}
\end{equation}
Here the zero-pivot term records both the missing zero-zero merger and
the truncated pair count. Since the parent deficit is of order $h$ in
the indicated degree range, extending the kernel sum to infinity costs
$O_p(h^{-2})$, using sufficiently high moments of $w_\ell$.
Summing increments from zero to $\mathbf d$ proves \eqref{eq:PF5}:
the exceptional increments at the creation of successive positive
entries are $3-p,4-p,\ldots,q+2-p$. The sum of the remainders is bounded
by $C_p(D+1)^4/h^2$.
\end{proof}

For $d>h^{1/4}$ in the dilaton identity, we use $0\le c\le1$, the
quadratic growth of $\delta_p$, and sufficiently high uniform absolute
moments of $b_d$; moments of order twelve already bound the
tail of $1+d^2/h$ by $O_p(h^{-2})$. This justifies summing the expansion,
although it gives no useful pointwise remainder at large degree.

\subsection{The differential equation for the leading term}

The three signed dilaton moments are exact:
\[
 \sum b_d=A,\qquad \sum d b_d=-A/2,\qquad
 \sum d^2 b_d=\bar H''(1)/3-7A/2.
\]
Apply \eqref{eq:PF5} to $(d,0^p)$ in \eqref{eq:TL12}. Changing the child
complexity $h+1$ to $h$ costs $O_p(h^{-2})$. We obtain
\begin{equation}
 h\frac{W_{g,p}}{W_{g,p+1}}
 =A\left(1-\frac{\beta_p}h+O_p(h^{-2})\right),\qquad
 \beta_p=A\left\{\frac{7-2p}{4}\Theta-M_1+2-p\right\}-\frac12.
 \label{eq:PF7}
\end{equation}
This holds for each fixed $p\ge0$: at $p=0$ its only profile input
has one marking. Inverting and undoing the normalization gives
\begin{equation}
 \frac{\mathrm{d}}{\mathrm{d}Y}\log\Phi_g^{(p)}(X(Y))
 =\frac{h+\beta_p+O_p(h^{-1})}{A r}.
 \label{eq:PF8}
\end{equation}

The correction in \eqref{eq:PF8} has an explicit primitive. Differentiate
$K_Y(r)=0$, using $\partial_YK=2\partial_zK+(K-a)/z$, to get
\begin{equation}
 r'=-(\Theta+2),\qquad
 (\log\mathcal Z)'=-\frac1{Ar}.
 \label{eq:PF9}
\end{equation}
The second identity also follows from $I'=-a\sqrt\Delta$.
Differentiating the root derivative and the corresponding equation
for $H$ gives
\begin{equation}
 \begin{split}
 (\log\Theta)'&=\frac{a'}a+
                  \frac{2\Theta+1-2M_1}{r},\\
 (\log A)'&=\frac{3\Theta/2+3-1/A}{r}-\frac{a'}a.
 \end{split}
 \label{eq:PF10}
\end{equation}
For the first formula use $\Theta=-a/(r\,\partial_zK_Y(r))$ and
$\partial_Y\partial_zK_Y=2\partial_z^2K_Y+\partial_zK_Y/r+a/r^2$ at the root $z=r$.
Equations \eqref{eq:PF7} and \eqref{eq:PF10} yield
\begin{equation}
 \frac{\beta_p}{Ar}
 =\frac12(\log(A\Theta))'+\frac p2(\log r)',\qquad
 \left|\frac{\mathrm{d}}{\mathrm{d}Y}\log
 \frac{\Phi_g^{(p)}(X(Y))}{\mathcal J_{g,p}(Y)}\right|
 \le\frac{C_p}{h A r}.
 \label{eq:PF11}
\end{equation}
Only exact identities have been differentiated here; the error term in
\eqref{eq:PF8} will be integrated, never differentiated.

\subsection{Fixing the constant at the critical endpoint}\label{sec:prefactor-endpoint}

Here we need the value at one pure-$\psi$ tuple, not only two-sided
bounds. In the normalization \eqref{eq:H4},
\begin{equation}
 C_g(2^{3g-3})=\frac1\pi+O(g^{-1}).
 \label{eq:PF12a}
\end{equation}
Indeed, the identity of Itzykson and Zuber~\cite[equation~(6.10)]{IZ}
expresses $C_g(2^{3g-3})$ through the coefficients of the formal
power-series solution of Painlev\'e~I.
Their large-order asymptotics follow from the Stokes data of Joshi--Kitaev
and Kapaev~\cite{JK,Kap}; see~\cite[Section~2.3]{GIKM}.
The constant itself is evaluated by Stirling's formula in
Lemma~\ref{rev:lem:endpoint-gevrey}(iii), whose leading case uses only
these two inputs. Subsection~\ref{sec:endpoint-expansion} proves a
quantitative version, \eqref{eq:FE10}, without using \eqref{eq:PF12a}. Estimate
\eqref{eq:PF12a} is also the case $\mathbf e=(2^{3g-3})$ of the uniform
theorem of Guo, Yang and Zagier~\cite[Theorem~1]{GYZ}; no other case of
that theorem is used.
For $m=3g-3+p$, string followed by dilaton gives the exact reduction
\begin{equation}
 C_g(0^p,2^m)=
 \prod_{i=0}^{p-1}
 \frac{5(m-i)}{3(2g-3+p+m-2i)}\ C_g(2^{3g-3})
 =\frac1\pi+O_p(g^{-1}).
 \label{eq:PF12}
\end{equation}
Each string step chooses a degree-two marking with normalized
coefficient $5/[3(2g-3+p+m)]$; the resulting degree-one marking is
removed by normalized dilaton with coefficient one. Since $p$ is fixed,
the displayed finite product is $1+O_p(g^{-1})$, and \eqref{eq:PF12a}
gives the final equality.

At the endpoint, \eqref{eq:H24} gives exactly
\begin{equation}
 \frac{\overline U_{g,p}(0^p;y_0)}
      {\mathcal M_{g,p}(0^p;y_0)}=C_g(0^p,2^m).
 \label{eq:PF13}
\end{equation}
Indeed the numerator is $3^m\langle\tau_0^p\tau_2^m\rangle_g/m!$
and the denominator is
$3^{h+m}\Gamma(h+m)/(4^g5^m m!)$.

We now show that the endpoint value persists on a short interval
$0\le y_0-Y\le h^{-2}$. The reason is that the amplitudes are polynomials
with nonnegative coefficients in the kernel coefficients, and these are
within a factor $1+O(y_0-Y)$ of their endpoint value one.
Here $\theta_\ell=[z^\ell]\bar K_Y(z)^{-1}$ are the rescaled kernel coefficients. Put $\delta=y_0-Y$. The normalized kernels are analytic in $\delta$
near zero on a common disk of radius greater than one. To see the
removable parameter singularity, apply the implicit-function theorem
to $K_{y_0-\delta}(\delta\rho)/\delta$ at $(\delta,\rho)=(0,3)$;
its $\rho$ derivative is nonzero. Division by the persistent root
$1-z$ then gives $\bar K_Y(z)/(1-z)=1+O(\delta)$ on a smaller common disk.
Together with coefficient positivity this proves
\begin{equation}
 1\le\theta_\ell(Y)\le1+C\delta\quad\hbox{for all }\ell\ge0.
 \label{eq:PF14}
\end{equation}
By the positive recursion, $\overline U$ is a polynomial with nonnegative
coefficients in the variables $\theta_\ell$, homogeneous of degree $h$. At complexity one the bases have degree one, with $\theta_0=1$;
each recursion adds one factor, while child complexities sum to $h-1$.
Thus
\begin{equation}
 1\le\frac{\overline U_{g,p}(0^p;Y)}
                {\overline U_{g,p}(0^p;y_0)}
 \le(1+C\delta)^h.
 \label{eq:PF15}
\end{equation}

Next, for every $Y$, we compare $\mathcal J_{g,p}$ with the harmonic
model at all degrees zero. Here $m=3g-3+p=3h/2-p/2$. Since
$\mathcal Z=Aar^{3/2}$, the powers of $r$ in $a^hr^m\mathcal J_{g,p}$ cancel
exactly. On the other hand, \eqref{eq:H12} and \eqref{eq:H21} give
$\mathcal M_{g,p}(0^p;Y)=3^h\Gamma(h)4^{-g}(3A)^{-h}\Pr(S_h=m)$, and the
central estimate \eqref{eq:H18} at displacement $-p/2$, with variance
$v_Y=1/(2A\Theta)$ from \eqref{eq:PF2}, gives
$\Pr(S_h=m)=\sqrt{A\Theta/(\pi h)}\,(1+O_p(h^{-1}))$. Hence
\begin{equation}
 a^h r^m\mathcal J_{g,p}(Y)
 =\Gamma(h)4^{-g}A^{-h}\sqrt{\frac{A\Theta}{\pi h}}
 =\mathcal M_{g,p}(0^p;Y)(1+O_p(h^{-1})),
 \label{eq:PF16}
\end{equation}
where the first equality is exact.
Its ratio to its endpoint value is $1+O_p(h\delta)$ when
$0\le\delta\le h^{-2}$, by analyticity and uniform positivity of
$A,\Theta$. Hence \eqref{eq:PF12}--\eqref{eq:PF16} prove
$\Phi_g^{(p)}(X(Y))/\mathcal J_{g,p}(Y)=1/\pi+O_p(h^{-1})$
throughout that interval.

\begin{proof}[Proof of Proposition~\ref{prop:prefactor}]
Take $Y_*=y_0-h^{-2}$. For $Y\le Y_*$ integrate \eqref{eq:PF11}:
by \eqref{eq:PF9}, its error is at most
\[
 \frac{C_p}{h}\log\frac{\mathcal Z(Y)}{\mathcal Z(Y_*)}
 \le C_p\frac{\log h}{h},
\]
since $\mathcal Z$ is decreasing, so that $\mathcal Z(Y)\le\mathcal Z(0)$, and
$\mathcal Z(Y_*)\asymp h^{-5}$. The value at $Y_*$ was fixed above.
Exponentiation gives \eqref{eq:PF4} for $Y\le Y_*$; the direct endpoint
comparison gives the stronger $O_p(h^{-1})$ error for $Y\ge Y_*$.
\end{proof}

At $Y=0$, $r=\pi^2/2$, $A=1/\pi^2$, $\Theta=2$, and
$\mathcal Z=\pi/(2\sqrt2)$. Taking $p=1$ in \eqref{eq:PF4} gives
\begin{equation}
 V_{g,1}=\frac{\Gamma(2g-1)(4\pi^2)^{2g-2}}{\sqrt{\pi g}}
 \left(1+O\!\left(\frac{\log(g+1)}g\right)\right).
 \label{eq:PF17}
\end{equation}
\begin{proof}[Proof of Corollary~\ref{cor:constant}]
Comparing \eqref{eq:PF17} with the expansion of
Mirzakhani--Zograf~\cite[Theorem~1.2]{MZ} at $n=1$
identifies their universal constant as $1/\sqrt\pi$. Since the constant
does not depend on $n$, their expansion then gives the asserted formula,
with relative error $O_n(g^{-1})$, for every fixed $n\ge0$.
\end{proof}

For completeness, we derive the boundary-length consequence stated in the
introduction. Fix $n\ge1$, put $x_i=L_i/2$, and set $S(x)=\sinh(x)/x$,
with $S(0)=1$. At $Y=0$, the normalization \eqref{eq:TL9} gives
$c_{g,n}(\mathbf d)=[\tau_{d_1}\cdots\tau_{d_n}]_{g,n}/V_{g,n}$.
Thus the volume polynomial, extended by zero outside its degree range, is
\[
 \frac{V_{g,n}(\mathbf L)}{V_{g,n}}
 =\sum_{\mathbf d\ge0}c_{g,n}(\mathbf d)
   \prod_{i=1}^n\frac{x_i^{2d_i}}{(2d_i+1)!}.
\]
Divide by $\prod_i S(x_i)$ and regard the product weights as the law of
independent nonnegative integers $d_i$. The identity
$\sum_{d\ge0}(2d)(2d+1)x^{2d}/(2d+1)!=x^2S(x)$ gives
$\mathbb E[d_i(d_i+1)]\le x_i^2/2$. Since
$D(D+1)\le n\sum_i d_i(d_i+1)$, where $D=\sum_i d_i$,
\eqref{eq:TL23}, including its one-marking case, yields
\begin{equation}
 0\le 1-\frac{V_{g,n}(\mathbf L)}{V_{g,n}\prod_i S(L_i/2)}
 \le \frac{C_n}{g}\sum_{i=1}^n L_i^2.
 \label{eq:boundary-factor}
\end{equation}
For $L_i=o(\sqrt g)$ the right side tends to zero. Corollary~\ref{cor:constant}
and Stirling's formula then give the asserted asymptotic formula.

The leading term $\mathcal J_{g,p}$ of Proposition~\ref{prop:prefactor}
is carried forward in two directions. Section~\ref{sec:higher-order}
refines the estimate \eqref{eq:PF4} to an expansion in $1/h$, and
Sections~\ref{extracting-fixed-cusp-counts}--\ref{sec:absolute} extract
from $\mathcal J_{g,1}$ the leading term $\widehat V_{g,n}$ of
Theorem~\ref{thm:intro-absolute}; the identities \eqref{eq:PF7},
\eqref{eq:PF9} and \eqref{eq:PF11} are used again there, in
\eqref{eq:CE1} and in the saddle-displacement estimate of
Subsection~\ref{subsec:quantitative-volume-ratios}.

\section{Higher-order corrections for marked generating functions}\label{sec:higher-order}

The estimate of Proposition~\ref{prop:prefactor} has relative error
$O_p(\log g/g)$. The logarithm comes from integrating the error
$O_p(1/(hAr))$ in~\eqref{eq:PF11} towards the critical endpoint, where
$Ar$ has a simple zero. This section states two refinements, still for a
fixed number $p$ of distinguished markings, obtained by iterating the
dilaton and KdV equations of Section~\ref{sec:prefactor} to higher
order: an explicit first correction
(Theorem~\ref{thm:first-correction}), which also gives the leading term
with relative error $O_p(g^{-1})$, and a complete expansion whose
truncation order may grow linearly with $g$
(Theorem~\ref{thm:full-expansion}). Both results are proved in
Appendix~\ref{app:marked}, with the endpoint constants again supplied by
Painlev\'e~I. They are the input of the extraction theorems of
Section~\ref{sec:sharp-extraction}, which make the passage to a
prescribed cusp count.

\subsection{The first correction}\label{subsec:first-correction}
The mechanism is that of Section~\ref{sec:prefactor}, carried one order
further. Lemma~\ref{lem:second-profile} extends the degree expansion
\eqref{eq:PF5} of the normalized profiles by a term $e_p(\mathbf d)/h^2$,
uniformly for $|\mathbf d|\le h^{1/8}$. Inserted into the dilaton identity,
it adds a term $\gamma_p/h$ to the numerator of the differential
equation~\eqref{eq:PF8}; see~\eqref{eq:ST8}. The coefficient $\gamma_p$,
defined in~\eqref{eq:ST7}, is a finite expression in kernel moments
(see~\eqref{eq:ST11}), and it vanishes at the pure-$\psi$ endpoint
(see~\eqref{eq:ST12}). Its primitive against $1/(Ar)$ is therefore bounded,
although $Ar$ vanishes there. The integration constant is fixed at the
endpoint by two inputs: the first correction to $C_g(2^{3g-3})$, which we
derive from Painlev\'e~I in Subsection~\ref{sec:endpoint-expansion}, and
the first correction to an explicit negative-binomial coefficient.
Accordingly, define
\begin{equation}
 t_p=\frac{-2p^2+12p-11}{20},\qquad
 \mathcal L_p(Y)=t_p-\int_Y^{y_0}\frac{\gamma_p(v)}{A(v)r(v)}\,\mathrm{d}v.
 \label{eq:ST13}
\end{equation}

\begin{theorem}[Uniform marked first correction]\label{thm:first-correction}
For every fixed $p\ge1$, uniformly for $0\le Y<y_0$,
\begin{equation}
 \Phi_g^{(p)}(X(Y))=\frac{\mathcal J_{g,p}(Y)}\pi
 \left\{1+\frac{\mathcal L_p(Y)}h+
            O_p\!\left(\frac{\log(g+1)}{g^2}\right)\right\}.
 \label{eq:ST14}
\end{equation}
In particular the leading relative error is $O_p(g^{-1})$.
\end{theorem}

At $Y=0$ the first-order coefficients of the cusp ratio, the genus
ratio and the degree profile that enter this expansion reduce to the
formulas of Mirzakhani and Zograf~\cite[Theorems~1.4 and~3.2(i)]{MZ};
this check is recorded at the end of Subsection~\ref{app:second-degree}.

\subsection{All orders}\label{subsec:all-orders}
The same scheme continues to every order. The recursion is triangular:
each new degree coefficient of the profiles uses earlier cusp and genus
coefficients, and the dilaton and KdV equations then give the new cusp and
genus coefficients. The difficulty is to control the constants as the
order grows.

\emph{Controlling the expansion order.} The aim is to lose only one
factorial per order. Polynomial degree estimates alone lead to two
factorials, because each recursion step both counts possible degree
splits and sums along a coordinate path.
Proposition~\ref{prop:gevrey-hierarchy} uses an exponential degree norm
instead. A small loss of degree exponent absorbs the path factor at a
cost proportional to the order, and in separating terms the dimension of
the bounded component controls its node degree independently of the
parent degree. Together these estimates give one factorial per expansion
order. A uniform Gamma envelope controls the omitted splittings through a
linear order range, so the same bound applies to the actual remainders
and not only to the formal coefficients. Endpoint transport and
coefficient extraction preserve this bound.

Appendix~\ref{app:marked} constructs, by the finite triangular rules
\eqref{eq:FE1}--\eqref{eq:FE6}, the formal coefficients of the profiles
and of the cusp and genus ratios, and from them functions $a_{p,j}(Y)$ and
\[
 \mathcal L_{p,j}(Y)=[t^j]\exp\!\Bigl(\sum_{i=1}^j a_{p,i}(Y)t^i\Bigr),
 \qquad \mathcal L_{p,1}=\mathcal L_p;
\]
see~\eqref{eq:FE11}. These functions are holomorphic on one complex
neighborhood of $[0,y_0]$ and satisfy
$\|a_{p,j}\|+\|\mathcal L_{p,j}\|\le C_p\mathfrak B_p^jj!$.

\begin{theorem}[Complete uniform marked expansion]\label{thm:full-expansion}
For each fixed $p\ge1$ there are $C_p,\mathfrak B_p,c_p,h_p>0$ such that,
for $h=2g-2+p\ge h_p$ and $0\le K\le c_ph$,
\begin{equation}
 \sup_{0\le Y<y_0}\left|
 \frac{\pi\Phi_g^{(p)}(X(Y))}{\mathcal J_{g,p}(Y)}
       -\sum_{j=0}^K\frac{\mathcal L_{p,j}(Y)}{h^j}\right|
 \le\frac{C_p\mathfrak B_p^{K+1}(K+1)!}{h^{K+1}},
 \qquad \mathcal L_{p,0}=1.
 \label{eq:FE12}
\end{equation}
The coefficients are specified by the finite rules and integrals of
Appendix~\ref{app:marked}, and have the Gevrey-1 norms stated above on
one complex neighborhood of $[0,y_0]$. The same error bound on the real parameter interval holds for $\log(\pi\Phi_g^{(p)}/\mathcal J_{g,p})$, with coefficients $a_{p,j}$ and zero constant term. A sufficiently small linear-order
truncation has error $O_p(\mathrm{e}^{-a_ph})$; no convergence of
the infinite series is required.
\end{theorem}

The proof has three steps. First, Proposition~\ref{prop:gevrey-hierarchy}
normalizes the profiles by $(p+Mj)^j\mathrm e^{\sigma D}$ and the
scalar ratios by $(p+Mj)^j$. After this normalization the
coefficient rules become coefficientwise inequalities between power
series with nonnegative coefficients, and a system of analytic majorants
(Lemma~\ref{lem:majorant}) bounds the normalized coefficients and the
actual remainders, uniformly in the order. Next, at the pure-$\psi$
endpoint, the integration constants are computed from the
Itzykson--Zuber identity and the large-order behavior of Painlev\'e~I;
the Borel representation of Costin, Costin and Huang~\cite{CCH} gives
these constants with a Gevrey-1 remainder. Finally, the logarithmic
derivative is integrated from a matching point at distance $h^{-K-2}$
from the endpoint, where the positive homogeneous structure of the
amplitudes controls the change from the endpoint value.

At $Y=0$ this is the fixed-cusp expansion mechanism of~\cite{MZ}, for
which effective order-dependent bounds are given by
Hide--Macera--Thomas~\cite[Theorem~4.2]{HMT} and single-factorial bounds
by Guillermit~\cite[Theorems~2 and~4]{Guillermit}; the comparison with
these bounds was made in the survey at the beginning of this part. What
is used below is the uniformity in the density parameter through the
pure-$\psi$ endpoint: Theorem~\ref{thm:full-expansion} is applied at a
fixed order in the proof of Theorem~\ref{thm:sharp-extraction} and over
the whole linear order range in that of
Theorem~\ref{thm:canonical-expansion}, and the linear-order truncation
in Theorem~\ref{thm:intro-absolute} is inherited from \eqref{eq:FE12}.

\section{Extracting a prescribed cusp count}\label{extracting-fixed-cusp-counts}

A generating-function estimate must be localized to recover a volume
with a prescribed cusp count. Positivity in $Y$
(Lemma~\ref{lem:kernel-positivity}) yields the exact convolution
factorization below and its Fourier bound \eqref{eq:TL35}. Combined
with the marked ratios of Proposition~\ref{prop:marked-ratios} and
the endpoint expansions of Lemma~\ref{lem:parameter}, it gives
Proposition~\ref{prop:local-limit}, uniformly across small, interior
and near-critical cusp weights.

The volume-ratio limits and bounded-boundary factorization that follow
use only the leading estimates of
Sections~\ref{ratios-with-a-fixed-number-of-distinguished-markings}--\ref{sec:prefactor}.
Higher-order marked expansions enter only when we seek inverse-genus
errors, in Sections~\ref{sec:sharp-extraction} and~\ref{sec:absolute}.

\subsection{A positive convolution factor}\label{a-positive-convolution-factor}

\begin{lemma}[Positive convolution factor]\label{lem:convolution}
For each stable marked pair \((g,p)\) with \(p\ge1\), the marked generating
function has the exact factorization
\begin{equation}
\Phi_g^{(p)}(x)=V(x)^h R_{g,p}(x),\qquad
 V(x)=Y'(x)=1/a(Y(x)),\qquad h=2g-2+p,
 \label{eq:TL33}
\end{equation} where \(R_{g,p}\) has nonnegative Taylor coefficients and
converges for \(0\le x<x_c\).
\end{lemma}

A relative approximation of \(\Phi_g^{(p)}\) on the real axis does not,
in general, give a relative approximation of its individual
coefficients. The factorization \eqref{eq:TL33} supplies the missing
structure: under the cusp weight \(x\), the number of additional cusps is
distributed as the sum of \(h\) independent variables with generating
function \(V(xz)/V(x)\) and one further independent nonnegative integer
variable. This yields the Fourier bound \eqref{eq:TL35} below, which is
the input a local limit theorem needs.

\begin{proof}
The proof uses positivity in the parameter: by
Lemma~\ref{lem:kernel-positivity}, the functions
\(\widehat\vartheta_\ell(Y)=a(Y)[z^\ell]K_Y(z)^{-1}\) and \(1/a(Y)\) have
nonnegative Taylor coefficients in \(Y\). Since this positivity in \(Y\) is
available for the unrescaled coefficients, we do not use the rescaling by
\(r(Y)\) of \eqref{eq:H7} here.

Undo the normalization in \eqref{eq:H9} and regard the unrescaled coefficients \(\vartheta_\ell=[z^\ell]K_Y(z)^{-1}\) as independent variables. In this convention the bases are \(U_{0,3}=\vartheta_0\), \(U_{1,1}(1)=\vartheta_0/8\), and \(U_{1,1}(0)=\vartheta_1/8\). Induction from these bases shows that every amplitude is a finite polynomial with nonnegative
coefficients in the variables \(\vartheta_\ell\); each monomial has total degree \(h\),
and weighted degree \(m\) when \(\vartheta_\ell\) is given weight \(\ell\). A
linear recursion term adds one factor \(\vartheta_\ell\) to a child of complexity \(h-1\); a
separating term adds one to children whose complexities sum to \(h-1\).
Thus factoring \(a^{-1}\) from each of the \(h\) factors \(\vartheta_\ell=a^{-1}\widehat\vartheta_\ell\) leaves a
positive polynomial in the \(\widehat\vartheta_\ell\).

Next, \(Y(x)\) has nonnegative coefficients and \(Y(0)=0\). To justify
the product used here, note that a nonzero zero \(y\) of \(X\) gives
\(u(s)=J_1(2\sqrt y\,s)\), with \(u(1)=0\) and regular behavior at zero.
Integration of \((s u')'+(4ys-1/s)u=0\) against its conjugate gives
\(4y\int_0^1s|u|^2\,\mathrm{d}s=\int_0^1(s|u'|^2+|u|^2/s)\,\mathrm{d}s>0\). Thus all
nonzero zeros are positive real. Since \(X\) has order one half and a
simple zero at zero with derivative one, its Hadamard product is
\(X(Y)=Y\prod_k(1-Y/\beta_k)\), \(\beta_k>0\). Lagrange inversion now
gives \[
 [x^k]Y(x)=\frac1k[Y^{k-1}](Y/X(Y))^k>0.
\] The inverse series has radius exactly \(x_c\): a smaller radius
would, by positivity and Pringsheim's theorem, give a positive real
singularity before \(x_c\), contradicting the analytic real inverse
there; the quadratic critical point at \(x_c\) supplies the opposite
inequality. Since
\(K_Y(z)/a(Y)=1+\sum_{\ell\ge1}X^{(\ell+1)}(Y)z^\ell/(a(Y)(2\ell+1)!!)\),
each \(\widehat\vartheta_\ell\) is a polynomial in the ratios
\(X^{(i)}(Y)/a(Y)\); by Lemma~\ref{lem:kernel-zeros} these ratios and
\(1/a\) are analytic on \(|Y|<y_0\), so their nonnegative Taylor series
converge on \([0,y_0)\). Composition therefore preserves coefficient positivity and
converges for \(0\le x<x_c\), since \(Y(x)<y_0\). Finally, by
\eqref{eq:H8}, \(\Phi_g^{(p)}(x)=U^{Y}_{g,p}(0^p)\) with \(Y=Y(x)\), and the
factor \(a(Y)^{-h}\) removed above equals \(V(x)^h\). This proves
\eqref{eq:TL33}.
\end{proof}

For the coefficient probability law with generating function
\(\Phi_g^{(p)}(xz)/\Phi_g^{(p)}(x)\), positivity in \eqref{eq:TL33} gives \begin{equation}
\left|\frac{\Phi_g^{(p)}(x\mathrm{e}^{i\theta})}{\Phi_g^{(p)}(x)}\right|
 \le\left|\frac{V(x\mathrm{e}^{i\theta})}{V(x)}\right|^h.
 \label{eq:TL35}
\end{equation} This Fourier estimate prevents a ratio of generating
functions on the real axis from concealing different local coefficient
behavior.

\subsection{Mean, variance and a central limit theorem}\label{mean-variance-and-a-real-central-limit-theorem}

We first identify the leading center and variance scale of the cusp count
under the weight \(x\) and prove a Gaussian limit when both the complexity
and the model mean tend to infinity. The next subsection upgrades this to
a local limit theorem. Recall
\(\mathscr D=x\,\mathrm{d}/\mathrm{d}x\) and the function \(\Psi\) of
Lemma~\ref{lem:parameter}. For convenience we repeat its endpoint expansions
\eqref{eq:SP16}, together with that of \(V=1/a(Y(x))\), which follows from
\(a(Y)=1-Y+O(Y^2)\) and \(Y(x)=x+O(x^2)\): \begin{equation}
\begin{aligned}
 x\downarrow0:\quad& V(x)=1+x+O(x^2),\quad \Psi(x)=2x+O(x^2),\\
 &\mathscr D\Psi(x)=2x+O(x^2),\quad \mathscr D^2\Psi(x)=2x+O(x^2);\\[1mm]
 x\uparrow x_c:\quad&\Psi(x)=\frac5{4\varepsilon}(1+O(\sqrt\varepsilon)),\\
 &\mathscr D\Psi(x)=\frac5{4\varepsilon^2}(1+O(\sqrt\varepsilon)),\quad
 \mathscr D^2\Psi(x)=O(\varepsilon^{-3}),
 \qquad\varepsilon=1-x/x_c.
 \end{aligned}
 \label{eq:TL36}
\end{equation} On compact subintervals of \((0,x_c)\), all needed
derivatives are bounded and \(\mathscr D\Psi\) is bounded away from
zero.

For \(F=\Phi_g^{(p)}\) at a positive cusp weight \(x\), let \begin{equation}
\mathbb P_x(N=k)=\frac{[z^k]F(z)x^k}{F(x)},\qquad
 \mu=h\Psi(x),\qquad \sigma^2=h\mathscr D\Psi(x).
 \label{eq:TL37}
\end{equation} The last two quantities are a model center and a model
variance; they are not claimed to be the exact finite-genus mean and
variance. Suppose \(h\to\infty\) and
\(\mu\to\infty\). For every fixed real \(u\), integrating \eqref{eq:TL30} gives
\begin{equation}
\log\mathbb E_x \mathrm{e}^{u(N-\mu)/\sigma}
 =\int_0^{u/\sigma}
 \{h\Psi(x\mathrm{e}^s)-h\Psi(x)+O_p(\Psi(x\mathrm{e}^s))\}\,\mathrm{d}s.
 \label{eq:TL38}
\end{equation} The interval of integration stays inside \((0,x_c)\). Expanding
\(h\Psi(x\mathrm e^s)-h\Psi(x)=sh\mathscr D\Psi(x)+O(s^2h\mathscr D^2\Psi)\),
the first part of the integral is \(u^2/2\) plus a Taylor remainder, and the
second part is an integrated error. Both errors are \(o(1)\). In the
critical regime and on an interior compact interval, the Taylor remainder
and the integrated error are \(O_u(h^{-1/2})\). At \(x\to0\), where
\(\sigma^2\sim\mu\), the Taylor remainder is \(O_u(\mu^{-1/2})\) and the
integrated error is \(O_u(x/\sqrt\mu)=O_u(\sqrt\mu/h)=o(1)\). These
assertions follow directly from \eqref{eq:TL36}; no estimate for the
derivative of the error in \eqref{eq:TL30} is required.

Thus the real moment generating functions converge to \(\exp(u^2/2)\)
along every such sequence. Exponential Markov bounds give tightness;
bounds at slightly larger real arguments give uniform integrability
along weakly convergent subsequences, whose limiting moment generating
functions then identify the standard Gaussian. Therefore \begin{equation}
(N-\mu)/\sigma\Longrightarrow\mathcal N(0,1).
 \label{eq:TL39}
\end{equation}

\subsection{Fourier bounds and the local limit theorem}\label{uniform-fourier-bounds-and-the-lattice-local-limit}

A central limit theorem alone does not determine ratios of individual
coefficients; for these we prove the corresponding local limit theorem.
Central and local limit theorems for coefficients of generating
functions go back to Bender~\cite{Bender}; see Hwang~\cite{Hwang} for
rates in the quasi-powers setting and~\cite[Chapter~IX]{FS} for an
account. Here the base function, the power and the evaluation point vary
together, up to both ends of the parameter interval, so the bounds are
proved directly.
The proof dominates the characteristic function by an integrable function
separately at small cusp weight, at interior cusp weights, and near the
critical point.

\begin{proposition}[Uniform local limit theorem]\label{prop:local-limit}
Fix \(p\ge1\) and use the notation \eqref{eq:TL37}. Along every sequence
with \(h\to\infty\), \(\mu\to\infty\) and \(x\in(0,x_c)\),
\begin{equation}
\mathbb P_x(N=k)=\frac{1+o(1)}{\sqrt{2\pi}\sigma}
 \mathrm{e}^{-u^2/2}\qquad\text{whenever }(k-\mu)/\sigma\longrightarrow u.
 \label{eq:TL47}
\end{equation}
No relation between \(h\) and the distance from \(x\) to \(x_c\) is required.
\end{proposition}

\begin{proof}
Passing to a subsequence, we may assume that \(x\) tends to zero, to an
interior point, or to \(x_c\). In each case we bound the characteristic
function of \((N-\mu)/\sigma\) by an integrable function, uniformly along
the sequence; Fourier inversion and the pointwise convergence implied by
\eqref{eq:TL39} then give \eqref{eq:TL47}.

On a compact subinterval of \((0,x_c)\), the coefficient probabilities
at zero and one of \(V(xz)/V(x)\) have positive uniform lower bounds. If
these probabilities are \(p_0,p_1\), then \begin{equation}
\left|\sum_jp_j\mathrm{e}^{ij\theta}\right|^2
 \le1-2p_0p_1(1-\cos\theta).
 \label{eq:TL40}
\end{equation} Equations \eqref{eq:TL35} and \eqref{eq:TL40} give an upper bound
\(\exp(-c h\theta^2)\) on \([-\pi,\pi]\). Since \(\sigma\) is comparable
to \(\sqrt h\), this is integrable Gaussian domination on the rescaled
Fourier axis.

At \(x\to0\), one has \(p_0\asymp1\), \(p_1\asymp x\). The same argument
gives \(\exp(-c h x\theta^2)\). Since \(hx\asymp\mu\) and
\(\sigma^2\sim\mu\to\infty\), it again gives integrable domination after
rescaling.

Near the critical point a different bound is needed, because the law
with generating function \(V(xz)/V(x)\) then has variance of order
\(\varepsilon^{-2}\) and its probabilities at zero and one tend to zero. Quadratic
analytic inversion at \(x_c\) gives \begin{equation}
V(z)=A_0(1-z/x_c)^{-1/2}
 \Xi\bigl(\sqrt{1-z/x_c}\bigr),\qquad
 A_0>0,\quad \Xi(0)=1,
 \label{eq:TL41}
\end{equation} where \(\Xi\) is analytic and nonzero near zero and has real
coefficients. Put \(\Lambda(w)=\log \Xi(\sqrt w)\) and
\(w(\theta)=1-(1-\varepsilon)\mathrm{e}^{i\theta}\). Then \begin{equation}
\log\left|\frac{V(x\mathrm{e}^{i\theta})}{V(x)}\right|
 =-\frac14\log\left(1+
 \frac{2(1-\varepsilon)(1-\cos\theta)}{\varepsilon^2}\right)
 +\Re \Lambda(w(\theta))-\Lambda(\varepsilon).
 \label{eq:TL42}
\end{equation} There are fixed \(\delta,\varepsilon_0>0\) such that \begin{equation}
\left|\frac{V(x\mathrm{e}^{i\theta})}{V(x)}\right|
 \le(1+\theta^2/(4\varepsilon^2))^{-1/8}
 \quad(|\theta|\le\delta,\ 0<\varepsilon\le\varepsilon_0).
 \label{eq:TL43}
\end{equation} Indeed, for \(|\theta|\le\varepsilon\), the first derivative at zero
of the real part of \(\Lambda(w(\theta))\) vanishes, while its second
derivative is \(O(\varepsilon^{-3/2})\). Its change is therefore
\(O(\theta^2\varepsilon^{-3/2})\), absorbed by half of the leading
logarithm in \eqref{eq:TL42} when \(\varepsilon\) is small. For
\(\varepsilon\le|\theta|\le\delta\), its change is bounded by
\(C\sqrt\delta\), while the logarithm is at least \(\log(5/4)\). Choose
\(\delta\) small enough for the same absorption. We used the elementary
bound \(2(1-\varepsilon)(1-\cos\theta)\ge\theta^2/4\) on the small arc.

On the outer arc \(\delta\le|\theta|\le\pi\), \begin{equation}
|V(x\mathrm{e}^{i\theta})/V(x)|\le C_\delta\sqrt\varepsilon.
 \label{eq:TL44}
\end{equation} To see that \(V\) is uniformly bounded there, note that nonnegative coefficients and
\(Y(x)\to y_0\) give absolute convergence of \(Y\) on \(|x|\le x_c\).
Its first two coefficients are positive, so strict triangle inequality
gives \(|Y(x_c \mathrm{e}^{i\theta})|<y_0\) away from \(\theta=0\), while
\(|Y(z)|\le Y(|z|)<y_0\) for \(|z|<x_c\). By continuity of \(Y\) on the
closed disk, \(|Y|\le y_1\) for some \(y_1<y_0\) on the closed sector
\(|z|\le x_c\), \(\delta\le|\arg z|\le\pi\). By Lemma~\ref{lem:kernel-zeros}, \(X'\) has no zero in \(|Y|<y_0\), so
\(V=1/a(Y)\) is uniformly bounded on this sector.
Since \eqref{eq:TL41} gives \(V(x)\asymp\varepsilon^{-1/2}\), \eqref{eq:TL44} follows.

On the critical Fourier scale \(\theta=\varepsilon t/\sqrt h\), \eqref{eq:TL35}
and \eqref{eq:TL43} are bounded by \begin{equation}
(1+t^2/(4h))^{-h/8}\le(1+t^2/64)^{-2},\qquad h\ge16.
 \label{eq:TL45}
\end{equation} The inequality follows because \((1+a/h)^h\) increases with \(h\).
The right side is integrable. The outer contribution after
multiplication by \(\sqrt h/\varepsilon\) is at most \begin{equation}
2\pi\sqrt h\,C_\delta^h\varepsilon^{h/2-1}\longrightarrow0.
 \label{eq:TL46}
\end{equation} This convergence is uniform for arbitrarily small \(\varepsilon\):
fix sufficiently small \(\varepsilon_0\) with
\(C_\delta\sqrt{\varepsilon_0}<1\), and bound \eqref{eq:TL46} by
\(2\pi\varepsilon_0^{-1}\sqrt h(C_\delta\sqrt{\varepsilon_0})^h\) for
\(h\ge2\).

These integrable bounds, in the three regimes, complete the proof.
\end{proof}

The same Fourier bounds also bound every point mass, not only those near
the center. We record this for later use in
Appendices~\ref{app:short-curves}, \ref{app:tight}, \ref{app:super-cuts}
and~\ref{app:sublinear-trace}.

\begin{lemma}[Uniform point-mass bound]\label{lem:point-mass}
For every $c_0>0$ there is a constant $C$ with the following property.
Let $x\in(0,x_c)$, let $h,h'$ be positive integers with $h'\ge\max(16,c_0h)$,
and let $G(z)=V(z)^{h'}R(z)$, where $R$ has nonnegative Taylor
coefficients, converges for $|z|<x_c$ and satisfies $R(x)>0$. Then,
with $\sigma^2=h\mathscr D\Psi(x)$,
\[
 \sup_{j\ge0}\frac{[z^j]G(z)\,x^j}{G(x)}\le\frac C\sigma.
\]
\end{lemma}

\begin{proof}
By Fourier inversion and positivity of the coefficients of $R$,
\[
 \frac{[z^j]G(z)\,x^j}{G(x)}
 \le\frac1{2\pi}\int_{-\pi}^{\pi}
       \left|\frac{V(x\mathrm e^{i\theta})}{V(x)}\right|^{h'}\mathrm d\theta.
\]
On a compact subinterval of $(0,x_c)$, \eqref{eq:TL40} bounds the
integrand by $\mathrm e^{-ch'\theta^2}$, so the integral is
$O(h^{-1/2})$, while $\sigma\asymp\sqrt h$. Near zero the same argument
with $p_1\asymp x$ gives $\mathrm e^{-ch'x\theta^2}$ and the bound
$O((hx)^{-1/2})$, while $\sigma^2\asymp hx$. Near $x_c$, put
$\theta=\varepsilon t/\sqrt{h'}$ on $|\theta|\le\delta$; \eqref{eq:TL43}
and the argument of \eqref{eq:TL45} bound this part by
$(\varepsilon/\sqrt{h'})\int(1+t^2/64)^{-2}\,\mathrm dt$. On the outer arc
\eqref{eq:TL44} gives at most
$2\pi(C_\delta\sqrt\varepsilon)^{h'}\le C\varepsilon h^{-1/2}$ once
$C_\delta\sqrt{\varepsilon_0}<1$. Since $\sigma\asymp\sqrt h/\varepsilon$
there, all three regimes give $O(1/\sigma)$.
\end{proof}

\subsection{Ratios at a prescribed number of cusps}\label{canonical-coefficient-ratios-and-all-cusp-count-uniformity}

We now prove the limits in Theorem~\ref{thm:volume}. Each volume ratio is
a ratio of coefficients, at the same index, of two marked generating
functions. We evaluate both at the same cusp weight, chosen so that the
index is, up to a negligible shift, the model center of both. The local limit theorem \eqref{eq:TL47} shows
that the two local probabilities are asymptotically equal, so the
coefficient ratio is asymptotic to the ratio of the generating functions
themselves, which Section~\ref{ratios-with-a-fixed-number-of-distinguished-markings}
determines.

\begin{proof}[Uniform convergence in Theorem~\ref{thm:volume}]
First suppose \(g,n\to\infty\). Put \begin{equation}
k=n-1,\qquad h=2g-1,
 \qquad h\Psi(x)=k.
 \label{eq:TL48}
\end{equation} The last equation has a unique solution in \((0,x_c)\) by \eqref{eq:TL32} and
the parameter bijection. The two families \(\Phi_g'\) and
\(\Phi_{g-1}'''\) have exactly the same complexity \(h\), model center,
and model variance. Thus their two local probabilities in \eqref{eq:TL47}, at the
common coefficient index \(k\), have ratio tending to one. Equations
\eqref{eq:TL31} and \eqref{eq:TL47} give \begin{equation}
\frac{[x^k]\Phi_{g-1}'''(x)}{[x^k]\Phi_g'(x)}
 =4r(Y)(1+o(1)).
 \label{eq:TL49}
\end{equation} By the definition \eqref{eq:intro-Phi}, the coefficients are \begin{equation}
[x^k]\Phi_g'(x)=
 \frac{V_{g,n}}{(2\pi^2)^{3g-3+n}k!},\qquad
 [x^k]\Phi_{g-1}'''(x)=
 \frac{V_{g-1,n+2}}{(2\pi^2)^{3g-4+n}k!}.
 \label{eq:TL50}
\end{equation} Consequently \begin{equation}
B_{g,n}=\frac{\pi^2}{2r(Y)}(1+o(1))
 =\mathcal B\left(\frac{2(n-1)}{2g-1}\right)(1+o(1)).
 \label{eq:TL51}
\end{equation} The factorial \(k!\) cancels, and the powers of \(2\pi^2\) differ by
exactly one; this explains the factor in \eqref{eq:TL51}.

For the cusp ratio, compare \(\Phi_g''\) with \(\Phi_g'\) at the same
index \(k\) and the same cusp weight \eqref{eq:TL48}. Their complexities are \(h+1\)
and \(h\). Their centers differ by \(\Psi(x)\), and their variances have
ratio \((h+1)/h\). Moreover \begin{equation}
\frac{\Psi(x)}{\sqrt{h\mathscr D\Psi(x)}}=O(h^{-1/2})
 \label{eq:TL52}
\end{equation} uniformly in \(x\), because \(\Psi^2/(\mathscr D\Psi)\) is continuous
in the interior and has finite endpoint limits zero and \(5/4\),
respectively. Thus the local probability ratio again tends to one.
Undoing \eqref{eq:TL24}, \begin{equation}
C_{g,n}=\frac{\pi^2}{g}
 \frac{[x^k]\Phi_g''}{[x^k]\Phi_g'}
 =\frac{\pi^2}{g}\frac{h}{Aar}(1+o(1))
 =\frac{\pi^2(n-1)}{gx}(1+o(1)).
 \label{eq:TL53}
\end{equation} Here, as in \eqref{eq:TL50}, the coefficients \([x^k]\Phi_g''\) and
\([x^k]\Phi_g'\) carry the normalizing powers \((2\pi^2)^{3g-2+n}\) and
\((2\pi^2)^{3g-3+n}\); division by \(2g\) gives the factor \(\pi^2/g\). With
\(\lambda_{\mathrm{eff}}=2(n-1)/(2g-1)\), the last expression equals \begin{equation}
\left(1-\frac1{2g}\right)\mathcal C(\lambda_{\rm eff})(1+o(1)).
 \label{eq:TL54}
\end{equation}

One may replace \(\lambda_{\mathrm{eff}}\) by \(n/g\) in both \eqref{eq:TL51} and
\eqref{eq:TL54}. At a zero limiting slope, use the continuous values
\(\mathcal B(0)=1\) and \(\mathcal C(0)=4\pi^2\) of \eqref{eq:SP3}; at a
finite positive slope use continuity; and at an infinite slope use the
endpoint equivalents \eqref{eq:SP15} of Lemma~\ref{lem:parameter}.
Since \(\lambda_{\mathrm{eff}}/(n/g)\to1\) whenever \(n\to\infty\),
these give the desired relative replacement in all cases.

Finally the exact identity \begin{equation}
D_{g,n}=C_{g,n-1}/C_{g,n}
 \label{eq:TL56}
\end{equation} and the adjacent continuity of the profile \(\mathcal C\) give
\(D_{g,n}\to1\) for every growing-cusp sequence. For bounded cusp
counts, the classical fixed-cusp theorems give \(B_{g,n}\to1\),
\(C_{g,n}\to4\pi^2\), and hence \(D_{g,n}\to1\); these are the only
bounded-cusp inputs used here.

If any supremum in Theorem~\ref{thm:volume} did not tend to zero, select a violating
sequence of cusp counts. A subsequence either has a fixed bounded cusp
count or tends to infinity. The preceding conclusions contradict a fixed
violation in both cases. This proves the uniform convergence assertions.
Subsection~\ref{subsec:quantitative-volume-ratios} supplies the stated
inverse-genus rate using Theorem~\ref{thm:sharp-extraction}.

\end{proof}

\subsection{Bounded geodesic boundaries}\label{subsec:bounded-boundaries}

The positive convolution factor also allows a fixed number of cusps to be
replaced by geodesic boundaries of bounded length. Each boundary then
contributes a factor that depends on the cusp density and connects the
familiar large-genus and large-cusp boundary profiles.
For $0\le Y\le y_0$ and $L\ge0$, put
\begin{equation}
 \mathfrak b_Y(L)=\int_0^1
 I_0\!\left(\frac L\pi\sqrt{T(Y)-\Delta(Y)s^2}\right)\,\mathrm ds,
 \label{eq:boundary-profile}
\end{equation}
where $I_0$ is the modified Bessel function of order zero, $T=Y+\Delta$
and $r=2\Delta$. At $Y=y_0$ use $T=y_0$ and $\Delta=0$.
Write $Y_\lambda$ for the inverse of $\lambda(Y)$; thus
$Y_0=0$ and $Y_\infty=y_0$.

\begin{proposition}[Boundary lengths at arbitrary cusp density]
\label{prop:bounded-boundaries}
Fix an integer $p\ge1$ and $L_0<\infty$. As $g\to\infty$,
\begin{equation}
 \sup_{\substack{n\ge p\\\mathbf L\in[0,L_0]^p}}
 \left|
 \frac{V_{g,n}(L_1,\ldots,L_p,0^{n-p})}
 {V_{g,n}\prod_{i=1}^p\mathfrak b_{Y_{n/g}}(L_i)}-1
 \right|\longrightarrow0.
 \label{eq:boundary-uniform}
\end{equation}
The endpoint factors are
\begin{equation}
 \mathfrak b_0(L)=\frac{\sinh(L/2)}{L/2},\qquad
 \mathfrak b_{y_0}(L)=I_0\!\left(\frac{j_{0,1}L}{2\pi}\right),
 \label{eq:boundary-endpoints}
\end{equation}
with value one at $L=0$.
\end{proposition}

The first endpoint factor occurs in the fixed-cusp estimate of Mirzakhani
and Petri~\cite[Proposition~3.1]{MP}; the second is the fixed-genus,
large-cusp factor of Hide and Thomas~\cite[Theorem~1.5]{HTLargeN}.
Here the genus tends to infinity and the cusp count is unrestricted,
while the number and lengths of the nonzero boundaries remain bounded.
For comparison, Budd and Lions obtain a boundary factor for tight volumes
with a random cusp count in a near-critical regime
\cite[Proposition~4.2]{BL}. Exact transformations between ordinary and
tight boundary generating functions were established in~\cite{BZ}.

\begin{proof}
Put $h=2g-2+p$, $F(x)=\Phi_g^{(p)}(x)$ and
\[
 G(x;\mathbf L)=\sum_{k\ge0}
 \frac{V_{g,p+k}(\mathbf L,0^k)}{(2\pi^2)^{3g-3+p+k}}
 \frac{x^k}{k!}.
\]
We first compare $G$ with $F$ on the real axis, and then show that the
same comparison holds for their coefficients at the saddle.

\emph{The boundary insertion.}
Use the string flow of \eqref{eq:H8}, holding its parameter $Y$ fixed
while differentiating the KW times. Its linear part gives
\[
 \partial_{t_d}=\sum_{b=0}^d\frac{Y^{d-b}}{(d-b)!}\partial_{t_b'}.
\]
For $g\ge2$ the free energy is invariant under the flow. Apply the
displayed identity to each external insertion before setting $x=X(Y)$.
In the volume normalization of \eqref{eq:H5}, an external degree $d$
has coefficient $(L^2/(4\pi^2))^d/d!$ times the unweighted derivative.
Consequently, with $t_i=L_i^2/(4\pi^2)$, the exact identity is
\begin{equation}
 \frac{G(X(Y);\mathbf L)}{F(X(Y))}
 =\sum_{\mathbf a,\mathbf b\ge0}c_{g,p}(\mathbf b;Y)
 \prod_{i=1}^p
 \frac{t_i^{a_i+b_i}Y^{a_i}r^{b_i}}
 {a_i!(a_i+b_i)!(2b_i+1)!!}.
 \label{eq:boundary-insertion}
\end{equation}
Indeed, \eqref{eq:TL9} gives
$U^Y_{g,p}(\mathbf b)/U^Y_{g,p}(\mathbf0)=r^{|\mathbf b|}c_{g,p}(\mathbf b;Y)$.
Terms with invalid $\mathbf b$ are zero; there is no bound on $\mathbf a$.

The identity
\[
 \int_0^1(1-s^2)^b\,\mathrm ds=\frac{2^b b!}{(2b+1)!!}
\]
and the series for $I_0$ show that replacing $c_{g,p}$ by one in
\eqref{eq:boundary-insertion} gives
$\prod_i\mathfrak b_Y(L_i)$. The same double series with an extra
factor $|\mathbf b|(|\mathbf b|+1)$ is bounded uniformly for
$Y\in[0,y_0]$ and $\mathbf L\in[0,L_0]^p$: the inequality
$(a+b)!\ge a!b!$ supplies a summable factorial majorant for every fixed
polynomial moment. Thus \eqref{eq:TL23}, including its one-marking case,
gives
\begin{equation}
 \frac{G(X(Y);\mathbf L)}{F(X(Y))}
 =\prod_i\mathfrak b_Y(L_i)\bigl(1+O_{p,L_0}(g^{-1})\bigr).
 \label{eq:boundary-real}
\end{equation}
All terms in \eqref{eq:boundary-insertion} are nonnegative. The
majorant also justifies convergence of the formal string-flow identity
for $0\le x<x_c$. The factors $\mathfrak b_Y(L)$ are at least one;
they and their first $Y$ derivatives are uniformly bounded on the
compact set under consideration.

\emph{Extraction at a prescribed cusp count.}
The proof of \eqref{eq:TL33} applies to every external degree tuple:
$a(Y)^hU^Y_{g,p}(\mathbf b)$ has nonnegative Taylor coefficients in $Y$.
Before division by $F$, the insertion identity is a sum of these
amplitudes with nonnegative powers of $Y$ and nonnegative constants.
Since $Y(x)$ also has nonnegative coefficients, it follows that
\[
 G(x;\mathbf L)=V(x)^h R_{g,p,\mathbf L}(x),\qquad
 [x^j]R_{g,p,\mathbf L}(x)\ge0.
\]
Convergence follows from \eqref{eq:boundary-real}. Hence
\begin{equation}
 \left|\frac{G(x\mathrm e^{i\theta};\mathbf L)}{G(x;\mathbf L)}\right|
 \le\left|\frac{V(x\mathrm e^{i\theta})}{V(x)}\right|^h.
 \label{eq:boundary-fourier}
\end{equation}

Suppose first that $k=n-p\to\infty$, and choose $x$ by
$h\Psi(x)=k$. Put $\sigma^2=h\mathscr D\Psi(x)$. For fixed real $u$,
the ratio of the centered moment generating functions associated with
$G$ and $F$ is
\[
 \frac{(G/F)(x\mathrm e^{u/\sigma};\mathbf L)}
 {(G/F)(x;\mathbf L)}.
\]
It tends to one, uniformly in $\mathbf L$. To see this, note that the endpoint
expansions \eqref{eq:TL36} and the quadratic inverse at $x_c$ show that
\[
 \sup_{0<x<x_c}\frac{xY'(x)}{\sqrt{\mathscr D\Psi(x)}}<\infty.
\]
On the tilted interval $x\mathrm e^{u/\sigma}$ the model variance is
comparable to its value at $x$, and the interval remains subcritical.
Therefore $Y(x\mathrm e^{u/\sigma})-Y(x)=O_u(h^{-1/2})$.
Now use \eqref{eq:boundary-real} and the bounded derivative of
$\mathfrak b_Y$.

The real Gaussian limit \eqref{eq:TL39} thus holds for $G$ with the same
center $k$ and scale $\sigma$. Its Fourier bound
\eqref{eq:boundary-fourier} is identical to the one used in
\eqref{eq:TL40}--\eqref{eq:TL46}. Fourier inversion consequently gives
the local limit \eqref{eq:TL47} for $G$ as well. At the common index
$k$, both local probabilities are asymptotic to
$(\sqrt{2\pi}\sigma)^{-1}$. Their ratio tends to one, and hence
\[
 \frac{[x^k]G(x;\mathbf L)}{[x^k]F(x)}
 =\prod_i\mathfrak b_{Y}(L_i)(1+o(1)),\qquad
 \lambda(Y)=\frac{2(n-p)}{2g-2+p}.
\]
All errors are uniform on the fixed length box, including along
sequences tending to either endpoint of the parameter interval.

For bounded $k$, pass to a subsequence on which $k$ is fixed.
Equation~\eqref{eq:TL23} at $Y=0$, with $p+k$ markings, gives the
fixed-cusp boundary estimate \eqref{eq:boundary-factor}. This yields the
same conclusion since the saddle then tends to zero.

Finally, replacing $2(n-p)/(2g-2+p)$ by $n/g$ changes the corresponding
$Y$ by $o(1)$ uniformly in $n\ge p$. Indeed, along a subsequence on
which $n/g$ tends to zero, a positive finite limit, or infinity, both
parameters have the same limit in $[0,y_0]$. Continuity of
$\mathfrak b_Y(L)$, uniformly on the length box, proves
\eqref{eq:boundary-uniform}. At $Y=0$, $T=\Delta=\pi^2/4$;
expanding $I_0$ and integrating gives the first formula in
\eqref{eq:boundary-endpoints}. At $Y=y_0$, the integrand is constant,
which gives the second.
\end{proof}

\section{Extraction with an inverse-genus error}
\label{sec:sharp-extraction}

This section completes the third idea of Subsection~\ref{subsec:strategy} with quantitative extraction estimates. The local limit theorem \eqref{eq:TL47} gives no rate. Its quantitative form, Proposition~\ref{prop:quantitative-local} in Appendix~\ref{app:extraction}, gives only a power rate in $\min(h,\mu)$, and its Gaussian approximation degenerates when the coefficient index stays bounded. To obtain an error of order $g^{-1}$ uniformly in the index, we retain an exact Poisson factor rather than its Gaussian approximation. Theorem~\ref{thm:sharp-extraction} supplies the estimate used for the leading absolute-volume formula and the volume-ratio rates in Section~\ref{sec:absolute}, and Theorem~\ref{thm:canonical-expansion} extends it to all orders. Their inputs are the marked expansions of Section~\ref{sec:higher-order}, the convolution factor of Lemma~\ref{lem:convolution} and the endpoint expansions \eqref{eq:TL36}; both theorems are proved in Appendix~\ref{app:extraction}.

\emph{Why a Poisson factor.} Near zero cusp weight the leading term of $\Phi_g^{(p)}$ behaves like $\mathrm e^{hS(x)}$ with $S(x)=2x+O(x^2)$; see~\eqref{eq:CE1} below. At a positive weight $x$, the normalized terms $\mathrm e^{-2hx}(2hx)^k/k!$ are exactly Poisson probabilities with mean $2hx$. Hence, when the coefficient index $k$ stays bounded, the number of additional cusps under the saddle weight is approximately Poisson, as anticipated in Subsection~\ref{subsec:parameter}. A Gaussian approximation of its central point mass then has a relative error that does not tend to zero for fixed positive $k$, and it is undefined at $k=0$. If the leading term is instead written with the exact Poisson mass $\mathfrak p_k=\mathrm e^{-k}k^k/k!$ at the mean, the ratio of the true coefficient to this model is $1+O(1/h)$ for every $k$. For large $k$ one has $\sqrt{2\pi k}\,\mathfrak p_k\to1$, which recovers the Gaussian formula.

Fix $p\ge1$, put $F_h=\Phi_g^{(p)}$, $h=2g-2+p$, and write
\begin{equation}
 \mathcal J_{g,p}(x)=\mathcal J_{g,p}(0)\mathrm{e}^{hS(x)}B_p(x),\qquad
 S=\log\frac{\mathcal Z(0)}{\mathcal Z(Y(x))},\qquad
 B_p=\frac{r^{p/2}\sqrt{A\Theta}}{r(0)^{p/2}\sqrt{A(0)\Theta(0)}}.
 \label{eq:CE1}
\end{equation}
We view $\mathcal J$ as a function of $x=X(Y)$. Then
$\mathscr DS=\Psi$ and $\mathscr D\log B_p=\beta_p\Psi$.
For each integer $k\ge1$, choose $x$ by $h\Psi(x)=k$, and put
$v=\mathscr D\Psi$, $\sigma^2=hv$, $M=\min(h,k)$.
Let $\mathfrak p_k=\mathrm{e}^{-k}k^k/k!$, with $\mathfrak p_0=1$.

\begin{theorem}[Uniform Poisson-corrected extraction]\label{thm:sharp-extraction}
For every fixed $p\ge1$ there is $h_p$ such that, for $h\ge h_p$ and
uniformly over integers $k\ge0$,
\begin{equation}
 [z^k]F_h(z)=\frac{\mathcal J_{g,p}(x)}\pi x^{-k}\mathfrak p_k
       \sqrt{\frac{k}{h\mathscr D\Psi(x)}}\,(1+O_p(h^{-1})).
 \label{eq:CE2}
\end{equation}
At $k=0$ take $x=0$, $x^{-k}=1$ and $\mathfrak p_0=1$; the
variance factor $\sqrt{k/(h\mathscr D\Psi(x))}$ is one by continuity
along $k=h\Psi(x)$.
\end{theorem}

Theorem~\ref{thm:sharp-extraction} is the case $K=0$ of Theorem~\ref{thm:canonical-expansion} below; it is proved separately because its proof uses Theorem~\ref{thm:full-expansion} only at a fixed order, whereas the proof of Theorem~\ref{thm:canonical-expansion} uses the full linear-order range. The proof, given in Subsections~\ref{app:local-complex}--\ref{app:small-indices}, has three parts. First, a sufficiently long marked expansion, namely Theorem~\ref{thm:full-expansion} at a fixed finite order, is transferred to small complex neighborhoods of the real saddle point. Second, the Poisson correction is identified. Third, the small indices are estimated directly. The argument uses no estimate for a growing number of markings and never differentiates a remainder; away from zero, the model is continued analytically only to small neighborhoods of the real axis.

\subsection{The full expansion at a prescribed cusp count}\label{sec:canonical-expansion}

Theorem~\ref{thm:sharp-extraction} recovers the leading term with an inverse-genus error. We now extract the whole marked expansion, keeping the remainder uniform in both the cusp count and the truncation order. The result is Theorem~\ref{thm:canonical-expansion}, which provides the remaining input for the full absolute-volume expansion in Section~\ref{sec:absolute}.

\emph{The idea.} Away from zero cusp weight, including near the critical point, the coefficient integral has a nondegenerate saddle point, and the Gaussian method of steepest descent gives an expansion in $1/h$. Near zero this fails for bounded $k$: the saddle point of $\mathrm e^{hS(z)}z^{-k}$ approaches the singularity of $z^{-k}$ at the origin. Following Temme, we change variables so that the phase $S(z)-(k/h)\log z$ becomes exactly $u-(k/h)\log u$. The coefficient is then a Poisson coefficient $[u^k]\mathrm e^{hu}f(u)$ with an analytic amplitude $f$. Expanding $f$ about $u=k/h$ produces the signed central polynomials $T_m(k)$ of the coefficient operator in~\eqref{eq:CG8}, divided by powers of $h$. These are not the central moments of a probability distribution: for example, $T_2(k)=-k$. The resulting expansion in $1/h$ is uniform in $k$ in this small-weight range, including $k=0$. The two coefficient formulas agree on an overlap of their ranges, which identifies them.

The proof is given in Subsections~\ref{app:phase-change}--\ref{app:canonical-proof}: a Poisson integral is used near zero cusp weight and a Gaussian integral away from zero, including near criticality. Near zero we follow Temme's method for uniform Laplace and loop integrals~\cite[Sections~2.3 and~3.7]{Temme}.

Retain $h=2g-2+p$, $k=h\Psi(x)$, $v=\mathscr D\Psi(x)$ and the
functions $S,B_p$ in \eqref{eq:CE1}. Regard $\mathcal L_{p,j}$
as functions of $x$ via $Y(x)$. Write
\begin{equation}
 \mathcal Q_{g,p,k}=
 \frac{\mathcal J_{g,p}(x)}\pi x^{-k}\mathfrak p_k
                   \sqrt{\frac{k}{hv}}.
 \label{eq:CG1}
\end{equation}
The conventions at $k=0$ are those of Theorem~\ref{thm:sharp-extraction}.

Put $\mathcal L_{p,0}=1$. Let $\mathfrak G_w$ be the linear functional
on polynomials which maps $w^{2a}$ to $(-1)^a(2a-1)!!$ and odd
powers to zero. Define the formal even series
\begin{equation}
\begin{split}
 \mathcal N_p(\tau,x)=\mathfrak G_w\Bigg[
 &\exp\!\left\{
 \sum_{m\ge3}\frac{\tau^{m-2}\mathscr D^mS}{m!v^{m/2}}w^m
 +\sum_{m\ge1}\frac{\tau^m\mathscr D^m\log B_p}{m!v^{m/2}}w^m
 \right\}\\[-2pt]
 &\hspace{20pt}\cdot\sum_{s\ge0}\tau^{2s}
       \mathcal L_{p,s}(x\mathrm{e}^{\tau w/\sqrt v})\Bigg].
\end{split}
 \label{eq:CG2}
\end{equation}
Formal composition in the last line means Taylor expansion at $x$.
At every order only finitely many derivatives are used. If
$\sqrt{2\pi k}\,\mathfrak p_k\sim\sum_{j\ge0}s_jk^{-j}$, $s_0=1$, set
\begin{equation}
 \mathcal A_{p,j}(Y)=[\tau^{2j}]
       \frac{\mathcal N_p(\tau,X(Y))}
            {\displaystyle\sum_{i\ge0}s_i\Psi(X(Y))^{-i}\tau^{2i}}.
 \label{eq:CG3}
\end{equation}
This initially defines the functions for $0<Y<y_0$.

\begin{theorem}[Uniform expansion at a prescribed cusp count]\label{thm:canonical-expansion}
For every fixed $p\ge1$ there are constants $C_p,\mathfrak B_p,c_p,h_p>0$
such that, for $h=2g-2+p\ge h_p$, every integer $k\ge0$ and
$0\le K\le c_ph$,
\begin{equation}
 \left|\frac{[z^k]\Phi_g^{(p)}}{\mathcal Q_{g,p,k}}
       -\sum_{j=0}^K\frac{\mathcal A_{p,j}(Y)}{h^j}\right|
 \le\frac{C_p\mathfrak B_p^{K+1}(K+1)!}{h^{K+1}},
 \qquad h\Psi(X(Y))=k,
 \label{eq:CG4}
\end{equation}
where $\mathcal A_{p,0}=1$.
The coefficients are analytic on one complex neighborhood $\Omega_p$
of $[0,y_0]$, with
$\|\mathcal A_{p,j}\|_{\Omega_p}\le C_p\mathfrak B_p^jj!$, and
$\mathcal A_{p,j}(0)=\mathcal L_{p,j}(0)$.
For a sufficiently small fixed $c'_p>0$, truncation at
$K=\lfloor c'_ph\rfloor$ has relative error
$O_p(\mathrm{e}^{-a_ph})$, uniformly in $k$, for some $a_p>0$.
The marking count is fixed. The infinite series need not converge.
\end{theorem}

The first coefficient $\mathcal A_{p,1}$ is given explicitly in \eqref{eq:CG12}, and its value at the critical endpoint in \eqref{eq:CG13}. Retaining it improves the relative remainder to $O_p(h^{-2})$, uniformly over all coefficient indices. Section~\ref{sec:absolute} applies Theorem~\ref{thm:sharp-extraction} with $p=1$, and with the pairs $(g,1),(g-1,3)$ and $(g,1),(g,2)$ for the two volume ratios, and applies Theorem~\ref{thm:canonical-expansion} with $p=1$ to obtain the coefficients $A_j(g,n)$ of Theorem~\ref{thm:intro-absolute}.

\section{Absolute volumes at every cusp count}\label{sec:absolute}

Specializing the extraction theorems to one distinguished marking,
$p=1$, gives the absolute-volume expansion and the inverse-genus
volume-ratio estimates. Theorem~\ref{thm:sharp-extraction} supplies the
explicit leading term with relative error $O(g^{-1})$;
Theorem~\ref{thm:canonical-expansion} supplies the higher corrections
and their uniform remainder. For $n=0$, the leading estimate also
follows from the fixed-cusp expansion of~\cite{MZ}, and the full
expansion is recovered from $n=1$ by the exact dilaton relation.

For $g\ge2$ and $n\ge1$, set $h=2g-1$, $k=n-1$, and choose $x=X(Y)$
by $h\Psi(x)=k$. At $k=0$ take $x=Y=0$. Define
\begin{equation}
 \begin{split}
 \widehat V_{g,n}&=
 (2\pi^2)^{3g-3+n}\frac{\mathcal J_{g,1}(Y)}{\pi}
 \left(\frac{k}{\mathrm{e} x}\right)^k
 \sqrt{\frac{k}{h\mathscr D\Psi(x)}}\quad(n\ge1),\\
 \widehat V_{g,0}&=
 \frac{\Gamma(2g-2)(4\pi^2)^{2g-3}}{\sqrt{\pi g}}.
 \end{split}
 \label{eq:AV1}
\end{equation}
At $k=0$, both $k$-dependent factors on the first line are defined by
their limits, which are one. All functions in this formula were defined
in \eqref{eq:TL32}, \eqref{eq:PF1}, and \eqref{eq:PF3}.
We have $r=2\Delta$, $Aar=I/\sqrt\Delta$ by \eqref{eq:TL32}, and
$\Theta=-1/\bar K_Y'(1)=-a/(r\,\partial_zK_Y(r))$ by \eqref{eq:PF1}.
Since $X''(v)=-X(v)/v$, differentiating $K_Y$ gives
$-2\partial_zK_Y(r)=\int_0^1(1-s^2)X(v)v^{-1}\,\mathrm ds$ with
$v=T-\Delta s^2$. Hence $rA\Theta=\varpi(Y)^2$, that is,
$r^{1/2}\sqrt{A\Theta}=\varpi(Y)$, and $\mathcal Z=\sqrt2\,I$, with
$\varpi$ as in the introduction. Hence \eqref{eq:AV1} is the formula
\eqref{eq:intro-leading}, with $\mathscr D\Psi(x)=x\Psi'(x)$.

\begin{theorem}[The leading term]\label{thm:absolute}
There is a constant $C>0$ such that, for all sufficiently large $g$,
\begin{equation}
 \sup_{n\ge0}\left|\frac{V_{g,n}}{\widehat V_{g,n}}-1\right|
       \le\frac Cg.
 \label{eq:AV2}
\end{equation}
The order $g^{-1}$ is sharp for the leading model \eqref{eq:AV1}:
\begin{equation}
 \liminf_{g\to\infty}g\sup_{n\ge0}
       \left|\frac{V_{g,n}}{\widehat V_{g,n}}-1\right|
       \ge\frac7{12}-\frac{17}{6\pi^2}>0.
 \label{eq:AV3}
\end{equation}
\end{theorem}
\begin{proof}
For $n\ge1$, apply Theorem~\ref{thm:sharp-extraction} with $p=1$.
The coefficient identity \eqref{eq:TL50} gives
\begin{equation}
 V_{g,n}=(2\pi^2)^{3g-3+n}(n-1)!x^{-k}
     \frac{\mathcal J_{g,1}(Y)}\pi\,
       \mathfrak p_k\sqrt{\frac{k}{h\mathscr D\Psi(x)}}
       (1+O(g^{-1})),
 \label{eq:AV4}
\end{equation}
uniformly, including $k=0$ by continuity.
Since $k!\mathfrak p_k=(k/\mathrm{e})^k$, this is precisely the first line of
\eqref{eq:AV1} with uniform relative error $O(g^{-1})$.
For $n=0$, the constant was independently identified in
Section~\ref{sec:prefactor}, so the fixed-cusp expansion of
Mirzakhani--Zograf~\cite[Theorem~1.2]{MZ} gives
$V_{g,0}/\widehat V_{g,0}=1+c^{(1)}_0g^{-1}+O(g^{-2})$, where
$c^{(1)}_n$ denotes their first coefficient. Its value is stated
in~\cite[Remark~1.3]{MZ}, and it also follows from two proved
statements. In the normalization of~\cite{MZ}, the first-correction
theorem of Liu and Xu~\cite[Theorem~1.3, at $t=1$]{LXConstant} gives
$c^{(1)}_2=1/12+1/(6\pi^2)$. Indeed, the conversion between the two
leading terms and the corresponding first coefficients is
\[
 \begin{aligned}
 \sqrt{2g}\,\frac{\Gamma(2g-\tfrac12)}{\Gamma(2g)}
     &=1+\frac3{16g}+O(g^{-2}),\\
 c^{(1)}_2=\mathfrak b_2(1)+\frac3{16}
     &=-\frac5{48}+\frac1{6\pi^2}+\frac3{16}
       =\frac1{12}+\frac1{6\pi^2},
 \end{aligned}
\]
where $\mathfrak b_2$ is the coefficient in~\cite{LXConstant}.
The cusp-ratio expansion
\cite[Theorem~1.4]{MZ} gives
$c^{(1)}_{n+1}-c^{(1)}_n=2/\pi^2-1/4-n/\pi^2$ for $n=0,1$. Hence
\begin{equation}
 \frac{V_{g,0}}{\widehat V_{g,0}}
 =1+\frac1g\left(\frac7{12}-\frac{17}{6\pi^2}\right)+O(g^{-2}).
 \label{eq:AV5}
\end{equation}
This closes \eqref{eq:AV2} and proves \eqref{eq:AV3}.

The order $g^{-1}$ is also attained at large cusp density, without
the value of $c^{(1)}_0$. By Theorem~\ref{thm:canonical-expansion}
with $p=1$ and $K=1$, together with the coefficient identity used in
the proof of Theorem~\ref{thm:intro-absolute} below,
$V_{g,n}/\widehat V_{g,n}=1+\mathcal A_{1,1}(Y_{g,n})/(2g-1)+O(g^{-2})$
uniformly in $n\ge1$, where $(2g-1)\Psi(X(Y_{g,n}))=n-1$. For fixed
$g$, $Y_{g,n}\to y_0$ as $n\to\infty$, and $\mathcal A_{1,1}(y_0)=-29/120$
by \eqref{eq:CG13}. Hence the left side of \eqref{eq:AV3} is also at
least $29/240$.
The latter concerns the uncorrected model chosen here, rather than an
obstruction to adding further asymptotic corrections.
\end{proof}

The model \eqref{eq:AV1} contains no unknown constant: its
normalization was fixed at the pure-$\psi$ endpoint in
Section~\ref{sec:prefactor}, and the fixed-cusp expansion of~\cite{MZ}
enters only afterwards, for $n=0$ in the proof above. The
matrix-model reading of its action $I(Y)$ and of the parameter $Y$ is
given in Subsection~\ref{subsec:jt-edge}.

\subsection{An inverse-genus error for the volume ratios}
\label{subsec:quantitative-volume-ratios}

The sharper extraction estimate also gives a uniform rate in the
volume-ratio theorem. For the genus-exchange ratio, the two generating
functions have the same saddle. For the cusp ratio, adding a cusp changes
the complexity by one, and the argument below controls the resulting
saddle displacement even near the endpoints.

\begin{proof}[Completion of the proof of Theorem~\ref{thm:volume}]
For $n\ge1$ put $h=2g-1$, $k=n-1$, and $h\Psi(x)=k$.
Apply Theorem~\ref{thm:sharp-extraction} at $(g,p)=(g,1)$ and
$(g-1,3)$. Their Poisson and variance factors cancel exactly, and
\eqref{eq:PF3} gives $\mathcal J_{g-1,3}/\mathcal J_{g,1}=4r$.
Consequently \eqref{eq:TL50} yields
\[
 B_{g,n}=\mathcal B(2k/h)(1+O(h^{-1})).
\]

For the cusp ratio assume first that $k\ge1$, and let $x'$ satisfy
$(h+1)\Psi(x')=k$. Write $q=\log x$, $q'=\log x'$ and
$v=\mathscr D\Psi$. The convergent endpoint expansions give
\[
 \sup_{0<x<x_c}\left\{
 \frac{\Psi(x)^2}{v(x)}+
 \left|\frac{\Psi(x)}{v(x)}\mathscr D\log v(x)\right|
 +|\beta_2(x)|\right\}<\infty.
\]
Indeed, near zero $\Psi\sim v\sim2x$; near criticality
$\Psi\sim5/(4\varepsilon)$ and $v\sim5/(4\varepsilon^2)$.
The function $\beta_2$ is bounded by \eqref{eq:PF7}.
Since $\Psi(x')/\Psi(x)=h/(h+1)$, integration with respect to
$\log\Psi$ shows that
\[
 \frac{v(x')}{v(x)}=1+O(h^{-1}),\qquad
 |q'-q|\le\frac{C\Psi(x)}{h v(x)}.
\]
For $L(q)=\log\mathcal J_{g,2}(\mathrm{e}^q)-kq$, the exact model identities give
\[
 L'(q)=(h+1)\Psi(\mathrm{e}^q)-k+\beta_2(\mathrm{e}^q)\Psi(\mathrm{e}^q).
\]
On $[q',q]$ its absolute value is at most $C\Psi(x)$, so
$|L(q')-L(q)|\le C\Psi(x)^2/(hv(x))=O(h^{-1})$.
Thus the saddle displacement and the variance factor in
\eqref{eq:CE2} each change the coefficient ratio by $1+O(h^{-1})$.
At the common point $x$, \eqref{eq:PF3} gives exactly
\[
 \frac{\mathcal J_{g,2}(x)}{\mathcal J_{g,1}(x)}
 =\frac{h}{Aar}\sqrt{\frac{h}{h+1}}.
\]
Using the coefficient normalization in \eqref{eq:TL53}, we obtain
\[
 C_{g,n}=\frac{h}{2g}\mathcal C(2k/h)(1+O(h^{-1})).
\]
At $k=0$ both saddles are zero. The conventions in
Theorem~\ref{thm:sharp-extraction} give the same two ratio formulas
directly, without logarithmic coordinates.

The positive profiles are analytic at zero, and their expansions at
infinity in \eqref{eq:SP13}--\eqref{eq:SP15} are convergent in the
inverse square root of the density. Hence
\[
 \sup_{\lambda\ge0}(1+\lambda)
 \bigl(|(\log\mathcal B)'(\lambda)|+
       |(\log\mathcal C)'(\lambda)|\bigr)<\infty.
\]
For $\lambda_*=2(n-1)/(2g-1)$ and $\lambda=n/g$,
$|\lambda_* -\lambda|/(1+\min(\lambda_*,\lambda))=O(g^{-1})$
uniformly in $n\ge1$. Integrating the preceding bound replaces
$\lambda_*$ by $\lambda$ with relative error $O(g^{-1})$.
For $n=0$, the fixed-cusp expansion of~\cite{MZ} at the finitely
many cusp counts $0,1,2$ gives the required rate directly.
Finally $D_{g,n}=C_{g,n-1}/C_{g,n}$ and the same logarithmic
derivative bound prove $\sup_{n\ge1}|D_{g,n}-1|=O(g^{-1})$.
\end{proof}

\subsection{Proof of Theorem~\ref{thm:intro-absolute}}

For $n\ge1$ the expansion follows from the extraction theorem,
Theorem~\ref{thm:canonical-expansion}, with $p=1$ and the same leading
term. The unmarked case $n=0$ is recovered from $n=1$ by the exact
dilaton relation, which for volume polynomials is the identity
$\partial_LV_{g,1}(2\pi i)=2\pi i(2g-2)V_{g,0}$ of Do and
Norbury~\cite[Theorem~2]{DoNorbury}; see also~\cite[equation~(3.5)]{Mir}.
Its higher coefficients therefore need no separate asymptotic
normalization. For $h=2g-1$, the scalar $\mathsf C_0$ of \eqref{eq:FE1} obeys
\begin{equation}
 V_{g,0}=\frac{V_{g,1}\mathsf C_0(h-1;0)}{4\pi^2(h-1)},\qquad
 \frac{\widehat V_{g,1}}{4\pi^2(h-1)\widehat V_{g,0}}
       =\sqrt{1+\frac1h}.
 \label{eq:AV6}
\end{equation}
Indeed \eqref{eq:TL12} holds at $p=0$; at $Y=0$, undo its
normalizations with $r=\pi^2/2$, $A=1/\pi^2$. Equivalently, the first
identity in \eqref{eq:AV6} is the case $n=0$
of~\cite[Theorem~2]{DoNorbury}.
The second identity follows exactly from \eqref{eq:AV1}.
The scalar expansion of $\mathsf C_0$ follows directly from the last
line of \eqref{eq:FE5}, using the proved one-marking profile and its
summable remainder; this step does not require an unmarked harmonic
comparison.

The coefficients in Theorem~\ref{thm:intro-absolute} are defined as
follows. For $n=0$, define constants $a_{0,j}$ by the finite-order formal
rule
\begin{equation}
 1+\sum_{j\ge1}a_{0,j}t^j
 =(1+t)^{1/2}\mathsf C_0\!\left(\frac t{1-t};0\right)
       \left(1+\sum_{j\ge1}\mathcal L_{1,j}(0)t^j\right).
 \label{eq:AV7}
\end{equation}
Here the second argument denotes $Y=0$, and $t=1/h$; thus the
cusp scalar uses complexity $h-1$. Put $A_j(g,0)=a_{0,j}$, and, for
$n\ge1$, put $A_j(g,n)=\mathcal A_{1,j}(Y_{g,n})$ with the coefficients
\eqref{eq:CG3} of Theorem~\ref{thm:canonical-expansion}, where
$(2g-1)\Psi(X(Y_{g,n}))=n-1$.

\begin{proof}[Proof of Theorem~\ref{thm:intro-absolute}]
For $n\ge1$, take $p=1,k=n-1$ in
Theorem~\ref{thm:canonical-expansion}. The coefficient identity \eqref{eq:TL50}
and $k!\mathfrak p_k=(k/\mathrm{e})^k$ identify $[z^k]\Phi_g^{(1)}/\mathcal Q_{g,1,k}$
with $V_{g,n}/\widehat V_{g,n}$, including $k=0$.
For $n=0$, the dilaton scalar in \eqref{eq:AV6} has a one-marking
child of complexity $h$, so exactly
\[
 \mathsf C_0(h-1;0)=\frac1{A(0)}
       \sum_{d\ge1}b_d(0)c_{g,1}(d;0).
\]
Apply \eqref{eq:gevrey-remainders} at this child with $L=K+1$.
Fix $0<\sigma<\log\rho$ within the absolute exponential moment
of $b_d$. Then
\[
 (1+ML)^L\sum_d\frac{|b_d(0)|}{A(0)}\mathrm{e}^{\sigma d}
 \le C D^LL!.
\]
This gives a Gevrey-1 remainder for the whole scalar sum, including
invalid degrees. Its coefficient series is precisely
$\mathsf C_0(t/(1-t);0)$: the child substitution in \eqref{eq:FE5}
becomes $(t/(1-t))/(1+t/(1-t))=t$.

The other two exact factors in \eqref{eq:AV6} are
$\sqrt{1+1/h}$, with a convergent expansion, and the $n=1$ ratio,
with the marked estimate at zero. All three exact factors are bounded.
The finite-product identity used in Proposition~\ref{prop:hierarchy}
contains one remainder per term; its factorial weights obey
$j!(L-j)!\le L!$. Applied twice, it gives the same
bound for their product and exactly the coefficients in \eqref{eq:AV7}.
Take common constants with the $n\ge1$ case. The coefficient bound
$|A_j(g,n)|\le CB^jj!$ follows from those of
Theorem~\ref{thm:canonical-expansion} and Proposition~\ref{prop:gevrey-hierarchy}.
Finally, for $K=\lfloor c'(2g-1)\rfloor$ with $c'\le c$ chosen after
the final base $B$ so that $Bc'<1$, the inequality
$B^LL!h^{-L}\le(BL/h)^L$ gives the exponentially small error.
\end{proof}

In particular, the explicit first correction \eqref{eq:CG12}
gives a relative error $O(g^{-2})$ uniformly over all cusp counts.
At the unmarked endpoint, by uniqueness of asymptotic expansions,
$a_{0,1}=7/6-17/(3\pi^2)$, the first coefficient of \eqref{eq:AV5}
once $h=2g-1$ is substituted.
Equivalently, by \eqref{eq:AV7} and $\beta_0(0)=4/\pi^2-1/2$ from
\eqref{eq:PF7}, $\mathcal L_{1,1}(0)=1/6-5/(3\pi^2)$, that is,
\[
 \int_0^{y_0}\frac{\gamma_1(v)}{A(v)r(v)}\,\mathrm{d}v
 =-\frac{13}{60}+\frac{5}{3\pi^2},
\]
a value that is not computed independently here.
The existence of an all-order expansion at a fixed cusp count is the
theorem of Mirzakhani--Zograf~\cite{MZ}; here the coefficients and
remainders remain controlled as the cusp count crosses every density
scale. This completes the proofs of Theorems~\ref{thm:intro-absolute}
and~\ref{thm:volume}. The remaining sections of this part treat the
diagonal chain and the coefficient ratios, which use the volume-ratio
limits and the harmonic comparison respectively; the volume ratios of
Theorem~\ref{thm:volume}, together with
Proposition~\ref{prop:bounded-boundaries}, are the inputs that
Parts~\ref{part:spectra} and~\ref{part:geometry} take from the volume
theory, beginning in Section~\ref{sec:geometry}.

\section{The full diagonal chain}\label{uniform-logarithmic-volume-rate-and-the-full-diagonal-chain}

This section proves Theorem~\ref{thm:chain}, the logarithmic rate along
the diagonal chain announced in the introduction, after
Theorem~\ref{thm:coefficient}.
The absolute expansion of Theorem~\ref{thm:intro-absolute} assumes
$g\to\infty$. Along the chain $V_{G,0},V_{G-1,2},\ldots,V_{0,2G}$, on
which the Euler characteristic is constant and genus is traded for pairs
of cusps, one also needs the bounded-genus end. To reach it we sum the
logarithms of the cusp ratios of Theorem~\ref{thm:volume}, whose profile
$\mathcal C$ is expressed through the parameter functions of
Lemma~\ref{lem:parameter}, and compare the result with the fixed-genus
asymptotics of Manin and Zograf~\cite{ManZ}. The result has only
logarithmic precision even though the cusp ratios are known with
relative errors, because summing $n$ of their logarithms accumulates an
error $n\varepsilon_g$ with $\varepsilon_g\to0$. The proof is given in
full here; apart from the limits in Theorem~\ref{thm:volume} and the
classical fixed-cusp and fixed-genus asymptotics, it uses only the
parameter identities of Lemma~\ref{lem:parameter}.

Let $Y(u)$ correspond to $\lambda=2u/(1-u)$ and define
\begin{equation}
 F(t)=\int_0^t\log\frac{4\Delta(Y(u))}{\pi^2}\,\mathrm{d}u
 \quad(0\le t<1),\qquad F(1)=-2\log(2x_c).
 \label{eq:WC1}
\end{equation}
The endpoint value agrees with the limit of the integral.
\Needspace{10\baselineskip}
\begin{theorem}[A logarithmic rate for the full chain]\label{thm:chain}
Uniformly over stable $(g,n)$ with $2g+n\to\infty$,
\begin{equation}
 \log\frac{V_{g,n}}{(2g+n-3)!(4\pi^2)^{2g+n-3}}
 =\frac{2g+n}{2}F\!\left(\frac n{2g+n}\right)+o(2g+n).
 \label{eq:WC2}
\end{equation}
Consequently,
\begin{equation}
 \sup_{0\le j\le G}
 \left|\frac1G\log\frac{V_{G-j,2j}}{V_{G,0}}-F(j/G)\right|\to0.
 \label{eq:WC3}
\end{equation}
\end{theorem}
The second statement determines each entry of the chain up to a factor
$\exp(o(G))$, uniformly in its position and including the planar
endpoint $V_{0,2G}$. When the genus tends to infinity,
Theorem~\ref{thm:intro-absolute} gives far more precise information;
the point of Theorem~\ref{thm:chain} is that it also covers bounded
genus.

We use the functions \(T,\Delta,I,J\) and the profile \(\mathcal C\)
of Section~\ref{the-stationary-parameter}, and the variable
\(t=\lambda/(2+\lambda)\in[0,1)\); thus \(Y(u)\) in \eqref{eq:WC1} is the
parameter with \(\lambda=2u/(1-u)\).

\begin{proof}[Proof of Theorem~\ref{thm:chain}]
We first give a closed form for \eqref{eq:WC1}. From \eqref{eq:SP11}
and differentiation under the integral, \[
 I'=-X'(Y)\sqrt\Delta,
 \qquad J'=-\frac{I\Delta'}{2\Delta^{3/2}}.
\] Writing \(t=\beta(Y)=X/J\), these identities show that \begin{equation}
F(t)=-2\log(2J(Y))
      -(1-t)\log\frac{4\Delta(Y)}{\pi^2}
 \label{eq:WC4}
\end{equation} has derivative \(\log(4\Delta/\pi^2)\). Its value at zero is zero,
because \(J(0)=1/2\) and \(\Delta(0)=\pi^2/4\); hence it is precisely
\eqref{eq:WC1}. As \(t\to1\), \(J\to x_c\) and \((1-t)\log\Delta\to0\), proving
its stated continuous endpoint.

Set \[
 Z_{g,n}=\frac{V_{g,n}}{(2g+n-3)!(4\pi^2)^{2g+n-3}},\qquad
 f(\lambda)=\log\frac{2\mathcal C(\lambda)}{4\pi^2(2+\lambda)}
            =-\log(2J(Y)).
\] Since \(J'>0\), \(f\) decreases from zero to \(-\log(2x_c)\) and has
finite total variation. Equation \eqref{eq:WC4} gives the exact identity \begin{equation}
\int_0^\lambda f(v)\,\mathrm{d}v
   =\left(1+\frac\lambda2\right)
       F\left(\frac\lambda{2+\lambda}\right).
 \label{eq:WC5}
\end{equation}

For \(g\ge2\), the exact adjacent-cusp relation is \[
 \log Z_{g,n}=\log Z_{g,0}
  +\sum_{j=0}^{n-1}
      \log\frac{2g C_{g,j}}{4\pi^2(2g+j-2)}.
\] Let \(\varepsilon_g\to0\) bound the logarithmic error in the uniform
\(C\) profile. Replacing \(C\) by \(\mathcal C\), and replacing
\(2g+j-2\) by \(2g+j\), gives total error
\(O(n\varepsilon_g+1+\log(1+n/g))\). The bounded variation of \(f\) gives
the uniform Riemann-sum estimate \[
 \sum_{j=0}^{n-1}f(j/g)
       =g\int_0^{n/g}f(v)\,\mathrm{d}v+O(1),
\] with a constant independent of \(n\). The fixed-cusp closed-volume
asymptotic gives \(\log Z_{g,0}=o(g)\). Using \eqref{eq:WC5}, we obtain \begin{equation}
 \log Z_{g,n}
  =\frac{2g+n}{2}F\left(\frac{n}{2g+n}\right)+o(2g+n),
 \label{eq:WC2a}
\end{equation} uniformly for all \(n\) as \(g\to\infty\).

The same statement holds uniformly over \emph{all} stable pairs with
\(2g+n\to\infty\). To check the remaining case, keep \(g\) fixed and use
the fixed-genus leading cusp asymptotic~\cite[Section~6]{ManZ}, which gives \[
 \log Z_{g,n}=-n\log(2x_c)+O_g(\log(n+1)).
\] Meanwhile \(F(1)-F(t)=O((1-t)|\log(1-t)|)\), because
\(\Delta\asymp\sqrt{1-t}\) by Lemma~\ref{lem:parameter}(iii), so the right
side of \eqref{eq:WC2} differs from the same leading term by only
\(O_g(\log(n+1))\). Every sequence of stable pairs has a subsequence on
which the genus is fixed or tends to infinity. This proves the full
uniformity claimed for \eqref{eq:WC2}.

Along \((g,n)=(G-j,2j)\), the value \(2g+n=2G\) and the factorial
normalization are constant. Subtracting the closed value
\(\log Z_{G,0}=o(G)\) from \eqref{eq:WC2} proves \eqref{eq:WC3}, including \(j=0\) and
\(j=G\). This full-chain statement includes bounded genera, whereas the
absolute equivalent in Theorem~\ref{thm:absolute} is uniform in cusp count
as the genus tends to infinity.\qedhere
\end{proof}

\section{Uniform one-boundary coefficient ratios}\label{uniform-one-boundary-coefficient-ratios}

This section proves the uniform equivalent \eqref{eq:C2} in
Theorem~\ref{thm:coefficient}, which answers Mirzakhani's question on
adjacent coefficients of $V_{g,1}(L)$; the quantitative error stated in the
theorem is deferred to Appendix~\ref{sec:quantitative-coefficients}.
The coefficient problem uses the first idea of
Subsection~\ref{subsec:strategy} only: the harmonic comparison of
Proposition~\ref{prop:harmonic} at $Y=0$, in the one-marking form
\eqref{eq:H44}, but neither the marked ratios nor the extraction
theorems. The comparison determines the coefficients only on the
exponential scale; log-concavity of the underlying mixed intersection
numbers, which follows from the nefness of the tautological classes,
converts this logarithmic information into adjacent ratios in the
interior. Near the top degree we use an exact intersection number, and
near degree zero the bracket difference estimates of Mirzakhani and
Zograf~\cite{MZ}. All proofs are given here.

Put \(N=3g-2\) and retain the defining normalization \[
 V_{g,1}(L)=\sum_{k=0}^{N}\frac{a_{g,k}}{(2k+1)!}L^{2k},
 \qquad b_{g,k}=[\tau_k]_{g,1}=4^k a_{g,k}.
\] For every valid adjacent index \(0\le k<N\), define
\(t_{g,k}\in(0,\pi)\) by \begin{equation}
\alpha_{g,k}:=\frac{2k+1}{2g-1}
 =\frac{t_{g,k}^{2}\sin t_{g,k}}
        {\sin t_{g,k}-t_{g,k}\cos t_{g,k}}.
 \label{eq:C1}
\end{equation} The assertion to prove is \begin{equation}
\sup_{0\le k<3g-2}
 \left|\frac{a_{g,k}/a_{g,k+1}}{4\pi^2/t_{g,k}^{2}}-1\right|
 \longrightarrow0.
 \label{eq:C2}
\end{equation} Equivalently, the comparison value for \(b_{g,k}/b_{g,k+1}\) is
\(\pi^2/t_{g,k}^2\). The ratio at \(k=N\) is undefined and is not
included. In the boundary-length normalization of \(V_{g,1}(L)\), the ratio tends to four when $k=o(g)$, while it diverges when
$1\le N-k=o(g)$.

\subsection{Coefficient identities and further inputs}\label{exact-dictionary-and-the-inputs-needed-here}

Let \[
 J_k=\int_{\overline{\mathcal M}_{g,1}}
       \kappa_1^{N-k}\psi_1^k,
 \qquad m=N-k,\quad h=2g-1.
\] The exact coefficient identities are \begin{equation}
a_{g,k}=(2k+1)!!\frac{(2\pi^2)^m}{m!}J_k,
 \qquad
 \frac{a_{g,k}}{a_{g,k+1}}
 =\frac{2\pi^2}{(2k+3)(N-k)}\frac{J_k}{J_{k+1}}.
 \label{eq:C3}
\end{equation} Here \(J_k\) is a mixed \(\kappa_1\)--\(\psi_1\)
intersection number, not a pure \(\psi\) intersection number.

We use the following consequences of Section~\ref{uniform-comparison}.
Define \[
 H(z)=\frac{3(\sin t-t\cos t)}{t^3},\quad t=\sqrt{2z},
 \qquad \widetilde A_m(h)=[z^m]H(z)^{-h},
\] with the removable value \(H(0)=1\), and \begin{equation}
M_{g,k}=4^k\frac{3^h\Gamma(h)}{4^g}
                  (2\pi^2)^m \widetilde A_m(h).
 \label{eq:C4}
\end{equation} Here \(H=H_0\) is the kernel of Section~\ref{the-positive-compact-family-of-kernels}
at \(Y=0\). The tilde indicates that \(\widetilde A_m(h)\) is a coefficient of the
unscaled kernel \(H_0\); thus \(\widetilde A_m(h)=r_0^{-m}A_m(h)\), with
\(A_m(h)\) from \eqref{eq:H12} at \(Y=0\) and \(r_0=\pi^2/2\).
There is an absolute \(C\) such that, for every sufficiently large
\(g\), uniformly in \(0\le k\le N\), \begin{equation}
\exp[-C\log^2(g+1)]
 \le\frac{b_{g,k}}{M_{g,k}}
 \le\exp(C\sqrt g).
 \label{eq:C5}
\end{equation} This is the one-marking case of
\eqref{eq:H42}--\eqref{eq:H44}. At \(Y=0\), the unscaled amplitude is
\(U^0=b/[4^k(2\pi^2)^m]\); undoing that scaling gives exactly
\eqref{eq:C4}, without an additional genus-dependent factor. In particular \begin{equation}
\sup_{0\le k\le N}\frac1g
       \left|\log\frac{b_{g,k}}{M_{g,k}}\right|\longrightarrow0.
 \label{eq:C6}
\end{equation} The lower estimate in \eqref{eq:C5} will also control the bottom endpoint.

We also use the following three facts, which are independent of the
harmonic comparison.

\begin{enumerate}[label=\arabic*.]
\item
  The cotangent classes are nef on the compactification. Indeed, the
  nef logarithmic relative dualizing class in~\cite[Theorem~0.4]{Keel}
  is the cotangent class at the additional mark under the identification
  of the universal curve with the moduli space with one additional mark.
  The logarithmic
  kappa class is nef as well: for the forgetful universal curve
  \(\pi:\overline{\mathcal M}_{g,n+1}\to\overline{\mathcal M}_{g,n}\),
  one has \(\kappa_1=\pi_*(\psi_{n+1}^2)\)~\cite{AC}. Its intersection against any
  complete curve is the square of the pulled-back nef cotangent class on
  the pulled-back universal surface, and hence is nonnegative. The
  argument can be carried out after passing to a finite cover supporting
  a family, and on its projective coarse surface; passage to the coarse
  space multiplies all mixed intersection numbers by the same positive
  factor and does not affect the inequality.
  Hence these nef classes, as rational divisor classes on the projective
  coarse moduli space, satisfy the Khovanskii--Teissier inequalities~\cite[Section~1.6]{Laz}; see also~\cite[Theorem~2.15]{BFJ}. This gives \begin{equation}
J_k^2\ge J_{k-1}J_{k+1},\qquad 1\le k<N.
  \label{eq:C7}
\end{equation} Consequently \(q_k:=J_k/J_{k+1}\) is nondecreasing. No log-concavity of the
  normalized bracket sequence is assumed.
\item
  Bracket monotonicity and the fixed-marking absolute difference estimate of Mirzakhani--Zograf~\cite[Theorem~3.2(ii)]{MZ} give, with \(c_k=b_{g,k}/V_{g,1}\), \begin{equation}
c_0=1,\qquad 0\le c_k-c_{k+1}\le C\frac{k+1}{g},
  \qquad 0\le k<N.
  \label{eq:C8}
\end{equation} For \(k\ge1\), their explicit bound is
  \(b_{g,k}-b_{g,k+1}\le16((k+1)V_{g,0}+kV_{g-1,2})\); the fixed-marking
  volume ratios turn it into \eqref{eq:C8}. The \(k=0\) case is the fixed-degree
  estimate \(1-c_1=O(1/g)\). Summing \eqref{eq:C8} gives \(1-c_k\le Ck(k+1)/g\).
  These are absolute, not relative, estimates, so they cannot be divided
  by an arbitrarily small \(c_k\).
\item
  The fixed-marking leading volume asymptotic gives \begin{equation}
V_{g,1}\asymp g^{-1/2}\Gamma(2g-1)
                      (4\pi^2)^{2g-2},
  \label{eq:C9}
\end{equation} and the one-point pure-$\psi$ identity is \(J_N=1/(24^g g!)\). Only
  the upper bound in \eqref{eq:C9} is needed below. The exact value of the top
  ratio is derived below from the ordinary DVV and string equations.
\end{enumerate}

\subsection{The saddle equation for the degree}\label{the-positive-saddle-and-its-effective-density}

The reciprocal of \(H\) has strictly positive Taylor coefficients and is
analytic on a neighborhood of \([0,\pi^2/2]\). These facts also follow
directly from \(2zH''+5H'+H=0\): writing \(-H'/H=\sum_{j\ge0}r_jz^j\)
gives \(r_0=1/5\) and \[
 (2j+5)r_j=2\sum_{u+v=j-1}r_ur_v>0\quad(j\ge1).
\] Hence \(-\log H\) and \(H^{-1}\) have positive nonconstant
coefficients. The function \(H\) has no positive real zero up to and
including \(\pi^2/2\), because \(\sin t-t\cos t>0\) for \(0<t\le\pi\); by
Pringsheim's theorem this excludes a smaller complex singularity of
\(H^{-1}\).

For \(0<z\le\pi^2/2\), tilt these coefficients to a probability law with
generating function \(H(z)/H(zw)\). Its mean is \begin{equation}
\mu(z)=-\frac{zH'(z)}{H(z)}
 =\frac12\left(3-
     \frac{t^2\sin t}{\sin t-t\cos t}\right),
 \qquad t=\sqrt{2z}.
 \label{eq:C10}
\end{equation} Its variance is \(z\mu'(z)>0\). Thus the right side of \eqref{eq:C1}, as a
function of \(t\), decreases strictly from three to zero on \((0,\pi)\),
proving existence and uniqueness in \eqref{eq:C1}. Note also the exact dimension
identity \begin{equation}
\frac mh=\frac{3-\alpha_{g,k}}2,
 \qquad \alpha_{g,k}=3-\frac{2m}{h}.
 \label{eq:C11}
\end{equation}

For \(m/h\) in a compact subset of \((0,3/2)\), there is a unique saddle
\(z\) with \(\mu(z)=m/h\). Fourier inversion for the tilted iid sum,
whose laws have uniform exponential moments and positive probabilities
at zero and one, gives \begin{equation}
\log \widetilde A_m(h)=-h\log H(z)-m\log z+O(\log h).
 \label{eq:C12}
\end{equation} The error is uniform on such compact subsets.

\subsection{Interior degrees: exponential comparison and concavity}\label{interior-degrees-exponential-comparison-plus-concavity}

The idea is that, for concave functions converging pointwise, slopes
also converge wherever the limit is differentiable. Suppose \(k/g\to\alpha\in(0,3)\) and let \(z=t^2/2\) solve
\(\mu(z)=(3-\alpha)/2\). Interpolate linearly on the degree mesh the
values \begin{equation}
f_g(k/g)=\frac1g\log
               \frac{J_k}{J_Ng^{2(N-k)}}.
 \label{eq:C13}
\end{equation} By \eqref{eq:C7}, these functions are concave. Equations \eqref{eq:C3}--\eqref{eq:C6}, \eqref{eq:C12},
and Stirling's formula give their pointwise limit on \((0,3)\): \begin{equation}
\begin{aligned}
 f(\alpha)={}&-2\log H(z)-(3-\alpha)\log z
 +(3-\alpha)\log(3-\alpha)\\
 &+3\log6-\alpha\log(2\alpha)-2(3-\alpha).
 \end{aligned}
 \label{eq:C14}
\end{equation} Explicitly, the model for the intersection number \(J_k\) is
\(J_k^{\mathrm{mod}}=3^h\Gamma(h)4^{-g}m!\widetilde A_m(h)/(2k+1)!!\). Its quotient
by the top model value is exactly \(m! \widetilde A_m(h)(2N+1)!!/(2k+1)!!\); the
two occurrences of \((N-k)\log g\) in its Stirling expansion cancel the
affine normalization in \eqref{eq:C13}. The logarithmic errors in \eqref{eq:C6}, including
the top one, vanish after division by \(g\).

Differentiation of \eqref{eq:C14} gives \begin{equation}
f'(\alpha)=\log\frac{z}{2\alpha(3-\alpha)}.
 \label{eq:C15}
\end{equation} The terms containing \(dz/d\alpha\) cancel by the saddle equation.
Concavity of $\log J_k$ now upgrades the exponential-scale information to
an adjacent ratio. For fixed small \(\varepsilon>0\), the slope on \([k/g,(k+1)/g]\)
is at most the secant slope on
\([\alpha-2\varepsilon,\alpha-\varepsilon]\) and at least the one on
\([\alpha+\varepsilon,\alpha+2\varepsilon]\), for all sufficiently large
\(g\). These four endpoints are fixed, so pointwise convergence
suffices. First let \(g\) tend to infinity and then let \(\varepsilon\)
decrease to zero. Since \eqref{eq:C14} is differentiable, the adjacent slope
converges to \eqref{eq:C15}. Consequently \begin{equation}
\frac{J_k}{g^2J_{k+1}}\longrightarrow
            \frac{2\alpha(3-\alpha)}z.
 \label{eq:C16}
\end{equation} Using the exact identity \eqref{eq:C3}, \begin{equation}
\frac{a_{g,k}}{a_{g,k+1}}\longrightarrow
      \frac{2\pi^2}{z}=\frac{4\pi^2}{t^2}.
 \label{eq:C17}
\end{equation}

\subsection{The top endpoint, including bounded codimension}\label{the-top-endpoint-including-bounded-codimension}

Near the top degree, log-concavity squeezes the ratios between the
interior limit and an exact value in codimension one. The one-kappa pushforward identity gives
\(J_{N-1}=\langle\tau_2\tau_{3g-3}\rangle_g\). Applying DVV at the
degree-two insertion and then the string equation twice gives \[
 15J_{N-1}=
 \{(6g-3)(6g-5)+12g\}J_N.
\] The nonseparating term is
\((1/2)\langle\tau_0^2\tau_{3g-3}\rangle_{g-1}=12gJ_N\); separating
terms vanish by dimension. At \(g=1\), the same identity follows
directly from \(\langle\tau_0^3\rangle_0=1\). Therefore \begin{equation}
q_{N-1}=\frac{12g^2-12g+5}{5},
 \qquad
 \frac{a_{g,N-1}}{a_{g,N}}
 =\frac{2\pi^2(12g^2-12g+5)}{5(6g-3)}.
 \label{eq:C18}
\end{equation}

Let \(m=N-k\ge1\) and \(m=o(g)\). For each fixed \(\alpha<3\),
eventually \(q_{\lfloor\alpha g\rfloor}\le q_{N-m}\le q_{N-1}\) by \eqref{eq:C7}.
As \(\alpha\) increases to three, the limit on the right side of \eqref{eq:C16}
tends to \(12/5\), because \begin{equation}
\frac{t^2\sin t}{\sin t-t\cos t}
       =3-\frac{t^2}{5}+O(t^4)\quad(t\to0).
 \label{eq:C19}
\end{equation} The exact top value \eqref{eq:C18} gives the same limit. Hence \begin{equation}
\frac{q_{N-m}}{g^2}\longrightarrow\frac{12}{5},
 \qquad
 \frac{a_{g,N-m}}{a_{g,N-m+1}}
       \sim\frac{4\pi^2g}{5m}\quad(m=o(g),\ m\ge1).
 \label{eq:C20}
\end{equation} In particular the statement includes every fixed positive
codimension.
For fixed $m$, this leading term agrees, after the length normalization
above, with the fixed-codimension expansion of
Liu--Xu~\cite[Theorem~1.2]{LX}. The same monotonicity argument also covers
unbounded $m=o(g)$ and joins this endpoint to the interior limit.

\subsection{The lowest degrees}\label{the-bottom-endpoint-without-a-relative-error-assumption}

Near degree zero the interior argument gives no information: the
comparison \eqref{eq:C5} determines \(\log b_{g,k}\) only up to an error
\(o(g)\), whereas for \(k=o(g)\) the quantities to be determined vary on a
smaller scale. Two other inputs replace it. For degrees up to a power of \(\log g\), the absolute estimate
\eqref{eq:C8} of Mirzakhani and Zograf gives the ratios directly. For
larger \(k=o(g)\), we argue by contradiction: by the monotonicity of the
ratios \(q_k\), a normalized bracket ratio $\rho_k=b_{g,k}/b_{g,k+1}$ bounded away from one would persist over a range
of degrees comparable to \(k\), and would force \(b_{g,k}\) to drop below
the lower bound \eqref{eq:C23}.

Set \(z_0=\pi^2/2\), and let \(S_h\) be the sum of \(h\) iid variables
with generating function \(H(z_0)/H(z_0w)\). Since \(H(z_0)=3/\pi^2\),
\eqref{eq:C4} becomes the exact identity \begin{equation}
M_{g,k}=\pi^2\Gamma(h)(4\pi^2)^{h-1}
                     \Pr(S_h=N-k).
 \label{eq:C21}
\end{equation} The summands have mean \(3/2\), nonzero variance, span one, and
exponential moments. Uniform exponential tilting in a fixed neighborhood
of that mean gives, for \(0\le k\le c_0g\) with one fixed small
\(c_0>0\), \begin{equation}
\Pr(S_h=N-k)\ge c g^{-1/2}
            \exp[-C(k+1)^2/g].
 \label{eq:C22}
\end{equation} Indeed the local limit theorem for the tilted sum gives the factor
\(g^{-1/2}\), while the rate function has a quadratic zero at the mean;
the displacement here is exactly \(k+1/2\). Combining \eqref{eq:C5}, \eqref{eq:C9}, and
\eqref{eq:C21}--\eqref{eq:C22} yields \begin{equation}
0\le\log\frac{V_{g,1}}{b_{g,k}}
 \le C\log^2(g+1)+C\frac{(k+1)^2}{g}.
 \label{eq:C23}
\end{equation}

Write \[
 \rho_k=\frac{b_{g,k}}{b_{g,k+1}}\ge1,
 \qquad \omega_k=(2k+3)(N-k),
 \qquad \rho_k=\frac{\pi^2q_k}{2\omega_k}.
\] Suppose \(k=o(g)\) but, along a subsequence, \(\rho_k\ge1+\eta\) for
a fixed \(\eta>0\). First consider \(k>\log^3(g+1)\); then \eqref{eq:C23} is
\(o(k)\) at every index comparable to \(k\). Put
\(\delta=\eta/[2(2+\eta)]\). For \(k\le j\le\lfloor(1+\delta)k\rfloor\),
\begin{equation}
\frac{\omega_j}{\omega_k}\le\frac{2j+3}{2k+3}\le1+\delta,
 \qquad
 \rho_j\ge\rho_k\frac{\omega_k}{\omega_j}
            \ge1+\eta/2.
 \label{eq:C24}
\end{equation} These inequalities follow from the monotonicity of the raw ratios $q_j$.
Multiplying the resulting lower bounds over this range forces
\(\log(V_{g,1}/b_{g,\lfloor(1+\delta)k\rfloor+1})\ge c_\eta k\),
contradicting \eqref{eq:C23}.

For the remaining indices \(k\le\log^3(g+1)\), summing the absolute
estimate \eqref{eq:C8} gives \(c_{k+1}=1-o(1)\) uniformly, so division in \eqref{eq:C8} is
legitimate. In fact there is a fixed small \(\delta_0>0\) for which \begin{equation}
\rho_k=1+O((k+1)/g)
       \quad(0\le k\le\delta_0\sqrt g).
 \label{eq:C25}
\end{equation} This proves \begin{equation}
\frac{a_{g,k}}{a_{g,k+1}}\longrightarrow4
                  \quad\text{whenever }k=o(g).
 \label{eq:C26}
\end{equation} No positive lower bound on \(c_k\) was assumed outside the small
square-root window.

\subsection{Combining the degree ranges}\label{uniformity-and-the-corrected-top-density}

\begin{proof}[Proof of the uniform equivalent \eqref{eq:C2}]
Take any sequence of valid indices and a subsequence on which \(k/g\)
converges in \([0,3]\). An interior limit is covered by \eqref{eq:C17}, and a
zero limit by \eqref{eq:C26}, while the comparison in \eqref{eq:C1} also tends to four. If
the limit is three, then \(m=N-k=o(g)\) and \eqref{eq:C20} applies. Equations
\eqref{eq:C11} and \eqref{eq:C19} give \[
 t_{g,k}^2\sim\frac{10m}{2g-1},\qquad
 \frac{4\pi^2}{t_{g,k}^2}\sim\frac{4\pi^2g}{5m}.
\] Thus the true and the proposed ratios are asymptotically equal on that
subsequence too. A sequence violating the supremum statement \eqref{eq:C2} would
contradict one of these three cases.
\end{proof}

The effective density \((2k+1)/(2g-1)\) in \eqref{eq:C1} is essential near the
top: it retains the exact codimension \(m\) in
\(3-\alpha_{g,k}=2m/(2g-1)\). Replacing it by \(k/g\) before taking the
endpoint limit would change the fixed-codimension denominator.

\subsection{Relation to the instanton prediction}\label{subsec:instanton-model}
The denominator \(H\) in \eqref{eq:C4} gives the effective potential of the instanton calculus for ordinary Weil--Petersson volumes. To see the correspondence, choose the positive sign convention
\[
 A_{\rm inst}(w)=\frac{\sin(2\pi w)-2\pi w\cos(2\pi w)}{4\pi^3}
     =\frac23 w^3 H(2\pi^2w^2).
\]
With $x=w^2$, this is the effective potential of~\cite[equation~(3.15)]{EGGLS}; its values at its stationary points are the ordinary Weil--Petersson instanton actions. The ordinary-volume polynomial prediction in~\cite[Section~3.3]{GPRS} and the one-boundary term in~\cite[equation~(4.22)]{JR} contain, up to an overall factor independent of the coefficient degree, the contour expression
\[
 \frac{\Gamma(h)}{L}\oint
       \frac{\mathrm{e}^{Lw}}{w A_{\rm inst}(w)^h}\,\mathrm{d}w,\qquad h=2g-1.
\]
The coefficient of \(L^{2k}\), multiplied by \((2k+1)!\), is proportional to
\[
 \Gamma(h)\operatorname*{Res}_{w=0} w^{2k}A_{\rm inst}(w)^{-h}\,\mathrm{d}w.
\]
Since \(3h-2k-1=2m\), expansion at zero gives the exact identity
\[
 \Gamma(h)\operatorname*{Res}_{w=0} w^{2k}A_{\rm inst}(w)^{-h}\,\mathrm{d}w
 =\Gamma(h)(3/2)^h(2\pi^2)^m \widetilde A_m(h)
 =2\,4^{-k}M_{g,k}.
\]
Thus the model for the coefficients \(a_{g,k}\) has precisely the normalization \(4^{-k}M_{g,k}\), up to the degree-independent overall factor in the prediction. Moreover, the exponent in the \(\mathrm{d}w/w\) representation has its stationary point at
\[
 \frac{2k+1}{h}=\frac{wA_{\rm inst}'(w)}{A_{\rm inst}(w)}
 =\frac{t^2\sin t}{\sin t-t\cos t},\qquad t=2\pi w.
\]
Thus both the comparison model and the saddle equation already appear in the instanton calculus. The estimates \eqref{eq:C5} and the concavity argument above justify the asymptotic formula for adjacent ratios of the true mixed intersection numbers, including both endpoint regimes, without using the contour expression.

\subsection{Cusp weights and the JT spectral edge}\label{subsec:jt-edge}

Saad, Shenker and Stanford related the genus expansion of
Jackiw--Teitelboim (JT) gravity to a matrix integral, with
Weil--Petersson volumes as its geometric input~\cite{SSS}; Stanford and
Witten extended this correspondence to the other random-matrix
ensembles~\cite{StanfordWitten}. See~\cite{MertensTuriaci} for a review.
The parameter $Y$ describes a spectral edge in its defect deformation,
which was introduced by Maxfield--Turiaci and by Witten~\cite{MT,WittenDef}.
We record the normalizations that relate our uniform estimates to that
model. Nothing in this subsection is used in the proofs; its purpose is
to explain why the endpoint $Y=y_0$, at which the constant of
Section~\ref{sec:prefactor} was fixed, is the point where
Painlev\'e~I appears. A defect of opening angle $2\pi\alpha$
is a cusp at $\alpha=0$~\cite[Section~2]{ET}.
Let $\lambda_{\mathrm{def}}$ be its weight, not to be confused with the
cusp density $\lambda=n/g$, and use the convention $2\gamma=1$ of
Maxfield--Turiaci~\cite{MT}. Their string equation at $\alpha=0$
and its branch connected to the undeformed edge are
\[
 \frac{\sqrt{E_0}}{2\pi}I_1(2\pi\sqrt{E_0})
       +\lambda_{\mathrm{def}}=0,\qquad
 E_0=-\frac{Y}{\pi^2},\qquad x=2\pi^2\lambda_{\mathrm{def}}.
\]
The Bessel expressions are interpreted by their entire power series.
Substituting gives exactly $X(Y)=x$. Moreover, the definition
\eqref{eq:H8} gives, for $g\ge2$,
\[
 \sum_{n\ge0}V_{g,n}\frac{\lambda_{\mathrm{def}}^n}{n!}
       =(2\pi^2)^{3g-3}\Phi_g(x).
\]
Thus a fixed cusp count is a coefficient of the defect-weighted
amplitude, whereas $Y$ parametrizes its spectral edge.

The genus-zero density, with its overall topological factor removed,
can be written using our kernel. Specializing the Abel integral in
\cite[equation~(3.50)]{MT} to $\alpha=0$ gives
\[
 \rho_0(E)=\frac1{4\pi}\int_{E_0}^{E}
        \frac{I_0(2\pi\sqrt u)}{\sqrt{E-u}}\,\mathrm du
   =\frac{\sqrt q}{2\pi}K_Y(-2\pi^2q),\qquad q=E-E_0\ge0.
\]
For the second equality set $u=E_0+q(1-s^2)$ and use
\eqref{eq:P1}. At $Y=0$ this is the JT density
$\sinh(2\pi\sqrt E)/(4\pi^2)$. At $Y=y_0$, put
$b=-X''(y_0)>0$. Since $a(y_0)=0$, expansion of $K_{y_0}$ gives
\[
 \rho_0(E_0+q)=\frac{\pi b}{3}q^{3/2}+O(q^{5/2}).
\]
The square-root edge has become a cubic critical edge, the local
geometry associated with the $(2,3)$ critical point, whose string equation
is Painlev\'e~I~\cite{BK,DS,GM}. This gives a
physical interpretation of the endpoint at which the normalized
intersection recursion reduces to pure-$\psi$ intersections.

The exponent of the large-genus action agrees with this critical
geometry. With $\delta=y_0-Y$, the endpoint
expansions of Section~\ref{the-stationary-parameter} give
\[
 x_c-x=\frac b2\delta^2(1+O(\delta)),\qquad
 \mathcal Z=\frac{2\sqrt3\,b}{5}\delta^{5/2}(1+O(\delta))
       =\frac{2^{9/4}\sqrt3}{5b^{1/4}}
          (x_c-x)^{5/4}\bigl(1+O(\sqrt{x_c-x})\bigr).
\]
Theorem~\ref{thm:full-expansion} retains the full action
$\mathcal Z(Y)$ uniformly as $Y\to y_0$. For a fixed marking count,
keeping only its leading critical power would introduce the factor
$\exp(O(g\sqrt{x_c-x}))$, and this replacement need not be uniform in
joint limits. Formal genus counting balances $\hbar$ against
$(x_c-x)^{5/4}$; here, however, $g$ is not the inverse string coupling.
A double-scaled sum over genera would need a separate summation
prescription and estimates beyond the coefficient asymptotics proved here.

\clearpage
\part{Spectral laws of random surfaces}\label{part:spectra}
\phantomsection\label{guide:spectral-laws}
Mirzakhani's integration formula turns volume ratios into expected
counts of short curves. The task in this part is to identify which of
these configurations produce the low eigenvalues and to exclude other
modes. We use the uniform ratios of Theorem~\ref{thm:volume}, the
bounded-boundary factorization of Proposition~\ref{prop:bounded-boundaries},
the positive convolution factor of
Section~\ref{a-positive-convolution-factor} and the all-marking
envelope \eqref{eq:gevrey-envelope}. Huang's trace expansion
\cite{HuangGap} supplies the principal external spectral input; its
precise hypotheses and the extensions needed here are given in
Appendices~\ref{app:critical-trace}, \ref{app:sublinear-trace}
and~\ref{app:moving-edge}.

The volume estimates feed two spectral questions: the law of the
few lowest eigenvalues, and the density of the many eigenvalues in a
fixed-energy window. The two questions require different normalizations
and different modes of convergence; Figure~\ref{fig:spectral-scales}
shows the two operations on the same counting measure. Throughout this
part, eigenvalue counting measures include multiplicities and exclude
the constant mode.
\begin{figure}[htbp]
\centering
\begin{tikzpicture}[
  box/.style={draw=citationblue,thin,rounded corners=2pt,
    text=black,align=center,text width=4.0cm,minimum height=1.95cm,
    inner sep=4pt,font=\small},
  route/.style={-{Latex[length=2mm]},draw=navigationred,line width=.65pt}]
\node[box,text width=3.3cm,minimum height=2.3cm] (spectrum) at (-5.2,0)
  {\textbf{Low spectrum}\\[2pt]
   $\displaystyle\Xi_X=\sum_{0<\lambda_j<1/4}\delta_{\lambda_j}$};
\node[box] (micro) at (0,1.35)
  {\textbf{Rescale energy}\\[2pt]
   $\displaystyle\sum_j\delta_{A_g\lambda_j}$\\[2pt]
   resolves the lowest levels};
\node[box] (fixed) at (0,-1.35)
  {\textbf{Normalize counts}\\[2pt]
   $d_g^{-1}\Xi_X$\\[2pt]
   fixed-energy windows};
\node[box] (poisson) at (5.2,1.35)
  {\textbf{Random limit}\\[2pt]
   $\operatorname{PPP}(2t\,dt)$\\[2pt]
   in distribution};
\node[box] (density) at (5.2,-1.35)
  {\textbf{Deterministic limit}\\[2pt]
   $\nu_*$, $\nu_\lambda$ or $\nu_\infty$\\[2pt]
   in $L^1$ against tests};
\draw[route] (spectrum.east)++(0,.45)--(micro.west);
\draw[route] (spectrum.east)++(0,-.45)--(fixed.west);
\draw[route] (micro.east)--(poisson.west);
\draw[route] (fixed.east)--(density.west);
\end{tikzpicture}
\caption{Two observables above the square-root cusp scale. The upper
branch applies in the ranges of Theorems~\ref{thm:supercritical-spectrum}
and~\ref{thm:dense-persistence}. The lower branch divides by
$d_g=n^2/g$, $g$ or $n$ in the regimes $\sqrt g\ll n\ll g$,
$n/g\to\lambda\in(0,\infty)$ and $g=o(n)$ of
Section~\ref{sec:fixed-energy-spectrum}. A deterministic limit of
normalized counts does not determine the law of the individual lowest
eigenvalues.}
\label{fig:spectral-scales}
\end{figure}

\subsection*{The lowest eigenvalues}
Section~\ref{sec:geometry} supplies the short-curve statistics and
their first variational consequences. At $n\asymp\sqrt g$,
Section~\ref{sec:critical-eigenvalues} identifies the complete spectrum
below $1/4$ by matching funnel trial states with Huang's positive trace.
Section~\ref{sec:moving-edge} refines this description near the quarter
threshold, including fluctuations and their moments.
Above the critical cusp scale and for bounded $n/g$,
Section~\ref{sec:supercritical-spectrum} proves the Poisson law at the
scale $A_g=2\pi^2\sqrt{q_gM_2(x_g)}$. Here $q_g$ is the two-cusp
short-curve intensity and $M_2$ is the area-weighted sum over planar
cusp types; see \eqref{eq:PS1} and~\eqref{eq:super-scale}.
For $n/g\to\infty$, Section~\ref{sec:dense-spectrum} proves that this
law persists in a range containing $n\asymp g^\alpha$ for every
$1<\alpha<11/9$, and that its normalization collapses when
$n\gg g^3/(\log g)^4$. Both ranges are only sufficient, so the
microscopic description has a genuine gap between them. The
residual-type conclusion extends to all sequences $g\to\infty$,
$n/g\to\infty$ (Corollary~\ref{cor:dense-residual}).

\subsection*{Spectral measures in fixed windows}
Let $\Xi_X=\sum_{0<\lambda_j(X)<1/4}\delta_{\lambda_j(X)}$.
Section~\ref{sec:fixed-energy-spectrum} studies $\Xi_X$ without rescaling
energy. Its three main theorems treat the regimes $\sqrt g\ll n\ll g$,
$n/g\to\lambda\in(0,\infty)$ and $g=o(n)$, the last of which includes
fixed genus and genus zero. After division by $n^2/g$, $g$ and $n$,
respectively, the measure converges in $L^1$ against continuous compactly
supported tests on $[0,1/4)$, and hence vaguely in probability. In the
sublinear regime the limit is the funnel pushforward
\eqref{eq:fixed-energy-measure}; at positive density all planar filling
types contribute.
Theorem~\ref{thm:regularized-spectrum} extends the orbital description
to a signed full-energy measure with scattering and threshold
corrections; Corollary~\ref{cor:cusp-density-determinant} is its
consequence for the determinant. The proofs of these results are
collected in Appendix~\ref{app:fixed-window-proofs}.

The distinction between the two kinds of limit is essential. A
fixed-window measure records the proportion of eigenvalues in an energy
range, and this proportion can be unchanged by a vanishing fraction of
much smaller eigenvalues; the microscopic process resolves those lowest
levels individually. The fixed-window theorems therefore coexist with
the unresolved dense microscopic law.

\subsection*{Previous work on the spectral transition}
This subsection places the results of this part next to earlier work,
in the order of the phase diagram of Table~\ref{tab:intro-phases}:
closed surfaces and the subcritical range $n=o(\sqrt g)$, the critical
scale $n\asymp\sqrt g$, the supercritical and dense ranges, counts and
spectral measures in fixed windows, and the tight length spectrum. The
statements themselves are in Sections~\ref{sec:geometry}
to~\ref{sec:fixed-energy-spectrum}; only the comparison is made here.

Random hyperbolic surfaces were first studied in a combinatorial
model, by Brooks and Makover~\cite{BrooksMakover}; the Weil--Petersson
model was first studied by Guth, Parlier and Young~\cite{GPY} and
by Mirzakhani~\cite{Mir}.
For Weil--Petersson random closed surfaces, Mirzakhani established a
positive spectral gap with high probability~\cite[Section~4.5]{Mir}.
Wu and Xue~\cite{WX} and, independently, Lipnowski and
Wright~\cite{LipWright} raised the lower bound to $3/16-\varepsilon$.
Anantharaman and Monk obtained $2/9-\varepsilon$, first in the original
arXiv version of the first part of their Friedman--Ramanujan
series~\cite{AMPartOne} and now in a separate
article~\cite[Theorem~1.1]{AMTwoNinths}, and then
$1/4-\varepsilon$~\cite[Theorem~1.6]{AMGap};
Hide, Macera and Thomas gave a polynomial rate of approach to $1/4$
\cite{HMT}. In the different model of uniformly random covers of a
fixed surface, Magee, Naud and Puder obtained the relative spectral gap
$3/16-\varepsilon$ for covers of a compact surface~\cite{MNP}, and
Magee, Puder and van Handel proved strong convergence of the
corresponding random permutation representations, which yields the
near-optimal relative gap $1/4-\varepsilon$~\cite{MPvH}. For
covers of a fixed cusped surface, Hide and Magee proved that no new
eigenvalues appear below $1/4-\varepsilon$, and deduced that there are
closed surfaces of genus tending to infinity with first eigenvalue
tending to $1/4$~\cite{HideMagee}. Growing cusp counts introduce the
separating pants mechanism studied by Shen and Wu~\cite{SW}.
See~\cite{MonkNaud,Petri} for surveys of spectral gaps of random
hyperbolic surfaces.

The transition occurs at the scale $\sqrt g$. Below it, Hide
proved~\cite[Theorem~1.3]{Hide} that for $n=O(g^\alpha)$ with
$0\le\alpha<1/2$ there is, with probability tending to one, no nonzero
eigenvalue below $1/4-((2\alpha+1)/4)^2-\varepsilon$. Shen and Wu
proved a uniform positive lower bound for the spectral gap throughout
the range $n=o(\sqrt g)$~\cite[Theorem~3]{SW}, and He, Wu and Xue
improved the bound for $n=O(g^\alpha)$ to
$1/4-(6(1-\alpha))^{-2}-\varepsilon$, which tends to $5/36-\varepsilon$
as $\alpha\to1/2$~\cite[Theorem~1.2]{HWX}. Huang subsequently proved
convergence to $1/4$ throughout $n=o(\sqrt g)$, together with matching
expected-count bounds and the necessity of this cusp
condition~\cite[Theorem~1.1, Proposition~6.14
and Corollary~6.15]{HuangGap}. We use his canonical trace coefficients and
positive exceptional measure; these subcritical results are not new
here.

When $n/\sqrt g\to c\in(0,\infty)$, Shen and Wu~\cite[Theorem~2]{SW} describe the
small-value distribution of the Cheeger constant, not a limiting law
for Laplace eigenvalues; for closed surfaces of large genus, Budzinski,
Curien and Petri bounded the Cheeger constant asymptotically by
$2/\pi$~\cite{BCPCheeger}. Theorem~\ref{thm:critical-eigenvalues}
determines the limiting Laplace-eigenvalue process. Its ingredients are
the funnel model of Section~\ref{sec:critical-eigenvalues}, spectral
regularization through pinching~\cite{Schulze,BJP}, and saturation of
the expected count by localized ground states. The spectral theory of
degenerating hyperbolic surfaces of finite area goes back to Hejhal,
Wolpert and Ji~\cite{HejhalDegen,WolpertSpecI,WolpertSpecII,JiGreen}
and, for the eigenvalues below $1/4$, to Colbois and
Courtois~\cite{ColboisCourtois};
Schulze's results cover surfaces with funnels, which is the case
needed for the model $F_s$. The moving-window
Theorem~\ref{thm:moving-edge} makes this comparison uniform for growing
waists and energies approaching $1/4$; in the stated near-critical
range it identifies the scale $n^2/g$ and the limiting inverse-square
Poisson intensity. Corollary~\ref{cor:edge-moments} explains why rare
surfaces give a logarithmic enhancement of the mean gap defect.
Above the critical scale the spectral gap vanishes in probability
(Corollary~\ref{cor:spectral}). This answers negatively the question
of Shen and Wu whether a uniform spectral gap persists when
$n/\sqrt g\to\infty$~\cite[Question~2]{SW}. Shen--Wu proved the
vanishing for $\sqrt g\ll n\ll g$ and observed that Zograf's
deterministic bound~\cite{ZogSpectral} covers $n/g\to\infty$; the
adjacent-volume limit supplies the remaining range $n\asymp g$.
The statement concerns the Weil--Petersson distribution: Guo, Hua
and Shen construct surfaces with $n/g$ tending to any prescribed
positive limit and a uniform spectral gap~\cite{GHS}.

Hide and Thomas obtained Poisson limits and small eigenvalues at fixed
genus~\cite{HTSphere,HTLargeN}; their deterministic linear count below
$1/4$~\cite[Theorem~1.3]{HTMany} applies whenever $g<an$, for each fixed $a>0$,
and is linear-optimal by the sharp topological upper bound $2g+n-2$ of
Otal and Rosas~\cite{OtalRosas}. Below
every fixed positive threshold, a deterministic linear total count for
$g=o(n)$ is expected to follow from the conformal capacitor bound
of~\cite{Hassannezhad} (Remark~\ref{rev:rem-hassannezhad});
Corollary~\ref{cor:many-spectrum} covers the full range $g\le an$ in
Weil--Petersson probability. Closed surfaces behave differently: Monk
proved that, with probability tending to one, a Weil--Petersson random
closed surface has $O(g(\log g)^{-3/4})$ eigenvalues below
$1/4$~\cite[Theorem~4]{Monk}. At high energy, Rudnick showed that the
Weil--Petersson variance of a smooth linear statistic counting
eigenvalues of closed surfaces in a short energy window has the
Gaussian orthogonal ensemble value in the double limit, first of large
genus and then of high energy~\cite{Rudnick}; Rudnick and Wigman proved
the corresponding central limit theorem~\cite{RudnickWigman}. In the range $\sqrt g\ll n\ll g$, the
second-moment volume inputs for our short-curve and small-eigenvalue
results already follow from~\cite{SW}. The uniform volume and boundary
theorems supply a single intensity across all cusp densities and the
growing-genus extension at positive and infinite density.

Gatignol and Thomas~\cite{GT} obtain primitive-geodesic Poisson
limits in the range $\sqrt g\ll n\ll g$ and, at the critical scale,
for lengths below $2\operatorname{arcsinh}1$ (their
Theorems~1.2--1.3). In the former range they also obtain
fixed-threshold lower eigenvalue counts of order
$n^2/g$~\cite[Theorem~1.10]{GT}. These results overlap our short-curve
and counting consequences at zero cusp density. Our critical length
process is unrestricted in length and retains the pants clusters, and
Theorem~\ref{thm:critical-eigenvalues} concerns the full subthreshold
Laplace spectrum rather than a short-length window or a lower
eigenvalue count.

Budd and Lions studied the tight length spectrum in a near-critical
Boltzmann model, in which the cusp count is random~\cite{BL}, and asked
whether their limit persists for other cusp weights and under the
fixed-count measures~\cite[Open Question~1.6]{BL}. Recovering a
prescribed cusp count requires a local coefficient estimate, not merely
concentration of the random count. Theorem~\ref{thm:tight-spectrum}
supplies the fixed-count extension: tight volume estimates and a
positive Fourier bound allow coefficient extraction at every cusp
density. Corollary~\ref{cor:tight-boltzmann} then treats every
subcritical cusp weight, including arbitrary approach to the critical
endpoint. The planar separating curves of
Theorem~\ref{thm:pants-poisson} lie outside the tight spectrum
(Subsection~\ref{subsec:tight-spectrum}), so the two length processes
describe disjoint families of curves.

\section{Short curves and small eigenvalues}\label{sec:geometry}

This section converts the volume theory of Part~\ref{part:volumes}
into statistics of short curves and into the eigenvalue upper bounds
they imply. Through Mirzakhani's integration formula~\cite{MirInv},
expected numbers of short curves are ratios of volumes of cut
surfaces, and these ratios are supplied by Theorem~\ref{thm:volume}
and, for fixed boundary lengths, by
Proposition~\ref{prop:bounded-boundaries}. We first apply the
adjacent-volume limit to Shen--Wu's vanishing-gap
argument~\cite{SW} (Corollary~\ref{cor:spectral},
Subsection~\ref{sec:spectral}). We then determine the short-curve
process and the systole when $n/\sqrt g\to\infty$ and $n/g$ has a
limit (Theorem~\ref{thm:pants-poisson},
Subsection~\ref{sec:pants-spectrum}), deduce many small eigenvalues
(Corollary~\ref{cor:many-spectrum},
Subsection~\ref{sec:pants-eigenvalues}) and locate them in the
residual spectrum (Theorem~\ref{thm:many-residual},
Subsection~\ref{sec:residual-count}). The last two subsections state
the limits of the tight length spectrum at every cusp density
(Theorem~\ref{thm:tight-spectrum}) and of the full primitive length
spectrum at the critical scale $n/\sqrt g\to c\in(0,\infty)$
(Theorem~\ref{thm:critical-cluster}); the longer proofs are in
Appendices~\ref{app:short-curves}, \ref{app:tight}
and~\ref{app:critical}. On the spectral side, everything here gives
eigenvalue upper bounds or lower counts; the matching upper counts,
which show that the short curves exhaust the low spectrum, are the
subject of Sections~\ref{sec:critical-eigenvalues}
to~\ref{sec:fixed-energy-spectrum}.

\subsection{The first positive eigenvalue}\label{sec:spectral}

The simplest spectral consequence, Corollary~\ref{cor:spectral}, uses
only the discrete cusp ratio $D_{g,n-1}$ of Theorem~\ref{thm:volume}
and holds throughout $n/\sqrt g\to\infty$. A short curve cutting off
two cusps bounds a two-cusp pair of pants of fixed area $2\pi$,
attached to the rest of the surface along a short boundary; the Cheeger
ratio, and hence the Rayleigh quotient, is therefore small. Shen--Wu's
second-moment argument makes such a curve likely. Shen and Wu asked for
the adjacent-volume limit needed to extend this argument beyond
$n=o(g)$~\cite[Question~3]{SW}, and Theorem~\ref{thm:volume} provides
that limit along every sequence with $g\to\infty$.

\begin{mainresult}
\begin{corollary}[Vanishing spectral gap]\label{cor:spectral}
Suppose that $g\to\infty$ and that the integer cusp counts satisfy $n(g)/\sqrt g\to\infty$. Then, for every $\varepsilon>0$,
\[
 \mathbb P^{\mathrm{WP}}_{g,n(g)}
 \{X:\Delta_X\text{ has an eigenvalue in }(0,\varepsilon)\}\longrightarrow1.
\]
\end{corollary}
\end{mainresult}

\begin{proof}[Proof of Corollary~\ref{cor:spectral}]
Put $n=n(g)$ and
\[
 L_g=\left(\frac{\sqrt g}{n}\right)^{1/2}.
\]
Let $N_g(X)$ count the simple closed geodesics of length at most $L_g$ whose complement has components of topological types $(0,3)$ and $(g,n-1)$. Proposition~35 of~\cite{SW}, in the notation of its arXiv version~4, states under the sole assumption $n/\sqrt g\to\infty$ that
\[
 \liminf_{g\to\infty}\mathbb P^{\mathrm{WP}}_{g,n}\{N_g\ge1\}
 \ge \liminf_{g\to\infty}
 \frac{V_{g,n-1}^{2}}{V_{g,n}V_{g,n-2}}.
\]
The ratio on the right is $D_{g,n-1}$, which tends to one by Theorem~\ref{thm:volume}. Thus $\mathbb P^{\mathrm{WP}}_{g,n}\{N_g\ge1\}\to1$. The same conclusion also follows within this paper: $N_g$ is the count $N_{L_g}$ of Subsection~\ref{sec:pants-spectrum}, with $q_g$ as in \eqref{eq:PS1}, and since $L_g\to0$ and $q_gL_g^2\asymp n\sqrt g/(g+n)\to\infty$, Proposition~\ref{prop:pants-concentration} gives $N_g/(q_gL_g^2)\to1$ in probability. Its proof is the same second-moment argument, with $D_{g,n-1}\to1$ as the volume input.

Write $\mathfrak h(X)$ for the Cheeger constant. On the event $N_g\ge1$, a counted geodesic bounds a hyperbolic pair of pants of area $2\pi$; the other component has area $2\pi(2g+n-3)\ge2\pi$ for all sufficiently large $g$. This separating curve is an admissible Cheeger cut, so
\[
 \mathfrak h(X)\le \frac{L_g}{2\pi}.
\]
Let $\mathcal R(X)$ be the infimum of the Rayleigh quotient over nonzero $H^1(X)$ functions orthogonal to constants. Buser's inequality~\cite[Theorem~7.1]{Buser}, in the form recorded in~\cite[Theorem~11]{SW} for complete noncompact hyperbolic surfaces, gives $\mathcal R(X)\le c\,\mathfrak h(X)$ with an absolute constant $c$; no smallness of $\mathfrak h$ is needed, because the Cheeger constant of a finite-area noncompact hyperbolic surface is bounded, so the quadratic term of the compact inequality is absorbed into the linear one. Since $L_g\to0$, this upper bound is eventually smaller than both $\varepsilon$ and $1/4$. The essential spectrum starts at $1/4$; the min--max principle, as recalled in~\cite[Theorem~7]{SW}, therefore yields a positive discrete eigenvalue $\lambda_1(X)=\mathcal R(X)<\varepsilon$ on $N_g\ge1$, for all sufficiently large $g$. This proves the asserted probability limit.
\end{proof}

The use of $D_{g,n-1}$ above is the only change to the final step of
Shen--Wu's proof; under their condition $n=o(g)$ this limit follows
from \cite[Lemma~5.1]{MZ}, as they note. The short-curve estimate, the
Cheeger cut and the passage to a discrete eigenvalue are their
geometric and spectral inputs. Corollary~\ref{cor:spectral} is the
supercritical half of the gap phase diagram
(Corollary~\ref{cor:gap-phase-diagram}); the rate at which the gap
vanishes, and the law of the lowest eigenvalues, are determined in
Section~\ref{sec:supercritical-spectrum}.

\subsection{Short-curve processes and the systole}
\label{sec:pants-spectrum}

This subsection determines the limiting process of short separating
curves and the systole above the critical scale
(Theorem~\ref{thm:pants-poisson}). The input is Mirzakhani's
integration formula~\eqref{eq:PS3}, whose volume ratios are evaluated
by Theorem~\ref{thm:volume} in the form~\eqref{eq:PS5}; the
uniformity in the cusp density comes from the coefficient-ratio bounds
of Section~\ref{extracting-fixed-cusp-counts}, and the proof is in
Appendix~\ref{app:short-curves}. The moment estimates are reused in
Subsection~\ref{sec:pants-eigenvalues} to count many disjoint short
pairs of pants. Throughout this and the next subsection,
$g\to\infty$, $n/\sqrt g\to\infty$, and $X$ has law
$\mathbb P^{\mathrm{WP}}_{g,n}$. Let $\mathcal P_2(X)$ be the set of
unoriented simple closed geodesics bounding a two-cusp pair of pants.
Define the deterministic intensity parameter
\begin{equation}
 q_g=\frac{n(n-1)}{8g\mathcal C(n/g)}.
 \label{eq:PS1}
\end{equation}
The profile endpoints in Theorem~\ref{thm:volume} imply
$q_g\asymp n^2/(g+n)\to\infty$. In particular
$q_g\sim n^2/(32\pi^2g)$ when $n/g\to0$, and
$q_g\sim x_cn/(8\pi^2)$ when $n/g\to\infty$.

We write $\operatorname{PPP}(\nu)$ for a Poisson point process with intensity measure $\nu$.

\begin{theorem}[Short-curve processes and the systole]\label{thm:pants-poisson}
Suppose $g\to\infty$ and $n/\sqrt g\to\infty$, and let $X$ have law
$\mathbb P^{\mathrm{WP}}_{g,n}$.

\noindent\emph{(i) Curves cutting off two cusps.} The random locally finite counting measures
\begin{equation}
 \sum_{\gamma\in\mathcal P_2(X)}
       \delta_{\sqrt{q_g}\ell_\gamma}
 \ \Longrightarrow\ \operatorname{PPP}(2t\,\mathrm{d}t)
 \quad\hbox{on }[0,\infty)
 \label{eq:PS2}
\end{equation}
converge in distribution for the vague topology. Thus counts in finitely
many disjoint bounded intervals $[a_i,b_i]$ converge jointly to independent
Poisson variables with means $b_i^2-a_i^2$. If $\ell_{\min}^{(2)}$ is the
shortest length in $\mathcal P_2(X)$, with value infinity when this set is
empty, then
\[
 \mathbb P^{\mathrm{WP}}_{g,n}
 \{\sqrt{q_g}\ell_{\min}^{(2)}>t\}\longrightarrow \mathrm{e}^{-t^2}
 \qquad(t\ge0).
\]
\smallskip
\noindent\emph{(ii) All planar cusp types.} Let $\mathcal P_k(X)$ consist of simple geodesics whose
genus-zero side contains exactly $k$ cusps. If $n/g\to\lambda\in[0,\infty]$,
put $x_\lambda=X(Y_\lambda)$, where $\lambda(Y_\lambda)=\lambda$, with
$x_0=0$ and $x_\infty=x_c$. Define
\begin{equation}
 w_k(x)=\frac{2V_{0,k+1}}{k!}
             \left(\frac{x}{2\pi^2}\right)^{k-2},\qquad
 W(x)=\sum_{k\ge2}w_k(x)=\frac{2\Phi_0'(x)}{x^2}.
 \label{eq:PC1}
\end{equation}
Here $w_2(0)=W(0)=1$, $w_k(0)=0$ for $k>2$, and $W(x_c)<\infty$.
Then
\begin{equation}
 \sum_{k\ge2}\ \sum_{\gamma\in\mathcal P_k(X)}
       \delta_{(\sqrt{q_g}\ell_\gamma,k)}
 \ \Longrightarrow\
 \operatorname{PPP}\!\left(2t\,\mathrm{d}t\otimes
                    \sum_{k\ge2}w_k(x_\lambda)\delta_k\right).
 \label{eq:PC2}
\end{equation}
Convergence is for the vague topology on the product with the discrete
type space. The projection to lengths also converges, to
$\operatorname{PPP}(2W(x_\lambda)t\,\mathrm{d}t)$.

\noindent\emph{(iii) Primitive geodesics and the systole.} Under the density-limit assumption of~(ii), its projected process is also the limit of
the process of all unoriented primitive closed geodesics, at the same
rescaling. Writing $\operatorname{sys}(X)$ for the shortest closed
geodesic length, we have
\[
 \mathbb P^{\mathrm{WP}}_{g,n}
 \{\sqrt{q_g}\operatorname{sys}(X)>t\}
 \longrightarrow \mathrm{e}^{-W(x_\lambda)t^2}\qquad(t\ge0).
\]
With probability tending to one every systolic curve is planar
separating. Select one by any measurable tie rule, and let $K_g$ be
its planar-side cusp count, with value zero if it is not planar. Then
$\bigl(\sqrt{q_g}\operatorname{sys}(X),K_g\bigr)$ converges jointly
to a pair $(Z,K)$ with
\begin{equation}
 \begin{aligned}
 \mathbb P\{Z\in \mathrm{d}t,K=k\}
 &=2t\,w_k(x_\lambda)\mathrm{e}^{-W(x_\lambda)t^2}\,\mathrm{d}t,\\
 \mathbb P\{K=k\}&=\frac{w_k(x_\lambda)}{W(x_\lambda)}
 \qquad(t>0,\ k\ge2).
 \end{aligned}
 \label{eq:FS1}
\end{equation}
\end{theorem}

The method of factorial moments for curves separating off a prescribed
number of cusps has a direct precedent in
Hide--Thomas~\cite[Theorem~1.7]{HTSphere} for punctured spheres, where
summed types and nested configurations are also
treated~\cite[Remark~1.9 and Proposition~3.6]{HTSphere}. Here the
uniform volume ratios give the entire family of type intensities when
$g\to\infty$, $n/\sqrt g\to\infty$ and $n/g\to\lambda\in[0,\infty]$.
These planar curves lie outside the tight spectrum
of~\cite[Definition~2.1]{BL}, which is treated separately in
Subsection~\ref{subsec:tight-spectrum}. The all-primitive process of
Mirzakhani--Petri~\cite[Theorem~4.1]{MP} concerns closed surfaces at
fixed physical lengths, whereas the full process here is observed on a
shrinking scale; a tight-spectrum limit at its longer natural scale is
a separate question. The same Poisson law, with the Mirzakhani--Petri
intensity, also governs short cycles on random unicellular maps and
random metric maps of large genus~\cite{JansonLouf,BGL}.

\emph{Factorial moments.} For sufficiently small $L$, the collar lemma
forces all curves counted by
$N_L=\#\{\gamma\in\mathcal P_2(X):\ell_\gamma\le L\}$ to be disjoint;
for closed surfaces, the tangle-free hypothesis of Monk and
Thomas~\cite{MonkThomas} gives embedded cylinders around short
geodesics up to lengths of order $\log g$, which is not needed here.
Their pants interiors are also disjoint: otherwise one pants boundary
would be an essential geodesic in the interior of the other pants,
hence boundary-parallel, contradicting uniqueness of geodesic
representatives. Thus their pairs of cusp labels are disjoint, and the
complement of $r$ such pants has type $(g,n-r)$. Since $g\ge2$
eventually, the exceptional double-pants decomposition of a
four-punctured sphere does not arise.

Mirzakhani's integration formula~\cite{MirInv,Mir}, with the same conventions as
\cite[Lemma~29]{SW}, gives, for any fixed $r$ and short length windows
$I_1,\ldots,I_r$, the ordered distinct-curve integral
\begin{equation}
 \mathbb E\!\sum_{\substack{\gamma_1,\ldots,\gamma_r\in\mathcal P_2(X)\\
                            \text{distinct}}}
       \prod_{j=1}^r\mathbf1_{\ell_{\gamma_j}\in I_j}
 =\frac{(n)_{2r}}{2^r V_{g,n}}
    \int_{I_1\times\cdots\times I_r}
       V_{g,n-r}(x_1,\ldots,x_r,0^{n-2r})
       \prod_{j=1}^r x_j\,\mathrm{d}x_j.
 \label{eq:PS3}
\end{equation}
Here $(n)_{2r}=n(n-1)\cdots(n-2r+1)$. The factor $2^{-r}$ forgets
the order of the two labels within each pair; since the tuple of curves
is ordered, there is no $r!$ divisor. Each pants volume is $V_{0,3}=1$.
Mirzakhani's boundary estimate~\cite[Lemma~3.2]{Mir}, in the form
of~\cite[Lemma~14]{SW},
\begin{equation}
 1\le\frac{V_{g,m}(\mathbf x)}{V_{g,m}}
       \le\exp\!\left(\frac12\sum_jx_j\right)
 \label{eq:PS4}
\end{equation}
is uniform in the topological type. For fixed $r$, it replaces the
boundary volume in \eqref{eq:PS3} by $V_{g,n-r}(1+O_r(L))$ when
all $I_j\subset[0,L]$.

Theorem~\ref{thm:volume} also implies, for every fixed $r$,
\begin{equation}
 \frac{V_{g,n-r}}{V_{g,n}}
  =(2g\mathcal C(n/g))^{-r}(1+o(1)).
 \label{eq:PS5}
\end{equation}
Indeed, the left side is a product of $r$ adjacent cusp ratios.
Replacing $(n-j)/g$ by $n/g$, for $1\le j\le r$, changes each profile
factor by $1+o(1)$: along any subsequence with
$n/g\to\lambda_*\in[0,\infty]$ this follows from the continuity and
positivity of $\mathcal C$ when $\lambda_*<\infty$, and from
$\mathcal C(\lambda)\sim\pi^2\lambda/x_c$ together with $(n-j)/n\to1$
when $\lambda_*=\infty$.

\smallskip\noindent\emph{Idea of the proof.} The proof, given in
Appendix~\ref{app:short-curves}, computes mixed factorial moments from
\eqref{eq:PS3}. Three inputs make the computation uniform in the cusp
density.
First, a coefficient-ratio bound obtained from the positive convolution
factor of Section~\ref{extracting-fixed-cusp-counts} shows that curves
cutting off many cusps contribute a summable tail, uniformly along the
sequence. Second, for finitely many fixed types the disjoint
(\emph{star}) configurations give the Poisson moments
$\prod_iw_{k_i}(x_\lambda)(b_i^2-a_i^2)$, while nested configurations
lose a factor $n^{-1}$ by a count of the root degree in the cut tree.
Third, the same positive factor bounds the total contribution of all
separating curves with positive genus on both sides by $O(1/g)$, and
nonseparating curves are rare on the scale $q_g^{-1/2}$; thus the
systole is planar separating with probability tending to one.

The intensity $q_g$ and the type weights $w_k$ are carried through the
rest of this part: they define the microscopic spectral scale
\eqref{eq:super-scale} of Section~\ref{sec:supercritical-spectrum},
the type process \eqref{eq:PC2} supplies the limiting rates in the
proof of Theorem~\ref{thm:supercritical-spectrum}, and the same
weights $w_k$ govern the cap process of
Section~\ref{sec:dense-spectrum} at growing density.

\subsection{Many small eigenvalues}\label{sec:pants-eigenvalues}

This subsection turns the short-curve counts into eigenvalue bounds.
Proposition~\ref{prop:pants-concentration} is a law of large numbers
for the number $N_L$ of short two-cusp pairs of pants, both at
thresholds $L\to0$, where \eqref{eq:PS3} and the boundary estimate
\eqref{eq:PS4} suffice, and at fixed $L$, where the boundary
factorization of Proposition~\ref{prop:bounded-boundaries} is needed.
Corollary~\ref{cor:many-spectrum} converts the count into many small
eigenvalues by half-collar trial functions, throughout
$n/\sqrt g\to\infty$ and, below every fixed threshold, throughout
$g\le an$. These are lower counts; the matching upper counts and the
limiting densities are proved in
Section~\ref{sec:fixed-energy-spectrum}.

\begin{proposition}[Many short pants]\label{prop:pants-concentration}
If $L=L_g\to0$ and $q_gL_g^2\to\infty$, then
\begin{equation}
 \frac{N_{L_g}}{q_gL_g^2}\longrightarrow1
 \quad\hbox{in probability}.
 \label{eq:PS6}
\end{equation}
Under these assumptions, for each fixed $\varepsilon\in(0,1)$, with probability
tending to one there are at least
$\lfloor(1-\varepsilon)q_gL_g^2\rfloor-1$ positive discrete eigenvalues
in $(0,6L_g]$, counted with multiplicity.
For each fixed $0<L<2\operatorname{arcsinh}1$, one also has
\begin{equation}\label{eq:fixed-pants-concentration}
 \frac{N_L}{q_g\mathfrak B_{Y_{n/g}}(L)}\longrightarrow1
 \quad\hbox{in probability},\qquad
 \mathfrak B_Y(L)=2\int_0^L t\mathfrak b_Y(t)\,\mathrm dt,
\end{equation}
where $\mathfrak b_Y$ is the boundary profile in
\eqref{eq:boundary-profile}.
\end{proposition}
\begin{proof}
The first two cases of \eqref{eq:PS3} give
\[
 m_L:=\mathbb E N_L=\frac{(n)_2L^2}{4}
       \frac{V_{g,n-1}}{V_{g,n}}(1+O(L))\sim q_gL^2,
\]
and
\[
 \frac{\mathbb E(N_L)_2}{m_L^2}
 =\frac{(n-2)(n-3)}{n(n-1)}D_{g,n-1}^{-1}(1+O(L))\longrightarrow1.
\]
Since $m_L\to\infty$, the relative variance tends to zero, proving
\eqref{eq:PS6}. This second-moment argument does not use the point-process
limit at a window growing with $g$.

Hide--Thomas~\cite[Lemma~5.3]{HTSphere} prove that $N$ such pants
boundaries of length at most $L$ give $N-1$ positive eigenvalues below
$6L$, provided $6L<1/4$. For the reader's convenience we recall the
variational argument.
Inside each pants choose a function equal to inward distance $\rho$ on
the half-collar $0\le\rho\le1$, equal to one beyond it, and zero outside
the pants. For all sufficiently small boundary lengths $\ell$, the
collar lemma supplies this embedded half-collar. The function lies in
$H^1(X)$, and
\[
 \int|\nabla f|^2=\ell\sinh1,\qquad
 \int|f|^2\ge2\pi-\ell\sinh1.
\]
Different pants have disjoint supports up to their boundaries, so the
span of these functions has dimension $N$ and Rayleigh quotient at most
$6L$ for small $L$. Its intersection with the orthogonal complement of
constants has dimension at least $N-1$. Since the essential spectrum
starts at $1/4$, the form version of min--max gives the claimed positive
discrete eigenvalues. Apply \eqref{eq:PS6} and let $g$ grow.

For fixed $L$, the $O(L)$ boundary error above is insufficient.
Instead apply Proposition~\ref{prop:bounded-boundaries} with one
and two nonzero boundaries to \eqref{eq:PS3}. Uniformly on $[0,L]^r$,
for $r=1,2$,
\[
 \frac{V_{g,n-r}(t_1,\ldots,t_r,0^{n-2r})}{V_{g,n-r}}
 =(1+o(1))\prod_{i=1}^r\mathfrak b_{Y_{n/g}}(t_i).
\]
The total end count here is $n-r$. Replacing $Y_{(n-r)/g}$ by
$Y_{n/g}$ costs $o(1)$ uniformly: both have the same subsequential
limit, including either endpoint of $[0,y_0]$, and the boundary
profile is uniformly continuous on the fixed length box.
Together with \eqref{eq:PS5}, this yields
\[
 \mathbb EN_L=(1+o(1))q_g\mathfrak B_{Y_{n/g}}(L),\qquad
 \mathbb E(N_L)_2=(1+o(1))q_g^2\mathfrak B_{Y_{n/g}}(L)^2.
\]
The combinatorial ratio for the second identity is
$(n-2)(n-3)/(n(n-1))\to1$. Since $\mathfrak b_Y\ge1$,
the mean is at least $(1+o(1))q_gL^2\to\infty$.
The relative variance tends to zero, proving
\eqref{eq:fixed-pants-concentration} without a density-limit assumption.
\end{proof}

\begin{corollary}[Small eigenvalues at shrinking and fixed thresholds]\label{cor:many-spectrum}
Suppose $g\to\infty$ and $n/\sqrt g\to\infty$. If $j=j(g)\to\infty$
and $j=o(q_g)$, then
\[
 \mathbb P^{\mathrm{WP}}_{g,n}
 \{\lambda_j(X)\text{ exists and }\lambda_j(X)\le12\sqrt{j/q_g}\}
 \longrightarrow1,
\]
where positive discrete eigenvalues are ordered with multiplicity.
For each fixed $0<\varepsilon<1/4$, let
$Z_\varepsilon(X)=\#\{j:0<\lambda_j(X)<\varepsilon\}$, with multiplicity.
Then also
\begin{equation}\label{eq:fixed-low-count}
 \mathbb P^{\mathrm{WP}}_{g,n}
 \left\{Z_\varepsilon(X)\ge\frac{\varepsilon^2}{576}q_g\right\}
 \longrightarrow1.
\end{equation}
Moreover, for every $a>0$ and $0<\varepsilon<1/4$ there is
$c_{\varepsilon,a}>0$ such that, along every sequence
$n\to\infty$ and $0\le g\le an$,
\begin{equation}\label{eq:linear-low-count}
 \mathbb P^{\mathrm{WP}}_{g,n}
 \{Z_\varepsilon(X)\ge c_{\varepsilon,a}(2g+n-2)\}
 \longrightarrow1.
\end{equation}
\end{corollary}

\begin{proof}[Proof of Corollary~\ref{cor:many-spectrum}]
For $j\to\infty$ with $j=o(q_g)$, choose $L_g=2\sqrt{j/q_g}$.
Then $L_g\to0$, $q_gL_g^2=4j$, and \eqref{eq:PS6} implies
$N_{L_g}\ge j+1$ with probability tending to one. The preceding
variational estimate gives $\lambda_j\le6L_g=12\sqrt{j/q_g}$.

For the fixed threshold choose $L=\varepsilon/12<1/48$.
The half-collar trial functions above have Rayleigh quotient at most
$6L=\varepsilon/2<\varepsilon$, so $Z_\varepsilon\ge N_L-1$.
By \eqref{eq:fixed-pants-concentration} and
$\mathfrak B_Y(L)\ge L^2$, with probability tending to one
$N_L\ge q_gL^2/2$. Since $q_gL^2\to\infty$, eventually
$N_L-1\ge q_gL^2/4=\varepsilon^2q_g/576$.
This proves \eqref{eq:fixed-low-count}.

If $g\to\infty$ and $g\le an$, put $\lambda=n/g\ge1/a$.
Uniformly in this range,
\[
 \frac{q_g}{2g+n-2}
 =(1+o(1))\frac{\lambda^2}{8\mathcal C(\lambda)(2+\lambda)}.
\]
The last expression is positive and continuous on $[1/a,\infty)$,
with positive limit $x_c/(8\pi^2)$ at infinity. Its infimum is
therefore positive, and \eqref{eq:fixed-low-count} implies
\eqref{eq:linear-low-count} in the growing-genus case.
For fixed genus, Hide--Thomas~\cite[Theorem~1.1]{HTLargeN}
give $C(\varepsilon)n$ eigenvalues below $\varepsilon$ with probability
tending to one, with $C(\varepsilon)>0$ independent of the fixed genus.
The order of quantifiers in \cite[Theorem~1.1]{HTLargeN} is the one
needed here: the constant is chosen first, for every $\varepsilon>0$,
and the statement then holds for every genus $g\ge0$ as $n\to\infty$.
A constant depending on the genus would not suffice, since the
subsequence argument below must produce a single $c_{\varepsilon,a}$
for all bounded genera.
Removing the zero eigenvalue and reducing the constant gives
\eqref{eq:linear-low-count} in that case as well. Choose
$c_{\varepsilon,a}$ smaller than both constants. Every subsequence
of a general genus sequence has a further subsequence with $g\to\infty$
or $g$ fixed, completing the proof.
\end{proof}

The spectral mechanism is established in the fixed-genus work of
Hide--Thomas~\cite[Theorem~1.10 and Remark~5.5]{HTSphere}; their later
bordered-volume asymptotics give linear counts below each fixed positive
threshold at fixed genus~\cite[Theorem~1.1]{HTLargeN}, and in the
range $\sqrt g\ll n\ll g$ the volume inputs of Shen--Wu~\cite{SW} and
the count of Gatignol--Thomas~\cite[Theorem~1.10]{GT} already give
the order $n^2/g$. Here the pants mechanism is combined with the
uniform boundary factorization to cover positive and unbounded cusp
density. This argument yields neither a matching spectral upper count
nor a limiting spectral density; both are obtained in
Section~\ref{sec:fixed-energy-spectrum}, where the count of order
$n^2/g$ is identified with the measure $\nu_*$ of
Theorem~\ref{thm:fixed-energy-spectrum}.

\begin{remark}[A deterministic total count via conformal capacitors]
\label{rev:rem-hassannezhad}
This remark sketches an alternative, deterministic route to a total
count of small eigenvalues when $g=o(n)$; it is not used anywhere else
in this paper, and the details indicated below are not carried out.
A deterministic form of the total-count statement for $g=o(n)$ expected
in~\cite[Section~1.1]{HTMany} should follow from the
conformal capacitor argument of Hassannezhad~\cite{Hassannezhad}:
for a finite-area surface $X$, its min--max levels, indexed from zero,
would satisfy
\[
 \lambda_j^{\mathrm{var}}(X)\operatorname{area}(X)\le C(g+j+1).
\]
To apply the classical argument one would compactify the cusps and
choose a smooth
conformal reference metric $h$ with curvature at least $-1$ and area
$O(g+1)$. The original area extends to a finite nonatomic measure on the
compactification. For $m=j+1$, the capacitor construction supplies $m$
disjoint plateau functions with mass at least $c\operatorname{area}(X)/m$
and energy at most $C(1+\operatorname{area}_h/m)$. One then has to check
that these bounded functions
belong to $H^1(X)$, using that the Dirichlet energy is conformally
invariant in dimension two while the cusp ends are noncompact, and
that cutting off cusp tails justifies the form domain. Min--max would
give the bound.
Since levels below $1/4$ are discrete, taking $j$ a sufficiently small
multiple of $\varepsilon n$ would prove the asserted total count when
$g=o(n)$.
\end{remark}

Equation~\eqref{eq:linear-low-count} also covers arbitrary $g\le an$ in
Weil--Petersson probability, below every prescribed $\varepsilon<1/4$;
the conformal bound need not do so at a fixed positive genus-to-cusp
ratio. Hide--Thomas~\cite[Theorem~1.3]{HTMany} give a deterministic
linear count below $1/4$ throughout that larger range. The residual
strengthening in Theorem~\ref{thm:many-residual} does not follow from
a total-count bound.

\subsection{Many small residual eigenvalues}\label{sec:residual-count}

The count of Corollary~\ref{cor:many-spectrum} can be obtained within
the residual spectrum (Theorem~\ref{thm:many-residual}), in the same
regimes $n/\sqrt g\to\infty$ and $g\le an$. Residual eigenfunctions
are the square-integrable residues of Eisenstein series, orthogonal to
the cuspidal spectrum~\cite{Hejhal,Iwaniec}; the constant zero mode is
excluded from our count. Whether a generic cofinite surface carries
cusp forms at all is a classical question going back to Phillips and
Sarnak~\cite{PhillipsSarnak}; see also Wolpert~\cite{WolpertCusp}. Let
$R_\varepsilon(X)$ be the total residual multiplicity in
$(0,\varepsilon)$.
Hide--Thomas~\cite[Section~1.1, after Remark~1.2]{HTLargeN} ask
whether their many small eigenvalues are cusp forms and anticipate a
linear residual count. The deterministic
Lemma~\ref{lem:residual-pants} gives that count from the short-pants
law of Proposition~\ref{prop:pants-concentration}, without assuming
the conjecture of Otal and Rosas~\cite[Section~6]{OtalRosas} that a
surface of genus $g\ge2$ with cusps has at most $2g-3$ small cuspidal
eigenvalues (see \cite[(1.1)]{BMP}); its main
external input is the topological eigenvalue bound of
Ballmann--Matthiesen--Mondal~\cite{BMM}. The idea is that cusp forms
have zero angular mean on cusp horocycles, so on a short two-cusp pair
of pants their mass is bounded by a fixed multiple of their energy.
Low-energy trial functions supported on the pants are therefore far
from the cuspidal spectral subspace, and project injectively onto the
residual eigenspace.

\begin{lemma}[Short pants force residual eigenvalues]
\label{lem:residual-pants}
Put $C_*=4(1+8\pi)$ and, for $0<\varepsilon<1/4$, set
\begin{equation}\label{eq:residual-length}
 L_\varepsilon=\frac{\varepsilon}{24(1+C_*\varepsilon)}.
\end{equation}
If a complete connected finite-area hyperbolic surface contains $N$
embedded two-cusp pairs of pants with disjoint interiors and geodesic
boundary lengths at most $L_\varepsilon$, then
$R_\varepsilon(X)\ge N-1$.
\end{lemma}

\begin{proof}
We first establish a uniform inequality on a two-cusp pants $P_s$
with waist $s\le1/48$. Double it across the geodesic boundary to
obtain a four-cusped sphere. The distance cutoff used in the proof of
Proposition~\ref{prop:pants-concentration} vanishes at the waist and
has Rayleigh quotient at most $6s\le1/8$. Its odd reflection therefore
produces a positive eigenvalue below $1/4$ in the odd sector of the
double. The even sector contains the constant zero mode.
By~\cite[Theorem~1.5 and equation~(1.4)]{BMM}, the double has at most
two eigenvalues at or below $1/4$, counting zero and multiplicity.
Hence its even sector has no positive spectrum below $1/4$.
Even reflection identifies this sector with the Neumann form on $P_s$;
both its norm and energy are doubled. Consequently
\begin{equation}\label{eq:residual-neumann}
 \int_{P_s}|v-\overline v_{P_s}|^2\le4\mathcal E_{P_s}(v),
 \qquad v\in H^1(P_s).
\end{equation}
This applies the topological eigenvalue bound to the boundaryless
double, not to a Neumann problem with boundary.

The standard cusp neighborhoods bounded by length-one horocycles
are embedded~\cite[Section~4.1]{HTMany}. In either pants cusp, the
subtail bounded by a horocycle of length $b=1/4$ has area $b$ and
lies inside $P_s$. Indeed, a closed geodesic meeting this subtail
must leave the length-one cusp neighborhood and so has length at
least $2\log4>s$. If $v$ has zero angular mean on the horocycles
of this subtail $H$, then $\int_Hv=0$, and
\[
 b|\overline v_{P_s}|^2
 \le\int_H|v-\overline v_{P_s}|^2
 \le4\mathcal E_{P_s}(v).
\]
Since $\operatorname{area}(P_s)=2\pi$, this and
\eqref{eq:residual-neumann} imply
\begin{equation}\label{eq:residual-anchor}
 \|v\|_{L^2(P_s)}^2\le C_*\mathcal E_{P_s}(v).
\end{equation}
Every finite linear combination of cusp-form eigenfunctions has the
required zero angular means. Their constant Fourier coefficients,
zero near infinity, vanish throughout the embedded cusp cylinder by
the Fourier equation. Summing \eqref{eq:residual-anchor} over the
pants union $U$ gives
\begin{equation}\label{eq:cuspidal-pants-mass}
 \|u\|_{L^2(U)}^2\le C_*\mathcal E_X(u)
\end{equation}
for every such combination, with no factor depending on $N$.

The $N$ distance cutoffs span a space supported in $U$ with
Rayleigh quotient at most
$\eta=6L_\varepsilon=\varepsilon/[4(1+C_*\varepsilon)]$.
Its intersection $W$ with the mean-zero space has dimension at least
$N-1$ and is still supported in $U$. If the residual spectral
subspace in $(0,\varepsilon)$ had dimension smaller than $\dim W$,
there would be a unit $f\in W$ orthogonal to it. The exact spectral
decomposition gives $f=u+v$, where $u$ belongs to the cuspidal
spectral subspace in $(0,\varepsilon)$ and $v$ has spectrum in
$[\varepsilon,\infty)$. The constant has already been removed;
the continuous spectrum starts at $1/4$. Spectral projections
contract the quadratic form, so
\[
 \mathcal E_X(u)\le\eta,\qquad \|v\|_2^2\le\eta/\varepsilon.
\]
As $f$ is supported in $U$, \eqref{eq:cuspidal-pants-mass} yields
\[
 1\le\|u\|_{L^2(U)}+\|v\|_{L^2(U)}
 \le\sqrt{C_*\eta}+\sqrt{\eta/\varepsilon}
 \le\sqrt{2(C_*+1/\varepsilon)\eta}=1/\sqrt2,
\]
a contradiction. Thus projection onto the actual residual eigenspace
is injective on $W$, proving the dimension bound. The argument allows
an eigenvalue to have both cuspidal and residual multiplicity and
counts only the latter.
\end{proof}

\begin{theorem}[Many small residual eigenvalues]\label{thm:many-residual}
For $g\to\infty$, $n/\sqrt g\to\infty$ and each fixed
$0<\varepsilon<1/4$,
\begin{equation}\label{eq:residual-supercritical}
 \mathbb P^{\mathrm{WP}}_{g,n}
 \{R_\varepsilon(X)\ge L_\varepsilon^2q_g/4\}\longrightarrow1.
\end{equation}
For every $a>0$ and $0<\varepsilon<1/4$ there is
$d_{\varepsilon,a}>0$ such that, along every sequence
$n\to\infty$ and $0\le g\le an$,
\begin{equation}\label{eq:residual-linear}
 \mathbb P^{\mathrm{WP}}_{g,n}
 \{R_\varepsilon(X)\ge d_{\varepsilon,a}(2g+n-2)\}
 \longrightarrow1.
\end{equation}
\end{theorem}

\begin{proof}
Since $L_\varepsilon<1/96$, the counted pants are interior-disjoint.
Proposition~\ref{prop:pants-concentration} and
$\mathfrak B_Y(L)\ge L^2$ give, with probability tending to one,
$N_{L_\varepsilon}-1\ge q_gL_\varepsilon^2/4$.
Lemma~\ref{lem:residual-pants} proves
\eqref{eq:residual-supercritical}. The positive lower bound for
$q_g/(2g+n-2)$ when $g\le an$, established in the proof of
Corollary~\ref{cor:many-spectrum}, gives
\eqref{eq:residual-linear} when $g\to\infty$.
For each fixed genus, use the short-pants count of
Hide--Thomas~\cite[Section~5, equation~(5.1)]{HTLargeN}:
$N_{L_\varepsilon}\ge c(L_\varepsilon)n$ with probability tending
to one, where the positive constant is independent of the fixed
genus. Applying Lemma~\ref{lem:residual-pants} gives the residual
count in this case as well. Reduce $d_{\varepsilon,a}$ to a common
constant and use the same fixed-or-diverging genus subsequence
argument as in Corollary~\ref{cor:many-spectrum}.
\end{proof}

This establishes the linear residual count anticipated
in~\cite[Section~1.1]{HTLargeN}, including every fixed genus, and
extends it to growing genus through the uniform volume theorem. The
count is far below the Otal--Rosas bound $2g+n-2$ on the number of
eigenvalues in $[0,1/4]$, counting the eigenvalue $0$~\cite{OtalRosas}: $q_g/(2g+n-2)$ is bounded
by \eqref{eq:PS1}, and $L_\varepsilon^2/4<10^{-4}$.
It does not assert that all small eigenvalues are residual, and it
does not bound the total number of small cuspidal eigenvalues.
The mechanism behind Lemma~\ref{lem:residual-pants}, that cusp forms
have vanishing cusp averages while localized pants states do not,
returns in sharper form at the critical scale
(Lemma~\ref{rev:lem-cusp-average-injective}) and in the dense regime
(Corollary~\ref{cor:dense-residual}), where it identifies the whole
low eigenspace as residual.

\subsection{The tight length spectrum}\label{subsec:tight-spectrum}

We now consider curves that remain noncontractible after all cusps are
filled, and prove their Poisson limit at every cusp density
(Theorem~\ref{thm:tight-spectrum}) and for every subcritical cusp
weight (Corollary~\ref{cor:tight-boltzmann}). Following Budd and
Lions~\cite{BL}, a closed geodesic is \emph{tight} if its class on the
filled surface is noncontractible and it is the unique shortest
geodesic on the original surface among those in that free homotopy
class after filling the cusps. The minimality condition matters, since
several ordinary geodesics can belong to the same class after the
cusps are filled. The corresponding tight Weil--Petersson volumes were
introduced by Budd and Zonneveld~\cite{BZ}; see also~\cite{Budd} for the
related irreducible metric maps. The tight integration formula and the
geometric argument used in the proof are due to~\cite{BL}; the new
input is a volume estimate at a prescribed cusp count, obtained from
the shifted amplitudes and the coefficient extraction of
Part~\ref{part:volumes}, and the proofs are in
Appendix~\ref{app:tight}.

Let $Y_\lambda$ satisfy $\lambda(Y_\lambda)=\lambda$,
as in Proposition~\ref{prop:bounded-boundaries}, and set
\begin{equation}
 \alpha(\lambda)=\frac{\sqrt{2r(Y_\lambda)}}\pi,
 \qquad \nu(\mathrm dt)=\frac{\cosh t-1}{t}\,\mathrm dt.
 \label{eq:tight-scale}
\end{equation}
Here $\alpha(0)=1$, and the density of $\nu$ is defined to be zero
at the origin.

\begin{theorem}[Tight lengths at arbitrary cusp density]
\label{thm:tight-spectrum}
Let $n=n(g)\ge0$ be any integer sequence, and let $X$ have law
$\mathbb P^{\mathrm{WP}}_{g,n}$. As $g\to\infty$, the point process
\begin{equation}
 \sum_{\substack{\gamma\text{ primitive, unoriented}\\
                  \gamma\text{ tight on }X}}
 \delta_{\alpha(n/g)\ell_X(\gamma)}
 \ \Longrightarrow\ \operatorname{PPP}(\nu)
 \label{eq:tight-process}
\end{equation}
converges in distribution for the vague topology on $[0,\infty)$.
For every fixed $L<\infty$, with probability tending to one all the
curves counted in $[0,L]$ are simple, pairwise disjoint, and have
connected complement after filling the cusps.
\end{theorem}

In particular, $n/g$ need not converge. At large density,
the endpoint expansion of $r$ gives
\[
 \alpha(\lambda)\sim
 \frac{\sqrt{3j_{0,1}\sqrt5}}\pi\lambda^{-1/4}.
\]
This agrees with the scale in~\cite[Corollary~1.3]{BL}.
Budd and Lions asked whether their limit persists for other cusp
weights and under the fixed-count measures
$\mathbb P^{\mathrm{WP}}_{g,n}$~\cite[Open Question~1.6]{BL}.
Theorem~\ref{thm:tight-spectrum} treats the fixed-count part of this
question, and Corollary~\ref{cor:tight-boltzmann} below treats the
remaining cusp-weight regimes. The planar curves in
Theorem~\ref{thm:pants-poisson} are contractible after the cusps are
filled, and therefore do not occur in \eqref{eq:tight-process}.

The proof follows the moment method of Budd and Lions: their tight
integration formula~\cite[Proposition~3.2]{BL} expresses factorial
moments of tight curve counts through tight volumes of the cut
surfaces. In Appendix~\ref{app:tight} we write the generating
functions of tight volumes through the shifted amplitudes of
Section~\ref{exact-amplitudes-and-the-harmonic-model}
(Lemma~\ref{lem:tight-gf}), extract the coefficient at the prescribed
cusp count as in Section~\ref{extracting-fixed-cusp-counts}
(Lemma~\ref{lem:tight-connected}), show that configurations with
disconnected complement are negligible
(Lemma~\ref{lem:tight-separating}), and identify the Poisson limit
(Subsection~\ref{subsec:tight-proof}).

\emph{Random cusp counts.}
For $0\le\mu<\mu_{\rm cr}=x_c/(2\pi^2)$, the Boltzmann measure
$\mathbb P^{\mathrm{WP}}_{g,\mu}$ first chooses the cusp count $N$
with probabilities proportional to $\mu^nV_{g,n}/n!$, then samples
$\mathbb P^{\mathrm{WP}}_{g,n}$ conditionally on $N=n$.
The probability generating function of the cusp count is
$\Phi_g(2\pi^2\mu z)/\Phi_g(2\pi^2\mu)$, as in~\cite{BL}.

\begin{corollary}[Arbitrary subcritical cusp weights]
\label{cor:tight-boltzmann}
For any sequence $0\le\mu_g<\mu_{\rm cr}$, put
\[
 a_g=\frac{\sqrt{2r(Y(2\pi^2\mu_g))}}\pi.
\]
As $g\to\infty$, under $\mathbb P^{\mathrm{WP}}_{g,\mu_g}$, the counting
measure of primitive unoriented tight lengths rescaled by $a_g$ converges
in distribution to $\operatorname{PPP}(\nu)$ for the vague topology
on $[0,\infty)$.
In particular, if $\mu_g\to\mu<\mu_{\rm cr}$, the unscaled lengths have
limiting intensity $(\cosh(a_\mu t)-1)\,\mathrm dt/t$, where
$a_\mu=\sqrt{2r(Y(2\pi^2\mu))}/\pi$. If $\mu_g\to\mu_{\rm cr}$, no condition on
the speed of convergence is required.
\end{corollary}
\noindent The proof is given in Subsection~\ref{subsec:tight-boltzmann-proof}.
The tight spectrum is not used in the spectral sections that follow;
its role here is to complete the description of short curves at every
cusp density, the planar separating curves having been treated in
Subsection~\ref{sec:pants-spectrum}. The next subsection returns to
those planar curves at the critical scale, where both families live
on the same length scale.

\subsection{The full length spectrum at the critical scale}
\label{sec:critical-cluster}

This subsection determines the full primitive length spectrum at the
critical scale $n/\sqrt g\to c\in(0,\infty)$
(Theorem~\ref{thm:critical-cluster}). The tight spectrum omits curves
that become contractible when the cusps are filled; at this scale,
these curves contribute on the same length scale as the nonseparating
curves, and their limiting contribution is a cluster process: each
two-cusp pair of pants contributes its whole primitive length
spectrum, not only its boundary. The volume inputs are the bounds of
Mirzakhani--Zograf and Shen--Wu~\cite{MZ,SW} together with the
zero-density case of Theorem~\ref{thm:volume} and
Proposition~\ref{prop:bounded-boundaries}; the proof is in
Appendix~\ref{app:critical}. The waist process $\mathcal Q_c$ obtained
here is the geometric input of the critical spectral law in
Section~\ref{sec:critical-eigenvalues}.

For $s>0$, let $P_s$ be the hyperbolic two-cusp pair of pants with
one geodesic boundary of length $s$. Write $K_s$ for its primitive
unoriented length counting measure, including the boundary once and
excluding parabolic classes. Equal lengths are counted with their
geometric multiplicities. Let
\begin{equation}
 \mu_c(\mathrm ds)=\frac{c^2}{8\pi^2}\sinh(s/2)\,\mathrm ds,
 \qquad \nu(\mathrm dt)=\frac{\cosh t-1}{t}\,\mathrm dt.
 \label{eq:cluster-intensities}
\end{equation}

\begin{theorem}[Critical primitive length spectrum]
\label{thm:critical-cluster}
Suppose $g\to\infty$ and $n/\sqrt g\to c\in(0,\infty)$, and let
$X$ have law $\mathbb P^{\mathrm{WP}}_{g,n}$. Then
\begin{equation}
 \sum_{\substack{\gamma\text{ primitive}\\\text{unoriented on }X}}
       \delta_{\ell_X(\gamma)}
 \ \Longrightarrow\ \Lambda_c:=\mathcal N+\sum_{s\in\mathcal Q_c}K_s,
 \label{eq:critical-cluster}
\end{equation}
in distribution in the space of locally finite counting measures on
$[0,\infty)$ with the vague topology. Here $\mathcal N$ and $\mathcal Q_c$ are independent Poisson
point processes of intensities $\nu$ and $\mu_c$. In particular, for
every nonnegative continuous compactly supported function $f$,
\begin{equation}
 \mathbb E\mathrm e^{-\langle f,\Lambda_c\rangle}
 =\exp\left\{-\int(1-\mathrm e^{-f(t)})\,\nu(\mathrm dt)
 -\int(1-\mathrm e^{-\langle f,K_s\rangle})\,\mu_c(\mathrm ds)\right\}.
 \label{eq:cluster-laplace}
\end{equation}
\end{theorem}

Shen and Wu obtained the factorial-moment asymptotics for short
two-cusp pairs of pants and the limiting probability that such a curve occurs
at this cusp scale~\cite[Lemmas~28 and~33]{SW}; on their stated short
window the corresponding Poisson mean is
$\mu_c([0,L])=c^2(\cosh(L/2)-1)/(4\pi^2)$.
Theorem~\ref{thm:critical-cluster} concerns the full primitive
process on every fixed length window, including the correlated lengths
inside each pair of pants. The background intensity and the geometric
moment method come from Mirzakhani and Petri~\cite{MP}. The bounds
of~\cite{MZ,SW} alone give an upper estimate and a positive lower bound
for the genus ratio $V_{g-1,n+2}/V_{g,n}$ at this scale, but not its
limit, which is why the zero-density
case of Theorem~\ref{thm:volume} is needed; the first-order expansion
of Gatignol and Thomas~\cite[Theorem~1.5]{GT} is another source of its
two-sided limit. The full uniform expansion is not used.

Nonsimple curves distinguish this regime from closed surfaces.
He, Shen, Wu and Xue proved that the nonsimple systole of a
Weil--Petersson random closed surface is asymptotic to $\log g$, and
more precisely concentrates about
$\log g-\log\log g$~\cite[Theorem~1.1]{HSWX}; in that model nonsimple
curves therefore leave every fixed length window. At the critical
scale, by contrast, two-cusp pairs of pants retain nonsimple curves of
bounded length.

\emph{Outline of the proof.} Short multicurves are organized by their
cut graphs. Using only the volume bounds of~\cite{MZ,SW},
Lemma~\ref{lem:critical-cuts} shows that a cut graph contributes in
the limit only if it consists of one component carrying the genus,
possibly with nonseparating loops, and two-cusp pairs of pants attached
to it as leaves. A bounded-multiplicity argument (Lemmas~\ref{lem:primitive-count}
and~\ref{lem:critical-complex}) reduces intersecting simple curves and
nonsimple curves to such cut graphs; hence, with probability tending to
one, every short nonsimple curve lies in a two-cusp pair of pants. The mixed
factorial moments of nonseparating lengths and pants waists converge to those
of independent Poisson processes (Lemma~\ref{lem:critical-skeleton}, where the
zero-density volume ratios enter), and attaching to each waist the full
spectrum of its pair of pants gives $\Lambda_c$. The details are in
Appendix~\ref{app:critical}.

The limit is not Poisson. If $m_L(s)=K_s([0,L])$, then
Lemma~\ref{lem:primitive-count} and the support restriction $s\le2L$
give finite moments and
\begin{equation}
 \operatorname{Var}\Lambda_c([0,L])-\mathbb E\Lambda_c([0,L])
 =\int m_L(s)(m_L(s)-1)\,\mu_c(\mathrm ds).
 \label{eq:cluster-overdispersion}
\end{equation}
Besides its waist, $P_s$ contains the primitive figure-eight around
its two cusps, of length $2\operatorname{arccosh}(2+\cosh(s/2))$
by \cite[Theorem~1.3]{Rivin}. For $L>2\operatorname{arccosh}3$,
both lengths lie below $L$ on an interval of positive $\mu_c$-measure.
The right side of \eqref{eq:cluster-overdispersion} is therefore
strictly positive. The theorem identifies these correlations in the
limit; no convergence of the full-spectrum second moments is asserted.
Of the cluster structure, only the waist process $\mathcal Q_c$, its
intensity $\mu_c$ and the factorial-moment bounds behind it are used
in Section~\ref{sec:critical-eigenvalues}; the cut-graph
Lemmas~\ref{lem:critical-cuts} and~\ref{lem:critical-complex} of
Appendix~\ref{app:critical} are reused there to isolate the pants.

\section{The critical low spectrum}\label{sec:critical-eigenvalues}

At the critical scale $n/\sqrt g\to c\in(0,\infty)$, the full positive
spectral point process below $1/4$ converges to a Poisson process of
model energies (Theorem~\ref{thm:critical-eigenvalues}). This gives the gap
phase diagram (Corollary~\ref{cor:gap-phase-diagram}). The geometric
input is the waist process $\mathcal Q_c$ of
Theorem~\ref{thm:critical-cluster}; the spectral input is the trace
expansion of Huang~\cite{HuangGap}, restated and extended to the
critical envelope in Appendix~\ref{app:critical-trace}; the model
input is the spectral theory of the funnel surface $F_s$ defined
below, through Borthwick--Judge--Perry~\cite{BJP},
Schulze~\cite{Schulze} and Ballmann--Matthiesen--Mondal~\cite{BMM}.
The argument is carried out in this section, with the cut-graph
lemmas of Appendix~\ref{app:critical} as the remaining input.

At the critical scale, a bounded number of two-cusp pairs of pants
occur in each bounded waist interval, and their exterior neighborhoods
become arbitrarily deep funnels as the genus grows. The spectral model
attached to a waist is therefore the complete surface obtained by gluing
a funnel to the pair of pants, not the pants with a Dirichlet condition
at its boundary. The distinction matters for bounded waist lengths,
since a cusped pair of pants need not have a Dirichlet eigenvalue below
$1/4$~\cite[Remark~2.2]{HTMany}.

Let $F_s$ be the complete hyperbolic surface with two cusps and one funnel
whose convex core is $P_s$, with boundary length $s>0$; we call $F_s$
the \emph{funnel model} of waist $s$. Write $\beta(s)$
for its lowest eigenvalue; see~\cite{Borthwick} for the spectral theory
of surfaces with funnels. The relation between the bottom of the spectrum
and the critical exponent goes back to Patterson~\cite{Patterson}. The
cusp theorem of Beardon~\cite{Beardon}
and the spectral relation of Sullivan~\cite[Theorems~2.17 and~2.21]{Sullivan}
give
\begin{equation}\label{eq:critical-model-energy}
 \beta(s)=\delta(s)(1-\delta(s))\in(0,1/4),
 \qquad 1/2<\delta(s)<1,
\end{equation}
where $\delta(s)$ is the critical exponent of its uniformizing group.
By~\cite[Theorem~1.5 and Corollary~1.6]{BMM}, there is at most one
eigenvalue below $1/4$, counted with multiplicity; since $\delta(s)>1/2$,
there is exactly one.
The function $\beta$ is continuous, including the limit $\beta(s)\to0$
as $s\downarrow0$, by spectral continuity through pinching
\cite[Theorems~2.14 and~3.21]{Schulze}; the pinched finite-area
component is the thrice-punctured sphere, whose constant eigenfunction
has eigenvalue zero. Lemma~\ref{lem:funnel-long-waist} below gives the
other endpoint, $\beta(s)\to1/4$ as $s\to\infty$.

Retain the pants process $\mathcal Q_c$ and its intensity $\mu_c$ from
\eqref{eq:cluster-intensities}, and define the energy intensity
\begin{equation}\label{eq:critical-energy-intensity}
 \nu_c(B)=\frac{c^2}{8\pi^2}\int_0^\infty
       \mathbf1_{\{\beta(s)\in B\}}\sinh(s/2)\,ds,
 \qquad B\subset[0,1/4).
\end{equation}
It is finite on compact subsets of $[0,1/4)$ and has no atom at zero.
We use $G(X)=\inf\operatorname{spec}(\Delta_X|\mathbf1^\perp)$ for the
gap, with $G(X)=1/4$ when no positive eigenvalue lies below $1/4$.

\begin{theorem}[Critical exceptional spectrum]\label{thm:critical-eigenvalues}
Suppose $g\to\infty$ and $n/\sqrt g\to c\in(0,\infty)$. For a
Weil--Petersson random surface $X\in\mathcal M_{g,n}$,
\begin{equation}\label{eq:critical-eigenvalue-process}
 \Xi_X:=\sum_{0<\lambda_j(X)<1/4}\delta_{\lambda_j(X)}
 \ \Longrightarrow\
 \Xi_c:=\sum_{s\in\mathcal Q_c}\delta_{\beta(s)}
\end{equation}
in the vague topology on $[0,1/4)$. All discrete eigenvalues, including
residual eigenvalues, are counted with multiplicity. The limit is a
Poisson random measure of intensity $\nu_c$, and
$\mathbb E\Xi_X\to\nu_c$ vaguely. In particular $G(X)$ converges in
distribution to $G_c\in(0,1/4)$ with
\begin{equation}\label{eq:critical-gap-law}
 \mathbb P\{G_c>t\}=\exp\{-\nu_c([0,t])\},\qquad 0\le t<1/4.
\end{equation}
Writing $\Xi_X^{\rm res}$ and $\Xi_X^{\rm cusp}$ for the residual and
cuspidal parts, respectively, we have jointly
\begin{equation}\label{eq:critical-residual-process}
 (\Xi_X^{\rm res},\Xi_X^{\rm cusp})
       \Longrightarrow(\Xi_c,0),\qquad
 \mathbb E\Xi_X^{\rm res}\longrightarrow\nu_c,
 \quad\mathbb E\Xi_X^{\rm cusp}\longrightarrow0
\end{equation}
in the same vague topologies. More generally, let $C_\varepsilon(X)$
be the cuspidal multiplicity in $(0,1/4-\varepsilon]$. For every fixed
$0<M<\infty$ and $0<\varepsilon<1/4$,
\begin{equation}\label{eq:critical-cuspidal-exclusion}
 \sup_{1\le n\le M\sqrt g}
       \mathbb E^{\rm WP}_{g,n}C_\varepsilon(X)\longrightarrow0.
\end{equation}
\end{theorem}

The primitive length process in Theorem~\ref{thm:critical-cluster}
retains all geodesics inside each pants and is a cluster process. In
contrast, each funnel model contributes exactly one
subthreshold eigenvalue to~\eqref{eq:critical-eigenvalue-process}.
Neither the pushforward description nor the gap law requires
monotonicity of $\beta$.

The proof combines the waist process $\mathcal Q_c$ with Huang's
orbital expansion and positive exceptional
measure~\cite[Sections~3.7 and~6.8]{HuangGap}, and adds three steps,
each treated in its own subsection after a preliminary lemma on large
waists.
First, we extend the trace estimate to each bounded critical envelope
and identify its quadratic correction with the model energies; a finite
contour deformation handles the continuous-spectrum contribution. The
resulting Proposition~\ref{prop:critical-mean}, together with
Proposition~\ref{prop:critical-identification}, is a two-sided limit
of the mean measure; in the count matching it is used as the upper
half, and the lower half is supplied by the localized lower count of
Lemma~\ref{rev:lem-localized-lower-count}.
Second, disjoint local ground states give lower spectral counts. Since
their mean exhausts the limiting expected count, convergence follows
without expanding joint eigenvalue traces. Third, compact model states
with small operator residuals have independent nonzero cusp averages.
Spectral projection and the matching dimension make the cusp-average
map injective on the actual low eigenspace, so this space is residual
(Lemma~\ref{rev:lem-cusp-average-injective}).
These arguments concern fixed compact windows below $1/4$, not windows
approaching the threshold; the latter are the subject of
Section~\ref{sec:moving-edge}.

\subsection{The model at large waist lengths}

The following lemma supplies the endpoint $\beta(s)\to1/4$ and the
integrability of the large-waist tail of the orbital sum. Both are
used in the identification of the mean measure
(Proposition~\ref{prop:critical-identification}); the endpoint also
gives the local finiteness of the fixed-window measure
\eqref{eq:fixed-energy-measure} in
Section~\ref{sec:fixed-energy-spectrum}, and the generators introduced
in the proof reappear in Section~\ref{sec:moving-edge}.

\begin{lemma}\label{lem:funnel-long-waist}
For every $z>1/2$, sufficiently large $S$ satisfies
\begin{equation}\label{eq:funnel-tail}
 \int_S^\infty\sinh(s/2)
   \sum_{\gamma\ {\rm primitive\ oriented}}\sum_{k\ge1}
       \ell_\gamma e^{-zk\ell_\gamma}\,ds<\infty.
\end{equation}
Moreover $\delta(s)\to1/2$ as $s\to\infty$. Thus for each $t<1/4$
there is $S_t<\infty$ such that $\beta(s)>t$ for $s>S_t$.
\end{lemma}

\begin{proof}
Use the two-parabolic marking
\[
 A=\begin{pmatrix}1&a\\0&1\end{pmatrix},\quad
 B=\begin{pmatrix}1&0\\-a&1\end{pmatrix},\quad
 a=2\cosh(s/4)=\kappa+\kappa^{-1},\quad\kappa=e^{s/4}.
\]
These generate the free group of $F_s$, and $AB$ represents the waist.
Every hyperbolic conjugacy class has an alternating cyclic representative
$W=A^{m_1}B^{n_1}\cdots A^{m_q}B^{n_q}$ with nonzero integer exponents.
For $C(u)=\left(\begin{smallmatrix}u&-1\\1&0\end{smallmatrix}\right)$,
$A^mB^n=-C(ma)C(na)$. The projective map of $C(ma)$ sends
$r\in[-\kappa^{-1},\kappa^{-1}]$ to $(ma-r)^{-1}$, where $r=x_2/x_1$ is
the ratio of the coordinates on which $C(ma)$ acts; it preserves this
interval and contracts it by at most $\kappa^{-2}$. Its first-coordinate
expansion is at least $|m|\kappa$. At the fixed direction of the product,
the expanding eigenvalue therefore gives
\begin{equation}\label{eq:funnel-word}
 \ell(W)\ge qs+2\sum_{j=1}^q\log|m_jn_j|.
\end{equation}
Consequently, for $z_0>1/2$, the sum of $e^{-z_0\ell(W)}$ over all
such words is bounded by
\[
 \sum_{q\ge1}\{4\zeta(2z_0)^2e^{-z_0s}\}^q.
\]
This overcounts every oriented conjugacy class, including powers.
Choose $1/2<z_0<z$ and then $S$ so that the ratio is at most $1/2$.
The inequality $ue^{-(z-z_0)u}\le[e(z-z_0)]^{-1}$ bounds the sum
in~\eqref{eq:funnel-tail} by $C_{z,z_0}e^{-z_0s}$, proving integrability.
The same word estimate makes the Selberg Euler product nonzero in
$\Re u\ge z_0$ for large $s$. Its eigenvalue zero at $\delta(s)$
\cite[Theorem~1.1]{BJP} implies $\delta(s)\le z_0$ there.
Let $z_0\downarrow1/2$ and use~\eqref{eq:critical-model-energy}.
\end{proof}

\subsection{The mean measure and its model identification}

This subsection proves the two-sided limit of the mean spectral
measure at the critical scale. It has two parts: a trace estimate that
expresses the mean spectrum through Huang's exceptional measure
(Proposition~\ref{prop:critical-mean}), and an identification of that
measure with the model energies
(Proposition~\ref{prop:critical-identification}). The first part
rests on the extension of Huang's orbital estimate to the critical
envelope $n\le M\sqrt g$, proved in Appendix~\ref{app:critical-trace};
the second rests on the regularized Selberg zeta function of the
funnel model through the pinching limit
(Lemma~\ref{rev:lem-schulze-pinching}), which uses
Schulze~\cite{Schulze} and Borthwick--Judge--Perry~\cite{BJP}.
We use the admissible multiplier
of Huang~\cite[Section~1.2]{HuangGap}:
$h=b_0*b_0$ is even, nonnegative, smooth, supported in $[-L,L]$ for a
fixed $L>0$, and has
integral one, $f(r)=\int h(u)e^{-iru}\,du$, and $a=f(i/2)>1$.
Thus $0\le f(r)\le1$ for real $r$, while $\tau\mapsto f(i\tau)$
is strictly increasing for $\tau>0$. The support length $L$ enters
the main text only through the cusp truncation height $Y_Q=3e^{LQ}$ below.
Subsection~\ref{rev:subsec-huang-inputs} of
Appendix~\ref{app:critical-trace} restates in our notation the results
of \cite{HuangGap} used in this subsection; in particular
Definition~\ref{rev:def-huang-orbital} introduces the orbital functional
$\mathcal O_X(p)$ and the canonical coefficients
$\omega_1$, $\omega_{j,n}$ and $W_i$ of the orbital expansion.

Write $\eta_2$ for Huang's positive exceptional measure and
$\mathscr T_2$ for his weighted distribution. Their construction
(Lemma~\ref{rev:lem-huang-exceptional}) is as follows. For each $n\ge1$
the functional $\mathscr K_n(q)=\omega_{2,n}(xq)$ on polynomials $q$
extends uniquely to a distribution supported in $[0,a]$ of order at most
$m=2D_0+2$, where $D_0$ is the exponent in the coefficient bound of
Lemma~\ref{rev:lem-huang-coefficients}, and
$|\mathscr K_n(\varphi)|\le K_f(n+1)^2\|\varphi\|_{C^m([0,a])}$. Taking
the second difference in $n$ of the quadratic identity
\eqref{rev:eq-huang-quadratic-p} gives
\[
 \mathscr T_2=\mathscr K_3-2\mathscr K_2+\mathscr K_1,\qquad
 \mathscr T_2(q)=W_2(xq),\qquad
 |\mathscr T_2(\varphi)|\le38K_f\|\varphi\|_{C^m([0,a])}.
\]
The restriction of $\mathscr T_2$ to $(1,\infty)$, divided by the
coordinate $x$, is a positive Radon measure $\eta_2$ on $(1,a]$:
\begin{equation}\label{eq:critical-weighted-distribution}
 \mathscr T_2(x^{j-1})=\Lambda_2(h^{*j}),\quad j\ge1,
 \qquad \mathscr T_2=x\eta_2\quad\hbox{on }(1,a].
\end{equation}
Here $\Lambda_2$ is the coefficient of $\binom n2$ in the first
orbital correction \eqref{rev:eq-huang-quadratic}, so that
$W_2(p)=\Lambda_2(H_p)$ and $H_{x^j}=h^{*j}$. Huang's construction
allows $\eta_2$ to have an atom at $a$, and his
\cite[Theorem~6.11]{HuangGap} is stated for tests that may be nonzero at
$a$; Proposition~\ref{prop:critical-identification} below shows that no
such atom is present; Huang's argument does not by itself exclude it.

\begin{proposition}[Critical mean spectrum]\label{prop:critical-mean}
If $n/\sqrt g\to c\in(0,\infty)$, then, for every
$\psi\in C_c^\infty((1,\infty))$,
\begin{equation}\label{eq:critical-mean-multiplier}
 \mathbb E\sum_{0<\lambda_j<1/4}
 \psi\bigl(f(i\sqrt{1/4-\lambda_j})\bigr)
 \longrightarrow \frac{c^2}{2}\int\psi\,d\eta_2.
\end{equation}
The tests may be nonzero at $a$; the constant eigenfunction is excluded
before taking the limit. The assertion includes residual eigenvalues.
\end{proposition}

\begin{proof}
Choose $M>c+1$, so $n+1\le M\sqrt g$ eventually.
Appendix~\ref{app:critical-trace} extends Huang's orbital estimate to
this envelope, with constants depending on $M$. In particular, for
$Q=\lfloor g^\sigma\rfloor$ with sufficiently small $\sigma>0$, and
$p=x\widetilde p$ of degree at most $Q$, it gives
\begin{equation}\label{eq:critical-polynomial-error}
 \left|\mathbb E\mathcal O_X(p)-\omega_1(p)
                   -g^{-1}\omega_{2,n}(p)\right|
 \le\left\{\frac{2(n+1)^3 A_Q^3}{g^2}
       +\frac{(n+1)^2}{g}2^{-Q}\right\}
       \|\widetilde p\|_{[0,a]}=o(1)\|\widetilde p\|_{[0,a]},
\end{equation}
where $A_Q=(3C_0Q)^{D_0}$ and the coefficients are exactly those of
Lemmas~\ref{rev:lem-huang-coefficients} and~\ref{rev:lem-huang-transfer}
(the coefficients of \cite[Lemma~6.7]{HuangGap}; the transfer is the
deterministic part of the proof of \cite[Proposition~6.9]{HuangGap}).
The same appendix
gives $\mathbb E N_X(1/4-\varepsilon)\le C_{M,\varepsilon}(n+1)^2/g$.

Approximate $\xi=\psi/x$, extended by zero near zero, by its shifted
Chebyshev truncation $r_Q$ of degree $Q-1$, and put $v_Q=xr_Q$.
For every fixed $K$, $\|r_Q-\xi\|_\infty=O(Q^{-K})$, and
$r_Q\to\xi$ in the fixed $C^m$ norm controlling
$\mathscr T_2$ and the three distributions $\mathscr K_1,\mathscr K_2,\mathscr K_3$
of Lemma~\ref{rev:lem-huang-exceptional}.
Let $\tau_{X,Q}(p)=\operatorname{Tr}(\chi_{Y_Q}p(A_X)\chi_{Y_Q})$ be the
positive trace of Lemma~\ref{rev:lem-huang-trace}, compressed to the
complement of the width-one cusp regions above the common height
$Y_Q=3e^{LQ}$ of \eqref{rev:eq-huang-YQ}. Its expected first moment is
$O_M(g+nQ)$ by Lemma~\ref{rev:lem-huang-first-moment}, whose hypothesis
is the base case of \eqref{eq:critical-monomial-error}. Since
$A_X\ge0$ and $|v_Q-\psi|\le x\|r_Q-\xi\|_\infty$ on $[0,a]$,
positivity, in the form of Lemma~\ref{rev:lem-huang-trace}(ii) applied
to $v_Q-\psi=x(r_Q-\xi)$, bounds the difference between the
constant-removed traces
$\tau^\circ_{X,Q}(v_Q)$ and
$\tau^\circ_{X,Q}(\psi)=\operatorname{Tr}(\chi_{Y_Q}\psi(A_X)(1-P_0)\chi_{Y_Q})$
of Lemma~\ref{rev:lem-huang-trace}(v) by
$\|r_Q-\xi\|_\infty\tau_{X,Q}(x)$, whose expectation is
$O_M((g+nQ)Q^{-K})$. As $\psi$ vanishes on $[0,1]$,
$\tau^\circ_{X,Q}(\psi)$ involves only the eigenvalues in $(0,1/4)$,
each weighted by its retained mass. In the exact
trace identity \eqref{rev:eq-huang-trace-identity}, the identity,
parabolic and leading real-spectrum terms, bounded by
\eqref{rev:eq-huang-deterministic} and by $|\mu_0(p)|\le5N_{\mathbb R}(p)$,
cost $O_M((g+nQ^2)Q^{-K})$ for this approximation, since $\xi$ vanishes
on $[0,1]$; choose $K$ large. Remove the constant projection as in
\eqref{rev:eq-huang-trace-compressed}: the removed mass of the constant
eigenfunction contributes exactly
$n v_Q(a)/(\operatorname{area}(X)Y_Q)=o(1)$, because $v_Q(a)$ is bounded
and $Y_Q$ grows exponentially in $Q$. Consequently the expected compressed
exceptional trace $\mathbb E\tau^\circ_{X,Q}(v_Q)$ equals
$g^{-1}\omega_{2,n}(v_Q)+o(1)$ by \eqref{eq:critical-polynomial-error}.

The exact quadratic identity \eqref{rev:eq-huang-quadratic-p} is
$\omega_{2,n}=W_0+nW_1+\binom n2 W_2$, where, on tests $xq$, each
$W_i(xq)$ is a fixed linear combination of
$\mathscr K_1(q),\mathscr K_2(q),\mathscr K_3(q)$.
The $C^m$ convergence just noted makes each $W_i(v_Q)$ converge;
only $\binom n2 W_2(v_Q)/g$ survives. By
\eqref{eq:critical-weighted-distribution}, its limit is the right
side of \eqref{eq:critical-mean-multiplier}.
Finally, the removed mass of an eigenfunction with
$\lambda\le1/4-\varepsilon$ is at most
$e^{-2\sqrt\varepsilon LQ}$ by Lemma~\ref{rev:lem-huang-retained-mass}.
The expected count bound removes these weights in expectation.
Smooth majorants and the same local count bound also extend the
conclusion to continuous compact tests on the physical interval
$[0,1/4)$.
\end{proof}

It remains to identify $\eta_2$. Put
$w(s)=\sinh(s/2)/(4\pi^2)$. Huang's explicit quadratic coefficient
\cite[Section~3.7]{HuangGap} gives the exact identity
\begin{equation}\label{eq:critical-orbital-identity}
 \Lambda_2(H)=\int_0^\infty w(s)\mathcal O_{F_s}(H)\,ds,
 \qquad
 \mathcal O_{F_s}(H)=\sum_{\gamma\ \mathrm{primitive\ oriented}}
 \sum_{k\ge1}\frac{\ell_\gamma H(k\ell_\gamma)}
                         {2\sinh(k\ell_\gamma/2)}.
\end{equation}
The sum includes the two oriented waist classes and all nonsimple
classes in the pants. The former give
$\sum_{k\ge1}sH(ks)/\sinh(ks/2)$, and the latter give Huang's
$Z_2=(8\pi^2)^{-1}\int 2\sinh(s/2)I_2(s)\,ds$; this verifies the
normalization in \eqref{eq:critical-orbital-identity}.

The analytic input for the identification is a regularized Selberg
zeta function of $F_s$ that is uniformly controlled through the pinching
limit $s\downarrow0$. Let $Z_s(z)$ be the Selberg zeta function of $F_s$
in the oriented primitive-class convention of \cite{BJP}, that is, the
Euler product used in Lemma~\ref{lem:funnel-long-waist}, and set
\[
 D_s(z)=\prod_{j\ge0}(1-e^{-(z+j)s})^2,\qquad
 N_s(z)=\beta(s)-z(1-z).
\]
Thus $D_s$ is the contribution of the two oriented waist classes to the
Euler product, and $N_s$ vanishes exactly at $z=\delta(s)$ and
$z=1-\delta(s)$.

\begin{lemma}[Regularized zeta function through pinching]
\label{rev:lem-schulze-pinching}
For every $s>0$ the function $Z_s/(D_sN_s)$ is holomorphic and zero-free
on $\{\Re z>1/2\}$. For every $S>0$ and every compact
$K\subset\{\Re z>1/2\}$ there are constants
$0<c_{K,S}\le C_{K,S}<\infty$ such that
\[
 c_{K,S}\le\Bigl|\frac{Z_s(z)}{D_s(z)N_s(z)}\Bigr|\le C_{K,S}
 \qquad(z\in K,\ 0<s\le S).
\]
\end{lemma}

\begin{proof}
We place $F_s$ in the framework of Schulze~\cite{Schulze}. There a
geometrically finite hyperbolic surface is assembled from pairs of pants
along an augmented admissible graph $G^*$ \cite[Section~1.2]{Schulze}:
the proper edges carry the gluing cylinders, and the remaining
\emph{phantom} edges carry the infinite ends, a phantom edge of positive
length $\ell$ giving a funnel and one of length zero a cusp. The
parameter space is
$\mathcal B=\{(\ell,\tau):\ell\colon G_1^*\to[0,\infty),\ \tau\colon G_1^*\to\mathbb R\}$,
in which length zero is admitted on every edge
\cite[Section~2.2]{Schulze}. The surface $F_s$ is the non-elementary
component of $Y_G(\lambda_s)$ for
the graph with one vertex and three phantom edges, two of length zero
(the cusps) and one, the funnel edge $d$, of length $\ell(d)=s$; the
twist parameters play no role on phantom edges. Each zero-length
phantom edge also contributes an isolated elementary parabolic
component, so $Y_G(\lambda_s)$ has two such components for $s>0$
\cite[Definition~1.5 and the discussion following it]{Schulze}.
These components have no $L^2$ eigenvalues and contribute the factor
$1$ to the Selberg zeta function \cite[Section~1.4]{Schulze}.
Shrinking the funnel
edge, $\ell(d)=s\downarrow0$, is therefore a degeneration inside
$\mathcal B$, and $\mathcal U=\{\lambda_s:0<s\le S\}$ has compact
closure $\{\lambda_s:0\le s\le S\}$ in $\mathcal B$. The limit surface
$Y_G(\lambda_0)$ is the disjoint union of a thrice-punctured sphere and
three elementary parabolic components, each isometric to
$\mathbb H/\langle z\mapsto z+1\rangle$.

In this notation Schulze's regularizing factor for the edge $d$ is
$\mathcal Z_{d,\lambda_s}(z)=\prod_{k\ge0}(1-e^{-(z+k)s})^2=D_s(z)$
\cite[Section~1.4 and Lemma~3.22]{Schulze}, and
$\mathcal N(\lambda_s,z)=\prod_{t<1/4}(t-z(1-z))$, the product over the
eigenvalues $t$ of $Y_G(\lambda_s)$ below $1/4$, equals $N_s(z)$: for
$s>0$ the only such eigenvalue is $\beta(s)$, by
\cite[Theorem~1.5 and Corollary~1.6]{BMM} together with
$\delta(s)>1/2$, while at $s=0$ it is the eigenvalue $0$ of the
constant function on the finite-area component, the elementary
components having no $L^2$ eigenvalues; this agrees with the limit
$\beta(s)\to0$ as $s\downarrow0$ recorded after
\eqref{eq:critical-model-energy}, that is, with the continuity of
$\beta$ through pinching \cite[Theorem~2.14]{Schulze}.
Theorem~3.21 of \cite{Schulze}, applied to the relatively compact sets
$K$ and $\mathcal U$, gives the two-sided bound. Theorem~3.18 there
gives continuity of $\lambda\mapsto\mathcal Z_\lambda/\prod_d\mathcal Z_{d,\lambda}$
on $\mathcal B$ for the topology of locally uniform convergence on
$\{\Re z>1/2\}$; in particular $Z_s/D_s$ is holomorphic on
$\{\Re z>1/2\}$ for each $s>0$. Finally, by \cite[Theorem~1.1]{BJP},
the zeros of $Z_s$ in $\{\Re z>1/2\}$ are exactly the points $z$ for
which $z(1-z)$ is an $L^2$ eigenvalue of $F_s$, with multiplicity equal
to that of the eigenvalue. By the uniqueness statement just quoted, the
only such zero is a simple zero at $z=\delta(s)$, which is cancelled by
the simple zero of $N_s$; and $D_s$ has no zeros in $\{\Re z>1/2\}$.
Hence $Z_s/(D_sN_s)$ is holomorphic and zero-free there.
\end{proof}

\begin{proposition}[Identification of the exceptional measure]
\label{prop:critical-identification}
With $x_s=f(i(\delta(s)-1/2))$,
\begin{equation}\label{eq:critical-eta-identification}
 \eta_2=\int_0^\infty w(s)\delta_{x_s}\,ds\quad\hbox{on }(1,a].
\end{equation}
In particular it has no atom at $a$, and the physical mean measure
in Proposition~\ref{prop:critical-mean} is $\nu_c$.
\end{proposition}

\begin{proof}
We keep the infinite contour tails on an absolutely convergent line.
This avoids a uniform high-frequency scattering estimate for pinching
surfaces. By Lemma~\ref{rev:lem-schulze-pinching}, for every $S>0$ the
function $Z_s/(D_sN_s)$ is holomorphic and zero-free in $\Re z>1/2$ and
bounded above and away from zero on each compact subset, uniformly for
$0<s\le S$. Cauchy's estimate for its derivative on slightly larger
compact sets, divided by the lower bound, therefore gives
\begin{equation}\label{eq:critical-zeta-regularization}
 \frac{Z_s'}{Z_s}(z)
  =\frac{2z-1}{\beta(s)-z(1-z)}+\frac{D_s'}{D_s}(z)+E_s(z),
\end{equation}
where $E_s$ is locally uniformly bounded and holomorphic, including
the pinching limit. The elementary series
$D_s'/D_s=2\sum_{k\ge1}s e^{-zks}/(1-e^{-ks})$ is
$O_K(1+|\log s|)$ on each compact $K\subset\{\Re z>1/2\}$:
split at $ks=1$ and bound the first part by a harmonic sum.
This is integrable against $w(s)\,ds$ at zero.

For $\Re z>1$ define
\[
 A(z)=\int_0^\infty w(s)\frac{Z_s'}{Z_s}(z)\,ds.
\]
The preceding bound and Lemma~\ref{lem:funnel-long-waist} give
absolute convergence. By the Euler product and
\eqref{eq:critical-orbital-identity}, $A$ is the Laplace transform
$\int_0^\infty e^{-(z-1/2)u}\,\mathcal M(du)$ of the positive
orbital measure defining $\Lambda_2$ on positive lengths. Hence
$|A(\sigma+iy)|\le A(\sigma)$ for $\sigma>1$.
For any $\sigma_0>1/2$, choose $S$ large using
Lemma~\ref{lem:funnel-long-waist}. Then
\begin{equation}\label{eq:critical-zeta-split}
 A(z)=M_S(z)+H(z),\qquad
 M_S(z)=\int_0^S w(s)\frac{2z-1}{\beta(s)-z(1-z)}\,ds,
\end{equation}
where $H$ is holomorphic on a half-plane slightly larger than
$\Re z\ge\sigma_0$. This follows from
\eqref{eq:critical-zeta-regularization} near zero, ordinary compact
family continuity for positive bounded waists, and the absolutely
convergent large-waist tail. Initially the equality holds for
$\Re z>1$; elsewhere it specifies the regularized continuation.
No growth assertion for $H$ in a whole vertical strip is needed.

Put $F(z)=f(i(z-1/2))$ and let
$B(v)=\mathscr T_2((v-x)^{-1})$, for $v\notin[0,a]$.
The moment identity \eqref{eq:critical-weighted-distribution} gives
$B(v)=\sum_{j\ge1}\Lambda_2(h^{*j})v^{-j}$ near infinity.
Fourier inversion shifted to $\Im r=\sigma-1/2$, followed by
integration against $\mathcal M$, yields, for real $v>F(\sigma)$,
\begin{equation}\label{eq:critical-cauchy-contour}
 B(v)=\frac1{2\pi i}\int_{\sigma-i\infty}^{\sigma+i\infty}
                  \frac{F(z)}{v-F(z)}A(z)\,dz.
\end{equation}
Indeed $h^{*j}(u)=(2\pi)^{-1}\int f(r+i(\sigma-1/2))^j
 e^{iru-(\sigma-1/2)u}\,dr$; set $z=\sigma-ir$ and reverse
orientation. The bound $|F(\sigma+iy)|\le F(\sigma)<v$ and
rapid vertical decay justify every interchange and the geometric sum.

The $M_S$ contribution to \eqref{eq:critical-cauchy-contour} is
\begin{equation}\label{eq:critical-model-cauchy}
 \int_0^S w(s)\frac{x_s}{v-x_s}\,ds.
\end{equation}
For fixed $s$, its resolvent factor is
$1/(z-\delta(s))+1/(z-(1-\delta(s)))$.
Shift to $\Re z=1/2$. Only the pole $\delta(s)$ is crossed, with
residue $x_s/(v-x_s)$; the new vertical integral is zero because
$F(1/2+iy)=f(y)$ is even, while the resolvent factor is odd in $y$.
Horizontal contributions vanish by rapid decay, and $v-F$ has no
zero in the strip. Integration over $0<s\le S$ is legitimate since
$\delta(s)$ stays uniformly above $1/2$ there, including its limit
$1$ at zero. This also checks the residue and orientation factors.

Fix $1<b<a$ and choose $\sigma_0>1/2$ so $F(\sigma_0)<b$.
Choose $S$ as in \eqref{eq:critical-zeta-split}, also ensuring
$\delta(s)<\sigma_0$ for $s>S$. Choose $R$ so large that $|F|<b$
on the horizontal segments at heights $\pm R$ in
$\sigma_0\le\Re z\le\sigma$, and on the two tails
$\Re z=\sigma$, $|\Im z|\ge R$.
Deform only the central segment of the $H$ integral to
$\Re z=\sigma_0$, retaining those two original infinite tails and
connecting them by the horizontal segments. On this contour $|F|<b$.
The resulting integral $D(v)$ is holomorphic for $|v|>b$:
the finite pieces are compact, and on the tails $H=A-M_S$ is bounded
and $F$ decays rapidly. Finite contour deformation for large real $v$
and then analytic continuation give
\[
 B(v)-\int_0^S w(s)\frac{x_s}{v-x_s}\,ds=D(v).
\]
Thus the Cauchy transform of
$\mathscr T_2-\int_0^S w(s)x_s\delta_{x_s}\,ds$ continues
holomorphically through $(b,\infty)$. The distributional jump formula
implies that this distribution vanishes there: the difference of its
upper and lower boundary values is $-2\pi i$ times the distribution,
by convolution with the Poisson approximate identity.
This argument includes a neighborhood of $a$, excluding endpoint
atoms and their derivatives. Divide by the nonzero coordinate $x$.
Since $s>S$ implies $x_s<b$, and $b>1$ was arbitrary,
\eqref{eq:critical-eta-identification} follows. Finally $\beta(s)>0$
for $s>0$, so there is no mass at zero energy; the factor $c^2/2$
in Proposition~\ref{prop:critical-mean} gives exactly $\nu_c$.
\end{proof}

The identification of $\eta_2$ is independent of the cusp scale and is
reused twice: in Section~\ref{sec:fixed-energy-spectrum}, where the
same pushforward gives the fixed-window density $\nu_*$ when
$\sqrt g\ll n\ll g$, and in Appendix~\ref{app:moving-edge}, where the
model energies near $1/4$ enter the moving-edge law of
Section~\ref{sec:moving-edge}.

\subsection{Localization and exhaustion of the mean}

This subsection supplies the lower half of the count matching and
completes the proof of the point-process and gap assertions of
Theorem~\ref{thm:critical-eigenvalues}. Lemma~\ref{lem:critical-funnels}
shows that short two-cusp pairs of pants have deep, disjoint funnel
neighborhoods in $X$, using the cut-graph
Lemmas~\ref{lem:critical-cuts} and~\ref{lem:critical-complex} of
Appendix~\ref{app:critical}; Lemma~\ref{rev:lem-localized-lower-count}
transplants model ground states into these neighborhoods to bound the
spectral count from below by the waist count. Together with the mean
limit of Propositions~\ref{prop:critical-mean}
and~\ref{prop:critical-identification}, this forces the actual count
and the waist count to agree in $L^1$.

\begin{lemma}[Deep disjoint model neighborhoods]\label{lem:critical-funnels}
Fix $0<d<S<\infty$ and $R<\infty$. If $n/\sqrt g\to c>0$,
with probability tending to one every two-cusp pants of waist length
in $[d,S]$ has an embedded exterior funnel neighborhood of radius $R$,
and all these enlarged pants are pairwise disjoint.
\end{lemma}

\begin{proof}
Lemma~\ref{lem:critical-complex} first excludes intersecting waists.
Expand the exterior normal collars simultaneously. Their metric is
$dr^2+s^2\cosh^2r\,d\theta^2$ with $\theta\in\mathbb R/\mathbb Z$, and
there are no focal points.
A first failure before radius $R$ supplies an embedded essential arc
of length at most $2R$ in the complement, joining one waist to itself
or two distinct waists. A boundary-parallel arc cannot cause a first
collision, since the corresponding cyclic funnel cover has embedded
normal coordinates.

A regular neighborhood of one pants and a returning arc has genus
zero, two cusps and two boundary components. Tighten its essential
boundaries. Filling a peripheral punctured disk instead gives a
genus-zero piece with at least three cusps and one essential boundary.
An unpunctured disk would make the arc boundary-parallel; both boundaries
cannot disappear when the ambient genus tends to infinity. If two
tightened boundaries coincide, retain the two incidences of the same
cut edge. An arc between distinct pants produces a genus-zero piece
with four cusps and one boundary. In all cases the total cut length
is bounded, for example by $4S+8R+1$.

There are finitely many cut-graph shapes with at most two edges.
The first type has a specified degree-two genus-zero vertex and
therefore $v\ge1$, $w\ge1$ in Lemma~\ref{lem:critical-cuts}; its
expected count is $O(g^{-1/2})$. The degree-one types have at least
three cusps and receive that lemma's extra $g^{-1/2}$ factor.
Markov's inequality excludes all of them. Only existence of such a
cut is used; no bound on the number of arcs realizing a cut is needed.
\end{proof}

We now show that the localized model states give a lower bound for the
actual eigenvalue count; together with the mean computed above, this
lower bound forces count matching. Write $N_X(t)=\Xi_X([0,t])$ and
$Z_g(t)=\#\{s\in\mathcal Q_X:\beta(s)\le t\}$, where
$\mathcal Q_X$ is the pants-waist process. At every continuity
threshold of $\nu_c$, the waist-process limit of
Lemma~\ref{lem:critical-skeleton} and the factorial-moment bounds in
its proof give convergence in distribution and uniform integrability of
$Z_g(t)$. Indeed the boundary of $\{s:\beta(s)\le t\}$ lies in
$\{0\}\cup\{\beta=t\}$, which is $\mu_c$-null; $\beta(s)\le t$ forces
$s\le S_t$, and the second moment is bounded by that of
$\mathcal Q_X([0,S_t])$.
The expected number of waists below $d$ is $O_c(d^2)$, uniformly for
large $g$, by the same integration formula.

\begin{lemma}[Localized lower count]\label{rev:lem-localized-lower-count}
Suppose $n/\sqrt g\to c\in(0,\infty)$. Then
\begin{equation}\label{eq:critical-lower-count}
 \mathbb P\{N_X(t)\ge Z_g(t-\eta)\}\longrightarrow1,
 \qquad 0<\eta<t<1/4.
\end{equation}
\end{lemma}

\begin{proof}
For each compact waist interval $[d,S]$ and $\epsilon>0$, choose
normalized compactly supported form-domain states $v_s$ on $F_s$ with
Rayleigh quotient at most $\beta(s)+\epsilon$. Their supports may
be taken in uniformly bounded funnel radius and cusp heights, and
$(\int|v_s|)^2\le C$ uniformly. To see this, approximate a ground
state at each fixed parameter in $H^1$ by a smooth compactly supported
function, transport it to neighboring marked metrics, and use continuity
of its norm, energy and $\beta$. A finite cover of $[d,S]$ supplies
common bounds; no uniformity down to $s=0$ is required.

On the event of Lemma~\ref{lem:critical-funnels}, transplant the states
of $k$ selected pants into $X$; they have disjoint supports. For
$u=\sum a_i v_i$ and the projection $P$ off the global constant,
\begin{equation}\label{eq:critical-constant-projection}
 \|Pu\|_2^2\ge
 \left(1-\frac{kC}{\operatorname{area}(X)}\right)\sum_i|a_i|^2,
 \qquad \mathcal E(Pu)=\mathcal E(u).
\end{equation}
Thus, if all selected energies are at most $t-\eta$, use trial error
$\eta/2$ and let $g\to\infty$. For bounded $k$, the projection is
injective and all projected Rayleigh quotients are below $t$, so
min--max gives $k$ positive eigenvalues, without losing a mode to the
constant. Tightness of pants counts removes the bound on $k$. Removing
the lower waist cutoff $d$ and using the deterministic upper cutoff
$S_t$ of Lemma~\ref{lem:funnel-long-waist} gives
\eqref{eq:critical-lower-count}.
More precisely, the failure probability first has limit superior
$O_c(d^2)$; let $d\downarrow0$ to obtain the displayed limit.
\end{proof}

\begin{proof}[Proof of the total-spectrum assertions in Theorem~\ref{thm:critical-eigenvalues}]
Propositions~\ref{prop:critical-mean} and
\ref{prop:critical-identification} give
$\mathbb E\Xi_X\to\nu_c$ vaguely on $[0,1/4)$. This is a two-sided
limit of the mean measure: at every continuity threshold $t$ of $\nu_c$
it gives $\mathbb EN_X(t)\to\nu_c([0,t])$. In the count matching below it
is used as the upper half, through the first term on the right of the
$L^1$ identity, while the lower half is supplied by the localized lower
count of Lemma~\ref{rev:lem-localized-lower-count}.
Take a continuity threshold $t$ and choose $t-\eta$ also to be one.
On the good event in \eqref{eq:critical-lower-count},
\[
 (Z_g(t)-N_X(t))_+\le Z_g(t)-Z_g(t-\eta).
\]
On its complement the left side is at most $Z_g(t)$, whose uniform
integrability makes that contribution vanish. First let $g\to\infty$
and then $\eta\downarrow0$. The limiting expected increment vanishes,
so $\mathbb E(Z_g(t)-N_X(t))_+\to0$. Since both expected counts
converge to $\nu_c([0,t])$, the identity
\[
 \mathbb E|N_X(t)-Z_g(t)|
 =\mathbb E N_X(t)-\mathbb E Z_g(t)
                    +2\mathbb E(Z_g(t)-N_X(t))_+
\]
proves $L^1$ matching. Finitely many thresholds give joint count
convergence; local mean bounds give tightness on $[0,1/4)$.
This proves the point-process assertion. It does not assume uniform
integrability of the actual spectral counts.

The pushforward of a Poisson random measure is Poisson, with intensity
\eqref{eq:critical-energy-intensity}. Its void probability is
\eqref{eq:critical-gap-law}. The measure is locally finite at zero
and has no atom there. Its total mass is infinite, since every
$\beta(s)<1/4$ while $\int_0^\infty\sinh(s/2)\,ds=\infty$.
It follows that $0<G_c<1/4$ almost surely. Count convergence at
continuity thresholds now gives convergence of the gap, including
the possible value $1/4$ for a finite surface without small positive
eigenvalues. The same argument gives convergence of every fixed
ordered small eigenvalue to the corresponding Poisson order statistic;
each such eigenvalue exists with probability tending to one.
\end{proof}

\subsection{The residual character of the critical spectrum}

The count matching just proved also identifies the spectral type, once
cusp averages are controlled on the whole low eigenspace; this gives
the residual assertions
\eqref{eq:critical-residual-process}--\eqref{eq:critical-cuspidal-exclusion}
of Theorem~\ref{thm:critical-eigenvalues}. The mechanism is the one of
Lemma~\ref{lem:residual-pants}, now applied to the exact low eigenspace
rather than to a trial space. For a compact annulus $A$ in a cusp,
bounded by two horocycles of length at most one, put
$h_A=\mathbf1_A/\sqrt{\operatorname{area}(A)}$. Every cuspidal
eigenfunction is orthogonal to $h_A$: its constant Fourier coefficient
vanishes throughout the cusp cylinder, by the Fourier equation.
The same holds for finite linear combinations.

\begin{lemma}[Injectivity of the cusp-average map]
\label{rev:lem-cusp-average-injective}
\textup{(a)} For every compact $[d,S]\subset(0,\infty)$ and $\zeta>0$,
there are smooth compactly supported functions $f_s$ on $F_s$ and
compact cusp annuli $A_s\subset P_s$, $s\in[d,S]$, such that
\begin{equation}\label{eq:critical-model-observation}
 \langle f_s,h_{A_s}\rangle=1,\qquad
 \|(\Delta_{F_s}-\beta(s))f_s\|_2\le\zeta.
\end{equation}
Their supports have uniformly bounded funnel radius and cusp heights.
No unit-norm condition on $f_s$ is imposed.

\textup{(b)} Let $0<b<t<1/4$. Suppose that $k$ models $F_{s_1},\dots,F_{s_k}$
with energies $\beta_i=\beta(s_i)\in[b,t-b]$, together with the functions
$f_i=f_{s_i}$ and annuli $A_i=A_{s_i}$ of \textup{(a)}, have been
transplanted isometrically into $X$ with pairwise disjoint supports, and
that $\dim E=k$ for $E=\mathbf1_{(0,t]}(\Delta_X)L^2(X)$. Let $P$ be the
orthogonal projection onto $E$, put $Fa=\sum_ia_if_i$ for
$a\in\mathbb C^k$, and $Tv=(\langle v,h_{A_i}\rangle)_{i\le k}$. Then
$TF=I$, $\|T\|\le1$, and
\begin{equation}\label{eq:critical-observation-injection}
 \|(I-P)f_i\|_2\le\zeta/b,\qquad
 \|TPF-I\|\le\sqrt{k}\,\zeta/b.
\end{equation}
If $\sqrt k\,\zeta/b<1$, then $T$ is injective on $E$, and $E$ consists
of residual eigenfunctions.
\end{lemma}

\begin{proof}
(a) To construct the functions at a fixed parameter, take the positive ground state
$u$ and divide it by $\langle u,h_A\rangle>0$. Multiply by smooth
compact cutoffs $\chi_R$ that equal one on $A$ and a growing compact
core. In logarithmic cusp height and funnel distance, choose the end
transitions with $|\nabla\chi_R|+|\Delta\chi_R|\le C/R$. Then
\[
 (\Delta-\beta)(\chi_Ru)
    =(\Delta\chi_R)u-2\langle\nabla\chi_R,\nabla u\rangle
       \longrightarrow0\quad\hbox{in }L^2,
\]
since $u$ has finite norm and energy. The observation remains one.
Transport a fixed such compact function to nearby marked metrics,
preserving the cusp coordinates on the annulus. The operator residual
depends continuously on the parameter, by continuity of $\beta$.
Renormalization by the
positive observation, if needed, preserves this continuity. A finite
cover of $[d,S]$ supplies the common support bounds in
\eqref{eq:critical-model-observation}. This argument requires neither
a globally chosen continuous ground state nor a uniform pinching limit.

(b) Since the transplanted supports are pairwise disjoint and the
transplantation is isometric, $\langle f_j,h_{A_i}\rangle=\delta_{ij}$
by \eqref{eq:critical-model-observation}, so $TF=I$; and $\|T\|\le1$
because the $h_{A_i}$ are orthonormal. Also
$\|(\Delta_X-\beta_i)f_i\|_2\le\zeta$. The spectral theorem and
\eqref{eq:critical-model-observation} give
\eqref{eq:critical-observation-injection}: the complement of $P$
includes zero and the continuous spectrum, and its distance from each
$\beta_i$ is at least $b$, so $\|(I-P)f_i\|_2\le\zeta/b$; then
$\|TPF-I\|=\|T(P-I)F\|\le\|(I-P)F\|\le\sqrt k\,\zeta/b$. If the last
bound is less than one, $TPF$ is invertible. Thus $PF$ maps onto $E$ and
$T$ is injective on $E$. Since cusp-form combinations have zero cusp
averages, no nonzero such combination can lie in $E$. The orthogonal
residual--cuspidal decomposition of each discrete
eigenspace therefore makes $E$ entirely residual, even when eigenvalues
have multiplicity. No separation between the individual $\beta_i$ is needed.
\end{proof}

\begin{proof}[Proof of \eqref{eq:critical-residual-process}--\eqref{eq:critical-cuspidal-exclusion}]
First suppose $n/\sqrt g\to c>0$. Fix a continuity threshold $t<1/4$.
The preceding total-spectrum proof gives
$\mathbb E|N_X(t)-Z_g(t)|\to0$, hence
$\mathbb P\{N_X(t)=Z_g(t)\}\to1$. Also $N_X(t)$ is uniformly
integrable: $Z_g(t)$ has bounded second moment and their difference
tends to zero in $L^1$.

Fix an arbitrarily small failure allowance $\alpha>0$. By count
tightness choose $K$ so that $Z_g(t)>K$ has limiting probability at
most $\alpha$. Choose $d>0$ so that a waist below $d$ occurs with
limiting probability at most $\alpha$. Finally choose $\eta>0$, with
$t-\eta$ a continuity point, so that a model energy in $(t-\eta,t]$
occurs with limiting probability at most $\alpha$. This follows from
mean convergence and $\nu_c(\{t\})=0$. All remaining waists lie in
$[d,S_t]$, and their energies lie in $[\beta_{\min},t-\eta]$, where
$\beta_{\min}=\min_{[d,S_t]}\beta>0$.

Put $b=\min(\beta_{\min},\eta)$ and choose
$\zeta<b/(2\sqrt{\max(K,1)})$. The common support radius in
\eqref{eq:critical-model-observation} is now fixed before $g\to\infty$.
Lemma~\ref{lem:critical-funnels} transplants these states and their
annuli into disjoint neighborhoods with probability tending to one.
On the count-matching event, $k=Z_g(t)=\dim E$, so
Lemma~\ref{rev:lem-cusp-average-injective}(b) makes $E$ residual; $k=0$
is immediate. Let $\alpha\downarrow0$. The cuspidal count in $(0,t]$
vanishes in probability, and its domination by the uniformly integrable
$N_X(t)$ gives convergence to zero in mean.

For any fixed $t_0<1/4$, choose a continuity threshold $t\in(t_0,1/4)$.
Domination gives the same mean conclusion at $t_0$, without assuming
that $\nu_c$ is atomless. Subtracting the cuspidal measure from the
total one proves the joint process and mean assertions in
\eqref{eq:critical-residual-process}.

If \eqref{eq:critical-cuspidal-exclusion} failed, there would be a
sequence $1\le n\le M\sqrt g$ along which the expected cuspidal count
stays bounded away from zero. Pass to $n/\sqrt g\to c\in[0,M]$.
For $c>0$ the preceding argument is a contradiction. For $c=0$, the local count bound of
Appendix~\ref{app:critical-trace}, applied on the envelope
$n+1\le(M+1)\sqrt g$, gives
$\mathbb EN_X(1/4-\varepsilon)\le C_{M,\varepsilon}(n+1)^2/g\to0$
(compare Huang's total-count bound~\cite[Theorem~1.1]{HuangGap}),
again a contradiction.
\end{proof}

\subsection{The three cusp regimes}

The corollary below assembles Huang's subcritical theorem,
Theorem~\ref{thm:critical-eigenvalues} and
Corollary~\ref{cor:spectral} into the gap phase diagram announced in
the introduction. Its two endpoints are refined in the next sections:
the approach to $1/4$ from below in Section~\ref{sec:moving-edge}, and
the rate of vanishing above the critical scale in
Section~\ref{sec:supercritical-spectrum}.

\begin{corollary}[Gap phase diagram]\label{cor:gap-phase-diagram}
For $n=n(g)\ge1$, as $g\to\infty$,
\[
 G(X)\xrightarrow{\mathbb P}1/4\ \Longleftrightarrow\ n/\sqrt g\to0,
 \qquad
 G(X)\xrightarrow{\mathbb P}0\ \Longleftrightarrow\ n/\sqrt g\to\infty.
\]
At $n/\sqrt g\to c\in(0,\infty)$ the limit is the nondegenerate
law \eqref{eq:critical-gap-law}. In particular there is no deterministic
limit in $(0,1/4)$, for any sequence of cusp counts.
\end{corollary}

\begin{proof}
The first equivalence is Huang's result
\cite[Theorem~1.1 and Corollary~6.15]{HuangGap}.
The forward cusp condition in the second implies gap convergence by
Corollary~\ref{cor:spectral}. Conversely, if $n/\sqrt g$ does not
tend to infinity, pass to a subsequence with a finite limit.
A zero limit gives gap limit $1/4$ by Huang; a positive limit gives
$G_c>0$ by Theorem~\ref{thm:critical-eigenvalues}. Both contradict
convergence to zero. For a proposed deterministic interior limit,
pass to a subsequence on which $n/\sqrt g$ converges in $[0,\infty]$.
The endpoints have already been treated. At a positive finite limit,
$\nu_c([0,t])$ is finite and strictly positive for every $t>0$:
positivity follows from $\beta(s)\to0$. Thus the void probability
lies strictly between zero and one, so the limiting law is not a point
mass. This also proves the claimed nondegeneracy.
\end{proof}

\section{Fluctuations near the quarter threshold}\label{sec:moving-edge}

Huang's theorem gives $G(X)\to1/4$ whenever $n=o(\sqrt g)$,
with bounds in each fixed spectral window~\cite[Theorem~1.1]{HuangGap}.
This section determines the fluctuations of the spectrum near $1/4$
in a range approaching the critical scale from below: the Fr\'echet
law of the gap defect and the Poisson law of all deficits on the scale
$n^2/g$ (Theorem~\ref{thm:moving-edge}), the residual type of the
corresponding eigenfunctions (Corollary~\ref{cor:edge-residual}) and
the moments of the defect (Corollary~\ref{cor:edge-moments}). The
inputs are the funnel model of Section~\ref{sec:critical-eigenvalues},
whose energies near $1/4$ are read off from the critical exponent of
Hecke groups~\cite{Soares}, and Huang's trace expansion, restated and
made uniform in a moving window in Appendix~\ref{app:moving-edge},
where the proofs are given. In this range the relevant pants waists
grow and the model eigenvalues approach the continuous spectrum, so
neither substituting a varying $c$ into
Theorem~\ref{thm:critical-eigenvalues} nor convergence of the geometric
counts alone gives the result below.

Write $d(w)$ for the critical exponent of the Hecke group generated
by $z\mapsto z+w$ and $z\mapsto-1/z$, $w>2$. Its index-two
torsion-free subgroup generated by $z\mapsto z+w$ and
$z\mapsto z/(1-wz)$ uniformizes a surface with two cusps and one
funnel~\cite[Sections~1 and~3]{Soares}; it is the funnel model $F_s$
when $w=2\cosh(s/4)$. Indeed, in terms of
the generators $A,B$ of Lemma~\ref{lem:funnel-long-waist}, with
$a=2\cosh(s/4)$, the waist is represented by $AB$ and
$\operatorname{tr}(AB)=2-a^2$, so
$|\operatorname{tr}(AB)|=a^2-2=2\cosh(s/2)$, the trace of a hyperbolic
element of translation length $s$; the translation $z\mapsto z+w$ of
the Hecke group is $A$ with $a=w$. Hence, with
$\kappa_0=(8\pi^2)^{-1}$,
\begin{equation}\label{eq:edge-model-coordinate}
 q(w)=(d(w)-1/2)^2,\qquad q(2\cosh(s/4))=1/4-\beta(s),\qquad
 h_*(y)=q(\sqrt{4+y}),\quad y\ge0.
\end{equation}
Continuity at $y=0$ is understood. The function $h_*$ decreases
strictly from $1/4$ to zero: $d(2)=1$, $d(w)\to1/2$ as $w\to\infty$,
and $w\mapsto d(w)$ is strictly decreasing on $[2,\infty)$, see
\cite[equation~(9) and the references therein]{Soares}. Soares's asymptotic
$d(w)=1/2+w^{-1}+O(w^{-2}\log w)$~\cite[Theorem~1.9]{Soares} gives
\begin{equation}\label{eq:edge-model-tail}
 h_*(y)=y^{-1}+O(y^{-3/2}\log y),\qquad
 h_*(y)\le\frac{C}{1+y}.
\end{equation}
In particular $h_*-(1+y)^{-1}$ is integrable. Define
$H(t)=|\{y\ge0:h_*(y)\ge t\}|$ for $0<t<1/4$.

\begin{theorem}[The moving spectral edge]\label{thm:moving-edge}
Let $X$ have normalized Weil--Petersson law on $\mathcal M_{g,n}$,
and suppose
\begin{equation}\label{eq:moving-edge-range}
 c=\frac n{\sqrt g}\longrightarrow0,\qquad c\log g\longrightarrow\infty.
\end{equation}
Counting every positive discrete eigenvalue with multiplicity, put
\[
 \mathcal D_X=\sum_{0<\lambda_j(X)<1/4}
       \delta_{(1/4-\lambda_j(X))/c^2}.
\]
Then $\mathcal D_X$ converges in distribution, vaguely on $(0,\infty)$,
to a Poisson random measure with intensity
\begin{equation}\label{eq:edge-poisson-intensity}
 \nu_{\rm edge}(dx)=\kappa_0 x^{-2}\,dx.
\end{equation}
The convergence also holds jointly for finitely many upper-tail counts
$[a,\infty)$, $a>0$, and their means converge to $\kappa_0/a$.
For the gap defect $D(X)=1/4-G(X)$,
\begin{equation}\label{eq:edge-gap-law}
 \frac{D(X)}{c^2}\Longrightarrow Z,\qquad
 \mathbb P(Z\le x)=\exp(-\kappa_0/x),\quad x>0.
\end{equation}
For each fixed $j$, the first $j$ positive eigenvalues below $1/4$
exist and are simple with probability tending to one. Their ordered
deficits, divided by $c^2$, converge jointly to
$(\kappa_0/T_1,\ldots,\kappa_0/T_j)$, where $T_i$ are the arrival
times of a rate-one Poisson process. Their normalized eigenfunctions
are, up to sign, $o_{\mathbb P}(1)$ in $L^2$ from the ground states of
the funnel models $F_s$ of the corresponding pairs of pants, cut off in
the funnel and transplanted to embedded neighborhoods of these pants in
$X$; the neighborhoods are mutually disjoint and have area at most a
deterministic $g^{o(1)}$.
\end{theorem}

The law is Fr\'echet with shape parameter one. Its geometric origin
is transparent in the coordinate $y=2(\cosh(s/2)-1)$:
the model intensity is $\kappa_0c^2\,dy$, while the deficit is
$h_*(y)\sim y^{-1}$. The substantial step is the transfer of this
model computation to the \emph{entire actual spectral count}.
Appendix~\ref{app:moving-edge} proves the moving mean, conditional
isolation and exact count matching needed for this transfer, as well
as the eigenfunction assertion.
The condition $c\log g\to\infty$ enters through localization: a model
ground state with deficit of order $c^2$ decays on the funnel only over
lengths of order $1/c$, so it must be cut off at a radius $R$ with
$cR\to\infty$, while our estimates for isolation of the embedded funnel
neighborhoods require $R=o(\log g)$.
Condition~\eqref{eq:moving-edge-range} includes
$n=\lfloor\sqrt g/(\log g)^\alpha\rfloor$ for $0<\alpha<1$, but
does not cover the whole subcritical range. The Poisson assertion
does not describe deficits of order $o(c^2)$. It does compute, for
instance, the probability that fewer than $j$ modes fall below a
prescribed moving energy: for each fixed $a>0$ and positive integer $j$,
\[
 \mathbb P\{\mathcal D_X([a,\infty))<j\}
 \longrightarrow e^{-\kappa_0/a}
       \sum_{r=0}^{j-1}\frac{(\kappa_0/a)^r}{r!}.
\]

\begin{corollary}[Spectral type at the moving edge]\label{cor:edge-residual}
Under the stronger condition
\begin{equation}\label{eq:edge-residual-range}
 c\longrightarrow0,\qquad
 \frac{\log(1/c)}{c\log g}\longrightarrow0,
\end{equation}
the expected cuspidal count in $(0,1/4-ac^2]$ tends to zero for
every fixed $a>0$. The residual and cuspidal deficit processes
converge jointly to $(\operatorname{PPP}(\nu_{\rm edge}),0)$, including their
upper-tail means. In particular each fixed initial list of small
eigenvalues is residual with probability tending to one. If
$G_{\rm cusp}$ is the smallest positive cuspidal eigenvalue below
$1/4$, or $1/4$ if none exists, then
$(1/4-G_{\rm cusp})/c^2\to0$ in probability.
\end{corollary}

The stronger condition controls the error relative to a cusp
observation that itself tends to zero; without this relative control,
localization does not distinguish residual eigenfunctions from cusp
forms. The proof is in Subsection~\ref{subsec:edge-residual-proof}.

\begin{corollary}[Rare gaps and moments]\label{cor:edge-moments}
For every sequence $n\to\infty$, $n=o(\sqrt g)$ and every fixed
$0<t<1/4$, let $N_X^{\rm phys}(t)$ count the positive eigenvalues
at most $1/4-t$. Then
\begin{equation}\label{eq:edge-rare-law}
 \mathbb P\{D(X)\ge t\}\sim\mathbb E N_X^{\rm phys}(t)
       \sim\kappa_0\frac{n^2}{g}H(t).
\end{equation}
Conditioned on this event, $N_X^{\rm phys}(t)=1$ with probability
tending to one. Under~\eqref{eq:moving-edge-range}, for each fixed
$p>0$ the moments of the actual gap defect satisfy
\begin{align}
 \frac{\mathbb E D^p}{c^{2p}}&\longrightarrow
       \kappa_0^p\Gamma(1-p),&&0<p<1,\label{eq:edge-small-moments}\\
 \mathbb E D&=\frac{c^2}{4\pi^2}\log(1/c)+O(c^2),
       \label{eq:edge-mean}\\
 \frac{\mathbb E D^p}{c^2}&\longrightarrow
       K_p:=\kappa_0\int_0^\infty h_*(y)^p\,dy,&&p>1.
       \label{eq:edge-large-moments}
\end{align}
Consequently $\operatorname{Var}(D)\sim K_2c^2$ and
$D/\mathbb E D\to0$ in probability.
\end{corollary}

Thus the typical defect has order $c^2$, the mean carries an
additional logarithm, and higher moments are governed by rare surfaces.
These conclusions require expectation-level estimates for discarded
geometries, since the weak limit in~\eqref{eq:edge-gap-law} has
infinite mean. Subsection~\ref{subsec:edge-moment-proof} proves the
moment statements directly, rather than exchanging a limit and an
expectation.

For comparison, Hide--Macera--Thomas discuss a Tracy--Widom prediction
for suitable random hyperbolic surface models after their
Problem~1.3~\cite[Introduction, Motivation]{HMT}; their main spectral
theorem concerns closed Weil--Petersson surfaces. The cusp-driven law
above has a different mechanism and parameter range. It resolves
neither that prediction nor the existence of closed surfaces of
arbitrarily large genus with gap at least $1/4$.

This completes the description of the lowest eigenvalues at and below
the critical scale. The remaining sections of this part treat
$n/\sqrt g\to\infty$, where the gap tends to zero
(Corollary~\ref{cor:spectral}) and the questions become the
microscopic law of the lowest eigenvalues
(Sections~\ref{sec:supercritical-spectrum}
and~\ref{sec:dense-spectrum}) and the density of eigenvalues in fixed
windows (Section~\ref{sec:fixed-energy-spectrum}).

\section{Supercritical spectral laws}
\label{sec:supercritical-spectrum}

When $n/\sqrt g\to\infty$, the short separating curves of
Theorem~\ref{thm:pants-poisson} produce eigenvalues tending to zero,
and Corollary~\ref{cor:many-spectrum} already gives their order. This
section identifies the law of the lowest eigenvalues when $n/g$ stays
bounded (Theorem~\ref{thm:supercritical-spectrum}): on the explicit
scale $A_g$ they form a Poisson process, all planar cusp types
contributing with weights given by the squares of their areas. The
work is to exclude every low mode not produced by a short curve,
including modes supported on a large component separated by several
curves; this is done by comparing the low spectrum with a weighted
star whose rates are collar conductances divided by enclosed areas.
The inputs are the type process of Theorem~\ref{thm:pants-poisson},
the all-marking volume envelope \eqref{eq:gevrey-envelope} of
Part~\ref{part:volumes}, and two appendices: the cut lemmas of
Appendix~\ref{app:super-cuts} and the surface--graph comparison of
Appendix~\ref{app:super-comparison}. The proof itself is given in
this section.

Retain $q_g$ and $w_k$ from \eqref{eq:PS1} and \eqref{eq:PC1}.
Let $x_g$ solve $(2g-1)\Psi(x_g)=n-1$, and put
\begin{equation}\label{eq:super-scale}
 M_2(x)=\sum_{k\ge2}(k-1)^2w_k(x),\qquad
 A_g=2\pi^2\sqrt{q_gM_2(x_g)}.
\end{equation}
The function $M_2$ is finite and continuous on $[0,x_c)$, with $M_2(0)=1$.
Indeed, writing $F=\Phi_0'$ and $D=x\partial_x$, it equals
$2(D-1)^2F/x^2=2/a(Y(x))-2Y(x)/x+W(x)$, extended continuously to zero.
Here $a(Y)=X'(Y)=J_0(2\sqrt Y)$ is the kernel coefficient of
Section~\ref{the-positive-compact-family-of-kernels}, and $Y(x)$
increases from $0$ to the first zero $y_0$ of $a$ as $x$ increases from
$0$ to $x_c$. Hence $a(Y(x))\to0$ and $M_2(x)\to\infty$ as
$x\uparrow x_c$. The hypothesis $\sup_gn/g<\infty$ in the theorem
below keeps $x_g$ in a compact subinterval of $[0,x_c)$, on which
$M_2$ is bounded and the cusp weights $w_k$ have a uniform exponential
tail; this compactness is used throughout the proof.

\begin{theorem}[Supercritical spectral process]\label{thm:supercritical-spectrum}
Suppose $g\to\infty$, $n/\sqrt g\to\infty$ and $\sup_g n/g<\infty$.
For $X$ with law $\mathbb P^{\mathrm{WP}}_{g,n}$,
\begin{equation}\label{eq:super-PPP}
 \sum_{0<\lambda_j(X)<1/4}\delta_{A_g\lambda_j(X)}
 \ \Longrightarrow\ \operatorname{PPP}(2t\,dt)
 \quad\text{on }[0,\infty)
\end{equation}
in distribution for the vague topology, counting multiplicities.
In particular, for every $t\ge0$,
\begin{equation}\label{eq:super-gap}
 \mathbb P\{A_gG(X)>t\}\longrightarrow e^{-t^2},
 \qquad \frac{G(X)}{\mathfrak h(X)}
                  \longrightarrow\frac1\pi\quad\text{in probability}.
\end{equation}
Here $G$ is the spectral gap, with the convention of
Section~\ref{sec:critical-eigenvalues}. If $\sqrt g\ll n\ll g$, then
$A_g\sim\pi n/(2\sqrt{2g})$.
\end{theorem}

Marked-geodesic and variational methods of this kind have precedents in
Hide--Thomas~\cite{HTSphere} and Gatignol--Thomas~\cite{GT}.
The additional conclusion here is a two-sided spectral comparison on
a shrinking energy scale. Two-sided estimates of small eigenvalues in
terms of short separating geodesics go back to Schoen, Wolpert and
Yau~\cite{SWY} for compact surfaces; see also Dodziuk, Pignataro,
Randol and Sullivan~\cite{DPRS}, and Wu and Xue~\cite{WXOptimal} for
lower bounds uniform in the genus. Burger's surface--graph comparison
\cite[Theorems~1.1--1.2]{Burger} concerns degeneration with constants
depending on the topological signature; we instead control the
pieces, traces and intermediate branches uniformly on events of
probability tending to one. The all-marking volume envelope
\eqref{eq:gevrey-envelope} is essential for excluding competing cuts.

The convergence \eqref{eq:super-PPP} is asserted only in distribution.
Unlike Theorem~\ref{thm:moving-edge}, whose statements about means
rest on a positive trace estimate
(Lemma~\ref{rev:lem:edge-trace-approximation}), no positive trace
input is available in the present range, so no convergence of spectral
first moments is claimed.
Theorem~\ref{thm:dense-persistence} extends the law to a specified
range with $n/g\to\infty$, at the cost of endpoint estimates and a more
general contraction argument.

\begin{proof}[Proof of Theorem~\ref{thm:supercritical-spectrum}]
Fix $M$ with $n\le Mg$, and write $P=2\pi^2$ and
$Q=nx_g/(4P)\sim q_g$. All constants may depend on $M$.
The proof has three steps: obtain uniform gaps on cut pieces, compare
the low spectrum with a weighted star, and transfer the marked-curve law.

\emph{The cut pieces.} Set
\begin{equation}\label{eq:super-cutoffs}
 \delta=Q^{-1/10},\quad b=Q^{-1/20},\quad
 K=C_M\log(Q+2),\quad R=Q^{1/20},\quad \tau=R/\sqrt Q.
\end{equation}
By Lemma~\ref{lem:super-geometry}, with probability tending to one,
all curves of length at most $\delta$ are planar separating. Cutting
along all of them gives a tree with a genus-$g$ root. Each root branch
has at most $K$ planar vertices and area at most $2\pi K$; every
geodesic piece has Neumann gap at least
\begin{equation}\label{eq:super-piece-gap}
 \gamma=c_M\delta^2/\log^2(Q+2).
\end{equation}

Trim the half-collars at circumference $b$. Write $m_v$ for the
remaining piece areas, $W_e$ for their total mass below an edge, and
$\kappa_e$ for the full collar conductance, as in
Lemma~\ref{lem:super-comparison}. Retain edges for which
$d_e=\kappa_e/W_e\le\tau$. If an edge encloses $k$ cusps, then
$W_e\le2\pi(k-1)$ and $\kappa_e\ge\ell_e/\pi$, so a retained edge
satisfies $\sqrt Q\ell_e/(k-1)\le PR$. The isolation assertion of
Lemma~\ref{lem:super-geometry}, applied with $PR$, makes each retained
edge an isolated root leaf. Other branches need not be unnested.

\emph{Spectral comparison.} Contract all nonretained vertices into the
root and use their total trimmed mass as its vertex mass. Let $\mu_j$
be the positive eigenvalues of this weighted star. For fixed $t$,
Lemma~\ref{lem:super-comparison} applies uniformly for
$0\le\Lambda\le t/\sqrt Q$. Its trace constant satisfies
$T\le C_MK/(b\gamma)=O(Q^{1/4}\log^3(Q+2))$. Thus
\[
 \frac{\Lambda}{\gamma}=O_t(Q^{-3/10}\log^2(Q+2)),\qquad
 \frac{K\Lambda}{\tau}=O_t(Q^{-1/20}\log(Q+2)),
\]
and $KT\Lambda=O_t(Q^{-1/4}\log^4(Q+2))$, while
$\tau T=O(Q^{-1/5}\log^3(Q+2))$. Every error tends to zero.
Consequently, with $c_{g,t}\to1$,
\begin{equation}\label{eq:super-count-sandwich}
 N_\star^+(\Lambda)\le N_X^+(\Lambda)
                  \le N_\star^+(c_{g,t}\Lambda)
 \quad(0\le\Lambda\le t/\sqrt Q).
\end{equation}
Counts are inclusive at the upper endpoint and exclude zero.

For a retained leaf with $k$ cusps, the area and conductance estimates give
\begin{equation}\label{eq:super-rate}
 d_e=(1+O(b+\delta/b))\frac{\ell_e}{P(k-1)}.
\end{equation}
The one-waist estimate \eqref{eq:super-one-waist} implies
\[
 \mathbb E\sum_{\text{retained caps}}\operatorname{area}(\text{cap})
 \le C Q\tau^2\sum_{k\ge2}(k-1)^3w_k(x_g)
                           e^{P\tau(k-1)}=O(R^2).
\]
It is enough to bound by all eligible waists; nesting only overcounts.
Since $Q=O(g)$, $R^2=o(g)$, and the full trimmed mass is
$(1+O(b))\operatorname{area}(X)$, the retained leaf mass fraction
$\rho_g$ tends to zero in probability. Ordering the retained rates,
\begin{equation}\label{eq:super-star-rates}
 d_j\le\mu_j\le d_j/(1-\rho_g).
\end{equation}
To see this, quotient the star by constants and use leaf--root
differences $z_i$. Its mass is
$\sum m_i z_i^2-(\sum m_i z_i)^2/\sum_v m_v$, between
$(1-\rho_g)\sum m_i z_i^2$ and $\sum m_i z_i^2$; its energy is
$\sum\kappa_i z_i^2$. Min--max gives \eqref{eq:super-star-rates}.

\emph{The limiting counts.} First take a subsequence with $n/g\to\lambda<\infty$.
For any fixed set of cusp types, Theorem~\ref{thm:pants-poisson}
maps under $(s,k)\mapsto s/[P(k-1)]$ to a Poisson process.
The expected contribution of types $k>K_0$ on a bounded score interval
is at most
\[
 C_t\sum_{k>K_0}(k-1)^2w_k(x_g e^{C_t/\sqrt Q}),
\]
by \eqref{eq:super-one-waist}. This tends uniformly to zero as
$K_0\to\infty$, since the arguments stay below a fixed $x_2<x_c$.
The estimate includes waists longer than $\delta$. For a fixed type
and bounded score the waist is eventually shorter than $\delta$, and
its rate is retained because $R\to\infty$. Thus the limiting rate
intensity is $2P^2M_2(x_\lambda)t\,dt$ in the $\sqrt Q$ scale.

Equations~\eqref{eq:super-count-sandwich}--\eqref{eq:super-star-rates}
transfer joint counts at any finite collection of thresholds. The
limiting intensity is atomless, so counts in intervals whose endpoints
differ by factors tending to one vanish in probability. Geometric
first-moment bounds and the comparison on events of probability tending
to one give spectral count tightness in probability on compact intervals.
This proves point-process convergence on $[0,\infty)$, including control
at zero. Continuity of $M_2$ and compactness of the bounded density
interval remove the subsequence assumption after the normalization
\eqref{eq:super-scale}.

The first-point law gives the gap assertion, since the probability of
having no positive subthreshold eigenvalue tends to zero by first taking
$g\to\infty$ and then $t\to\infty$. For the ratio assertion, the
uniformity in \eqref{eq:super-count-sandwich} permits evaluation at
$\Lambda=\lambda_j(X)$ whenever $\mu_j\le t/\sqrt Q$; hence
$\mu_j/c_{g,t}\le\lambda_j(X)\le\mu_j$.
The first rate is tight and bounded away from zero in probability in
this scale. Lemma~\ref{lem:super-geometry} identifies the geodesic
Cheeger minimum with the same minimum of $\ell/[2\pi(k-1)]$.
Here the geodesic Cheeger constant $H$ of~\cite{SW} is an infimum over
separating geodesic multicurves that cut the surface into two connected
pieces, which are the cuts treated by Lemma~\ref{lem:super-cuts};
allowing disconnected sides does not lower this infimum, by the
graph-averaging reduction in the proof of
Lemma~\ref{lem:relative-geodesic} (Appendix~\ref{app:super-comparison}),
which selects, among the same geodesics, a cut with both sides
connected and no larger ratio.
Finally $H/(1+H)\le\mathfrak h\le H$ for the geodesic Cheeger
constant $H$~\cite[Proposition~13]{SW}, and $H\to0$ in probability.
This proves the coupled ratio $G/\mathfrak h\to1/\pi$.
The stated sparse-density constant follows from
$q_g\sim n^2/(32\pi^2g)$ and $M_2(x_g)\to1$.
\end{proof}

The scale $A_g$ and the star comparison are carried into
Section~\ref{sec:dense-spectrum}, where the law is extended to a range
with $n/g\to\infty$ and shown to fail at very large density. The
same regime $\sqrt g\ll n\ll g$, observed in a fixed-energy window
rather than on the scale $A_g$, is treated in
Section~\ref{sec:fixed-energy-spectrum}.

\section{The low spectrum at growing cusp density}
\label{sec:dense-spectrum}

This section follows the microscopic law of
Theorem~\ref{thm:supercritical-spectrum} into the dense regime
$n/g\to\infty$. The law persists in a range specified in
Theorem~\ref{thm:dense-persistence}, which contains $n\asymp g^\alpha$
for every $1<\alpha<11/9$, but not at all densities: once the planar
components themselves grow, their internal geometry can produce
eigenvalues smaller than the conductance of their attaching collars
divided by their areas, and Theorem~\ref{thm:dense-collapse} proves
that this collapse occurs whenever $n\gg g^3/(\log g)^4$, with no
upper restriction on the cusp count. Independently of either range,
Corollary~\ref{cor:dense-residual} shows that the lowest eigenspaces
are residual throughout the dense regime. The inputs are the
large-density endpoint expansions of Part~\ref{part:volumes}, which
give the growing-cap law (Lemma~\ref{lem:dense-caps},
Appendix~\ref{app:dense-collapse-details}), the star comparison of
Section~\ref{sec:supercritical-spectrum} in the more general form of
Appendix~\ref{app:dense-comparison}, and, for the collapse, the
low-energy cusp cutoffs of Appendices~\ref{app:cusp-cutoffs}
and~\ref{app:cusp-opening}. Below $g^3$, the geometric cap process
still has its expected Poisson limit.

\begin{theorem}[Dense spectral collapse]\label{thm:dense-collapse}
Let $g\to\infty$ and let $n=n(g)$ be any integer sequence satisfying
\begin{equation}\label{eq:dense-window}
 \frac{g^3}{(\log g)^4}\ll n.
\end{equation}
For $X$ with law $\mathbb P^{\mathrm{WP}}_{g,n}$,
\begin{equation}\label{eq:dense-collapse}
 A_g\lambda_r^+(X)\longrightarrow0
 \quad\text{in probability for every fixed }r\ge1.
\end{equation}
Here $A_g$ is defined by \eqref{eq:super-scale}, and $\lambda_r^+$ is
the $r$th positive eigenvalue below $1/4$, set to $1/4$ when fewer
than $r$ such eigenvalues exist. Consequently, the spectral counting
measures in \eqref{eq:super-PPP} are not tight on $[0,\infty)$ along
any sequence satisfying \eqref{eq:dense-window}.
\end{theorem}

This is a shrinking-scale conclusion. Hide--Thomas~\cite[Theorem~1.3]{HTMany}
prove a deterministic linear lower bound on the number of eigenvalues
below $1/4$ whenever $g<an$, for each fixed $a>0$. The classical genus
bound of Grigor'yan--Netrusov--Yau~\cite[Corollary~5.6]{GNY}
already gives \eqref{eq:dense-collapse} when $n/g^3\to\infty$; the
additional content is collapse below and through the scale $g^3$,
where that genus bound alone is insufficient. The small-eigenvalue
estimates of Gatignol--Thomas~\cite[Theorem~1.10]{GT} apply when
$n/g\to0$.

Theorem~\ref{thm:dense-collapse} neither gives a replacement law for all
sequences with $n/g\to\infty$ nor locates a sharp transition. In
particular, the collapse argument does not address a fixed power
$n=g^\alpha$ with $1<\alpha<3$; Theorem~\ref{thm:dense-persistence}
gives the matching Poisson law for $1<\alpha<11/9$.
The logarithmic exponent in \eqref{eq:dense-window} is a sufficient
condition, not a claimed transition exponent. Finally, non-tightness
on $[0,\infty)$ does not rule out a vague limit on $(0,\infty)$, since
points may escape through zero in the latter space.

The proof places low-energy cutoffs around a typical cusp of a large
WP sphere. Appendix~\ref{app:cusp-cutoffs} constructs them
quantitatively, by conformal regularization and the finite planar
resistance estimate of Gurel-Gurevich--Nachmias~\cite[Lemma~3.6]{GGN},
and Appendix~\ref{app:cusp-opening} transfers them to the actual
bordered cap law using the positive tree encoding of
Budd--Meeusen--Zonneveld~\cite{BMZ}. We first establish the
growing-cap law from the uniform volume estimates.

\subsection{Growing planar caps}
At growing density a large planar cap cut off by a short curve
contains about $\epsilon^{-1}\asymp n/g$ cusps, and the spectral scale
is no longer set by the two-cusp pairs of pants of
Section~\ref{sec:geometry} alone. This subsection records the parameters of this regime and the joint law of
cap waists and cap sizes (Lemma~\ref{lem:dense-caps}), which is the
geometric input of both the collapse and the persistence theorems.
Put $P=2\pi^2$, let $(2g-1)\Psi(x)=n-1$, and set
\begin{equation}\label{eq:dense-parameters}
 \epsilon=1-x/x_c,\quad Q=\frac{nx}{4P},\quad
 M_a(x)=\sum_{k\ge2}(k-1)^aw_k(x),\quad
 S=\sqrt{QM_2(x)},\quad \eta=\frac1{n\epsilon^{3/2}}.
\end{equation}
The endpoint expansions, which require only $n/g\to\infty$, give
\begin{equation}\label{eq:dense-scales}
 \epsilon\sim\frac{5g}{2n},\qquad
 M_2(x)\sim\sqrt{\frac{2}{bx_c}}\epsilon^{-1/2},\qquad
 S\sim\frac{(y_0/5)^{1/4}}{2\pi}n^{3/4}g^{-1/4},
\end{equation}
where $b=-X''(y_0)=x_c/y_0$ as in Lemma~\ref{lem:parameter}(iii). Also $Q/q_g\to1$, so $A_g\sim PS$,
and $M_a(x)=O(\epsilon^{-(a-3/2)})$ for $a\ge2$.

\begin{lemma}[Cap process with growing cusp count]\label{lem:dense-caps}
Suppose $g\ll n\ll g^3$.
For each simple geodesic bounding a genus-zero cap with $k$ cusps,
assign score $S\ell/(k-1)$ and mark $\epsilon(k-1)$.
On bounded score intervals and compact mark intervals in $(0,\infty)$,
their counting measure converges in distribution to a Poisson process
with intensity
\begin{equation}\label{eq:dense-cap-intensity}
 2t\,dt\,\frac{v^{-1/2}e^{-v}}{\sqrt\pi}\,dv.
\end{equation}
The expected number of nested pairs in any such window is $O(\eta)$.
\end{lemma}

The mark records the size of a cap on its natural scale
$\epsilon^{-1}$, and the score normalizes its waist length by the
enclosed area. The nested-pair bound permits simultaneous use of
disjoint caps in the spectral comparison. The factorial-moment proof
is given in Appendix~\ref{app:dense-collapse-details}.

\subsection{Cutoffs and the spectral consequence}
The new mechanism is that a large planar cap can support a low-energy
function whose cost is controlled inside the cap rather than by the
conductance of its attaching collar. The argument has three cases.
Below $g^3$, the cap process above supplies the locations
and sizes of the caps, and a cutoff at a typical cusp, transported to
the bordered cap, gives the required trial function. There the
decisive comparison is between the cutoff energy $O(1/\log K)$ inside
a large planar component and the scale $\ell_*$ of its attaching
waist. Here $K=\epsilon^{-1}$ is the typical cusp count of a large cap
and $\ell_*=K/S$ its typical waist length, with $\epsilon$ and $S$ as in
\eqref{eq:dense-parameters}. When $n/g^3\to0$ in the range of
Theorem~\ref{thm:dense-collapse}, one has $\ell_*\log K\to\infty$, so the
cutoff energy is much smaller than the waist scale. At the cubic scale we count smaller caps to recover a diverging
supply of disjoint supports; above it the classical genus bound
applies. The complete argument, including the exceptional-cap
estimates and the centering of the trial space, is in
Appendix~\ref{app:dense-collapse-details}.

This construction is variational. Determining the actual dense
spectral scale requires matching lower estimates and control of the
other low modes; neither follows from the construction.

\subsection{Persistence of the Poisson spectrum}

There is nevertheless an unbounded-density range in which collar
conductances still determine all low eigenvalues. Write
$L=\log(n+2)$ and retain the exact parameters
\eqref{eq:dense-parameters}.

\begin{theorem}[Poisson spectrum at growing density]\label{thm:dense-persistence}
Suppose $g\to\infty$, $n/g\to\infty$ and
\begin{equation}\label{eq:dense-persistence-range}
 \frac{L^{10}}{\epsilon^3S}\longrightarrow0.
\end{equation}
For $X$ with law $\mathbb P^{\mathrm{WP}}_{g,n}$,
\begin{equation}\label{eq:dense-persistence-PPP}
 \sum_{0<\lambda_j(X)<1/4}\delta_{PS\lambda_j(X)}
       \ \Longrightarrow\ \operatorname{PPP}(2t\,dt)
       \quad\text{on }[0,\infty)
\end{equation}
in distribution for the vague topology, with multiplicities and
residual eigenvalues included. The normalization $PS$ may be replaced
by $A_g$. In particular, for every $t\ge0$,
\[
 \mathbb P\{A_gG(X)>t\}\longrightarrow e^{-t^2},\qquad
 \frac{G(X)}{\mathfrak h(X)}\longrightarrow\frac1\pi
       \quad\text{in probability}.
\]
Every fixed initial list of rescaled positive eigenvalues converges
jointly to the corresponding ordered points of this Poisson process.
\end{theorem}

By \eqref{eq:dense-scales}, condition
\eqref{eq:dense-persistence-range} is equivalent to
$n^{9/4}\log^{10}(n+2)/g^{11/4}\to0$; in particular, it includes
$n\asymp g^\alpha$ for every fixed $1<\alpha<11/9$.
The exponent $11/9$ has a simple origin in the proof. The spectral
comparison requires the low eigenvalues, of size $1/S$ by
\eqref{eq:dense-scales}, to be small relative to the Neumann gaps
$\gamma$ of the pieces cut out by short planar curves; these gaps are
obtained from a Cheeger inequality and carry its square loss,
$\gamma=cu^2$ with $u=c\delta/K$, where $\delta=1/L$ is the cut length
and $K=C_1L/\epsilon$ is the maximal cusp number of a $\delta$-short
cap, so that $\gamma\asymp\epsilon^2L^{-4}$. The square accounts for
$\epsilon^{2}$; the attachment errors carry one further factor $K$,
the area scale of a branch, which accounts for $\epsilon^{-1}$; and $S$
enters through $\Lambda\le T/S$ and $\tau=R/S$. The critical quantity is
therefore $\epsilon^3S\asymp g^{11/4}n^{-9/4}$. A gap estimate for
the root piece that is linear rather than quadratic in its Cheeger
constant would improve the exponent, which we do not claim to be sharp.
Neither this exponent nor the logarithmic power is a proposed phase
boundary, and the theorem leaves unresolved a substantial interval
below the collapse regime of Theorem~\ref{thm:dense-collapse}.

The proof, given in Appendix~\ref{app:dense-comparison}, combines the
growing-cap law with a comparison on the whole low spectral subspace.
Unlike in the bounded-density case, the retained waists need not be
isolated leaves of the original cut tree. We therefore contract the
intervening branches and control their joint attachment errors,
including the loss incurred in passing from surface energy to piece
means. This yields a continuum conclusion beyond the geometric point
process, without invoking a fixed-topology degeneration estimate
uniformly in the number of cusps.

\subsection{Residual type at every dense growth rate}

Identifying the spectral type of the lowest eigenspaces does not
require the persistence condition \eqref{eq:dense-persistence-range}.
The following result concerns these actual eigenspaces throughout the
dense regime.

\begin{corollary}[Residual microscopic spectrum]\label{cor:dense-residual}
Suppose $g\to\infty$ and $n/g\to\infty$. For every fixed $T>0$,
\[
 \mathbb P\{X\text{ has a cuspidal eigenfunction with }
                 0<\lambda\le T/A_g\}\longrightarrow0.
\]
For each fixed $j\ge1$, the first $j$ positive eigenvalues, counted with
multiplicity, exist below $1/4$ and their full eigenspaces are residual
with probability tending to one.
If, in addition, \eqref{eq:dense-persistence-range} holds, the positive
residual eigenvalue process rescaled by $A_g$ converges in distribution
to $\operatorname{PPP}(2t\,dt)$ vaguely on $[0,\infty)$.
\end{corollary}

\begin{proof}
Let $\eta_b$ be the infimum of the Rayleigh quotient on the $H^1$ functions
whose angular average vanishes throughout every fixed-circumference
cusp tail; the precise form domain is given in
Appendix~\ref{app:cuspidal-exclusion}.
Proposition~\ref{prop:dense-cusp-zero} of that appendix proves
$S\eta_b\to\infty$ in probability without an upper restriction on
$n/g$; here $\eta_b$ is unrelated to the parameter $\eta$ of
\eqref{eq:dense-parameters}.
Every cusp-form eigenfunction belongs to this form domain, hence its
eigenvalue is at least $\eta_b$. Since $A_g\sim PS$, this proves the
first assertion. The cuspidal/residual decomposition below $1/4$
\cite[Section~3.3]{BMP} then makes every eigenspace in this window
entirely residual, including in the presence of multiplicities.

We separately check that $S\lambda_j^+=O_{\mathbb P}(1)$ for each fixed
$j$. If $n/g^3\to0$, fix a mark interval $[a_0,a_1]\subset(0,\infty)$
in Lemma~\ref{lem:dense-caps}. By choosing a large fixed score bound $R$,
with probability arbitrarily close to one there are at least $j+1$
interior-disjoint caps in this window. Their waist lengths tend to zero;
the collar lemma and the vanishing nested-pair bound ensure disjointness.
On each cap take a function equal to one off the inward width-one
waist half-collar, equal to inward distance on that half-collar, and
zero outside. For a cap with $k$ cusps and waist $\ell$, its energy is
at most $\ell\sinh1$ and its squared norm is at least
$2\pi(k-1)-\ell\sinh1$. Its Rayleigh quotient is therefore at most
$C\ell/(k-1)\le CR/S$. Disjoint supports give the same bound on the
entire $(j+1)$-dimensional span. Min--max, with the constant eigenvalue
counted once, yields $\lambda_j^+\le CR/S$.

If $n/g^3$ stays bounded below by a positive constant, then
Theorem~\ref{thm:dense-collapse} gives $S\lambda_j^+\to0$ in probability.
Every subsequence has a further subsequence of one of these two types,
which proves the asserted upper tightness for every dense sequence.
First choose a large fixed $T$ and then apply the first assertion to
$(0,T/A_g]$. Tightness shows that the first $j$ positive levels lie in
this window with probability arbitrarily close to one, so their full
eigenspaces are residual with probability tending to one. Under the
additional persistence condition, the residual point-process limit
now follows from Theorem~\ref{thm:dense-persistence}.
\end{proof}

This settles the spectral type of any fixed initial list of
eigenvalues at all dense growth rates, and complements the fixed-genus
residual-count question in \cite[Section~1.1, after Remark~1.2]{HTLargeN}.
It does not bound the total number of small cusp forms or settle the
Otal--Rosas conjecture discussed in \cite[Introduction]{BMP}. The
correct dense normalization and the joint law outside the persistence
range remain open. The fixed-window counterpart of this regime is
deterministic and has no such gap: Theorem~\ref{thm:cusp-dominated-spectrum}
identifies the limiting spectral measure for every sequence with
$g=o(n)$.

\section{Spectral measures in fixed-energy windows}
\label{sec:fixed-energy-spectrum}

The shrinking-scale process of Theorem~\ref{thm:supercritical-spectrum}
describes only the lowest eigenvalues. This section studies the
counting measure $\Xi_X$ of the spectrum in $(0,1/4)$ without
rescaling energy and proves deterministic limits in the three regimes
above the critical scale: $\nu_*$ when $\sqrt g\ll n\ll g$
(Theorem~\ref{thm:fixed-energy-spectrum}), $\nu_\lambda$ when
$n/g\to\lambda\in(0,\infty)$
(Theorem~\ref{thm:positive-density-spectrum}) and $\nu_\infty$ when
$g=o(n)$ (Theorem~\ref{thm:cusp-dominated-spectrum}), followed by the
regularized full-energy measure (Theorem~\ref{thm:regularized-spectrum})
and its determinant consequence
(Corollary~\ref{cor:cusp-density-determinant}). The common method is
the one announced in Subsection~\ref{subsec:strategy}: concentration
of orbital statistics, obtained from the volume theory through the
integration formula, is transferred to the spectrum by Huang's trace
approximation~\cite{HuangGap}. The trace and orbital inputs are
proved in Appendices~\ref{app:sublinear-trace},
\ref{app:positive-density-orbits} and~\ref{app:cusp-dominated-orbits},
and the spectral transfers in Appendix~\ref{app:fixed-window-proofs}.

When $\sqrt g\ll n\ll g$, a fixed-energy window below $1/4$ contains
order $n^2/g$ eigenvalues, and we first identify their deterministic
density, including matching expected counts.
Retain the funnel model $F_s$ and its ground-state energy $\beta(s)$
from Section~\ref{sec:critical-eigenvalues}, and set
\begin{equation}\label{eq:fixed-energy-measure}
 \mu_*(ds)=\frac{\sinh(s/2)}{8\pi^2}\,ds,
 \qquad \nu_*=\beta_*\mu_*.
\end{equation}
This measure is locally finite on $[0,1/4)$ by
Lemma~\ref{lem:funnel-long-waist}, and has no atom at zero.

\begin{theorem}[Fixed-energy spectral density]\label{thm:fixed-energy-spectrum}
Suppose $g\to\infty$ and $\sqrt g\ll n\ll g$, and put $d_g=n^2/g$.
For every $\varphi\in C_c([0,1/4))$,
\begin{equation}\label{eq:fixed-energy-L1}
 \mathbb E\left|
 \frac1{d_g}\sum_{0<\lambda_j(X)<1/4}\varphi(\lambda_j(X))
                  -\int\varphi\,d\nu_*\right|\longrightarrow0.
\end{equation}
In particular the normalized counting measures converge vaguely in
probability to $\nu_*$, and their mean measures converge vaguely.
Eigenvalues are counted with multiplicity, including residual eigenvalues.
\end{theorem}

Gatignol--Thomas~\cite[Theorem~1.10]{GT} prove a lower count of
order $n^2/g$ below each fixed positive threshold in this regime;
Theorem~\ref{thm:fixed-energy-spectrum} supplies the matching density.
The upper bound uses Huang's canonical trace coefficients and positive
truncation method~\cite{HuangGap}, together with a new
complementary-volume tail, proved in Appendix~\ref{app:sublinear-trace},
which removes the quadratic cusp restriction on the trace remainder.
Positive limiting $n/g$ is treated separately in
Theorem~\ref{thm:positive-density-spectrum}; neither result concerns
moving windows at $1/4$.

\smallskip\noindent\emph{Structure of the proof.} The proof, given in
Appendix~\ref{app:fixed-window-proofs}, has two halves. The upper half
identifies the mean measure: after the cusp labels are summed at the
saddle, the trace formula leaves only the term of the two-cusp pairs of pants, whose
distribution is $\nu_*$ by Proposition~\ref{prop:critical-identification}.
The lower half places disjoint funnel trial states on short two-cusp
pants. A law of large numbers for their waists, together with the
removal of overlapping configurations, shows that these states account
for the whole mean. Matching the two halves gives convergence in $L^1$.
In the sublinear regime only two-cusp pairs of pants contribute; the
next two subsections treat the densities at which all planar filling
types enter.

\subsection{Positive cusp density}
\label{sec:positive-density-spectrum}

When $n/g\to\lambda\in(0,\infty)$, planar subsurfaces with more
than two cusps contribute at the same order as two-cusp pairs of
pants, so the limiting low spectrum involves all planar filling types.
Theorem~\ref{thm:positive-density-spectrum} identifies the limit of
$g^{-1}\Xi_X$ as a measure $\nu_\lambda$ characterized through the
Weil--Petersson integrals of these types and the boundary profile of
Proposition~\ref{prop:bounded-boundaries}; this is the regime in which
the density parameter $Y_\lambda$ at an interior point, rather than
at an endpoint of $[0,y_0]$, enters the spectrum. The orbital concentration
is proved in Appendix~\ref{app:positive-density-orbits} and the
spectral transfer in Appendix~\ref{app:fixed-window-proofs}.

\begin{theorem}[Low spectrum at positive density]
\label{thm:positive-density-spectrum}
There is a vaguely continuous family of Radon measures
$\{\nu_\lambda:0<\lambda<\infty\}$ on $[0,1/4)$ such that, whenever
$g\to\infty$ and $n/g\to\lambda\in(0,\infty)$,
\begin{equation}\label{eq:positive-density-L1}
 \mathbb E\left|
 \frac1g\sum_{0<\lambda_j(X)<1/4}\varphi(\lambda_j(X))
                   -\int\varphi\,d\nu_\lambda\right|\longrightarrow0,
 \qquad \varphi\in C_c([0,1/4)).
\end{equation}
Eigenvalues are counted with multiplicity, including residual modes.
The measure $\nu_\lambda$ has no atom at zero.
The planar orbital integrals \eqref{eq:positive-orbital-functional}
and the limit \eqref{eq:positive-spectral-characterization} determine
$\nu_\lambda$ uniquely.
\end{theorem}

Thus the normalized low spectrum converges vaguely in probability,
and its mean converges as well. No absolute continuity of $\nu_\lambda$
is assumed. Hide--Thomas~\cite[Remark~1.8]{HTLargeN} ask about spectral
measures in the large-cusp regime; the theorem treats their low-energy
part when genus and cusp count grow proportionally. The endpoint
$g=o(n)$, including fixed genus, is treated in
Theorem~\ref{thm:cusp-dominated-spectrum}. Neither theorem treats
the continuous spectrum above $1/4$.

The proof differs from the sublinear case in its order of limits.
After division by $g$, one may let $g$ tend to infinity with the
polynomial degree fixed, and only then let the degree tend to
infinity. No trace expansion uniform in a growing degree is needed.

\smallskip\noindent\emph{Structure of the proof.}
Appendix~\ref{app:positive-density-orbits} constructs a functional
$\Omega_\lambda$ from planar orbital integrals and proves
\begin{equation}\label{eq:positive-orbital-L1}
 \mathbb E|g^{-1}\mathcal O_X(H)-\Omega_\lambda(H)|\longrightarrow0
 \quad(H\in C_c^\infty(\mathbb R)\text{ even}).
\end{equation}
The functional is continuous in $\lambda$ on compact density intervals.
Huang's positive trace~\cite[Sections~5--6]{HuangGap} then transfers
this fixed-support assertion to the actual spectrum, as carried out in
Appendix~\ref{app:fixed-window-proofs}. The functional
$\Omega_\lambda$ and the measures $\nu_\lambda$ are carried forward:
their $\lambda\to\infty$ limits are identified in the next subsection
\eqref{eq:cusp-dominated-continuity}, and they reappear as the
low-energy part of the regularized measure in
Theorem~\ref{thm:regularized-spectrum}.

\subsection{The cusp-dominated limit}
\label{sec:cusp-dominated-spectrum}

The endpoint $n/g\to\infty$ has a universal fixed-energy limit
$\nu_\infty$, which agrees with the large-cusp limit at every fixed
genus, including genus zero
(Theorem~\ref{thm:cusp-dominated-spectrum}). The input is the
critical-weight endpoint of the volume theory, $x=x_c$ and $Y=y_0$,
at which the planar orbital functional of the previous subsection has
a limit; the orbital concentration is proved in
Appendix~\ref{app:cusp-dominated-orbits} and the spectral transfer in
Appendix~\ref{app:fixed-window-proofs}. The statement concerns
energies held fixed below $1/4$; the shrinking scale of the first few
eigenvalues in the same regime is the unresolved question of
Section~\ref{sec:dense-spectrum}.

\begin{theorem}[Universal cusp-dominated low spectrum]
\label{thm:cusp-dominated-spectrum}
There is a Radon measure $\nu_\infty$ on $[0,1/4)$ such that, whenever
$n\to\infty$ and $g=g(n)\ge0$ is any integer sequence with $g=o(n)$,
\begin{equation}\label{eq:cusp-dominated-L1}
 \mathbb E\left|\frac1n\sum_{0<\lambda_j(X)<1/4}\varphi(\lambda_j(X))
                  -\int\varphi\,d\nu_\infty\right|\longrightarrow0,
 \qquad \varphi\in C_c([0,1/4)).
\end{equation}
This includes multiplicities and residual eigenvalues. For the measures
of Theorem~\ref{thm:positive-density-spectrum},
\begin{equation}\label{eq:cusp-dominated-continuity}
 \lambda^{-1}\nu_\lambda\longrightarrow\nu_\infty
       \quad\text{vaguely as }\lambda\to\infty.
\end{equation}
The measure is characterized by
\eqref{eq:positive-spectral-characterization} with $\Omega_\lambda$
replaced by $\Omega_\infty$ defined below.
It has no atom at zero.
\end{theorem}

Here $\Omega_\infty$ is obtained by dividing the right side of
\eqref{eq:positive-orbital-functional} by $\lambda$ and setting
\[
 x=x_c,\qquad K_x=K_{x_c},\qquad
 \mathfrak b_Y(L)=\mathfrak b_{y_0}(L)
                  =I_0\!\left(\frac{j_{0,1}L}{2\pi}\right).
\]
The series $K_{x_c}$ converges, since its coefficients at critical
weight decay as $j^{-5/2}$ by the genus-zero case of
Manin--Zograf~\cite[Theorem~6.1]{ManZ}; see also
\cite[Theorem~2.4]{HTLargeN}. The boundary bound and fixed-support
fiber estimates therefore give
\begin{equation}\label{eq:cusp-dominated-functional}
 \lambda^{-1}\Omega_\lambda(H)\longrightarrow\Omega_\infty(H)
 \quad(H\in C_c^\infty(\mathbb R)\text{ even}).
\end{equation}

Hide--Thomas~\cite[Theorem~1.1]{HTLargeN} prove a linear lower count
of small eigenvalues at fixed genus;
Theorem~\ref{thm:cusp-dominated-spectrum} identifies the full limiting
low spectral measure in that regime and, more generally, throughout
$g=o(n)$. This treats the exceptional discrete part of the
large-cusp spectral-measure question in~\cite[Remark~1.8]{HTLargeN};
the full regularized measure, including scattering, is treated in
Subsection~\ref{sec:regularized-spectrum}.

\smallskip\noindent\emph{Structure of the proof.}
Appendix~\ref{app:cusp-dominated-orbits} proves
\begin{equation}\label{eq:cusp-dominated-orbital-L1}
 \mathbb E|n^{-1}\mathcal O_X(H)-\Omega_\infty(H)|\longrightarrow0
 \qquad(g=o(n),\ n\to\infty).
\end{equation}
It treats diverging genus, fixed positive genus and genus zero
separately and obtains the same functional in each case. The
distinction is necessary:
the positive-convolution local limit requires diverging genus, whereas
on a sphere the leading complementary component is selected by a
summable coefficient convolution. The transfer to the spectrum is that
of Theorem~\ref{thm:positive-density-spectrum}, with normalization $n$;
see Appendix~\ref{app:fixed-window-proofs}. The functional
$\Omega_\infty$ is the $\theta=0$ endpoint of the family used in the
regularized law of the next subsection.

\subsection{The full regularized spectral measure}
\label{sec:regularized-spectrum}

The three preceding theorems concern the discrete spectrum below
$1/4$. To include the continuous spectrum, which starts at $1/4$, we
use the classical geometric regularization of the heat trace
\cite{JLHeat,HJLCounting}, and prove that the resulting signed
spectral measure, divided by $n$, has a deterministic full-energy
limit whenever $g/n\to\theta\in[0,\infty)$
(Theorem~\ref{thm:regularized-spectrum}). For the trace formula,
Eisenstein series and the scattering determinant on cofinite surfaces
see~\cite{Hejhal,Iwaniec}. The regularized measure combines all
discrete eigenvalues with the logarithmic derivative of the scattering
determinant, and the parabolic corrections must be retained when the
number of cusps is comparable to the area. The inputs are the orbital
limits \eqref{eq:positive-orbital-L1} and
\eqref{eq:cusp-dominated-orbital-L1} of the two preceding subsections;
the proof is in Appendix~\ref{app:fixed-window-proofs}.

In width-one cusp coordinates, let $\Phi_X(s)$ be the
scattering matrix and $\phi_X(s)=\det\Phi_X(s)$. For
$F\in C_c([0,\infty))$, define
\begin{align}\label{eq:regularized-measure}
 \mathfrak m_X(F)
 &=\sum_j F(\lambda_j(X))
   -\frac1{4\pi}\int_{\mathbb R}F(1/4+r^2)
                 \frac{\phi_X'}{\phi_X}(1/2+ir)\,dr\notag\\
 &\quad+\frac n{2\pi}\int_{\mathbb R}F(1/4+r^2)
             \bigl(\operatorname{Re}\psi(1+ir)+\log2\bigr)\,dr
       -\frac{n-\operatorname{tr}\Phi_X(1/2)}4F(1/4).
\end{align}
Here $\psi=\Gamma'/\Gamma$, the prime on $\phi_X$ is the complex
$s$ derivative, and the sum includes every discrete $L^2$ eigenvalue,
with multiplicity, including zero and embedded eigenvalues.
The scattering matrix is regular and unitary on the critical line;
its value at $1/2$ is an involution. Formula
\eqref{eq:regularized-measure} defines a signed locally finite measure,
not a positive eigenvalue-counting measure, and fixes the
regularization used throughout this subsection.

For $\operatorname{Re}z>0$, set
\begin{align}\label{eq:regularized-heat-functions}
 G_z(u)&=\frac{e^{-z/4-u^2/(4z)}}{\sqrt{4\pi z}},&
 k_{\mathbb H}(z)&=\frac1{4\pi}\int_{\mathbb R}
                   e^{-z(1/4+r^2)}r\tanh(\pi r)\,dr,\notag\\
 \mathcal G_\theta(H)&=\theta\Omega_{1/\theta}(H)\quad(\theta>0),&
 \mathcal G_0(H)&=\Omega_\infty(H).
\end{align}
The square root is positive on the positive real axis.
Initially $H$ is smooth, even and compactly supported; the proof in
Appendix~\ref{app:fixed-window-proofs} extends $\mathcal G_\theta$ to
$G_z$ by length cutoffs.

\begin{theorem}[Full regularized spectral law]
\label{thm:regularized-spectrum}
There is a vaguely continuous family of signed Radon measures
$\mathfrak m_\theta$, $0\le\theta<\infty$, on $[0,\infty)$ such that
whenever $n\to\infty$ and $g/n\to\theta$,
\begin{equation}\label{eq:regularized-L1}
 \mathbb E\left|\frac{\mathfrak m_X(F)}n-\mathfrak m_\theta(F)\right|
       \longrightarrow0,
       \qquad F\in C_c([0,\infty)).
\end{equation}
The limit is characterized by its heat transform
\begin{equation}\label{eq:regularized-limit-heat}
 \mathcal R_\theta(z):=\int_0^\infty e^{-zE}\,d\mathfrak m_\theta(E)
 =2\pi(1+2\theta)k_{\mathbb H}(z)+\mathcal G_\theta(G_z),
       \qquad \operatorname{Re}z>0.
\end{equation}
Writing $\mathcal R_X(z)=\mathfrak m_X(e^{-zE})$, one also has
\begin{equation}\label{eq:regularized-analytic-L1}
 \mathbb E\sup_{z\in K}
   |n^{-1}\mathcal R_X^{(q)}(z)-\mathcal R_\theta^{(q)}(z)|\longrightarrow0
\end{equation}
for every compact $K\subset\{\operatorname{Re}z>0\}$ and fixed
$q\ge0$. On $[0,1/4)$ the limit equals
$\theta\nu_{1/\theta}$ if $\theta>0$ and $\nu_\infty$ if $\theta=0$.
In particular it has no atom at zero.
\end{theorem}

Thus, for the standard regularization \eqref{eq:regularized-measure},
the large-cusp spectral-measure question in \cite[Remark~1.8]{HTLargeN}
has a full-energy limit, including the threshold.
For uniformly discrete sequences Benjamini--Schramm converging to
$\mathbb H$, the normalized spectral measures converge to the Plancherel
measure~\cite{ABBGNRS}; for hyperbolic surfaces see Le Masson and
Sahlsten~\cite{LMSQE} and Monk~\cite{Monk}. Here the local limit is
not $\mathbb H$.
Le Masson--Sahlsten~\cite[Theorem~5.1]{LMSScattering} obtain spectral
convergence when the cusp count is sublinear in the genus, under
Benjamini--Schramm convergence to $\mathbb H$ and a uniform systole
bound; the present regime has a different orbital limit and imposes
no systole bound.

The theorem does not separate cusp forms from scattering, identify the
atom at $1/4$, or give a moving-window law. Nor does signed weak
convergence alone imply sharp interval counts, even at continuity
points of the signed limit.

\smallskip\noindent\emph{Structure of the proof.} The starting point is the
exact regularized trace formula \cite[Proposition~5.1]{HJLCounting},
\begin{equation}\label{eq:regularized-exact-heat}
 \mathcal R_X(z)=\operatorname{area}(X)k_{\mathbb H}(z)
                         +\mathcal O_X(G_z).
\end{equation}
This uses the unshifted Laplacian and the oriented convention for
$\mathcal O_X$ of Definition~\ref{rev:def-huang-orbital} (a sum over
oriented primitive classes and their repetitions); every term in the parabolic
part retains the factor $e^{-z/4}$. The proof then has three steps.
A packing bound in the thick part controls long repetitions without
any systole bound. After these are cut off, the orbital term converges
by \eqref{eq:positive-orbital-L1} and
\eqref{eq:cusp-dominated-orbital-L1}. Finally, continuous energy tests
are recovered from heat transforms through a bound on the signed total
variation, which uses the residual-pole decomposition of the
scattering determinant. The details are in
Appendix~\ref{app:fixed-window-proofs}. The heat transform
\eqref{eq:regularized-limit-heat} is the input of the determinant
limit in the next subsection.

\subsection{A determinant consequence}
\label{sec:spectral-determinant}

The full spectral law also controls a singular global statistic, the
log determinant of the Laplacian: Corollary~\ref{cor:cusp-density-determinant}
gives the $L^1$ limit of the normalized log determinant, and of
$\log Z_X'(1)$, whenever $g/n\to\theta\in[0,\infty)$. Its inputs are
the heat transform \eqref{eq:regularized-limit-heat} of
Theorem~\ref{thm:regularized-spectrum} on compact time intervals and,
at both ends of the time axis, the short-geodesic intensity bound and
low-energy tail of Proposition~\ref{prop:low-energy-tail}
(Appendix~\ref{app:cusp-dominated-orbits}); the proof is in
Appendix~\ref{app:fixed-window-proofs}. On compact hyperbolic surfaces the
zeta-regularized determinant is expressed through the Selberg zeta
function~\cite{DHokerPhong,SarnakDet}, and Efrat treated surfaces of
finite area~\cite{Efrat}. Naud~\cite{NaudDet} proved
convergence of the normalized log determinant for several closed
random-surface models, and He--Wu~\cite{HeWuDet} obtained $L^1$
convergence in the closed Weil--Petersson model.
Jorgenson, Smajlovi\'c and Spilioti~\cite[Theorem~5.3]{JSSDet} obtain
deterministic large-volume limits for weighted Laplacian determinants
on compact surfaces, under a weak spectral gap, uniform discreteness
and a non-accumulation hypothesis for bounded-length geodesics.
Gross--Lachman--Nachmias~\cite[Theorem~1]{GLNDet} treat arbitrary local
limits of compact hyperbolic surfaces under uniform integrability of
short inverse lengths, without a uniform spectral gap. In the present
noncompact setting, the growing-cusp law has a nonzero orbital
contribution, and the scattering and parabolic terms in
\eqref{eq:regularized-measure} are retained.

Write $A_X=\operatorname{area}(X)$ and $a_\theta=2\pi(1+2\theta)$, and put
\begin{equation}\label{eq:determinant-orbital-limit}
 \mathcal B_X(t)=\mathcal O_X(G_t),\qquad
 \mathcal B_\theta(t)=\mathcal G_\theta(G_t),\qquad
 I_\theta=\int_0^\infty\mathcal B_\theta(t)\,\frac{dt}{t}.
\end{equation}
Use the Selberg product
\[
 Z_X(s)=\prod_{\gamma\text{ primitive oriented}}\prod_{j\ge0}
                 (1-e^{-(s+j)\ell_\gamma}),\qquad \operatorname{Re}s>1,
\]
with the same oriented convention as \eqref{eq:regularized-exact-heat}.
Its continuation has a simple zero at $s=1$; see~\cite{Hejhal,Efrat,BJP}.\footnote{The
product runs over oriented primitive classes, equivalently over the
primitive hyperbolic conjugacy classes of the covering group, since an
element and its inverse are not conjugate in a torsion-free Fuchsian
group; each unoriented primitive closed geodesic therefore contributes
two factors. This is the convention forced by the hyperbolic terms of
the trace formula in \eqref{eq:regularized-exact-heat}, which count
both orientations, and it is the one in which the Mellin identity
$\int_0^\infty e^{-s(s-1)t}\mathcal B_X(t)\,dt/t=-\log Z_X(s)$ used in
Appendix~\ref{app:fixed-window-proofs} holds without an extra factor.
If a reference defines the Selberg zeta function as a product over
unoriented primitive geodesics, that function is the square root of
ours; all statements about $Z_X$ in this paper refer to the oriented
product.}
To specify the determinant on the noncompact surface, define
\begin{equation}\label{eq:geometric-determinant-definition}
 \zeta_X^{\mathrm{geom}}(w)
 =\frac1{\Gamma(w)}\int_0^\infty
       t^{w-1}(\mathcal R_X(t)-1)\,dt,\qquad
 \log\det{}'_{\mathrm{geom}}\Delta_X
       =-(\zeta_X^{\mathrm{geom}})'(0).
\end{equation}
The integral initially converges for $\operatorname{Re}w>1$;
continuation to zero is justified in the proof. The hyperbolic-plane density is
\[
 E_{\mathbb H}=
 \frac{4\zeta_{\mathrm R}'(-1)-1/2+\log(2\pi)}{4\pi},
\]
where $\zeta_{\mathrm R}$ is the Riemann zeta function
\cite[Section~2]{NaudDet}.

\begin{corollary}[Log determinants at positive cusp density]
\label{cor:cusp-density-determinant}
The function $I_\theta$ in \eqref{eq:determinant-orbital-limit} is finite,
strictly positive and continuous for $0\le\theta<\infty$.
If $n\to\infty$ and $g/n\to\theta$ under normalized Weil--Petersson law,
then
\begin{align}\label{eq:cusp-density-determinant-law}
 \mathbb E\left|\frac{\log Z_X'(1)}n+I_\theta\right|&\longrightarrow0,
       \notag\\
 \mathbb E\left|\frac{\log\det{}'_{\mathrm{geom}}\Delta_X}{A_X}
            -\left(E_{\mathbb H}-\frac{I_\theta}{a_\theta}\right)\right|
       &\longrightarrow0.
\end{align}
In particular the latter limit is strictly below $E_{\mathbb H}$.
\end{corollary}

The two limits are equivalent: for every surface $X$ one has the exact
identity
$\log\det{}'_{\mathrm{geom}}\Delta_X=A_XE_{\mathbb H}+\log Z_X'(1)$
(Appendix~\ref{app:fixed-window-proofs}), and $A_X/n\to a_\theta$.
The corollary includes $g=o(n)$ and every fixed genus. The determinant
is defined precisely by \eqref{eq:geometric-determinant-definition};
alternative noncompact regularizations may introduce cusp-dependent
constants, whereas the Selberg-derivative statement fixes the
normalization independently of this terminology. Neither statement
asserts convergence of exponential moments.

\smallskip\noindent\emph{Structure of the proof.} The classical
heat-to-determinant method~\cite{NaudDet,HeWuDet,GLNDet} expresses
$\log Z_X'(1)$ through the integral of the orbital heat term
$\mathcal B_X(t)$ against $dt/t$, with the zero mode subtracted at
large times. Theorem~\ref{thm:regularized-spectrum} gives convergence on
compact time intervals, so the work lies at the two ends, which require
averaged bounds uniform in the surface. At small times we use an
intensity bound for short simple geodesics of every separating
topology; at large times, the low-energy tail of
Proposition~\ref{prop:low-energy-tail} replaces a uniform spectral gap.
Short two-cusp pairs of pants make $I_\theta$ strictly positive. The
details are in Appendix~\ref{app:fixed-window-proofs}. This corollary
closes the spectral part; the regularized measure and the determinant
are not used in Part~\ref{part:geometry}.

\clearpage
\part{Local geometry and scaling limits}\label{part:geometry}
\phantomsection\label{guide:geometric-context}
At positive cusp density $n/g\to\lambda\in(0,\infty)$, the
volume-ratio profiles of Theorem~\ref{thm:volume} enter the carrier
integrals of Appendix~\ref{app:positive-density-orbits}. Retaining
based holonomies in these integrals determines whole pointed balls and
constructs the local family. This part uses no spectral theorem of
Part~\ref{part:spectra}; it also uses the genus-zero convolution
identities of Appendix~\ref{app:cusp-dominated-orbits} and the cusp
cutoffs of Appendix~\ref{app:cusp-cutoffs}. The family's planar endpoint
has a Brownian-plane scaling limit, proved through a separate
growing-radius tree coupling.

In genus zero, Budd and Curien showed that the Weil--Petersson random sphere with $n+1$ cusps, with its metric rescaled by $n^{-1/4}$, converges to the Brownian sphere~\cite{BC}, the universal scaling limit of
uniform random planar maps~\cite{LeGallUniq,Miermont}.
In the discrete setting, uniform triangulations of high genus converge
locally to the planar stochastic hyperbolic triangulations of
Curien~\cite{CurienPSHT,BudzinskiLouf}, whose near-critical scaling limit
is the hyperbolic Brownian plane of Budzinski~\cite{BudzinskiHBP}. Theorem~\ref{thm:wp-local-family} is a
continuum Weil--Petersson analogue of the high-genus local-limit
statement. For closed Weil--Petersson surfaces, Monk proved
quantitative Benjamini--Schramm convergence to the hyperbolic
plane~\cite{Monk}; the limit $\lambda\downarrow0$ in
\eqref{eq:wp-local-endpoints} recovers the plane from positive density.
Theorems~\ref{thm:wp-local-family} and~\ref{thm:wp-conformal-phase}
construct the local family at every finite positive cusp density and
show that the cusp-filled surfaces are conformally discs, with the
recurrent infinite WP sphere as the infinite-density endpoint.
Budd and Curien asked for a hyperbolic analogue of their infinite
sphere, as in~\cite{CurienPSHT}, and conjectured that the infinite
sphere is parabolic, with recurrent Brownian motion and ergodic
geodesic flow~\cite[Section~7.3]{BC}. The family $S_\lambda$ provides
such analogues, and Theorem~\ref{thm:wp-conformal-phase} with
Corollary~\ref{cor:wp-geodesic-flow} proves the conjecture. The
conformal-type conclusion concerns cusp-filled surfaces and does not
construct a canonical embedding or describe the finer conformal
structure; the near-critical scaling question remains open. A related peeling approach by Budd, Budzinski, Curien and Petri
was announced in \cite[reference~19]{BC}; our local construction
instead retains based holonomies in the carrier integration formula.
Theorem~\ref{thm:wp-brownian-plane} answers their question about the
Brownian-plane scaling limit of the infinite area-rooted WP surface
\cite[Section~7.3]{BC}, with cusp-trimmed measured GHP convergence,
full-area GP convergence and a finite mesoscopic version. Its
growing-radius coupling, proved in Appendix~\ref{app:brownian-plane},
describes the planar geometry but does not supply the energy
comparison required for the dense microscopic spectrum.

\section{Local geometry and conformal type at positive cusp density}
\label{sec:local-conformal-phase}

We construct the local WP family at
$n/g\to\lambda\in(0,\infty)$ and determine the conformal type after
filling its cusps. Theorem~\ref{thm:wp-local-family} is proved in
Appendix~\ref{app:rooted-local-law}; the finite-density and
recurrent-endpoint parts of Theorem~\ref{thm:wp-conformal-phase} are
proved in Appendices~\ref{app:conformal-phase} and~\ref{app:cusp-cutoffs},
respectively.

Throughout this section, $B_r(S,o)$ carries the restricted ambient
metric and the hyperbolic area, and $r$ is fixed before the large-genus
limit. Total variation distances are taken between laws of pointed
measured balls up to isometry. Write $A_{g,n}=2\pi(2g-2+n)$.

\begin{theorem}[The local WP family]\label{thm:wp-local-family}
For every $\lambda\in(0,\infty)$ there is a unique law $\mathcal L_\lambda$
of complete pointed hyperbolic surfaces $(S_\lambda,o_\lambda)$ with the
following properties. If $g\to\infty$, $n/g\to\lambda$, $X$ has
normalized WP law on $\mathcal M_{g,n}$, and $o$ is conditionally uniform
for its area, then, for every fixed $r<\infty$,
\begin{equation}\label{eq:wp-local-tv}
 d_{\rm TV}\bigl(\operatorname{Law}(B_r(X,o)),
               \operatorname{Law}(B_r(S_\lambda,o_\lambda))\bigr)
       \longrightarrow0.
\end{equation}
Moreover, for every bounded Borel function $F$ of a pointed measured
$r$-ball,
\begin{equation}\label{eq:wp-local-quenched}
 \frac1{A_{g,n}}\int_X F(B_r(X,p))\,dA(p)
       \longrightarrow\mathbb E_{\mathcal L_\lambda}F(B_r)
       \quad\text{in }L^2(\mathrm{WP}).
\end{equation}
The law $\mathcal L_\lambda$ is unimodular. It depends continuously on
$\lambda$ in fixed-radius total variation, with endpoints
\begin{equation}\label{eq:wp-local-endpoints}
 \mathcal L_\lambda\longrightarrow(\mathbb H^2,o)
       \quad(\lambda\downarrow0),\qquad
 \mathcal L_\lambda\longrightarrow(S_\infty^*,o_\infty)
       \quad(\lambda\uparrow\infty)
\end{equation}
in the same sense, where the second endpoint is area-rooted as
in~\cite[Theorem~34]{BC}.
\end{theorem}

Here unimodularity means that, for every nonnegative isometry-invariant
Borel function $m$ of a surface with two marked points,
\begin{equation}\label{eq:local-mass-transport}
 \mathbb E\int_S m(S,o,z)\,dA(z)
       =\mathbb E\int_S m(S,z,o)\,dA(z).
\end{equation}
This is the manifold mass-transport identity of~\cite{AbertBiringer};
for unimodular random networks see~\cite{AldousLyons}.
The theorem requires no ergodicity assertion. The family is
nonconstant: for every fixed circumference $0<b\le2$ of
the standard cusp horocycles, its cusp occupation probability is
\begin{equation}\label{eq:local-cusp-occupation}
 \mathbb P_{\mathcal L_\lambda}
       \{o_\lambda\text{ lies in a }b\text{-cusp tail of }S_\lambda\}
       =\frac{b\lambda}{2\pi(2+\lambda)}.
\end{equation}

\begin{theorem}[Conformal transition]\label{thm:wp-conformal-phase}
For $0<\lambda<\infty$, almost surely $S_\lambda$ has genus zero and
infinite area, and every simple closed geodesic bounds a unique finite
planar cap. After filling its isolated cusps, the resulting Riemann
surface is conformally equivalent to $\mathbb D$. Brownian motion for
the original complete hyperbolic metric is transient.

At the endpoint $\lambda=\infty$, the Riemann surface obtained from
$S_\infty^*$ by filling its cusps is conformally equivalent to
$\mathbb C$, and Brownian motion for the hyperbolic metric of
$S_\infty^*$ is recurrent.
\end{theorem}

The endpoint assertion is the parabolicity conjectured by Budd and
Curien~\cite[Section~7.3]{BC}, who also conjectured that the geodesic
flow is ergodic. The latter follows from the classical equivalence
between recurrence of Brownian motion and ergodicity of the geodesic
flow.

\begin{corollary}[Geodesic flow]\label{cor:wp-geodesic-flow}
For $0<\lambda<\infty$, almost surely the geodesic flow on the unit
tangent bundle of $S_\lambda$ is not ergodic for the Liouville measure.
At the endpoint, almost surely the geodesic flow on the unit tangent
bundle of $S_\infty^*$ is ergodic. In particular, the area-rooted
infinite Weil--Petersson sphere is parabolic, has recurrent Brownian
motion and has ergodic geodesic flow, as conjectured
in~\cite[Section~7.3]{BC}.
\end{corollary}

\begin{proof}
Each surface is complete and hyperbolic, hence of the form
$\mathbb H^2/\Gamma$ for a torsion-free Fuchsian group $\Gamma$. For such
a surface, Brownian motion is recurrent if and only if the Poincar\'e
series $\sum_{\gamma\in\Gamma}e^{-d(z,\gamma z)}$ diverges, because the
Green function of the quotient is the $\Gamma$-sum of that of
$\mathbb H^2$, which decays like $e^{-d}$; and divergence of this series
is equivalent to ergodicity of the geodesic flow for the Liouville
measure. This is the theorem of Hopf, Tsuji and
Sullivan~\cite{SullivanErgodic}; see~\cite[Section~1.1]{BHS} for this
list of equivalent conditions. Theorem~\ref{thm:wp-conformal-phase}
gives transience for $0<\lambda<\infty$ and recurrence at the endpoint.
\end{proof}

Appendix~\ref{app:rooted-local-law} proves
Theorem~\ref{thm:wp-local-family} by retaining the primitive holonomies
of short based loops. Their number in a fixed displacement window has
a deterministic bound, even in a thin cusp, and a support lemma for
generators and their pairwise products permits exact carrier unfolding.
The one- and two-carrier coefficient estimates of
Appendix~\ref{app:positive-density-orbits}, applied to uniformly bounded
rooted fibers, then determine all local statistics and their variances.

The finite-density assertion of Theorem~\ref{thm:wp-conformal-phase}
rests on the subcritical cap-size tail, which implies that only
finitely many caps with small boundary-to-area ratio can contain a
fixed non-elementary core.
This gives an eventual isoperimetric inequality for connected domains
containing a fixed point, and the Green-level argument in
Appendix~\ref{app:conformal-phase} proves transience directly for the
original metric. The recurrent endpoint is the conformal conclusion
proved in Appendix~\ref{app:cusp-cutoffs}. The dichotomy parallels the
hyperbolic--parabolic dichotomy for unimodular random planar
maps~\cite{AHNR}, the recurrence of local limits of finite planar
graphs~\cite{BenjaminiSchramm,GGN} and the parabolicity of the Riemann
surface of the uniform infinite planar triangulation~\cite{GillRohde};
we use these only as analogies. For
recurrence and transience of Brownian motion on Riemannian manifolds
see~\cite{Grigoryan}.

The limits of these statements are as follows. The convergence
\eqref{eq:wp-local-tv} keeps the radius fixed, and its estimates may
depend on $\lambda$; the results give neither a growing-radius
approximation uniform as $\lambda\to\infty$ nor the inverse-energy
estimate needed for the dense microscopic spectrum. The near-critical
scaling limit of a hyperbolic family requested
in~\cite[Section~7.3]{BC} therefore remains a separate question
(Section~\ref{sec:conclusions}). The next section treats a metric
scaling limit at the planar endpoint.

\section{The planar endpoint and the Brownian plane}
\label{sec:brownian-plane}

The area-rooted infinite WP sphere, the $\lambda=\infty$ endpoint of
\eqref{eq:wp-local-endpoints}, has the Brownian plane as its scaling
limit. The same limit holds for finite WP spheres at every mesoscopic
radius $1\ll r_n\ll n^{1/4}$. We prove measured
Gromov--Hausdorff--Prokhorov convergence after trimming the deep cusps
and projecting their full area onto the trim boundaries, and
Gromov--Prokhorov convergence for the untrimmed surface with its full
area. The proof in Appendix~\ref{app:brownian-plane} couples growing
balls through the continuous-type trees of Budd and Curien~\cite{BC},
then transfers their finite Brownian-sphere limit~\cite[Theorem~2]{BC}
by the intermediate-size strategy of Curien and Le Gall~\cite{CLGPlane}.

Write $S_n$ for a WP sphere with $n+1$ cusps, $A_n=2\pi(n-1)$ for its
area, and $o_n$ for a conditionally uniform area point. Let
$(S_\infty^*,o_\infty)$ be the area-rooted infinite WP sphere
of~\cite[Theorem~34]{BC}. For either surface, remove the open cusp
neighborhoods bounded by horocycles of length one and write $S^\circ$
for the remaining closed subset. Set
\[
 X(S,o)=S^\circ\cup\{o\}.
\]
This space carries the original ambient hyperbolic distance, not its
intrinsic path distance. The root is retained even when it lies in a
removed cusp. Let $P:S\to S^\circ$ be vertical projection in each cusp
and the identity elsewhere, and put $\bar A=P_*dA$. In particular,
\begin{equation}\label{eq:bp-projected-area}
 \bar A=dA|_{S^\circ}+\sum_C d\ell_{\partial C}.
\end{equation}
Thus each cusp contributes its entire unit area as boundary arclength.

Let $(\mathcal P,D_\infty,\mu_\infty,\rho_\infty)$ be the Brownian
plane, with the natural volume obtained by pushing Lebesgue measure
through the coding of Curien and Le Gall~\cite{CLGPlane}. With the
Brownian-sphere distance convention of~\cite{BC,CLGPlane}, put
\begin{equation}\label{eq:bp-normalization}
 c_{\rm WP}=\frac{2\pi}{\sqrt{3j_{0,1}}},
 \qquad v_{\rm WP}=\frac{2\pi}{c_{\rm WP}^4}.
\end{equation}

We use local pointed Gromov--Hausdorff--Prokhorov convergence on proper
spaces, allowing the marked root to lie outside the measure support.
Explicitly, the local distance is the integral, against $e^{-r}dr$, of
the compact pointed GHP distance between closed $r$-balls, capped at
one. Balls carry restricted ambient distances and finite measures,
without renormalization. This is the formula by which Abraham, Delmas
and Hoscheit define the Gromov--Hausdorff--Prokhorov distance on
locally compact length spaces~\cite{ADH}; we use it on proper spaces.
The pointed GP convention omits the Hausdorff
term and identifies spaces on the support of the measure together with
the marked root, as in the Gromov--Prokhorov metric of Greven,
Pfaffelhuber and Winter~\cite{GPW} and its pointed locally finite
form, the Gromov-vague topology of~\cite{ALW}. These definitions also apply to the possibly nonlength
spaces $X(S,o)$, whose root may lie outside the support of the measure.
Indeed, for a proper pointed space with a locally finite measure every
closed ball is compact and carries a finite restricted measure, so the
compact pointed GHP distance between $r$-balls is finite, with the
roots compared as marked points whether or not they lie in the
supports; as $r'\downarrow r$ the closed $r'$-balls decrease to the
closed $r$-ball in the Hausdorff sense and the restricted measures
converge weakly by monotone convergence, so the integrand is
right-continuous in $r$, hence Borel, and the integral is defined. No
length-space property is used, and if two spaces have isometric
measured $r$-balls for almost every $r$, then for every $r$ (restrict
an isometry of a larger ball) and hence, by a diagonal argument on
proper spaces, the spaces are isometric. Convergence in this local
distance follows from convergence of the compact measured $r$-balls
for almost every $r$, by dominated convergence of the capped integrand
against $e^{-r}dr$; this is the form in which the convergences below
are established in Appendix~\ref{app:brownian-plane}, through an almost
sure representation for the finite spheres and through bounded
Lipschitz test functions of the restricted balls for the infinite
surface. The limiting spaces there are
length spaces, and only the prelimit spaces lack this property.

\begin{theorem}[Brownian-plane limit]\label{thm:wp-brownian-plane}
As $L\to\infty$, the following convergences hold in distribution:
\begin{align}
 &(X(S_\infty^*,o_\infty),L^{-1}d,
       v_{\rm WP}^{-1}L^{-4}\bar A,o_\infty)
   \ \Longrightarrow
   (\mathcal P,D_\infty,\mu_\infty,\rho_\infty),
       \label{eq:bp-ghp}\\
 &(S_\infty^*,L^{-1}d,
       v_{\rm WP}^{-1}L^{-4}dA,o_\infty)
   \ \Longrightarrow
   (\mathcal P,D_\infty,\mu_\infty,\rho_\infty).
       \label{eq:bp-gp}
\end{align}
The topology is local pointed GHP in \eqref{eq:bp-ghp} and local
pointed GP in \eqref{eq:bp-gp}. Both statements also hold with
$(S_n,o_n)$ and $L=r_n$ for every deterministic sequence
$1\ll r_n\ll n^{1/4}$.
\end{theorem}

The main step of the proof, given in Appendix~\ref{app:brownian-plane},
is a coupling of finite and infinite WP balls at every radius $o(n^{1/4})$.
An exact likelihood for a tree prefix marked by a uniform leaf,
together with mixing of the continuous angle types, permits the prefix
depth to grow. A functional central limit theorem for the spine labels
then gives the necessary geometric separation. Passing through a finite
surface of intermediate size transfers the Brownian-sphere tangent
coupling of~\cite{CLGPlane}; its coding construction preserves the full
volume measure as well as distances.

\begin{corollary}[Radial volume profile]\label{cor:wp-bp-volume}
Put $V_\infty(t)=\mu_\infty(B_t(\rho_\infty))$. Then
\begin{equation}\label{eq:bp-volume-profile}
 \left(\frac{A(B_{Lt}(S_\infty^*,o_\infty))}
                   {v_{\rm WP}L^4}\right)_{t\ge0}
 \Longrightarrow (V_\infty(t))_{t\ge0}
 \quad\text{in }C_{\rm loc}([0,\infty)).
\end{equation}
In particular, $\log A(B_L(S_\infty^*,o_\infty))/\log L\to4$
in probability. The same conclusions hold for $(S_n,o_n)$ and $L=r_n$
in the mesoscopic range of Theorem~\ref{thm:wp-brownian-plane}.
\end{corollary}

\begin{proof}
The Brownian-plane radial measure has a continuous positive density
on $(0,\infty)$~\cite[Proposition~1]{LeGallArea}. Thus every positive
radius has zero boundary mass, and $V_\infty$ is continuous and finite,
with $V_\infty(t)>0$ for $t>0$. Theorem~\ref{thm:wp-brownian-plane}
gives joint convergence of the ball masses at any finite list of radii.
The prelimit profiles are continuous and nondecreasing: Riemannian
distance spheres have area zero. On a finite mesh, monotonicity bounds
each intermediate value between its two endpoint values. Uniform
continuity of the limiting profile therefore upgrades these joint
limits to convergence in $C_{\rm loc}$. Finally the logarithm of the
normalized volume at radius one is tight, proving the stated exponent.
\end{proof}

The cusp trimming cannot be dropped from the Hausdorff statement: the
untrimmed metric spaces do not converge locally in GH. Indeed,
choose any $K$ distinct cusp tails in $S_\infty^*$ and entrance points
at distances at most $a_K<\infty$ from the root. Points at depth $L/2$
in these tails lie in $B_{3L/4}(o_\infty)$ when $L\ge4a_K$, and
are pairwise at distance at least $L$. The rescaled unit balls thus
fail GH tightness. Since the deep cusps carry exponentially small
mass, this is compatible with \eqref{eq:bp-gp}. The theorem describes
geometry and area; the energy comparison needed for the dense
microscopic spectrum remains a separate problem.

\part*{Concluding remarks}
\section{Perspectives and open problems}\label{sec:conclusions}

The uniform expansion separates the changing leading growth of
$V_{g,n}$ from its inverse-genus corrections. The saddle parameter
absorbs the former; the positive string-flow recursion controls the
latter through both endpoints. Painlev\'e~I fixes the normalization at
the pure-$\psi$ endpoint, while the positive convolution factor and
Poisson extraction make the passage to individual volumes uniform
from bounded cusp counts to the critical singularity. The resulting
volume ratios provide the counting input for the spectral laws and
the local geometry.

The two applications use this input differently. When
$n/\sqrt g\to c\in(0,\infty)$, rare two-cusp pairs of pants support
funnel ground states. A trace comparison proves that they
asymptotically exhaust the positive spectrum on each fixed compact
energy window below $1/4$. In the supercritical Poisson regimes, the
lowest eigenvalues are controlled by collar conductances divided by
enclosed areas. At sufficiently large cusp density, low-energy
functions inside growing planar caps invalidate that normalization.
Fixed-energy windows still admit deterministic limits, because orbital
statistics see the aggregate contribution of the many local pieces.
The distinction between these limits leaves the dense microscopic
problem below unresolved.

At $n\asymp g$, all planar cusp types contribute to the same carrier
integrals. Retaining based holonomies makes these integrals determine
whole pointed neighborhoods, yielding the local WP family. Its
cusp-filled surfaces are conformally discs at each finite positive
density; the infinite-density endpoint is the recurrent infinite WP
sphere, whose cusp-filling is conformally the plane, as Budd and
Curien conjectured. A separate
growing-prefix coupling gives the Brownian-plane metric and area
limit. Thus the common volume estimates support two further arguments:
spectral exhaustion for the low eigenvalues and holonomy control for
the local family. Passing from the latter's metric limit to an energy
limit remains an additional problem.

\subsection{Open problems}
We list the problems that the present methods leave open. Each lies
beyond a specific obstruction identified in the proofs.

\textbf{1. The dense microscopic spectrum.}
Theorem~\ref{thm:dense-persistence} proves the Poisson law at the scale
$A_g$ under \eqref{eq:dense-persistence-range}, which holds for
$n\le g^{11/9-\delta}$, and Theorem~\ref{thm:dense-collapse} shows that
this normalization fails for $n\gg g^3/(\log g)^4$. Neither range is
asserted to be sharp. Between them one must determine the correct scale
and the limit law of the smallest eigenvalues. Upper bounds come from
trial functions on large planar caps; the missing ingredient is a
matching lower bound on the whole low spectral subspace, which requires
control of the interaction between growing planar caps and the rest of
the surface. Corollary~\ref{cor:dense-residual} already identifies the
spectral type of these eigenvalues at every dense growth rate.

\textbf{2. Edge fluctuations without short pants.}
Theorem~\ref{thm:moving-edge} gives the Fr\'echet law of $1/4-G(X)$ when
$n/\sqrt g\to0$ and $(n/\sqrt g)\log g\to\infty$; it does not cover the
whole subcritical range. When short two-cusp pants are too rare to
dominate, and in particular for closed surfaces, the fluctuations of the
spectral gap around $1/4$ remain open; Huang's theorem~\cite{HuangGap}
and the polynomial rate of Hide, Macera and Thomas~\cite{HMT} give
convergence but no limit law.

\textbf{3. Small cusp forms.}
With probability tending to one, there are no cusp forms with eigenvalue
in $(0,1/4-\varepsilon]$ when $n=O(\sqrt g)$
(Theorem~\ref{thm:critical-eigenvalues}), and none below $T/A_g$ when
$n/g\to\infty$ (Corollary~\ref{cor:dense-residual}). Whether Weil--Petersson
random surfaces carry any cusp form below $1/4$ in the other regimes,
and how many, is open; this is the probabilistic counterpart of the
question of Phillips and Sarnak~\cite{PhillipsSarnak} and of the
conjecture of Otal and Rosas~\cite[Section~6]{OtalRosas} that, for
$g\ge2$, at most $2g-3$ small eigenvalues are cuspidal; see~\cite{BMP}, and~\cite{Mondal}
for partial results. Relatedly, the
regularized measure of Theorem~\ref{thm:regularized-spectrum} does not
separate cusp forms from scattering, and the fixed-window limits are
compatible with a vanishing proportion of eigenvalues collapsing to
zero.

\textbf{4. Resurgence uniformly in the cusp density.}
The single-factorial bounds in Theorem~\ref{thm:intro-absolute} give
Borel transforms convergent near the origin. Uniform analytic
continuation, growth bounds for Borel summation and the corresponding
Stokes data remain to be related to the nonperturbative predictions
of~\cite{EGGLS}, uniformly in $n/g$.

\textbf{5. Near-critical scaling of the local family.}
The local family of Theorem~\ref{thm:wp-local-family} is obtained at
fixed radius and fixed positive density. The near-critical scaling
question of Budd and Curien~\cite[Section~7.3]{BC} requires estimates
uniform as the radius and the density grow together. In the discrete
setting the corresponding limit is Budzinski's hyperbolic Brownian
plane~\cite{BudzinskiHBP}; a continuum Weil--Petersson analogue is not
known.

\textbf{6. Energy on the Brownian plane.}
Theorem~\ref{thm:wp-brownian-plane} gives a metric and area limit at the
planar endpoint, but a limiting Dirichlet form or diffusion clock needs
further energy control. This is also the obstacle to deriving a dense
microscopic spectral law from the Brownian-plane geometry alone.

\textbf{7. Mixed intersection numbers.}
Relative asymptotics for mixed intersection numbers when both the number
of nonzero degrees and their sum grow with the genus, uniformly over
their allocation, are not supplied by the present methods: the bounded
boundary factorization of Proposition~\ref{prop:bounded-boundaries}
does not give such individual coefficient estimates.

The uniform volume expansion determines the weights of finite geometric
configurations across cusp densities. The spectral and local-limit
results show what can be recovered from these weights, while the dense
microscopic spectrum and joint geometric scaling require further control
of energy and growing neighborhoods.

\clearpage
\appendix
\addtocontents{toc}{\protect\newpage}
\part*{Proof appendices}
\section{Proof locations and notation}\label{app:reader-guide}\label{subsec:notation}

\subsection{Proof locations}
For each main result, the following table gives the place of its proof and
the principal results and methods used there. References to Sections~\ref{the-positive-compact-family-of-kernels}--\ref{uniform-comparison} include the supporting Appendices~\ref{app:kernels} and~\ref{app:primitive-bounds}.
In the table, Thm., Cor., Prop., Sec.\ and App.\ abbreviate Theorem,
Corollary, Proposition, Section and Appendix.
{\small
\renewcommand{\arraystretch}{1.18}
\hyphenpenalty=10000\exhyphenpenalty=10000
\begin{longtable}{@{}>{\raggedright\arraybackslash}p{0.27\textwidth}>{\raggedright\arraybackslash}p{0.25\textwidth}>{\raggedright\arraybackslash}p{0.42\textwidth}@{}}
\toprule
Result & Proved in & Main ingredients\\
\midrule
\endhead
\multicolumn{3}{@{}l}{\emph{Part~\ref{part:volumes}: volume asymptotics}}\\*[1pt]
Thm.~\ref{thm:intro-absolute} & Sec.~\ref{sec:absolute} & Thms.~\ref{thm:full-expansion} and~\ref{thm:canonical-expansion} (App.~\ref{app:marked} and~\ref{app:extraction})\\
Cor.~\ref{cor:constant} (Zograf's constant) & Sec.~\ref{sec:prefactor} & Secs.~\ref{the-positive-compact-family-of-kernels}--\ref{ratios-with-a-fixed-number-of-distinguished-markings}; Painlev\'e~I at the pure-$\psi$ endpoint; \cite{MZ}\\
Thm.~\ref{thm:volume}, limits & Sec.~\ref{extracting-fixed-cusp-counts} & Secs.~\ref{the-positive-compact-family-of-kernels}--\ref{ratios-with-a-fixed-number-of-distinguished-markings}; fixed-cusp limits~\cite{Mir}\\
Thm.~\ref{thm:volume}, rate $O(g^{-1})$ & Sec.~\ref{sec:absolute} & Thm.~\ref{thm:sharp-extraction} (App.~\ref{app:extraction}), which uses Thm.~\ref{thm:full-expansion} at a fixed order (App.~\ref{app:marked})\\
Thm.~\ref{thm:coefficient}, equivalent & Sec.~\ref{uniform-one-boundary-coefficient-ratios} & Prop.~\ref{prop:harmonic}; log-concavity; \cite{MZ}\\
Thm.~\ref{thm:coefficient}, rate & App.~\ref{sec:quantitative-coefficients} & as above\\
Thm.~\ref{thm:chain} & Sec.~\ref{uniform-logarithmic-volume-rate-and-the-full-diagonal-chain} & limits in Thm.~\ref{thm:volume}; \cite{ManZ}\\[3pt]
\multicolumn{3}{@{}l}{\emph{Part~\ref{part:spectra}: spectral laws}}\\*[1pt]
Cors.~\ref{cor:spectral}, \ref{cor:many-spectrum}; Thm.~\ref{thm:pants-poisson} & Sec.~\ref{sec:geometry}; App.~\ref{app:short-curves} & limits in Thm.~\ref{thm:volume}; Secs.~\ref{ratios-with-a-fixed-number-of-distinguished-markings}, \ref{sec:prefactor} and~\ref{extracting-fixed-cusp-counts}\\
Thm.~\ref{thm:many-residual} & Sec.~\ref{sec:residual-count} & fixed-length pants concentration; doubling and cusp averages; residual spectral projection\\
Thm.~\ref{thm:tight-spectrum}, Cor.~\ref{cor:tight-boltzmann} & App.~\ref{app:tight} & Secs.~\ref{ratios-with-a-fixed-number-of-distinguished-markings}, \ref{sec:prefactor} and~\ref{extracting-fixed-cusp-counts}; \cite{BL,BZ}\\
Thm.~\ref{thm:critical-cluster} & App.~\ref{app:critical} & volume bounds of~\cite{MZ,SW}; zero-density limits in Thm.~\ref{thm:volume}; Prop.~\ref{prop:bounded-boundaries}\\
Thm.~\ref{thm:critical-eigenvalues}, Cor.~\ref{cor:gap-phase-diagram} & Sec.~\ref{sec:critical-eigenvalues}, using App.~\ref{app:critical-trace} & pants counts; Huang trace expansion~\cite{HuangGap}; model spectrum~\cite{BMM,Schulze,BJP}; count matching\\
Thm.~\ref{thm:moving-edge}, Cors.~\ref{cor:edge-residual}, \ref{cor:edge-moments} & App.~\ref{app:moving-edge} & moving expected count; conditional isolation; count matching; rare-event moments\\
Thm.~\ref{thm:supercritical-spectrum} & Sec.~\ref{sec:supercritical-spectrum}, using App.~\ref{app:super-cuts}, \ref{app:super-comparison} & all-marking envelope; coefficient ratios; marked curves; relative Cheeger and collar comparisons\\
Thm.~\ref{thm:dense-persistence} & App.~\ref{app:dense-comparison} & growing-cap law; branch contraction; comparison on the whole low spectral subspace\\
Thm.~\ref{thm:dense-collapse} & App.~\ref{app:dense-collapse-details}, using App.~\ref{app:cusp-cutoffs}--\ref{app:cusp-opening} & growing caps; finite planar resistance; energy transport; genus bound\\
Cor.~\ref{cor:dense-residual} & Sec.~\ref{sec:dense-spectrum}; App.~\ref{app:cuspidal-exclusion} & cusp-penalized Cheeger comparison; all-topology volume sum; fixed-index upper tightness\\
Thm.~\ref{thm:fixed-energy-spectrum} & App.~\ref{app:sublinear-trace}, \ref{app:fixed-window-proofs} & saddle-summed complementary tails; Huang's trace coefficients; model count exhaustion\\
Thm.~\ref{thm:positive-density-spectrum} & App.~\ref{app:positive-density-orbits}, \ref{app:fixed-window-proofs} & planar carrier sums; two-carrier concentration; positive truncated trace\\
Thm.~\ref{thm:cusp-dominated-spectrum} & App.~\ref{app:cusp-dominated-orbits}, \ref{app:fixed-window-proofs} & endpoint graph bounds; fixed-genus and planar convolutions; positive trace\\
Thm.~\ref{thm:regularized-spectrum}, Cor.~\ref{cor:cusp-density-determinant} & App.~\ref{app:fixed-window-proofs} & orbital limits; uniform long-repetition tail; classical regularized trace and integrated residual-pole bound; heat-to-determinant integral\\[3pt]
\multicolumn{3}{@{}l}{\emph{Part~\ref{part:geometry}: local geometry}}\\*[1pt]
Thm.~\ref{thm:wp-local-family} & App.~\ref{app:rooted-local-law} & short-word support; bounded rooted fibers; one- and two-carrier coefficients of App.~\ref{app:positive-density-orbits}\\
Thm.~\ref{thm:wp-conformal-phase} & App.~\ref{app:conformal-phase}, \ref{app:cusp-cutoffs} & finite-cap multiplicity; topological anchor; recurrent endpoint in App.~\ref{app:cusp-cutoffs}\\
Cor.~\ref{cor:wp-geodesic-flow} & Sec.~\ref{sec:local-conformal-phase} & Hopf--Tsuji--Sullivan theorem\\
Thm.~\ref{thm:wp-brownian-plane} & App.~\ref{app:brownian-plane} & growing-prefix coupling; joint finite measured limit; measured tangent coupling\\
\bottomrule
\end{longtable}}
The volume-ratio limits and the short-pants consequences in Section~\ref{sec:geometry} use the main-text estimates and Appendices~\ref{app:kernels}--\ref{app:primitive-bounds}. The subsequent length and spectral results use the additional arguments in Appendices~\ref{app:short-curves}--\ref{app:fixed-window-proofs}, and the geometric results those in Appendices~\ref{app:rooted-local-law}--\ref{app:brownian-plane}, as indicated above. Within each part the appendices follow the order in which they are used. Appendices~\ref{app:marked}--\ref{sec:quantitative-coefficients} are needed only for the inverse-genus rates and for the expansions to higher order, apart from the all-marking envelope \eqref{eq:gevrey-envelope}, which Part~\ref{part:spectra} uses. They can be read independently of one another: Appendix~\ref{app:marked} after Section~\ref{sec:higher-order}, Appendix~\ref{app:extraction} after Section~\ref{sec:sharp-extraction} (it uses Theorem~\ref{thm:full-expansion} only through its statement), and Appendix~\ref{sec:quantitative-coefficients} after Section~\ref{uniform-one-boundary-coefficient-ratios}.

\Needspace{10\baselineskip}
\subsection{Notation} The first table lists the principal symbols of Part~\ref{part:volumes} and where they are defined; the second does the same for the recurring symbols of Parts~\ref{part:spectra} and~\ref{part:geometry}.
{\small
\renewcommand{\arraystretch}{1.15}
\begin{longtable}{@{}lll@{}}
\toprule
Symbol & Meaning & Defined in\\
\midrule
\endhead
$X$, $a=X'$, $y_0$, $x_c$ & genus-zero relation, $X'(y)=J_0(2\sqrt y)$ & Subsection~\ref{subsec:parameter}, Section~\ref{the-positive-compact-family-of-kernels}\\
$x=X(Y)$, $Y\in[0,y_0)$ & cusp weight and parameter & Subsection~\ref{subsec:parameter}\\
$\Phi_g$ & normalized cusp generating function & \eqref{eq:intro-Phi}\\
$K_Y$, $H_Y$ & positive kernels & \eqref{eq:P1}\\
$\bar K_Y$, $\bar H_Y$ & kernels rescaled by $r(Y)$ & \eqref{eq:P9}\\
$\theta_\ell$, $w_\ell$ & coefficients of $1/\bar K_Y$, increments & \eqref{eq:P11}\\
$\vartheta_\ell$ & coefficients of $1/K_Y$ & Lemma~\ref{lem:kernel-positivity}\\
$r(Y)=2\Delta(Y)$ & first positive zero of $K_Y$ & Sections~\ref{the-positive-compact-family-of-kernels}, \ref{the-stationary-parameter}\\
$T$, $\Delta=T-Y$ & stationary root and half-zero & \eqref{eq:SP1}, \eqref{eq:SP2}\\
$I$, $\lambda$ & action and cusp density & \eqref{eq:SP2}\\
$\Psi=\lambda/2$ & leading cusp-count scale per unit complexity & Section~\ref{the-stationary-parameter}\\
$\mathscr D=x\,\mathrm d/\mathrm dx$ & logarithmic derivative operator & Section~\ref{the-stationary-parameter}\\
$h$ & complexity, $2g-2+n$ or $2g-2+p$ & see below\\
$m$ & deficit $3g-3+n-|\mathbf d|$ & Section~\ref{exact-amplitudes-and-the-harmonic-model}\\
$U^Y$, $\mathcal U$ & shifted amplitudes, normalized & \eqref{eq:H6}, \eqref{eq:H7}\\
$\mathcal M$, $A_m(h)$ & harmonic model & \eqref{eq:H12}\\
$\overline U$, $W_{g,p}$, $c_{g,p}$ & marked amplitudes and degree profiles & \eqref{eq:TL9}\\
$A$, $\Theta$, $M_1$, $\mathcal Z=\sqrt2\,I$ & kernel scalars and action & \eqref{eq:P9}, \eqref{eq:PF1}\\
$\mathcal J_{g,p}$ & leading term for $\Phi_g^{(p)}$ & \eqref{eq:PF3}\\
$\widehat V_{g,n}$, $\varpi$ & leading term for $V_{g,n}$ & \eqref{eq:intro-leading}, \eqref{eq:intro-varpi}\\
$B_{g,n}$, $C_{g,n}$, $D_{g,n}$ & volume ratios & Subsection~\ref{subsec:intro-ratios}\\
$\mathcal B$, $\mathcal C$ & volume-ratio profiles & Subsection~\ref{subsec:intro-ratios}\\
$C_g(\mathbf e)$ & normalized pure $\psi$ intersection number & \eqref{eq:H4}\\*
$\mathfrak p_k=\mathrm e^{-k}k^k/k!$ & Poisson point mass at its mean & Section~\ref{sec:sharp-extraction}\\
\bottomrule
\end{longtable}}

\Needspace{8\baselineskip}
\noindent In Parts~\ref{part:spectra} and~\ref{part:geometry}, $X$ denotes a random surface, not the genus-zero relation of Part~\ref{part:volumes}, and $\lambda$ is the cusp density $n/g\to\lambda$.
{\small
\renewcommand{\arraystretch}{1.15}
\hyphenpenalty=10000\exhyphenpenalty=10000
\begin{longtable}{@{}>{\raggedright\arraybackslash}p{0.28\textwidth}>{\raggedright\arraybackslash}p{0.40\textwidth}>{\raggedright\arraybackslash}p{0.26\textwidth}@{}}
\toprule
Symbol & Meaning & Defined in\\
\midrule
\endhead
\multicolumn{3}{@{}l}{\emph{Part~\ref{part:spectra}: spectral laws}}\\*[1pt]
$\mathbb P^{\mathrm{WP}}_{g,n}$, $\Delta_X$ & normalized WP law on $\mathcal M_{g,n}$; nonnegative Laplacian of $X$ & Subsection~\ref{subsec:intro-geometry}\\
$\lambda_j(X)$, $G(X)$ & eigenvalues below $1/4$ with multiplicity; spectral gap, equal to $1/4$ if none & Subsection~\ref{subsec:intro-geometry}, Section~\ref{sec:critical-eigenvalues}\\
$\Xi_X$ & counting measure of the spectrum in $(0,1/4)$ & Theorem~\ref{thm:critical-eigenvalues}\\
$\lambda_r^+(X)$ & $r$th positive eigenvalue & Theorem~\ref{thm:dense-collapse}\\
$\mathfrak h(X)$ & Cheeger constant & Subsection~\ref{sec:spectral}\\
$\operatorname{PPP}(\nu)$ & Poisson point process with intensity $\nu$ & Subsection~\ref{sec:pants-spectrum}\\
$\mathcal P_2(X)$, $\mathcal P_k(X)$ & simple geodesics bounding a two-cusp pants; simple geodesics whose planar side has $k$ cusps & Subsection~\ref{sec:pants-spectrum}, Theorem~\ref{thm:pants-poisson}\\
$q_g$ & two-cusp short-curve intensity & \eqref{eq:PS1}\\
$N_L$ & number of curves in $\mathcal P_2(X)$ of length at most $L$ & Subsection~\ref{sec:pants-spectrum}\\
$Y_\lambda$, $x_\lambda$ & density parameter with $\lambda(Y_\lambda)=\lambda$; $x_\lambda=X(Y_\lambda)$, $x_0=0$, $x_\infty=x_c$ & Subsection~\ref{subsec:bounded-boundaries}, Theorem~\ref{thm:pants-poisson}\\
$w_k(x)$, $W(x)$ & weights of the planar cusp types; their sum & \eqref{eq:PC1}\\
$\operatorname{sys}(X)$ & systole & Theorem~\ref{thm:pants-poisson}\\
$\alpha(\lambda)$, $\nu$ & tight length scale; tight intensity $(\cosh t-1)\,\mathrm dt/t$ & \eqref{eq:tight-scale}\\
$\mathbb P^{\mathrm{WP}}_{g,\mu}$ & Boltzmann measure with cusp weight $\mu$ & Subsection~\ref{subsec:tight-spectrum}\\
$P_s$, $K_s$ & two-cusp pants with waist $s$; its primitive length counting measure & Subsection~\ref{sec:critical-cluster}\\
$\mu_c$, $\mathcal Q_c$, $\mathcal N$, $\Lambda_c$ & waist intensity $c^2\sinh(s/2)\,\mathrm ds/(8\pi^2)$; independent Poisson processes of intensities $\mu_c$ and $\nu$; critical cluster limit & \eqref{eq:cluster-intensities}, Theorem~\ref{thm:critical-cluster}\\
$F_s$, $\beta(s)$ & complete two-cusp funnel surface with convex core $P_s$; its eigenvalue below $1/4$ & Section~\ref{sec:critical-eigenvalues}, \eqref{eq:critical-model-energy}\\
$\nu_c$ & critical energy intensity, $\beta_*\mu_c$ & \eqref{eq:critical-energy-intensity}\\
$\kappa_0=1/(8\pi^2)$, $h_*$, $\mathcal D_X$ & edge constant; model deficit function; deficit point process & \eqref{eq:edge-model-coordinate}, Theorem~\ref{thm:moving-edge}\\
$x_g$, $M_2(x)$, $A_g$ & supercritical saddle $(2g-1)\Psi(x_g)=n-1$; area-weighted type sum; microscopic spectral scale & \eqref{eq:super-scale}\\
$\epsilon$, $Q$, $S$, $\eta$ & dense-regime parameters, $A_g\sim2\pi^2S$ & \eqref{eq:dense-parameters}, \eqref{eq:dense-scales}\\
$\mu_*$, $\nu_*$ & sublinear waist intensity and fixed-window limit & \eqref{eq:fixed-energy-measure}\\
$\nu_\lambda$, $\nu_\infty$ & fixed-window limits at positive and infinite density & Theorems~\ref{thm:positive-density-spectrum}, \ref{thm:cusp-dominated-spectrum}\\
$\mathcal O_X(H)$, $\Omega_\lambda(H)$ & orbital functional of Huang's trace expansion; its planar limit & Definition~\ref{rev:def-huang-orbital}, \eqref{eq:positive-orbital-functional}\\
$\mathfrak m_X$, $\mathfrak m_\theta$ & regularized spectral measure; its limit when $g/n\to\theta$ & \eqref{eq:regularized-measure}, Theorem~\ref{thm:regularized-spectrum}\\[3pt]
\multicolumn{3}{@{}l}{\emph{Part~\ref{part:geometry}: local geometry}}\\*[1pt]
$A_{g,n}$, $B_r(S,o)$ & hyperbolic area $2\pi(2g-2+n)$; pointed measured $r$-ball & Section~\ref{sec:local-conformal-phase}\\
$\mathcal L_\lambda$, $(S_\lambda,o_\lambda)$ & law and pointed surface of the local WP family & Theorem~\ref{thm:wp-local-family}\\
$(S_\infty^*,o_\infty)$ & area-rooted infinite WP sphere & \eqref{eq:wp-local-endpoints}\\
$S_n$, $A_n$, $o_n$ & WP sphere with $n+1$ cusps, its area, a uniform area point & Section~\ref{sec:brownian-plane}\\
$S^\circ$, $X(S,o)$, $P$, $\bar A$ & cusp-trimmed surface; rooted trimmed space; cusp projection; projected area & \eqref{eq:bp-projected-area}\\
$c_{\rm WP}$, $v_{\rm WP}$ & distance and volume normalizations & \eqref{eq:bp-normalization}\\
$(\mathcal P,D_\infty,\mu_\infty,\rho_\infty)$, $V_\infty$ & Brownian plane; its radial volume profile & Section~\ref{sec:brownian-plane}, Corollary~\ref{cor:wp-bp-volume}\\
\bottomrule
\end{longtable}}

The complexity parameter $h$ is specified locally: we use $h=2g-2+n$ for volumes with $n$ boundary components and $h=2g-2+p$ for generating functions with $p$ distinguished markings. In particular, $h=2g-1$ in \eqref{eq:intro-leading}, where $p=1$ and the remaining cusps are counted by $k$. Positive constants $C$, $c$ may change from line to line, and a subscript, as in $O_p$, indicates the parameters on which an implied constant depends. 

A few letters in Part~\ref{part:volumes} have an additional local meaning, stated where it occurs: $D$ is the degree sum $|\mathbf d|$ in Sections~\ref{exact-amplitudes-and-the-harmonic-model}--\ref{sec:prefactor} and in Appendix~\ref{app:marked}; $h_0=2g-2$ is the complexity of a closed surface (Appendices~\ref{app:marked}, \ref{app:tight} and~\ref{app:critical}); the harmonic coefficients $A_m(h)$ are unrelated to the kernel scalar $A(Y)$ and to the coefficients $A_j(g,n)$; the Laplacian $\Delta_X$ of a surface $X$ is unrelated to $\Delta(Y)$; and $T$ is a length window in Appendix~\ref{app:short-curves} and an energy window $T/A_g$ in Corollary~\ref{cor:dense-residual}. 

In Parts~\ref{part:spectra} and~\ref{part:geometry} the following letters also carry more than one meaning. The cusp density $\lambda$ of \eqref{eq:SP2} is unrelated to the eigenvalues $\lambda_j(X)$ and $\lambda_r^+(X)$; $s$ is a waist length in $P_s$, $F_s$, $\beta(s)$ and $\mathcal Q_c$, but the spectral parameter of $\Phi_X(s)$ and $Z_X(s)$ in Subsections~\ref{sec:regularized-spectrum} and~\ref{sec:spectral-determinant}; $c$ is the critical ratio $n/\sqrt g\to c$ in Sections~\ref{sec:critical-eigenvalues} and~\ref{sec:moving-edge} and in $\mathcal Q_c$, $\mu_c$, $\nu_c$ and $\Lambda_c$, apart from its use as a generic constant and in $c_{\rm WP}$.

$Q$ is the parameter $nx/(4P)$ in \eqref{eq:dense-parameters} and a polynomial degree in Section~\ref{sec:critical-eigenvalues} and Appendix~\ref{app:critical-trace}, and neither is the pants process $\mathcal Q_c$; $K$ denotes the kernel $K_Y$ in Part~\ref{part:volumes}, the length measure $K_s$ in Subsection~\ref{sec:critical-cluster}, the systolic cusp count $K_g$ in Theorem~\ref{thm:pants-poisson}, the moment constants $K_p$ in Corollary~\ref{cor:edge-moments} and the planar series $K_x(L)$ of \eqref{eq:positive-orbital-functional} in Section~\ref{sec:fixed-energy-spectrum} and Appendix~\ref{app:positive-density-orbits}; $H$ denotes the kernel $H_Y$ in Part~\ref{part:volumes}, an even test function in $\mathcal O_X(H)$ and $\Omega_\lambda(H)$, and the level-set function $H(t)$ of Section~\ref{sec:moving-edge}.

$D$ is also the gap defect $D(X)=1/4-G(X)$ and the deficit process $\mathcal D_X$ in Section~\ref{sec:moving-edge}, the operator $x\partial_x$ in Section~\ref{sec:supercritical-spectrum}, and the Brownian-plane distance $D_\infty$; $P$ is the constant $2\pi^2$ in \eqref{eq:dense-parameters}, the pants $P_s$, the cusp projection in Section~\ref{sec:brownian-plane} and, as $\mathcal P$, the Brownian plane; $S$ is the dense spectral scale in \eqref{eq:dense-parameters} and a surface in $S_\lambda$, $S_n$ and $S_\infty^*$; and $\nu$ without subscript is the tight intensity \eqref{eq:tight-scale}, whereas $\nu_*$, $\nu_\lambda$, $\nu_\infty$, $\nu_c$ and $\nu_{\rm edge}$ are spectral limit measures.

\part*{Appendices to Part~1: volume asymptotics}
\section{Proofs for the kernels and the parameter}\label{app:kernels}

This appendix proves the three lemmas on which the kernel family and
the density parameter of Part~\ref{part:volumes} rest:
Lemma~\ref{lem:kernel-zeros} (real zeros of \(K_Y\)),
Lemma~\ref{lem:kernel-positivity} (positivity of the Taylor
coefficients in \(Y\), which supplies the convolution factor of
Lemma~\ref{lem:convolution}) and Lemma~\ref{lem:parameter} (the
functions \(T\), \(\Delta\), \(I\), \(\lambda\) and \(\Psi\), whose
endpoint expansions are used in
Sections~\ref{extracting-fixed-cusp-counts} and~\ref{sec:absolute}).
The proofs use only the genus-zero relation \(X'(y)=J_0(2\sqrt y)\),
the explicit kernels \eqref{eq:P1} and classical facts about entire
functions of order \(1/2\); no volume estimate enters, and the only
dependence among them is that the later two proofs use the sign
properties of \(X\) and the genus-zero product for \(a\) established in
the proof of Lemma~\ref{lem:kernel-zeros}. The notation is that of
Sections~\ref{the-positive-compact-family-of-kernels}
and~\ref{the-stationary-parameter}, and the subsections follow the order
of the three lemmas.

\subsection{Proof of Lemma~\ref{lem:kernel-zeros}}\label{app:kernel-zeros}

The proof has three stages, following the sketch given after the
statement of the lemma in Subsection~\ref{real-zeros-of-the-kernel}. We
first verify the sign properties of \(X\) and its derivatives and locate
the zeros of \(a\). We then write \(K_Y\) along characteristics as
\(\sin t\) minus a positive correction \eqref{eq:P3} and count its real
zeros by concavity. Finally a Rouch\'e comparison with \(\sin t/t\) on
large circles shows that the real zeros exhaust all zeros, which gives
the product \eqref{eq:P5}.

We use the elementary identities \(X''(y)=-X(y)/y\) and
\(X'''(y)=\sum_{j\ge0}(-1)^jy^j/(j!(j+2)!)\). On
\([0,\pi^2/4]\), \(X(y)>0\) for \(y>0\), \(X''<0\), and \(X'''>0\). For
the first assertion we use the alternating-series
lower bound \[
 X(y)/y\ge1-y/2+y^2/12-y^3/144>0\quad(0\le y\le5/2);
\] the cubic on the right is decreasing and is positive at \(5/2\). The
other assertions follow from the displayed identities and
\(X'''(y)\ge1/2-y/6>0\). Similarly, the partial sum of the series for
\(a(3/2)\) through \(j=4\) is an upper bound by the alternating-series
test, and it is negative; this proves \(y_0<3/2\). In particular
\(X(Y)\le Y<3/2\).

All zeros of \(a(y)\) are positive real. Indeed, if \(a(y)=0\), the
function \(u(s)=J_0(2\sqrt y\,s)\) solves \((su')'+4ysu=0\), is regular
at zero, and vanishes at one. Integration against \(\bar u\) gives \[
 4y\int_0^1s|u|^2\,\mathrm{d}s=\int_0^1s|u'|^2\,\mathrm{d}s>0.
\] Since the zeros are positive real and \(y_0\) is the smallest of
them, \(a\) has no zeros in \(|y|<y_0\). Moreover the series is entire
of order \(1/2\), so its Hadamard product has genus zero and no
nonconstant exponential factor; this product is used in the proof of
Lemma~\ref{lem:kernel-positivity} below.

We next prove that every zero of \(K_Y\) is positive and simple. At
\(Y=0\), \[
 K_0(z)=\frac{\sin\sqrt{2z}}{\sqrt{2z}}.
\] Coefficient comparison in \eqref{eq:P1} gives \[
 \partial_YK_Y=2\partial_zK_Y+(K_Y-a(Y))/z.
\] Integration along the characteristics \(z+2Y=\text{constant}\) yields
\begin{equation}
\sqrt{t^2-4Y}\ K_Y((t^2-4Y)/2)=\sin t-R_Y(t),\qquad
 R_Y(t)=2\int_0^Y\frac{a(u)}{\sqrt{t^2-4u}}\,\mathrm{d}u.       \label{eq:P3}
\end{equation} For \(0<Y<y_0\) and real \(t>2\sqrt Y\), \(R_Y>0\), \(R_Y'<0\), \(R_Y''>0\), and \[
 R_Y(t)\le\frac{2X(Y)}{\sqrt{t^2-4Y}}.
\] At \(Y=0\), \(R_0\) vanishes identically, as in the sinc formula above. For \(Y>0\), \(\sin t-R_Y(t)\) is strictly concave on each positive
sine interval. On the first interval \((2\sqrt Y,\pi)\), it tends to
zero at the left endpoint, is positive immediately afterwards by
\(K_Y(0)=a(Y)>0\), and is negative at \(\pi\); hence it has exactly one
simple interior zero. On \((2j\pi,(2j+1)\pi)\), \(j\ge1\), it is
negative at both endpoints and positive at the midpoint, since \[
 R_Y((2j+\tfrac12)\pi)
 <\frac3{\sqrt{(5\pi/2)^2-6}}<1.
\] Strict concavity gives exactly two simple zeros there. Negative sine
intervals contain no zero. For real \(z\le0\), every argument of \(X'\)
in \eqref{eq:P1} is less than \(y_0\), where \(X'>0\); for negative arguments
this follows from the positive modified-Bessel series. Thus \(K_Y\) has
no nonpositive real zero.

To exclude nonreal zeros, set \(G_Y(t)=K_Y((t^2-4Y)/2)\). It is entire
and even. Formula \eqref{eq:P3}, with the square root asymptotic to \(t\), holds
outside the segment \([-2\sqrt Y,2\sqrt Y]\). On \(|t|=(N+1/2)\pi\), \begin{equation}
\frac{G_Y(t)}{\sin t/t}
 =\frac t{\sqrt{t^2-4Y}}
    \left(1-\frac{R_Y(t)}{\sin t}\right)=1+O(N^{-1}).   \label{eq:P4}
\end{equation} The error is uniform: \(R_Y=O(N^{-1})\), and \(|\sin t|\) is bounded
below on these circles. For the latter assertion use
\(|\sin t|\ge\sinh1\) when \(|\Im t|\ge1\); otherwise \(|\Re t|\)
differs from \((N+1/2)\pi\) by \(O(N^{-1})\). Rouch\'e's theorem gives
exactly \(2N\) zeros of \(G_Y\) inside each sufficiently large circle.
The real zeros already counted give exactly \(N\) positive and \(N\)
negative zeros there: for even \(N\) the circle ends at the midpoint of
a positive sine interval, after its first zero, and for odd \(N\) it
ends in the subsequent negative interval. They exhaust the count,
proving the assertion. Consequently the product formula \eqref{eq:P5}
holds. This is again a genus-zero product, since \eqref{eq:P1} has order \(1/2\).

\subsection{Proof of Lemma~\ref{lem:kernel-positivity}}\label{app:kernel-positivity}

The coefficients of \(\log(a/H_Y)\) satisfy a system of differential
equations in \(Y\) whose right-hand sides have nonnegative coefficients.
Since these coefficients are positive at \(Y=0\), where the zeros of
\(H_0\) are positive real, induction on the order of differentiation
shows that all their Taylor coefficients in \(Y\) are nonnegative. The
same argument applies to \(\log(a/K_Y)\).

We first prove \eqref{eq:P6}. By \eqref{eq:P1}, \[
 \partial_YH_Y=2\partial_zH_Y+3(H_Y-a(Y))/z.
\] Writing \(\mathcal T=\log(a/H_Y)\) and \(\gamma=-a'/a\), this becomes \[
 \partial_Y\mathcal T=-\gamma+2\partial_z\mathcal T+3(\mathrm{e}^{\mathcal T}-1)/z.
\] The constant identity is \(\gamma=5t_1\); for \(m\ge1\), \begin{equation}
t_m'=(2m+5)t_{m+1}
  +3\sum_{\ell=2}^{m+1}\frac1{\ell!}
       \sum_{\substack{i_1+\cdots+i_\ell=m+1\\ i_j\ge1}}
          t_{i_1}\cdots t_{i_\ell}.                 \label{eq:P7}
\end{equation} At \(Y=0\), \[
 H_0(z)=3\frac{\sin\sqrt{2z}-\sqrt{2z}\cos\sqrt{2z}}
                    {(\sqrt{2z})^3}.
\] Its zeros \(\rho_j\) are positive real. For a direct proof, a nonzero root \(t\)
of \(\sin t-t\cos t\) gives the regular function \(u(s)=j_1(ts)\),
satisfying \((s^2u')'+(t^2s^2-2)u=0\) and \(u(1)=0\). Integration
against \(\bar u\) shows \[
 t^2\int_0^1s^2|u|^2\,\mathrm{d}s
   =\int_0^1(s^2|u'|^2+2|u|^2)\,\mathrm{d}s>0.
\] Here \(j_1(w)=\sin w/w^2-\cos w/w\) is regular at zero. The root is
simple because \((\sin t-t\cos t)'=t\sin t\) does not vanish at a
nonzero root. The order-\(1/2\) product therefore gives
\(t_m(0)=m^{-1}\sum_j\rho_j^{-m}>0\).

Induct on the derivative order \(q\), simultaneously for all indices
\(m\), in \eqref{eq:P7}. Every differentiated term is a product of derivatives of
order at most \(q\) with nonnegative coefficients, so the induction is
valid; each fixed equation has only finitely many terms. For each \(m\),
\(t_m(Y)\) is a polynomial in the ratios \(X^{(j)}(Y)/a(Y)\), so it is
analytic on \(|Y|<y_0\), and its Taylor series converges to
\(t_m(Y)\) there. This proves \eqref{eq:P6}, and gives strictly positive
coefficients of \(H_Y^{-h}\) for every real \(h>0\), and nonnegative
coefficients of \(-H_Y'/H_Y\).

The same argument applied to \(S=\log(a/K_Y)=\sum_{m\ge1}s_m(Y)z^m\)
gives \begin{equation}
s_m'=(2m+3)s_{m+1}
  +\sum_{\ell=2}^{m+1}\frac1{\ell!}
       \sum_{\substack{i_1+\cdots+i_\ell=m+1\\ i_j\ge1}}
           s_{i_1}\cdots s_{i_\ell}.               \label{eq:P8}
\end{equation} The initial sinc product and the same analytic argument imply that
every Taylor coefficient in \(Y\) of \[
 \widehat\vartheta_\ell(Y):=a(Y)[z^\ell]K_Y(z)^{-1}
\] is nonnegative. The product for \(a\) implies the same for \(a^{-1}\)
and for \(\vartheta_\ell=[z^\ell]K_Y^{-1}\). This completes the proof.

\subsection{Proof of Lemma~\ref{lem:parameter}}\label{app:parameter-proof}

We rewrite the equation for the stationary root in a normalized form
that is strictly monotone in \(T\); the implicit function theorem then
gives analyticity. The remaining assertions follow from integral
formulas for \(I\), \(J\) and \(\beta\) and from expansions at the two
endpoints of the parameter interval.

We use from Lemma~\ref{lem:kernel-zeros} that $X>0$,
$X''<0$ and $X'''>0$ on the relevant interval, and that $0<y_0<3/2$.
The four parts of the lemma are proved in order.

\emph{Part~(i): the root \(T\).}
Normalize the left side of \eqref{eq:SP1} by \(2\sqrt{T-Y}\). Its root equation
becomes \begin{equation}
Q(Y,T)=\int_0^1X'(T(1-s^2)+Ys^2)\,\mathrm{d}s=0.
 \label{eq:SP8}
\end{equation} It is strictly decreasing in \(T\). At \(T=y_0\) it is positive when
\(Y<y_0\). The series and the beta integral give the Sonine identity \begin{equation}
\int_0^T\frac{X'(v)}{\sqrt{T-v}}\,\mathrm{d}v
       =\sin(2\sqrt T).
 \label{eq:SP9}
\end{equation} At \(T_0=\pi^2/4\), the full integral is zero, and removing \([0,Y]\)
makes it negative when \(Y>0\). Since \(Q>0\) for \(T\le y_0\) and
\(Q\) is strictly decreasing in \(T\), this proves existence and
uniqueness of the root of \eqref{eq:SP1} in \((y_0,\pi^2/4]\), including
\(T(0)=T_0\). The derivative \(Q_T\) never vanishes
there, so the implicit function theorem gives analyticity. Also
\(T'<0\), and \(T\to y_0\) as \(Y\to y_0\).

Twice differentiating \eqref{eq:SP8}, with \(v=T(1-s^2)+Ys^2\), gives \[
 T''(Y)=
 -\frac{\int_0^1X'''(v)[T'(1-s^2)+s^2]^2\,\mathrm{d}s}
        {\int_0^1X''(v)(1-s^2)\,\mathrm{d}s}>0.
\] At \((Y,T)=(0,T_0)\), \eqref{eq:SP9} gives \(Q_T=-1/(2T_0)\), and
differentiation of the lower limit \(Y\) in
\(\int_Y^TX'(v)(T-v)^{-1/2}\,\mathrm dv=2\sqrt{T-Y}\,Q\) gives
\(Q_Y=-X'(0)/(2T_0)=-1/(2T_0)\); hence \(T'(0)=-1\). Consequently \(-1\le T'<0\).

The integral \(I\) is positive. Indeed, subtracting \((T-y_0)\) times
\eqref{eq:SP1} from its definition gives \begin{equation}
I(Y)=\int_Y^T
 \frac{X'(v)(y_0-v)}{\sqrt{T-v}}\,\mathrm{d}v>0.
 \label{eq:SP10}
\end{equation}
The integrand is nonnegative and not identically zero, because \(X'\) is
decreasing on \([0,\pi^2/4]\) and vanishes at \(y_0\).

\emph{Part~(ii): the density \(\lambda\).}
Introduce \(J\) and \(\beta\) as in \eqref{eq:SP11}.
The equality for \(J\) is integration by parts, and the last equality
for \(\beta\) follows from the definition of \(\lambda\). If
\(Q_2=\int_0^1s^2X'(T-\Delta s^2)\,\mathrm{d}s\), then \(Q_2>0\): the function
\(X'(T-\Delta s^2)\) is strictly increasing in \(s\) and has mean zero
by \eqref{eq:SP8}. Moreover \(Q_2\le X'(Y)/3\). Differentiating \(J\) gives \begin{equation}
J'=(1-T')Q_2\le\frac23X'(Y),
 \qquad
 \beta'=\frac{X'J-XJ'}{J^2}\ge\frac{X'}{3J}>0.
 \label{eq:SP12}
\end{equation} At zero, integrating \eqref{eq:SP9} with respect to \(T\) gives
\[
 \int_0^TX'(v)\sqrt{T-v}\,\mathrm{d}v
 =\tfrac14\bigl(\sin2\sqrt T-2\sqrt T\cos2\sqrt T\bigr),
\]
whose value at \(T_0=\pi^2/4\) is \(I(0)=\pi/4\); then \(J(0)=0+(\pi/4)/(\pi/2)=1/2\), and
\(\lambda(Y)=4Y+O(Y^2)\). At the other endpoint \(J\to x_c\) and
\(X(Y)\to x_c\), so \(\beta\to1\). This proves the asserted analytic
bijection and the values in \eqref{eq:SP3} at zero.

\emph{Part~(iii): expansions at the critical endpoint.}
The endpoint constants are explicit. Put \(\delta=y_0-Y\) and
\(b=-X''(y_0)=x_c/y_0>0\). The implicit equation at \((y_0,y_0)\) has
\(T'=-1/2\); this gives the first formula in \eqref{eq:SP13}. In
particular \(T-y_0=\delta/2+O(\delta^2)=\Delta/3+O(\Delta^2)\), so
that \(v-y_0=\Delta/3-\Delta s^2+O(\Delta^2)\) for \(v=T-\Delta s^2\). For the
second formula, substitute
\(X'(T-\Delta s^2)=-b(v-y_0)+O(\Delta^2)=b\Delta(s^2-1/3)+O(\Delta^2)\) into
\(I=2\Delta^{3/2}\int_0^1s^2X'(T-\Delta s^2)\mathrm{d}s\). Together with
\(X(Y)\to x_c\), this proves \eqref{eq:SP14}, and \eqref{eq:SP15} follows
from \eqref{eq:SP3}. Since \(1-\beta=2/(2+\lambda)\), \eqref{eq:SP14}
also gives \(\Delta\asymp\sqrt{1-\beta}\).

\emph{Part~(iv): the function \(\Psi\).}
By part~(ii), \(\Psi\) is analytic and strictly increasing for
\(0<x<x_c\). The first line of \eqref{eq:SP16} follows from
\(\lambda(Y)=4Y+O(Y^2)\) and \(x=X(Y)=Y+O(Y^2)\). For the second line, quadratic inversion gives
\(x_c-x=(b/2)\delta^2(1+O(\delta))\). Combining this with \eqref{eq:SP13} and
\eqref{eq:SP14} gives the constant \(5/4\). All expansions here arise from
convergent analytic expansions in \(x\) or \(\sqrt\varepsilon\), so
their asserted derivatives follow by termwise differentiation. In
particular \(\Psi^2/(\mathscr D\Psi)\) has finite endpoint limits \(0\)
and \(5/4\), and is bounded on the full interval.

\section{Uniform bounds for positive intersection numbers}
\label{app:primitive-bounds}

This appendix proves Lemma~\ref{lem:primitive-bounds}, the two-sided
bound $1/6\le C_g(\mathbf e)\le1$ for the normalized pure-$\psi$
intersection numbers of \eqref{eq:H4}. The bound is the only
intersection estimate that the harmonic comparison of
Section~\ref{uniform-comparison} imports: it enters the terminal bounds
of Lemma~\ref{lem:terminal}, and through them
Proposition~\ref{prop:harmonic}. The lemma
is a consequence of the stronger bounds in~\cite[Lemmas~2--6]{GYZ};
the proof below uses their separating-term estimate, which we reproduce,
and replaces their two-variable majorant by a scalar one, so that a
single induction on the complexity $h$ gives both bounds. The only
external inputs are the classical one-point formula
$\langle\tau_{3g-2}\rangle_g=1/(24^gg!)$ and the two-point bound of
Delecroix, Goujard, Zograf and Zorich, proved in~\cite{DGZZ} and quoted
as~\cite[equation~(6)]{DGZZBounds}. After the statement we record the
normalized string, dilaton and DVV identities, outline the induction,
and then give the proof in five marked parts.

\begin{lemma}\label{lem:primitive-bounds}
For $g\ge1$ and positive integers $d_1,\ldots,d_n$ with
$\sum d_i=3g-3+n$, the normalization in~\eqref{eq:H4} satisfies
\[
 \frac16\le C_g(d_1,\ldots,d_n)\le1.
\]
\end{lemma}

Write $C(\mathbf d)=C_g(\mathbf d)$ and $h=2g-2+n$ when the
genus is determined by the degrees. The normalized string and dilaton
identities are
\begin{equation}
 C(0,\mathbf d)=\sum_{d_j>0}\frac{2d_j+1}{3(Z-1)}
                    C(\mathbf d-\mathbf e_j),\qquad
 C(1,\mathbf d)=C(\mathbf d),
 \label{eq:primitive-string}
\end{equation}
where $Z$ is the complexity of $(0,\mathbf d)$; the terminal value
$C_1(1)=1/6$ is kept separately. For a pivot $d_1\ge2$, DVV gives
\begin{align}
 C(\mathbf d)&=\sum_{j=2}^n\frac{2d_j+1}{3(h-1)}
 C(d_2,\ldots,d_j+d_1-1,\ldots,d_n)\notag\\
 &\quad+\frac{2}{3(h-1)}\sum_{a+b=d_1-2}
 C_{g-1}(a,b,d_2,\ldots,d_n)\notag\\
 &\quad+\sum_{\substack{a+b=d_1-2\\I\sqcup J=\{2,\ldots,n\}}}
 \frac{\Gamma(h_1)\Gamma(h_2)}{6\Gamma(h)}
 C(a,\mathbf d_I)C(b,\mathbf d_J).
 \label{eq:primitive-dvv}
\end{align}
Here $h_1,h_2$ are the complexities of the two factors, and only stable,
dimension-compatible terms are included. The weights in the first two
lines sum to one.

Both bounds are proved by induction on complexity. The lower bound
uses the DVV recursion with its nonnegative separating terms
discarded. For the upper bound we prove the stronger estimate
$C(\mathbf d)\le F(h)$ for tuples with positive entries, where $F\le1$
is an explicit nondecreasing function. In the induction step the
separating terms contribute at most $2/((h-1)(h-2))$; the cases
$h\le7$ are treated directly.

\begin{proof}[Proof of Lemma~\ref{lem:primitive-bounds}]
\emph{The lower bound.}
We prove the same lower bound for positive-genus tuples with at most two
zero entries, by induction on complexity. For a tuple with $z\in\{1,2\}$
zeros, the total string weight is $(3h-z)/(3h-3)\ge1$.
String lowers complexity without increasing the number of zeros.
For a positive tuple, remove the ones by dilaton. Unless the tuple
reduces to the genus-one base, its genus is at least two. Its merger
children are positive and its handle children have positive genus and
at most two zeros. Dropping the nonnegative separating terms in
\eqref{eq:primitive-dvv} therefore proves the lower bound.

\emph{The upper majorant and zero entries.}
Define the nondecreasing function
\[
 F(h)=\begin{cases}1/3,&1\le h\le7,\\
                   1-4/(h-1),&h\ge7.
       \end{cases}
\]
We will prove $C(\mathbf d)\le F(h)$ for positive tuples.
Its increment for $h\ge8$ is $4/((h-1)(h-2))$.
The corresponding bound for a tuple with one zero is $C\le1$:
remove all ones, and let $Z$ be the remaining complexity. String bounds
the result by
\[
 \left(1+\frac{2}{3(Z-1)}\right)F(Z-1).
\]
For $Z\le8$ this is at most $5/9$; for $Z\ge9$, substitution of
$F(Z-1)=1-4/(Z-2)$ gives a bound below one. This observation uses only
positive tuples of smaller complexity and so is available in the induction.

\emph{The separating contribution.}
Suppose all $d_i\ge2$. Both separating genera are positive by dimension.
For a fixed subset $I$ of size $m$ and genus $g_1$, the degree $a$ is
determined. The total separating weight is thus at most
\[
 \sum_{m=0}^{n-1}\binom{n-1}{m}
 \sum_{\substack{g_1+g_2=g\\g_1,g_2\ge1}}
 \frac{(2g_1+m-2)!(2g_2+n-m-3)!}{6(h-1)!}.
\]
As in~\cite[Lemma~3]{GYZ}, the terms with at least one genus equal to
one, allowing double counting, are bounded by
\[
 \frac{1}{3(h-1)(h-2)}
 \sum_{m=0}^{n-1}\prod_{j=1}^m\frac{n-j}{h-2-j}
 \le\frac{3}{2(h-1)(h-2)}.
\]
Indeed $n\le3h/5$ and $n\le h-2$. The terms with $m\le2$ are at most
one, and for $m\ge3$ each factor with $j\ge3$ is at most $3/5$, so
the sum is at most $1+1+\sum_{m\ge2}(3/5)^{m-2}=9/2$.
For $g_1,g_2\ge2$, the inequality
$\binom ab\binom cd\le\binom{a+c}{b+d}$ gives the bound
\[
 \frac{n(g-3)}{6(h-1)(h-2)\binom{2g-4}{2}}
 \le\frac{3}{8(h-1)(h-2)}.
\]
This part vanishes for $g\le3$; otherwise use $n\le3g-3$.
The total weight is therefore at most $2/((h-1)(h-2))$.
Each factor is positive or has one zero, so the same bound applies to
the full separating contribution once the smaller-complexity bounds hold.

\emph{Initial cases.}
Put $a_g=C_g(3g-2)$ and $b_g=C_g(0,3g-1)$.
The one-point formula and the two-point bound of~\cite{DGZZ}, quoted
in~\cite[equation~(6)]{DGZZBounds}, give
\[
 C_g(d_1,d_2)\le b_g=\frac{6g-1}{6g-3}a_g<\frac1\pi<\frac8{25}.
\]
For completeness, $b_{g+1}/b_g=(6g+5)(6g+1)(6g-3)/
(108g(g+1)(2g-1))>1$, and Stirling's formula gives $b_g\to1/\pi$.
All positive genus-two tuples reduce by dilaton to $(4)$, $(2,3)$ or
$(2,2,2)$. String gives
\[
 C_1(0,2)=\frac5{18},\quad C_1(0,0,3)=\frac{35}{108},\quad
 C_1(0,0,2,2)=\frac{25}{81}.
\]
Using these values in~\eqref{eq:primitive-dvv} yields
\[
 C_2(4)=\frac{35}{144},\quad C_2(2,3)=\frac{1015}{3888},\quad
 C_2(2,2,2)=\frac{175}{648}<\frac3{11}.
\]
The other two values are also below $3/11$. Hence a genus-two tuple
with one zero is at most $(11/9)(3/11)=1/3$ after removing ones;
its remaining complexity is at least four. The genus-one one-zero
case reduces to $(0,2)$ and has the same bound.

For $h\le7$, the only additional positive cores with at least three
markings are the genus-three tuples $(2,2,5)$, $(2,3,4)$ and $(3,3,3)$.
For the first two, two string steps bound the handle term by
$(16/15)(11/9)(3/11)=16/45$. Each separating factor has genus one
or two and is at most $1/3$; the separating contribution is therefore at
most the general separating weight bound $2/((h-1)(h-2))=2/30$ at
$h=7$, multiplied by the bound $1/3$ for each of the two factors, that
is, by $1/9$. Thus
\[
 C\le\frac89\frac8{25}+\frac19\frac{16}{45}
              +\frac19\frac2{30}=\frac{671}{2025}<\frac13.
\]
For $(3,3,3)$, both handle terms reduce to $(0,3,3)$ in genus two, so
\[
 C\le\frac79\frac8{25}+\frac29\frac13
              +\frac19\frac2{30}=\frac{223}{675}<\frac13.
\]
Together with dilaton, the one- and two-point bounds and the genus-one
base, this proves every initial case $h\le7$.

\emph{Induction step.}
Let $\mathbf d$ have no entry equal to one (ones having been removed by
dilaton) and complexity $h\ge8$; assume $n\ge3$ and choose a
minimal degree $d=d_1\ge2$. If $d\ge3$, all handle children except
the two endpoint terms have positive degrees. Each endpoint has one
zero and at most one degree-one entry. Removing that entry if necessary
and applying string bounds it by
$(1+2/(3(h-3)))F(h-1)$. The total endpoint weight is $4/(3(h-1))$.
Since $F\le1$, the excess over $F(h-1)$ is at most
\[
 \frac{8}{9(h-1)(h-3)}+\frac2{(h-1)(h-2)}
 \le\frac{46}{15(h-1)(h-2)}<F(h)-F(h-1).
\]
If $d=2$, the handle has two zeros. Two string steps, allowing an
intervening dilaton step, give the factor
\[
 R_h=\frac{(3h-5)(3h-10)}{(3h-6)(3h-12)}.
\]
The resulting positive tuples have complexity at most $h-3$.
The handle weight is $2/(3(h-1))$, and
\[
 \frac{2(R_h-1)}{3(h-1)}
 =\frac{2(9h-22)}{27(h-4)(h-1)(h-2)}
 \le\frac{25}{27(h-1)(h-2)}.
\]
Adding separation gives $79/(27(h-1)(h-2))$, again below the increment
of $F$. This closes the induction and proves the upper bound.
\end{proof}

\section{The marked expansion beyond the leading term}\label{app:marked}

This appendix proves Theorems~\ref{thm:first-correction}
and~\ref{thm:full-expansion}, the higher-order refinements of
Proposition~\ref{prop:prefactor} stated in Section~\ref{sec:higher-order}.
The final marked asymptotic is stated for fixed $p$; the intermediate
recursive estimates also control larger marking ranges. The strategy
is the one of Section~\ref{sec:prefactor} carried to higher
order: the adjacent-degree subtraction in the recursion \eqref{eq:H9}
expands the normalized profiles $c_{g,p}$ in powers of $1/h$; the
dilaton and KdV equations convert these into expansions of the cusp
and genus ratios and hence of the logarithmic derivative in $Y$; and the
integration constants are fixed at the pure-$\psi$ endpoint, where
Painlev\'e~I supplies the asymptotics of $\langle\tau_2^{3g-3}\rangle_g$.
The inputs are the marked ratios of Proposition~\ref{prop:marked-ratios}
together with the estimates \eqref{eq:TL14}--\eqref{eq:TL27} from its
proof, the leading term of Proposition~\ref{prop:prefactor}, the Itzykson--Zuber
identity~\cite{IZ}, and the Painlev\'e~I asymptotics in the form
of~\cite{GIKM} and~\cite{CCH}.
Subsections~\ref{app:second-degree}--\ref{app:critical-matching} compute
one further term; the coefficient to be integrated vanishes at the
endpoint, so its primitive stays bounded, and this gives both the
$O_p(g^{-1})$ relative error for the leading term and its first explicit
correction (Theorem~\ref{thm:first-correction}).
Subsections~\ref{sec:full-expansion}--\ref{app:remainder-bounds}
construct the full formal series and prove polynomial bounds on its
coefficients and remainders (Proposition~\ref{prop:hierarchy});
Subsections~\ref{app:gevrey-one}--\ref{app:gevrey-remainders} sharpen
these to single-factorial bounds (Proposition~\ref{prop:gevrey-hierarchy});
and Subsection~\ref{sec:endpoint-expansion} determines the endpoint
constants to all orders and completes the proof of
Theorem~\ref{thm:full-expansion}.

\subsection{The second degree coefficient}\label{app:second-degree}

We first extend the degree expansion \eqref{eq:PF5} of the normalized
profiles by one order (Lemma~\ref{lem:second-profile}). The argument is
the adjacent-degree subtraction used for \eqref{eq:PF5} in
Section~\ref{sec:prefactor}, now
keeping the $h^{-2}$ terms of the prefactors and the two smallest
separating configurations.

We retain $\delta_p$ and $\beta_p$ from \eqref{eq:PF5} and
\eqref{eq:PF7}. The kernel coefficients $\theta_\ell$ are those of \eqref{eq:H2}. Set
\begin{equation}
 \alpha_p=A\{(p-3)\Theta+2M_1+2p-3\},\qquad
 \mathcal T_j=\sum_{\ell\ge\max(0,j)}w_\ell.
 \label{eq:ST1}
\end{equation}
For a pivot $i$, put $t=d_i$ and let $z$ count zero entries among the
other $p-1$ markings. Let $\mathbf d^{ij,\ell}$ be the merger child,
with new degree $t+d_j+\ell-1$, and let $\mathbf d^{i;u,v}$ be the
nonseparating child with the pivot replaced by $u,v$. Define
\begin{equation}
\begin{split}
 B_{p,i}(\mathbf d)={}&
 A\sum_{j\ne i}(2d_j+1)
 \sum_{\ell\ge\max(0,1-t-d_j)}w_\ell
 [1-\beta_{p-1}-\delta_{p-1}(\mathbf d^{ij,\ell})]\\
 &+2A\sum_{\ell\ge\max(0,2-t)}w_\ell
 \sum_{\substack{u+v=t+\ell-2\\u,v\ge0}}
 [1-\beta_{p+1}-\alpha_p-\delta_{p+1}(\mathbf d^{i;u,v})]\\
 &+A^2\left\{\left[\binom z2+\frac{\theta_1}2\right]
                    \mathcal T_{2-t}+\frac12\mathcal T_{3-t}\right\}.
\end{split}
 \label{eq:ST2}
\end{equation}
The first line is absent at $p=1$. Starting at zero along a fixed
coordinate-increment path, define $e_p$ by
\begin{equation}
 e_p(\mathbf0)=0,\qquad
 e_p(\mathbf d+\mathbf e_i)-e_p(\mathbf d)=-B_{p,i}(\mathbf d).
 \label{eq:ST3}
\end{equation}
The formulas converge absolutely by \eqref{eq:TL10}; they give
$B_{p,i}=O_p((D+1)^3)$ and $e_p=O_p((D+1)^4)$.

\begin{lemma}[Second profile correction]\label{lem:second-profile}
For fixed $p\ge1$, uniformly in $Y$ and $D\le h^{1/8}$,
\begin{equation}
 c_{g,p}(\mathbf d)=1-\frac{\delta_p(\mathbf d)}h
       +\frac{e_p(\mathbf d)}{h^2}
       +O_p\!\left(\frac{(D+1)^6}{h^3}\right).
 \label{eq:ST4}
\end{equation}
The function $e_p$ is independent of the chosen path.
\end{lemma}
\begin{proof}
We first need the linear term in the genus exchange. For $p\ge2$,
the two genus-zero three-marking contributions to \eqref{eq:TL18}
give $T_{g,p}=2(p-2)A/h+O_p(h^{-2})$. Every other genus-zero term
loses at least two cusp factors, and the positive-genus portion is
$O_p(h^{-2})$ by \eqref{eq:TL19}. Thus \eqref{eq:TL14} and
\eqref{eq:PF5} imply
\begin{equation}
 \frac{W_{g-1,p+2}}{W_{g,p}}=1-\frac{\alpha_p}h+O_p(h^{-2}).
 \label{eq:ST5}
\end{equation}
Indeed the coefficient is
$\delta_p(\mathbf e_1)+3(p-2)A=\alpha_p$.
For $p=0,1$, descend from $p+1$ using \eqref{eq:PF7}; the cusp ratio
changes its first coefficient by
$\beta_{p+2}-\beta_p=-A(\Theta+2)$, while
$\alpha_{p+1}-\alpha_p=A(\Theta+2)$.
Consequently the two adjacent-difference prefactors are
\begin{equation}
\begin{split}
 \frac{W_{g,p-1}}{W_{g,p}}
 &=\frac Ah\left(1+\frac{1-\beta_{p-1}}h\right)+O_p(h^{-3}),\\
 2\frac{W_{g-1,p+1}}{W_{g,p}}
 &=\frac{2A}h\left(1+\frac{1-\beta_{p+1}-\alpha_p}h\right)
     +O_p(h^{-3}).
\end{split}
 \label{eq:ST6}
\end{equation}
The first line is used only at $p\ge2$. The merger and nonseparating children both have complexity
$h-1$, which accounts for the two constants $+1$.

We next retain the leading separating terms of the adjacent difference.
If the smaller component has genus $j$ and contains $s$ of the other
external labels, the complementary component has complexity $h-(2j+s)$.
Stability forces $2j+s\ge2$.
Equality occurs only for $(j,s)=(0,2)$ and $(1,0)$. The first component
is a three-marked sphere, so its three degrees must be zero; it gives
$\binom z2\mathcal T_{2-t}$. The second is a one-marked torus, whose
$W$ amplitudes are $\theta_1/2$ and $1/2$ at degrees zero and one.
It gives $(\theta_1/2)\mathcal T_{2-t}+(1/2)\mathcal T_{3-t}$.
The two orientations cancel the factor $1/2$ in the recursion.
In both cases the large all-zero component divided by the parent is
$A^2/h^2+O_p(h^{-3})$, by fixed iterations of \eqref{eq:PF7} and
\eqref{eq:ST5}. Its degrees can be replaced by zero at cost
$O_p((D+\ell+1)^2/h)$.

After removing these two types the all-zero separating mass is
$O_p(h^{-3})$. To verify this sharper bound, repeat the proof of
\eqref{eq:TL20}: the remaining genus-zero terms have $s\ge3$, the
remaining genus-one terms have $s\ge1$, and $j\ge2$ starts at
$h^{-4}$. The same factorial-geometric series sums these bounds,
with any fixed polynomial in $j$; \eqref{eq:TL21} handles the
macroscopic splits. At $p=1$, insert $2j+1$ cusps in the larger
component before raising its genus to marking count two, then use
$W_{g,2}/W_{g,1}=O(g)$. This bounds its ratio by
$O((C/g)^{2j})$. After the two genus-one orientations are removed,
the factorial series starts at $j=2$ and is $O(h^{-4})$.
This argument uses the marking-independent constant in \eqref{eq:TL16}.
Pair counts and the moments of $w_\ell$ give an adjacent-difference
error $O_p((D+1)^3/h^3)$ for the separating part.

Now substitute \eqref{eq:PF5} and \eqref{eq:ST6} into the exact
adjacent subtraction of \eqref{eq:H9}. The two retained separating
types give the last line of \eqref{eq:ST2}. The result is
\[
 c(\mathbf d)-c(\mathbf d+\mathbf e_i)
 =\frac{\delta_p(\mathbf d+\mathbf e_i)-\delta_p(\mathbf d)}h
  +\frac{B_{p,i}(\mathbf d)}{h^2}
  +O_p((D+1)^5/h^3).
\]
To justify the substitution, cut the kernel index $\ell$ at $h^{1/8}$. Below this cutoff, child
degrees are at most $2h^{1/8}+O_p(1)$, inside the range of
\eqref{eq:PF5} for all sufficiently large $h$. Its remainder, multiplied
by the linear prefactor and degree count, is
$O_p((D+\ell+1)^5/h^3)$. Above the cutoff use $0\le c\le1$ and
polynomial growth of the displayed coefficients. Since the uniform
moments \eqref{eq:TL10} of $w_\ell$ are finite to every order, these
tails are $O_p(h^{-3})$, even after summing at most $D$ increments. The dimension cutoff is farther away and is
covered by the same argument. Summing increments proves \eqref{eq:ST4}.
For a fixed tuple, expansions along two paths must have the same
$h^{-2}$ coefficient as $g\to\infty$, proving path independence.
\end{proof}

\emph{Comparison with the fixed-cusp expansion.}
At $Y=0$ one has $A=1/\pi^2$, $\Theta=2$ and $M_1=3/2$, and the
normalization \eqref{eq:TL9} gives $W_{g,p}/W_{g,p+1}=4V_{g,p}/V_{g,p+1}$
and $W_{g-1,p+2}/W_{g,p}=V_{g-1,p+2}/V_{g,p}$. With $h=2g-2+p$,
equations \eqref{eq:PF7} and \eqref{eq:ST5} become
\[
 \frac{V_{g,p+1}}{8\pi^2gV_{g,p}}
 =1+\frac{(\frac12-\pi^{-2})p-\frac54+2\pi^{-2}}{g}+O_p(g^{-2}),
 \qquad
 \frac{V_{g-1,p+2}}{V_{g,p}}=1+\frac{3-2p}{\pi^2g}+O_p(g^{-2}),
\]
and \eqref{eq:PF5} gives, for fixed $k\ge1$,
\[
 \frac{[\tau_k\tau_0^{p-1}]_{g,p}}{V_{g,p}}
 =1-\frac{k^2+(p-\frac52)k-\frac p2+\frac32}{\pi^2g}+O_p(g^{-2}).
\]
These are the first-order formulas of Mirzakhani and
Zograf~\cite[Theorems~1.4 and~3.2(i)]{MZ}, recovered here as the
$Y=0$ case of the parameter-dependent coefficients.

\subsection{The coefficient in the differential equation}\label{app:transport-coefficient}

Inserting Lemma~\ref{lem:second-profile} into the dilaton identity gives
the next coefficient $\gamma_p$ of the differential equation~\eqref{eq:PF8}.
We then evaluate it with finitely many kernel moments and check that it
vanishes at the pure-$\psi$ endpoint.

Apply \eqref{eq:ST4} to the dilaton child of complexity $h+1$ and set
\begin{equation}
 E_p=\sum_{d\ge1}b_d e_{p+1}(d,0^p),\qquad
 \gamma_p=\beta_p^2-\beta_p-E_p/A.
 \label{eq:ST7}
\end{equation}
The moments of $b_d$ justify this sum and yield
\begin{equation}
 \frac{\mathrm{d}}{\mathrm{d}Y}\log\Phi_g^{(p)}(X(Y))
 =\frac{h+\beta_p+\gamma_p/h+O_p(h^{-2})}{Ar}.
 \label{eq:ST8}
\end{equation}
Indeed, replacing $h+1$ by $h$ in the child profile gives
$1-\delta/h+(\delta+e)/h^2$.
Thus the dilaton ratio is
$A\{1-\beta_p/h+(\beta_p+E_p/A)/h^2+O_p(h^{-3})\}$,
and inversion proves \eqref{eq:ST8}. To justify the sum, truncate the
dilaton degree at $h^{1/8}$; the tail is bounded using
$1+d^2/h+d^4/h^2$ and a sufficiently high absolute moment of $b_d$
from \eqref{eq:TL13}.

The coefficient in \eqref{eq:ST7} can be evaluated using finitely many
kernel moments, as follows.
Put $P=p+1$, and write $B_P(t)=B_{P,1}(t,0^{P-1})$. With
$k=t+\ell-1$ and $s=t+\ell-2$, \eqref{eq:ST2} becomes
\[
\begin{split}
 R_P(k)&=1-\beta_{P-1}
       -A\{\Theta(k^2+(P-5)k)+2M_1k+(4-P)\mathbf1_{k\ge1}\},\\
 S_P(s)&=(s+1)(1-\beta_{P+1}-\alpha_P)\\
 &\quad-A\{(s+1)(\Theta(s^2+(P-3)s)+2M_1s)
                   +(5-2P)s-1+\mathbf1_{s=0}\}.
\end{split}
\]
Then
\begin{equation}
\begin{split}
 B_P(t)={}&A(P-1)\sum_{\ell\ge\max(0,1-t)}w_\ell R_P(t+\ell-1)\\
 &+2A\sum_{\ell\ge\max(0,2-t)}w_\ell S_P(t+\ell-2)\\
 &+A^2\left\{\left[\binom{P-1}2+\frac{\theta_1}2\right]
                     \mathcal T_{2-t}+\frac12\mathcal T_{3-t}\right\}.
\end{split}
 \label{eq:ST9}
\end{equation}
The zero-degree correction in $S_P$ follows by summing the
$q$ term of $\delta_{P+1}$: for $s\ge1$ the two endpoint pairs have
$q=1$ and the $s-1$ interior pairs have $q=2$.

For $t\ge3$, remove both lower cutoffs, replace
$\mathbf1_{k\ge1}$ by one, $\mathbf1_{s=0}$ by zero, and both tails
by $\Theta$. This gives a cubic polynomial
$B_\infty(t)=a_3t^3+a_2t^2+a_1t+a_0$.
Its coefficients use only $M_j=\sum\ell^jw_\ell$, $0\le j\le3$.
Put $\Delta_j=B_P(j)-B_\infty(j)$ for $j=0,1,2$, and define
\begin{equation}
 \begin{split}
 \Sigma(d)&=\frac{a_3}4d^2(d-1)^2+
       \frac{a_2}6d(d-1)(2d-1)+\frac{a_1}2d(d-1)+a_0d,\\
 e_P(d,0^{P-1})&=-\Sigma(d)-\sum_{j=0}^2\Delta_j\mathbf1_{d\ge j+1}.
 \end{split}
 \label{eq:ST10}
\end{equation}
The three exceptional values use the same moments and $w_0,w_1,w_2$;
there is no remaining infinite recursion. Set
\[
 L_j=\frac13[(z\partial_z+1)^j\bar H_Y(z)]_{z=1}\quad(1\le j\le4),
 \qquad b_1=\frac13,\quad b_2=-\frac{\theta_1}5,
\]
and
\[
 \Sigma_b=\frac{a_3}4(L_4-2L_3+L_2)
   +\frac{a_2}6(2L_3-3L_2+L_1)
   +\frac{a_1}2(L_2-L_1)+a_0L_1.
\]
Then the finite formula for \eqref{eq:ST7} is
\begin{equation}
 \gamma_p=\beta_p^2-\beta_p+
 \frac{\Sigma_b+\Delta_0A+\Delta_1(A-b_1)+\Delta_2(A-b_1-b_2)}A.
 \label{eq:ST11}
\end{equation}
All ingredients extend analytically to $Y=y_0$, by the removable
normalized-kernel extension used in \eqref{eq:PF14}.

At that endpoint $A=2/15$, $\Theta=1$, $M_j=0$ for $j\ge1$,
$w_0=\theta_1=1$, and $\beta_p=-p/5$. Equation~\eqref{eq:ST9}
gives $B_P(0)=0$ and $B_P(1)=A(P-1)(P+4)/5$.
Only $b_1=1/3$ and $b_2=-1/5$ remain, so
$E_p/A=p(p+5)/25=\beta_p^2-\beta_p$. In particular,
\begin{equation}
 \gamma_p(y_0)=0,\qquad \gamma_p(Y)=O_p(y_0-Y),\qquad
 \gamma_p/(Ar)=O_p(1).
 \label{eq:ST12}
\end{equation}
This agrees with an exact computation: string and dilaton at the pure-$\psi$
endpoint give $W_{g,p+1}/W_{g,p}=3(h+m)=(15h-3p)/2$,
so that $A\,W_{g,p+1}/(hW_{g,p})=1-p/(5h)$ exactly, with $A=2/15$;
there is no term of order $h^{-2}$, and $\gamma_p=0$.

\subsection{Critical matching to the next order}\label{app:critical-matching}

This subsection proves Theorem~\ref{thm:first-correction}.
Recall $t_p$ and $\mathcal L_p$ from \eqref{eq:ST13}.
By \eqref{eq:ST9}--\eqref{eq:ST12}, $\mathcal L_p$ is an explicit
bounded function, given by finitely many kernel derivatives and a
convergent one-dimensional integral.
The proof integrates the next-order transport equation \eqref{eq:ST8}
down from the pure-$\psi$ endpoint, using the positive comparison
\eqref{eq:PF15}--\eqref{eq:PF16} near the endpoint as in
Section~\ref{sec:prefactor}. There the first correction is fixed
by two inputs: a sharper asymptotic formula for the critical
intersection numbers, whose derivation in
Subsection~\ref{sec:endpoint-expansion} is independent of the transport
argument, and the first correction to the critical harmonic coefficient
of \eqref{eq:PF13}.

\begin{proof}[Proof of Theorem~\ref{thm:first-correction}]
The sharper critical input is
$C_g(2^{3g-3})=\pi^{-1}(1-2/(9g)+O(g^{-2}))$.
It is derived in Subsection~\ref{sec:endpoint-expansion}, as the
coefficient $d_1=-2/9$ in \eqref{eq:FE10}, from the Itzykson--Zuber
identity, the large-order relation of Garoufalidis--Its--Kapaev--Mari\~no
and the first one-instanton coefficient $v_1=-1/8$; that derivation is
independent of the transport argument of
Subsections~\ref{app:second-degree}--\ref{app:critical-matching}, and the value
agrees with \cite[equation~(25)]{GYZ}. A leading estimate alone would not identify
this correction. Expanding the finite product \eqref{eq:PF12} gives
\[
 \pi C_g(0^p,2^m)=
 1+\frac{p(5-p)/15-4/9}{h}+O_p(h^{-2}).
\]
The critical harmonic coefficient also has a first correction. Its
summand law is geometric with mean $3/2$, variance $v=15/4$,
third cumulant $15$, and fourth cumulant $705/8$. At displacement
$d=-p/2$ from $3h/2$, expansion of its exact negative-binomial
coefficient by Stirling's formula gives
\[
 \frac{\mathcal M_{g,p}(0^p;y_0)}
 {\Gamma(h)4^{-g}A_*^{-h}/\sqrt{2\pi hv}}
 =1+\frac{-p^2/30+4p/15-19/180}{h}+O_p(h^{-2}),
 \qquad A_*=2/15.
\]
Explicitly, the coefficient before substituting $d=-p/2$ is
$-d^2/(2v)-15d/(2v^2)+(705/8)/(8v^2)-5(15)^2/(24v^3)$.
Combining these two expansions with \eqref{eq:PF13} proves
\[
 \lim_{Y\uparrow y_0}
 \frac{\pi\Phi_g^{(p)}(X(Y))}{\mathcal J_{g,p}(Y)}
 =1+\frac{t_p}{h}+O_p(h^{-2}).
\]
Take $Y_*=y_0-h^{-3}$. The positive homogeneous-polynomial comparison
\eqref{eq:PF15} and the analytic normalized factor in \eqref{eq:PF16}
change this ratio by $O_p(h(y_0-Y))=O_p(h^{-2})$ for $Y\ge Y_*$.
Moreover $\mathcal L_p(Y)=t_p+O_p(y_0-Y)$ there.
For $Y\le Y_*$, subtract the exact logarithmic derivative of
$\mathcal J_{g,p}$ in \eqref{eq:PF11} from \eqref{eq:ST8}, and integrate:
\[
 \log\frac{\pi\Phi_g^{(p)}(X(Y))}{\mathcal J_{g,p}(Y)}
 =\frac{\mathcal L_p(Y)}h+
 O_p\!\left(h^{-2}+\frac1{h^2}
             \log\frac{\mathcal Z(Y)}{\mathcal Z(Y_*)}\right).
\]
Here \eqref{eq:ST12} controls the terminal integral, and
$\mathcal Z(Y_*)\asymp h^{-15/2}$. The remainder is
$O_p(\log h/h^2)$ uniformly. Boundedness of $\mathcal L_p$ permits
exponentiation, proving \eqref{eq:ST14} in both regions.
\end{proof}

\subsection{The formal hierarchy}\label{sec:full-expansion}

We now construct the full series and control its remainder as the truncation order grows. These are still estimates for generating functions: the passage to a prescribed cusp count is made in Subsection~\ref{sec:canonical-expansion} and Appendix~\ref{app:extraction}. The recursion is triangular, in the sense described in Subsection~\ref{subsec:all-orders}. This subsection writes down the formal rules \eqref{eq:FE1}--\eqref{eq:FE6} and states Proposition~\ref{prop:hierarchy}, which bounds the coefficients and the actual remainders polynomially in the order; its inputs are the ratio estimates \eqref{eq:TL14}--\eqref{eq:TL27} of Section~\ref{ratios-with-a-fixed-number-of-distinguished-markings} and the kernel bounds of Section~\ref{the-positive-compact-family-of-kernels}. The proof divides each order-$j$ coefficient by $S_j^{2j}$, where $S_j=D+p+Mj$; after this division each rule becomes a coefficientwise inequality between series with nonnegative coefficients, an analytic majorant (Lemma~\ref{lem:majorant}) bounds the normalized coefficients by $CB^j$, and the true remainders are controlled by a second, linear majorant system with the same normalization. Matching at the pure-$\psi$ endpoint in Subsection~\ref{sec:endpoint-expansion} then completes the proof of Theorem~\ref{thm:full-expansion}.

Use a formal variable $t$ for inverse complexity, and put
\begin{equation}
 \mathsf P_p(\mathbf d;t)=\sum_{j\ge0}e_{p,j}(\mathbf d)t^j,
 \quad \mathsf C_p(h)=\frac hA\frac{W_{g,p}}{W_{g,p+1}},
 \quad \mathsf G_p(h)=\frac{W_{g-1,p+2}}{W_{g,p}}.
 \label{eq:FE1}
\end{equation}
The formal series for the last two quantities will also be denoted by
$\mathsf C_p(t),\mathsf G_p(t)$. Their constant terms, and $e_{p,0}$,
are one. Here $e_{p,1}=-\delta_p$ and $e_{p,2}=e_p$ from
\eqref{eq:PF5} and \eqref{eq:ST3}. All series operations below are
truncated at a fixed finite order; convergence of a formal series is
never asserted. A complexity change from $h$ to $h-c$ means replacing
$t$ by $t/(1-ct)$.

For a bounded component of genus $j$ containing $s$ other external
labels, set $r_*=2j+s$. Thus $r_*$ is the decrease in complexity
from the parent to the complementary component. Its all-zero ratio has the exact
factorization
\begin{equation}
\begin{split}
 \frac{W_{g-j,p-s}}{W_{g,p}}
 ={}&\prod_{a=0}^{r_*-1}
   \frac{A\mathsf C_{p-s+a}(h-r_*+a)}{h-r_*+a}
      \prod_{b=0}^{j-1}\mathsf G_{p+2b}(h),\\
 \mathsf R_{p;j,s}(t)
 ={}&(At)^{r_*}\prod_{a=0}^{r_*-1}
 \frac{\mathsf C_{p-s+a}(t/(1-(r_*-a)t))}{1-(r_*-a)t}
       \prod_{b=0}^{j-1}\mathsf G_{p+2b}(t).
\end{split}
 \label{eq:FE2}
\end{equation}
The genus steps have unchanged complexity; the cusp steps remove
$r_*$ markings at genus $g-j$. All their marking counts are positive
in the applications below. To obtain $\mathsf R$ through order $t^K$
requires scalar coefficients only through order $K-r_*$.

The coefficients are computed as follows. Use the child tuples of
\eqref{eq:ST2}, and set $\mathcal W_{j,q}(\mathbf d)=4^j\overline
U_{j,q}(\mathbf d)$ for the exact bounded-complexity amplitudes.
At each order $K\ge1$, form
\begin{equation}
\begin{split}
 \mathsf A_{p,i}={}&\frac{At}{1-t}\mathsf C_{p-1}\!\left(\frac t{1-t}\right)
 \sum_{j\ne i}(2d_j+1)
 \sum_{\ell\ge\max(0,1-d_i-d_j)}w_\ell
 \mathsf P_{p-1}\!\left(\mathbf d^{ij,\ell};\frac t{1-t}\right),\\
 \mathsf N_{p,i}={}&\frac{2At}{1-t}\mathsf C_{p+1}\!\left(\frac t{1-t}\right)
 \mathsf G_p(t)\sum_{\ell\ge\max(0,2-d_i)}w_\ell
 \sum_{u+v=d_i+\ell-2}\mathsf P_{p+1}\!\left(\mathbf d^{i;u,v};\frac t{1-t}\right).
\end{split}
 \label{eq:FE3}
\end{equation}
Nonnegative node degrees are understood, and $\mathsf A$ is absent
at $p=1$. The retained separating contribution is
\begin{equation}
 \mathsf S_{p,i}=
 \sum_{\substack{j\ge0,\ I\subset\{1,\ldots,p\}\setminus\{i\}\\
        s=|I|,\ (j,s+1)\text{ stable}\\2\le2j+s\le K}}
 \mathsf R_{p;j,s}(t)\sum_{\ell\ge0}w_\ell
 \sum_{u+v=d_i+\ell-2}
 \mathcal W_{j,s+1}(u,\mathbf d_I)
 \mathsf P_{p-s}\!\left(v,\mathbf d_J;\frac t{1-(2j+s)t}\right),
 \label{eq:FE4}
\end{equation}
where $J$ is the remaining set of labels. Only bounded $j$ occur;
for large $g$ these orient the smaller component uniquely. The two
orientations cancel the separating factor $1/2$. Define
\begin{equation}
 \begin{split}
 e_{p,K}(\mathbf0)&=0,\\
 e_{p,K}(\mathbf d+\mathbf e_i)-e_{p,K}(\mathbf d)
   &=-[t^K](\mathsf A_{p,i}+\mathsf N_{p,i}+\mathsf S_{p,i}),\\
 \mathsf C_p(t)&=\frac1A\sum_{d\ge1}b_d
       \mathsf P_{p+1}\!\left(d,0^p;\frac t{1+t}\right).
 \end{split}
 \label{eq:FE5}
\end{equation}
Finally, for $p\ge2$, KdV gives
\begin{equation}
\begin{split}
 \mathsf G_p(t)={}&\mathsf P_p(\mathbf e_1;t)
 -3\sum_{\substack{0\le s\le p-2,\ j\ge0\\
             (j,s+2)\text{ stable}\\1\le2j+s\le K}}
  \binom{p-2}s\mathcal W_{j,s+2}(\mathbf0)\mathsf R_{p;j,s}(t),\\
 \mathsf G_p(h)&=\mathsf G_{p+1}(h+1)
                   \mathsf C_{p+2}(h)/\mathsf C_p(h).
\end{split}
 \label{eq:FE6}
\end{equation}
The last line descends from $p=2$ to $p=1$, and optionally to $p=0$.
It is an exact identity; its shift $h+1$ must be retained.

The coefficient and remainder bounds below are obtained by comparing
coefficientwise inequalities with an analytic system. We isolate the
elementary principle used for this comparison; it is applied in the
proof of Lemma~\ref{lem:coefficient-bounds} and again, with the
exponential degree norm, in the proof of
Proposition~\ref{prop:gevrey-hierarchy} in Subsection~\ref{app:gevrey-one}.

\begin{lemma}[Analytic majorants]\label{lem:majorant}
Let $\Phi=(\Phi_1,\ldots,\Phi_m)$ be power series in
$(t,\mathsf v_1,\ldots,\mathsf v_m)$ with nonnegative coefficients, convergent near the
origin, such that $\Phi(0,0)=0$ and the matrix $\partial_v\Phi(0,0)$ is
strictly lower triangular. Then $U=\Phi(t,U)$ has a unique formal
solution with $U(0)=0$. Its coefficients are nonnegative, and it
converges near $t=0$. Let $e_{i,j}\in[0,\infty]$, for $1\le i\le m$ and
$j\ge1$, and put $E_i=\sum_{j\ge1}e_{i,j}t^j$. If
$e_{i,j}\le[t^j]\Phi_i(t,E)$ for all $i,j$, then
$e_{i,j}\le[t^j]U_i$. In particular every $e_{i,j}$ is finite, and
$e_{i,j}\le C\rho^{-j}$ with $C,\rho>0$ depending only on $\Phi$.
\end{lemma}
\begin{proof}
Order the pairs $(j,i)$ lexicographically. Since $E$ has no constant
term, $[t^j]\Phi_i(t,E)$ is a polynomial with nonnegative coefficients
in the numbers $e_{k,l}$. A monomial $t^av^\beta$ with $a\ge1$ or
$|\beta|\ge2$ contributes only coefficients with $l<j$, since each
factor has order at least one. The remaining terms are
$\partial_{v_k}\Phi_i(0,0)\,e_{k,j}$, and they occur only for $k<i$.
Thus $[t^j]\Phi_i(t,E)$ involves only pairs preceding $(j,i)$. The
equations $U=\Phi(t,U)$ therefore determine the coefficients of $U$
recursively and uniquely, and these are nonnegative. The same induction
and the monotonicity of polynomials with nonnegative coefficients give
$e_{i,j}\le[t^j]U_i$. The derivative of $(t,v)\mapsto v-\Phi(t,v)$
with respect to $v$ at the origin is $I-\partial_v\Phi(0,0)$, which is
invertible. The analytic implicit-function theorem gives a solution
holomorphic near $t=0$ and vanishing at $t=0$; its Taylor series is the
formal solution. Cauchy's estimate on a small circle gives the final
bound.
\end{proof}

\begin{proposition}[Uniform recursive expansions]\label{prop:hierarchy}
There are $C,B,M>0$ and one complex neighborhood $\Omega$ of
$[0,y_0]$ such that the preceding rules determine path-independent
holomorphic coefficients with
\[
 \|e_{p,j}(\mathbf d)\|_\Omega\le CB^j(D+p+Mj)^{2j},\qquad
 \|\mathsf C_{p,j}\|_\Omega+\|\mathsf G_{p,j}\|_\Omega
       \le CB^j(p+Mj)^{2j}.
\]
For every stable state $p\ge1$, $h=2g-2+p\ge1$, every $L\ge1$,
and every nonnegative degree tuple, uniformly for $0\le Y\le y_0$,
\begin{equation}
 \begin{split}
 \left|c_{g,p}(\mathbf d)-\sum_{j<L}e_{p,j}(\mathbf d)h^{-j}\right|
       &\le CB^L(D+p+ML)^{2L}h^{-L},\\
 \left|\mathsf C_p(h)-\sum_{j<L}\mathsf C_{p,j}h^{-j}\right|
 +\left|\mathsf G_p(h)-\sum_{j<L}\mathsf G_{p,j}h^{-j}\right|
       &\le CB^L(p+ML)^{2L}h^{-L}.
 \end{split}
 \label{eq:FE7}
\end{equation}
Here $D=\sum d_i$, $\mathsf C_{p,j}=[t^j]\mathsf C_p(t)$ and
$\mathsf G_{p,j}=[t^j]\mathsf G_p(t)$; invalid-degree profiles and invalid genus numerators
are zero. For fixed $p,j$, the rules use finitely many kernel moments,
low kernel coefficients and exact bounded-complexity amplitudes.
\end{proposition}

The proof of Proposition~\ref{prop:hierarchy} occupies
Subsections~\ref{app:weights}--\ref{app:remainder-bounds}. At each order the
construction treats the profile, the cusp ratio, the genus ratio at
$p\ge2$, and the descent to $p=1$, in this order; every profile on the right of
\eqref{eq:FE3}--\eqref{eq:FE4} has smaller order, and the dilaton and KdV
rules use the current profile only after it has been constructed.
Subsection~\ref{app:weights} fixes the analytic neighborhood and the
weights $S_j$ and records four elementary bounds.
Lemma~\ref{lem:coefficient-bounds} bounds the formal coefficients,
uniformly in the order, by a system of analytic majorants; this is the
first assertion of the proposition. Lemmas~\ref{lem:tails}
and~\ref{lem:first-order} control the omitted topological types and give
first-order estimates for the true quantities. Finally,
Subsection~\ref{app:remainder-bounds} proves the remainder bounds
\eqref{eq:FE7} by induction on the order, using a linear majorant system
for the normalized errors.

\subsection{Weights and elementary bounds}\label{app:weights}
We fix a common analytic neighborhood and the normalizing weights, and
record four elementary bounds used in both the coefficient and the
remainder estimates.

Shrink one complex neighborhood $\Omega$ of $[0,y_0]$ so that the
normalized kernels $\bar K_Y,\bar H_Y$ and their nonvanishing quotients
from Section~\ref{the-positive-compact-family-of-kernels} are holomorphic
and bounded on a common disk $|w|\le1+2\eta$, $\eta>0$, for every
$Y\in\Omega$, and so that $A(Y)$ is bounded above and away from zero on
$\Omega$. Since $\sum_\ell w_\ell w^\ell=(1-w)/\bar K_Y(w)$ and
$b_d=\frac13[w^{d-1}]\bar H_Y(w)$, Cauchy's estimate on the circle
$|w|=1+\eta$ gives constants $C_0$ and $\rho=1+\eta>1$ with
\begin{equation}
 |w_\ell(Y)|\le C_0\rho^{-\ell},\qquad |b_d(Y)|\le C_0\rho^{-d}
 \qquad(Y\in\Omega,\ \ell\ge0,\ d\ge1).
 \label{eq:FE-geom}
\end{equation}
The absolute-series argument in the proof of Lemma~\ref{lem:terminal}(iii) also works on $\Omega$:
the background majorant $B_Y$ of
Subsection~\ref{elementary-terminal-bounds-and-a-global-crude-estimate}
has a uniformly bounded germ, so one small radius $z_*$ makes
$B_Y(z_*)\le1/2$ for all $Y\in\Omega$. A stable component of complexity
$r$ has at most $r+2$ markings and deficit $m\le3r/2$, so the primitive
bound of Lemma~\ref{lem:terminal}(i) and this majorant give
\begin{equation}
 |\mathcal W_{j,q}(\mathbf d)|\le C_0^{\,r}\,\Gamma(r)\qquad(r=2j-2+q\ge1),
 \label{eq:FE-small}
\end{equation}
uniformly in the markings, the degrees and $Y\in\Omega$; enlarging
$C_0$ covers the components of complexity one.

Fix an integer $M\ge2$ with $\mathrm{e}^{2/M}<\rho$. For a profile at a
state with $p$ markings and degree sum $D$, and for a scalar with $p$
markings, the order-$j$ weights are
\begin{equation}
 S_j=D+p+Mj\qquad\text{and}\qquad S_j=p+Mj
 \label{eq:FE-weights}
\end{equation}
respectively, so that $S_j\ge Mj\ge j$ and $S_j\ge p$. The following four
bounds are used repeatedly.

\begin{enumerate}[label=(\alph*),leftmargin=*]
\item\emph{Child weights.} Let $k\le j$. The merger child
$\mathbf d^{ij',\ell}$ has $p-1$ markings and degree sum $D+\ell-1$, the
nonseparating child $\mathbf d^{i;u,v}$ has $p+1$ markings and degree
sum $D+\ell-2$, and the large child $(v,\mathbf d_J)$ of a separating
term has $p-s$ markings and degree sum at most $D+\ell-2$. In all three
cases the order-$k$ weight of the child is at most $S_j+\ell$. A scalar
factor occurring at loss $r$ has at most $p+r$ markings; if it enters at
order $k$ with $k+r\le j$, its weight is at most $p+r+Mk\le p+Mj\le S_j$,
because $r\le j-k\le M(j-k)$. The dilaton child $(d,0^p)$ of a scalar
with $p$ markings has weight at most $S_j+d+1$, and the children in the
KdV and descent rules have weights at most $S_j+2$.
\item\emph{Moment bound.} For $0\le k\le2j$ and $a\ge0$,
$(S_j+a)^k\le S_j^k(1+a/S_j)^{2j}\le S_j^k\mathrm{e}^{2a/M}$, since
$S_j\ge Mj$. Hence, with
$C_1=C_0\sum_{n\ge0}\rho^{-n}\mathrm{e}^{2(n+1)/M}<\infty$,
\begin{equation}
 \sum_{\ell\ge0}|w_\ell|(S_j+\ell+1)^k\le C_1S_j^k,\qquad
 \sum_{d\ge1}|b_d|(S_j+d+1)^k\le C_1S_j^k\qquad(0\le k\le2j).
 \label{eq:FE-moment}
\end{equation}
\item\emph{Shifts.} For integers $a\ge1$, $q\ge0$ and $0\le c\le c_0$,
\[
 [t^q](1-ct)^{-a}=\frac{a(a+1)\cdots(a+q-1)c^q}{q!}
 \le\frac{((a+q)c_0)^q}{q!},
\]
with the empty product equal to one. For a product
$\prod_k(1-c_kt)^{-a_k}$, taking absolute values if necessary, the
coefficient is bounded by that of $(1-c_0t)^{-a}$, where
$a=\sum_k a_k$ and $c_0\ge\max_k|c_k|$. In an order-$j$ term, a
prefactor of order $r$ and coefficients of total order $i$ leave shift
degree $q=j-r-i$. The total denominator exponent is at most $r+i$,
so $a+q\le j$; also $c_0\le\max(1,r)\le j$. Thus the coefficient is
at most $j^{2q}/q!\le S_j^{2q}/q!$. The case $a=0$ contributes only
$q=0$. After normalization all shifts are therefore majorized by
$\mathrm{e}^{t}$.
\item\emph{Counts.} At a profile state, $\sum_{j'\ne i}(2d_{j'}+1)\le
2D+p\le2S_j$, the number of pairs $(u,v)$ with $u+v=d_i+\ell-2$ is at
most $D+\ell\le S_j+\ell$, and a coordinate path from $\mathbf0$ to
$\mathbf d$ has $D\le S_j$ steps. In a separating term of loss
$r=2j_0+s$ the number of label subsets is $\binom{p-1}s\le p^s/s!$,
there are at most $r+1\le2^r$ pairs $(j_0,s)$, and the small component
has complexity $r-1$ in the profile rule and $r$ in the KdV rule, so
\eqref{eq:FE-small} gives
\begin{equation}
 \binom{p-1}{s}\,|\mathcal W|\le C_0^{\,r}\frac{p^s\Gamma(r)}{s!}
 \le C_0^{\,r}(p+r)^r\le C_0^{\,r}S_j^{\,r}.
 \label{eq:FE-count}
\end{equation}
Indeed $\Gamma(r)/s!\le r^{r-1-s}$ when $s\le r-1$, and $\Gamma(r)/r!\le1$
when $s=r$; and $p+r\le p+Mj\le S_j$ because $r\le j$.
\end{enumerate}

\subsection{Bounds for the formal coefficients}\label{app:coefficient-bounds}

We bound the formal coefficients uniformly in the order. After division
by the weights, each coefficient rule becomes a coefficientwise
inequality between power series with nonnegative coefficients, and the
normalized coefficient sequences are dominated by the solution of a convergent
analytic system.

\begin{lemma}[Coefficient bounds]\label{lem:coefficient-bounds}
The rules \eqref{eq:FE3}--\eqref{eq:FE6} determine coefficients
$e_{p,j}(\mathbf d)$, $\mathsf C_{p,j}$ and $\mathsf G_{p,j}$ that are
holomorphic on $\Omega$, and there are $C,B_0$ such that, on $\Omega$, for
all $p\ge1$, all degree tuples and all $j\ge0$,
\[
 |e_{p,j}(\mathbf d)|\le CB_0^{\,j}S_j^{2j},\qquad
 |\mathsf C_{p,j}|+|\mathsf G_{p,j}|\le CB_0^{\,j}S_j^{2j},
\]
with the weights \eqref{eq:FE-weights}.
\end{lemma}

\begin{proof}
Define the normalized coefficient sequences
\[
 \varepsilon_j=\sup\frac{|e_{p,j}(\mathbf d)|}{S_j^{2j}},\quad
 \gamma_j=\sup\frac{|\mathsf C_{p,j}|}{S_j^{2j}},\quad
 \eta_j=\sup_{p\ge2}\frac{|\mathsf G_{p,j}|}{S_j^{2j}},\quad
 \zeta_j=\sup\frac{|\mathsf G_{1,j}|}{S_j^{2j}},
\]
the suprema being taken over $Y\in\Omega$, all marking counts $p\ge1$ and
all degree tuples; $\varepsilon_0=\gamma_0=\eta_0=\zeta_0=1$. Write
$\varepsilon(t)=\sum_j\varepsilon_jt^j$, and similarly $\gamma,\eta,\zeta$,
and put $\bar\eta=\max(\eta,\zeta)$ coefficientwise, which bounds every
genus factor. Each coefficient of order $j$ is a finite combination of
coefficients of lower order, summed against $w_\ell$ or $b_d$; by
\eqref{eq:FE-geom} and the polynomial growth in $\ell$ or $d$ of the
lower-order weights these sums converge absolutely and uniformly on
$\Omega$, so each coefficient is holomorphic on $\Omega$ and each
supremum is finite, by induction on $j$. We bound the order-$j$
coefficients term by term.

\textbf{Profile.} Fix a state $(\mathbf d,p)$ and a direction $i$. By
\eqref{eq:FE3}, $[t^j]\mathsf A_{p,i}$ is the sum over $j'\ne i$,
$\ell\ge0$ and $i_1+i_2+q=j-1$ of
\[
 A(2d_{j'}+1)\,w_\ell\;\mathsf C_{p-1,i_1}\;
 e_{p-1,i_2}(\mathbf d^{ij',\ell})\;[t^q](1-t)^{-(1+i_1+i_2)}.
\]
By (a) of Subsection~\ref{app:weights} the two coefficients are at most
$\gamma_{i_1}S_j^{2i_1}$ and $\varepsilon_{i_2}(S_j+\ell)^{2i_2}$; by
(b) the sum over $\ell$ is at most $C_1S_j^{2(i_1+i_2)}$; by (c) the
shift coefficient is at most $S_j^{2q}/q!$; and by (d) the merger count is
at most $2S_j$. Hence
\[
 |[t^j]\mathsf A_{p,i}|\le2|A|C_1\,S_j^{\,2j-1}\,
 [t^{j-1}]\bigl(\mathrm{e}^{t}\gamma(t)\varepsilon(t)\bigr).
\]
The term $\mathsf N_{p,i}$ has in addition one genus coefficient and the
pair count $S_j+\ell$, which \eqref{eq:FE-moment} absorbs with
$k=2(i_1+i_2+i_3)+1\le2j$; it is at most
$2|A|C_1S_j^{2j-1}[t^{j-1}](\mathrm{e}^{t}\gamma\bar\eta\varepsilon)$.
A separating term of loss $r\ge2$ in \eqref{eq:FE4} carries the
prefactor $(At)^r$, $r$ cusp coefficients, $j_0\le r/2$ genus
coefficients, one profile coefficient, the count \eqref{eq:FE-count}, the
pair count and the shifts. Since $\gamma$ and $\bar\eta$ have constant
term one and nonnegative coefficients,
$\gamma^r\bar\eta^{j_0}\preceq\gamma^{2r}\bar\eta^{2r}$, and the term is
at most
$(2|A|C_0)^rC_1S_j^{\,r+1}S_j^{\,2(j-r)}[t^{j-r}](\mathrm{e}^t\gamma^{2r}\bar\eta^{2r}\varepsilon)$,
where $S_j^{r+1}\le S_j^{2r-1}$ because $r\ge2$. Summing the increments
along a coordinate path of at most $S_j$ steps supplies one more factor
$S_j$. After division by $S_j^{2j}$ we obtain, for a constant $b_1$
depending only on $\sup_\Omega|A|$, $C_0$ and $C_1$,
\begin{equation}
 \varepsilon_j\le\Bigl[b_1\mathrm{e}^{t}\varepsilon(t)
 \Bigl(t\gamma\bar\eta+\sum_{r\ge2}(b_1t)^r\gamma^{2r}\bar\eta^{2r}\Bigr)\Bigr]_j
 \qquad(j\ge1).
 \label{eq:FE-maj-E}
\end{equation}

\textbf{Cusp ratio.} By \eqref{eq:FE5}, $[t^j]\mathsf C_p$ is the sum over
$d\ge1$ and $i+q=j$ of $A^{-1}b_d\,e_{p+1,i}(d,0^p)\,[t^q](1+t)^{-i}$.
The child weight is at most $S_j+d+1$ by (a), so \eqref{eq:FE-moment}
and (c) give
\begin{equation}
 \gamma_j\le\bigl[b_2\,\mathrm{e}^{t}\varepsilon(t)\bigr]_j\qquad(j\ge1),
 \qquad b_2=C_1/\inf_\Omega|A|.
 \label{eq:FE-maj-C}
\end{equation}

\textbf{Genus ratio at $p\ge2$.} In \eqref{eq:FE6} the term
$\mathsf P_p(\mathbf e_1;t)$ has weight at most $S_j+1$ and contributes
at most $\mathrm{e}^{2/M}\varepsilon_j$. A separating term of loss $r\ge1$
carries $3\binom{p-2}s|\mathcal W_{j_0,s+2}(\mathbf0)|\le3C_0^rS_j^r\le3C_0^rS_j^{2r}$
by \eqref{eq:FE-count}, the prefactor $(At)^r$, at most $2r$ scalar
coefficients and the shifts. Hence
\begin{equation}
 \eta_j\le\Bigl[b_3\varepsilon(t)+b_3\mathrm{e}^{t}
 \sum_{r\ge1}(b_3t)^r\gamma^{2r}\bar\eta^{2r}\Bigr]_j\qquad(j\ge1),
 \label{eq:FE-maj-H}
\end{equation}
for a constant $b_3$ depending only on $M$, $\sup_\Omega|A|$ and $C_0$.

\textbf{Descent.} The identity
$\mathsf G_1(h)=\mathsf G_2(h+1)\mathsf C_3(h)/\mathsf C_1(h)$ is a
product of a genus series with the shift $t/(1+t)$, a cusp series and
the formal inverse of a cusp series; the three weights are at most
$S_j+2$. If $|\mathsf C_{1,i}|\le\gamma_iS_j^{2i}$ for $i\le j$, then the
coefficients of $1/\mathsf C_1$ through order $j$ are bounded by
$S_j^{2i}$ times those of $1/(2-\gamma(t))=\sum_{k\ge0}(\gamma-1)^k$.
Therefore
\begin{equation}
 \zeta_j\le\Bigl[b_4\,\mathrm{e}^{t}\frac{\eta(t)\gamma(t)}{2-\gamma(t)}\Bigr]_j
 \qquad(j\ge1),\qquad b_4=\mathrm{e}^{12/M}.
 \label{eq:FE-maj-G}
\end{equation}

\textbf{The majorant system.} Put $b=\max(1,b_1,b_2,b_3,b_4)$ and let
$E_*,C_*,H_*,G_*$ be the formal solution, with constant terms one, of
\begin{equation}
\begin{aligned}
 E_*&=1+b\,\mathrm{e}^{t}E_*\Bigl(tC_*G_*+\frac{(btC_*^2G_*^2)^2}{1-btC_*^2G_*^2}\Bigr),&
 C_*&=1+b(\mathrm{e}^{t}E_*-1),\\
 H_*&=1+b(E_*-1)+b\,\mathrm{e}^{t}\frac{btC_*^2G_*^2}{1-btC_*^2G_*^2},&
 G_*&=1+b\Bigl(\mathrm{e}^{t}\frac{H_*C_*}{2-C_*}-1\Bigr).
\end{aligned}
\label{eq:FE-system}
\end{equation}
Here $E_*$ bounds the profile coefficients after division by their
weights, $C_*$ the cusp coefficients, $H_*$ the genus coefficients at
$p\ge2$ and $G_*$ the descent. Since $b\ge1$, $C_*/(2-C_*)\succeq1$ and
$\mathrm{e}^t\succeq1$, the last equation gives $G_*\succeq H_*$; so
$G_*$ majorizes $\bar\eta$ as soon as $H_*,G_*$ majorize $\eta,\zeta$.
Writing $E_*=1+\mathsf v_1$, $C_*=1+\mathsf v_2$, $H_*=1+\mathsf v_3$, $G_*=1+\mathsf v_4$, the
right-hand sides of \eqref{eq:FE-system} are power series in $(t,v)$ with
nonnegative coefficients, convergent near the origin and vanishing there.
At $t=0$ the first equation has no linear term, the second and third are
linear in $\mathsf v_1$ only, and the fourth depends only on $\mathsf v_2,\mathsf v_3$; so the
Jacobian at the origin is strictly lower triangular, and
Lemma~\ref{lem:majorant} shows that the formal solution of
\eqref{eq:FE-system} converges near $t=0$. The inequalities
\eqref{eq:FE-maj-E}--\eqref{eq:FE-maj-G} now give
$\varepsilon_j\le[E_*]_j$, $\gamma_j\le[C_*]_j$, $\eta_j\le[H_*]_j$ and
$\zeta_j\le[G_*]_j$ for all $j$, by induction on $j$ in the stage order:
at order $j$ the right-hand sides involve $\bar\eta$ only through orders
below $j$, or through $\eta_j$ in the descent, and
$\max(\eta_i,\zeta_i)\le[G_*]_i$ as soon as $\eta_i\le[H_*]_i$ and
$\zeta_i\le[G_*]_i$. Cauchy's estimate for the convergent solution gives
$\varepsilon_j,\gamma_j,\eta_j,\zeta_j\le CB_0^{\,j}$ for all $j$, with
$C,B_0$ depending only on the system.
\end{proof}

\subsection{Tails over topological types and first-order estimates}\label{app:tails}

We next bound the tails of the sums over topological types, before
truncating any expansion, and deduce first-order estimates for the
true quantities. In the sums $S_{g,p}$ and $T_{g,p}$ of
\eqref{eq:TL18}, orient the smaller-complexity component; its loss is
its complexity plus one in $S_{g,p}$ and its complexity in $T_{g,p}$. Write
$S_{g,p}^{[L]}$ and $T_{g,p}^{[L]}$ for the parts of these sums in which
the loss is at least $L$; for $p=1$, $S_{g,1}$ is the one-marking
separating sum of \eqref{eq:TL26}.

\begin{lemma}[Tails over topological types]\label{lem:tails}
Fix $M_0>0$. There are $C_{M_0},B_{M_0}$ such that, for all sufficiently
large $h$ and all $1\le p\le M_0\sqrt h$,
\begin{equation}
 \begin{aligned}
 S_{g,p}^{[L]}&\le C_{M_0}\left(\frac{B_{M_0}(p+L)}h\right)^L,
                  &&2\le L\le\sqrt h,\\
 T_{g,p}^{[L]}&\le C_{M_0}\left(\frac{B_{M_0}(p+L)}h\right)^L,
                  &&1\le L\le\sqrt h,
 \end{aligned}
 \label{eq:FE13}
\end{equation}
the second bound being used only for $p\ge2$. The ranges of $L$ start at
the stability minima two and one respectively.
\end{lemma}

\begin{proof}
Note first that, before specializing to fixed markings,
Lemmas~\ref{lem:lower} and~\ref{lem:upper}, with zero count and marked count
both $p$, give, using the common local limit \eqref{eq:H17},
\[
 \log W_{g,p}=\log\Gamma(h)-h\log A
 +O_{M_0}(\sqrt h\log(h+1)+\log^2(h+1)).
\]
Indeed the harmonic probability has displacement $-p/2=O(\sqrt h)$
and lies between $c_{M_0}/\sqrt h$ and one; the terminal factors cost
$O(p\log(p+1)+(p+\log(h+1))\log(h+1))$.
For losses between $\varepsilon h$ and half the parent complexity,
Stirling entropy dominates this error. Both children have comparable
complexity and $O(\sqrt h)$ markings. The $2^p$ subsets cost only
$\mathrm{e}^{O(\sqrt h)}$, so their total is $O(\mathrm{e}^{-c_\varepsilon h})$.
For $r\le\varepsilon h$, \eqref{eq:TL16}--\eqref{eq:TL17} give the
complementary ratio at most $(C/h)^r$. If its marking count is one,
insert $2j_0+1$ cusps before raising genus $j_0$ times, then the
remaining $s-1$ cusps; this uses $p\ge2$ and avoids a one-marking
genus comparison. The same subset/Gamma bound gives
$a_r=[C(p+r)/h]^r$, after absorbing the number of types into $C^r$.
Their consecutive ratio satisfies
\[
 \frac{a_{r+1}}{a_r}
 \le \frac{eC(p+r+1)}h
 \le eC\left(\frac{M_0}{\sqrt h}+\varepsilon+\frac1h\right).
\]
Choose $\varepsilon$ first, and then a lower bound on $h$ depending on
$M_0$, so this is at most $1/2$. Thus the sum over $L\le r\le\varepsilon h$
is at most $2a_L$. The macroscopic tail also fits \eqref{eq:FE13}:
for $1\le L\le\sqrt h$,
\[
 \left(\frac{p+L}{h}\right)^L
 \ge \exp\{-\sqrt h\log h\},
\]
which eventually exceeds $\mathrm{e}^{-c_\varepsilon h}$. Enlarging the constants
covers the bounded remaining complexities and proves the asserted bound.
At $p=1$, compare to $T_{g,2}$: two uses of \eqref{eq:TL16} and the
proved fixed-one-marking cusp ratio bound multiply its summand by
at most $C/j_0$ for $j_0\le g/2$. The loss is $2j_0$ in both sums.
This proves \eqref{eq:FE13}, including $p=1$.
\end{proof}

\begin{lemma}[First-order estimates]\label{lem:first-order}
Fix $M_0>0$. For all sufficiently large $h$, all $1\le p\le M_0\sqrt h$
and all degree tuples,
\[
 0\le1-c_{g,p}(\mathbf d)
 \le\min\{1,C_{M_0}D(D+p+1)/h\},\qquad
 |\mathsf C_p(h)-1|+|\mathsf G_p(h)-1|\le C_{M_0}(p+1)/h.
\]
\end{lemma}

\begin{proof}
Subtract the recursion at adjacent degrees, as in the proof of
\eqref{eq:TL23}, and bound its separating part by Lemma~\ref{lem:tails} at
the stability minima. For the profile, each increment costs
$C[(D+p+1)/h+(D+1)(p+2)^2/h^2]$; sum the path.
Dilaton sums the child bound against absolute $b_d$ moments, KdV uses
degree one and $T^{[1]}$, and the exact descent gives $p=1$.
These estimates concern the true quantities, in a range of marking counts
that grows with $h$.
\end{proof}

\subsection{Remainder bounds: proof of Proposition~\ref{prop:hierarchy}}\label{app:remainder-bounds}

We now compare the true quantities with their truncated formal series.
With the same normalization as for the coefficients, the remainders
satisfy a linear majorant system whose kernels have no constant term,
and induction on the order bounds them geometrically.

Before the proof we record that the coefficient rules, which contain
infinite sums over the kernel index $\ell$, over node degrees and over
dilaton degrees, reduce at each order to finitely many kernel data, and
that the coefficients grow polynomially in the degrees. For an integer
$b\ge0$, a \emph{stratum of threshold $b$} in $\mathbb N^p$ is the set
of tuples obtained by prescribing, for each coordinate, either a value
in $\{0,\ldots,b-1\}$ or membership in $[b,\infty)$; the coordinates of
the second kind are called \emph{free}. A function on $\mathbb N^p$ is
\emph{polynomial on the strata of threshold $b$} if its restriction to
each stratum is a polynomial in the free coordinates. The path in
\eqref{eq:FE5} is taken to be the one that raises $d_1$ from $0$ to its
value, then $d_2$, and so on; path independence is proved in
Proposition~\ref{prop:hierarchy} below.

\begin{lemma}[Polynomial stratification]\label{rev:lem:stratification}
Define $b_0=0$ and $b_K=2b_{K-1}+2K+4$ for $K\ge1$. For every $K\ge0$
and $p\ge1$:
\begin{enumerate}[label=\textup{(\roman*)},leftmargin=*]
\item the coefficient $\mathbf d\mapsto e_{p,K}(\mathbf d)$ is
polynomial on the strata of threshold $b_K$, of total degree at most
$2K$ in the free coordinates;
\item for $K\ge1$ and each direction $i$, the adjacent coefficient
$\mathbf d\mapsto[t^K](\mathsf A_{p,i}+\mathsf N_{p,i}+\mathsf S_{p,i})(\mathbf d)$
is polynomial on the strata of threshold $b_K$, of total degree at most
$2K-1$;
\item the coefficients of the polynomials in \textup{(i)} and
\textup{(ii)}, and the scalar coefficients $\mathsf C_{p,K}$,
$\mathsf G_{p,K}$, are finite combinations of the kernel moments
$\sum_{\ell\ge\ell_0}w_\ell\ell^m$ and $\sum_{d\ge d_0}b_dd^m$ with
$\ell_0,d_0,m\le C(K,p)$, of finitely many $w_\ell$ and $b_d$, of exact
amplitudes of complexity at most $K$, of $A$, and of scalar
coefficients of order less than $K$ (for the descent rule in
\eqref{eq:FE6}, also of the current-order scalar coefficients
$\mathsf C_{1,K},\mathsf C_{3,K},\mathsf G_{2,K}$, which precede
$\mathsf G_{1,K}$ in the stage order).
\end{enumerate}
In particular, for fixed $p$ and $K$ the rules
\eqref{eq:FE3}--\eqref{eq:FE6} use finitely many kernel data, and
$|e_{p,K}(\mathbf d)|\le C_{p,K}(D+1)^{2K}$.
\end{lemma}

\begin{proof}
Induction on $K$. For $K=0$ all coefficients equal one. Let $K\ge1$ and
assume (i)--(iii) for all orders below $K$. Write $b=b_K$ and
$b'=b_{K-1}$. The recursion gives
\begin{equation}\label{rev:eq:strat-threshold}
 b\ge2b'+1,\qquad b\ge b'+\tfrac32K,\qquad b\ge\tfrac32K-1,
 \qquad b\ge b_k\ \ (k<K),
\end{equation}
the first three because $b-2b'=2K+4$ and $b'\ge0$, the last because
$(b_k)$ is increasing. Fix a stratum $\Sigma$ of threshold $b$ and a
direction $i$; a coordinate is \emph{low} if it is prescribed on
$\Sigma$ and free otherwise. Every free coordinate is at least
$b\ge b_k$ for all $k<K$, and every low coordinate has a fixed value.
Hence, for a child tuple obtained from $\mathbf d\in\Sigma$ by one of the
rules, the stratum of threshold $b_k$ in which it lies is determined by
$\Sigma$ as soon as the modified coordinates are controlled; we check
these coordinates rule by rule.

\emph{Merger.} In $\mathsf A_{p,i}$ the child $\mathbf d^{ij,\ell}$
replaces $d_i,d_j$ by $d_i+d_j+\ell-1$ and keeps the other coordinates,
and the child coefficient has order $k\le K-1$. If $d_i$ or $d_j$ is
free, then $d_i+d_j+\ell-1\ge b-1\ge b_k$ for every $\ell\ge0$ by the
first inequality of \eqref{rev:eq:strat-threshold}, the lower summation
limit $\max(0,1-d_i-d_j)$ is zero, and by induction
$e_{p-1,k}(\mathbf d^{ij,\ell})$ is a polynomial of degree at most $2k$
in $d_i+d_j+\ell-1$ and in the other free coordinates; the sum over
$\ell$ against $w_\ell$ converges absolutely by \eqref{eq:FE-geom} and
yields a polynomial in $d_i,d_j$ and the free coordinates whose
coefficients are the moments $\sum_{\ell\ge0}w_\ell\ell^m$, $m\le2k$.
The factor $2d_j+1$ raises the degree to at most $2k+1\le2K-1$. If
$d_i$ and $d_j$ are both low, $d_i+d_j$ is a fixed value; the finitely
many $\ell$ with $d_i+d_j+\ell-1<b_k$ give fixed child tuples, and the
remaining $\ell$ give the polynomial form with the tail moments
$\sum_{\ell\ge\ell_0}w_\ell\ell^m$.

\emph{Nonseparating pairs.} In $\mathsf N_{p,i}$ the child
$\mathbf d^{i;u,v}$ replaces $d_i$ by the pair $(u,v)$ with
$u+v=d_i+\ell-2$. Suppose $d_i$ is free. Then $d_i\ge b\ge2b_k+1$, so
for every $\ell\ge0$ the pairs fall into three ranges: $u<b_k$, in which
case $v=d_i+\ell-2-u\ge(2b_k+1)-2-(b_k-1)=b_k$; $v<b_k$, in which case
$u\ge b_k$; and $u,v\ge b_k$. The first two ranges consist of $b_k$
fixed low values each and do not overlap, since $u,v\le b_k-1$ would
force $u+v\le2b_k-2<d_i+\ell-2$. On each of them the child coefficient
is, by induction, a polynomial of degree at most $2k$ in the free
coordinate $v$ (respectively $u$) and the other free coordinates. On
the third range it is a polynomial $P$ of degree at most $2k$ in
$(u,v)$, and
\[
 \sum_{u=b_k}^{d_i+\ell-2-b_k}P(u,\,d_i+\ell-2-u)
\]
is, by Faulhaber's formula, a polynomial of degree at most $2k+1$ in
$d_i+\ell-2$; the summation range is consistent with this polynomial
form because $d_i+\ell-2\ge2b_k-1$. Summation against $w_\ell$ over
$\ell\ge0$ gives a polynomial of degree at most $2k+1\le2K-1$ in $d_i$
and the other free coordinates, with moment coefficients. If $d_i$ is
low, the finitely many $\ell$ with $d_i+\ell-2<2b_k+1$ give finite sums
of polynomials of degree at most $2K-1$ in the other free coordinates,
and the remaining $\ell$ are treated as before with
tail moments.

\emph{Retained splits.} In $\mathsf S_{p,i}$ the bounded component has
genus $j$, label set $I$ with $|I|=s$, $2\le2j+s\le K$, node degree $u$
and amplitude $\mathcal W_{j,s+1}(u,\mathbf d_I)$, which vanishes unless
$u+|\mathbf d_I|\le3j+s-2\le\tfrac32K-2$. By the third inequality of
\eqref{rev:eq:strat-threshold} a free coordinate cannot belong to $I$,
so the admissible sets $I$ and the values $\mathbf d_I$ are fixed on
$\Sigma$, and $u$ ranges over at most $\tfrac32K-1$ values. The large
child $(v,\mathbf d_J)$ has $v=d_i+\ell-2-u$. If $d_i$ is free, the
second inequality of \eqref{rev:eq:strat-threshold} gives
$v\ge b-2-(\tfrac32K-2)\ge b_k$ for all $\ell\ge0$, so the child
coefficient is a polynomial of degree at most $2k$ in $v$ and in the
free coordinates of $\mathbf d_J$, and summation against $w_\ell$ gives
a polynomial of degree at most $2k\le2K-2$ in $d_i$. If $d_i$ is low,
finitely many $\ell$ are exceptional and the rest are polynomial with
tail moments. The prefactor $\mathsf R_{p;j,s}$ involves $A$ and scalar
coefficients of order at most $K-2j-s<K$ only.

\emph{Scalars and shifts.} The scalar and shift coefficients in
\eqref{eq:FE3}--\eqref{eq:FE4} do not depend on the degrees. This
proves (ii), with degree at most $2K-1$.

\emph{Path sum.} By \eqref{eq:FE5}, $e_{p,K}(\mathbf d)$ is minus the
sum of the adjacent coefficients along the chosen path. While $d_i$ is
raised, the states visited are $(d_1,\ldots,d_{i-1},m,0,\ldots,0)$,
$0\le m<d_i$, and their strata of threshold $b$ are determined by
$\Sigma$ and by whether $m<b$. If $d_i$ is free, the sum over $m$ is,
for $m<b$, a finite sum of polynomials of degree at most $2K-1$ in the
free coordinates among $d_1,\ldots,d_{i-1}$, plus, by (ii) and Faulhaber's
formula, a polynomial of degree at most $2K$ in $d_i$ and in the free
coordinates among $d_1,\ldots,d_{i-1}$ for $b\le m<d_i$; if $d_i$ is low
it is a finite sum. Hence $e_{p,K}$ is polynomial on the strata of
threshold $b$, of total degree at most $2K$, which is (i).

\emph{Dilaton, KdV and descent.} In \eqref{eq:FE5}, $\mathsf C_{p,K}$
is a finite combination of the values $e_{p+1,k}(d,0^p)$, $k\le K$,
against $b_d$ and shift coefficients. By (i) at orders $k\le K$,
$d\mapsto e_{p+1,k}(d,0^p)$ is a polynomial of degree at most $2k$ for
$d\ge b_k$ and takes finitely many values for $d<b_k$, so the sum over
$d$ uses the signed moments $\sum_{d\ge b_k}b_dd^m$, $m\le2k$, and
finitely many $b_d$. In \eqref{eq:FE6}, $\mathsf P_p(\mathbf e_1;t)$ is
evaluated at a single tuple, the amplitudes $\mathcal W_{j,s+2}(\mathbf0)$
have complexity $2j+s\le K$, and the descent is a finite product of
scalar series. This proves (iii). The final assertion follows from (i),
since there are finitely many strata and each free coordinate is at
most $D$.
\end{proof}

\begin{proof}[Proof of Proposition~\ref{prop:hierarchy}]
The coefficient bounds are Lemma~\ref{lem:coefficient-bounds}, with
$S_j^{2j}\le(D+p+Mj)^{2j}$. It remains to prove the remainder bounds
\eqref{eq:FE7}, the path independence and the finiteness of the
coefficient rules.

For $L\ge1$ define the normalized remainders
\[
 \mathscr E_L=\sup\frac{\bigl|c_{g,p}(\mathbf d)-\sum_{j<L}e_{p,j}(\mathbf d)h^{-j}\bigr|}{S_L^{2L}h^{-L}},
 \qquad
 \mathscr C_L=\sup\frac{\bigl|\mathsf C_p(h)-\sum_{j<L}\mathsf C_{p,j}h^{-j}\bigr|}{S_L^{2L}h^{-L}},
\]
and $\mathscr H_L,\mathscr G_L$ in the same way for $\mathsf G_p$ with
$p\ge2$ and with $p=1$; the suprema are over real $Y\in[0,y_0]$, all
stable states with $h\ge1$, $p\ge1$, and all degree tuples, invalid
tuples having profile zero. We prove by induction on $L$, in the stage
order profile, cusp, genus, descent, that these numbers are finite and
at most $CB^L$; this is \eqref{eq:FE7}.

\textbf{Large weights.} Put $\delta=S_L^2/h$ and fix $\delta_0\in(0,1]$
with $B_0\delta_0<1/2$. The exact quantities are bounded: $0\le c\le1$,
$|\mathsf C_p(h)|\le C$ by \eqref{eq:TL16}, $\mathsf G_p(h)\le1$ for
$p\ge2$ by \eqref{eq:TL17}, and $\mathsf G_1(h)\le C$ by the case $p=1$ of \eqref{eq:TL27}.
If $\delta\ge\delta_0$, then Lemma~\ref{lem:coefficient-bounds} and $S_j^{2j}h^{-j}\le\delta^j$
($j\le L$) bound the normalized remainder by
\[
 \delta^{-L}\Bigl(C+\sum_{j<L}CB_0^{\,j}\delta^{j}\Bigr)
 \le C\delta_0^{-L}\Bigl(1+\sum_{j<L}(B_0\delta_0)^j\Bigr)\le3C\delta_0^{-L},
\]
since $\delta^{j-L}\le\delta_0^{j-L}$ for $j<L$. Invalid-degree parents
satisfy $D\ge h$, hence $\delta\ge h\ge1$, and lie in this region. In the remaining region
$\delta<\delta_0$ one has $p+ML\le S_L<\sqrt h$, so every marking count
that occurs, including those of children, is at most $p+L+2\le4\sqrt h$:
Lemmas~\ref{lem:tails} and~\ref{lem:first-order} apply with $M_0=4$. Also $h>M^2L^2/\delta_0$, and
shrinking $\delta_0$ once makes Lemma~\ref{lem:first-order} give $\mathsf C_p(h)\ge1/2$
wherever a cusp ratio is inverted. From now on $\delta<\delta_0$.

\textbf{Products of exact factors.} Every term of the exact adjacent
difference, of the exact dilaton identity, of the exact KdV identity and
of the exact descent is a product of exact factors: a prefactor
$A^r\prod(h-c_k)^{-1}$, cusp and genus ratios at shifted complexities,
and profiles of children. Each exact factor $F$ has, for every $q\ge1$,
a decomposition $F=T_{<q}F^{\mathrm{form}}+R_F(q)$ into its truncated
formal series and a remainder, and the formal rule at order $L$ is the
truncation of the product of the formal series. The product of exact
factors and the truncated formal product differ by the finite identity
\begin{equation}
 \prod_kF_k-T_{<q}\prod_kF_k^{\rm form}
 =R_{F_1}(q)\prod_{k>1}F_k
 +\sum_{j<q}f_{1,j}h^{-j}\Bigl(\prod_{k>1}F_k-T_{<q-j}\prod_{k>1}F_k^{\rm form}\Bigr),
 \label{eq:FE-product}
\end{equation}
applied again to the bracket on the right, and so on. Each resulting term
contains exactly one remainder $R_{F_k}(q_k)$, multiplied by known
coefficients of the factors before it and by exact factors after it. For
an inverse, with known formal inverse coefficients $d_j$,
\begin{equation}
 \mathsf C\,T_{<q}(1/\mathsf C)-1
 =\sum_{j<q}d_jh^{-j}\bigl(\mathsf C-T_{<q-j}\mathsf C\bigr),
 \label{eq:FE-inverse}
\end{equation}
so the remainder of $1/\mathsf C$ is $1/\mathsf C$ times a sum of known
coefficients against remainders of $\mathsf C$.

\textbf{Bounding one term.} Normalize by $S_L^{2L}h^{-L}$. A known
coefficient of order $j'$ of a factor with weight $S'$ is at most
$CB_0^{j'}(S')^{2j'}$ by Lemma~\ref{lem:coefficient-bounds}, where $S'\le S_L+\ell$, $S_L+d+1$,
$S_L+2$ or $S_L$ by (a) of Subsection~\ref{app:weights}, applied with $j'\le L$. A remainder
$R_{F_k}(q_k)$ of a factor at complexity $h-c$, $|c|\le L$, is at most
$\mathscr X_{q_k}(S')^{2q_k}(h-c)^{-q_k}$, where $\mathscr X$ is the error
sequence of its type, and
$(h-c)^{-q_k}\le h^{-q_k}\mathrm{e}^{2|c|q_k/h}\le\mathrm{e}^{2\delta_0/M^2}h^{-q_k}$
because $q_k\le L$ and $L^2/h\le\delta_0/M^2$. The exact factors after
the remainder are bounded by the constants above. The known shift tails
are bounded as follows. If the orders of the known coefficients and of
the prefactor sum to $a=L-q$ and the first omitted shift degree is $q$,
the elementary inequality
\[
 \sum_{n\ge q}\frac{(a)_n}{n!}x^n\le\frac{(a)_q}{q!}x^q(1-x)^{-a-q}
 \qquad(a\ge1,\ 0\le x<1),
\]
which follows from $(a)_{q+k}=(a)_q(a+q)_k$ and $(q+k)!\ge q!k!$, bounds
the tail by $\binom{L-1}{q}(c/h)^q(1-c/h)^{-L}$ with $c\le L$; after
division by $S_L^{2q}h^{-q}$ this is at most
$\mathrm{e}^{2\delta_0/M^2}/q!$, since $c(L-1)\le L^2\le S_L^2$
and $(1-c/h)^{-L}\le\mathrm e^{2cL/h}$ in the small-weight region.
The counts, kernel moments and path sums are
handled exactly as in the proof of Lemma~\ref{lem:coefficient-bounds}, with the same powers of $S_L$: a linear
term supplies $S_L^2$ against the prefactor order one, and a separating
term of loss $r$ supplies $S_L^{r+2}\le S_L^{2r}$ ($r\ge2$) in the profile
rule and $S_L^{r}\le S_L^{2r}$ ($r\ge1$) in the KdV rule.

\textbf{Forcing terms.} Three contributions do not involve the error
sequences. The exact topological types of loss $r\ge L$, which the formal
rule truncated below order $L$ omits, have total mass bounded by
\eqref{eq:FE13} at $L$; after the kernel sum and the coordinate path their profile
contribution is at most
$CD(D+1)(B(p+L)/h)^{L}\le CB^LS_L^{L+2}h^{-L}\le CB^LS_L^{2L}h^{-L}$ for
$L\ge2$, and the KdV contribution has no path or pair factor. The
difference between the exact kernel sum over $\ell\le m$ and the formal
sum over all $\ell$, and the macroscopic cutoffs in the proof of Lemma~\ref{lem:tails}, are
$O(\mathrm{e}^{-ch})$, and
$\mathrm{e}^{-ch}\le C^LL!\,h^{-L}\le C^LS_L^{2L}h^{-L}$. At $L=1$ the
normalized remainders are bounded by Lemma~\ref{lem:first-order}.

\textbf{The linear majorant system.} Let $F_0=C/(1-B_0t)$ be the majorant
of the normalized known coefficients from Lemma~\ref{lem:coefficient-bounds}, and
$V=C_{\rm ex}+F_0$, where $C_{\rm ex}\ge1$ bounds all exact factors and
$1/\mathsf C$ in the region $\delta<\delta_0$. Collecting the bounds
above, the error sequences $\mathscr E,\mathscr C,\mathscr H,\mathscr G$,
regarded as formal series in $t$, satisfy coefficientwise
\begin{equation}
\begin{aligned}
 \mathscr E&\preceq F+K_E\mathscr R,&
 \mathscr C&\preceq F+b\mathscr E,\\
 \mathscr H&\preceq F+b\mathscr E+K_H\mathscr R,&
 \mathscr G&\preceq F+b(\mathscr H+\mathscr C)+K_G\mathscr R,
 \qquad \mathscr R=\mathscr E+\mathscr C+\mathscr H+\mathscr G,
\end{aligned}
\label{eq:FE-error-system}
\end{equation}
with the zero-constant kernels
\[
 K_E=b\,\mathrm{e}^{bt}\Bigl[tV^2+\sum_{r\ge2}(2r+1)(bt)^rV^{2r}\Bigr],\quad
 K_H=b\,\mathrm{e}^{bt}\sum_{r\ge1}2r(bt)^rV^{2r},\quad
 K_G=b\bigl[\mathrm{e}^{bt}V^3-V(0)^3\bigr],
\]
and the forcing
\[
 F=\frac{bV^3\mathrm{e}^{bt}}{(1-btV^2)^3}+\frac{b}{1-\delta_0^{-1}t},
\]
for a constant $b$ depending only on the constants already fixed. The
factors $2r+1$ and $2r$ count the possible positions of the single
remainder among the factors of a product with $r$ cusp, at most $r$ genus
and at most one profile factor; the term $b\mathscr E$ in the second and
third lines is the current-order remainder of the profile, which enters
the dilaton and KdV rules without a prefactor $t$, and $b(\mathscr H+\mathscr C)$
plays the same role in the descent. The first summand of $F$ dominates
the shift tails, the omitted topology and the cutoff terms, each of which
is a product of at most $2r+1$ known or exact factors times $C^q/q!$ or
$C^LL!/S_L^{2L}$; the second summand dominates the region
$\delta\ge\delta_0$.

\textbf{Induction.} Substituting the first line of
\eqref{eq:FE-error-system} into the others gives
$\mathscr R\preceq A_0F+A_0(K_E+K_H+K_G)\mathscr R$ for a constant
$A_0$. Since $K_E+K_H+K_G$ has no constant term, the coefficient of
$t^L$ on the right involves $\mathscr R_j$ only for $j<L$. Let
$\mathscr R_*=A_0F/(1-A_0(K_E+K_H+K_G))$, a series with nonnegative
coefficients, analytic near $t=0$. If $\mathscr R_j\le[\mathscr R_*]_j$
for all $j<L$, then
$\mathscr R_L\le[A_0F+A_0(K_E+K_H+K_G)\mathscr R_*]_L=[\mathscr R_*]_L$
by monotonicity of series with nonnegative coefficients; the case $L=1$
is the base estimate. Cauchy's estimate for $\mathscr R_*$ gives
$\mathscr R_L\le CB^L$ and proves \eqref{eq:FE7}; no infinite
inverse-complexity series is evaluated. For each fixed degree tuple, two
coordinate paths give asymptotic expansions of the same ratio, and
uniqueness of asymptotic coefficients as $g\to\infty$ proves the asserted
path independence.

\textbf{Finiteness of the coefficient rules.} The finite coefficient
rule and the polynomial degree bounds used above are the content of
Lemma~\ref{rev:lem:stratification}: at each order the degree lattice
admits a finite stratification by the thresholds $b_K$ (an admissible,
not optimal, choice), the coefficients are polynomial on every stratum,
the adjacent coefficient has degree at most $2K-1$ and its path sum
degree at most $2K$, and every displayed infinite sum in the coefficient
algorithm reduces to finitely many kernel moments, low kernel
coefficients and exact small amplitudes, the latter being finite
polynomials in the $\theta_\ell$. The common analytic disks and
removable critical normalization give analyticity of that data
throughout $[0,y_0]$.
\end{proof}

\subsection{Exponential degree norms and single-factorial bounds}
\label{app:gevrey-one}

This subsection and the next prove
Proposition~\ref{prop:gevrey-hierarchy}, the single-factorial form of
Proposition~\ref{prop:hierarchy}; it is the version used in the proof of
Theorem~\ref{thm:full-expansion} and, together with that theorem, in
the coefficient bounds of Section~\ref{sec:absolute}. The polynomial bounds of
Proposition~\ref{prop:hierarchy} retain useful information at every
degree. For summing the dilaton equation, however, an exponential degree
norm gives a sharper dependence on the expansion order: in this norm the
order-$j$ bounds grow like $(p+Mj)^j$ rather than $(p+Mj)^{2j}$, with a
factor $\mathrm{e}^{\sigma D}$ for profiles. The argument has three
parts. Lemma~\ref{lem:linear-tails} enlarges the range of the
topological-tail estimate of Lemma~\ref{lem:tails} to linearly many
markings and a linear order range, using the all-marking envelope
\eqref{eq:gevrey-envelope}. Lemma~\ref{rev:lem:path-sum} isolates the
mechanism that saves one factorial per order. The coefficient bounds \eqref{eq:gevrey-coefficients} then
follow from the majorant system \eqref{eq:gevrey-majorant} and
Lemma~\ref{lem:majorant}; the actual remainders
\eqref{eq:gevrey-remainders} are proved in
Subsection~\ref{app:gevrey-remainders}.
All constants in this subsection are uniform in the real parameter
$Y\in[0,y_0]$ unless a complex neighborhood is specified.

\emph{Dependence on earlier results.} The chain of dependencies is
linear. The marked ratio estimates of
Section~\ref{ratios-with-a-fixed-number-of-distinguished-markings} feed
Proposition~\ref{prop:prefactor} in Section~\ref{sec:prefactor}, whose
endpoint constant is fixed by \eqref{eq:PF12a}, that is, by the
Itzykson--Zuber identity and the large-order asymptotics of
Painlev\'e~I in the form recorded by
Garoufalidis--Its--Kapaev--Mari\~no~\cite[Section~2.3]{GIKM}. Of
Subsection~\ref{sec:endpoint-expansion} it uses only the leading case
of Lemma~\ref{rev:lem:endpoint-gevrey}(iii), a Stirling computation that
rests on these two inputs alone. The leading asymptotic of
Proposition~\ref{prop:prefactor} gives the envelope
\eqref{eq:gevrey-envelope} used in Lemma~\ref{lem:linear-tails};
Proposition~\ref{prop:gevrey-hierarchy} uses
Lemma~\ref{lem:linear-tails} and the formal rules
\eqref{eq:FE3}--\eqref{eq:FE6}; and Theorem~\ref{thm:full-expansion}
follows from Proposition~\ref{prop:gevrey-hierarchy} together with
Subsection~\ref{sec:endpoint-expansion}, which has two independent
parts. The endpoint Gevrey expansion
(Lemma~\ref{rev:lem:endpoint-gevrey}) uses only the Painlev\'e~I input
and none of the results of this subsection or the next. The vanishing
\eqref{eq:FE9} of the higher cusp coefficients at $y_0$ uses
\eqref{eq:gevrey-remainders} and the exact string--dilaton identity at
the endpoint, but not the Painlev\'e input. Neither part is used in the
proof of Proposition~\ref{prop:gevrey-hierarchy}. Thus
Section~\ref{ratios-with-a-fixed-number-of-distinguished-markings}
$\to$ Section~\ref{sec:prefactor} $\to$ Lemma~\ref{lem:linear-tails}
$\to$ Proposition~\ref{prop:gevrey-hierarchy} $\to$
Theorem~\ref{thm:full-expansion} contains no cycle.

\begin{lemma}[Linear marking and order ranges]\label{lem:linear-tails}
There are constants $\eta,c\in(0,1)$ and $C,B>0$ with the following
properties, uniformly in $Y\in[0,y_0]$. For every stable $(g,p)$ with
$p\ge1$ and $h=2g-2+p$,
\begin{equation}\label{eq:gevrey-envelope}
 W_{g,p}\le C\,\Gamma(h)A^{-h}/\sqrt h .
\end{equation}
If moreover $p\le\eta h$, the bounds \eqref{eq:FE13} hold, with the
constants $C,B$, for all $L\le ch$ above the same stability minima
(the bound for $T_{g,p}^{[L]}$ when $p\ge2$), and, for all degree
tuples $\mathbf d$ with $D=|\mathbf d|$,
\[
 |\mathsf C_p(h)-1|+|\mathsf G_p(h)-1|\le C(p+1)/h,
 \qquad
 0\le1-c_{g,p}(\mathbf d)
 \le\min\{1,CD(D+p+1)/h\}.
\]
The envelope \eqref{eq:gevrey-envelope} uses only
Proposition~\ref{prop:prefactor} at $p\in\{1,2\}$, the exact
inequalities \eqref{eq:TL16}--\eqref{eq:TL17} and \eqref{eq:TL27} at
$p=1$; no expansion in a growing number of markings enters.
\end{lemma}

\begin{proof}
For $p=1,2$, the leading asymptotic of
Proposition~\ref{prop:prefactor}, after the normalization
\eqref{eq:TL9}, gives
\[
 W_{g,p}\asymp\Gamma(h)A^{-h}/\sqrt h.
\]
Indeed its leading factor is
$\Gamma(h)A^{-h}\sqrt{A\Theta}/(\pi\sqrt{\pi h})$; all powers
of $r$ cancel. For each fixed $(g,p)$ the normalized amplitude $W_{g,p}(Y)$ and
$A(Y)$ are continuous and positive on the closed interval $[0,y_0]$,
the former being a polynomial with nonnegative coefficients in the
kernel coefficients $\theta_\ell(Y)$; hence the two-sided bound, which
Proposition~\ref{prop:prefactor} gives uniformly on $[0,y_0)$ for
large $g$, extends to $Y=y_0$, and the finitely many remaining stable
pairs are covered by compactness. Genus exchange
\eqref{eq:TL17} gives $W_{g,p}\le W_{g+1,p-2}$ for $p\ge4$.
Iterating it reduces even $p$ to two and odd $p$ to three. In the
odd case use the bounded one-marking genus ratio \eqref{eq:TL27}
once to reach $p=1$. Each other exchange has constant one. This
proves \eqref{eq:gevrey-envelope}.
The complexity-one bases are $W_{0,3}=1$ and
$W_{1,1}=\theta_1/2$. Conversely, for $g\ge2$, iteration of
\eqref{eq:TL16} from marking count one gives
\[
 W_{g,p}\ge C_0^{-(p-1)}
       \frac{\Gamma(h)}{\Gamma(2g-1)}W_{g,1}
 \ge c_0\mathrm{e}^{-C_1p}\Gamma(h)A^{-h}/\sqrt h.
\]
The residual power $A^{p-1}$ is absorbed into $\mathrm{e}^{-C_1p}$.
These estimates use only the leading fixed-marking result.

Orient the smaller component as in Lemma~\ref{lem:tails}. For loss
$r\le\varepsilon h$, its complementary ratio is bounded by
$(C/h)^r$ using \eqref{eq:TL16}--\eqref{eq:TL17}; the one-marking
complement is handled by the same cusp-insertion order as there.
The small-component and subset bound \eqref{eq:FE-count} gives
$a_r=[C(p+r)/h]^r$, with
$a_{r+1}/a_r\le\mathrm{e}\,C(p+r+1)/h$. Choose $\varepsilon$, then $\eta$,
so that this is at most $1/2$.

For $\varepsilon h\le r\le h/2+O(1)$ both child complexities are
proportional to $h$. The upper and lower envelopes above and Stirling's
formula give a bound
$h^{C_2}\exp\{-hH(\varepsilon/2)+C_3p\}$, where
$H(t)=-t\log t-(1-t)\log(1-t)$.
The label subsets cost at most $2^p$, and the number of remaining
types is polynomial in $h$. Shrinking $\eta$ therefore bounds this
whole sector by $C\mathrm{e}^{-a h}$ for a fixed $a>0$.
Choose $c<\min(\varepsilon,\mathrm{e}^{-1})$ with $c\log(1/c)<a/2$.
For $L\le ch$,
\[
 ((p+L)/h)^L\ge(L/h)^L\ge \mathrm{e}^{-c\log(1/c)h},
\]
so the macroscopic tail fits the same bound as $\sum_{r\ge L}a_r$.
The comparison with $T_{g,2}$ in Lemma~\ref{lem:tails} also treats
$S_{g,1}$. Enlarging constants includes bounded complexities.

The adjacent-degree argument of Lemma~\ref{lem:first-order} now costs
$C[(D+p+1)/h+(D+1)(p+2)^2/h^2]$ per increment. Put
$x=D(D+p+1)/h$. When $x\ge1$, the profile bound follows from
$0\le c_{g,p}\le1$. When $x<1$ and $D\ge1$, one has
$D(p+2)/h\le2$.
Consequently the second path contribution is bounded by the first:
\[
 \frac{D(D+1)(p+2)^2}{h^2}
 \le C\frac{D(D+p+1)}h.
\]
This proves the profile estimate in the linear marking range;
substituting $p\le\eta h$ in each increment would only give
$O(D^2p/h)$.
Dilaton sums this bound against its geometric kernel, KdV uses degree
one and the loss-one tail, and descent gives the remaining scalar case.
\end{proof}

\begin{proposition}[Single-factorial recursive bounds]
\label{prop:gevrey-hierarchy}
Let $\rho>1$ be as in \eqref{eq:FE-geom}, and fix
$0<\sigma<\log\rho$. There are $C,B,M>0$ and a complex
neighborhood $\Omega$ of $[0,y_0]$ such that
\begin{equation}\label{eq:gevrey-coefficients}
 \begin{split}
 \|e_{p,j}(\mathbf d)\|_\Omega
    &\le CB^j(p+Mj)^j\mathrm{e}^{\sigma D},\\
 \|\mathsf C_{p,j}\|_\Omega+\|\mathsf G_{p,j}\|_\Omega
    &\le CB^j(p+Mj)^j.
 \end{split}
\end{equation}
For every stable state with $p\ge1$, $h\ge1$, and every $L\ge1$,
uniformly for real $Y\in[0,y_0]$,
\begin{equation}\label{eq:gevrey-remainders}
 \begin{split}
 \left|c_{g,p}(\mathbf d)-\sum_{j<L}e_{p,j}(\mathbf d)h^{-j}\right|
    &\le CB^L(p+ML)^Lh^{-L}\mathrm{e}^{\sigma D},\\
 \left|\mathsf C_p(h)-\sum_{j<L}\mathsf C_{p,j}h^{-j}\right|
 +\left|\mathsf G_p(h)-\sum_{j<L}\mathsf G_{p,j}h^{-j}\right|
    &\le CB^L(p+ML)^Lh^{-L}.
 \end{split}
\end{equation}
Invalid-degree profiles and invalid genus numerators are zero, as in
Proposition~\ref{prop:hierarchy}.
\end{proposition}

We prove the coefficient estimates here and the actual remainders in
Subsection~\ref{app:gevrey-remainders}. The latter are needed to pass
from formal factorial bounds to an exponentially accurate approximation.

\emph{Weights and normalized suprema.} Fix $0<\sigma_0<\sigma_1<\log\rho$
and an integer $M\ge\max\{2,2(\sigma_1-\sigma_0)\}$. For
$0<a\le\sigma_1-\sigma_0$, use the degree exponent $\sigma_0+a$ and the
weights
\begin{equation}\label{rev:eq:gevrey-weights}
 Z_j=Z_j(p,a)=p+\frac{Mj}a,\qquad B_j=Z_j^{\,j}\ (j\ge1),\qquad B_0=1,
\end{equation}
so that $Z_j\ge Mj/a\ge2j$, $Z_j\ge p$ and $j/Z_j\le a/M\le1/2$.
An order-$j$ profile coefficient at a state with $p$ markings and degree
sum $D$ is normalized by $B_j\mathrm{e}^{(\sigma_0+a)D}$, a scalar
coefficient by $B_j$. Let $\varepsilon_j$, $\gamma_j$, $\eta_j$ and
$\zeta_j$ be the suprema of the normalized profile, cusp, genus ($p\ge2$)
and descent coefficients of order $j$, taken over all $p$, all admissible
$a$, all degree tuples and all $Y\in\Omega$; put
$\varepsilon_0=\gamma_0=\eta_0=\zeta_0=1$ and
$\bar\eta_j=\max(\eta_j,\zeta_j)$. At each fixed $p,a$, the polynomial
estimates of Lemma~\ref{lem:coefficient-bounds} ensure convergence of
the kernel sums and finiteness of each coefficient. Uniform finiteness
in $p$ will follow by induction; it is not a consequence of the
polynomial $p^{2j}$ bound. The single-factorial mechanism is isolated in
the following lemma: the exponential degree norm absorbs the sum along a
coordinate path at the cost of one power of $Z_j$, and in separating
terms the dimension of the bounded component removes the path count
altogether. Concretely, a child of positive order $0<k<j$ is measured at the
smaller exponent $a'=ak/j$; since $Mk/a'=Mj/a$, its weight is still
comparable to $Z_j$, while the exponential gap $\mathrm e^{-(a-a')D}$
pays for a path of length $D$ at the price of one factor $Z_j$. Over
$j$ orders this costs $Z_j^j\approx(Mj/a)^j$, a single factorial.

\begin{lemma}[Path-sum estimate]\label{rev:lem:path-sum}
Let $j\ge1$, let $C_0,\rho$ be as in \eqref{eq:FE-geom}, and let $C$,
$C_1$ denote constants depending only on $\sigma_0,\sigma_1,\rho,C_0$
and $\sup_\Omega|A|$.
\begin{enumerate}[label=\textup{(\roman*)},leftmargin=*]
\item\textup{(Linear terms.)} Consider a term of
$[t^j](\mathsf A_{p,i}+\mathsf N_{p,i})$ in \eqref{eq:FE3} whose child
profile coefficient has order $k$, whose scalar coefficients have total
order $i$ (one cusp coefficient in $\mathsf A$; one cusp and one genus
coefficient in $\mathsf N$) and whose shift has degree $q$, so that
$i+k+q=j-1$. Bound the child profile coefficient in the norm with
exponent $a'=ak/j$, so that $Mk/a'=Mj/a$ (for $k=0$ the child profile
is the constant one and no norm is needed), and write $\Gamma_i$ for the
product of the normalized scalar suprema involved ($\gamma_i$, or
$\gamma_{i_1}\bar\eta_{i_2}$ with $i_1+i_2=i$). At a parent state with
$p$ markings and degree sum $m$ the term is bounded by
\begin{equation}\label{rev:eq:path-increment}
 C(2m+p+1)\;\varepsilon_kZ_j^{\,k}\mathrm{e}^{a/M}\cdot
 \Gamma_iZ_j^{\,i}\cdot\frac{(a/M)^qZ_j^{\,q}}{q!}\;
 \mathrm{e}^{(\sigma_0+a')m}.
\end{equation}
Summing \eqref{rev:eq:path-increment} along a coordinate path from
$\mathbf0$ to $\mathbf d$, with $D=|\mathbf d|$, gives at most
\begin{equation}\label{rev:eq:path-sum}
 C(D+p+1)\,\mathrm{e}^{(\sigma_0+a')D}\;
 \varepsilon_k\Gamma_i\,\mathrm{e}^{a/M}\frac{(a/M)^q}{q!}\,Z_j^{\,j-1}.
\end{equation}
Division by the parent normalization
$Z_j^{\,j}\mathrm{e}^{(\sigma_0+a)D}$ leaves the factor
$Z_j^{-1}(D+p+1)\mathrm{e}^{-(a-a')D}$, and, since $a-a'=a(j-k)/j$,
\begin{equation}\label{rev:eq:path-normalized}
 \sup_{D\ge0}\,(D+p+1)\,\mathrm{e}^{-a(j-k)D/j}
 \le p+1+\frac{j}{\mathrm{e}\,a(j-k)}\le3Z_j .
\end{equation}
Hence the normalized contribution of all linear terms of order $j$ is at
most
\[
 C\mathrm{e}^{a/M}\sum_{i+k+q=j-1}\varepsilon_k\,\Gamma_i\,\frac{(a/M)^q}{q!},
\]
the power count being $1+i+k+q=j$.
\item\textup{(Separating terms.)} In a term of \eqref{eq:FE4} of loss
$r=2u+s\ge2$, the bounded component has complexity $r-1$ and dimension
\begin{equation}\label{rev:eq:path-dimension}
 3u+s-2=u+r-2\le\tfrac32r-2,
\end{equation}
so its node degree takes at most $2r+1$ values, independently of $D$ and
$\ell$. Its amplitude, the label subsets and the node choice contribute
at most
\begin{equation}\label{rev:eq:path-separating-count}
 C^r\binom{p-1}s(2r+1)\Gamma(r-1)\le C_1^{\,r}(p+r)^r\le C_1^{\,r}Z_j^{\,r}.
\end{equation}
The degree path contributes the geometric sum
$\sum_{m<D}\mathrm{e}^{(\sigma_0+a)m}\le C\mathrm{e}^{(\sigma_0+a)D}$, with
no factor $D$, so the large child may be bounded at the parent exponent
$a$. A shift of degree $q$ whose total denominator exponent is at most
$r+i+k$, each shift being at most $r$, costs at most
\begin{equation}\label{rev:eq:path-shift}
 \frac{(jr)^q}{q!}\le\Bigl(\frac aM\Bigr)^{q}Z_j^{\,q}\,\mathrm{e}^{r}.
\end{equation}
Every scalar factor of order $v$ in such a term has marking count at
most $p+r$ and $r+v\le j$, hence weight at most $Z_j^{\,v}$.
\end{enumerate}
\end{lemma}

\begin{proof}
(i) The term consists of the prefactor $A$, which carries the power $t$,
the count $\sum_{j'\ne i}(2d_{j'}+1)\le2m+p$ for a merger or the number
$\#\{(u,v):u+v=d_i+\ell-2\}\le m+\ell+1$ of nonseparating pairs, the
kernel coefficient $w_\ell$, the scalar coefficients, the child profile
coefficient and the shift coefficient. The child has $p\mp1$ markings and
degree sum at most $m+\ell$. In the norm with exponent $a'=ak/j$ its
weight is $Z_k(p\mp1,a')=p\mp1+Mj/a\le Z_j+1$, and
\[
 (Z_j+1)^k\le Z_j^{\,k}\Bigl(1+\frac1{Z_j}\Bigr)^{j}
 \le Z_j^{\,k}\mathrm{e}^{j/Z_j}\le Z_j^{\,k}\mathrm{e}^{a/M};
\]
so the child profile coefficient is at most
$\varepsilon_kZ_j^{\,k}\mathrm{e}^{a/M}\mathrm{e}^{(\sigma_0+a')(m+\ell)}$.
A scalar coefficient of order $i'$ at a marking count $p'\le p+1$ is at
most its normalized supremum times $Z_{i'}(p',a)^{i'}\le(Z_j+1)^{i'}$,
and the same inequality bounds the product of the scalar weights by
$\mathrm{e}^{a/M}Z_j^{\,i}$; the constant $C$ absorbs this second factor
$\mathrm{e}^{a/M}\le\mathrm{e}^{1/2}$. Every factor of \eqref{eq:FE3} other than
$\mathsf G_p(t)$ carries the shift $t/(1-t)$, so the total denominator
exponent is at most $1+i+k\le j$; by (c) of
Subsection~\ref{app:weights} with $c_0=1$, and since
$[t^q](1-t)^{-a}$ is nondecreasing in the exponent, the shift
coefficient is at most
$(1+i+k+q)^q/q!\le j^q/q!\le(a/M)^qZ_j^{\,q}/q!$, because
$j\le(a/M)Z_j$. Finally, since $\sigma_0+a'\le\sigma_1<\log\rho$,
\eqref{eq:FE-geom} gives
$\sum_{\ell\ge0}|w_\ell|(m+\ell+1)\mathrm{e}^{(\sigma_0+a')\ell}
\le C_0\sum_{\ell\ge0}(m+\ell+1)\rho^{-\ell}\mathrm{e}^{\sigma_1\ell}\le C(m+1)$.
Collecting these bounds proves \eqref{rev:eq:path-increment}.

A coordinate path from $\mathbf0$ to $\mathbf d$ visits each total
degree $m=0,\ldots,D-1$ exactly once, and
\[
 \sum_{m=0}^{D-1}(2m+p+1)\mathrm{e}^{(\sigma_0+a')m}
 \le(2D+p+1)\frac{\mathrm{e}^{(\sigma_0+a')D}}{\mathrm{e}^{\sigma_0+a'}-1}
 \le\frac{2(D+p+1)}{\mathrm{e}^{\sigma_0}-1}\,\mathrm{e}^{(\sigma_0+a')D},
\]
which gives \eqref{rev:eq:path-sum}. Division by
$Z_j^{\,j}\mathrm{e}^{(\sigma_0+a)D}$ leaves
$Z_j^{-1}(D+p+1)\mathrm{e}^{-(a-a')D}$ with $a-a'=a(j-k)/j>0$.
Since $\sup_{D\ge0}D\mathrm{e}^{-\alpha D}=1/(\mathrm{e}\alpha)$ for $\alpha>0$,
the first inequality of \eqref{rev:eq:path-normalized} follows; for the
second, $p+1\le Z_j+1\le2Z_j$ and
$j/(a(j-k))\le j/a\le Z_j/M\le Z_j/2$.

(ii) The bounded component has genus $u$ and $s+1$ markings, hence
complexity $2u-2+s+1=r-1$ and dimension $3u-3+s+1=3u+s-2$; since
$2u+s=r$ and $u\le r/2$, this gives \eqref{rev:eq:path-dimension}. Its
amplitude vanishes when the sum of its degrees exceeds its dimension, so its node
degree is at most $\tfrac32r-2$ and takes at most $\tfrac32r-1\le2r+1$
values; the other node degree $v$ is then fixed by $u_0+v=d_i+\ell-2$.
By \eqref{eq:FE-small} the amplitude is at most $C_0^{\,r-1}\Gamma(r-1)$,
the number of label subsets is $\binom{p-1}s\le p^s/s!$, and
$\Gamma(r-1)/s!\le r^{\,r-2-s}$ for $s\le r-2$, while $s=r$ gives
$\Gamma(r-1)/r!\le1$; the case $s=r-1$ does not occur since $r-s=2u$ is
even. Hence $\binom{p-1}s(2r+1)\Gamma(r-1)\le3^r(p+r)^r$, which gives
\eqref{rev:eq:path-separating-count} with $C_1=3C$, and
$p+r\le p+2j\le Z_j$ because $r\le j$ and $M/a\ge2$. None of these factors, nor the kernel sum
$\sum_\ell|w_\ell|\mathrm{e}^{(\sigma_0+a)\ell}\le C$, grows with $m$, so the
increments of a separating term sum along the path to at most a constant
times $\sum_{m<D}\mathrm{e}^{(\sigma_0+a)m}\le\mathrm{e}^{(\sigma_0+a)D}/(\mathrm{e}^{\sigma_0}-1)$,
and after normalization by $\mathrm{e}^{(\sigma_0+a)D}$ nothing depending
on $D$ remains; in particular the large child, of degree sum at most
$m+\ell$, is bounded at the parent exponent $a$ with weight
$Z_k(p-s,a)\le Z_j$. For the shift, (c) of
Subsection~\ref{app:weights} with total exponent at most $r+i+k$ and
$c_0\le r$ gives the coefficient bound $((r+i+k+q)r)^q/q!\le(jr)^q/q!$,
because $r+i+k+q\le j$; then $j^q\le(a/M)^qZ_j^{\,q}$ and
$r^q/q!\le\mathrm{e}^r$ prove \eqref{rev:eq:path-shift}. Finally a scalar
factor of order $v$ at marking count at most $p+r$ has weight
$p+r+Mv/a\le p+M(r+v)/a\le p+Mj/a=Z_j$, since $M/a\ge2$ and $r+v\le j$.
\end{proof}

\begin{proof}[Proof of \eqref{eq:gevrey-coefficients}]
Use the weights \eqref{rev:eq:gevrey-weights} and the normalized suprema
defined before Lemma~\ref{rev:lem:path-sum}. By
Lemma~\ref{rev:lem:path-sum}(i), the linear terms of order $j$
contribute to $\varepsilon_j$ at most
\[
 C\mathrm{e}^{a/M}\sum_{i+k+q=j-1}\varepsilon_k\,\Gamma_i\,\frac{(a/M)^q}{q!};
\]
one factor $Z_j$ is spent on the prefactor $t$ together with the path sum,
by \eqref{rev:eq:path-normalized}, and the remaining $Z_j^{\,i+k+q}$ on
the coefficients and the shift. The exponent $a'=ak/j$ is used only in
the child profile, and $a'=0$ when $k=0$.

By Lemma~\ref{rev:lem:path-sum}(ii), a separating term of loss $r\ge2$
contributes at most $C_1^{\,r}$ times the normalized product of its $r$
cusp coefficients, of its at most $r/2$ genus coefficients and of one
profile coefficient, times $(a/M)^q\mathrm{e}^r$ for its shift; the
degree path is a geometric sum without an extra factor $D$, and the
power $Z_j^{\,r}$ against the prefactor $(At)^r$ is supplied by
\eqref{rev:eq:path-separating-count}.
The KdV split starts at $r=1$, with bounded-component complexity $r$;
the bound $\binom{p-2}s\Gamma(r)\le C^r(p+r)^r$ applies.

Dilaton is a bounded functional because
$\sum_d|b_d|\mathrm{e}^{\sigma_1d}<\infty$. The current profile in KdV is
evaluated at degree one. Descent changes the marking count by at most
two and costs at most $\mathrm{e}^{2a/M}$ in the weights.

Here is an explicit majorant. Let $Q(t)=(1-bt)^{-1}$, with $b\ge1$
large enough for all the preceding constants, and give $E,C,H,G$
constant terms one. Set
\begin{equation}\label{eq:gevrey-majorant}
\begin{aligned}
 E&=1+bQE\left(tCG+\frac{(btC^2G^2)^2}{1-btC^2G^2}\right),\\
 C&=1+b(QE-1),\\
 H&=1+b(E-1)+bQ\frac{btC^2G^2}{1-btC^2G^2},\\
 G&=1+b\left(Q\frac{HC}{2-C}-1\right).
\end{aligned}
\end{equation}
They bound respectively profiles, cusp ratios, genus ratios at $p\ge2$,
and descent; the last equation also gives $G\succeq H$.
A loss-$r$ term has $r$ cusp factors and at most $r/2$ genus factors,
so the powers in \eqref{eq:gevrey-majorant} dominate their products.
$Q$ bounds the normalized shifts. With the constant terms subtracted,
the Jacobian at the origin is
\[
 \begin{pmatrix}0&0&0&0\\ b&0&0&0\\ b&0&0&0\\0&2b&b&0\end{pmatrix}.
\]
It is strictly lower triangular. Lemma~\ref{lem:majorant} gives an
analytic solution with nonnegative coefficients bounded by $CB^j$.
Induction in order and then in the stage order $E,C,H,G$ proves the
same bound for the normalized suprema. Taking $\sigma_1=\sigma$ and
one fixed positive $a$, and absorbing $1/a$ into $M$, proves
\eqref{eq:gevrey-coefficients}.
\end{proof}

\subsection{Actual remainders in the exponential norm}
\label{app:gevrey-remainders}

This subsection proves the remainder bounds \eqref{eq:gevrey-remainders}
and thereby completes the proof of
Proposition~\ref{prop:gevrey-hierarchy}; the coefficient bounds
\eqref{eq:gevrey-coefficients} were proved in
Subsection~\ref{app:gevrey-one}. The argument has the same four stages as
the remainder proof for polynomial weights in
Subsection~\ref{app:remainder-bounds}, and the proof below is marked
accordingly: a region of large weights, where the exact bounds and the
coefficient bounds suffice; retained exact products, decomposed by the
finite identities \eqref{eq:FE-product}--\eqref{eq:FE-inverse} into terms
containing one unknown remainder each, with
Lemma~\ref{rev:lem:path-sum} supplying the weight count; forcing terms
from the omitted topological types, bounded by
Lemma~\ref{lem:linear-tails}; and a convergent linear majorant system for
the normalized errors, solved by induction on the order. The difference
from Subsection~\ref{app:remainder-bounds} is that the marking count and
the order may now be linear in $h$, so no smallness of $L^2/h$ is
available, and only the ratios $r/h$ and $L/h$ are used.

\begin{proof}[Proof of \eqref{eq:gevrey-remainders}]
Choose $0<\varsigma_0<\varsigma_1<\varsigma_2<\log\rho$.
Use the preceding formal bounds with base exponent $\varsigma_0$;
for actual remainders use $\varsigma_1+a$, with
$0<a\le\varsigma_2-\varsigma_1$. Put $Z_L=p+ML/a$, with
$M\ge\max\{2,2(\varsigma_2-\varsigma_1)\}$ as in
\eqref{rev:eq:gevrey-weights}, and $\delta=Z_L/h$.
We prove geometric bounds for the remainder divided by $\delta^L$,
with the additional factor $\mathrm{e}^{(\varsigma_1+a)D}$ for a profile.
Precisely, for $L\ge1$ put
\[
 R_E(L)=\sup\frac{\bigl|c_{g,p}(\mathbf d)-\sum_{j<L}e_{p,j}(\mathbf d)h^{-j}\bigr|}
        {\delta^L\,\mathrm e^{(\varsigma_1+a)D}},\qquad
 R_C(L)=\sup\frac{\bigl|\mathsf C_p(h)-\sum_{j<L}\mathsf C_{p,j}h^{-j}\bigr|}
        {\delta^L},
\]
and $R_H(L)$, $R_G(L)$ in the same way for $\mathsf G_p$ with $p\ge2$
and with $p=1$. The suprema are over real $Y\in[0,y_0]$, all stable
states with $p\ge1$, all degree tuples (invalid tuples having profile
zero) and all $a\in(0,\varsigma_2-\varsigma_1]$; the supremum over $a$
is what allows a child remainder to be estimated below at the smaller
exponent $a'=ak/L$.
The constants are independent of the allowed $a,p,L,h$.

\textbf{Large weights and degrees.}
The exact profiles and scalar ratios are bounded by
\eqref{eq:TL15}--\eqref{eq:TL17} and \eqref{eq:TL27}.
The formal coefficient bounds are at most $CB_0^jZ_L^j$ for $j<L$,
with the profile exponential included: if $M_0$ is the constant of
\eqref{eq:gevrey-coefficients} at the exponent $\varsigma_0$, then
$(p+M_0j)^j\le C^j(p+ML/a)^j$ for $j<L$ with
$C=\max\{1,M_0(\varsigma_2-\varsigma_1)/M\}$, because
$a\le\varsigma_2-\varsigma_1$ gives $ML/a\ge ML/(\varsigma_2-\varsigma_1)$;
this comparison uses that $a$ is bounded above, and would fail for
unbounded $a$. Choose $\delta_0$ so that
$B_0\delta_0<1/2$. When $\delta\ge\delta_0$, the normalized error is
at most $C\delta_0^{-L}$ by the large-weight calculation in the proof
of Proposition~\ref{prop:hierarchy}.

Suppose $\delta<\delta_0$. Then $p\le\delta_0h$ and
$L/h\le\delta_0(\varsigma_2-\varsigma_1)/M$. Shrink $\delta_0$ to fit
Lemma~\ref{lem:linear-tails} and to bound the cusp ratios inverted
below away from zero. Every state occurring below has
marking count at most $p+L+2$ and complexity at least $h-L$, so for
$\delta_0$ small the lemma applies at all of them. To treat
$D\ge h/4$, note that $B_0\delta<1/2$ and the formal bounds at the
exponent $\varsigma_0$ give
$\sum_{j<L}|e_{p,j}(\mathbf d)|h^{-j}\le2C\mathrm e^{\varsigma_0D}$,
while $0\le c_{g,p}\le1$. After division by
$\mathrm e^{(\varsigma_1+a)D}$ the error is therefore at most
$(1+2C)\mathrm{e}^{-b_0h}$ with $b_0=(\varsigma_1-\varsigma_0)/4$.
Since
\[
 \mathrm{e}^{-b_0h}h^L\le(L/(b_0\mathrm{e}))^L,\qquad
 Z_L\ge ML/(\varsigma_2-\varsigma_1),
\]
this contributes only a geometric forcing sequence. Invalid parents
are included: their dimension is $(3h-p)/2>h/4$ in this region.
We may use the exact adjacent recursion at $D<h/4$.

\textbf{Retained products and shifts.}
Retain losses $r<L$. The factorization \eqref{eq:FE2} is an exact
finite product of shifted cusp and genus ratios. Apply
\eqref{eq:FE-product} with the large profile first. Each term contains
one unknown remainder, either the profile remainder or a scalar
remainder preceded by known profile and scalar coefficients. All exact
factors after it are scalar and bounded. There are $O(r)$ such factors,
costing at most $C^r$. For descent use \eqref{eq:FE-inverse}; its
exact inverse is bounded in the small-weight region.

The exact kernel sums can be extended by zero to all $\ell\ge0$.
Indeed, with parent dimension $N=3g-3+p$ and deficit $m=N-D$,
a merger has child dimension $N-1$ and degree $D+\ell-1$;
a nonseparating child has dimension $N-2$ and degree $D+\ell-2$.
Both deficits are $m-\ell$. In a split the sum of child deficits
is also $m-\ell$, so a valid bounded child and $\ell>m$ force the
large child to be invalid. The expanded sums converge absolutely by
$\varsigma_2<\log\rho$.

In each resulting term, the loss $r$, total known coefficient order $i$,
unique remainder order $q\ge1$, and shift order $v$ satisfy
$r+i+q+v=L$. Thus the weights are allocated to one total order $L$.
Lemma~\ref{rev:lem:path-sum} applies verbatim with $j$ replaced by $L$
and $(\sigma_0,\sigma_1)$ by $(\varsigma_1,\varsigma_2)$, the normalized
supremum of the one unknown factor being its error sequence and those
of the known factors being the normalized formal coefficients.
For a linear branch let $k\le L-1$ be the order of the profile factor,
whether it is the remainder or a known coefficient. If it is the
remainder and $k>0$, bound it in the norm with exponent $\varsigma_1+a'$,
$a'=ak/L$; then \eqref{rev:eq:path-normalized} gives the degree-path cost
\[
 \sup_{D\ge0}(D+p+1)\mathrm{e}^{-a(L-k)D/L}
 \le p+1+\frac{L}{\mathrm{e}\,a(L-k)}\le3Z_L .
\]
If it is a known coefficient, its formal bound carries the exponential
$\mathrm{e}^{\varsigma_0D'}$ at the child degree sum $D'\le D+\ell$, and
$\varsigma_0<\varsigma_1\le\varsigma_1+a'$; hence
\[
 \mathrm{e}^{\varsigma_0D'}\le\mathrm{e}^{\varsigma_0\ell}\,\mathrm{e}^{(\varsigma_1+a')D},
\]
so the known coefficient enters the path sum exactly as a remainder of
the same order would, with the margin $\varsigma_1-\varsigma_0$ to spare.
This is how the formal-coefficient bound at the exponent $\varsigma_0$
is absorbed into the norm $\mathrm{e}^{(\varsigma_1+a)D}$ of the actual
remainders. At $k=0$ the profile is one.
All other weights have the same product bounds as in the formal proof.
For a split, \eqref{rev:eq:path-dimension}--\eqref{rev:eq:path-separating-count}
remove the extra degree count, giving $C^rZ_L^r$ against its prefactor
of order $r$.

A remainder of order $q\le L$ at complexity $h-c$, $|c|\le r$,
costs the additional factor
\[
 \Bigl(\frac h{h-c}\Bigr)^q\le \mathrm{e}^{2rq/h}\le\mathrm{e}^{2rL/h}\le C^r,
\]
because $|c|/h\le r/h\le1/2$ in the small-weight region, where
$-\log(1-x)\le2x$ for $0\le x\le1/2$, and
$L/h\le\delta_0(\varsigma_2-\varsigma_1)/M$ is bounded. Only $r/h$ and
$L/h$ enter here; no smallness of $L^2/h$ is needed, in contrast to the
proof of Proposition~\ref{prop:hierarchy}. Unit shifts in
dilaton and descent cost a fixed constant. For a known shift tail
whose first omitted degree is $q$, the finite negative-binomial bound
used in Proposition~\ref{prop:hierarchy} gives
\[
 \frac{(Lr)^q}{q!}h^{-q}(1-r/h)^{-L}.
\]
After division by $Z_L^qh^{-q}$, the inequalities
$L/Z_L\le a/M$ and $r^q/q!\le \mathrm{e}^r$ leave a geometric $q$ bound
times $C^r$; again no condition on $L^2/h$ is needed.

\textbf{Tail forcing and the base case.}
For $L\ge2$, Lemma~\ref{lem:linear-tails} bounds the omitted exact
topology by $C[B(p+L)/h]^L$. In a profile it is multiplied by at
most a kernel-weighted quadratic polynomial in $D$. The norm
$\mathrm{e}^{-(\varsigma_1+a)D}$ absorbs that polynomial uniformly. This gives
$CB^L\delta^L$. At $L=1$ the first-order estimates in the same
lemma give the bound directly. The large-weight and large-degree
regions have already supplied geometric forcing.

\textbf{A convergent error majorant.}
Regard the four normalized error sequences $R_E,R_C,R_H,R_G$ defined
above as series in $z$, with zero constant terms and positive orders
starting at $L=1$.
Known formal coefficients, inverse-cusp coefficients and exact scalar
factors are bounded by one series
$V(z)=\upsilon_0/(1-\upsilon_1z)$, with $\upsilon_0,\upsilon_1\ge1$ (the symbol $\upsilon_1$ is unrelated to the instanton coefficient $v_1=-1/8$ of Subsection~\ref{sec:endpoint-expansion}).
Choose a fixed $b\ge1$ large enough for the retained-product constants.
Put $x(z)=bzV(z)^8$, and set
\[
 k(z)=\frac{bV(z)^8}{1-bz}
     \left(z+\frac{1+x(z)}{(1-x(z))^3}-1\right).
\]
This analytic series has nonnegative coefficients and $k(0)=0$.
Its second term is the sum, over $r\ge1$, of
$b(r+1)^2(bz)^rV(z)^{8r+8}/(1-bz)$.
A separating loss-$r$ product, with bounded-component genus $u$,
has at most $r+u+1\le3r/2+1$ expandable factors, hence $O(r)$ possible remainder positions. Together with the
$O(r)$ node choices, these give the displayed quadratic factor;
known-order allocations are already counted by the powers of $V$.
Exponential costs
in $r$ are absorbed in $b^r$. The first term covers linear branches
and positive-order unit shifts. Terms whose only error is a truncated
known shifted power contain no unknown remainder. Their same analytic
loss sum is included in the forcing.

Let $J_4$ be the four-by-four matrix with all entries one, and define
\[
 N_0=\begin{pmatrix}0&0&0&0\\b&0&0&0\\b&0&0&0\\0&b&b&0\end{pmatrix},
 \qquad F(z)=\frac{f_0z}{1-f_1z}(1,1,1,1)^{\mathsf T}.
\]
Choose $f_0,f_1$ after $b,\upsilon_0,\upsilon_1$ to bound all forcing sequences.
The only current-order dependencies are profile to cusp, profile to
genus at $p\ge2$, and cusp/genus to descent. They are covered by
$N_0$; all remaining errors have lower order. Thus, componentwise,
\[
 R\preceq F+N_0R+kJ_4R.
\]
Since $N_0^4=0$, $(I-N_0)^{-1}$ has nonnegative entries. The series
\[
 U(z)=(I-N_0-k(z)J_4)^{-1}F(z)
\]
is analytic near zero, with nonnegative coefficients bounded by $CB^L$.
Induction in $L$, then in stage order, proves $R\preceq U$ and proves
finiteness of the suprema at the same time. Finally take
$\varsigma_1+a=\sigma$ and absorb the fixed $1/a$ into $M$.
This proves \eqref{eq:gevrey-remainders}.
\end{proof}

For fixed $p$, $(p+ML)^L\le C_p^LL!$. Thus the scalar ratios have
Gevrey-1 remainders. The same holds for $\sum_{\mathbf d}\lambda_{\mathbf d}c_{g,p}(\mathbf d)$
whenever $\sum_{\mathbf d}|\lambda_{\mathbf d}|\mathrm e^{\sigma|\mathbf d|}<\infty$
for the exponent $\sigma$ of \eqref{eq:gevrey-remainders}; the dilaton
functional $f\mapsto\sum_{d\ge1}b_df(d,0^p)$ is of this kind.

\subsection{Determining the constants at the pure-\texorpdfstring{$\psi$}{psi} endpoint}\label{sec:endpoint-expansion}

To integrate the transport equation to all orders, we need the value of
every coefficient at the pure-$\psi$ endpoint, with constants controlled
as the order grows. This subsection supplies these values and then
proves Theorem~\ref{thm:full-expansion}. It proceeds in four stages.
First, the exact endpoint identity following \eqref{eq:ST12} shows
that the higher cusp-ratio coefficients vanish at $y_0$, see
\eqref{eq:FE9}. Second, the large-order asymptotics of the
Painlev\'e~I zero-instanton coefficients are made quantitative: the
Borel representation of Costin, Costin and Huang~\cite{CCH} gives a
Gevrey-1 remainder for the jump between adjacent tritronqu\'ee
solutions, whose first coefficient $v_1=-1/8$ is derived directly, and
the contour formula of Garoufalidis, Its, Kapaev and
Mari\~no~\cite{GIKM} transfers this remainder to the coefficients
themselves. Third, the Itzykson--Zuber identity~\cite{IZ} converts
these into the expansion \eqref{eq:FE10} of $\pi C_g(2^{3g-3})$ and,
through string and dilaton, into the exact normalized endpoint values
\eqref{eq:FE14}; Lemma~\ref{rev:lem:endpoint-gevrey} records the
resulting Gevrey-1 bounds. Finally, the coefficients $a_{p,j}$ and
$\mathcal L_{p,j}$ of Theorem~\ref{thm:full-expansion} are defined in
\eqref{eq:FE11}, and the theorem is proved by integrating the
logarithmic derivative from a matching point near $y_0$, using
Proposition~\ref{prop:gevrey-hierarchy}, the inverse identity
\eqref{eq:FE-inverse} and the positive comparison
\eqref{eq:PF15}--\eqref{eq:PF16}.

Invert the cusp expansion and write, at each fixed order,
\begin{equation}
 \frac1{\mathsf C_p(h)}
 =1+\frac{\beta_p}h+\sum_{j=1}^K\frac{\gamma_{p,j}}{h^{j+1}}
      +O_{p,K}(h^{-K-2}).
 \label{eq:FE8}
\end{equation}
Here $\gamma_{p,1}=\gamma_p$. Three facts determine the endpoint values
of these coefficients. First, the remainder estimate
\eqref{eq:gevrey-remainders} holds uniformly for real $Y\in[0,y_0]$, in
particular at $Y=y_0$, and the finite inverse identity
\eqref{eq:FE-inverse} transfers it to $1/\mathsf C_p(h;y_0)$; thus
\eqref{eq:FE8} is a valid asymptotic expansion at $Y=y_0$. Second, the
formal coefficients $\gamma_{p,j}$ are holomorphic on the neighborhood
$\Omega$ of $[0,y_0]$, hence continuous at $y_0$, so the coefficients of
that expansion are the values $\beta_p(y_0)$ and $\gamma_{p,j}(y_0)$.
Third, the exact endpoint identity following \eqref{eq:ST12} gives
$1/\mathsf C_p(h;y_0)=1-p/(5h)$, a finite expansion whose only term is
$\beta_p(y_0)/h=-p/(5h)$. Uniqueness of asymptotic expansions in
$h^{-1}$ therefore gives
\begin{equation}
 \gamma_{p,j}(y_0)=0\quad(j\ge1),\qquad
 \gamma_{p,j}/(Ar)\ \hbox{is analytic at }y_0.
 \label{eq:FE9}
\end{equation}
Analyticity follows from the finite coefficient rules, not from remainder
estimates. Division by the simple zero of $Ar$ is a bounded operation
on a fixed smaller complex neighborhood: write a vanishing function
as $(Y-y_0)\int_0^1 f'(y_0+t(Y-y_0))\,\mathrm{d}t$ and apply Cauchy's bound
once on the larger neighborhood. No successive shrinkage in order occurs.

We need a quantitative version of the pure-$\psi$ endpoint input.
Use the conventions of Garoufalidis--Its--Kapaev--Mari\~no~\cite{GIKM}:
$-u''/6+u^2=Z$, zero-instanton coefficients $u_{g,0}$,
$\mathfrak a=8\sqrt3/5$, $\nu_g=2g-1/2$ and
$S_1=-i3^{1/4}/(2\sqrt\pi)$. Put
\[
 Z_g=\frac{\pi i}{S_1}
       \frac{\mathfrak a^{\nu_g}u_{g,0}}{\Gamma(\nu_g)}.
\]
The Borel representation of Costin--Costin--Huang~\cite[Theorem~2]{CCH}
makes the coefficient formula in \cite[Section~2.3]{GIKM} quantitative,
as follows. That formula gives the expansion of $u_{g,0}$ to each fixed
order; what is needed here, and proved below, is a remainder whose
factorial growth in the order is uniform. Set $\varrho=6^{-1/5}\mathrm{e}^{i\pi/5}$;
$y(Z_c)=\varrho^2u(\varrho Z_c)$ satisfies $y''=6y^2+Z_c$.
Their Boutroux coordinate is $\xi=\mathfrak a Z^{5/4}$, and the
normalized remainder is $u(Z)/\sqrt Z-1+4/(25\xi^2)$. The CCH variable
$x$ is our $\xi$, and their unknown $h(x)$ is this normalized remainder.
The adjacent GIKM solutions $u_1,u_0$ have sectors
$(-2\pi/5,6\pi/5)$ and $(-6\pi/5,2\pi/5)$ respectively.
Their normalized remainders are $O(\xi^{-4})$ on $+i\infty$ for $u_1$
and on $-i\infty$ for $u_0$, each ray lying inside the corresponding
sector. Thus the respective lateral constants in the exact Borel
representation vanish: multiplying by $\mathrm{e}^\xi\xi^{1/2}$ gives zero
on those rays. This is the characterization of the tritronqu\'ee
solutions as the tronqu\'ee solutions whose constant beyond all orders
vanishes, \cite[Theorem~7 and Note~8(ii)]{CCH}. It identifies the
lateral Borel sums with the solutions themselves, not merely with their
formal expansions.

The Borel transform at its first positive singularity has the form
$S(1-s)^{-1/2}\alpha(s)+\beta(s)$, where $\alpha,\beta$ are analytic
near one and $\alpha(1)=1$~\cite[Proposition~3(i)]{CCH}. We keep the
constant $S$ symbolic: only the structure of this singularity is used.
The jump between the two solutions is normalized by the Stokes
multiplier $S_1$ of Kapaev~\cite{Kap}, in the form recorded
in~\cite[Section~2.1]{GIKM}; with this value of $S_1$,
\[
 \mathfrak J(\xi)=
 \frac{u_1(Z)-u_0(Z)}{-S_1\mathrm{e}^{-\xi}Z^{-1/8}},\qquad
 \mathfrak J(\xi)\longrightarrow1.
\]
Choose a fixed small $r>0$ and angles $\pm\phi$ so
$(1+r)\cos\phi>1+r/2$. Deform only the finite initial parts of the
two Laplace rays into a local Hankel contour about $[1,1+r]$ and
compact joining arcs on $|s|=1+r$. The infinite ray tails are kept
unchanged. The analytic continuation and weighted $L^1$ estimates
of~\cite[Lemma~17(c),(f)]{CCH} bound these outer parts by $C\mathrm{e}^{-(1+\varepsilon)\xi}$,
with constants independent of Taylor order.
The local analytic term has no jump. Cauchy's estimate on $|s-1|\le2r$
bounds the degree-$L-1$ remainder of $\alpha(1+t)$ by
$C(2r)^{-L}t^L$. Its normalized Laplace integral is bounded by
$CB^L\Gamma(L+1/2)\xi^{-L}$. Extending the retained polynomial
integrals from $r$ to infinity obeys the same bound; the outer term
$O(\sqrt\xi \mathrm{e}^{-\varepsilon\xi})$ is at most $CB^LL!\xi^{-L}$.
Consequently, for all $L\ge1$ and real $\xi\ge\xi_0$,
\[
 \left|\mathfrak J(\xi)-\sum_{j<L}v_j\xi^{-j}\right|
       \le CB^LL!\xi^{-L},\qquad v_j=\mathfrak a^ju_{j,1}.
\]
The value of the Borel constant $S$ of~\cite{CCH} is not used; the only
numerical Stokes datum is Kapaev's $S_1$.

The first coefficient is $v_1=-1/8$, and we derive it here because it
determines the first correction at the endpoint. The difference
$u_1-u_0$ is exponentially small and satisfies the linearization
$(u_1-u_0)''=12u_0(u_1-u_0)$ of $-u''/6+u^2=Z$ up to a term quadratic in
$u_1-u_0$, which is $O(\mathrm{e}^{-2\xi})$ and does not affect the
coefficients $v_j$. Here $12u_0=12Z^{1/2}-Z^{-2}/4+O(Z^{-9/2})$, because
$u_{1,0}=-1/48$. Put $f=Z^{-1/8}\mathrm{e}^{-\xi}$. Since
$\xi'=2\sqrt3Z^{1/4}$ and $\xi''=(\sqrt3/2)Z^{-3/4}$,
\[
 \frac{f''}f=\Bigl(\frac1{8Z}+\xi'\Bigr)^2+\frac1{8Z^2}-\xi''
 =12Z^{1/2}+\frac9{64}Z^{-2};
\]
the terms of order $Z^{-3/4}$ cancel, which is why the prefactor is
$Z^{-1/8}$. Writing $u_1-u_0=-S_1f\mathfrak J$ and dividing the
linearized equation by $-S_1f$ gives
$(f''/f)\mathfrak J+2(f'/f)\mathfrak J'+\mathfrak J''=12u_0\mathfrak J$.
The local Laplace representation above gives the expansion of
$\mathfrak J$ in a sector about the positive ray; Cauchy's estimate
therefore permits differentiation of its remainder on a smaller sector.
With $\mathfrak J=1+v_1\xi^{-1}+O(\xi^{-2})$ one has
$\mathfrak J'=-v_1\xi'\xi^{-2}+O(Z^{-7/2})=-(2\sqrt3v_1/\mathfrak a^2)Z^{-9/4}+O(Z^{-7/2})$,
so $2(f'/f)\mathfrak J'=-2\xi'\mathfrak J'+O(Z^{-13/4})=(24v_1/\mathfrak a^2)Z^{-2}+O(Z^{-13/4})$,
while $\mathfrak J''=O(Z^{-13/4})$. Comparing the coefficients of
$Z^{-2}$ gives $9/64+24v_1/\mathfrak a^2=-1/4$, that is
$v_1=-25\mathfrak a^2/1536=-1/8$, since $\mathfrak a^2=192/25$.
Equivalently $u_{1,1}=-5/(64\sqrt3)$.

The exact coefficient contour in~\cite[Section~2.3]{GIKM} expresses
$u_{g,0}$ as a circle integral and deforms it onto the five rays
$\arg t=2\pi k/5$, $k=0,\pm1,\pm2$; the symmetry of the tritronqu\'ee
solutions reduces their sum to five times the positive ray
\cite[equations~(2.25) and~(2.29)]{GIKM}:
\[
 u_{g,0}=-\frac5{2\pi i}\int_{R_0}^\infty t^{5g-2}
        \bigl(u_1(t^2)-u_0(t^2)\bigr)\,\mathrm{d}t+r_g(R_0),
 \qquad |r_g(R_0)|\le CR_0^{5g},
\]
for a fixed radius $R_0$, large enough for~\cite{GIKM} and such that
$\xi_0=\mathfrak aR_0^{5/2}$ is at least the threshold in the jump
expansion above. Put $\xi=\mathfrak at^{5/2}$, so that
$t^{5g-9/4}\,\mathrm{d}t=\tfrac25\mathfrak a^{-\nu_g}\xi^{\nu_g-1}\,\mathrm{d}\xi$,
and write $u_1-u_0=-S_1\mathrm{e}^{-\xi}t^{-1/4}\mathfrak J(\xi)$. The
factor $5\cdot\tfrac25$ cancels the $2$ in $2\pi i$, and the definition
of $Z_g$ gives exactly
\[
 Z_g=\frac1{\Gamma(\nu_g)}\int_{\xi_0}^{\infty}
        \mathrm{e}^{-\xi}\xi^{\nu_g-1}\mathfrak J(\xi)\,\mathrm{d}\xi
       +O\!\left(\frac{D_1^{\,g}}{\Gamma(\nu_g)}\right),
 \qquad D_1=\mathfrak a^2R_0^{5};
\]
the error term is $(\pi i/S_1)\mathfrak a^{\nu_g}r_g(R_0)/\Gamma(\nu_g)$,
of modulus at most $C|\pi/S_1|\mathfrak a^{-1/2}D_1^{\,g}/\Gamma(\nu_g)$
since $\nu_g=2g-1/2$. Applying the jump remainder
with $L=K+1\le g/2$ gives
$CB^LL!\Gamma(\nu_g-L)/\Gamma(\nu_g)\le CB^LL!g^{-L}$.
The omitted low integrals and inner circle are
$O(\mathrm{e}^{-c g\log g})$ uniformly for $K\le c_0g$, with $c_0$ small.
For the low integrals, bound the finite sum by
$C\xi_0^{\nu_g}\sum_{j<L}B^jj!\xi_0^{-j}/\Gamma(\nu_g)$
and use Stirling's lower bound.
Thus the finite falling-factorial expansion of $Z_g$ has a Gevrey-1
remainder. Re-expansion in $g^{-1}$ preserves it:
\[
 \frac{\Gamma(\nu_g-j)}{\Gamma(\nu_g)}
 =(2g)^{-j}\prod_{a=1}^j\left(1-\frac{a+1/2}{2g}\right)^{-1}.
\]
Cauchy on radius $1/(D_0(K+1))$ for the product, with a fixed $D_0\ge4$,
on which each factor has modulus at most two and which contains
$1/(2g)$ in its half-radius disk once the linear-order constant is
small, followed by
$j!(K+1)^{K+1-j}\le \mathrm{e}^{K+1}(K+1)!$, bounds the total omitted
inverse-genus powers by $CB^{K+1}(K+1)!g^{-K-1}$.

The exact identity of Itzykson--Zuber~\cite[equation~(6.10)]{IZ} is
\[
 \langle\tau_2^{3g-3}\rangle_g
 =\frac{2^g(3g-3)!}{(5g-5)(5g-3)}\,\phi_g,\qquad g\ge2,
\]
where $\phi(t)=\sum_{g\ge0}\phi_g t^{(1-5g)/2}$, $\phi_0=-1$,
satisfies $\phi''/3+\phi^2=t$. The substitution
$\phi(t)=-2^{1/5}u(2^{-2/5}t)$ gives exactly
$\phi_g=-2^g u_{g,0}$.
Combined with the expansion of $Z_g$ above, this gives, in the normalization of \eqref{eq:H4} and for every fixed $K$,
\begin{equation}
 \pi C_g(2^{3g-3})=1+\sum_{j=1}^K d_jg^{-j}+O_K(g^{-K-1}).
 \label{eq:FE10}
\end{equation}
Lemma~\ref{rev:lem:endpoint-gevrey}(iii) below gives the Gevrey-1 form
of this remainder, uniform for $K\le c_0g$.
The leading constant is also independently given by
\cite[Theorem~1]{GYZ}. To include the markings, put $h_0=2g-2$, $h=h_0+p$,
$a_*=2/15$, and define the exact normalized endpoint
\[
 E_{g,p}=\frac{\pi\overline U_{g,p}(0^p;y_0)}
 {\Gamma(h)4^{-g}a_*^{-h}\sqrt{a_*/(\pi h)}}.
\]
For $p\ge1$ this is the limit of $\pi\Phi_g^{(p)}/\mathcal J_{g,p}$;
at $p=0$ it is defined by the same exact pure-$\psi$ expression.
String and dilaton give
$E_{g,p+1}/E_{g,p}=(h_0+4p/5)(h_0+p)^{-1}
\sqrt{(h_0+p+1)/(h_0+p)}$. The base conversion gives exactly
$\overline U_{g,0}(y_0)=-108^gu_{g,0}/[27(5g-5)(5g-3)]$.
Substitution therefore proves
\begin{equation}
 \begin{split}
 E_{g,p}&=Z_g\,\frac{\Gamma(h_0+3/2)}{\Gamma(h_0)h_0^{3/2}}
             \frac{R_p(h_0)}{1+4/(5h_0)},\\
 R_p(h_0)&=\sqrt{1+p/h_0}\prod_{i=0}^{p-1}
                       \frac{1+4i/(5h_0)}{1+i/h_0}.
 \end{split}
 \label{eq:FE14}
\end{equation}
This identity is exact for $g\ge2,p\ge0$, with the empty product one.

The quantitative endpoint input is the following lemma. It is used at
the matching point in the proof of Theorem~\ref{thm:full-expansion}, and
it gives the Gevrey-1 form of \eqref{eq:FE10}.

\begin{lemma}[Gevrey-1 endpoint expansion]\label{rev:lem:endpoint-gevrey}
Let $h_0=2g-2$ and $h=h_0+p$. There are constants $c_0,C,B>0$ and, for
each $p\ge0$, constants $C_p,B_p>0$ such that the following hold for all
$g\ge2$ and $0\le K\le c_0g$.
\begin{enumerate}[label=\textup{(\roman*)},leftmargin=*]
\item With $v_0=1$ and $v_j=\mathfrak a^ju_{j,1}$,
\[
 \Bigl|Z_g-\sum_{j=0}^{K}v_j\frac{\Gamma(\nu_g-j)}{\Gamma(\nu_g)}\Bigr|
 \le CB^{K+1}(K+1)!\,g^{-K-1},
\]
and the same bound holds after re-expanding the sum in powers of
$g^{-1}$ or of $h_0^{-1}$ and truncating at order $K$.
For each fixed $p$, re-expansion in powers of $h^{-1}$ has remainder
at most $C_pB_p^{K+1}(K+1)!h^{-K-1}$.
\item The function $\log E_{g,p}$ has an asymptotic expansion
$\sum_{j\ge1}a^*_{p,j}h^{-j}$ as $g\to\infty$, whose coefficients are
obtained by re-expanding \eqref{eq:FE14} and \textup{(i)} in powers of
$h^{-1}$, and
\[
 \Bigl|\log E_{g,p}-\sum_{j=1}^Ka^*_{p,j}h^{-j}\Bigr|
 +\Bigl|E_{g,p}-\sum_{j=0}^K\mathcal E_{p,j}h^{-j}\Bigr|
 \le C_pB_p^{K+1}(K+1)!\,h^{-K-1},
\]
where $\mathcal E_{p,j}=[t^j]\exp\bigl(\sum_{i=1}^ja^*_{p,i}t^i\bigr)$.
\item For the coefficients $d_j$ of \eqref{eq:FE10}, with $d_0=1$,
\[
 \Bigl|\pi C_g(2^{3g-3})-\sum_{j=0}^Kd_jg^{-j}\Bigr|
 \le CB^{K+1}(K+1)!\,g^{-K-1}.
\]
\end{enumerate}
\end{lemma}

\begin{proof}
(i) This is the Gevrey-1 remainder for $Z_g$ established above: the jump
remainder in the integral representation of $Z_g$ gives the bound in
falling factorials for $L=K+1\le g/2$, and Cauchy's estimate on the
radius $1/(D_0(K+1))$ for the product
$\prod_{a\le j}(1-(a+1/2)/(2g))^{-1}$, followed by
$j!(K+1)^{K+1-j}\le\mathrm{e}^{K+1}(K+1)!$, converts it into a bound in
powers of $g^{-1}$. Since $g^{-1}=2h^{-1}(1+(2-p)h^{-1})^{-1}$ is a
rational function of $h^{-1}$, analytic on a disk about zero depending
on the fixed value of $p$, the
same Cauchy estimate applied to the degree-$K$ polynomial on the radius
$1/(D_p(K+1))$ gives the re-expansion in $h^{-1}$ with constants
depending on $p$; taking $p=0$ gives the re-expansion in $h_0^{-1}$.

(ii) By \eqref{eq:FE14},
\[
 \log E_{g,p}=\log Z_g+\log\frac{\Gamma(h_0+3/2)}{\Gamma(h_0)h_0^{3/2}}
 +\log R_p(h_0)-\log\Bigl(1+\frac4{5h_0}\Bigr).
\]
The first term is covered by (i), since $Z_g=1+O(g^{-1})$ stays in a
fixed disk on which the logarithm is analytic. For the Gamma factor,
for fixed $\alpha>0$,
\[
 \log\frac{\Gamma(h_0+\alpha)}{\Gamma(h_0)h_0^\alpha}
 =\int_0^\infty \mathrm{e}^{-h_0t}
 \left[\alpha-\frac{1-\mathrm{e}^{-\alpha t}}{1-\mathrm{e}^{-t}}\right]\frac{\mathrm{d}t}{t}.
\]
The kernel is analytic at zero on a fixed disk and bounded outside a
small positive interval. Its global Taylor remainder is bounded by
$CB^Lt^{L-1}$, giving $CB^L(L-1)!h_0^{-L}$ after integration.
The finite insertion factors $R_p(h_0)$ and $(1+4/(5h_0))^{-1}$ have
convergent expansions at $h_0^{-1}=0$.
Finite shifts from $h_0$ to $h=h_0+p$ preserve the factorial error by the
same Cauchy estimate for rational products used in (i).
For logarithms, exponentials and inverses, take the degree-$K$ polynomial
on radius $1/(D_0(K+1))$; its unit or zero-constant part stays bounded
away from the relevant singularity. Cauchy's tail and the bounded
real derivative of the composed function preserve the remainder bound.
This proves (ii); uniqueness of asymptotic expansions identifies
$a^*_{p,j}$ with the coefficients of the formal re-expansion.

(iii) By the Itzykson--Zuber identity, the conversion
$\phi_g=-2^gu_{g,0}$ and the definition of $Z_g$, in the normalization
\eqref{eq:H4},
\[
 \pi C_g(2^{3g-3})=Z_g\,\Xi(g),\qquad
 \Xi(g)=\frac{3^{1/4}}{2\sqrt\pi}\cdot
 \frac{4^g\,15^{3g-3}\,2^{2g}\,\Gamma(3g-2)\,\Gamma(\nu_g)}
      {3^{5g-5}\,\Gamma(5g-5)\,(5g-5)(5g-3)\,\mathfrak a^{\nu_g}},
\]
where $3^{1/4}/(2\sqrt\pi)=-S_1/i$; this is an explicit product of
powers $c^g$, of Gamma functions at arguments $\alpha g+\beta$ and of a
rational function of $g$. Stirling's formula
$\Gamma(\alpha g+\beta)\sim\sqrt{2\pi}\,(\alpha g)^{\alpha g+\beta-1/2}\mathrm{e}^{-\alpha g}$
gives $\Xi(g)=1+O(g^{-1})$ directly: the exponentials and the powers
of $g$ and of $5$ cancel, the powers of $2$ and $3$ reduce to
$2^{1/2}3^{-1/4}$, and
$\sqrt{2\pi}\,2^{1/2}3^{-1/4}\cdot3^{1/4}/(2\sqrt\pi)=1$. Together
with $Z_g=1+O(g^{-1})$, which is the leading case of
\cite[Section~2.3]{GIKM}, this proves \eqref{eq:PF12a}. Binet's formula
$\log\Gamma(x)=(x-\tfrac12)\log x-x+\tfrac12\log2\pi+\mu(x)$, with
$\mu(x)=\int_0^\infty\bigl(\frac1{\mathrm{e}^t-1}-\frac1t+\frac12\bigr)\frac{\mathrm{e}^{-xt}}t\,\mathrm{d}t$,
expresses $\log\Xi(g)$ as a sum of elementary terms and of integrals
$\pm\mu(\alpha g+\beta)$. The elementary terms are convergent series in
$g^{-1}$ once their contributions of orders $g\log g$, $g$, $\log g$
and $g^0$ are collected, and these must cancel because $\log\Xi(g)\to0$.
The kernel of $\mu$ is analytic on $|t|<2\pi$ and bounded on the positive
axis, so the Taylor-remainder argument just used for the Gamma factor
gives for each $\mu(\alpha g+\beta)$ an expansion in $(\alpha g+\beta)^{-1}$
with remainder $CB^{L}(L-1)!(\alpha g+\beta)^{-L}$, and the rational
re-expansion in $g^{-1}$ is handled as in (i). Hence $\log\Xi(g)$, and
by exponentiation $\Xi(g)$, have Gevrey-1 expansions in $g^{-1}$ for
$K\le c_0g$. Multiplying by (i), as in (ii), proves (iii).
\end{proof}

At first order, \eqref{eq:FE14}, the value $v_1=-1/8$ derived above,
and $\nu_g-1=h_0+1/2$ give
\[
 a^*_{p,1}=-\frac18+\frac38-\frac45+\frac p2-\frac{p(p-1)}{10}
 =\frac{-2p^2+12p-11}{20}=t_p,
\]
the constant in \eqref{eq:ST13}; replacing $h_0$ by $h=h_0+p$ does not
change the first-order coefficient.

In particular the first coefficient in \eqref{eq:FE10} is $d_1=-2/9$.
Indeed, at $p=0$ the definition of $E_{g,0}$ and \eqref{eq:PF13} give
\[
 \pi C_g(2^{3g-3})=\frac{E_{g,0}}{N_{h_0}},\qquad
 N_{h_0}=\frac{\mathcal M_{g,0}(y_0)}
          {\Gamma(h_0)4^{-g}a_*^{-h_0}/\sqrt{2\pi h_0v}}
     =\sqrt{2\pi h_0v}\,\Pr(S_{h_0}=3h_0/2),\qquad v=\frac{15}4,
\]
where $S_{h_0}$ is the sum of $h_0$ independent copies of the geometric law
with generating function $(2/5)/(1-3w/5)$; this is the normalized
negative-binomial coefficient of Theorem~\ref{thm:first-correction} at
displacement zero, so $N_{h_0}=1-19/(180h_0)+O(h_0^{-2})$. On the other side,
$Z_g=1+v_1\Gamma(\nu_g-1)/\Gamma(\nu_g)+O(g^{-2})=1-1/(8h_0)+O(h_0^{-2})$,
because $\nu_g-1=h_0+1/2$ and $v_1=-1/8$;
$\Gamma(h_0+3/2)/(\Gamma(h_0)h_0^{3/2})=1+3/(8h_0)+O(h_0^{-2})$; and
$(1+4/(5h_0))^{-1}=1-4/(5h_0)+O(h_0^{-2})$. Hence, by \eqref{eq:FE14},
\[
 \pi C_g(2^{3g-3})
 =1+\frac1{h_0}\Bigl(-\frac18+\frac38-\frac45+\frac{19}{180}\Bigr)+O(h_0^{-2})
 =1-\frac4{9h_0}+O(h_0^{-2})
 =1-\frac2{9g}+O(g^{-2}),
\]
since $h_0=2g-2$. This value agrees with the expansion recorded by Guo,
Yang and Zagier~\cite[equation~(25)]{GYZ}, which rests on the same
Painlev\'e~I input; it is used in the proof of
Theorem~\ref{thm:first-correction}.

Define
\begin{equation}
 a_{p,j}(Y)=a^*_{p,j}-\int_Y^{y_0}\frac{\gamma_{p,j}(v)}{A(v)r(v)}\,\mathrm{d}v,
 \qquad
 \mathcal L_{p,j}(Y)=[t^j]\exp\!\left(\sum_{i=1}^j a_{p,i}(Y)t^i\right).
 \label{eq:FE11}
\end{equation}
Proposition~\ref{prop:gevrey-hierarchy}, inversion and the bounded division
and integration above give common-domain bounds
$\|a_{p,j}\|+\|\mathcal L_{p,j}\|\le C_p\mathfrak B_p^jj!$.
For formal exponentiation this follows also by Cauchy on the finite
log polynomial at radius $1/(D_0j)$. The first coefficient is
$\mathcal L_{p,1}=\mathcal L_p$. These are the coefficients in
Theorem~\ref{thm:full-expansion}.

\begin{proof}[Proof of Theorem~\ref{thm:full-expansion}]
We first expand the logarithmic derivative, then fix its integration
constant near the endpoint. After integration, exponentiation gives
the claimed expansion. The extra orders retained below compensate for
the error accumulated near the zero of $Ar$.

Put $F_g(Y)=\log(\pi\Phi_g^{(p)}/\mathcal J_{g,p})$.
For fixed $p$ the inverse cusp ratio is bounded, and the finite
inverse identity \eqref{eq:FE-inverse}, with
Proposition~\ref{prop:gevrey-hierarchy}, transfers the remainder to
$1/\mathsf C_p$ with error $CB^NN!h^{-N}$, for $N\le ch$.
The exact dilaton and \eqref{eq:PF11} give
\[
 F'_g=\frac{h}{Ar}(1/\mathsf C_p-1-\beta_p/h).
\]
For $L=K+1$, use $N=L+2=K+3$, retaining $\gamma_{p,j}$ through
$j=L$. Its unknown derivative error is
$CB^{L+2}(L+2)!h^{-L-1}/(Ar)$.
Choose $\delta_*=h^{-L-1}$.
The positive comparison \eqref{eq:PF15} and the exact first equality
of \eqref{eq:PF16} change $F_g$ from its endpoint value by at most
$Ch\delta_*=Ch^{-L}$ on this strip. The leading harmonic
approximation in the second equality of \eqref{eq:PF16} is not used.
Lemma~\ref{rev:lem:endpoint-gevrey}(ii) supplies the value at the
matching point through order $L$, with remainder
$C_pB_p^{L+1}(L+1)!h^{-L-1}$; relative to $L!h^{-L}$ this is
$C_pB_p^{L+1}(L+1)/h\le C_p(c+h^{-1})B_p^{L+1}\le2C_pB_p^{L+1}$, a fixed
exponential in $L$.
The bounded coefficient derivatives from \eqref{eq:FE9} show that the
whole retained coefficient sum changes by $O(\delta_*/h)$ there,
after choosing the linear-order constant small.

Integrate backwards from $y_0-\delta_*$. Since $Ar$ has a simple zero,
the unknown error costs at most
\[
 CB^{L+2}(L+2)!h^{-L-1}[1+(L+1)\log h].
\]
Relative to $L!h^{-L}$, its extra factor is a fixed exponential
base times $(L+1)(L+2)[1+(L+1)\log h]/h$, bounded by another fixed
exponential in $L$. Indeed $L\le ch$ implies $h\ge L/c$;
the decrease of $\log h/h$ for $h\ge e$ bounds this factor by
$C(L+2)^2\log(L+2)\le C\,D_2^{\,L}$ for a fixed $D_2>1$.
The bounded exceptional values of $h,L$ are absorbed into $C$.
The known coefficients integrate to
\eqref{eq:FE11}. Discarding the last term $a_{p,L}h^{-L}$ proves
the asserted log error, on both the strip and its complement.
Choose $c_p$ smaller and $h_p$ larger once so that the inverse-cusp
remainder estimate remains valid through order $K+3$. The endpoint range is linear in $h$.

Finally the finite degree-$K$ log polynomial is uniformly bounded on
$|t|\le1/(D_0(K+1))$. Cauchy's exponential tail at $t=1/h$ is
$C[D_0(K+1)/h]^{K+1}\le C(\mathrm{e}D_0)^{K+1}(K+1)!h^{-K-1}$.
At the real evaluation point the polynomial is $O(h^{-1})$ and the
function $F_g$ itself is bounded by Proposition~\ref{prop:prefactor}.
Thus the exponential is Lipschitz on the required bounded interval
and transfers the error bound for the logarithm to \eqref{eq:FE12}.
Choosing $K=\lfloor c'_ph\rfloor$ after the final error base,
with $\mathfrak B_pc'_p<1$, gives the exponential upper error.
\end{proof}

\section{Coefficient extraction with inverse-genus errors}\label{app:extraction}

This appendix proves Theorems~\ref{thm:sharp-extraction}
and~\ref{thm:canonical-expansion} of Section~\ref{sec:sharp-extraction},
which extract a single coefficient of $\Phi_g^{(p)}$ with an error of
order $1/h$, and to all orders, uniformly in the coefficient index. The
strategy is the third idea of Subsection~\ref{subsec:strategy}: the
coefficient is a contour integral around the real saddle point, the
true generating function is replaced there by a finite marked expansion
transferred to a small complex neighborhood, and the resulting model
integral is evaluated against an exact Poisson factor rather than a
Gaussian. The inputs are the marked expansion of
Theorem~\ref{thm:full-expansion}, used only through its statement; the
positivity of Lemma~\ref{lem:convolution}, which bounds the modulus of
$\Phi_g^{(p)}$ on a circle by its value on the real axis and supplies the
Fourier bounds of Subsection~\ref{uniform-fourier-bounds-and-the-lattice-local-limit};
the endpoint expansions \eqref{eq:TL36} of Lemma~\ref{lem:parameter};
and, near zero cusp weight, Temme's change of variables~\cite{Temme}.
Subsection~\ref{sec:quantitative-local} makes the local limit theorem of
Section~\ref{extracting-fixed-cusp-counts} quantitative; its complex
comparison and Fourier estimates are reused afterwards.
Subsections~\ref{app:local-complex}--\ref{app:small-indices} prove
Theorem~\ref{thm:sharp-extraction}, and
Subsections~\ref{app:phase-change}--\ref{app:canonical-proof} prove
Theorem~\ref{thm:canonical-expansion} and compute its first coefficient.

\subsection{A quantitative local limit}\label{sec:quantitative-local}

We first make the local limit theorem \eqref{eq:TL47} quantitative. The argument transfers a real-axis error to a complex neighborhood and applies Fourier inversion; the positive convolution factor controls the remaining part of the integral. This gives an error bound by a negative power of $\min(h,\mu)$. The complex comparison and Fourier estimates established in the proof will also be used in Subsection~\ref{app:local-complex}, where retaining a Poisson factor gives an inverse-genus error uniformly over all coefficient indices.

\begin{proposition}\label{prop:quantitative-local}
Fix $p\ge1$, and use $h,\mu,\sigma$ from \eqref{eq:TL37}. There are
constants $C_p,\kappa_p>0$ such that, with $M=\min(h,\mu)$,
\begin{equation}
 \sup_{k\in\mathbb Z}
 \left|\sigma\mathbb P_x(N=k)-
 \frac{\mathrm{e}^{-(k-\mu)^2/(2\sigma^2)}}{\sqrt{2\pi}}\right|
 \le C_p M^{-\kappa_p},
 \label{eq:QL1}
\end{equation}
whenever $M$ is sufficiently large, uniformly for $0<x<x_c$.
\end{proposition}

\begin{proof}
Write $Z=(N-\mu)/\sigma$ and $\mathcal L(z)=\mathbb E_x \mathrm{e}^{zZ}$.
The integral identity \eqref{eq:TL38} and the convergent endpoint
expansions \eqref{eq:TL36} give constants $c,C_p>0$ for which
\begin{equation}
 \left|\log\mathcal L(u)-\frac{u^2}{2}\right|
 \le C_p M^{-1/2}(|u|+|u|^3),\qquad |u|\le c\sqrt M.
 \label{eq:QL2}
\end{equation}
To check the uniformity, note that near zero $\Psi,\mathscr D\Psi$ and
$\mathscr D^2\Psi$ are comparable to $x$, so the model Taylor remainder
is $O(|u|^3/\sqrt\mu)$ and the integrated error is
$O(|u|x/\sqrt\mu)$. On an interior compact interval both are
$O_p(h^{-1/2}(|u|+|u|^3))$. Near criticality,
$\sigma\asymp\sqrt h/\varepsilon$; for $|u|\le c\sqrt h$ the tilted
distance to the critical point stays comparable to $\varepsilon$.
The bounds $\Psi=O(\varepsilon^{-1})$ and
$\mathscr D^2\Psi=O(\varepsilon^{-3})$ give the same estimate.
Shrinking $c$ ensures the real moment generating function is defined
on the whole interval. As before, the error in \eqref{eq:TL30} is never
differentiated.

Put $T=\sqrt{\log M/80}$ and $f(z)=\mathcal L(z)-\mathrm{e}^{z^2/2}$.
The function is analytic on $|z|\le4T$, because $4T\le c\sqrt M$ for
all sufficiently large $M$ and the moment generating function is
defined for $|\Re z|\le c\sqrt M$, the range of \eqref{eq:QL2}. Positivity gives
$|\mathcal L(z)|\le\mathcal L(\Re z)$, whence \eqref{eq:QL2} implies
\begin{equation}
 \sup_{|z|\le4T}|f(z)|\le2\mathrm{e}^{9T^2}=B,\qquad
 \sup_{u\in[-T,T]}|f(u)|
 \le C_p M^{-1/2}(1+T)^3\mathrm{e}^{T^2}=b.
 \label{eq:QL3}
\end{equation}
For sufficiently large $M$, $b\le B$. In the disk $|z|<4T$ slit along
$[-T,T]$, the harmonic measure of the slit is at least $1/2$ for
$|z|\le T$. To check this, rescale $T$ to one and use
$G(w)=\log|w+\sqrt{w^2-1}|$, with the branch at infinity.
It vanishes on the slit, is at most $\log(1+\sqrt2)$ on $|w|\le1$,
and is at least $\log(4+\sqrt{15})$ on $|w|=4$.
The maximum principle compares harmonic measure to
$1-G/\log(4+\sqrt{15})$, which is at least $1/2$ on the inner disk.
Applying the two-constants principle to $\log|f|$ therefore yields
\begin{equation}
 |f(it)|\le\sqrt{bB}
 \le C_p M^{-1/4}(1+T)^{3/2}\mathrm{e}^{5T^2},\qquad |t|\le T.
 \label{eq:QL4}
\end{equation}
Its integral over this interval is $O_p(M^{-1/8})$.

The Fourier estimates of Subsection~\ref{uniform-fourier-bounds-and-the-lattice-local-limit}
give constants $c_0,c_1,C>0$ such that
\begin{equation}
 \int_{T<|t|\le\pi\sigma}|\mathbb E_x \mathrm{e}^{itZ}|\,\mathrm{d}t
 \le C\mathrm{e}^{-c_0T^2}+C\mathrm{e}^{-c_0h},\qquad 1\le T\le c_1\sqrt M.
 \label{eq:QL5}
\end{equation}
Away from a fixed critical neighborhood this follows from
\eqref{eq:TL40}: the product of the zero and one probabilities of
$V(xz)/V(x)$, divided by $\mathscr D\Psi(x)$, is bounded below.
Near zero both numerator and denominator are comparable to $x$.
In the critical neighborhood \eqref{eq:TL43} becomes
$(1+c t^2/h)^{-h/8}$. For $|t|\le\sqrt h$ this is bounded by
$\mathrm{e}^{-c't^2}$. Beyond $\sqrt h$, split the exponent in half, bound
one half at $\sqrt h$, and integrate the other against a fixed
integrable power. The integral is $O(\mathrm{e}^{-c'h})$.
The outer arc is exponentially small uniformly in $\varepsilon$ by
\eqref{eq:TL46}. These observations prove \eqref{eq:QL5}.

Fourier inversion, \eqref{eq:QL4}, \eqref{eq:QL5}, and the Gaussian
Fourier tail prove \eqref{eq:QL1}. The three contributions are
$O_p(M^{-1/8})$, $O(M^{-c_0/80})+O(\mathrm{e}^{-c_0h})$, and
$O(\mathrm{e}^{-T^2/2})=O(M^{-1/160})$, so one may take
$\kappa_p=\min(c_0,1/2)/80$. This exponent is not intended to be sharp.
\end{proof}

As a first application, and only as an illustration of
Proposition~\ref{prop:quantitative-local}, we make the volume-ratio
limits quantitative; the rates actually used in
Section~\ref{sec:absolute} are the sharper ones that follow from
Theorem~\ref{thm:sharp-extraction}, and \eqref{eq:QL6} is not used
elsewhere. Put $\kappa=\min(\kappa_1,\kappa_2,\kappa_3)$, so that the
exponent serves for each of the marking counts $p=1,2,3$. At the common
saddle \eqref{eq:TL48}, put $M=\min(2g-1,n-1)$. Repeating the
coefficient comparisons in \eqref{eq:TL49}--\eqref{eq:TL54} now gives
\begin{equation}
 B_{g,n}=\frac{\pi^2}{2r(Y)}(1+O(M^{-\kappa})),\qquad
 C_{g,n}=\frac{\pi^2(n-1)}{gx}(1+O(M^{-\kappa}))
 \label{eq:QL6}
\end{equation}
when $g,n\to\infty$. The center shift in \eqref{eq:TL52} contributes
only $O(g^{-1})$, since its Gaussian argument is $O(g^{-1/2})$.

The next three subsections prove Theorem~\ref{thm:sharp-extraction}. We
first compare the true and the model generating functions on small
complex neighborhoods of the real saddle point, and then treat large
and small coefficient indices separately.

\subsection{Local complex comparison}\label{app:local-complex}

To extract a coefficient by a contour integral we need the marked
expansion on complex neighborhoods of the real saddle, not only on the
real axis. We transfer a finite-order expansion to such neighborhoods and
compare the true and model Fourier integrals.

Choose a fixed finite order $K$, to be enlarged below, and define
\begin{equation}
 F_{h,[K]}(z)=\frac{\mathcal J_{g,p}(z)}\pi
       \left(1+\sum_{j=1}^K\frac{\mathcal L_{p,j}(z)}{h^j}\right).
 \label{eq:CE3}
\end{equation}
This is initially a local analytic function near the positive real
interval. We do not assume that it is analytic on a whole cusp-weight
disk. Theorem~\ref{thm:full-expansion} gives its real relative error
$O_{p,K}(h^{-K-1})$.

The convergent endpoint expansions and compactness give a fixed $c>0$
such that all functions in \eqref{eq:CE3} are analytic at
$z=x\exp(w/\sigma)$ for $|w|\le c\sqrt M$. Near zero this is an
ordinary inverse-function neighborhood. Near criticality the radius
in $z$ is proportional to $x_c-x$, and
$\sigma(x_c-x)/x\asymp\sqrt h$; the square-root inverse and the
analytic $\delta$ coefficients give the required branch. Each fixed
collection of $\mathcal L_{p,j}$ is bounded on a smaller such
neighborhood. The following lemma makes the analytic radius and the
resulting derivative bounds explicit, with constants that are geometric
in the order of differentiation. It is used for \eqref{eq:CE4}, for the
disk bound following \eqref{eq:CE4}, for \eqref{eq:CE5}, for the
remainder in \eqref{eq:CE8}, and for the disk bound in the proof of
Theorem~\ref{thm:canonical-expansion}.

\begin{lemma}[Normalized analytic radius]\label{rev:lem:normalized-radius}
Fix $p\ge1$. There are constants $c_0,c,C>0$ such that for every
$0<x<x_c$, with $\varepsilon=1-x/x_c$ and
\[
 \rho_x=c_0\min(1,\varepsilon),
\]
the functions
\[
 \phi_x(s)=S(x\mathrm{e}^{s})-S(x)-\Psi(x)s,\qquad
 \psi_x(s)=\log B_p(x\mathrm{e}^{s})-\log B_p(x)
\]
are analytic on $|s|\le\rho_x$ and satisfy
\begin{equation}\label{rev:eq:radius-sup}
 \sup_{|s|\le\rho_x}|\phi_x(s)|+\sup_{|s|\le\rho_x}|\psi_x(s)|\le Cx .
\end{equation}
Consequently, by Cauchy's estimate,
\begin{equation}\label{rev:eq:radius-cauchy}
 |\mathscr D^mS(x)|\le C\,m!\,x\,\rho_x^{-m}\quad(m\ge2),\qquad
 |\mathscr D^m\log B_p(x)|\le C\,m!\,x\,\rho_x^{-m}\quad(m\ge1).
\end{equation}
Moreover, with $k=h\Psi(x)$, $M=\min(h,k)$ and $\sigma^2=h\mathscr D\Psi(x)$,
\begin{equation}\label{rev:eq:radius-scale}
 \sigma\rho_x\ge c\sqrt M,\qquad
 c\,\sigma^2\le hx\rho_x^{-2}\le C\sigma^2 ,
\end{equation}
so that for $|w|\le\sigma\rho_x$
\begin{equation}\label{rev:eq:radius-disk}
 h\,|\phi_x(w/\sigma)|\le C|w|^2,\qquad |\psi_x(w/\sigma)|\le C|w| ,
\end{equation}
and, for $m\ge3$ and $m\ge1$ respectively,
\begin{equation}\label{rev:eq:radius-derivatives}
\begin{aligned}
 \frac{h|\mathscr D^mS(x)|}{\sigma^m}&\le C\,m!\,(\sigma\rho_x)^{2-m}
 \le C\,m!\,c^{2-m}M^{1-m/2},\\
 \frac{|\mathscr D^m\log B_p(x)|}{\sigma^m}&\le C\,m!\,(\sigma\rho_x)^{-m}
 \le C\,m!\,c^{-m}M^{-m/2}.
\end{aligned}
\end{equation}
\end{lemma}

\begin{proof}
Near zero, $S$ and $\log B_p$ are analytic on a fixed disk $|z|\le z_0$
and vanish at $z=0$, so $|S(z)|+|\log B_p(z)|\le C|z|$ there; for
$x\le z_0/\mathrm{e}$ and $|s|\le1$ one has $|x\mathrm{e}^s|\le\mathrm{e}x\le z_0$,
and $\Psi(x)\le Cx$, which gives \eqref{rev:eq:radius-sup} with
$\rho_x=1$. On a compact subinterval of $(0,x_c)$ the same holds with a
fixed radius, by compactness. Near $x_c$, write $\delta=y_0-Y$ and use
the convergent local factors
\[
 x_c-X(Y)=\delta^2\alpha(\delta),\qquad
 \mathcal Z(Y)=\delta^{5/2}\zeta(\delta),\qquad
 B_p(X(Y))=\delta^{p/2}\chi_p(\delta),
\]
with $\alpha,\zeta,\chi_p$ analytic and nonzero at zero; they
follow from the analytic endpoint parametrization constructed in the
proof of Lemma~\ref{lem:parameter}(iii) in
Appendix~\ref{app:parameter-proof}, $\mathcal Z=\sqrt2I$ and
\eqref{eq:CE1}, and are recorded again in
Subsection~\ref{app:canonical-proof}. For $z=x\mathrm{e}^s$ the parameter
$\delta_z=\delta q$ solves
$q^2\alpha(\delta q)=\alpha(\delta)-x(\mathrm{e}^s-1)\delta^{-2}$; since
$\partial_q[q^2\alpha(\delta q)]=2\alpha(0)\ne0$ at $(q,\delta)=(1,0)$,
the implicit function theorem gives $q$ analytic in $\delta$ and in
$\tau=x(\mathrm{e}^s-1)\delta^{-2}$ for $|\tau|\le c_1$, hence for
$|s|\le c_2\delta^2/x$; as $\varepsilon=\delta^2\alpha(\delta)/x_c$, this
is a disk $|s|\le c_0\varepsilon$ on which $q$ stays in a fixed
neighborhood of one. There
\[
 S(z)-S(x)=-\tfrac52\log q-\log\frac{\zeta(\delta q)}{\zeta(\delta)},\qquad
 \log B_p(z)-\log B_p(x)=\tfrac p2\log q+\log\frac{\chi_p(\delta q)}{\chi_p(\delta)}
\]
are bounded, and $\Psi(x)|s|\le C\varepsilon^{-1}\cdot c_0\varepsilon\le C$
by \eqref{eq:SP16}. This proves \eqref{rev:eq:radius-sup} near $x_c$,
where $Cx\asymp C$.

Since $\frac{\mathrm{d}^m}{\mathrm{d}s^m}S(x\mathrm{e}^s)\big|_{s=0}=\mathscr D^mS(x)$
and $\mathscr DS=\Psi$, the Taylor coefficients of $\phi_x$ at $s=0$ are
$\mathscr D^mS(x)/m!$ for $m\ge2$, and those of $\psi_x$ are
$\mathscr D^m\log B_p(x)/m!$ for $m\ge1$; Cauchy's estimate on
$|s|=\rho_x$ gives \eqref{rev:eq:radius-cauchy}. For
\eqref{rev:eq:radius-scale}, $\sigma^2=h\mathscr D\Psi$ with
$\mathscr D\Psi\asymp x$ near zero, $\asymp1$ on compact subintervals
and $\asymp\varepsilon^{-2}$ near $x_c$, by \eqref{eq:SP16}; hence
$hx\rho_x^{-2}\asymp hx\asymp\sigma^2$ near zero and
$hx\rho_x^{-2}\asymp h\varepsilon^{-2}\asymp\sigma^2$ near $x_c$, while
$\sigma\rho_x\asymp\sqrt{hx}\asymp\sqrt{h\Psi(x)}=\sqrt k\ge\sqrt M$ near
zero and $\sigma\rho_x\asymp\sqrt h\ge\sqrt M$ otherwise. Since $\phi_x$
vanishes to second order and $\psi_x$ to first order at $s=0$,
Schwarz's lemma on $|s|\le\rho_x$ gives
$|\phi_x(s)|\le Cx|s|^2\rho_x^{-2}$ and $|\psi_x(s)|\le Cx|s|\rho_x^{-1}$;
at $s=w/\sigma$ this is \eqref{rev:eq:radius-disk}, by
\eqref{rev:eq:radius-scale} and
$x/(\sigma\rho_x)\le C\sqrt{x/h}\le C$ for $h\ge1$. Finally
\eqref{rev:eq:radius-derivatives} follows from
\eqref{rev:eq:radius-cauchy} and \eqref{rev:eq:radius-scale}, since
$hx\rho_x^{-m}\sigma^{-m}=(hx\rho_x^{-2}\sigma^{-2})(\sigma\rho_x)^{2-m}$
and $x\rho_x^{-m}\sigma^{-m}\le x_c(\sigma\rho_x)^{-m}$.
\end{proof}

The model derivatives obey
\begin{equation}
 \frac{h\mathscr D^mS}{\sigma^m}=O_m(M^{1-m/2})\quad(m\ge3),
 \qquad
 \frac{\mathscr D^m\log B_p}{\sigma^m}=O_{p,m}(M^{-m/2})\quad(m\ge1),
 \label{eq:CE4}
\end{equation}
with constants of the form $C\,m!\,c^{-m}$, by
\eqref{rev:eq:radius-derivatives}; after division by the factor $m!$
of the Taylor coefficients this dependence on $m$ is geometric, which is
what the disk bounds below require.
At zero use $\mathscr D^mS=2x+O_m(x^2)$ and
$\mathscr D^m\log B_p=O_{p,m}(x)$. At criticality, with
$\varepsilon=1-x/x_c$, both numerator derivatives before their
$h$ factor are $O(\varepsilon^{-m})$, whereas
$\sigma\asymp\sqrt h/\varepsilon$.
These bounds also hold locally in the complex neighborhoods by the
same convergent expansions.

Put
\[
 \mathcal L_h(w)=\mathrm{e}^{-kw/\sigma}
   \frac{F_h(x\mathrm{e}^{w/\sigma})}{F_h(x)},\qquad
 \mathcal L_{h,[K]}(w)=\mathrm{e}^{-kw/\sigma}
   \frac{F_{h,[K]}(x\mathrm{e}^{w/\sigma})}{F_{h,[K]}(x)}.
\]
Since $k=h\Psi(x)$, the logarithm of the finite analytic model is
\[
 \log\mathcal L_{h,[K]}(w)=h\phi_x(w/\sigma)+\psi_x(w/\sigma)
 +\log\frac{P_K(x\mathrm{e}^{w/\sigma})}{P_K(x)},
 \qquad P_K=1+\sum_{j=1}^K\mathcal L_{p,j}h^{-j},
\]
in the notation of Lemma~\ref{rev:lem:normalized-radius}; the last
ratio lies between $1/2$ and $2$ in modulus on the chosen neighborhoods
for large $h$. Hence \eqref{rev:eq:radius-disk} gives the disk bound
$C_{p,K}\mathrm{e}^{C|w|^2+C|w|}$, for $|w|\le c\sqrt M$, for the centered,
normalized finite analytic model, with $C$ independent of $K$; the
centered, normalized true generating function obeys the same bound,
because positivity bounds its modulus by its value at $\Re w$, where the
real expansion applies. No error term is differentiated. The explicit model also satisfies
\begin{equation}
 \left|\frac{\mathcal J_{g,p}(x\mathrm{e}^{it/\sigma})}
                    {\mathcal J_{g,p}(x)}\right|
       \le C_p \mathrm{e}^{-c't^2},\qquad |t|\le c_*\sqrt M,
 \label{eq:CE5}
\end{equation}
for a suitable $c_*>0$ and all sufficiently large $M$. Indeed the real action is
\[
 \Re\,h\phi_x(it/\sigma)=-\frac{t^2}2
 +\Re\sum_{m\ge4}\frac{h\mathscr D^mS\,(it)^m}{m!\,\sigma^m},
\]
the odd terms being purely imaginary, and by
\eqref{rev:eq:radius-derivatives} the sum is at most
$CM\sum_{m\ge4}(|t|/(c\sqrt M))^m\le C_1t^4/M$ for $|t|\le c\sqrt M/2$,
with $C_1=2C/c^4$;
the real amplitude change $\Re\,\psi_x(it/\sigma)$ is at most
$C\sum_{m\ge2}(|t|/(c\sqrt M))^m\le C_1t^2/M$ on the same range.
For $|t|\le c_*\sqrt M$, with $c_*\le c/2$ and $C_1c_*^2\le1/8$, the
former is at most $t^2/8$, and for $M\ge8C_1$ so is the latter; this
gives \eqref{eq:CE5} with $c'=1/4$.

Take $T^2=A_0\log h$, with $A_0$ large enough that the Fourier
tail \eqref{eq:QL5} of the true distribution, the polynomial-Gaussian tails below, and the
sparse-model Fourier tail below are $O(h^{-3})$. These tail
constants depend on the fixed action $S$ and amplitude $B_p$, not
on the finite truncation order.
The central truncation order will first be fixed below; then choose $C_0$
so that $M\ge C_0\log h$ implies both $4T\le c\sqrt M$ for its
analytic neighborhood and $T\le\min(c_1,c_*)\sqrt M$, the ranges of
\eqref{eq:QL5} and \eqref{eq:CE5}.

\begin{lemma}[Local comparison]\label{lem:local-comparison}
There are a fixed order $K_c$ and a constant $C_0$ such that, uniformly
for $M\ge C_0\log h$,
\[
 \sigma\mathbb P_x(N=k)=\frac1{2\pi}\int_{-T}^{T}\mathcal L_{h,[K_c]}(it)\,\mathrm dt
 +O(h^{-2}),\qquad T^2=A_0\log h.
\]
\end{lemma}

\begin{proof}
Here $k\ge1$ is an integer and $x$ is its saddle, $h\Psi(x)=k$, so
$\mathcal L_h(it)$ is the characteristic function of $(N-k)/\sigma$.
Fourier inversion gives
$\sigma\mathbb P_x(N=k)=(2\pi)^{-1}\int_{-\pi\sigma}^{\pi\sigma}\mathcal L_h(it)\,\mathrm dt$,
and by \eqref{eq:QL5} and the choice of $A_0$ the part
$T<|t|\le\pi\sigma$ is $O(h^{-3})$. It remains to compare the two
integrands on $[-T,T]$.
The difference between the two centered, normalized functions has real-slit bound
$C h^{-K-1}\mathrm{e}^{CT^2}$ on $[-T,T]$ and disk bound $C\mathrm{e}^{C'T^2}$.
The harmonic-measure argument of \eqref{eq:QL4} then gives
\begin{equation}
 \sup_{|t|\le T}|\mathcal L_h(it)-\mathcal L_{h,[K]}(it)|
       \le C_{p,K}h^{-(K+1)/2+C''A_0}.
 \label{eq:CE6}
\end{equation}
The exponential constants here are independent of $K$: for each fixed
$K$, the modulus of its finite correction factor lies between $1/2$ and $2$ on the
chosen neighborhoods for all sufficiently large $h$. Only the threshold
and multiplicative constants depend on $K$. Fix the central order $K=K_c$ sufficiently large
that the integral of \eqref{eq:CE6} is $O(h^{-2})$, and then fix $C_0$.
Thus $\sigma\mathbb P_x(N=k)$ equals the centered finite-model integral
on $[-T,T]$, divided by $2\pi$, up to $O(h^{-2})$.
The normalization $F_{h,[K]}(x)/F_h(x)=1+O(h^{-K-1})$ is retained
when recovering coefficients. All orders chosen here are fixed
independently of $g,x,k$.
\end{proof}

\subsection{Cancellation of the Poisson singularity}

For indices $k\ge h^{1/4}$, so that $M=\min(h,k)\ge h^{1/4}$, we compute
the model coefficient by the saddle-point method. At order $h^{-j}$ its Gaussian coefficient has a leading pole
$s_j\Psi^{-j}$ as the cusp weight tends to zero, matching the
corresponding Poisson coefficient. Division by the Poisson factor removes
these leading poles; the weaker residual bounds in~\eqref{eq:CE10}
suffice for the error in~\eqref{eq:CE2}.

\begin{lemma}[Large indices]\label{lem:large-indices}
The estimate \eqref{eq:CE2} holds uniformly for integers $k\ge h^{1/4}$.
\end{lemma}

\begin{proof}
Here $M=\min(h,k)\ge h^{1/4}$. Put $b=\log B_p$. For a formal
variable $u$ define the model polynomials through the exponential
\begin{equation}
 \exp\!\left\{
 \sum_{m\ge3}u^{m-2}\frac{\mathscr D^{m-1}\Psi}{m!v^{m/2}}w^m
 +\sum_{m\ge1}u^m\frac{\mathscr D^mb}{m!v^{m/2}}w^m\right\}.
 \label{eq:CE7}
\end{equation}
Let $E_j(x)$ be the Gaussian Fourier integral of its $u^{2j}$
coefficient: replace $w^{2a}$ by $(-1)^a(2a-1)!!$ and odd powers
by zero. Each coefficient uses only finitely many derivatives.
By Lemma~\ref{lem:local-comparison} it suffices to evaluate the central
model integral. Expansion under the Gaussian bound \eqref{eq:CE5} gives
\begin{equation}
 \frac1{2\pi}\int_{-T}^T \mathrm{e}^{-ikt/\sigma}
       \frac{\mathcal J_{g,p}(x\mathrm{e}^{it/\sigma})}{\mathcal J_{g,p}(x)}\,\mathrm{d}t
 =\frac1{\sqrt{2\pi}}\left(1+\sum_{j=1}^3\frac{E_j(x)}{h^j}
                                  +O_p(M^{-4}+h^{-2})\right).
 \label{eq:CE8}
\end{equation}
To justify the stated remainder, retain weights through seven in
\eqref{eq:CE7}; a term of weight $q$ has parity $q$ in $w$, so odd
weights integrate to zero. By \eqref{rev:eq:radius-derivatives}, a
factor of weight $m-2$ in the action sum of \eqref{eq:CE7} is at most
$Cc^{2-m}M^{1-m/2}|t|^m$ and a factor of weight $m$ in the amplitude sum
is at most $Cc^{-m}M^{-m/2}|t|^m$, geometrically in $m$; so the terms of
total weight $q$ are at most $C^qM^{-q/2}(1+|t|)^{3q}$, and their sum
over $q\ge8$ is at most $2(C(1+|t|)^3/\sqrt M)^8$ as long as
$C(1+|t|)^3\le\sqrt M/2$. Thus the remaining terms are bounded by
$M^{-4}$ times a fixed polynomial in $|t|$ times $\mathrm{e}^{-c't^2}$.
Their integral is $O(M^{-4})$. The cutoff is $O(\sqrt{\log h})$, which is smaller than
any fixed positive power of $M$ in this range; the Taylor remainder
and exponential remainder bounds therefore apply uniformly. The
removed polynomial-Gaussian tails are $O(h^{-2})$ by the initial
choice of $A_0$. Replacing the finite model by $\mathcal J$ in its central
integral costs $O_p(h^{-1})$, by boundedness of its corrections and
\eqref{eq:CE5}.

The coefficients $E_j$ are bounded at criticality by
\eqref{eq:CE4}. Near zero, $v=\Psi+O(\Psi^2)$ and
$\mathscr D^m\Psi=\Psi+O_m(\Psi^2)$. In a monomial of weight $2j$
with $a$ action factors and $b_*$ amplitude factors, the total degree
in $w$ is $2j+2a$. Its leading power in $\Psi$ is therefore
$\Psi^{b_*-j}$. Every amplitude factor, or higher term in an action
coefficient, improves the pole by at least one. Consequently
\begin{equation}
 E_j(x)=s_j\Psi(x)^{-j}+O_{p,j}(\Psi(x)^{1-j})
       \quad(x\downarrow0),
 \label{eq:CE9}
\end{equation}
where $s_j$ is the coefficient of $k^{-j}$ in
$\sqrt{2\pi k}\,\mathfrak p_k$. Indeed the leading action terms in
\eqref{eq:CE7} are exactly those of a centered Poisson law: for the
exact Poisson generating function $\mathrm e^{hz}$ one has $\Psi(z)=z$,
$v=\Psi$ and $b=0$, the coefficient is $h^k/k!$, each $E_j$ equals
$e_j\Psi^{-j}$ by homogeneity, and \eqref{eq:CE8}, carried to any fixed
order, reads
$\sqrt{2\pi k}\,\mathfrak p_k=1+\sum_{j\le J}e_jk^{-j}+O(k^{-J-1})$,
so that $e_j=s_j$.
Parity makes the remaining Laurent powers integral.
For example, direct Gaussian integration gives
\[
 E_1=-\frac{\mathscr D^2b+(\mathscr Db)^2}{2v}
      +\frac{(\mathscr Db)\mathscr D^2\Psi}{2v^2}
      +\frac{\mathscr D^3\Psi}{8v^2}
      -\frac{5(\mathscr D^2\Psi)^2}{24v^3}.
\]
The amplitude terms are needed because $k=h\Psi$ is a model center,
not an assertion about the exact finite-genus mean.

Divide \eqref{eq:CE8} by the central Poisson expansion
$\sqrt{2\pi k}\,\mathfrak p_k=1+\sum_{j=1}^3s_j/(h\Psi)^j+O(k^{-4})$.
The principal pole cancels at each order. The quotient coefficients
$R_j$ satisfy
\begin{equation}
 |R_j(x)|\le C_{p,j}(1+\Psi(x)^{1-j}),\qquad j=1,2,3.
 \label{eq:CE10}
\end{equation}
This follows inductively: a product of a Poisson pole of order $i$
and a previous quotient pole of order at most $j-i-1$ has order
at most $j-1$. The quotient is bounded at the critical end as well.
Since $k=h\Psi\ge1$, every $h^{-j}R_j$ is $O_p(h^{-1})$.
Also $M^{-4}\le h^{-1}$. Combining this with \eqref{eq:CE6}
and the bounded finite amplitude corrections proves
\eqref{eq:CE2} in this range.
\end{proof}

\subsection{Small coefficient indices}\label{app:small-indices}

For $k\le h^{1/4}$ the preceding finite-order Gaussian estimate does
not give the required uniform error $O(h^{-1})$. We expand the model coefficient directly around the Poisson
coefficient of $\mathrm e^{\mathfrak s_1hz}$.

\begin{lemma}[Small indices]\label{lem:small-indices}
The estimate \eqref{eq:CE2} holds uniformly for integers $0\le k\le h^{1/4}$.
\end{lemma}

\begin{proof}
Near zero write $S(z)=\mathfrak s_1z+\mathfrak s_2z^2+O(z^3)$, $\mathfrak s_1=2$, and
$B_p(z)=1+d_1z+O(z^2)$. The full model \eqref{eq:CE3} is analytic
on a fixed neighborhood of zero. For $k\le h^{1/4}$, coefficient
majorants give
\begin{equation}
 \frac{k![z^k]F_{h,[K]}(z)}
      {(\mathcal J_{g,p}(0)/\pi)(\mathfrak s_1h)^k}
 =1+\frac{(\mathfrak s_2/\mathfrak s_1^2)k(k-1)+(d_1/\mathfrak s_1)k+\mathcal L_{p,1}(0)}h
       +O_{p,K}\!\left(\frac{(k+1)^4}{h^2}\right).
 \label{eq:CE11}
\end{equation}
Factor $\mathrm{e}^{\mathfrak s_1hz}$ from the model. Its retained multiplier is
$1+h\mathfrak s_2z^2+d_1z+\mathcal L_{p,1}(0)/h$.
Analytic absolute-coefficient majorants bound the omitted terms by
convergent products starting with
$hz^3,h^2z^4,z^2,z/h,h^{-2}$. In convolution with the exponential,
a degree-$a$ term contributes the factor $(k)_a/(\mathfrak s_1h)^a$, bounded
by $(k/(\mathfrak s_1h))^a$. Evaluation of the majorants at that last argument
gives the remainder in \eqref{eq:CE11}; their exponent is
$O(k^2/h)=o(1)$. This avoids a spurious $\sqrt k$ loss from a
supremum Cauchy estimate.

At the saddle, direct expansion gives the same polynomial without
$\mathcal L_{p,1}(0)$ for the right-hand side of \eqref{eq:CE2},
relative to this Poisson baseline. Indeed,
$hS-k+k\log(k/(\mathfrak s_1hx))$ contributes $\mathfrak s_2k^2/(\mathfrak s_1^2h)$,
$\log B_p$ contributes $d_1k/(\mathfrak s_1h)$, and
$\tfrac12\log(k/(hv))$ contributes $-\mathfrak s_2k/(\mathfrak s_1^2h)$.
Their error is $O_p((k+1)^4/h^2)$. Thus the model coefficient divided
by the proposed saddle expression is $1+O_p(h^{-1})$ for
$k\le h^{1/4}$.

For $C_0\log h\le k\le h^{1/4}$, comparison
\eqref{eq:CE6} applies. Here the entire circle $|z|=x$ is inside
the fixed zero-neighborhood. The convergent series for $S$ and
$|1-\cos(m\theta)|\le m^2(1-\cos\theta)$ give
\[
 \Re(S(x\mathrm{e}^{i\theta})-S(x))\le-cx(1-\cos\theta),
 \qquad |\theta|\le\pi.
\]
The model amplitude is bounded on this circle, and $hx\asymp k$,
so its standardized Fourier tail beyond $T$ is $O(h^{-2})$ by the
initial choice of $A_0$. Complete the central model integral to its full
coefficient integral and apply \eqref{eq:CE11}. This proves
\eqref{eq:CE2} in the middle sparse range. The scaled model central mass is bounded below by a
positive constant here, by \eqref{eq:CE11} and the saddle comparison;
therefore the additive central-integral errors remain relative errors.

For $0\le k\le C_0\log h$ use a different contour and, if needed,
a different fixed truncation order $K=K_\ell$, chosen after $C_0$. Put
$r_k=(k+1)/(2h)$ and consider
$(F_h(r_k z)-F_{h,[K]}(r_k z))/(\mathcal J_{g,p}(0)/\pi)$
in the disk $|z-1|<16$ slit along $[0,2]$.
Both functions are analytic there for all large $h$; only a fixed
neighborhood of zero is used. Positivity, the real expansion, and
the explicit model give disk bound $C\mathrm{e}^{C'(k+1)}$ and slit bound
$C_{p,K}h^{-K-1}\mathrm{e}^{C'(k+1)}$.
On $|z|=1$ the slit harmonic measure is greater than $1/2$:
translate $w=z-1$, so $|w|\le2$, and compare the Green function
in the proof of \eqref{eq:QL4} on radii two and sixteen.
Its lower bound is
\[
 1-\frac{\log(2+\sqrt5)}{\log(16+\sqrt{255})}>\frac12.
\]
Thus the unit-circle error is at most
$C_{p,K}h^{-(K+1)/2}\mathrm{e}^{C''(k+1)}=O(h^{-2})$ upon fixing $K$
sufficiently large. Cauchy extraction changes the coefficient,
relative to the baseline of \eqref{eq:CE11}, by at most
\[
 Ch^{-2}\frac{k!}{(2hr_k)^k}
       =Ch^{-2}\frac{k!}{(k+1)^k}\le Ch^{-2}.
\]
Equation~\eqref{eq:CE11} and its saddle comparison finish this
range, including $k=0$. The two finite orders $K_c,K_\ell$ are fixed before $h$ tends to
infinity, and each is used only in its own range. Enlarging the
low-mean order does not change the already fixed central analytic
neighborhood or $C_0$.
\end{proof}

\begin{proof}[Proof of Theorem~\ref{thm:sharp-extraction}]
Lemmas~\ref{lem:large-indices} and~\ref{lem:small-indices} cover all
$k\ge0$; taking common constants and a common threshold $h_p$ proves
the uniform assertion.
\end{proof}

\subsection{The analytic phase change at zero}\label{app:phase-change}

To obtain an expansion that is uniform down to $k=0$, we change variables
so that the phase becomes exactly that of a Poisson coefficient, following
Temme~\cite{Temme}. This subsection and the next supply the two lemmas
used in the small-density part of the proof of
Theorem~\ref{thm:canonical-expansion} in
Subsection~\ref{app:canonical-proof}: Lemma~\ref{lem:phase-change}
constructs the change of variables and records the resulting coefficient
identity \eqref{eq:CG7}, and Lemma~\ref{lem:poisson-central} expands the
transformed amplitude with a remainder that is uniform in the index $k$
and factorial in the order. Only the local analyticity of the action $S$
and the amplitude near zero is used.

\begin{lemma}[Analytic phase change]\label{lem:phase-change}
Let $S$ be analytic near zero with $S(0)=0$ and $S'(0)=s>0$, and let
$x=x(\nu)$ be the solution of $xS'(x)=\nu$ near zero, so that
$x=\nu/s+O(\nu^2)$. There is a jointly analytic map $z=Z(u,\nu)$ on a fixed
polydisk, with $Z(0,\nu)=0$, such that
\begin{equation}
 S(z)-\nu\log z=u-\nu\log u+\mathcal H(\nu),\qquad
 \mathcal H(\nu)=S(x)-\nu+\nu\log(\nu/x),\qquad Z(\nu,\nu)=x.
 \label{eq:CG5}
\end{equation}
If $H$ is analytic near zero and
$f(u,\nu)=H(Z(u,\nu))\,uZ_u(u,\nu)/Z(u,\nu)$, then $f$ is jointly analytic,
including at $u=0$, and for all sufficiently small $\nu=k/h$,
\begin{equation}
 [z^k]\mathrm{e}^{hS(z)}H(z)
       =\mathrm{e}^{h\mathcal H(\nu)}[u^k]\mathrm{e}^{hu}f(u,\nu).
 \label{eq:CG7}
\end{equation}
\end{lemma}

\begin{proof}
The logarithmic terms in \eqref{eq:CG5} are to be read as
$-\nu\log z+\nu\log u=-\nu\log(Z/u)$: the quotient $Z(u,\nu)/u$
constructed below is analytic and nonvanishing on the polydisk, so
$\log(Z/u)$ is a single-valued analytic function there, and no branch
of $\log z$ or of $\log u$ separately is needed. We construct $Z$
directly. Put $\mathfrak c=S'(x)$ at the saddle $x=x(\nu)$,
$xS'(x)=\nu$; thus $\mathfrak c=\nu/x$ for $\nu\ne0$, the function
$\mathfrak c$ is analytic near $\nu=0$ with $\mathfrak c(0)=S'(0)$, and
$\mathcal H(\nu)=S(x)-\nu+\nu\log\mathfrak c$ is analytic with
$\mathcal H(0)=0$. Put $d=z-x$ and
$\phi=S(z)-\mathfrak cz+\nu\log\mathfrak c-\mathcal H$.
Then $\phi(x,\nu)=\phi_z(x,\nu)=0$, and its double-zero quotient is
$\int_0^1(1-t)S''(x+td)\,\mathrm{d}t$, jointly analytic.
Seek $u=z(\mathfrak c+d\eta)$. The phase equation becomes
\begin{equation}
 \frac{\phi}{d^2}-\eta+
 \nu\sum_{m\ge2}\frac{(-1)^{m-1}}m \mathfrak c^{-m}d^{m-2}\eta^m=0.
 \label{eq:CG6}
\end{equation}
At $z=\nu=0$ it is $S''(0)/2-\eta=0$, with derivative $-1$ in
$\eta$. The analytic implicit-function theorem gives $\eta(z,\nu)$.
Since $u_z(0,0)=s$, the inverse $Z$ is jointly analytic on a fixed
polydisk. The construction gives $u(x,\nu)=\nu$ and $u(z,0)=S(z)$.
Only local analyticity of $S$ is required.

The function $f$ is jointly analytic because $u/Z=\mathfrak c+d\eta$ is
analytic and nonvanishing, its value at $z=\nu=0$ being $s$. For \eqref{eq:CG7}, write the coefficient as a residue and
change variables $z=Z(u,\nu)$. For $\nu>0$ sufficiently small, the
circle $|u|=\nu$ lies in the fixed analytic neighborhood, and its image
surrounds zero once. At $k=\nu=0$, both sides equal $H(0)$, since
$S(0)=\mathcal H(0)=0$ and $f(0,0)=H(0)$.
\end{proof}

At the saddle, differentiating \eqref{eq:CG5} twice in $u$ gives
$(S''(x)+\nu/x^2)Z_u(\nu,\nu)^2=1/\nu$. With
$v=\mathscr D^2S(x)=\nu+x^2S''(x)$ this is
$(\nu Z_u(\nu,\nu)/x)^2=\nu/v$, and $Z_u>0$ near the origin, so
$f(\nu,\nu)=H(x)\sqrt{\nu/v}$, with the positive branch.
For $H=B_p$ the normalization on the right is exactly
$\mathrm{e}^{h\mathcal H}h^k B_p(x)\sqrt{\nu/v}/k!$,
which equals the $x$-dependent part of \eqref{eq:CG1}.

\subsection{Poisson central polynomials with a uniform remainder}

In the new variable the coefficient is a Poisson coefficient of an
analytic amplitude. Expanding the amplitude about $u=k/h$ expresses it
through the central polynomials $T_m(k)$ below; we bound the remainder
uniformly in $k$.

Define the normalized coefficient operator by
\[
 \mathscr P_{h,k}f=k!h^{-k}[u^k]\mathrm{e}^{hu}f(u,\nu),\qquad \nu=k/h.
\]
Its central polynomial moments are
\begin{equation}
 \mathscr P_{h,k}(u-\nu)^m=h^{-m}T_m(k),\qquad
 \sum_{m\ge0}T_m(k)\frac{t^m}{m!}=\mathrm{e}^{-kt}(1+t)^k.
 \label{eq:CG8}
\end{equation}
Thus $T_0=1$, $T_1=0$; for $m\ge2$, $T_m$ has degree at most
$\lfloor m/2\rfloor$ and zero constant term. Write
$T_m(k)=\sum_{j=\lceil m/2\rceil}^{m-1}t_{m,j}k^{m-j}$.

\begin{lemma}[Poisson central expansion]\label{lem:poisson-central}
Let $f(u,\nu)$ be jointly analytic and bounded by one on a fixed polydisk.
There are $\nu_0,C>0$ such that, uniformly for $0\le\nu=k/h\le\nu_0$ and
$0\le K\le h/2$,
\begin{equation}
 \begin{split}
 \mathscr P_{h,k}f&=\sum_{j=0}^Kh^{-j}\mathscr T_j[f](\nu)
       +O\!\left(C^{K+1}(K+1)!h^{-K-1}\right),\\
 \mathscr T_0[f]&=f(\nu,\nu),\qquad
 \mathscr T_j[f]=\sum_{m=j+1}^{2j}
 t_{m,j}\nu^{m-j}\frac{\partial_u^mf(\nu,\nu)}{m!}\quad(j\ge1).
 \end{split}
 \label{eq:CG10}
\end{equation}
Moreover $\|\mathscr T_j[f]\|\le C^{j+1}j!$ on a smaller complex disk in
$\nu$, each $\mathscr T_j[f]$ with $j\ge1$ vanishes at $\nu=0$, and
$|\mathscr P_{h,k}f|\le C$.
\end{lemma}

\begin{proof}
Choose $\nu_0>0$ small and a fixed Cauchy radius $r_0$ about every
$u=\nu$, with $2\nu_0<r_0/2$. The Taylor remainder through degree
$2K+1$ on $|u|=\nu$ is at most
$2r_0^{-2K-2}|u-\nu|^{2K+2}$.
For $k\ge1$, normalized coefficient extraction therefore bounds it by
\begin{equation}
 \begin{split}
 &Cr_0^{-2K-2}\mathfrak p_k^{-1}\nu^{2K+2}
   \int_{-\pi}^{\pi}|\theta|^{2K+2}\mathrm{e}^{-ck\theta^2}\,\mathrm{d}\theta\\
 &\hspace{20pt}\le C^{K+1}(K+1)!\nu^{K+1}h^{-K-1}
       \le C_1^{K+1}(K+1)!h^{-K-1}.
 \end{split}
 \label{eq:CG9}
\end{equation}
Here $\mathfrak p_k^{-1}\le C\sqrt{k+1}$ and
$1-\cos\theta\ge c\theta^2$ on $[-\pi,\pi]$.
The Gaussian integral is an upper bound at every order, even when
$K$ exceeds $k$. At $k=0$ the operator is evaluation at zero and
this remainder vanishes.

The retained moment terms above inverse order $K$ must also be bounded.
The generating function in \eqref{eq:CG8} gives
\[
 |[k^a]T_m(k)|\le\frac{m!}{a!}
 [t^m]\left(\sum_{r\ge2}\frac{t^r}{r}\right)^a
 \le2^m\frac{m!}{a!},\qquad 1\le a\le m/2.
\]
For $\ell=m-a=K+1+d>K$, the retained range has
$0\le d\le K-1$, $m\le2K+1$, and
\[
 \frac{m!}{(m-\ell)!}\le\frac{(2K+1)!}{(K-d)!}
       \le C^K(K+1)!K^d.
\]
The sum of $(K/h)^d$ is bounded when $K\le h/2$; the polynomial
number of terms is absorbed into $C^K$. At $K=0$ this sum is empty.
This proves \eqref{eq:CG10}. At order $j$, $m!/(m-j)!\le(2j)!/j!\le4^jj!$.
Cauchy's bound thus gives
$\|\mathscr T_j[f]\|\le C^{j+1}j!$ on one smaller complex
$\nu$ disk. Each positive-order coefficient vanishes at zero.
The same contour estimate with no vanishing factor gives
$|\mathscr P_{h,k}f|\le C\|f\|$.
\end{proof}

Apply Lemmas~\ref{lem:phase-change} and~\ref{lem:poisson-central} to
$f_s(u,\nu)=B_p(Z(u,\nu))\mathcal L_{p,s}(Z(u,\nu))uZ_u/Z$.
Put $f_0^*(\nu)=B_p(x)\sqrt{\nu/v}$, which is bounded away from
zero on a fixed small complex disk, with $f_0^*(0)=1$.
The extracted coefficients there are
\begin{equation}
 \mathcal A_{p,j}(Y)=\frac1{f_0^*(\nu)}
       \sum_{s=0}^j\mathscr T_{j-s}[f_s](\nu),\qquad \nu=\Psi(X(Y)).
 \label{eq:CG11}
\end{equation}
This proves their zero-end analyticity and endpoint values.
The common norms of Theorem~\ref{thm:full-expansion} give
$\|f_s\|\le C_p(\mathfrak B_p^{\mathrm{FE}})^ss!$ on one polydisk. Applying
\eqref{eq:CG10} at order $K-s$ uses
$s!(K-s+1)!\le(K+1)!$.
For a finite marked model of order $J=\lfloor\kappa h\rfloor>K$,
the bounded operator controls all unused amplitudes by their first
geometric tail term, since $\mathfrak B_p^{\mathrm{FE}}(s+1)/h<1/2$ through $J$.
Thus both model errors and coefficient norms have the required Gevrey
bound. The functions \eqref{eq:CG11} agree with \eqref{eq:CG3} where
both are defined; this is proved at the end of the proof of
Theorem~\ref{thm:canonical-expansion}, from the two expansions on a
common range of densities.

\Needspace{7\baselineskip}
\subsection{Proof of the uniform coefficient expansion}\label{app:canonical-proof}

This subsection proves Theorem~\ref{thm:canonical-expansion} and then
computes the first extracted coefficient $\mathcal A_{p,1}$, see
\eqref{eq:CG12}, together with its value \eqref{eq:CG13} at the critical
endpoint. The inputs are Theorem~\ref{thm:full-expansion} over its whole
linear order range, which provides one finite marked model with
exponentially small relative error on the real axis; the Fourier bounds
\eqref{eq:QL5} and the two-constants argument of \eqref{eq:QL4}, which
transfer this error to the contour; Lemma~\ref{rev:lem:normalized-radius}
for the disk bounds; and Lemmas~\ref{lem:phase-change}
and~\ref{lem:poisson-central} for the Poisson range. The proof is marked
in stages: the order in which the constants are fixed, the comparison of
the true coefficients with those of the model, the Poisson range of small
densities, uniform Gaussian coordinates through the critical point, the
Gaussian expansion, the exact Poisson normalization, the identification
of the two coefficient families on their overlap, and the absorption of
the residual exponential error.

\begin{proof}[Proof of Theorem~\ref{thm:canonical-expansion}]
We first replace the true generating function by one finite marked
expansion, with exponentially small error after coefficient extraction.
For small cusp density the Poisson calculation above applies. For the
remaining densities we use uniform Gaussian coordinates and divide by
the exact Poisson normalization. The last step absorbs the exponential
comparison error into the stated remainder at every allowed order.

\emph{Order of the constants.} The constants below are chosen in the
order
\[
 \kappa\ \to\ \nu_0\ \to\ \nu_*\ \to\ d\ \to\ \eta\ \to\ \tau\ \to\ a
 \ \to\ c_p\ \to\ c'_p ,
\]
each depending only on $p$ and on those before it. The proportion
$\kappa$ fixes the model order $J=\lfloor\kappa h\rfloor$ and satisfies
$\kappa\mathfrak B_p^{\mathrm{FE}}<1/2$, where we write
$\mathfrak B_p^{\mathrm{FE}}$ for the Gevrey base of
Theorem~\ref{thm:full-expansion} throughout this proof, to distinguish
it from the error base $\mathfrak B_p$ of the present theorem; we also
take $\kappa$ at most the linear-order constant of that theorem, so that
it applies at order $J$. The proportion $\kappa$ determines the exponential rate
$a_0$ of the model. The radius $\nu_0$ is that of
Lemma~\ref{lem:poisson-central} for the amplitudes $f_s$ of
\eqref{eq:CG11}, and $0<\nu_*<\nu_0$ is the switching density. The
small-disk radius $d$ is chosen so that the two-constants bound on
$|z|\le d$ is exponentially small, then $\eta$ so that the saddle for
$k<\eta h$ lies in $|z|\le d/4$ and
\[
 \eta\le\nu_*/2,
\]
so that the Gaussian comparison \eqref{eq:CG14}, proved for
$k\ge\eta h$, covers every density $\Psi\ge\nu_*/2$. The Fourier
cutoff $T=\tau\sqrt h$ is fixed after $\eta$ and $a_0$, and $a>0$
denotes the smallest of the exponential rates obtained. Two further
smallness conditions are imposed where they arise in the Gaussian
range: $\tau$ is small compared with the radius of the common disk of
the Gaussian coordinates, and $c_p$ is small compared with that radius
and with $\tau$; neither depends on a constant chosen later. Finally
$c_p$ satisfies
\[
 c_p<\kappa/2,\qquad c_p\le c_R\,\nu_*/2,\qquad 2c_p<1,\qquad
 2c_p\log\frac{\mathrm{e}}{2c_p}<a,
\]
where $c_R$ is the constant in the Gevrey remainder of $R(k)$ below,
which requires $K\le c_Rk$ while $k\ge(\nu_*/2)h$ on the Gaussian range;
and $c'_p<c_p$ satisfies $\mathfrak B_pc'_p<1$ for the error base
$\mathfrak B_p$ of the present theorem.

\textbf{The finite model.} Write $F_h=\Phi_g^{(p)}$ and $J_h=\mathcal J_{g,p}/\pi$.
Choose a sufficiently small fixed $\kappa>0$, set
$J=\lfloor\kappa h\rfloor$, and use the single finite model
\[
 M_h(z)=J_h(z)P_J(z),\qquad
 P_J=\sum_{s=0}^J\mathcal L_{p,s}h^{-s}.
\]
Theorem~\ref{thm:full-expansion} gives
$|F_h/M_h-1|\le C\mathrm{e}^{-a_0h}$ on the entire real interval.
Its common coefficient norm also gives $\|P_J-1\|\le C/h$ on one
complex parameter neighborhood, since the consecutive majorant ratio
is $\mathfrak B_p^{\mathrm{FE}}(s+1)/h<1/2$. All choices below use this same model.
The constants in the two coefficient ranges can be chosen in the order
$d$, then $\eta$, then $\tau$; the small-disk argument below specifies
the first two choices.

\textbf{Comparison with the model: large indices.} We first compare the true coefficients with those of this model. Put
$m_* =\min(h,k)$ and, for $k\ge1$,
\[
 \mathscr F_x(w)=\mathrm{e}^{-kw/\sigma}
       \frac{F_h(x\mathrm{e}^{w/\sigma})}{F_h(x)},\qquad
 \mathscr M_x(w)=\mathrm{e}^{-kw/\sigma}
       \frac{M_h(x\mathrm{e}^{w/\sigma})}{M_h(x)}.
\]
There is a fixed $c>0$ on whose disks $|w|\le c\sqrt{m_*}$
the complex inverse parameter is defined and stays in the common
coefficient neighborhood. At criticality,
$\varepsilon(x\mathrm{e}^{w/\sigma})/\varepsilon(x)=1+O(w/\sqrt h)$;
choosing $c$ small fixes its square-root branch and keeps the tilted
real radius below $x_c$. Near zero use ordinary analytic inversion,
and use compactness between the endpoints.
The centered model action has zero first derivative and bounded second
derivative: its second derivative is
$v(x\mathrm{e}^{w/\sigma})/v(x)$. In the notation of
Lemma~\ref{rev:lem:normalized-radius}, with $k=h\Psi(x)$,
\[
 \log\mathscr M_x(w)=h\phi_x(w/\sigma)+\psi_x(w/\sigma)
 +\log\frac{P_J(x\mathrm{e}^{w/\sigma})}{P_J(x)},
\]
and $\|P_J-1\|\le C/h$ on the common neighborhood. Consequently, by
\eqref{rev:eq:radius-disk}, for $|w|\le c\sqrt{m_*}$,
\[
 |\mathscr M_x(w)|\le C\mathrm{e}^{C|w|^2+C|w|}.
\]
Positivity bounds $|\mathscr F_x(w)|$ by $\mathscr F_x(\Re w)$.
On the real diameter the exact identity
\[
 \mathscr F_x(u)=\mathscr M_x(u)
 \frac{1+e_h(x\mathrm{e}^{u/\sigma})}{1+e_h(x)},\qquad
 e_h=F_h/M_h-1,
\]
therefore gives the same disk bound and an additional factor
$\mathrm{e}^{-a_0h}$ for their difference on that diameter.

For a fixed small $\eta>0$, take $k\ge\eta h$ and
$T=\tau\sqrt h$. Choose $\tau>0$ after $\eta,a_0$, so that
$4T\le c\sqrt{m_*}$, $T\le c_1\sqrt{m_*}$ in \eqref{eq:QL5},
$T\le c_*\sqrt{m_*}$ in \eqref{eq:CE5},
and $T<\pi\sigma$. The scale-invariant harmonic measure used in
\eqref{eq:QL4} gives
\[
 \sup_{|t|\le T}|\mathscr F_x(it)-\mathscr M_x(it)|
 \le C\exp\{-a_0h/2+C T^2+C T\}
 \le C\mathrm{e}^{-a_1h}.
\]
Its integral has the same bound after decreasing $a_1$.
The integrated tail \eqref{eq:QL5} of the true distribution is also exponentially small
in $h$; the factor $\sigma\asymp\sqrt h/\varepsilon$ from the outer-arc
estimate is absorbed into this exponential bound.
The model central integral is bounded by \eqref{eq:CE5} and the
bounded $P_J$ ratio. Restoring the real normalization thus shows
\begin{equation}
 \left|[z^k]F_h-
 \frac{x^{-k}}{2\pi\sigma}\int_{-T}^T
             \mathrm{e}^{-ikt/\sigma}M_h(x\mathrm{e}^{it/\sigma})\,\mathrm{d}t\right|
 \le C\mathcal Q_{g,p,k}\mathrm{e}^{-a_2h}.
 \label{eq:CG14}
\end{equation}
Indeed the error before division is at most
$CF_h(x)x^{-k}\sigma^{-1}\mathrm{e}^{-a_2h}$, and
$\mathfrak p_k\sqrt k$ is bounded above and below for integer $k\ge1$.
The real matching factor $F_h(x)/M_h(x)=1+O(\mathrm e^{-a_2h})$
has been absorbed into the error. Here $k\ge\eta h$ with
$\eta\le\nu_*/2$, so that \eqref{eq:CG14} covers every density
$\Psi(x)\ge\nu_*/2$ used in the Gaussian range below.

\textbf{Comparison with the model: small indices.} For the remaining $k$, put $\varrho=d$.
On $|z|\le4\varrho$, positivity and the near-zero model give
$|(F_h-M_h)/J_h(0)|\le C\mathrm{e}^{Ch\varrho}$; on $[0,\varrho]$ this improves to
$C\mathrm{e}^{-a_0h+Ch\varrho}$. In the disk $|z-\varrho/2|<2\varrho$ slit along
$[0,\varrho]$, the slit harmonic measure is greater than $1/4$ on
$|z|\le\varrho/4$. To see this, rescale to $|w|<4$ with slit $[-1,1]$.
On the target $|w|\le3/2$ its Green function is less than $\log4$,
and on $|w|=4$ it is greater than $\log7$; comparison gives
$1-\log4/\log7>1/4$.
Choose $d$ small, then $\eta$ small, so the two-constants bound is
$C\mathrm{e}^{-a_3h}$ and the saddle for $1\le k<\eta h$
lies in $|z|\le\varrho/4$. Cauchy's bound on $|z|=x$ is on the
$\mathcal Q$ scale with a cost at most $C\sqrt{k+1}$:
$v/\Psi$ is bounded and $J_h(x)/J_h(0)\ge c$ near zero.
This polynomial cost is absorbed into a smaller exponential.
At $k=0$ use the marked estimate at zero directly.

\textbf{The Poisson range.} Choose a switching density $\nu_*>0$ strictly inside the fixed
Poisson neighborhood. For $\Psi(x)\le\nu_*$, the full model circle
is analytic and
\[
 \Re(S(x\mathrm{e}^{i\theta})-S(x))\le-cx(1-\cos\theta).
\]
This follows from $S(z)=2z+\sum_{m\ge2}\mathfrak s_mz^m$ and
$\sum_{m\ge2}m^2|\mathfrak s_m|x^m=O(x^2)$, using
$1-\cos(m\theta)\le m^2(1-\cos\theta)$.
The amplitudes are bounded, and $hv\asymp k$, so the standardized
model tail is Gaussian on the whole circle. For $k\ge\eta h$
it completes the integral in \eqref{eq:CG14} to $[z^k]M_h$ at
exponential precision. Together with the small disk, this compares the
true coefficients with those of the full model for every
$\Psi(x)\le\nu_*$.
The preceding Poisson argument now gives \eqref{eq:CG4}, with one
additional fixed error $C\mathrm{e}^{-ah}$, throughout this range.

\textbf{Uniform coordinates near criticality.} For the other range we justify the growing-order Gaussian expansion
uniformly through the critical point. In $\delta=y_0-Y$, write the
convergent local factors as
\[
 x_c-X(Y)=\delta^2\alpha(\delta),\quad
 \mathcal Z(Y)=\delta^{5/2}\zeta(\delta),\quad
 B_p(X(Y))=\delta^{p/2}\chi_p(\delta),\quad
 v(X(Y))=\delta^{-4}v_0(\delta),
\]
where all four factors are analytic and nonzero at zero.
For $z=X(Y)\mathrm{e}^{iu/\sqrt v}$ the inverse has
$\delta_z=\delta q(\delta,u)$, jointly analytic and close to one:
divide its equation $x_c-X(y_0-\delta q)=x_c-z$ by $\delta^2$
and apply the implicit-function theorem at $q=1,u=\delta=0$.
Its $q$ derivative is $2\alpha(0)\ne0$.
Thus $S(z)-S(x)$ and $B_p(z)/B_p(x)$ are analytic joint families,
expressed through $\log q$, $\zeta(\delta q)/\zeta(\delta)$ and
$q^{p/2}\chi_p(\delta q)/\chi_p(\delta)$.
Also $\Psi/\sqrt v$ extends analytically. Compactness gives common
neighborhoods for $\Psi\ge\nu_*/2$.

Set
\[
 \psi_Y(u)=S(x\mathrm{e}^{iu/\sqrt v})-S(x)-i\Psi u/\sqrt v,
 \qquad W_Y(u)=u\sqrt{-2\psi_Y(u)/u^2}.
\]
The removable quotient is one at zero; $W'_Y(0)=1$.
On one fixed disk $W_Y$ has a jointly analytic inverse $U_Y$,
with $\psi_Y(U_Y(w))=-w^2/2$.
Put $u_0=T/\sqrt h=\tau$.
The image endpoints are $\pm u_0+O(u_0^2)$. Straight connectors
to $\pm u_0$ have length $O(u_0^2)$ and
$\Re(w^2)\ge u_0^2/2$. A finite contour deformation therefore
costs at most $C\sqrt h\,u_0^2\mathrm{e}^{-hu_0^2/4}$.
After factoring $J_h(x)x^{-k}/\sigma$, the central integral becomes
\[
 \frac{\sqrt h}{2\pi}\int_{-u_0}^{u_0}\mathrm{e}^{-hw^2/2}
       \sum_{s=0}^J\widetilde f_s(w,Y)h^{-s}\,\mathrm{d}w
       +O(\mathrm{e}^{-ah}),
\]
where
$\widetilde f_s=[B_p(z(U_Y(w)))/B_p(x)]
 \mathcal L_{p,s}(z(U_Y(w)))U'_Y(w)$ has common norm
$C_p(\mathfrak B_p^{\mathrm{FE}})^ss!$ and $\widetilde f_0(0,Y)=1$.

\textbf{Gaussian expansion.} For $s\le K$, Taylor expand through degree $2(K-s)+1$ on that
fixed disk. The absolute Gaussian moment bounds the remainder by
$C^{K-s+1}(K-s+1)!h^{-K+s-1}\|\widetilde f_s\|$.
Odd terms vanish. To extend the polynomial moments to the whole line,
use, for $|w|\ge u_0$,
\[
 \mathrm{e}^{-hw^2/2}\le \mathrm{e}^{-hu_0^2/4}\mathrm{e}^{-hw^2/4}.
\]
The remaining moment majorants $C^j j!h^{-j}$ have bounded sum
for $j\le K\le c_ph$, after choosing $c_p$ sufficiently small
relative to the fixed disk and $u_0$. The sum over amplitude orders
is bounded by the same choice. All polynomial tails are therefore
$O(\mathrm{e}^{-ah})$, even at growing Taylor degree.
The unused $s>K$ amplitudes again form a geometric tail.
The factorial convolution used above yields a numerator remainder
$CB^{K+1}(K+1)!h^{-K-1}+C\mathrm{e}^{-ah}$.
After factoring $1/\sqrt{2\pi}$, its coefficients have the same
common-domain bounds and agree with \eqref{eq:CG2}. Indeed, for fixed
$Y$ and each fixed $s$, the model integral with the single amplitude
$\mathcal L_{p,s}$ is defined for every real $h$, and both families
are the coefficients of its asymptotic expansion as $h\to\infty$
through real values: \eqref{eq:CG2} before and the present ones after
the analytic change of variable $u=U_Y(w)$. Such expansion coefficients
are unique. The argument is applied separately to each $s$, and the
coefficients of the full model are then assembled by the formal product
in \eqref{eq:CG2}; no uniqueness statement is needed for the model with
the $h$-dependent order $J$.

\textbf{The exact Poisson normalization.} The divisor in \eqref{eq:CG1} is exact. Put
$R(k)=\sqrt{2\pi k}\,\mathfrak p_k$. Its logarithm satisfies
\[
 \log R(k)=-\int_0^\infty \mathrm{e}^{-kt}
 \left(\frac1{\mathrm{e}^t-1}-\frac1t+\frac12\right)\frac{\mathrm{d}t}{t}.
\]
The kernel has a removable origin, analytic radius $2\pi$, and is
bounded on the positive axis. Its Taylor remainder through degree
$L-2$ is bounded globally there by $CD^Lt^{L-1}$: use Cauchy near
zero and bound each retained power outside a fixed small interval.
Integration gives $CD^L(L-1)!k^{-L}$.
For exponentiation and inversion, the degree-$K$ log polynomial is
bounded on $|t|\le1/(D_0(K+1))$ for a fixed $D_0>0$, as in the proof
of Lemma~\ref{rev:lem:endpoint-gevrey}. Cauchy's tail at $t=1/k$, for
$K\le c_Rk$, is at most
$C[D_0(K+1)/k]^{K+1}\le C(\mathrm{e}D_0)^{K+1}(K+1)!k^{-K-1}$.
Since the exponential is bounded on the relevant real interval, the
remainder for $\log R$ transfers to $R$ and to $R^{-1}$, which
therefore have Gevrey-1 remainders. Here $k\ge(\nu_*/2)h$ and $\Psi^{-1}$ is
bounded analytic through criticality. By \eqref{eq:CG1},
$\mathcal Q_{g,p,k}/(\mathcal J_{g,p}(x)x^{-k}/(\pi\sigma))
=R(k)/\sqrt{2\pi}$. The two factors $1/\sqrt{2\pi}$ cancel.
Dividing the normalized numerator by $R(k)$ consequently preserves
the Gevrey-1 error and gives exactly
\eqref{eq:CG3}.

\textbf{Identification on the overlap.} For $\nu_*/2\le\Psi\le\nu_*$
both expansions hold; write $\mathcal A^{\mathrm P}_{p,j}$ for the
Poisson coefficients \eqref{eq:CG11} and $\mathcal A^{\mathrm G}_{p,j}$
for the Gaussian coefficients \eqref{eq:CG3}. Fixed-$Y$ uniqueness is
not available directly, since $k=h\Psi(X(Y))$ must be an integer, so we
argue by density. Suppose the two families agree below order $j$.
Subtracting the expansions at order $K=j$ gives
$\mathcal A^{\mathrm P}_{p,j}(Y_{h,k})-\mathcal A^{\mathrm G}_{p,j}(Y_{h,k})=O(h^{-1})$,
uniformly over all admissible $(h,k)$ with $Y_{h,k}$ in the overlap,
where $h\Psi(X(Y_{h,k}))=k$. Given $Y$ in the open overlap, take
$k=\lfloor h\Psi(X(Y))\rfloor$; since $\Psi\circ X$ is continuous and
strictly increasing, $Y_{h,k}\to Y$ as $h\to\infty$. Both coefficients
are continuous, so they agree at $Y$. By induction the two families
coincide on the open overlap. Both are analytic on
$0<\Psi(X(Y))<\nu_*$, so the identity theorem extends the equality to
this whole interval. Hence \eqref{eq:CG11} is the analytic continuation
of \eqref{eq:CG3} to $Y=0$, and the coefficients are analytic on one
complex neighborhood of $[0,y_0]$.

\textbf{Absorbing the residual error.} It remains to absorb the fixed exponential error at every allowed
order. After fixing the model, neighborhoods and cutoff, choose
$c_p<\kappa/2$ so small that $2c_p<1$ and
$2c_p\log(\mathrm{e}/(2c_p))<a$. For $L=K+1\le2c_ph$ and $\mathfrak B_p\ge1$,
\[
 \mathfrak B_p^LL!h^{-L}
 \ge \mathrm{e}^{-L\log(\mathrm{e}h/L)}
 \ge \mathrm{e}^{-2c_p\log(\mathrm{e}/(2c_p))h}.
\]
This bounds the fixed residual by the stated remainder simultaneously
for every $K\le c_ph$, including zero, once $h$ is large enough.
Finally choose $c'_p<c_p$ after the final $\mathfrak B_p$, with
$\mathfrak B_pc'_p<1$. The upper bound
$\mathfrak B_p^LL!h^{-L}\le(\mathfrak B_pL/h)^L$ proves the exponential
approximation at $K=\lfloor c'_ph\rfloor$.
\end{proof}

The first coefficient after extraction follows from
\eqref{eq:CG2}--\eqref{eq:CG3}, in which $s_1=-1/12$: with $E_1$ the
Gaussian coefficient displayed in the proof of
Lemma~\ref{lem:large-indices}, one has
$\mathcal A_{p,1}=\mathcal L_{p,1}+E_1+1/(12\Psi)$, that is,
\begin{equation}
 \begin{split}
 \mathcal A_{p,1}={}&\mathcal L_{p,1}
 -\frac{\mathscr D^2b+(\mathscr Db)^2}{2v}
 +\frac{(\mathscr Db)\mathscr D^2\Psi}{2v^2}
 +\frac{\mathscr D^3\Psi}{8v^2}
 -\frac{5(\mathscr D^2\Psi)^2}{24v^3}+\frac1{12\Psi},
 \qquad b=\log B_p.
 \end{split}
 \label{eq:CG12}
\end{equation}
Its apparent pole at zero is removable by \eqref{eq:CG11}.
By Theorem~\ref{thm:canonical-expansion} with $K=1$, retaining this term improves the relative remainder
to $O_p(h^{-2})$ over every coefficient index.

At the critical endpoint its value is
\begin{equation}
 \mathcal A_{p,1}(y_0)=-\frac{p^2}{8}+\frac p2-\frac{37}{60}.
 \label{eq:CG13}
\end{equation}
Indeed, writing $\varepsilon=1-x/x_c$, the convergent factorizations
recorded above in the proof of Theorem~\ref{thm:canonical-expansion}, together with
the exact relation $\varepsilon=\delta^2\alpha(\delta)/x_c$, give $S=-(5/4)\log\varepsilon+\alpha_S(\delta)$ and
$\log B_p=(p/4)\log\varepsilon+\alpha_B(\delta)$ with $\alpha_S,\alpha_B$
analytic at $\delta=0$. Since $\mathscr D\varepsilon=-(1-\varepsilon)$, one
has $\mathscr D^m\log\varepsilon=-(m-1)!\,\varepsilon^{-m}(1+O(\varepsilon))$,
while $\mathscr D^m$ applied to an analytic function of $\delta$ is
$O(\delta^{1-2m})=O(\varepsilon^{1/2-m})$, since
$\mathscr D=-(x/a(Y))\,\mathrm d/\mathrm d\delta$ and $a(Y)\asymp\delta$;
relative to $\varepsilon^{-m}$ this is $O(\sqrt\varepsilon)$. With
$v=\mathscr D^2S=\tfrac54\varepsilon^{-2}(1+O(\sqrt\varepsilon))$, the
normalized derivatives have the endpoint limits
\[
 \frac{\mathscr D^mS}{v^{m/2}}\to\Bigl(\frac54\Bigr)^{1-m/2}(m-1)!\quad(m\ge2),
 \qquad
 \frac{\mathscr D^mb}{v^{m/2}}\to-\frac p4\,(m-1)!\Bigl(\frac54\Bigr)^{-m/2}\quad(m\ge1),
\]
and the four terms of $E_1$ in \eqref{eq:CG12} tend respectively to
\[
 -\frac{\mathscr D^2b+(\mathscr Db)^2}{2v}\to-\frac{p^2}{40}+\frac p{10},\quad
 \frac{(\mathscr Db)\mathscr D^2\Psi}{2v^2}\to-\frac p5,\quad
 \frac{\mathscr D^3\Psi}{8v^2}\to\frac35,\quad
 -\frac{5(\mathscr D^2\Psi)^2}{24v^3}\to-\frac23,
\]
whose sum is $E_1(y_0)=-(p^2+4p)/40-1/15$.
Combine this with $\mathcal L_{p,1}(y_0)=t_p=(-2p^2+12p-11)/20$ from
\eqref{eq:ST13} and $\Psi^{-1}\to0$.

\section{The rate for adjacent coefficient ratios}
\label{sec:quantitative-coefficients}

This appendix proves Theorem~\ref{thm:quantitative-coefficient}, which
contains the rate in Theorem~\ref{thm:coefficient}. We quantify the three
degree ranges treated in Section~\ref{uniform-one-boundary-coefficient-ratios}.
The logarithmic comparison and concavity give the rate away from the
lowest degrees, including bounded positive codimension at the top.
Bracket difference estimates supply the remaining range. Balancing the
bounds gives $O(g^{-1/8})$ uniformly; we do not claim that this exponent
is optimal. The coefficients and the exact density are those of
Theorem~\ref{thm:coefficient}.

\begin{theorem}[Quantitative adjacent ratios]\label{thm:quantitative-coefficient}
There is a constant $C$ such that, for all sufficiently large $g$,
\begin{equation}
 \sup_{0\le k<3g-2}
 \left|\frac{a_{g,k}/a_{g,k+1}}{4\pi^2/t_{g,k}^2}-1\right|
 \le Cg^{-1/8}.
 \label{eq:QC1}
\end{equation}
For every fixed $\varepsilon\in(0,3)$ there is $C_\varepsilon$ such that
\begin{equation}
 \sup_{\varepsilon g\le k<3g-2}
 \left|\frac{a_{g,k}/a_{g,k+1}}{4\pi^2/t_{g,k}^2}-1\right|
 \le C_\varepsilon g^{-1/4}.
 \label{eq:QC2}
\end{equation}
In the second supremum $k$ is an integer. Both bounds include bounded
positive codimension from the top coefficient.
\end{theorem}

\subsection{The logarithmic model up to its endpoints}

We first make the logarithmic comparison of
Section~\ref{uniform-one-boundary-coefficient-ratios} quantitative, up to
both ends of the degree range. The result is
Lemma~\ref{lem:log-comparison}: the logarithm of the normalized
coefficient sequence agrees with $gf(k/g)$ to within $O(\sqrt g)$ for
every degree $0\le k\le N$. Its inputs are the exact model \eqref{eq:C4},
the comparison \eqref{eq:C5} and the saddle-point description
\eqref{eq:C10}--\eqref{eq:C15} of the function $f$, all from
Section~\ref{uniform-one-boundary-coefficient-ratios}. The new point is
the degenerate end $m/h\to0$, where the tilted summands depend on $h$;
Lemma~\ref{rev:lem:triangular-local-limit} supplies the local limit
theorem for such a triangular array, and Lemma~\ref{lem:uniform-saddle}
uses it to bound the model coefficients uniformly in the codimension $m$.

Use $H,\widetilde A_m,\mu$ from \eqref{eq:C4} and \eqref{eq:C10}.
Let $\zeta(u)$ be the inverse of $\mu(z)$ on $0\le u\le3/2$,
with $\zeta(0)=0$. Since $\mu(z)=z/5+O(z^2)$ at zero and
$\mu'(z)>0$ throughout $[0,\pi^2/2]$, the function
\[
 \chi(u)=u/\zeta(u),\qquad \chi(0)=1/5,
\]
is positive and analytic in a neighborhood of the closed real interval.
Write $f$ for the function in \eqref{eq:C14}, extended continuously to
$[0,3]$, and set
\begin{equation}
 \mathfrak q(\alpha)=4\alpha\chi((3-\alpha)/2),\qquad
 f'(\alpha)=-\log\mathfrak q(\alpha),\qquad
 |f''(\alpha)|\le C(1+\alpha^{-1})\quad(0<\alpha\le3).
 \label{eq:QC3}
\end{equation}
These follow from \eqref{eq:C15}; in particular
$\mathfrak q(3)=12/5$ and $f(3)=0$.

We need the endpoint-uniform version of the logarithmic saddle estimate
\eqref{eq:C12}. At the degenerate end $m/h\to0$ the tilted summands
depend on $h$, so the local limit theorem required is one for a
triangular array; we isolate it.

\begin{lemma}[Triangular local limit]\label{rev:lem:triangular-local-limit}
For each $h$ let $(p_j)_{j\ge0}=(p_j^{(h)})_{j\ge0}$ be a probability
law on $\mathbb N$ with characteristic function
$\varphi_1(\theta)=\sum_jp_j\mathrm{e}^{ij\theta}$, mean $\mu_1$ and variance
$\operatorname{Var}_1$, and let $S_h$ be the sum of $h$ independent
copies. Suppose that there are $z=z_h>0$ and constants $c,C>0$ such
that $\varphi_1$ does not vanish on $[-\pi,\pi]$ and
\begin{equation}\label{rev:eq:tll-hypotheses}
 p_0p_1\ge cz,\qquad \operatorname{Var}_1\le Cz,\qquad
 \sup_{|\theta|\le\pi}\bigl|\partial_\theta^3\log\varphi_1(\theta)\bigr|\le Cz .
\end{equation}
Then $\operatorname{Var}_1\ge p_0p_1\ge cz$, and for every integer $m$
with $|m-h\mu_1|\le1$, as $hz\to\infty$,
\begin{equation}\label{rev:eq:tll-conclusion}
 \Pr(S_h=m)=\bigl(2\pi h\operatorname{Var}_1\bigr)^{-1/2}(1+o(1)),
\end{equation}
uniformly over all arrays satisfying \eqref{rev:eq:tll-hypotheses} with
the same $c,C$.
\end{lemma}

\begin{proof}
First, for every $\theta$,
\begin{equation}\label{rev:eq:tll-fourier}
 |\varphi_1(\theta)|^2=\sum_{j,k}p_jp_k\cos((j-k)\theta)
 =1-\sum_{j,k}p_jp_k\bigl(1-\cos((j-k)\theta)\bigr)
 \le1-2p_0p_1(1-\cos\theta),
\end{equation}
since $\sum_jp_j=1$ and every term of the middle sum is nonnegative;
this is the bound \eqref{eq:TL40}, reproved here. Also
$\operatorname{Var}_1=\tfrac12\sum_{j,k}p_jp_k(j-k)^2\ge p_0p_1$.
Write $\log\varphi_1(\theta)=i\mu_1\theta-\tfrac12\operatorname{Var}_1\theta^2+R(\theta)$;
by Taylor's formula and the third hypothesis,
$|R(\theta)|\le Cz|\theta|^3/6$ on $[-\pi,\pi]$. Fourier inversion gives
\[
 \sqrt{h\operatorname{Var}_1}\,\Pr(S_h=m)
 =\frac1{2\pi}\int_{-\pi\sqrt{h\operatorname{Var}_1}}^{\pi\sqrt{h\operatorname{Var}_1}}
 \exp\Bigl\{-\frac{\vartheta^2}2+hR\Bigl(\frac\vartheta{\sqrt{h\operatorname{Var}_1}}\Bigr)
 -i(m-h\mu_1)\frac{\vartheta}{\sqrt{h\operatorname{Var}_1}}\Bigr\}\,\mathrm{d}\vartheta ,
\]
after the substitution $\theta=\vartheta/\sqrt{h\operatorname{Var}_1}$.
For fixed $\vartheta$ the exponent tends to $-\vartheta^2/2$, because
$h|R|\le Chz|\vartheta|^3(h\operatorname{Var}_1)^{-3/2}\le C'|\vartheta|^3(hz)^{-1/2}\to0$
by $\operatorname{Var}_1\ge cz$, and $|m-h\mu_1|\,|\vartheta|(h\operatorname{Var}_1)^{-1/2}\to0$.
For domination, \eqref{rev:eq:tll-fourier} and $1-\cos\theta\ge2\theta^2/\pi^2$
on $[-\pi,\pi]$ give
$|\varphi_1(\theta)|^h\le\exp\{-hp_0p_1(1-\cos\theta)\}\le\exp\{-2chz\theta^2/\pi^2\}
\le\exp\{-2c\vartheta^2/(\pi^2C)\}$, using $\operatorname{Var}_1\le Cz$.
Dominated convergence gives the limit $(2\pi)^{-1/2}$ of the right side,
which is \eqref{rev:eq:tll-conclusion}. The convergence is uniform over
the class: on $|\vartheta|\le(hz)^{1/8}$ the exponent differs from
$-\vartheta^2/2$ by at most $C''(hz)^{-1/8}$ by the two displayed
bounds, and the complement is controlled by the Gaussian majorant;
every constant depends only on $c,C$.
\end{proof}

\begin{lemma}[Endpoint-uniform saddle estimate]\label{lem:uniform-saddle}
For a positive integer $h$ and integers $0\le m\le3h/2$, let $z=\zeta(m/h)$, interpreting the
$m=0$ expressions by continuity. Then, uniformly in $m$,
\begin{equation}
 \log \widetilde A_m(h)=-h\log H(z)-m\log z+O(\log(h+1)).
 \label{eq:QC4}
\end{equation}
\end{lemma}

\begin{proof}
Tilt $H^{-1}$ at $z$ and let
$S_h$ be the sum of $h$ independent variables with probability
generating function $H(z)/H(zw)$. Exactly,
\[
 \widetilde A_m(h)=H(z)^{-h}z^{-m}\Pr(S_h=m).
\]
There is $c>0$ such that $\Pr(S_h=m)\ge c/\sqrt{m+1}$ uniformly.
We give the details at the degenerate endpoint $m/h\to0$, where
Lemma~\ref{rev:lem:triangular-local-limit} applies with the parameter
$z=\zeta(m/h)$. Indeed, the summand law has
$p_j=H(z)z^j[w^j]H(w)^{-1}=H(z)\mathsf h_jz^j$ with the positive coefficients
$\mathsf h_j=[z^j]H^{-1}$ ($\mathsf h_0=1$, $\mathsf h_1=1/5$), so $p_0=H(z)$ and
$p_1=H(z)z/5$, and $p_0p_1=H(z)^2z/5\ge cz$ on $0<z\le\pi^2/2$. Its
characteristic function is $\varphi_1(\theta)=H(z)/H(z\mathrm{e}^{i\theta})$,
which does not vanish for $z\le\pi^2/2$, and
\[
 \log\varphi_1(\theta)=\log H(z)+\sum_{j\ge1}c_jz^j\mathrm{e}^{ij\theta},
 \qquad -\log H(z)=\sum_{j\ge1}c_jz^j ,
\]
where the series converges for $|z|$ below the first zero of $H$,
which exceeds $\pi^2/2$. Hence
$\partial_\theta^3\log\varphi_1(\theta)=\sum_j(ij)^3c_jz^j\mathrm{e}^{ij\theta}$
has modulus at most $\sum_jj^3|c_j|z^j\le Cz$, and
$\operatorname{Var}_1=\sum_jj^2c_jz^j\le Cz$, uniformly for
$0<z\le\pi^2/2$. When $m\to\infty$ and $m/h\to0$, one has $z\asymp m/h$,
so $hz\asymp m\to\infty$, and the mean of $S_h$ is exactly
$h\mu(z)=m$. Lemma~\ref{rev:lem:triangular-local-limit} therefore gives
$\Pr(S_h=m)\sim(2\pi h\operatorname{Var}_1)^{-1/2}\ge c/\sqrt{hz}\ge c'/\sqrt{m+1}$.
When $m/h$ stays away from zero, the summand law lies in a compact
tilted family with positive zero and one probabilities and variance
bounded above and below; hence the sum has variance comparable to $h$.
If $m$ stays bounded, pass to a fixed integer $m$;
then $z=5m/h+O_m(h^{-2})$ and
$\widetilde A_m(h)=(h/5)^m/m!\,(1+O_m(h^{-1}))$, so the central mass tends to
$\mathrm{e}^{-m}m^m/m!>0$. At $m=0$ it is one.
These cases and a subsequence contradiction prove the asserted uniform
lower bound. The upper bound is one, proving \eqref{eq:QC4}.
\end{proof}

Put
\[
 G_g(k)=\log\frac{J_k}{J_Ng^{2(N-k)}},\qquad N=3g-2.
\]

\begin{lemma}[Logarithmic comparison up to the endpoints]\label{lem:log-comparison}
There is a constant $C$ such that
\begin{equation}
 \max_{0\le k\le N}|G_g(k)-g f(k/g)|\le C\sqrt g.
 \label{eq:QC5}
\end{equation}
\end{lemma}

\begin{proof}
This follows from the model \eqref{eq:C4}, the comparison \eqref{eq:C5} and
Lemma~\ref{lem:uniform-saddle}. The passage to $f$ is uniform, even at
$m=0$. To see this, use the
exact model for the intersection number,
$J^{\mathrm{mod}}_k=3^h\Gamma(h)4^{-g}m!\widetilde A_m(h)/(2k+1)!!$ from
Section~\ref{uniform-one-boundary-coefficient-ratios}, with $h=2g-1$ and
$m=N-k$. Its quotient by the top model value is
$m!\,\widetilde A_m(h)\,(2N+1)!!/(2k+1)!!$, so by
Lemma~\ref{lem:uniform-saddle}
\[
 G_g^{\mathrm{mod}}(k):=\log\frac{J^{\mathrm{mod}}_k}{J^{\mathrm{mod}}_Ng^{2(N-k)}}
 =\log m!+\Phi(h,m)+\sum_{j=k+1}^{N}\log(2j+1)-2m\log g+O(\log(g+1)),
\]
where the saddle action is
\[
 \Phi(h,m)=-h\log H(\zeta(m/h))-m\log\zeta(m/h).
\]
Since $\zeta(m/h)$ is the stationary point of $z\mapsto-h\log H(z)-m\log z$,
whose derivative is $(h\mu(z)-m)/z$ by \eqref{eq:C10}, the envelope
identity gives
\[
 \partial_h\Phi(h,m)=-\log H(\zeta(m/h)),\qquad
 \partial_m\Phi(h,m)=-\log\zeta(m/h).
\]
The first is bounded, since $H$ decreases from $1$ to $3/\pi^2$ on
$[0,\pi^2/2]$; the second has only the integrable singularity
$\log\zeta(u)=\log u+O(1)$ at $u=0$, since $\zeta(u)\asymp u$. Hence
replacing $(h,m)=(2g-1,3g-2-k)$ by $(2g,3g-k)$ along a straight path,
on which $m'/h'$ stays in $[0,3/2]$, changes $\Phi$ by
$O(1)+\int|\log\zeta(m'/h')|\,\mathrm{d}m'=O(\log(g+1))$, and
$\Phi(2g,3g-k)=g[-2\log H(z)-(3-\alpha)\log z]$ with $\alpha=k/g$ and
$z=\zeta((3-\alpha)/2)$, the first two terms of $gf(\alpha)$ in
\eqref{eq:C14}. For the remaining terms, Stirling's bound with
logarithmic error and the convention $0\log0=0$ give, uniformly in
$0\le k\le N$,
\[
\begin{aligned}
 \log m!&=g(3-\alpha)\log g+g(3-\alpha)\log(3-\alpha)-g(3-\alpha)+O(\log(g+1)),\\
 \sum_{j=k+1}^{N}\log(2j+1)
 &=\bigl[3g\log g+3g\log6-3g\bigr]-\bigl[g\alpha\log g+g\alpha\log(2\alpha)-g\alpha\bigr]
   +O(\log(g+1)),\\
 -2m\log g&=-2g(3-\alpha)\log g+O(\log(g+1)),
\end{aligned}
\]
the second line using $\log(2n+1)!!=n\log(2n)-n+O(\log(n+1))$ at
$n=N$ and $n=k$. The coefficients of $g\log g$ add up to
$(3-\alpha)+(3-\alpha)-2(3-\alpha)=0$, and the remaining terms are
$g[(3-\alpha)\log(3-\alpha)+3\log6-\alpha\log(2\alpha)-2(3-\alpha)]$,
which completes $gf(\alpha)$. Thus
$G^{\mathrm{mod}}_g(k)=gf(k/g)+O(\log(g+1))$ uniformly. The two logarithmic errors between
the true and the model values are at most $C\sqrt g$ by \eqref{eq:C5}, including the top
coefficient, which proves \eqref{eq:QC5}.
\end{proof}

\subsection{Quantitative secants and the top endpoint}

Concavity converts the logarithmic comparison into two-sided secant
bounds for the adjacent ratios. Together with the exact top anchor,
these give the rate for all $k\ge g^{3/4}$; the remaining range
$k\le g^{3/4}$ is treated in the final part of the same proof.

\begin{proof}[Proof of Theorem~\ref{thm:quantitative-coefficient}]
\textbf{Secant bounds.}
The sequence $G_g$ is concave by \eqref{eq:C7}. Write
$q_k=J_k/J_{k+1}$. For an integer $\ell\ge1$ with
$\ell\le k/2$ and $k+\ell\le N$, its secants imply
\[
 -\frac{G_g(k)-G_g(k-\ell)}\ell
 \le\log(q_k/g^2)
 \le-\frac{G_g(k+\ell)-G_g(k)}\ell.
\]
Equations \eqref{eq:QC3} and \eqref{eq:QC5} yield
\begin{equation}
 \left|\log\frac{q_k}{g^2\mathfrak q(k/g)}\right|
 \le C\left(\frac{\sqrt g}\ell+\frac\ell k\right),
 \label{eq:QC6}
\end{equation}
where $\ell\le k/2$ ensures $\alpha^{-1}\le2g/k$ on both secant
intervals $[(k-\ell)/g,k/g]$ and $[k/g,(k+\ell)/g]$, so that the
second-derivative bound in \eqref{eq:QC3} gives the secant error
$O(\ell/k)$.
If $k+\ell>N$, retain the lower secant and use
$q_k\le q_{N-1}$ for the upper bound. The exact top anchor
\eqref{eq:C18} gives $q_{N-1}/g^2=12/5+O(g^{-1})$.
Since $\log\mathfrak q$ is Lipschitz near three and
$3-k/g\le(\ell+2)/g$, this gives the same bound
\eqref{eq:QC6}, with a changed constant.

Take $\ell=\lfloor g^{5/8}\rfloor$. For $k\ge g^{3/4}$ it is
at most $k/2$ for large $g$, and \eqref{eq:QC6} is $O(g^{-1/8})$.
For $k\ge\varepsilon g$, instead take
$\ell=\lfloor g^{3/4}\rfloor$ to obtain $O_\varepsilon(g^{-1/4})$.

\textbf{From intersection numbers to coefficient ratios.}
These estimates for the intersection numbers translate into the
statement of the theorem without losing the bounded-codimension range.
Indeed, put
$m=N-k$, $h=2g-1$, and $z_{g,k}=t_{g,k}^2/2$. By \eqref{eq:C11},
$\mu(z_{g,k})=m/h$. Therefore the value of $q_k=J_k/J_{k+1}$ that
corresponds, through \eqref{eq:C3}, to the proposed equivalent
$4\pi^2/t_{g,k}^2$ of $a_{g,k}/a_{g,k+1}$ is
\begin{equation}
 \frac{(2k+3)m}{z_{g,k}}
 =(2k+3)h\chi(m/h)
 =g^2\mathfrak q(k/g)\left(1+O(k^{-1}+g^{-1})\right).
 \label{eq:QC7}
\end{equation}
The last error is uniform for $k\ge1$: the two arguments of $\chi$
differ by $O(g^{-1})$, while $\chi$ is positive with bounded
logarithmic derivative on the closed interval. This calculation
retains $m$ exactly, including $m=1$. By \eqref{eq:C3}, the two
secant estimates prove \eqref{eq:QC1} for $k\ge g^{3/4}$ and all
of \eqref{eq:QC2}.

\textbf{The small-degree range.}
It remains to consider $k\le g^{3/4}$. Put
$\rho_k=b_{g,k}/b_{g,k+1}\ge1$ and
$\omega_k=(2k+3)(N-k)$ as before. Equation~\eqref{eq:C25} gives
$\rho_k=1+O(g^{-1/2})$ for $k\le\delta_0\sqrt g$ with a fixed
small $\delta_0>0$. For the remaining indices we claim
\begin{equation}
 0\le\rho_k-1\le C\left(
       \frac{\log(g+1)}{\sqrt k}+\sqrt{\frac kg}+\frac1k\right)
 \quad(\delta_0\sqrt g\le k\le g^{3/4}).
 \label{eq:QC8}
\end{equation}
First $\rho_k<2$ uniformly in this range for large $g$. Otherwise
monotonicity of $q_j$ and $\omega_j/\omega_k\le1+\delta$ propagate a fixed
ratio gap over a block of length comparable to $k$: with $\delta=1/3$,
for $k\le j\le4k/3$ one has
$\omega_j/\omega_k\le(2j+3)/(2k+3)\le4/3$, hence
$\rho_j\ge\rho_k\omega_k/\omega_j\ge3/2$ on a block of at least $k/3$
indices. Its product
would force $\log(V_{g,1}/b_{g,J+1})\ge\tfrac k3\log\tfrac32$ at
$J=\lfloor4k/3\rfloor$,
contradicting \eqref{eq:C23} at $J+1\le2k$, whose upper bound
$C\log^2(g+1)+Ck^2/g$ is $o(k)$ uniformly
on this range.

Now put $\eta=\rho_k-1\in[0,1]$ and $\delta=\eta/6$.
For $k\le j\le\lfloor(1+\delta)k\rfloor$, the same monotonicity
and the exact factor $\omega_j$ give
$\rho_j\ge(1+\eta)/(1+\delta)\ge1+\eta/2$.
This block has $\lfloor\eta k/6\rfloor+1\ge\eta k/6$ indices.
Summing logarithms over the block and using \eqref{eq:C23} at its
last index plus one gives
\[
 c\eta^2 k\le C\log^2(g+1)+Ck^2/g.
\]
All these indices lie below the fixed $c_0g$ range of \eqref{eq:C23}
for large $g$.
This proves \eqref{eq:QC8}. Its right side is $O(g^{-1/8})$ on
the displayed range. Thus $\rho_k=1+O(g^{-1/8})$ uniformly for
$0\le k\le g^{3/4}$.

Finally the saddle equation at $t=\pi$ gives
$\pi^2/t_{g,k}^2=1+O((k+1)/g)$ uniformly on this small-degree
range: by \eqref{eq:C10}--\eqref{eq:C11}, $t=t_{g,k}$ solves
$R(t)=\alpha_{g,k}=(2k+1)/(2g-1)$ with
$R(t)=t^2\sin t/(\sin t-t\cos t)$, and $R(\pi)=0$, $R'(\pi)=-\pi$, so
$\pi-t_{g,k}=\alpha_{g,k}/\pi+O(\alpha_{g,k}^2)$ and
$\pi^2/t_{g,k}^2=1+2\alpha_{g,k}/\pi^2+O(\alpha_{g,k}^2)$. Since $a_{g,k}/a_{g,k+1}=4\rho_k$, the desired
relative error is $O(g^{-1/8})$ there as well. Together with the
secant range this proves both assertions.
\end{proof}

The argument is a quantitative version of the one used for the leading
equivalent, based on the same logarithmic comparison and concavity. Liu--Xu's fixed-codimension expansion
has finer errors at each fixed codimension~\cite{LX}; the bounds here
hold with a single constant across all valid degrees.

\part*{Appendices to Part~2: spectral laws}
\section{Proof of the short-curve theorem}\label{app:short-curves}

This appendix proves Theorem~\ref{thm:pants-poisson} by the method of
moments, following the outline at the end of
Subsection~\ref{sec:pants-spectrum}. The proof is a single argument
in three steps. Step~1 derives a coefficient-ratio bound,
\eqref{eq:PC3}, from the positive convolution factor \eqref{eq:TL33}
and Lemma~\ref{lem:point-mass}; it shows that curves cutting off many
cusps form a uniformly summable tail, and it is reused in the proofs
of Lemmas~\ref{lem:super-cuts} and~\ref{lem:super-geometry}
(Appendix~\ref{app:super-cuts}). Step~2 computes the mixed factorial
moments for finitely many fixed cusp types: the disjoint
configurations give the Poisson moments, while nested configurations
are negligible. Step~3 shows that, on the same length scale, all
short curves other than planar separating ones are rare, which
identifies the full primitive process and the systole. The volume
inputs are Mirzakhani's integration formula \eqref{eq:PS3}, the
boundary estimate \eqref{eq:PS4}, the volume ratio \eqref{eq:PS5}
obtained from Theorem~\ref{thm:volume}, and the coefficient identity
\eqref{eq:TL50}; we use the notation of
Subsection~\ref{sec:pants-spectrum}, in particular the intensity
$q_g$ in \eqref{eq:PS1} and the weights \eqref{eq:PC1}.

\begin{proof}[Proof of Theorem~\ref{thm:pants-poisson}]
Write $P=2\pi^2$, $h=2g-1$, $m=n-1$, and $F_h=\Phi_g'$.
Choose $x\in(0,x_c)$ by $h\Psi(x)=m$ and put $q_g^*=nx/(4P)$.
By \eqref{eq:SP3}, $\mathcal C(\lambda)=\pi^2\lambda/X(Y_\lambda)$, so
$q_g=(n-1)X(Y_{n/g})/(8\pi^2)$ and
$q_g^*=nX(Y_{\lambda_{\rm eff}})/(8\pi^2)$ with the effective density
$\lambda_{\rm eff}=2(n-1)/(2g-1)$. Since
$\lambda_{\rm eff}g/n\to1$ and
$\lambda\mapsto X(Y_\lambda)$ is continuous and positive on
$(0,\infty)$, with $X(Y_\lambda)\sim\lambda/4$ as $\lambda\to0$ and
$X(Y_\lambda)\to x_c$ as $\lambda\to\infty$, a subsequence argument gives
$q_g^*/q_g\to1$; and $q_g\asymp n^2/(g+n)\to\infty$ by
Subsection~\ref{sec:pants-spectrum}.

\emph{Step 1: a uniform tail bound over cusp types.}
We first show that curves cutting off a large number of cusps contribute
negligibly, uniformly along the sequence. The input is the
coefficient-ratio bound \eqref{eq:PC3}, which compares volumes whose cusp
counts differ by $j$.
By the positive convolution \eqref{eq:TL33} and
Lemma~\ref{lem:point-mass}, every point mass is at most $C/\sigma$,
where $\sigma^2=h\mathscr D\Psi(x)$. The local limit
\eqref{eq:TL47} gives a lower bound $c/\sigma$ at $m$.
Consequently, for all sufficiently large $h,m$,
\begin{equation}
 \frac{[z^{m-j}]F_h(z)}{[z^m]F_h(z)}
 =x^j\frac{\mathbb P_x(N=m-j)}{\mathbb P_x(N=m)}
 \le Cx^j\qquad(0\le j\le m).
 \label{eq:PC3}
\end{equation}
Only the denominator uses the local limit; the numerator may be at any
index, including zero.

Let $N_{k,L}$ count curves of type $k$ with length at most $L$, and put
$a_k=V_{0,k+1}/(k!P^{k-2})$. Mirzakhani integration gives
\begin{equation}
 \begin{aligned}
 \mathbb E N_{k,L}=\frac{\binom nk}{V_{g,n}}
 \int_0^L\ell\,V_{0,k+1}(\ell,0^k)
                  V_{g,n-k+1}(\ell,0^{n-k})\,\mathrm{d}\ell,\\
 \binom nk\frac{V_{0,k+1}V_{g,n-k+1}}{V_{g,n}}
 =\frac nP a_k\frac{[z^{n-k}]F_h}{[z^{n-1}]F_h}.
 \end{aligned}
 \label{eq:PC4}
\end{equation}
The two sides have different genera, so there is no side-exchange
factor. The coefficient identity follows from \eqref{eq:TL50}; it
holds up to $k=n$. Bound \eqref{eq:PS4} changes the integrand's two
volumes by a factor between one and $\mathrm{e}^L$, independently of $k$.
Using \eqref{eq:PC3} with $j=k-1$ therefore gives, for each fixed $T$,
\begin{equation}
 \mathbb E\sum_{k\ge K}N_{k,T/\sqrt{q_g^*}}
 \le CT^2\sum_{k\ge K}a_kx_c^{k-2}\longrightarrow0
 \quad(K\to\infty),
 \label{eq:PC5}
\end{equation}
uniformly along the sequence. To justify summability directly, apply
the string-shift identity \eqref{eq:H8} at $(g,n)=(0,3)$: by
\eqref{eq:H6} with $h=1$ and $m=0$,
$\Phi_0'''(x)=U^{Y}_{0,3}(0,0,0)=\langle\tau_0^3\rangle_0/a(Y)=1/a(Y(x))=Y'(x)$.
Since $\Phi_0(x)=x^3/6+O(x^4)$ has no unstable terms,
$\Phi_0''(0)=\Phi_0'(0)=0$; with $Y(0)=0$ this gives
$\Phi_0''=Y$ and $\Phi_0'(x)=\int_0^xY(u)\,\mathrm{d}u$. This is finite at
$x_c$; positivity gives convergence of the coefficient series there.
Its initial term $x^2/2$ also gives $W(0)=1$.

\emph{Step 2: fixed types and nested configurations.}
Next we compute the mixed factorial moments for finitely many fixed cusp
types. The main contribution comes from configurations whose planar
sides are pairwise disjoint, which we call \emph{stars}; nested
configurations are shown to be negligible.
Fix $r$ ordered distinct curves of prescribed types $k_1,\ldots,k_r$,
and write $K=\sum k_i$. For large $g$ all rescaled windows lie below
$2\operatorname{arcsinh}1$, so the collar lemma makes their boundaries
disjoint. For a star, the ordered
cusp-label multiplicity is $(n)_K/\prod_i k_i!$, and the complement
has type $(g,n-K+r)$. There is no $r!$ divisor. The boundary error is
$1+O_r(L)$, as each cut occurs twice in \eqref{eq:PS4}.
For fixed $d$, \eqref{eq:PS5} and $q_g^*/q_g\to1$ give
$V_{g,n-d}/V_{g,n}=(x/(Pn))^d(1+o(1))$.
Thus for windows $[a_i,b_i]/\sqrt{q_g^*}$ the star contribution tends to
\begin{equation}
 \prod_{i=1}^r
       w_{k_i}(x_\lambda)(b_i^2-a_i^2).
 \label{eq:PC6}
\end{equation}
Before the limit is taken, the star contribution equals the right side
at $x=x(g,n)$ up to a factor $1+o(1)$; the error is relative, so the
statement keeps its meaning when $x\to0$ and $w_k(x)\to0$ for $k>2$.

Filling cusps turns each planar side into a disk. Since each curve separates and $g\ge1$, disjoint boundaries
therefore have nested or disjoint disks; their cut graph is a tree with
one genus-$g$ root and $r$ planar vertices. Let $d$ be the root degree,
$a_v$ the number of cusps at a planar vertex and $d_v$ its degree;
put $a=\sum_v a_v$. Prescribed types bound $a$ by $K$.
There are finitely many such trees with ordered edges and cusp
allocations, depending only on the fixed types. Classification of
surfaces gives the corresponding mapping-class orbits with bounded
symmetry factors; assigning cusp labels costs at most $Cn^a$.
The root volume ratio is $O((x/(Pn))^{a-d})$, the planar volumes
are fixed, and the $r$ length integrals cost $O((nx)^{-r})$.
Each cut type contributes at most
\begin{equation}
 Cn^{d-r}x^{a-d-r}=O(n^{-1})\quad\hbox{if it is not a star}.
 \label{eq:PC7}
\end{equation}
Indeed stability gives $a_v+d_v\ge3$ and
$\sum_vd_v=2r-d$, hence $a\ge r+d$. The power of $x$ is
nonnegative, so no loss occurs as $x\to0$. A nonstar has $d\le r-1$;
if $d=r$, every planar vertex is a leaf and the configuration is a star.
This accounts for every fixed-type mixed factorial moment.

The limits \eqref{eq:PC6} are the mixed factorial moments of independent
Poisson counts. Higher fixed moments give uniform integrability, and
the product Poisson law is moment determinate. First moments give
tightness on compact sets in the marked space. This proves
\eqref{eq:PC2}. Its projection onto the length coordinate requires more than vague
convergence: truncate the types, pass to the finite-type limit, then
let the cutoff grow using \eqref{eq:PC5} and the convergent limiting
series. This proves the stated unmarked process. Finally replace
$q_g^*$ by $q_g$. The type-two restriction gives \eqref{eq:PS2}
without assuming a limit for $n/g$, by subsequences and compactness of
$[0,\infty]$. Its zero-count probabilities give the pants minimum law.

\emph{Step 3: the remaining curves and the systole.}
It remains to show that, on the scale $q_g^{-1/2}$, there are no short
curves other than planar separating ones: nonseparating curves, and
separating curves with positive genus on both sides, are rare.
Nonseparating simple curves of length at most $L\le L_0$ have expected
count at most $C\mathrm{e}^LL^2V_{g-1,n+2}/V_{g,n}\le CL^2$, by integration,
\eqref{eq:PS4}, and the uniform lower bound for $B_{g,n}$.
To sum every positive-genus separating cut, put
$S_g(z)=\sum_{i=1}^{g-1}\Phi_i'(z)\Phi_{g-i}'(z)$. Exact positivity
\eqref{eq:TL33} includes genus one, so the child powers add to $h-1$:
\begin{equation}
 S_g(z)=V(z)^{h-1}Q_g(z),\qquad Q_g\ge0
 \ \hbox{coefficientwise},\qquad
 \frac{S_g(x)}{F_h(x)}\le\frac{Ca(Y)r(Y)^2}{g^2}.
 \label{eq:FS2}
\end{equation}
The last bound is \eqref{eq:TL26}: undoing \eqref{eq:TL9} leaves
exactly $ar^2$, as the child complexities sum to $2g-2$ and their
deficits to $3g-4$. In particular bounded child genera are retained.
Lemma~\ref{lem:point-mass}, applied to the positive factor $V^{h-1}$,
gives $\sup_j\mathbb P_x^{S}(N=j)\le C/\sigma$. Against the parent lower
bound $\mathbb P_x^{F}(N=n-1)\ge c/\sigma$, the standard deviations
cancel. The coefficient identity \eqref{eq:TL50} therefore yields
\begin{equation}
 \begin{aligned}
 \mathcal A_{g,n}
 &:=\sum_{i=1}^{g-1}\sum_{k=0}^n\binom nk
       \frac{V_{i,k+1}V_{g-i,n-k+1}}{V_{g,n}}
   =\frac1P\frac{[z^n]S_g(z)}{[z^n]\Phi_g(z)}\\
 &=\frac{n}{Px}\frac{S_g(x)}{F_h(x)}
       \frac{\mathbb P_x^S(N=n)}{\mathbb P_x^F(N=n-1)}
   \le\frac{Cnar^2}{xg^2}\le\frac Cg.
 \end{aligned}
 \label{eq:FS3}
\end{equation}
Indeed $(n-1)/x=h/(Aar)$, while $A$ is bounded below and $r$ above.
Thus no small-density factor is lost. Each separating orbit with two
positive genera occurs in this ordered sum. Its integration symmetry
and one-handle factors are at most one; applying \eqref{eq:PS4} to
both components bounds its total expected count by
$C\mathrm{e}^LL^2\mathcal A_{g,n}\le CL^2/g$. This includes all cusp partitions.

For fixed $T$, take $L=T/\sqrt{q_g}$. Markov's inequality shows that
with probability $1-O(T^2/q_g)$ every simple geodesic in this window
is planar. Primitive geodesics shorter than $2\operatorname{arcsinh}1$
are simple, by the standard collar theorem used in
\cite[Section~1.1]{HTSphere}. Since $q_g\to\infty$, the full primitive
process and planar process agree on each bounded rescaled window
with probability tending to one. This transfers the process limit
and its zero-count law to the systole.

Finally restrict the first point to $[\eta,T]$ and finitely many cusp
types. The first-point map is continuous when its first length is
unique and avoids the endpoints, as happens almost surely in the
limiting Poisson process. Let the type cutoff grow using
\eqref{eq:PC5}, then $\eta\downarrow0$ and $T\to\infty$.
The omitted length probabilities tend to zero by the zero-count law;
$W\ge1$. Nonplanar competitors are excluded on each fixed window.
This also makes any tie rule immaterial and proves planarity of every
systolic curve with probability tending to one. The first-arrival
density of the limiting marked process is \eqref{eq:FS1}.
\end{proof}

\section{Proofs for the tight length spectrum}\label{app:tight}

This appendix proves Theorem~\ref{thm:tight-spectrum} and
Corollary~\ref{cor:tight-boltzmann}, following the outline in
Subsection~\ref{subsec:tight-spectrum}. The scale $\alpha$ and the
intensity $\nu$ are those of \eqref{eq:tight-scale}. The tight
integration formula of Budd and Lions~\cite[Proposition~3.2]{BL}
expresses factorial moments of tight curve counts through tight
volumes of the cut surfaces; the inputs from
Part~\ref{part:volumes} are the shifted amplitudes of
Section~\ref{exact-amplitudes-and-the-harmonic-model}, the coefficient
extraction of Section~\ref{extracting-fixed-cusp-counts}
(Lemma~\ref{lem:point-mass} and the estimates
\eqref{eq:TL24}--\eqref{eq:TL31}), and the prefactor of
Proposition~\ref{prop:prefactor}. Subsection~\ref{subsec:tight-gf}
expresses the generating functions of tight volumes through the shifted
amplitudes (Lemma~\ref{lem:tight-gf}); Subsection~\ref{subsec:tight-count}
extracts the coefficient at the prescribed cusp count for connected
complements (Lemma~\ref{lem:tight-connected});
Subsection~\ref{subsec:tight-sep} shows that disconnected complements
are negligible (Lemma~\ref{lem:tight-separating});
Subsection~\ref{subsec:tight-proof} identifies the Poisson limit, and
Subsection~\ref{subsec:tight-boltzmann-proof} averages over the random
cusp count to obtain the corollary.

\subsection{Tight volume generating functions}\label{subsec:tight-gf}

This subsection proves Lemma~\ref{lem:tight-gf}, which identifies the
tight generating functions with sums of the shifted amplitudes
$U^Y_{g,p}$ of \eqref{eq:H6}, compares them with the ordinary marked
generating functions, and records their positive convolution
structure. The identity comes from the tight recursion of Budd
and Zonneveld~\cite{BZ}; the comparison and positivity come from the
factorial series \eqref{eq:TL15}, \eqref{eq:TL23} and the positive
factor \eqref{eq:TL33}. Both are used in the next two subsections:
the comparison transfers local limits from the marked functions to the
tight ones (proof of Lemma~\ref{lem:tight-connected}), and positivity
supplies the point-mass bounds (proof of
Lemma~\ref{lem:tight-separating}).
Put $P=2\pi^2$. Let $T_{g,p,k}(\mathbf L)$ be the tight volume with
$p$ protected geodesic boundaries and $k$ cusps that may be filled,
in the conventions of~\cite{Budd,BZ,BL}. For stable $(g,p)$ define
\[
 \mathscr T_{g,p}(x;\mathbf L)
 =\sum_{k\ge0}\frac{T_{g,p,k}(\mathbf L)}{P^{3g-3+p+k}}
                 \frac{x^k}{k!},\qquad F_{g,p}=\Phi_g^{(p)}.
\]
Unmarked functions are used only for $g\ge2$. Set
$\varrho(Y)=\sqrt{r(Y)/P}$ and $S(t)=\sinh(t)/t$, with $S(0)=1$.

\begin{lemma}[Tight generating functions]\label{lem:tight-gf}
Let $(g,p)$ be stable and $0\le Y<y_0$. The tight generating function is
expressed through the shifted amplitudes by the exact identity
\begin{equation}
 \mathscr T_{g,p}(X(Y);\mathbf L)
 =\sum_{\mathbf d}U^Y_{g,p}(\mathbf d)
       \prod_{i=1}^p\frac{L_i^{2d_i}}{P^{d_i}(2d_i+1)!},
 \label{eq:tight-insertion}
\end{equation}
which fixes both the normalization and the meaning of the protected
boundaries. Its ratio to the ordinary marked generating function satisfies
\begin{equation}
 \begin{split}
 1\le\frac{\mathscr T_{g,p}(X(Y);\mathbf L)}{F_{g,p}(X(Y))}
   &\le\prod_{i=1}^p S(\varrho(Y)L_i),\\
 \frac{\mathscr T_{g,p}(X(Y);\mathbf L)}{F_{g,p}(X(Y))}
   &=\prod_{i=1}^p S(\varrho(Y)L_i)
        \bigl(1+O_{p,B}(g^{-1})\bigr)
       \quad\text{if }\max_i\varrho(Y)L_i\le B,
 \end{split}
 \label{eq:tight-real}
\end{equation}
where the inequality holds at every stable genus. Finally, for each fixed
physical length vector,
\begin{equation}
 \mathscr T_{g,p}(x;\mathbf L)=V(x)^{2g-2+p}R_{g,p,\mathbf L}(x),
 \qquad [x^j]R_{g,p,\mathbf L}(x)\ge0,
 \label{eq:tight-positive}
\end{equation}
and the series converge for $0\le x<x_c$.
\end{lemma}

\begin{proof}
The identity \eqref{eq:tight-insertion} follows from~\cite[Proposition~12]{BZ}. Indeed, their root and
moments are $R=Y/P$ and $M_j=P^jX^{(j+1)}(Y)$. With $q=P/2$,
their shifted times are $q^{k-1}t'_k$, where the $t'_k$ are the
times in \eqref{eq:H8}. The dimension constraint gives the factor
$q^{3g-3+p-|\mathbf d|}$. Combining their external insertion
$2L^{2d}/(4^d d!)$ with the free-energy factor $2^{3g-3}$ and
the double factorial in \eqref{eq:H6} gives
\eqref{eq:tight-insertion}. This also covers genus zero with $p\ge3$.

Three special cases check the normalization in \eqref{eq:tight-insertion}.
The dependence on $\mathbf L$ is tested at $Y=0$. There $x=0$, no cusp
is present, tightness is vacuous and
$T_{g,p,0}(\mathbf L)=V_{g,p}(\mathbf L)$. By \eqref{eq:H5},
$U^0_{g,p}(\mathbf d)=[\tau_{d_1}\cdots\tau_{d_p}]_{g,p}/(4^{|\mathbf d|}P^{m})$
with $m=3g-3+p-|\mathbf d|$, so the right side of
\eqref{eq:tight-insertion} is
\[
 P^{-(3g-3+p)}\sum_{\mathbf d}[\tau_{d_1}\cdots\tau_{d_p}]_{g,p}
   \prod_{i=1}^p\frac{(L_i/2)^{2d_i}}{(2d_i+1)!}
 =\frac{V_{g,p}(\mathbf L)}{P^{3g-3+p}},
\]
by the expansion of $V_{g,p}(2\mathbf L)$ in
Section~\ref{exact-amplitudes-and-the-harmonic-model}.
Next let $\mathbf L=\mathbf 0$. Only $\mathbf d=\mathbf 0$ survives on
the right, and \eqref{eq:H8} gives
\[
 \mathscr T_{g,p}(X(Y);\mathbf 0)=U^Y_{g,p}(\mathbf 0)
 =\Phi_g^{(p)}(X(Y))
 =\sum_{k\ge0}\frac{V_{g,p+k}}{P^{3g-3+p+k}}\frac{X(Y)^k}{k!},
\]
that is, $T_{g,p,k}(\mathbf 0)=V_{g,p+k}$ for every $k$. This is the
correct value: a protected boundary of length zero is a cusp, every
other closed geodesic in its extended free homotopy class has positive
length, so the tightness condition is vacuous and the tight moduli
space is the full moduli space; Budd and Lions record
$T_{g,n,0}=T_{g,0,n}=V_{g,n}$ and
$\mathcal M^{\rm tight}_{g,n,p}(\mathbf 0,\mathbf 0)=\mathcal M_{g,n+p}(\mathbf 0,\mathbf 0)$
in \cite[Section~2]{BL}.
Finally let $g=0$ and $p=3$. Then $m=3g-3+p-|\mathbf d|=-|\mathbf d|$ in
\eqref{eq:H6}, so again only $\mathbf d=\mathbf 0$ contributes, with
$h=2g-2+p=1$, no $\tau_k$ insertions and
$\langle\tau_0^3\rangle_0=1$:
\[
 \sum_{\mathbf d}U^Y_{0,3}(\mathbf d)
 \prod_{i=1}^3\frac{L_i^{2d_i}}{P^{d_i}(2d_i+1)!}
 =U^Y_{0,3}(\mathbf 0)=\frac1{a(Y)}=\frac1{X'(Y)},
\]
independently of $\mathbf L$; here $a=X'$ as in
Section~\ref{exact-amplitudes-and-the-harmonic-model}. On the other
side, in the normalization of \cite{BZ,BL} with cusp weight $\mu=x/P$,
\[
 \mathscr T_{0,3}(x;\mathbf L)
 =\sum_{k\ge0}\frac{\mu^k}{k!}T_{0,3,k}(\mathbf L,\mathbf 0^k)
 =T_{0,3}(\mathbf L,\mu)=\frac1{M_0(\mu)},
\]
the last identity being the genus-zero initial case
$\mathcal P_{0,3}\equiv1$ of the tight recursion \cite[Theorem~4]{BZ},
as recalled in \cite[Theorem~2.2]{BL}. Their root $R(\mu)$ solves
$\mu=\frac{\sqrt R}{\sqrt2\pi}J_1(2\pi\sqrt{2R})$ and
$M_0(\mu)=\frac{d}{dr}\bigl[\frac{\sqrt r}{\sqrt2\pi}J_1(2\pi\sqrt{2r})\bigr]_{r=R(\mu)}$.
Since $X(y)=\sqrt y\,J_1(2\sqrt y)$ and $P=2\pi^2$, the function
differentiated is $X(Pr)/P$; hence $R(\mu)=Y/P$ when $\mu=X(Y)/P$,
and $M_0(\mu)=X'(Y)=a(Y)$, in agreement with the root and moments
quoted above. Both sides of \eqref{eq:tight-insertion} therefore equal
$1/a(Y)$ for $g=0$, $p=3$. In particular $T_{0,3,k}(\mathbf L)=V_{0,3+k}$
for all $\mathbf L$ and $k$: the tight volumes with three protected
boundaries in genus zero do not depend on the boundary lengths, which is
the content of $\mathcal P_{0,3}\equiv1$.

Equations~\eqref{eq:TL15} and~\eqref{eq:TL23} give \eqref{eq:tight-real}. For
the asymptotic estimate,
the factorial series with the extra weight
$|\mathbf d|(|\mathbf d|+1)$ is uniformly summable on the indicated
box. Thus the lengths may grow as $r(Y)^{-1/2}$.

For each fixed physical length vector, the proof of
\eqref{eq:TL33} applies to every amplitude in
\eqref{eq:tight-insertion}, which gives \eqref{eq:tight-positive}. The
finite sum in \eqref{eq:tight-insertion} and \eqref{eq:tight-real} justify
convergence throughout $0\le x<x_c$.
\end{proof}

\subsection{Prescribing the cusp count}\label{subsec:tight-count}

This subsection proves Lemma~\ref{lem:tight-connected}, the volume
ratio for a cut surface with connected complement, as it occurs in the
tight integration formula. We extract the coefficient at the
prescribed cusp count by a local limit theorem at a suitably chosen
cusp weight, or by a Poisson limit when the cusp count stays bounded;
the inputs are the real and Fourier estimates
\eqref{eq:TL24}--\eqref{eq:TL31} and \eqref{eq:TL40}--\eqref{eq:TL46}
of Section~\ref{extracting-fixed-cusp-counts}, transferred from the
marked generating functions to the tight ones by
Lemma~\ref{lem:tight-gf}. Fix a number $k$ of
cuts and put
$h_0=2g-2$. When $n\to\infty$, choose $x$ by
\[
 (h_0+1)\Psi(x)=n-1,\qquad
 \sigma^2=(h_0+1)\mathscr D\Psi(x),\qquad x=X(Y),
\]
and write $\alpha_x=2\varrho(Y)$. If $n$ stays bounded, we pass to a
subsequence on which $n$ is fixed and take instead $x=c/h_0$ with $c>0$.

\begin{lemma}[Connected cut configurations]\label{lem:tight-connected}
Fix $k$ and a compact box of scaled lengths $\mathbf t\in[0,\infty)^{2k}$.
With $x$ chosen as above, uniformly for $\mathbf t$ in the box,
\begin{equation}
 \frac{[z^n]\mathscr T_{g-k,2k}(z;\mathbf t/\alpha_x)}
      {[z^n]\Phi_g(z)}
 =(4r(Y))^k\prod_{i=1}^{2k}S(t_i/2)(1+o(1)).
 \label{eq:tight-connected}
\end{equation}
Moreover $\alpha_x/\alpha(n/g)\to1$.
\end{lemma}

\begin{proof}
Suppose first that $n\to\infty$. The denominator
$F_{g,1}=\Phi_g'$ has its local limit at $n-1$.
For the connected-complement numerator
$G(x)=\mathscr T_{g-k,2k}(x;\mathbf t/\alpha_x)$, the physical
lengths are held fixed when its argument is varied. Its model center
is $h_0\Psi(x)$, and
\[
 \frac{n-h_0\Psi(x)}\sigma
 =\frac{1+\Psi(x)}\sigma\longrightarrow0.
\]
The endpoint expansions give
\[
 \sup_{0<x<x_c}
 \frac{|\mathscr D\log r(Y(x))|}{\sqrt{\mathscr D\Psi(x)}}<\infty.
\]
Consequently a real tilt $x\mathrm e^{u/\sigma}$ changes $r$ by a
relative $O_u(g^{-1/2})$, uniformly up to either endpoint.
Equation~\eqref{eq:tight-real} therefore transfers the real Gaussian
limit from $F_{g-k,2k}$ to $G$. The positive factor
\eqref{eq:tight-positive} transfers the Fourier estimates
\eqref{eq:TL40}--\eqref{eq:TL46}. The two local probabilities at
$n$ and $n-1$ are asymptotic to $(\sqrt{2\pi}\sigma)^{-1}$.

We use the derivative only in the denominator: adding a protected cusp
to a tight surface would change which homotopies are allowed.
For a generating function $J$ with nonnegative coefficients write
$\mathbb P_x^J(N=m)=[z^m]J(z)x^m/J(x)$. The exact identity is
\begin{equation}
 \frac{[z^n]G(z)}{[z^n]\Phi_g(z)}
 =\frac n x\frac{G(x)}{\Phi_g'(x)}
       \frac{\mathbb P_x^G(N=n)}
            {\mathbb P_x^{\Phi_g'}(N=n-1)}.
 \label{eq:tight-coefficient}
\end{equation}
Equation~\eqref{eq:TL24} at $p=0$ shows that
$n\Phi_g(x)/(x\Phi_g'(x))\to1$. Iterating
\eqref{eq:TL31} now yields \eqref{eq:tight-connected}.

For bounded $n$, with $x=c/h_0$, equation~\eqref{eq:TL24} gives
\eqref{eq:TL30} also at $p=0$. Integrating that estimate between
$x$ and $xz$, and using $\Psi(x)=2x+O(x^2)$, proves convergence of
the probability generating functions to $\exp(2c(z-1))$.
The same holds for $G$ by \eqref{eq:tight-real}. The probabilities at
the fixed index $n$ thus converge to the same positive Poisson
probability $\mathrm e^{-2c}(2c)^n/n!$. In this case we use, instead of
\eqref{eq:tight-coefficient}, the undifferentiated identity
\[
 \frac{[z^n]G(z)}{[z^n]\Phi_g(z)}
 =\frac{G(x)}{\Phi_g(x)}\,
   \frac{\mathbb P_x^{G}(N=n)}{\mathbb P_x^{\Phi_g}(N=n)},
\]
which is meaningful also at $n=0$. This proves
\eqref{eq:tight-connected} in this case too.

Finally $\alpha_x/\alpha(n/g)\to1$. At finite density this follows
by continuity; at infinite density use $r(Y_\lambda)\sim
(3\sqrt5j_{0,1}/2)\lambda^{-1/2}$. This relative comparison is
needed because both scales tend to zero at infinite density.
\end{proof}

\subsection{Separating configurations}\label{subsec:tight-sep}

This subsection proves Lemma~\ref{lem:tight-separating}: tuples whose
union disconnects the filled surface contribute negligibly to the
tight integration formula. The mechanism is that, compared with the
parent, the product of the component generating functions loses one
inverse power of the genus for each additional component. The inputs
are the two-sided prefactor of Proposition~\ref{prop:prefactor}, a
convolution bound for the resulting Gamma-function profiles, and the
point-mass bound of Lemma~\ref{lem:point-mass} applied to the positive
factor of Lemma~\ref{lem:tight-gf}.

\begin{lemma}[Separating configurations]\label{lem:tight-separating}
Fix $k$ and $L$. Uniformly in $n$, with lengths scaled by $\alpha_x$,
\begin{equation}
 \mathbb E^{\mathrm{WP}}_{g,n}
 \bigl[\text{ordered separating tight $k$-tuples with scaled
 lengths at most $L$}\bigr]=o(1),
 \label{eq:tight-separating-count}
\end{equation}
where only simple, pairwise disjoint tuples are counted, and separating
means that their union disconnects the filled surface.
\end{lemma}

\begin{proof}
For a stable connected cut graph with $k$ edges and $Q$ vertices, let
$(g_i,p_i)$ be the component types after filling the cusps. No filled
component is a disk or an annulus: the first would make a cut curve
contractible on the filled surface, and the second would make two
distinct tight curves homotopic there, contrary to the uniqueness in
the definition of tightness. Hence every component is stable, and
\[
 p_i\ge1,\quad h_i=2g_i-2+p_i\ge1,\quad
 \sum_i h_i=h_0,\quad \sum_i p_i=2k,\quad
 \sum_i g_i=g-k+Q-1.
\]
For fixed positive $p$, Proposition~\ref{prop:prefactor} implies the
two-sided bound
\begin{equation}
 F_{g,p}(X(Y))\asymp_p
 4^{-g}\mathcal Z(Y)^{-h}r(Y)^{p/2}f(h),
 \qquad f(h)=\frac{\Gamma(h)}{\sqrt h},\quad h=2g-2+p.
 \label{eq:tight-gamma}
\end{equation}
The constants are uniform in $Y$. Bounded stable genera are included
by positivity and continuity of the normalized finite amplitudes on
$[0,y_0]$. For $p=0$, $g\ge2$, the same bound follows from $p=1$
and \eqref{eq:TL24}.

The identities for the component types cancel all powers of
$\mathcal Z$ in the product-to-parent ratio. They leave
\[
 r(Y)^{-k}\frac{\prod_i F_{g_i,p_i}(X(Y))}{\Phi_g(X(Y))}
 \le C_k\frac{\prod_i f(h_i)}{f(h_0)}.
\]
To sum over genera, observe that
$\sum_{a=1}^{\ell-1}f(a)f(\ell-a)\le C f(\ell-1)$.
Indeed,
\[
 \frac{\Gamma(a)\Gamma(\ell-a)}{\Gamma(\ell)}
 =\frac{1}{(\ell-1)\binom{\ell-2}{a-1}},
\]
the reciprocal binomial coefficients have uniformly bounded sum,
and $\sqrt{\ell}/\sqrt{a(\ell-a)}$ is bounded. Iteration gives
$f^{*Q}(\ell)\le C_Q f(\ell-Q+1)$. Since the marking allocations and
gluings have bounded multiplicity for fixed $k$, it follows that
\begin{equation}
 r(Y)^{-k}\sum_{\substack{\text{stable cut types}\\Q\text{ components}}}
 \frac{\prod_i F_{g_i,p_i}(X(Y))}{\Phi_g(X(Y))}
 \le C_k g^{-(Q-1)},\qquad 2\le Q\le k+1.
 \label{eq:tight-separating-real}
\end{equation}

For scaled lengths in a fixed compact box,
\eqref{eq:tight-real} bounds each boundary factor uniformly, including
on small-genus components. Every product in this positive sum has
the factor $V^{h_0}$ by \eqref{eq:tight-positive}. Its largest tilted
point probability is at most $C/\sigma$, by Lemma~\ref{lem:point-mass}. Applying \eqref{eq:tight-coefficient}, against the
denominator probability at least $c/\sigma$, therefore transfers
\eqref{eq:tight-separating-real} to coefficients. For bounded $n$,
use $x=c/h_0$ and the undifferentiated identity above, whose
denominator probability has a positive limit;
the numerator probability is at most one.

In the tight integration formula~\cite[Proposition~3.2]{BL}, the
physical normalization costs $P^{-k}$ and the length Jacobian is
$\alpha_x^{-2k}=(P/(4r))^k$. These cancel the remaining $r^k$.
The symmetry factors are at most one. This proves
\eqref{eq:tight-separating-count}.
\end{proof}

\subsection{Proof of Theorem~\ref{thm:tight-spectrum}}\label{subsec:tight-proof}

Lemmas~\ref{lem:tight-connected} and~\ref{lem:tight-separating} give
the Poisson limits of the moments of disjoint tight tuples. The proof
below relates these to the moments of the actual counts, using the
deterministic multicurve constructions of Budd and Lions~\cite{BL} to
show that intersecting tight curves are excluded at the level of
expectations, and then passes from moments to the point process.

\begin{proof}[Proof of Theorem~\ref{thm:tight-spectrum}]
For fixed disjoint bounded intervals in $[0,\infty)$, count ordered tuples
of simple, pairwise disjoint tight curves in those intervals.
The connected-complement symmetry is $2^{-k}$ in the tight integration
formula. Lemma~\ref{lem:tight-connected}, with each cut length
appearing on two boundary components, gives the limiting integrand
\[
 \prod_{i=1}^k\frac12 t_i S(t_i/2)^2\,\mathrm dt_i
 =\prod_{i=1}^k\frac{\cosh t_i-1}{t_i}\,\mathrm dt_i.
\]
Lemma~\ref{lem:tight-separating} removes all other cut types.
Thus the disjoint-tuple moments have the Poisson limits of every
fixed order.

Here the disjoint-tuple moments must be related to the actual count
moments: a small probability of intersections alone would not suffice.
Choose $L$ above all interval endpoints, and let $A_L$ be the event that
scaled lengths at most $L$ include a
self-intersecting tight curve or two intersecting tight curves. Write
$\widehat N_{\mathbf r}$ for the disjoint-tuple count in the prescribed
intervals, with multiplicities $\mathbf r=(r_1,\ldots,r_j)$.
We use the deterministic constructions in
\cite[proofs of Propositions~5.2--5.3 and Theorem~1.2, Section~5.2]{BL}, rather than their probabilistic
conclusions in the near-critical model. These constructions use only
curve surgery at intersection points, the subadditivity of length
under surgery and a lower bound on the scaled tight systole; they
therefore apply after rescaling by any positive length scale, with the
explicit bounds on the number and lengths of the added curves given
in~\cite[Section~5.1]{BL}.

Fix $\varepsilon>0$ below the interval endpoints. On $A_L$, outside
the event of a tight systole shorter than $\varepsilon$, the construction
extends a counted tuple to a separating tight multicurve. Its number of
components and scaled lengths are bounded in terms of
$\varepsilon,L,\mathbf r$, independently of $g,n$. If there is a shorter
tight systole, either it adds a disjoint curve in $[0,\varepsilon]$, or
the tuple again extends to a separating multicurve. The original ordered
tuple is retained in these extensions. Thus the separating terms vanish
by \eqref{eq:tight-separating-count}, and the mixed disjoint-tuple limit gives
\[
 \limsup_{g\to\infty}\mathbb E[
   \widehat N_{\mathbf r}\mathbf1_{A_L}]
 \le \nu([0,\varepsilon])\prod_{i=1}^j\nu(I_i)^{r_i},
\]
where $I_i$ are the prescribed intervals. Letting $\varepsilon\downarrow0$
removes the right side. The same construction without a marked tuple
gives $\mathbb P(A_L)\to0$. Conditional on $A_L^c$, factorial counts
are exactly disjoint-tuple counts. Their limiting moments therefore
remain unchanged, and the moment-determinate Poisson law gives joint
count convergence.

The first-moment bound for simple tight curves in $[0,\varepsilon]$
tends to $\nu([0,\varepsilon])\to0$. A shortest tight curve is
simple: resolving a self-intersection yields two closed curves, at
least one of which is noncontractible on the filled surface, and its
tight representative is strictly shorter. So the same bound controls the event of any tight curve in
this window. This extends the convergence to $[0,\infty)$ and gives
tightness on compact intervals. Finally the disjoint-tuple moments
bound the number of curves in a fixed window in probability.
For any fixed upper bound on that number,
\eqref{eq:tight-separating-count} rules out every disconnected
complement. Letting the upper bound increase proves the last claim.
\end{proof}

\subsection{Proof of Corollary~\ref{cor:tight-boltzmann}}\label{subsec:tight-boltzmann-proof}

Under the Boltzmann measure $\mathbb P^{\mathrm{WP}}_{g,\mu}$ of
Subsection~\ref{subsec:tight-spectrum} the cusp count is random. We
show, using the real central limit theorem of
Section~\ref{extracting-fixed-cusp-counts}, that the scale evaluated at
the random density is asymptotic in probability to the deterministic
scale $a_g$, and then average the conditional limits of
Theorem~\ref{thm:tight-spectrum} over the cusp count.

\begin{proof}[Proof of Corollary~\ref{cor:tight-boltzmann}]
Write $x=P\mu_g$ and $h=2g-2$. The case $x=0$ has $N=0$.
Otherwise \eqref{eq:TL24} at $p=0$ gives
\[
 \mathbb E N=h\Psi(x)+O(\Psi(x)).
\]
If $h\Psi(x)$ stays bounded, $N$ is tight and both $N/g$ and
$2\Psi(x)$ tend to zero. If $h\Psi(x)\to\infty$, integrating
\eqref{eq:TL30} at $p=0$ gives the same real central limit theorem
as \eqref{eq:TL39}, hence
\[
 N=h\Psi(x)+O_{\mathbb P}\!\left(\sqrt{h\mathscr D\Psi(x)}\right).
\]
The endpoint estimates \eqref{eq:TL36} imply
$\mathscr D\Psi/\Psi^2\le C(1+1/\Psi)$. It follows that
$N/(h\Psi(x))\to1$ in probability in this second case.
Continuity of $\alpha$ at finite density and its
$\lambda^{-1/4}$ asymptotic at infinity therefore give, in both cases,
\[
 \frac{\alpha(N/g)}{\alpha(2\Psi(x))}
 =\frac{\alpha(N/g)}{a_g}\longrightarrow1
 \quad\text{in probability}.
\]
Subsequences cover arbitrary choices of $x$.

For each fixed nonnegative compactly supported continuous test function,
the conditional Laplace-functional convergence in
Theorem~\ref{thm:tight-spectrum} is uniform in $n$: a sequence violating
uniformity would contradict that theorem. Averaging over $N$ gives the
Poisson limit with random normalization $\alpha(N/g)$.
The displayed scale comparison and continuity of dilation of point
measures replace this by $a_g$. For $\mu_g\to\mu<\mu_{\rm cr}$, a change of
scale gives the stated unscaled intensity. At the critical endpoint,
$r\sim3(y_0-Y)$ and $x_c-x\sim(-X''(y_0))(y_0-Y)^2/2$ show that
$a_g$ is asymptotic to the scale in~\cite[Theorem~1.2]{BL}.
\end{proof}

\section{Proofs for the critical length spectrum}\label{app:critical}

This appendix proves Theorem~\ref{thm:critical-cluster}, following the
outline in Subsection~\ref{sec:critical-cluster}. Throughout,
$n/\sqrt g\to c\in(0,\infty)$, and $P_s$, $K_s$, $\mu_c$ and $\nu$ are as
in \eqref{eq:cluster-intensities}. The volume inputs are the bounds of
Mirzakhani--Zograf and Shen--Wu~\cite{MZ,SW}, the zero-density case
of Theorem~\ref{thm:volume} and the boundary factorization of
Proposition~\ref{prop:bounded-boundaries}; the geometric inputs are
the intersection and filling arguments of \cite{MP,Hide}. The four
subsections bound the contribution of each cut graph
(Lemma~\ref{lem:critical-cuts}), reduce intersecting simple tuples and
nonsimple curves to cut graphs (Lemmas~\ref{lem:primitive-count}
and~\ref{lem:critical-complex}), compute the joint Poisson limit of
nonseparating lengths and waists of two-cusp pairs of pants
(Lemma~\ref{lem:critical-skeleton}), and finally attach to each waist
the primitive length spectrum of its pair of pants. The first three
lemmas are also used outside this appendix: Lemmas~\ref{lem:critical-cuts}
and~\ref{lem:critical-complex} isolate the pants in
Section~\ref{sec:critical-eigenvalues}, and
Lemma~\ref{lem:primitive-count} is used in
Subsection~\ref{sec:critical-cluster} and in
Appendices~\ref{app:positive-density-orbits},
\ref{app:fixed-window-proofs} and~\ref{app:rooted-local-law}.

\subsection{Ordinary cuts and cusp allocations}

Short multicurves are organized by their cut graphs: the vertices are the
components of the complement, labeled by their genera, and the edges are
the curves. Lemma~\ref{lem:critical-cuts} bounds the expected number of
multicurves with a given cut graph, using only Mirzakhani integration
and the volume bounds of \cite{SW}. It shows that, in the limit, only
nonseparating curves and boundaries of two-cusp pairs of pants
contribute; the later subsections compute these contributions jointly.

Fix a connected cut graph with $k\ge1$ edges and $Q$ vertices;
loops and multiple edges are allowed. At vertex $i$ let $a_i$ be the
genus and $p_i\ge1$ its degree. Then
\[
 \sum_i a_i=g-k+Q-1,\qquad \sum_i p_i=2k.
\]
Partition the vertices into $u$ of type $(0,1)$, $v$ of type $(0,2)$,
and $w$ with $h_i=2a_i-2+p_i\ge1$. Put $h_0=2g-2$. For large
$g$, $w\ge1$, and the latter complexities sum to $H=h_0+u$.
A component of type $(0,1)$ needs at least two cusps, and one of
type $(0,2)$ at least one.

\begin{lemma}[Critical cut estimate]\label{lem:critical-cuts}
Fix $k$ and a bound on the cut lengths. At $n/\sqrt g\to c>0$,
the expected number of cut multicurves of the specified graph and
vertex partition, summed over all admissible genera and cusp allocations,
is $O(g^{1-w-v/2})$. If a specified $(0,1)$ component has at least
three cusps, the bound gains a factor $g^{-1/2}$.
\end{lemma}

\begin{proof}
Write $h_n=2g-2+n$ for the complexity of the cusped surface. The volume bounds in \cite[Lemmas~18--19]{SW} give
\[
 V_{a,b}\le C\frac{(4\pi^2)^{d-1}\Gamma(d)}{\sqrt d},
 \quad d=2a-2+b\ge1,\qquad
 V_{g,n}\ge c_0\frac{(4\pi^2)^{h_n-1}\Gamma(h_n)}{\sqrt{h_n}}.
\]
The first bound includes $d=1$ after enlarging $C$ for $(0,3)$ and
$(1,1)$; the second uses $n=O(\sqrt g)$. If vertex $i$ receives
$q_i$ cusps, put $d_i=h_i+q_i$. Then $\sum_i d_i=h_n$ and
$\sqrt{h_n}/\prod_i\sqrt{d_i}\le\sqrt Q$. Mirzakhani integration,
the bounded-boundary estimate \cite[Lemma~3.2]{Mir}, and the
multinomial allocation of cusp labels therefore bound the expected
cut count $E_{g,n}$ by
\[
 E_{g,n}\le\frac{C_{k,c}\,n!}{\Gamma(h_n)}
       \sum_{\substack{\text{admissible genera}\\q_1+\cdots+q_Q=n}}
       \prod_i\frac{\Gamma(h_i+q_i)}{q_i!}.
\]
Here and below the constants may also depend on the length bound.
The powers of $4\pi^2$ leave only the fixed factor $(4\pi^2)^{1-Q}$.

The cusp sums at the two unstable vertex types are bounded using
\begin{equation}
 \phi_{0,1}(z)=\sum_{q\ge2}\frac{z^q}{q(q-1)}=\frac{z^2}{2}+O(z^3),
 \qquad \phi_{0,2}(z)=\sum_{q\ge1}\frac{z^q}{q}=z+O(z^2).
 \label{eq:cluster-unstable}
\end{equation}
At a stable vertex their generating function is
$\Gamma(h)(1-z)^{-h}$. Moreover
\[
 \sum_{\substack{h_1+\cdots+h_w=H\\h_i\ge1}}
       \prod_i\Gamma(h_i)\le C_w\Gamma(H-w+1).
\]
For $w=2$ this follows by dividing by $\Gamma(H-1)$ and summing
$\binom{H-2}{h_1-1}^{-1}$; iteration gives the general bound.
Dropping genus parity constraints only increases the sum.
Let $a_m=[z^m]\phi_{0,1}(z)^u\phi_{0,2}(z)^v$. Expanding $(1-z)^{-H}$ gives
\[
 E_{g,n}\le C_{k,c}\frac{\Gamma(H-w+1)}{\Gamma(H)}
       \sum_{m\le n}a_m(n)_m\frac{\Gamma(h_n+u-m)}{\Gamma(h_n)}.
\]
Uniformly for $m\le n=O(\sqrt{h_n})$, the finite-product formula for
the Gamma ratio gives
$\Gamma(h_n+u-m)/\Gamma(h_n)\le C_{k,c}h_n^{u-m}$: the logarithmic
error is $O((m+1)^2/h_n)$. Since $(n)_m\le n^m$, positivity yields
\begin{equation}
 E_{g,n}\le C_{k,c}h_n^{u-w+1}\phi_{0,1}(n/h_n)^u\phi_{0,2}(n/h_n)^v
          =O(g^{1-w-v/2}).
 \label{eq:critical-graph-bound}
\end{equation}
Requiring three cusps at a specified $(0,1)$ vertex replaces one
$\phi_{0,1}(z)$ by $\phi_{0,1}(z)-z^2/2$ and gains $n/h_n=O(g^{-1/2})$.
\end{proof}

Only $w=1$, $v=0$ can survive: the graph then has one large-genus
vertex, loops at that vertex, and two-cusp pairs of pants attached as
leaves. All other graphs, and all leaves with more than two cusps, vanish.

\subsection{Excluding complex short subsurfaces}\label{subsec:critical-exclusion}

We now reduce intersecting simple tuples, and nonsimple curves that do
not fill two-cusp pairs of pants, to the cut graphs just estimated
(Lemma~\ref{lem:critical-complex}). Since we need expectation bounds,
not only a bound on the probability of a bad configuration, we give
the multiplicity argument explicitly. It rests on a uniform bound for
the number of short primitive geodesics on a surface of bounded
topology (Lemma~\ref{lem:primitive-count}), on the local intersection
argument of Mirzakhani--Petri~\cite{MP} extended to cusp
neighborhoods, and on the filling construction of Hide~\cite{Hide}.

\begin{lemma}[Bounded primitive count]\label{lem:primitive-count}
For fixed $K,L$, a finite-area hyperbolic surface $S$ with geodesic
boundary and $0<-\chi(S)\le K$ has at most $C(K,L)$ primitive
unoriented geodesics of length at most $L$, uniformly in its metric.
\end{lemma}

\begin{proof}
Double $S$ along its boundaries. The double has area at
most $4\pi K$. Fix a Margulis constant $\varepsilon>0$. A primitive
geodesic wholly in a thin component is its core; there are boundedly
many such cores. Every other geodesic meets the thick part. A maximal
$\varepsilon/4$-separated net there has bounded cardinality by area
packing. At each net point, orbit points of a lift are uniformly
separated, so hyperbolic ball packing bounds the number of based
loops of length at most $L+\varepsilon$. Connecting a short geodesic
to a net point proves the assertion. Primitivity is preserved upon
inclusion in the double or in an ambient surface: a root would have
the same geodesic image and would already be a closed curve in $S$.
\end{proof}

\begin{lemma}[Complex short configurations]\label{lem:critical-complex}
Fix $r$ and $L$. As $g\to\infty$ with $n/\sqrt g\to c\in(0,\infty)$,
\begin{equation}
 \mathbb E[\text{intersecting ordered simple $r$-tuples of lengths
 at most $L$}]=O_{r,L,c}(g^{-1/2}),
 \label{eq:critical-intersections}
\end{equation}
\begin{equation}
 \mathbb E[\text{nonsimple primitive curves of length at most $L$
 not filling two-cusp pants}]=O_{L,c}(g^{-1/2}),
 \label{eq:critical-bad-curves}
\end{equation}
and the expected number of simple separating curves of length at most $L$
that are not waists of two-cusp pairs of pants tends to zero.
\end{lemma}

\begin{proof}
For fixed $r,L$, a tuple of $r$ distinct simple geodesics of length
at most $L$ has bounded total intersection number, on every finite-area
hyperbolic surface. The segment argument in the proof of
\cite[Proposition~4.5]{MP} rests on \cite[Lemma~4.4]{MP}: if two
geodesic segments $I,J$ have common endpoints, meet transversally at
at least one of them, and $I$ lies on a closed geodesic of length at
most $L$, then $\ell(I)+\ell(J)$ is bounded below by a positive
constant depending only on $L$. On a cusped surface this holds
whenever the loop $I\cup J$ is not peripheral; here $I$ and $J$ lie on
two distinct geodesics of the tuple, which meet transversally. Fix $Y=4$ and write $\{y>Y\}$ for
the cusp neighborhoods in width-one cusp coordinates. The larger
regions $\{y>2\}$ are embedded and pairwise disjoint by
\cite[Lemma~5.2]{HuangGap}. A
peripheral loop through a point of height $y\le Y$ has length at least
$2\operatorname{arcsinh}(1/(2Y))$, so the argument of~\cite{MP} bounds
the number of intersection points outside the cusp neighborhoods. For
those inside, note that no closed geodesic lies entirely in a cusp
neighborhood, and that two consecutive excursions of a geodesic into
$\bigcup\{y>Y\}$ are separated by length at least
$2\log(Y/2)=2\log2$; a
closed geodesic of length at most $L$ therefore makes at most
$1+L/(2\log2)$ excursions. Each excursion lifts to an arc of a
Euclidean semicircle centered on the real axis, with endpoints on
$\{y=Y\}$ and radius at most $Ye^{L/2}$; two distinct semicircles meet
at most once, and a semicircle meets at most $4Ye^{L/2}+1$ of the
integer translates of another. Hence the total intersection number of
the tuple is bounded by a constant depending only on $r$ and $L$. Thicken each nonannular
component of their union, fill complementary disks and peripheral
once-punctured disks, and tighten the essential boundaries. The
resulting pieces have bounded topology and total boundary length at
most $2rL$. Keep one ambient geodesic for each parallel boundary class.
Empty annuli between parallel copies collapse; their two incidences
are the two sides of a single edge, possibly a loop. An annulus with
cusps does not collapse and is a $(0,2)$ component counted among the
$v$ vertices.
Isolated tuple curves are retained as additional cut edges, unless
they are already one of the tightened boundaries. Minimal intersection
ensures that these boundaries are disjoint from the original tuple.

There are finitely many resulting graph shapes and markings at fixed
$r,L$. For a fixed cut multicurve, Lemma~\ref{lem:primitive-count}
bounds the choices of original curves in each marked piece. Hence the
tuple count is bounded pointwise by a constant times a finite sum of
marked cut-multicurve counts. An intersecting simple tuple cannot fill
a two-cusp pair of pants. Its marked component either has $h_i\ge1$, forcing
a second stable vertex to carry the unbounded genus, or has type
$(0,2)$, or has type $(0,1)$ with at least three cusps. The bounds in
Lemma~\ref{lem:critical-cuts} give \eqref{eq:critical-intersections}.

For a single nonsimple primitive curve of length at most $L$, use
the filling construction of \cite[Definition~4.10 and Lemma~4.11]{Hide}.
The filled piece has boundary length at most $2L$ and area at most
$4L$, so its topology is bounded. The same recording retains formal
boundary sides that tighten to the same ambient curve as sides of a
cut edge. Lemma~\ref{lem:primitive-count} controls the number of curves
filling each piece. The preceding classification therefore gives
\eqref{eq:critical-bad-curves}.
A filled piece cannot be the whole ambient surface for large $g$,
since its topology is bounded. Finally, the one-edge case of
Lemma~\ref{lem:critical-cuts} shows that simple separating curves
of bounded length other than waists of two-cusp pairs of pants vanish
in expectation.
\end{proof}

\subsection{Joint nonseparating and pants counts}

With intersecting tuples and the other cut graphs excluded, it remains
to compute the mixed factorial moments of the two surviving counts
(Lemma~\ref{lem:critical-skeleton}). This is where the two-sided
volume ratios enter: the cusp ratios come from \cite{MZ}, the genus
ratios from the zero-density case of Theorem~\ref{thm:volume}, and
the boundary factors from Proposition~\ref{prop:bounded-boundaries}.
Let $\mathcal N_g$ count simple nonseparating lengths and let
$\mathcal Q_X$ count the waist lengths of two-cusp pairs of pants.

\begin{lemma}[Joint Poisson limit]\label{lem:critical-skeleton}
As $g\to\infty$ with $n/\sqrt g\to c\in(0,\infty)$,
\begin{equation}
 (\mathcal N_g,\mathcal Q_X)\Longrightarrow
 (\mathcal N,\mathcal Q_c)
 \label{eq:critical-skeleton}
\end{equation}
jointly in the vague topology on $[0,\infty)$, where $\mathcal N$ and
$\mathcal Q_c$ are independent Poisson point processes of intensities $\nu$
and $\mu_c$.
\end{lemma}

\begin{proof}
In the mixed
factorial moments of $(\mathcal N_g,\mathcal Q_X)$, \eqref{eq:critical-intersections} removes intersecting
tuples. Lemma~\ref{lem:critical-cuts} removes every disjoint cut type
except the graph with a single large component, $k$ nonseparating
loops and $r$ two-cusp leaves. For prescribed bounded intervals for
the $k+r$ lengths, its exact contribution is
\begin{equation}
 \frac{(n)_{2r}}{2^{r+k}V_{g,n}}
 \int V_{g-k,n-r+2k}(\mathbf s,\mathbf t,\mathbf t,0^{n-2r})
        \prod_{j=1}^r s_j\,\mathrm ds_j
        \prod_{i=1}^k t_i\,\mathrm dt_i.
 \label{eq:critical-mixed-moment}
\end{equation}
Here $(n)_{2r}$ is the falling factorial. The factor $2^{-r}$ counts
unordered pairs of cusp labels at the ordered pants leaves, as in
\cite[Lemma~28]{SW}; $2^{-k}$ exchanges the sides of each
nonseparating cut, as in \cite[Proposition~4.2]{MP}. There is no
factor $r!$ or $k!$, since both tuples are ordered.

The volume ratio is a product of $r$ cusp ratios and $k$ genus ratios.
The cusp ratios satisfy $V_{g,j+1}/V_{g,j}\sim8\pi^2g$ for $j\le n$ by
\cite[Lemma~5.1(ii)]{MZ}, since $n=o(g)$.
For the genus ratios, \cite[Lemma~5.1(iii)]{MZ} gives only an upper bound,
and the two-sided limit at $n\asymp\sqrt g$ is the zero-density case
$\mathcal B(0)=1$ of Theorem~\ref{thm:volume}. Together,
\[
 \frac{V_{g-k,n-r+2k}}{V_{g,n}}\sim(8\pi^2g)^{-r}.
\]
Proposition~\ref{prop:bounded-boundaries}, applied at zero limiting
density, where $\mathfrak b_0(L)=\sinh(L/2)/(L/2)$, supplies, uniformly
on the fixed length box, the boundary factor
$\prod_j S(s_j/2)\prod_i S(t_i/2)^2$, with $S(t)=\sinh(t)/t$; compare
Shen--Wu's bounded-boundary estimate~\cite[Proposition~25]{SW}, which
they attribute essentially to \cite[Proposition~3.1]{MP}
and~\cite[Lemma~20]{NWX}.
Consequently the limiting factors in \eqref{eq:critical-mixed-moment}
are
\[
 \frac12 tS(t/2)^2\,\mathrm dt=\nu(\mathrm dt),\qquad
 \frac{c^2}{16\pi^2}sS(s/2)\,\mathrm ds=\mu_c(\mathrm ds).
\]
Hence, for disjoint bounded intervals, the mixed factorial moments of
the counts of $\mathcal N_g$ and $\mathcal Q_X$ converge to those of
independent Poisson variables with the corresponding $\nu$- and
$\mu_c$-masses. A multivariate Poisson law is determined by its
moments, so these counts converge jointly in distribution. The
intensities $\nu$ and $\mu_c$ are diffuse and finite on bounded
intervals, including neighborhoods of zero, where their masses are
$O(\varepsilon^2)$; convergence of the counts on finite families of
disjoint bounded intervals is therefore equivalent to vague convergence
in distribution of the pair of point
processes~\cite[Chapter~4]{Kallenberg}. This proves the lemma.
\end{proof}

\subsection{Proof of Theorem~\ref{thm:critical-cluster}}

By Lemmas~\ref{lem:critical-cuts} and~\ref{lem:critical-complex}, with
probability tending to one the primitive length spectrum below a
fixed bound consists of the simple nonseparating lengths together with
the primitive length spectra of the two-cusp pairs of pants with short
waists. The proof combines this with the joint Poisson limit of
Lemma~\ref{lem:critical-skeleton} and the continuity of the pants
spectra in the waist length.

\begin{proof}[Completion of the proof of Theorem~\ref{thm:critical-cluster}]
Fix $L<\infty$. By Lemma~\ref{lem:critical-complex}, with probability
tending to one every nonsimple primitive curve of length at most $L$
fills a two-cusp pair of pants with waist at most $2L$. By
\eqref{eq:critical-intersections}, the pants waists up to $2L$ are
pairwise disjoint with probability tending to one. The pants interiors
are then disjoint as well. Other simple separating curves up to $L$
are absent with probability tending to one. On this event, the full
primitive counting measure restricted to $[0,L]$ is exactly
\[
 \left.\left(\mathcal N_g+
       \sum_{s\in\mathcal Q_X\cap[0,2L]}K_s\right)\right|_{[0,L]}.
\]
Every interior primitive geodesic of $P_s$ is nonsimple and fills it.
The filling boundary bound implies $s\le2\ell$ for such a geodesic
of length $\ell$. Thus pants with waist greater than $2L$ cannot
contribute, and the limiting sum is locally finite.

For $\varepsilon>0$, the family $P_s$, $\varepsilon\le s\le2L$,
is compact in its fixed-topology moduli space. Lengths of marked
classes vary continuously, and only finitely many such classes have
length at most $L$ somewhere on this compact family. For the latter
assertion, use a finite cover by uniformly bilipschitz markings,
with standard neighborhoods at the two cusps, and the finiteness of
each center metric's length spectrum below a fixed bound. Hence
$s\mapsto\langle f,K_s\rangle$ is continuous away from zero for
continuous compactly supported $f$. Write $F(s)=\langle f,K_s\rangle$
and choose $L$ larger than the support of $f$; then $F$ is bounded by
Lemma~\ref{lem:primitive-count} and vanishes for $s\ge2L$. For
$\varepsilon>0$ the map
$(\eta,\xi)\mapsto\langle f,\eta\rangle+\int_{[\varepsilon,\infty)}F\,\mathrm d\xi$
is vaguely continuous at every pair with $\xi(\{\varepsilon\})=0$,
hence almost surely at $(\mathcal N,\mathcal Q_c)$, because $\mu_c$
has no atoms. By \eqref{eq:critical-skeleton} and the continuous
mapping theorem,
\[
 \langle f,\mathcal N_g\rangle
 +\sum_{s\in\mathcal Q_X,\ s\ge\varepsilon}F(s)
 \ \Longrightarrow\
 \langle f,\mathcal N\rangle
 +\sum_{s\in\mathcal Q_c,\ s\ge\varepsilon}F(s).
\]
On the event described above, the left side differs from the pairing
of $f$ with the primitive length spectrum only if a waist is shorter
than $\varepsilon$, which has probability at most
$\mathbb E\mathcal Q_X([0,\varepsilon])$. By the case $k=0$, $r=1$ of
\eqref{eq:critical-mixed-moment} this expectation converges to
$\mu_c([0,\varepsilon])=O(\varepsilon^2)$, and the right side
converges almost surely to $\langle f,\Lambda_c\rangle$ as
$\varepsilon\downarrow0$. Hence the Laplace functionals converge,
which proves \eqref{eq:critical-cluster}.
The independent Poisson Laplace functionals give
\eqref{eq:cluster-laplace}.
\end{proof}

\section{Trace estimates on bounded critical envelopes}
\label{app:critical-trace}

This appendix extends Huang's orbital trace estimate~\cite{HuangGap}
to cusp counts of order at most $\sqrt g$. Its two outputs are the
polynomial trace estimate \eqref{eq:critical-polynomial-error} and the
local expected-count bound for eigenvalues below a fixed threshold,
both used in the proof of Proposition~\ref{prop:critical-mean};
Appendix~\ref{app:sublinear-trace} restates the corresponding inputs
from \cite{HuangGap} for the sublinear range, with tail bounds proved
there. Huang's
stated additive genus threshold contains $C(n+1)^2$, so it cannot
simply be invoked for an arbitrary fixed limit $n/\sqrt g=c$. We
retain his canonical coefficients and modify only the
complementary-volume tail estimates.
Subsection~\ref{rev:subsec-huang-inputs} restates, in the notation of
this paper, every result of \cite{HuangGap} that the appendix uses, and
records which two of his estimates are replaced.
Subsection~\ref{rev:subsec-critical-tails} proves the replaced tail
bounds on the new envelope, Subsection~\ref{rev:subsec-critical-records}
checks the retained records and integrates over orbital types, and
Subsection~\ref{rev:subsec-critical-polynomial} passes to polynomial
tests of growing degree and to eigenvalue counts.
All numbered references to \cite{HuangGap} below refer to version~1.

Fix $M\ge1$ and write $\mathsf s=n+1\le M\sqrt g$; the letter $s$ is
reserved for waist lengths, as in Section~\ref{sec:critical-eigenvalues}.
Constants may depend on
$M$ and on the fixed multiplier, but not on $g,n$ or the expansion
order. There are $\theta>0$ and $g_M$ such that for $g\ge g_M$ and
$\rho=R+\lceil T\rceil+1\le g^\theta$, the orbital expansion obeys
\begin{equation}\label{eq:critical-orbital-remainder}
 \left|\mathbb E\mathcal O_H-\Omega(H)
                -\sum_{j=1}^R g^{-j}\kappa_{j,n}(H)\right|
 \le g^{-R-1}\mathsf s^{R+2}(C_M\rho)^{C_M\rho}
                  e^{3T}\|H\|_\infty,
\end{equation}
for bounded even Borel $H$ supported in $[-T,T]$. The functionals $\mathcal O_H$,
$\Omega$ and $\kappa_{j,n}$ are recalled in
Definition~\ref{rev:def-huang-orbital}; the coefficients are precisely
those in \cite[Theorem~4.8]{HuangGap}. We justify the changed parameter
domain, including the retained kernels, before transferring this estimate
to the positive spectral trace.

\subsection{Inputs from Huang's expansion}\label{rev:subsec-huang-inputs}

This subsection contains no new mathematics. It fixes notation,
introduces the orbital functional and the carrier bookkeeping of
\cite{HuangGap} (Definitions~\ref{rev:def-huang-orbital}
and~\ref{rev:def-huang-carriers}), and restates as
Lemmas~\ref{rev:lem-huang-envelope}--\ref{rev:lem-huang-detector} the
estimates of that paper that the later subsections use, each with the
hypotheses under which it is proved there or, for
Lemmas~\ref{rev:lem-huang-first-moment} and~\ref{rev:lem-huang-transfer},
under the hypotheses of the deterministic parts of its proofs that we
use. The paragraph at the end of
the subsection lists the two estimates whose genus thresholds are
quadratic in the cusp count and which are therefore replaced on the
envelope $\mathsf s\le M\sqrt g$.

Throughout, $g\ge3$ and $n\ge0$ are integers, $N=2g+n-3$ is the
factorial index of $V_{g,n}$, and
$\mathcal A_{g,n}=2\pi(2g-2+n)=\operatorname{area}(X)$. The fixed
multiplier is that of Proposition~\ref{prop:critical-mean}: $h=b_0*b_0$
with $\operatorname{supp}h\subset[-L,L]$ and $\|h\|_1=1$,
$f=\widehat h$ with $\widehat H(r)=\int_{\mathbb R}H(u)e^{-iru}\,du$,
and $a=f(i/2)>1$; thus $0\le f\le1$ on $\mathbb R$ and
$1\le f(i\tau)\le a$ for $0\le\tau\le1/2$. For a real polynomial
$p=x\widetilde p$ we write $N_*(p)=\|\widetilde p\|_{[0,a]}$ and
$N_{\mathbb R}(p)=\|\widetilde p\|_{[0,1]}$. Put
$\phi(x)=\sinh(x/2)/(x/2)$, $\phi(0)=1$. Volumes $V_{u,v}(\boldsymbol\ell)$
are Weil--Petersson orbifold volumes with $V_{0,3}=1$ and
$V_{1,1}(\ell)=(\ell^2+4\pi^2)/48$; a zero length denotes a cusp, and a
type $(u,v)$ is stable when its complexity $c=2u+v-2$ is at least one.

\begin{definition}[Orbital functional and canonical coefficients]
\label{rev:def-huang-orbital}
Let $H$ be a bounded even Borel function supported in $[-T,T]$. The
orbital functional is
\begin{equation}\label{rev:eq-orbital-functional}
 \mathcal O_H(X)=\sum_{\gamma\in\mathcal P_{\rm or}(X)}\sum_{\nu\ge1}
 \frac{\ell_\gamma H(\nu\ell_\gamma)}{2\sinh(\nu\ell_\gamma/2)},
\end{equation}
where $\mathcal P_{\rm or}(X)$ is the set of primitive oriented closed
geodesics, both orientations are counted, and all repetitions $\nu\ge1$
occur \cite[(3.1)]{HuangGap}; it is the functional $\mathcal O_{F_s}(H)$
of \eqref{eq:critical-orbital-identity} evaluated on $X$. The leading
functional is
$\Omega(H)=\int_0^\infty\sum_{\nu\ge1}\frac{2\sinh^2(\ell/2)}{\sinh(\nu\ell/2)}
H(\nu\ell)\,d\ell$ \cite[(3.20)]{HuangGap}. By
\cite[Theorems~3.5 and~4.8]{HuangGap} there are linear functionals
$\kappa_{j,n}$, $j\ge1$, on such tests, independent of $g$, of the
truncation order $R$ and of the support cutoff $T$, such that
$\mathbb E_{g,n}\mathcal O_H=\Omega(H)+\sum_{j=1}^Rg^{-j}\kappa_{j,n}(H)
+(\text{remainder})$ under Huang's genus thresholds;
\eqref{eq:critical-orbital-remainder} asserts the same expansion, with
the same $\kappa_{j,n}$, on the present envelope. By
\cite[Proposition~3.8]{HuangGap},
\begin{equation}\label{rev:eq-huang-quadratic}
 \kappa_{1,n}=\Lambda_0+n\Lambda_1+\binom n2\Lambda_2
\end{equation}
with linear functionals $\Lambda_i$ independent of $n$; $\Lambda_2$ is
the functional \eqref{eq:critical-orbital-identity}. For a polynomial
$p(x)=\sum_{t\ge1}a_tx^t$ put
\begin{equation}\label{rev:eq-huang-polynomial-test}
\begin{gathered}
 H_p=\sum_{t\ge1}a_th^{*t},\qquad
 \mathcal O_X(p)=\mathcal O_{H_p}(X),\qquad
 \omega_{1,n}(p)=\omega_1(p)=\Omega(H_p),\\
 \omega_{j,n}(p)=\kappa_{j-1,n}(H_p)\quad(j\ge2),\qquad
 W_i(p)=\Lambda_i(H_p)\quad(0\le i\le2),
\end{gathered}
\end{equation}
so that $\widehat{H_p}=p\circ f$, $\operatorname{supp}H_p\subset[-LQ,LQ]$
when $\deg p\le Q$, $\|h^{*t}\|_\infty\le\|h\|_\infty$, and
\eqref{rev:eq-huang-quadratic} reads
\begin{equation}\label{rev:eq-huang-quadratic-p}
 \omega_{2,n}=W_0+nW_1+\binom n2W_2 .
\end{equation}
These are \cite[(3.43), (3.56) and~(6.55)]{HuangGap}. The leading
coefficient $\omega_1$ does not depend on $n$
\cite[Proposition~5.1]{HuangGap}.
\end{definition}

\begin{definition}[Carriers, patterns, tails and records]
\label{rev:def-huang-carriers}
\emph{(i) Carriers.} The normalized carrier $Y$ of a nonsimple primitive
geodesic $\gamma$ on $X$ is the subsurface filled by $\gamma$: thicken
the geodesic graph, fill disks, complete peripheral punctured disks to the
inherited cusps, and tighten the remaining boundary curves to geodesics,
keeping two boundary slots when two of them tighten to the same ambient
geodesic \cite[Section~3.3]{HuangGap}. Its type $(h_0,a_0,k)$ records
the genus, the number of inherited cusps and the number $k\ge1$ of
positive boundary slots; its complexity is $m=2h_0+a_0+k-2\ge1$. By
\cite[Lemma~3.3]{HuangGap}, $\Lambda:=\ell(\partial Y)\le2\ell_\gamma$
and $\ell_\gamma\ge\pi m/2$. A carrier without stable complementary
component has $m=2g+n-2$, hence $\ell_\gamma\ge\pi(g-1+n/2)$
\cite[(3.38)]{HuangGap}.

\emph{(ii) Patterns and complementary data.} A gluing pattern $P$ on
the $k$ slots consists of a matching of $2\varsigma$ slots into pairs,
glued to one another, and a partition of the remaining slots into
$b\ge1$ nonempty blocks $I_1,\dots,I_b$, the block $I_i$ being attached
to the $i$th complementary component; $k_i=|I_i|$. The number of
patterns is at most $(4t)^{4t}$ whenever $k\le t+1$. Complementary data
are integers $(g_i,a_i)_{i\le b}$ with
\[
 \sum_ia_i=n-a_0,\qquad \sum_ig_i=g-h_0-k+\varsigma+b,\qquad
 c_i=2g_i+a_i+k_i-2\ge1 .
\]
We write $R_i=c_i-1$, so that $\sum_iR_i=N-d$ with $d=m+b-1$, and
define the \emph{deficit}
\[
 \mathfrak d=N-\max_iR_i\ \ge\ d\ \ge\ m,
\]
the loss of factorial index on passing to a largest complementary
component \cite[(2.5)--(2.6), (2.11)]{HuangGap}. The labeled weight and
the complementary tail of order $R$ are
\begin{equation}\label{rev:eq-huang-tail-def}
 W=\frac{n!}{a_0!\prod_ia_i!}
   \frac{V_{h_0,a_0+k}\prod_iV_{g_i,a_i+k_i}}{V_{g,n}},\qquad
 \mathcal T^{h_0,a_0,k}_{g,n}(R)=\sum_P\ \sum_{\mathfrak d\ge R+1}W
\end{equation}
\cite[(2.7)]{HuangGap}. For a fixed pattern the boundary version assigns
a length $y_e\ge0$ to each of the $e_P=k-\varsigma$ actual edges,
induces slot lengths $x_j$ (paired slots share one length), replaces the
numerator volumes by boundary volumes with these slot lengths and zero
cusp lengths, and multiplies by the Jacobian $J_P=\prod_ey_e$; then
$\Lambda=\sum_{j=1}^kx_j$.

\emph{(iii) Separating tail.} For a simple separating cut with $a_0$
cusps and genus $i$ on a distinguished side put $c_1=2i+a_0-1$,
$c_2=2(g-i)+n-a_0-1$ and $\mathfrak d=\min(c_1,c_2)$. Then
\begin{equation}\label{rev:eq-huang-sep-tail-def}
 \mathcal S^{\ge r}_{g,n}(x)=\sum_{\mathfrak d\ge r}\binom n{a_0}
 \frac{V_{i,a_0+1}(\mathbf 0_{a_0},x)\,V_{g-i,n-a_0+1}(\mathbf 0_{n-a_0},x)}
      {V_{g,n}}
\end{equation}
is the ordered separating kernel restricted to deficit at least $r$
\cite[(3.3), (3.9), (3.13)]{HuangGap}.

\emph{(iv) Records.} A retained record of a pattern consists of a
distinguished block $*$ and stable data $(g_i,a_i)_{i\ne*}$; the
component attached to $I_*$ is reconstructed from the ambient genus. For
a record put
\begin{equation}\label{rev:eq-huang-record}
\begin{gathered}
 \mathfrak d=m+\sum_{i\ne*}c_i,\qquad
 A_{\rm out}=a_0+\sum_{i\ne*}a_i,\qquad
 H_0=h_0+k-\varsigma-b+\sum_{i\ne*}g_i,\\
 k_*=|I_*|,\qquad
 p=n-A_{\rm out}+k_*,\qquad
 X_*=\sum_{j\in I_*}x_j .
\end{gathered}
\end{equation}
Thus $A_{\rm out}$ is the number of cusps outside the large component,
$H_0$ its genus shift, $p$ its number of marked ends, $k_*$ its number of
gluing boundaries and $X_*$ their total length. One has
$2H_0+A_{\rm out}-k_*=\mathfrak d$, $A_{\rm out}\le\mathfrak d+1$,
$p\ge k_*\ge1$, the large component has genus $g-H_0$ and factorial
index $N-\mathfrak d$, and there are at most $b(\mathfrak d+3)^{2(b-1)}$
records of deficit $\mathfrak d$ \cite[(2.40)--(2.42)]{HuangGap}. For a
truncation order $R\ge\mathfrak d$, $J_{\rm loc}=R-\mathfrak d$ denotes
the order through which the shifted volume ratio of the large component
is expanded.
\end{definition}

\begin{lemma}[Volume envelope {\cite[Lemma~2.1]{HuangGap}}]
\label{rev:lem-huang-envelope}
Put $A_0=4\pi^2$. There are absolute constants $U\ge1$,
$\beta_{\rm v}>0$ and $\kappa_{\rm v}>0$ such that, for every stable
type $(u,v)$ of complexity $c$, every $g\ge2$, $n\ge0$, and every
length vector $\boldsymbol\ell$,
\begin{equation}\label{rev:eq-huang-envelope}
 V_{u,v}\le UA_0^{c-1}\frac{(c-1)!}{\sqrt c},\qquad
 V_{g,n}\ge\beta_{\rm v}A_0^{N}\frac{N!}{\sqrt g}
   \exp\Bigl(-\frac{\kappa_{\rm v}n(n+1)}{2g}\Bigr),\qquad
 V_{u,v}(\boldsymbol\ell)\le V_{u,v}\exp\Bigl(\tfrac12\sum_j\ell_j\Bigr).
\end{equation}
\end{lemma}

The lower bound is the only volume estimate in which the cusp count
enters exponentially; the cusp-label counts $(b+1)^n$ and $2^n$ of the
deep sector are the second place where it does, and they are absorbed
in \eqref{rev:eq-critical-deep-check}. Huang imposes $g\ge\kappa_{\rm v}\mathsf s^2$,
which makes the exponential at least $e^{-1/2}$; on the envelope
$\mathsf s\le M\sqrt g$ it is at least $e^{-\kappa_{\rm v}M^2/2}$, see
\eqref{rev:eq-critical-envelope-loss}.

\begin{lemma}[Tail skeleton {\cite[proofs of Theorem~2.2 and Lemma~4.6]{HuangGap}}]
\label{rev:lem-huang-tail-skeleton}
Let $(h_0,a_0,k)$ be a carrier type of complexity $m\ge1$, let $R\ge0$,
$r=R+1$, $t=R+m+1\ge2$, $D=\max(m,r)$, $J=\lfloor N/(64\mathsf s)\rfloor$
and $L_{\mathsf s}=\log(64\mathsf s)$, and let $\Theta\ge1$ be such that
the exponential in \eqref{rev:eq-huang-envelope} is at least
$\Theta^{-1}$. Assume
\begin{equation}\label{rev:eq-huang-tail-hypotheses}
 g\ge128t\mathsf s,\qquad g\le N\le3g,\qquad N-D+1\ge g,
 \qquad J\ge t,\qquad N/(128\mathsf s)\le J\le N/2 .
\end{equation}
Fix a pattern $P$ with $b$ blocks, so that $b\le k\le m+2\le t+1$, and
put $P_b=\Theta B^{b+2}\sqrt b\,(m-1)!$ with an absolute constant
$B\ge2$. Then:
\begin{enumerate}[label=\textup{(\roman*)},leftmargin=*]
\item every labeled weight satisfies
 $W\le P_b\,\dfrac{\prod_iR_i!}{N!}\,\dfrac{n!}{a_0!\prod_ia_i!}$;
\item $A_{\rm out}\le\mathfrak d+1$, the multinomial in \textup{(i)} is
 at most $\min\{\mathsf s^{\mathfrak d+1},(b+1)^n\}$, and the total weight at
 deficit $\mathfrak d$ is at most $P_b\,b\,x_{\mathfrak d}$, where
 $x_{\mathfrak d}=\mathsf s^{\mathfrak d+1}(\mathfrak d+3)^{2b}/\binom N{\mathfrak d}$;
\item \textup{(near sector)} for $D\le\mathfrak d<J$ one has
 $x_{\mathfrak d+1}/x_{\mathfrak d}\le e^2/63<1/2$, and the sum over
 $D\le\mathfrak d\le J$ is at most
 $2P_b\,b\,(D+3)^{2b}D!\,g^{-D}\mathsf s^{D+1}$;
\item \textup{(deep sector)} for $b\ge2$ the sum over $\mathfrak d>J$ is
 at most $P_b(b+1)^n(N+1)^{2b}\binom NJ^{-1}$, and
 $\binom NJ\ge(64\mathsf s)^J\ge\exp\bigl(NL_{\mathsf s}/(128\mathsf s)\bigr)$;
 for $b=1$ the deep sector is empty;
\item $r\log g+n\log(b+1)+2b\log(N+1)\le F_t(g)
 :=4t\log g+\mathsf s\log(2t)+3t\log4$;
\item there are at most $(4t)^{4t}$ patterns, and
 $B^{b+2}\sqrt b\,(m-1)!\le(Ct)^{3t}$,
 $2b(D+3)^{2b}D!\le(2t)^{5t}$;
\item in the boundary version each of the bounds above acquires the
 factor $(1+\Lambda)^ke^{\Lambda}$.
\end{enumerate}
For the separating kernel, under \eqref{rev:eq-huang-tail-hypotheses}
with $t=r\ge1$: each zero-length summand of
\eqref{rev:eq-huang-sep-tail-def} at deficit $\mathfrak d$ is at most
$B\Theta\binom n{a_0}/(\mathfrak d\binom N{\mathfrak d})$; the sum at fixed
deficit is at most $6B\Theta\mathsf s^{\mathfrak d+1}/\binom N{\mathfrak d}$,
with successive ratio at most $1/63$ for $r\le\mathfrak d<J$, so the near
sector $r\le\mathfrak d\le J$ is at most $12B\Theta\,r!\,g^{-r}\mathsf s^{r+1}$;
the deep sector $\mathfrak d>J$ is at most
$B\Theta(g+1)2^n\binom NJ^{-1}$, and
$\log\bigl(g^r(g+1)2^n\bigr)\le F_r(g)$; boundary lengths cost the factor
$e^{x}$.
\end{lemma}

\begin{proof}
These are the displayed steps (2.12)--(2.19) and (4.19)--(4.24) in
Huang's proofs of Theorem~2.2 and Lemma~4.6, and (3.10)--(3.12) for the
separating kernel, with the following bookkeeping of hypotheses.
(i) uses \eqref{rev:eq-huang-envelope} for every factor of $W$: the
largest complementary complexity is at least $(N+1-m)/b\ge g/b$ because
$g\ge m+2$, which cancels the ambient $\sqrt g$ up to $\sqrt b$; the
powers of $A_0$ total $-b$; and the exponential of the lower bound is
absorbed into $\Theta$. (ii) is combinatorial. (iii) uses
$\mathfrak d+1\le J\le N/(64\mathsf s)$, hence
$\mathsf s(\mathfrak d+1)/(N-\mathfrak d)\le1/63$, together with
$(1+1/(\mathfrak d+3))^{2b}\le e^2$ from $b\le\mathfrak d+2$, and
$\binom ND=(N)_D/D!\ge g^D/D!$ from $N-D+1\ge g$. (iv) uses the
capped factorial inequality $\prod_iz_i!\le(N-J)!\,J!$ for nonnegative
integers with $\sum_iz_i\le N$, $z_i\le N-J$, $1\le J\le N/2$
\cite[(2.18)]{HuangGap}, applied to the $R_i$ (all of which are smaller
than $N-J$ when $\mathfrak d>J$), the count $(N+1)^{2b}$ of deep
allocations, and $N/J\ge64\mathsf s$, $J\ge N/(128\mathsf s)$. (v) uses
$N+1\le4g$, $r+2b\le4t$, $b+1\le2t$ and $2b\le3t$. (vi) is elementary
and (vii) is the boundary estimate in \eqref{rev:eq-huang-envelope}
together with $J_P\le(1+\Lambda)^k$ \cite[(2.10)]{HuangGap}. The
separating statements use \eqref{rev:eq-huang-envelope} in the form
\cite[(3.10)]{HuangGap}, the bound $\binom n{a_0}\le\mathsf s^{\mathfrak d+1}$
for the smaller side, at most $2(\mathfrak d+2)$ ordered records at fixed
deficit, binomial unimodality for $\mathfrak d>J$, and
$\sum_{a_0}\binom n{a_0}=2^n$ with at most $g+1$ genus choices.

In \cite{HuangGap} the single hypothesis $g\ge K(\mathsf s^2+t^2\mathsf s)$,
$K\ge\max\{4096,\kappa_{\rm v}\}$, serves two purposes: it yields
\eqref{rev:eq-huang-tail-hypotheses} with $\Theta=e^{1/2}$
\cite[(4.20)]{HuangGap}, and it closes the deep sector through the
absorption inequality $F_t(g)<gL_{\mathsf s}/(128\mathsf s)$
\cite[(4.25)]{HuangGap}. Only these two uses depend on a quadratic genus
condition; they are replaced in Subsection~\ref{rev:subsec-critical-tails}.
\end{proof}

\begin{lemma}[Retained records and orbital integration {\cite[Lemmas~4.1--4.2, Propositions~3.4 and~4.7, Theorem~4.8]{HuangGap}}]
\label{rev:lem-huang-records}
\textup{(a)} \emph{Graded shifted-volume error.} There is an absolute
$C\ge2$ such that, for integers $J_{\rm loc},H_0\ge0$ and
$\mathfrak d\ge1$ with $p=n+2H_0-\mathfrak d\ge1$, for
$L_{\rm sh}=J_{\rm loc}+H_0+\mathfrak d+1$, $0\le\ell\le p$,
$\mathbf x\in[0,\infty)^\ell$, $X=\sum_ix_i$ and
$g\ge[CL_{\rm sh}]^C\mathsf s$,
\begin{equation}\label{rev:eq-huang-shifted}
 \Bigl|\frac{V_{g-H_0,p}(\mathbf 0_{p-\ell},\mathbf x)}{V_{g,n}}
 -(8\pi^2)^{-\mathfrak d}g^{-\mathfrak d}
  \sum_{j=0}^{J_{\rm loc}}g^{-j}Q_j(\mathbf x)\Bigr|
 \le g^{-\mathfrak d-J_{\rm loc}-1}\mathsf s^{J_{\rm loc}+1}
     [CL_{\rm sh}]^{C(J_{\rm loc}+1)}(1+X)^{2J_{\rm loc}+2}e^{X/2},
\end{equation}
where the $Q_j=Q_j^{n,H_0,\mathfrak d;\ell}$ are canonical functions
independent of $g$ and $J_{\rm loc}$, with
$|Q_j(\mathbf x)|\le\mathsf s^j[CL_{\rm sh}]^{Cj}(1+X)^{2j}e^{X/2}$. The
boundary expansion of $V_{G,p}(\mathbf 0_{p-\ell},\mathbf x)/V_{G,p}$
through order $J$ holds for $G\ge[C(J+1)]^C(p+1)$
\cite[Lemma~4.1]{HuangGap}, and the handle-ratio and source thresholds
of \cite[Proposition~2.3]{HuangGap} used inside these lemmas are bounded
by $[Ct]^C\mathsf s$ whenever $p\le n+3t$ and the order is at most $t$
\cite[proof of Proposition~4.7]{HuangGap}. In particular every retained-record
genus threshold used here from \cite[Section~4.1]{HuangGap} is of the form
$[Ct]^C\mathsf s$.

\textup{(b)} \emph{Small-volume prefactor.} For a record $e$ of deficit
$\mathfrak d$, without any genus condition,
\begin{equation}\label{rev:eq-huang-record-prefactor}
\begin{gathered}
 A_e(\mathbf y)=J_P\,V_{h_0,a_0+k}(\mathbf 0_{a_0},\mathbf x)
 \frac{n!}{a_0!\,(n-A_{\rm out})!\prod_{i\ne*}a_i!}
 \prod_{i\ne*}V_{g_i,a_i+k_i}(\mathbf 0_{a_i},\mathbf x_{I_i}),\\
 |A_e(\mathbf y)|\le\mathsf s^{A_{\rm out}}
   [C(\mathfrak d+1)]^{C(\mathfrak d+1)}(1+\Lambda)^ke^{\Lambda-X_*/2}
\end{gathered}
\end{equation}
\cite[(2.43), (4.13)]{HuangGap}.

\textup{(c)} \emph{Kernel errors.} Let $t=R+m+1$. For a record of
deficit $\mathfrak d\le R$ take $J_{\rm loc}=R-\mathfrak d$; then
$H_0\le t$, $L_{\rm sh}=R+H_0+1\le2t$ and $p\ge k_*\ge1$. The product
of \eqref{rev:eq-huang-shifted} and \eqref{rev:eq-huang-record-prefactor}
has cusp power $\mathsf s^{A_{\rm out}+J_{\rm loc}+1}\le\mathsf s^{R+2}$,
exponential factor $e^{\Lambda}$ and length power at most
$k+2(R-\mathfrak d)+2\le2t$. Reconstruction of the large component
requires $N>2R$ and $g>H_0$; unavailable cusp allocations are omitted
before any volume is evaluated. The sum over records of the retained
terms is Huang's canonical coefficient
$\mathcal F^{h_0,a_0,k;P}_{j,n}$ \cite[(2.45)]{HuangGap}, and the full
nonnegative kernel
\begin{equation}\label{rev:eq-huang-kernel}
 F^{h_0,a_0,k;P}_{g,n}(\mathbf y)
 =J_P\,V_{h_0,a_0+k}(\mathbf 0_{a_0},\mathbf x)
  \sum\frac{n!}{a_0!\prod_ia_i!}
  \frac{\prod_iV_{g_i,a_i+k_i}(\mathbf 0_{a_i},\mathbf x_{I_i})}{V_{g,n}}
\end{equation}
differs from $\sum_{j=1}^Rg^{-j}\mathcal F^{h_0,a_0,k;P}_{j,n}$ by at
most the sum of the record errors plus the boundary version of
$\mathcal T^{h_0,a_0,k}_{g,n}(R)$ for the pattern. The nonseparating
kernel $U_{g,n}(x)=V_{g-1,n+2}(\mathbf 0_n,x,x)/V_{g,n}$ is the product
of a scalar handle ratio and a boundary ratio and has no tail; its
expansion through order $R$ has error
$g^{-R-1}\mathsf s^{R+1}[Cr]^{Cr}(1+x)^{2R+2}e^x$ under
$g\ge[Cr]^C\mathsf s$. The separating kernel splits into retained
endpoints $\mathfrak d\le R$, treated by \eqref{rev:eq-huang-shifted}
with $H_0=h_0\le(\mathfrak d+1)/2$ and the prefactor
$\binom n{a_0}V_{h_0,a_0+1}(\mathbf 0_{a_0},x)
\le\mathsf s^{a_0}[C(\mathfrak d+1)]^{C(\mathfrak d+1)}e^{x/2}$, and the
omitted tail $\mathcal S^{\ge R+1}_{g,n}(x)$.

\textup{(d)} \emph{Integration.} With $K_{g,n}=U_{g,n}+S_{g,n}$, where
$S_{g,n}$ is the full ordered separating kernel,
\begin{equation}\label{rev:eq-huang-simple-integration}
 \mathbb E_{g,n}\mathcal O^{\rm simp}_H
 =\int_0^\infty\sum_{\nu\ge1}\frac{\ell^2H(\nu\ell)}{2\sinh(\nu\ell/2)}
   K_{g,n}(\ell)\,d\ell ,
\end{equation}
and a kernel bound $|k(\ell)|\le A(1+\ell)^de^{\ell}$ gives
$\bigl|\int_0^\infty\sum_{\nu\ge1}\frac{\ell^2k(\ell)H(\nu\ell)}
{2\sinh(\nu\ell/2)}\,d\ell\bigr|\le\frac{\zeta(3)}2AT^2(1+T)^de^T\|H\|_\infty$
\cite[Lemma~3.2 and~(3.19)]{HuangGap}. For a carrier type with at least
one stable complementary component,
\begin{equation}\label{rev:eq-huang-carrier-integration}
 \mathbb E_{g,n}\mathcal O^{h_0,a_0,k}_H
 =\frac1{k!}\sum_P\int_{(0,\infty)^{e_P}}
   \Psi^{h_0,a_0,k}_H(\mathbf x_P(\mathbf y))\,
   F^{h_0,a_0,k;P}_{g,n}(\mathbf y)\,d\mathbf y ,
\end{equation}
where $\Psi_H$ is the fiber integral, over the carrier moduli space with
the prescribed boundary lengths, of the sum
\eqref{rev:eq-orbital-functional} restricted to filling geodesics,
divided by the carrier volume; it satisfies
$|\Psi_H|\le Cm(1+T)e^T\|H\|_\infty$ and vanishes unless $\Lambda\le2T$
and $m\le2T/\pi$, and $\operatorname{Leb}\{\Lambda_P\le2T\}\le(1+2T)^k$
\cite[(3.23), (3.25), (3.29), (3.33)]{HuangGap}. At fixed $m$ there are
at most $(m+2)^2$ carrier types and $k\le m+2$. The vanishing rules
\cite[(3.32)]{HuangGap} identify the sum of the retained coefficients
over all types with $\kappa_{j,n}$ for every $R$ and $T$. Apart from the
exclusion of carriers without stable complementary component, which
requires $\pi(g-1+n/2)>T$, no statement in \textup{(d)} involves a genus
condition.
\end{lemma}

\begin{lemma}[Truncated positive trace {\cite[Theorem~5.3, Propositions~5.1 and~5.4, (5.44)--(5.45)]{HuangGap}}]
\label{rev:lem-huang-trace}
Let $X$ be a complete hyperbolic surface of finite area with $n\ge1$
cusps and no orbifold points, let $A_X=f(\sqrt{\Delta_X-1/4})$ be the
bounded self-adjoint operator with Laplace multiplier $f(r(\lambda))$,
$r(\lambda)=\sqrt{\lambda-1/4}$ for $\lambda\ge1/4$ and
$r(\lambda)=i\sqrt{1/4-\lambda}$ for $0\le\lambda<1/4$, so that
$0\le A_X\le a$, the constant eigenfunction has $A_X$-eigenvalue $a$, and
the spectrum of $\Delta_X$ in $[1/4,\infty)$ is mapped into $[0,1]$.
For $Y>1$ let $C_b(Y)$ be the width-one cusp regions above height $Y$,
which are embedded and pairwise disjoint, let
$X(Y)=X\setminus\bigcup_bC_b(Y)$, $\chi_Y=\mathbf 1_{X(Y)}$, and fix
\begin{equation}\label{rev:eq-huang-YQ}
 Y_Q=3e^{LQ}\qquad(Q\ge1)
\end{equation}
\cite[(5.23), (6.1)]{HuangGap}. Let $p$ be a real polynomial with
$p(0)=0$ and $\deg p\le Q$.
\begin{enumerate}[label=\textup{(\roman*)},leftmargin=*]
\item The compression $\chi_{Y_Q}p(A_X)\chi_{Y_Q}$ is trace class and
\begin{equation}\label{rev:eq-huang-trace-identity}
 \tau_{X,Q}(p):=\operatorname{Tr}\bigl(\chi_{Y_Q}p(A_X)\chi_{Y_Q}\bigr)
 =\Bigl(\mathcal A_{g,n}-\frac n{Y_Q}\Bigr)\mathcal I(p)
  +n\,\mathcal P_{Y_Q}(p)+\mathcal O_X(p),
\end{equation}
where $\mathcal I(p)=\frac1{4\pi}\int_{\mathbb R}r\tanh(\pi r)p(f(r))\,dr$
and $\mathcal P_Y(p)=2\int_0^{LQ}K_{H_p}(\varrho)\cosh(\varrho/2)
\mathsf H_{\lfloor2Y\sinh(\varrho/2)\rfloor}\,d\varrho$, with $K_H$ the Abel
transform $K_H(\varrho)=-\frac1{\sqrt2\pi}\int_\varrho^\infty
\frac{H'(u)\,du}{\sqrt{\cosh u-\cosh\varrho}}$ and
$\mathsf H_m=\sum_{k=1}^mk^{-1}$.
\item If $p\ge0$ on the spectrum of $A_X$ then $\tau_{X,Q}(p)\ge0$;
 in general $|\tau_{X,Q}(p)|\le N_*(p)\tau_{X,Q}(x)$. More generally,
 if $\psi=x\xi$ with $\xi$ bounded and Borel on $[0,a]$, then
 $\chi_{Y_Q}\psi(A_X)\chi_{Y_Q}$ is trace class and its trace, also
 after insertion of $1-P_0$ with $P_0$ as in \textup{(v)}, has absolute
 value at most $\|\xi\|_{[0,a]}\tau_{X,Q}(x)$ \cite[(5.29)]{HuangGap}.
\item There are constants $I_f,J_f,b_f$ depending only on $f$ with
\begin{equation}\label{rev:eq-huang-deterministic}
 |\mathcal I(p)|\le I_fN_{\mathbb R}(p),\qquad
 |\mathcal P_{Y_Q}(p)|\le J_fQ^2N_{\mathbb R}(p),\qquad
 |\mathcal P_{Y_Q}(x)|\le b_fQ .
\end{equation}
\item There is a finite signed measure $\mu_0$ on $[0,1]$, independent
 of $n$, with $\|\mu_0\|_{\mathrm{TV}}<5$ and
 $\omega_1(p)=p(a)+\int_{[0,1]}p\,d\mu_0$; in particular
 $|\omega_1(x)|\le6a$.
\item Let $P_0$ be the orthogonal projection onto constants. The
 compressed trace
 $\tau^\circ_{X,Q}(p):=\operatorname{Tr}\bigl(\chi_{Y_Q}p(A_X)(1-P_0)\chi_{Y_Q}\bigr)
 =\tau_{X,Q}(p)-\bigl(1-\frac n{\mathcal A_{g,n}Y_Q}\bigr)p(a)$ is
 nonnegative when $p\ge0$ on $[0,a]$, and
\begin{equation}\label{rev:eq-huang-trace-compressed}
 \tau^\circ_{X,Q}(p)
 =\Bigl(\mathcal A_{g,n}-\frac n{Y_Q}\Bigr)\mathcal I(p)
  +n\,\mathcal P_{Y_Q}(p)+\mu_0(p)
  +\bigl(\mathcal O_X(p)-\omega_1(p)\bigr)
  +\frac{n\,p(a)}{\mathcal A_{g,n}Y_Q}.
\end{equation}
The last term is the contribution of the mass of the constant
eigenfunction removed by the truncation.
\end{enumerate}
\end{lemma}

\begin{lemma}[Retained mass {\cite[Lemma~6.5]{HuangGap}}]
\label{rev:lem-huang-retained-mass}
If $\psi$ is a normalized $L^2$ eigenfunction of $\Delta_X$ with
eigenvalue $0<\lambda\le1/4-\varepsilon$, including residual
eigenfunctions, then $\|\chi_{Y_Q}\psi\|_2^2\ge1-e^{-2\sqrt\varepsilon LQ}$,
uniformly in the surface and in its number of cusps.
\end{lemma}

\begin{lemma}[Coefficient bounds {\cite[Lemmas~6.1 and~6.7]{HuangGap}}]
\label{rev:lem-huang-coefficients}
\textup{(a)} If $p(x)=\sum_{t=1}^Qa_tx^t$, then
$\sum_{t=1}^Q|a_t|\le7^QN_*(p)$.

\textup{(b)} Let $c_f\ge2$ be the integer of
\cite[Proposition~3.6]{HuangGap} and $D_0=2c_f+4$. There is $C_0\ge2$,
depending only on $f$, such that for all integers $j\ge2$, $Q,n\ge1$ and
all real $p=x\widetilde p$ of degree at most $Q$,
\begin{equation}\label{rev:eq-huang-coefficient-bound}
 |\omega_{j,n}(p)|\le\mathsf s^{j}[C_0(Q+j)]^{D_0j}N_*(p),
\end{equation}
with no condition on the genus. With $A_Q=(3C_0Q)^{D_0}$ this gives
$|\omega_{j,n}(p)|\le\mathsf s^jA_Q^jN_*(p)$ for $2\le j\le2Q$.
\end{lemma}

\begin{lemma}[First moment of the truncated trace {\cite[proof of Lemma~6.8]{HuangGap}}]
\label{rev:lem-huang-first-moment}
Suppose $1\le n\le g$ and $|\mathbb E_{g,n}\mathcal O_X(x)-\omega_1(x)|\le g$.
Then, for every $Q\ge1$,
\begin{equation}\label{rev:eq-huang-first-moment}
 \mathbb E_{g,n}\tau_{X,Q}(x)\le D_f(g+nQ),\qquad D_f=1+6a+6\pi I_f+b_f,
\end{equation}
and every polynomial trace and orbital statistic of degree at most $Q$
with zero constant term is absolutely integrable, with
$\mathbb E_{g,n}|\mathcal O_X(p)|\le N_*(p)\{D_f(g+nQ)+(4\pi g+2\pi n)I_f+nJ_fQ^2\}$.
\end{lemma}

\begin{proof}
Take expectations in \eqref{rev:eq-huang-trace-identity} at $p=x$,
where $H_x=h$ and $\mathcal O_X(x)\ge0$. By
\eqref{rev:eq-huang-deterministic}, Lemma~\ref{rev:lem-huang-trace}(iv),
the hypothesis and $\mathcal A_{g,n}-n/Y_Q\le4\pi g+2\pi n\le6\pi g$,
\[
 \mathbb E\tau_{X,Q}(x)\le(4\pi g+2\pi n)I_f+nb_fQ+6a+g\le D_f(g+nQ).
\]
Lemma~\ref{rev:lem-huang-trace}(ii) and the exact identity give the
remaining statements \cite[(6.5)--(6.6)]{HuangGap}.
\end{proof}

\begin{lemma}[Polynomial transfer {\cite[proof of Proposition~6.9]{HuangGap}}]
\label{rev:lem-huang-transfer}
Let $Q\ge1$, $q=Q+2$, and let $p=x\widetilde p$ be real of degree at most
$Q$. Suppose that for every $1\le t\le Q$
\begin{equation}\label{rev:eq-huang-monomial-hypothesis}
 \Bigl|g^{-1}\mathbb E_{g,n}\mathcal O_X(x^t)
  -\sum_{j=1}^{q-1}g^{-j}\omega_{j,n}(x^t)\Bigr|\le E ,
\end{equation}
and that $\mathsf sA_Q/g\le1/2$. Then
\begin{equation}\label{rev:eq-huang-transfer}
 \bigl|\mathbb E_{g,n}\mathcal O_X(p)-\omega_1(p)-g^{-1}\omega_{2,n}(p)\bigr|
 \le\Bigl\{\frac{2\mathsf s^3A_Q^3}{g^2}+7^QgE\Bigr\}N_*(p).
\end{equation}
\end{lemma}

\begin{proof}
Linearity, Lemma~\ref{rev:lem-huang-coefficients}(a) and
\eqref{rev:eq-huang-monomial-hypothesis} give
$\mathbb E\mathcal O_X(p)=\omega_1(p)+\sum_{j=2}^{Q+1}g^{1-j}\omega_{j,n}(p)
+\mathcal R$ with $|\mathcal R|\le7^QgE\,N_*(p)$ \cite[(6.51)]{HuangGap}.
For $3\le j\le Q+1$, Lemma~\ref{rev:lem-huang-coefficients}(b) gives
\[
 \sum_{j=3}^{Q+1}g^{1-j}|\omega_{j,n}(p)|
 \le\frac{\mathsf s^3A_Q^3}{g^2}N_*(p)\sum_{j\ge3}(\mathsf sA_Q/g)^{j-3}
 \le\frac{2\mathsf s^3A_Q^3}{g^2}N_*(p).
\]
In \cite{HuangGap} the hypothesis is supplied by Lemma~6.8 with
$E=g^{-q}\mathsf s^q[c_3q]^{c_3q}$ under the additive threshold (6.46),
and $7^QgE\le(\mathsf s^2/g)2^{-Q}$ there by (6.52); here it is supplied by
\eqref{eq:critical-monomial-error}.
\end{proof}

\begin{lemma}[Weighted distributions and the exceptional measure {\cite[Lemma~6.10, Theorem~6.11]{HuangGap}}]
\label{rev:lem-huang-exceptional}
Put $m=2D_0+2$ (here $m$ is an order of differentiation, not the carrier complexity of Subsection~\ref{rev:subsec-huang-inputs}). For each $n\ge1$ the functional
$\mathscr K_n(q)=\omega_{2,n}(xq)$ on polynomials $q$ extends uniquely
to a distribution on $\mathbb R$ supported in $[0,a]$ and of order at
most $m$, with $|\mathscr K_n(\varphi)|\le K_f(n+1)^2\|\varphi\|_{C^m([0,a])}$
for a constant $K_f$ depending only on $f$. The unweighted extension
$\mathscr V_n(\psi)=\mathscr K_n\bigl((\psi-\psi(0))/x\bigr)$ restricts on
$(1,\infty)$ to a positive Radon measure $\mu_n$ supported in $(1,a]$,
with $\mu_n([b,a])\le C_{f,b}(n+1)^2$ for $1<b\le a$
\cite[(6.57)--(6.58)]{HuangGap}. The second difference in
$n$ of \eqref{rev:eq-huang-quadratic-p} gives the weighted extension
\begin{equation}\label{rev:eq-huang-T2}
 \mathscr T_2=\mathscr K_3-2\mathscr K_2+\mathscr K_1,\qquad
 \mathscr T_2(q)=W_2(xq),\qquad
 |\mathscr T_2(\varphi)|\le38K_f\|\varphi\|_{C^m([0,a])}.
\end{equation}
There is a positive Radon measure $\eta_2$ on $(1,a]$ such that
\begin{equation}\label{rev:eq-huang-eta2}
 \mathscr T_2(\psi/x)=\int\psi\,d\eta_2,\qquad
 \eta_2([b,a])\le C_{f,b}\quad(1<b\le a),
\end{equation}
for all $\psi\in C_c^\infty((1,\infty))$; these tests may be nonzero at
$a$, and the statement allows $\eta_2$ to have an atom at $a$.
\end{lemma}

\begin{lemma}[Detector and count bound {\cite[Lemma~6.4, (6.26), Proposition~6.12]{HuangGap}}]
\label{rev:lem-huang-detector}
Fix $0<\varepsilon<1/4$ and put $b_\varepsilon=f(i\sqrt\varepsilon)\in(1,a)$.
There are $\kappa_\varepsilon>0$, $B_\varepsilon<\infty$, an integer
$Q_\varepsilon$ and, for every $Q\ge1$, a real polynomial $p_Q$ of degree
$Q$ with $p_Q(0)=0$, $p_Q(a)=1$ and $0\le p_Q\le1$ on $[0,a]$, such that
for $Q\ge Q_\varepsilon$
\begin{equation}\label{rev:eq-huang-detector}
 N_*(p_Q)\le B_\varepsilon,\qquad
 N_{\mathbb R}(p_Q)\le B_\varepsilon e^{-\kappa_\varepsilon Q},\qquad
 \inf_{[b_\varepsilon,a]}p_Q\ge1-e^{-8\kappa_\varepsilon Q}.
\end{equation}
Let $N_\varepsilon(X)$ be the number of eigenvalues of $\Delta_X$ in
$(0,1/4-\varepsilon]$, counted with multiplicity and including residual
eigenvalues. Moreover $e^{-8\kappa_\varepsilon Q}\le\tfrac12$ and
$e^{-2\sqrt\varepsilon LQ}\le\tfrac12$ for $Q\ge Q_\varepsilon$
\cite[(6.14)]{HuangGap}. Then, by \eqref{rev:eq-huang-detector} and
Lemma~\ref{rev:lem-huang-retained-mass},
\begin{equation}\label{rev:eq-huang-count-from-trace}
 N_\varepsilon(X)\le4\,\tau^\circ_{X,Q}(p_Q)\qquad(Q\ge Q_\varepsilon).
\end{equation}
Finally, there is $C_{f,\varepsilon}$ such that, for all $n\ge1$ and
$Q\ge Q_\varepsilon$, without any genus condition,
\begin{equation}\label{rev:eq-huang-degree-independent}
 |\omega_{2,n}(p_Q)|\le C_{f,\varepsilon}(n+1)^2 .
\end{equation}
\end{lemma}

The estimate \eqref{rev:eq-huang-degree-independent} is deduced in
\cite{HuangGap} from Lemma~\ref{rev:lem-huang-exceptional} by splitting
$p_Q$ at a cutoff strictly below the detector transition; it uses only the
local mass bound for the measures $\mu_n$ associated with $\mathscr K_n$
and the exponential smallness of $p_Q/x$ below the transition, and hence
holds for every $n$.

\medskip\noindent
\emph{Verbatim inputs and replaced estimates.}
Lemmas~\ref{rev:lem-huang-envelope}--\ref{rev:lem-huang-detector} are
used below as stated. Lemmas~\ref{rev:lem-huang-first-moment}
and~\ref{rev:lem-huang-transfer} record the deterministic parts of the
proofs of \cite[Lemma~6.8 and Proposition~6.9]{HuangGap} under the
hypotheses displayed there, which replace his additive genus
thresholds. Two estimates of
\cite{HuangGap} are not available on the envelope $\mathsf s\le M\sqrt g$
and are replaced:
\begin{enumerate}[label=\textup{(R\arabic*)},leftmargin=*]
\item the use of the volume lower bound in \eqref{rev:eq-huang-envelope}
 with exponential factor at least $e^{-1/2}$, which Huang obtains from
 $g\ge\kappa_{\rm v}\mathsf s^2$; it is replaced by the uniform loss
 $e^{\kappa_{\rm v}M^2/2}$ in \eqref{rev:eq-critical-envelope-loss};
\item the deep-sector absorption \cite[(4.25)]{HuangGap}, which closes
 Lemma~\ref{rev:lem-huang-tail-skeleton}(iv)--(v) under
 $g\ge K\mathsf s(\mathsf s+t^2)$; it is replaced by
 \eqref{rev:eq-critical-deep-chain}--\eqref{rev:eq-critical-deep-check},
 which use $t\le g^\theta$ and $\mathsf s\le M\sqrt g$ instead.
\end{enumerate}
Besides these, the linear thresholds $[Ct]^C\mathsf s$ of
Lemma~\ref{rev:lem-huang-records}(a) and the carrier exclusion
$\pi(g-1+n/2)>T$ are verified on the envelope in
Subsection~\ref{rev:subsec-critical-records}, and Huang's additive
threshold for polynomial transfer is replaced by inserting
\eqref{eq:critical-monomial-error} into Lemma~\ref{rev:lem-huang-transfer}
in Subsection~\ref{rev:subsec-critical-polynomial}. Every inequality that
is not displayed in Subsections~\ref{rev:subsec-critical-tails}--\ref{rev:subsec-critical-polynomial}
is taken verbatim from the lemmas above.

\subsection{Near and deep complementary-volume tails}\label{rev:subsec-critical-tails}

This subsection proves the graded tail bound
\eqref{eq:critical-graded-tail} for carrier types and the separating
tail bound \eqref{eq:critical-separating-tail} on the envelope
$\mathsf s\le M\sqrt g$, $t\le g^\theta$; these replace the two
estimates \textup{(R1)}--\textup{(R2)} listed at the end of
Subsection~\ref{rev:subsec-huang-inputs}. The input is the tail
skeleton of Lemma~\ref{rev:lem-huang-tail-skeleton}. The tails are
split into a near sector, summed as a geometric series, and a deep
sector, bounded by the capped factorial inequality; the only new
step is the comparison \eqref{rev:eq-critical-deep-check}, which closes
the deep sector using the envelope in place of a quadratic genus
condition. Use the
notation of Lemma~\ref{rev:lem-huang-tail-skeleton}: put
$N=2g+n-3$, $r=R+1$, $t=R+m+1$, $D=\max(m,r)$, and
$J=\lfloor N/(64\mathsf s)\rfloor$. Choose initially
\[
 0<\theta<\min\{1/4,1/(512M^2)\}.
\]
For $1\le t\le g^\theta$, uniformly at sufficiently large $g$,
\begin{equation}\label{eq:critical-tail-domain}
 g\ge128t\mathsf s,\quad g\le N\le3g,\quad N-D+1\ge g,
 \quad J\ge t,\quad N/(128\mathsf s)\le J\le N/2 .
\end{equation}
Indeed $g\ge128t\mathsf s$ follows from $g^{1/2-\theta}\ge128M$;
$g\le N\le3g$ from $g\ge3$ and $n\le g+3$; $N-D+1\ge g$ from
$D\le t\le g^\theta\le g-2$; and $N\ge g\ge128t\mathsf s$ gives
$N/(64\mathsf s)\ge2t$, hence $J\ge N/(128\mathsf s)\ge t$. Thus
\eqref{rev:eq-huang-tail-hypotheses} holds. The inverse volume lower
bound in \eqref{rev:eq-huang-envelope} loses at most a bounded factor on
the envelope: since $n(n+1)\le\mathsf s^2\le M^2g$,
\begin{equation}\label{rev:eq-critical-envelope-loss}
 V_{g,n}\ge\beta_{\rm v}A_0^N\frac{N!}{\sqrt g}
   \exp\Bigl(-\frac{\kappa_{\rm v}n(n+1)}{2g}\Bigr)
 \ge e^{-\kappa_{\rm v}M^2/2}\,\beta_{\rm v}A_0^N\frac{N!}{\sqrt g},
\end{equation}
so Lemma~\ref{rev:lem-huang-tail-skeleton} applies with
$\Theta=e^{\kappa_{\rm v}M^2/2}$. Its factorial envelope therefore
remains valid with prefactor
\[
 P_{b,M}=e^{\kappa_{\rm v}M^2/2}B^{b+2}\sqrt b\,(m-1)!.
\]
Here $b$ is the number of complementary components. The largest
complementary complexity is at least $g/b$, by the same Euler identity;
the terminal types $(0,3)$ and $(1,1)$ are included in the envelope.

\emph{Near and deep sectors.}
The exact cusp budget is $A_{\rm out}\le\mathfrak d+1$, and the number
of fixed-deficit allocations is at most
$b(\mathfrak d+3)^{2(b-1)}$ (Lemma~\ref{rev:lem-huang-tail-skeleton}(ii)).
Neither statement changes with the genus
threshold. In the near sector, the successive ratio up to $J$ is at
most $e^2/63<1/2$ by \eqref{eq:critical-tail-domain}
(Lemma~\ref{rev:lem-huang-tail-skeleton}(iii)). Summing gives
\begin{equation}\label{eq:critical-near-tail}
 2P_{b,M}b(D+3)^{2b}D!g^{-D}\mathsf s^{D+1}
 \le2P_{b,M}b(D+3)^{2b}D!g^{-r}\mathsf s^{r+1},
\end{equation}
where the inequality uses $D\ge r$ and $\mathsf s/g\le1$.
For $b=1$ this includes every allocation. For $b\ge2$, the capped
factorial inequality and multinomial sum of
Lemma~\ref{rev:lem-huang-tail-skeleton}(iv), before Huang's final
absorption step, give the deep-sector bound
\begin{equation}\label{eq:critical-deep-tail}
 P_{b,M}(b+1)^n(N+1)^{2b}
               \exp\{-N\log(64\mathsf s)/(128\mathsf s)\}.
\end{equation}
\emph{Replacing the absorption step.}
Write $L_{\mathsf s}=\log(64\mathsf s)$
and $F_t(g)=4t\log g+\mathsf s\log(2t)+3t\log4$. By
Lemma~\ref{rev:lem-huang-tail-skeleton}(iv)--(v) and $N\ge g$,
\begin{equation}\label{rev:eq-critical-deep-chain}
 g^{r}\,P_{b,M}(b+1)^n(N+1)^{2b}\binom NJ^{-1}
 \le P_{b,M}\exp\Bigl\{F_t(g)-\frac{NL_{\mathsf s}}{128\mathsf s}\Bigr\}
 \le P_{b,M}\exp\Bigl\{F_t(g)-\frac{gL_{\mathsf s}}{128\mathsf s}\Bigr\},
\end{equation}
so the deep sector is at most $P_{b,M}g^{-r}$ as soon as
$F_t(g)<gL_{\mathsf s}/(128\mathsf s)$. Uniformly on the envelope
$1\le t\le g^\theta$, $1\le\mathsf s\le M\sqrt g$,
\begin{equation}\label{rev:eq-critical-deep-check}
\begin{aligned}
 \frac{\mathsf sF_t(g)}{gL_{\mathsf s}}
 &=\frac{\mathsf s(4t\log g+3t\log4)}{gL_{\mathsf s}}
   +\frac{\mathsf s^2\log(2t)}{gL_{\mathsf s}},\\
 \frac{\mathsf s(4t\log g+3t\log4)}{gL_{\mathsf s}}
 &\le\frac{M\sqrt g\,g^{\theta}(4\log g+3\log4)}{g\log64}
   =O\bigl(g^{\theta-1/2}\log g\bigr)=o(1),\\
 \frac{\mathsf s^2\log(2t)}{gL_{\mathsf s}}
 &\le\frac{M^2(\theta\log g+\log2)}{\log(64M\sqrt g)}
   =2M^2\theta+o(1)<\frac1{256}+o(1).
\end{aligned}
\end{equation}
For the last inequality, use the monotonicity of $x^2/\log(64x)$ for
$x\ge1$ to replace $\mathsf s$ by $M\sqrt g$, and $\log(2t)\le\theta\log g+\log2$.
Thus eventually $\mathsf sF_t(g)/(gL_{\mathsf s})<1/128$, that is,
$F_t(g)<gL_{\mathsf s}/(128\mathsf s)$, and
\eqref{eq:critical-deep-tail} is at most
$P_{b,M}g^{-r}$. Summing the at most $(4t)^{4t}$ patterns with
\eqref{eq:critical-near-tail}, and using
Lemma~\ref{rev:lem-huang-tail-skeleton}(vi), gives the graded tail
\begin{equation}\label{eq:critical-graded-tail}
 \mathcal T_{g,n}^{h_0,a_0,k}(R)
       \le g^{-R-1}\mathsf s^{R+2}(C_Mt)^{C_Mt}.
\end{equation}
The boundary version has the same extra factor
$(1+\Lambda)^k e^\Lambda$ as in Huang's lemma, since that factor
comes from the numerator boundary estimate
(Lemma~\ref{rev:lem-huang-tail-skeleton}(vii)).

\emph{The separating kernel.}
For the ordered separating kernel, repeat the split with $t=r\ge1$,
using the separating part of Lemma~\ref{rev:lem-huang-tail-skeleton}
with $\Theta=e^{\kappa_{\rm v}M^2/2}$ and $B_M=Be^{\kappa_{\rm v}M^2/2}$.
The near sector is at most $12B_Mr!g^{-r}\mathsf s^{r+1}$, and the deep
sector before absorption is
$B_M(g+1)2^n\exp\{-N\log(64\mathsf s)/(128\mathsf s)\}$.
The logarithm of $g^r(g+1)2^n$ is bounded by $F_r(g)$.
The same comparison \eqref{rev:eq-critical-deep-check}, including $r=1$,
gives
\begin{equation}\label{eq:critical-separating-tail}
 \mathcal S_{g,n}^{\ge r}(x)
       \le g^{-r}\mathsf s^{r+1}(C_Mr)^{C_Mr}e^x.
\end{equation}
The original allocation sums include the central tie and empty sums,
so no endpoint configuration is discarded.

\subsection{Retained records and orbital integration}\label{rev:subsec-critical-records}

This subsection completes the proof of the orbital estimate
\eqref{eq:critical-orbital-remainder}. We verify that the remaining
inputs of Huang's orbital expansion, namely the linear genus thresholds
of Lemma~\ref{rev:lem-huang-records}(a) and the whole-carrier exclusion
of Definition~\ref{rev:def-huang-carriers}(i), hold on the new
envelope; we then combine the record prefactors of
Lemma~\ref{rev:lem-huang-records}(b)--(c) with the tail bound
\eqref{eq:critical-graded-tail} into a graded carrier error, and
integrate it over orbital types using
Lemma~\ref{rev:lem-huang-records}(d).

\emph{Linear thresholds.}
The graded boundary and shifted-volume estimates of
Lemma~\ref{rev:lem-huang-records}(a) require only $g\ge[Ct]^C\mathsf s$.
Decrease $\theta$ so that $C\theta<1/2$ for each of their finitely many
fixed constants; then
\[
 [Ct]^C\mathsf s\le C^Cg^{C\theta}M\sqrt g=o(g)
 \qquad(t\le g^\theta,\ \mathsf s\le M\sqrt g),
\]
so these requirements, including fixed multiples of
$t$, hold uniformly on $t\le g^\theta$, $\mathsf s\le M\sqrt g$.
\emph{Retained records.}
For a retained carrier record of deficit $\mathfrak d\le R$, its exact
identities \eqref{rev:eq-huang-record} and
Lemma~\ref{rev:lem-huang-records}(c) give
\[
 H_0\le t,\qquad p=n-A_{\rm out}+k_*\ge k_*\ge1,
 \qquad J_{\rm loc}=R-\mathfrak d,
 \qquad A_{\rm out}+J_{\rm loc}+1\le R+2 .
\]
The small-volume prefactor is bounded by
$\mathsf s^{A_{\rm out}}[C(\mathfrak d+1)]^{C(\mathfrak d+1)}
 (1+\Lambda)^k e^{\Lambda-X_*/2}$ by
\eqref{rev:eq-huang-record-prefactor}, without a genus condition.
Multiplication by the local shifted-volume error
\eqref{rev:eq-huang-shifted} therefore has cusp
power at most $\mathsf s^{R+2}$, length power at most $2t$ and exponential
$e^\Lambda$. There are at most $b(\mathfrak d+3)^{2(b-1)}$ records.
Together with \eqref{eq:critical-graded-tail}, this proves the same
graded carrier error as \cite[Proposition~4.7]{HuangGap}, on the new
envelope:
\[
 \Bigl|F^{h_0,a_0,k;P}_{g,n}(\mathbf y)
  -\sum_{j=1}^Rg^{-j}\mathcal F^{h_0,a_0,k;P}_{j,n}(\mathbf y)\Bigr|
 \le g^{-R-1}\mathsf s^{R+2}(C_Mt)^{C_Mt}(1+\Lambda)^{2t}e^{\Lambda}.
\]
For a nonseparating kernel, the scalar handle-ratio expansion
times the two-boundary expansion has the smaller cusp power $\mathsf s^{R+1}$.
For a separating kernel, the retained endpoints use the same shifted
estimate and the omitted ones use \eqref{eq:critical-separating-tail}.

The inequalities $N>2R$ and $g>H_0$, which hold since $R,H_0\le t\le g^\theta$,
justify reconstruction of retained
records. Unavailable cusp allocations are omitted before evaluating a
volume. When $R=0$ or $m>R$, apply the tail estimate directly rather
than a negative-order shifted expansion. The canonical finite sums
are consequently unchanged; this argument does not interpolate in $n$.

\emph{Integration over orbital types.}
For a length support $[-T,T]$, the whole-carrier bound
$\pi(g-1+n/2)>T$ of Definition~\ref{rev:def-huang-carriers}(i), which
holds since $T\le\rho\le g^\theta$, excludes filled surfaces and entirely
self-paired carriers. The remaining types have
$1\le m\le\lfloor2T/\pi\rfloor$, $k\le m+2$ and
$t=R+m+1\le R+\lceil T\rceil+1=\rho$.
The exact rooted integration formula
\eqref{rev:eq-huang-carrier-integration} and its fiber and
boundary-simplex bounds (Lemma~\ref{rev:lem-huang-records}(d)) give, per
pattern, an error at most
\[
 g^{-R-1}\mathsf s^{R+2}(C_Mt)^{C_Mt}
 Cm(1+T)(1+2T)^{2t+k}e^{3T}\|H\|_\infty .
\]
There are at most $(k+1)^{2k}$ patterns and $(m+2)^2$ types at each
$m$. Absorbing these factors, and adding the simple part obtained from
\eqref{rev:eq-huang-simple-integration} with the kernel errors of
Lemma~\ref{rev:lem-huang-records}(c), proves
\eqref{eq:critical-orbital-remainder}. Simple repetitions have the
same absolutely integrable bound. The zero-support case is zero and
an unavailable nonsimple type contributes an empty sum. Absolute
bounds also cover signed tests. The coefficient vanishing rules
\cite[(3.32)]{HuangGap}, recalled in
Lemma~\ref{rev:lem-huang-records}(d), identify the retained finite sum
with the canonical coefficients at every $R,T$.

\subsection{Growing-degree polynomial tests and eigenvalue counts}\label{rev:subsec-critical-polynomial}

Finally we apply the orbital estimate
\eqref{eq:critical-orbital-remainder} to polynomial tests whose degree
grows with the genus, and deduce the local expected-count estimate for
eigenvalues. There are three stages: the monomial estimate
\eqref{eq:critical-monomial-error}, which also supplies the first-moment
hypothesis of Lemma~\ref{rev:lem-huang-first-moment}; the transfer to
polynomials of degree $Q=\lfloor g^\sigma\rfloor$ through
Lemma~\ref{rev:lem-huang-transfer}, which proves
\eqref{eq:critical-polynomial-error}; and the count bound, obtained
from the detector polynomial and the positive trace of
Lemmas~\ref{rev:lem-huang-trace}--\ref{rev:lem-huang-retained-mass}
and~\ref{rev:lem-huang-detector}.

\emph{Monomial tests.}
Apply \eqref{eq:critical-orbital-remainder} to $H=h^{*j}$ with
$1\le j<q$, $R=q-2$, and $T=jL$, where $\operatorname{supp}h\subset[-L,L]$.
Its supremum is at most $\|h\|_\infty$, and
$\rho=q-1+\lceil jL\rceil\le(1+L)q$. For
$q\le g^\theta/(2+L)$, division by $g$ gives
\begin{equation}\label{eq:critical-monomial-error}
 \left|g^{-1}\mathbb E\mathcal O_X(x^j)
       -\sum_{i=1}^{q-1}g^{-i}\omega_{i,n}(x^j)\right|
       \le g^{-q}\mathsf s^q(C_{M,f}q)^{C_{M,f}q},
\end{equation}
since $(C_M\rho)^{C_M\rho}e^{3jL}\|h\|_\infty\le(C_{M,f}q)^{C_{M,f}q}$.
The base case $q=2,j=1$ gives
$|\mathbb E\mathcal O_X(x)-\omega_1(x)|\le\mathsf s^2(2C_{M,f})^{2C_{M,f}}/g
\le M^2(2C_{M,f})^{2C_{M,f}}\le g$ for large $g$, so
Lemma~\ref{rev:lem-huang-first-moment} bounds the first moment of the
positive truncated trace by $D_f(g+nQ)$ for every cutoff $Q$; a
large-degree genus threshold is not required for this base-moment bound.

\emph{Polynomial transfer.}
Take $Q=\lfloor g^\sigma\rfloor$ with
\[
 0<\sigma<\theta,\qquad \sigma C_{M,f}<1/4,
 \qquad 3\sigma D_0<1/2 .
\]
For $p=x\widetilde p$ of degree at most $Q$, put $q=Q+2$ in
\eqref{eq:critical-monomial-error}, so that
\eqref{rev:eq-huang-monomial-hypothesis} holds with
$E=g^{-q}\mathsf s^q(C_{M,f}q)^{C_{M,f}q}$. The monomial coefficient bound
of Lemma~\ref{rev:lem-huang-coefficients}(a)
contributes $7^Q$. Relative to $(\mathsf s^2/g)\|\widetilde p\|_{[0,a]}$,
the terminal error $7^QgE\,N_*(p)$ is bounded by
\[
 7^Q[C_{M,f}(Q+2)]^{C_{M,f}(Q+2)}(\mathsf s/g)^Q\le2^{-Q}
\]
for sufficiently large $g$: its logarithm divided by $Q\log g$ has
limit superior at most $\sigma C_{M,f}-1/2<0$, because
$\log(\mathsf s/g)\le-\frac12\log g+\log M$.
The genus-independent canonical coefficient estimate of
Lemma~\ref{rev:lem-huang-coefficients}(b), with $A_Q=(3C_0Q)^{D_0}$, gives
\[
 \sum_{i=3}^{Q+1}g^{1-i}|\omega_{i,n}(p)|
 \le\frac{2\mathsf s^3A_Q^3}{g^2}\|\widetilde p\|_{[0,a]},
 \qquad \mathsf sA_Q/g=O_M\bigl(g^{\sigma D_0-1/2}\bigr)=o(1).
\]
Lemma~\ref{rev:lem-huang-transfer} therefore
proves \eqref{eq:critical-polynomial-error}, with an $o(1)$
right side since $\mathsf s^3A_Q^3/g^2=O_M(g^{3\sigma D_0-1/2})$.
The same reasoning holds for $Q=\lceil K\log g\rceil$
for any fixed $K$.

\emph{Eigenvalue counts.}
Finally apply the binomial detector, exact constant cancellation,
parabolic bound and retained-mass estimate of
Lemmas~\ref{rev:lem-huang-trace}, \ref{rev:lem-huang-retained-mass}
and~\ref{rev:lem-huang-detector} at
$Q=\lceil K_\varepsilon\log g\rceil$ with
$K_\varepsilon\ge\max\{3/\kappa_\varepsilon,2/L\}$. Taking expectations in
\eqref{rev:eq-huang-trace-compressed} at $p=p_Q$, using
\eqref{rev:eq-huang-count-from-trace}, the bounds
\eqref{rev:eq-huang-deterministic} and $|\mu_0(p_Q)|\le5N_{\mathbb R}(p_Q)$,
the detector norms \eqref{rev:eq-huang-detector}, the
degree-independent bound \eqref{rev:eq-huang-degree-independent}, and
the polynomial error \eqref{eq:critical-polynomial-error} just
established at this $Q$, we obtain, exactly as in
\cite[proof of Theorem~6.13, (6.74)]{HuangGap},
\[
 \mathbb E N_X(1/4-\varepsilon)\le4\,\mathbb E\tau^\circ_{X,Q}(p_Q)\le
 C_{M,f,\varepsilon}\left\{
 (g+nQ^2)e^{-\kappa_\varepsilon Q}+\frac{\mathsf s^2}{g}
 +\frac{\mathsf s^3Q^{3D_0}}{g^2}+\frac{\mathsf s^2}{g}2^{-Q}
 +\frac{n}{\operatorname{area}(X)Y_Q}\right\}.
\]
Here $N_X(1/4-\varepsilon)=N_\varepsilon(X)$ is the count of
Lemma~\ref{rev:lem-huang-detector}. The exponential term is
$O(g^{-2}\log^2g)$ by the choice of $K_\varepsilon$, since $n\le g$ and
$e^{-\kappa_\varepsilon Q}\le g^{-3}$; the third term divided by
$\mathsf s^2/g$ is $O_M(Q^{3D_0}/\sqrt g)=o(1)$; and the truncation
height $Y_Q=3e^{LQ}\ge3g^2$ makes the last term $O(g^{-2})$, because
$\operatorname{area}(X)\ge2\pi n$. Thus all terms are
$O_{M,f,\varepsilon}(\mathsf s^2/g)$.
This is the local expected-count estimate used to remove the
eigenfunction-mass weights in Proposition~\ref{prop:critical-mean}.

\section{Proof of the moving spectral law}\label{app:moving-edge}

This appendix proves the moving spectral law,
Theorem~\ref{thm:moving-edge}; Subsections~\ref{subsec:edge-residual-proof}
and~\ref{subsec:edge-moment-proof} then prove
Corollaries~\ref{cor:edge-residual} and~\ref{cor:edge-moments}.
Throughout, $c=n/\sqrt g$, $\kappa_0=(8\pi^2)^{-1}$,
and the model functions are those of~\eqref{eq:edge-model-coordinate}.
For $a>0$, let $N_g(a)=\mathcal D_X([a,\infty))$ and let $Z_g(a)$
count two-cusp pairs of pants whose model deficits are at least $ac^2$.
The inputs are Huang's trace expansion~\cite{HuangGap}, used here
under his original genus thresholds, the model identification of
Proposition~\ref{prop:critical-identification}, the cut-graph bound
\eqref{eq:critical-graph-bound} of Appendix~\ref{app:critical}, and the
funnel model $F_s$ of Section~\ref{sec:critical-eigenvalues}.
The proof has three parts, kept separate so that their ranges are
explicit: a positive trace estimate for the mean of the eigenvalue
counts (first subsection), Poisson convergence and conditional
isolation of the pants counts (second subsection), and a spectral
comparison that matches the two counts and identifies the
eigenfunctions (third subsection).

\subsection{The moving expected count}

This subsection proves the limit \eqref{eq:edge-mean-count} of the
mean of the eigenvalue counts, using Huang's trace expansion with test
functions adapted to the moving scale. The main step is
Lemma~\ref{rev:lem:edge-trace-approximation}, a trace approximation
whose error terms are explicit in $c$; the inputs
\ref{rev:input-exact}--\ref{rev:input-retained} from \cite{HuangGap}
are listed before it, and the identification of the limit with the
model intensity uses Proposition~\ref{prop:critical-identification}.
In this subsection we assume only $c\to0$ and $\log(1/c)=o(\log g)$.
This is implied by~\eqref{eq:moving-edge-range}: if $c\log g\to\infty$
then $c\ge1/\log g$ eventually, so $\log(1/c)\le\log\log g$; it also
forces $n=g^{1/2-o(1)}\to\infty$.
Keep the fixed multiplier $f$, its endpoint $a_0=f(i/2)$ and its
support constant $L$ from Section~\ref{sec:critical-eigenvalues}.
The function $F(t)=f(i\sqrt t)$ is analytic on $[0,1/4]$ and
$F'(t)>0$ there. At zero this follows from the positive second
moment of $h=b_0*b_0$; at positive $t$ it follows by differentiating
the even integral against $\cosh(\sqrt t\,u)$. Its inverse $J$ is
smooth on $[1,a_0]$.

Let $\phi$ be smooth on $[0,\infty)$, zero below some $a>0$,
and either compactly supported or constant beyond a finite point.
Put $\psi_c(x)=\phi(J(x)/c^2)$ for $x\ge1$ and zero for $x\le1$.
Extend smoothly beyond $a_0$ when necessary. For each fixed $r$,
\begin{equation}\label{eq:edge-test-norm}
 \|\psi_c/x\|_{C^r([0,a_0])}\le C_{\phi,r}c^{-2r}.
\end{equation}
Choose $Q=\lfloor g^\sigma\rfloor$ with
$0<\sigma<\min\{(2D_3)^{-1},(6D_0)^{-1},1/4\}$, using Huang's
constants~\cite[Proposition~6.9]{HuangGap}. His original threshold
applies: $(n+1)^2/g\to0$ and $Q^{D_3}(n+1)/g\to0$.
No critical-envelope extension is needed in this range.

Let $r_Q$ be the shifted Chebyshev truncation of $\psi_c/x$ of
degree $Q-1$ and $v_Q=xr_Q$. For fixed $K,m$, integration by parts
after $x=a_0(1+\cos\theta)/2$ and \eqref{eq:edge-test-norm} with
$r=K+2m+2$ give
\begin{equation}\label{rev:eq:edge-chebyshev}
 e_m:=\|r_Q-\psi_c/x\|_{C^m([0,a_0])}
 \le C_{K,m}Q^{-K}c^{-2(K+2m+2)}=Q^{-K}g^{o(1)}.
\end{equation}
Indeed the $j$th Chebyshev coefficient costs $j^{-K-2m-2}$ and
its $m$th derivative costs at most $C_mj^{2m}$.

\emph{Inputs from Huang's trace expansion.}
We now record, in the notation of~\cite{HuangGap} and of
Proposition~\ref{prop:critical-mean}, the statements of that paper
that enter the trace estimate below; each is quoted with its own
hypotheses, and none of them is modified here. Let $A_X$ be the
bounded self-adjoint operator with multiplier $f$: a discrete
eigenvalue $0<\lambda_j<1/4$ with eigenfunction $u_j$ has
$A_X$-eigenvalue $x_j=F(1/4-\lambda_j)=f(i\sqrt{1/4-\lambda_j})\in(1,a_0]$,
the constant eigenfunction has $A_X$-eigenvalue $a_0$, and the part of
the spectrum of $\Delta_X$ in $[1/4,\infty)$ is mapped into $[0,1]$.
Let $\chi_Q$ be the indicator of the complement of the $n$ width-one
cusp regions above height $Y_Q=3e^{LQ}$, and let $P_0$ be the
projection onto constants. For a real polynomial $p=x\widetilde p$
with $p(0)=0$ and $\deg p\le Q$, write
$N_*(p)=\|\widetilde p\|_{[0,a_0]}$,
$N_{\mathbb R}(p)=\|\widetilde p\|_{[0,1]}$, and
$\tau_{X,Q}(p)=\operatorname{Tr}(\chi_Qp(A_X)\chi_Q)$. We write
$\mathcal I(p)$ and $\mathcal P_{Y_Q}(p)$ for Huang's identity and
parabolic terms, $\mathcal O_X(p)$ for his orbital term, and
$\omega_1,\omega_{2,n}$ for the first two canonical coefficients of
the expansion of $\mathbb E\mathcal O_X(p)$ in powers of $1/g$, as in
\eqref{eq:critical-polynomial-error}. Put $\mathsf s=n+1$, as in Appendix~\ref{app:critical-trace}.
\begin{enumerate}[label=\textup{(I\arabic*)},leftmargin=*]
\item\label{rev:input-exact}
\textup{(\cite[Theorem~5.3, Proposition~5.1 and equations~(5.44)--(5.45)]{HuangGap}.)}
For every such $X$, $Q$ and $p$, the compression $\chi_Qp(A_X)\chi_Q$
is trace class and
\[
 \operatorname{Tr}\bigl(\chi_Qp(A_X)(I-P_0)\chi_Q\bigr)
 =\Bigl(\operatorname{area}(X)-\frac n{Y_Q}\Bigr)\mathcal I(p)
  +n\mathcal P_{Y_Q}(p)+\mu_0(p)+\frac{n\,p(a_0)}{\operatorname{area}(X)Y_Q}
  +\bigl(\mathcal O_X(p)-\omega_1(p)\bigr),
\]
where $\mu_0$ is a fixed signed measure supported in $[0,1]$ with
total variation less than $5$, and $\omega_1(p)=p(a_0)+\mu_0(p)$ is
independent of $n$. Moreover $\tau_{X,Q}(p)\ge0$ when $p\ge0$ on the
spectrum of $A_X$, and $|\tau_{X,Q}(p)|\le N_*(p)\tau_{X,Q}(x)$.
More generally, if $\psi=x\xi$ with $\xi$ bounded and Borel on
$[0,a_0]$, then $\chi_Q\psi(A_X)(I-P_0)\chi_Q$ is trace class and
$|\operatorname{Tr}(\chi_Q\psi(A_X)(I-P_0)\chi_Q)|
\le\|\xi\|_{[0,a_0]}\tau_{X,Q}(x)$
\cite[(5.29) and (6.63)]{HuangGap}.
\item\label{rev:input-deterministic}
\textup{(\cite[Proposition~5.4]{HuangGap}.)}
There are constants $I_f,J_f$, depending only on $f$, such that for
every integer $Q\ge1$ and every such $p$,
$|\mathcal I(p)|\le I_fN_{\mathbb R}(p)$ and
$|\mathcal P_{Y_Q}(p)|\le J_fQ^2N_{\mathbb R}(p)$. These bounds use
only the values of $\widetilde p$ on the real spectral image $[0,1]$.
\item\label{rev:input-first-moment}
\textup{(\cite[Lemma~6.8, equation~(6.43)]{HuangGap}.)}
There is an integer $c_3\ge2$, depending only on $f$, such that
$\mathbb E\tau_{X,Q}(x)\le D_f(g+nQ)$ for every $Q\ge1$ whenever
$g\ge c_3\mathsf s^2+[2c_3]^{c_3}\mathsf s$; here $D_f$ depends only on $f$.
At this threshold every degree-at-most-$Q$ polynomial trace is
absolutely integrable.
\item\label{rev:input-orbital}
\textup{(\cite[Proposition~6.9, equations~(6.46)--(6.47)]{HuangGap}.)}
With $D_0,C_0$ from \cite[Lemma~6.7]{HuangGap}, $D_3=\max\{D_0,3c_3\}$
and $A_Q=(3C_0Q)^{D_0}$, there is $C_{\rm add}\ge2$, depending only on
$f$, such that for all integers $n,Q\ge1$ with
$g\ge C_{\rm add}\mathsf s^2+C_{\rm add}Q^{D_3}\mathsf s$ and every such $p$,
\[
 \bigl|\mathbb E\mathcal O_X(p)-\omega_1(p)-g^{-1}\omega_{2,n}(p)\bigr|
 \le\Bigl\{\frac{2\mathsf s^3A_Q^3}{g^2}+\frac{\mathsf s^2}{g}2^{-Q}\Bigr\}N_*(p).
\]
This threshold implies the one in \ref{rev:input-first-moment}.
\item\label{rev:input-distribution}
\textup{(\cite[Lemma~6.10 and Theorem~6.11]{HuangGap}, with the quadratic identity \cite[(3.56)]{HuangGap}.)}
Put $m=2D_0+2$. For each $n\ge1$ the functional
$\mathscr K_n(q)=\omega_{2,n}(xq)$ on polynomials extends uniquely to a
distribution supported in $[0,a_0]$ of order at most $m$, with
$|\mathscr K_n(\varphi)|\le K_f(n+1)^2\|\varphi\|_{C^m([0,a_0])}$ for a
constant $K_f$ depending only on $f$. The exact identity
$\omega_{2,n}=W_0+nW_1+\binom n2W_2$ gives
$\mathscr K_n=\mathscr K^{(0)}+n\mathscr K^{(1)}+\binom n2\mathscr T_2$,
where $\mathscr T_2=\mathscr K_3-2\mathscr K_2+\mathscr K_1$,
$\mathscr K^{(1)}=3\mathscr K_2-\mathscr K_3-2\mathscr K_1$ and
$\mathscr K^{(0)}=3\mathscr K_1-3\mathscr K_2+\mathscr K_3$ are the
weighted extensions of $W_2$, $W_1$ and $W_0$; each has order at most
$m$ and norm at most $C_f$ on $C^m([0,a_0])$. On $(1,a_0]$ the
distribution $\mathscr T_2$ is the positive measure $x\eta_2$ of
\eqref{eq:critical-weighted-distribution}, and
$\mathscr T_2(\psi/x)=\int\psi\,d\eta_2$ for smooth $\psi$ vanishing
near $[0,1]$, including tests with $\psi(a_0)\ne0$.
\item\label{rev:input-retained}
\textup{(\cite[Lemma~6.5]{HuangGap}.)}
If $u$ is a normalized eigenfunction with eigenvalue
$0<\lambda\le1/4-\varepsilon$, then
$\|\chi_Qu\|_2^2\ge1-e^{-2\sqrt\varepsilon LQ}$, uniformly in the
surface and in its cusp count.
\end{enumerate}
Inputs \ref{rev:input-exact}--\ref{rev:input-distribution} are used
in the following lemma; \ref{rev:input-retained} is used afterwards
to remove the truncation weights.

\begin{lemma}[Trace approximation for $c$-dependent tests]
\label{rev:lem:edge-trace-approximation}
Let $\psi_c$ be smooth on $[0,a_0]$, vanishing on $[0,1]$, with
\[
 \|\psi_c/x\|_{C^r([0,a_0])}\le M_rc^{-2r}
 \qquad\text{for every fixed }r\ge0,
\]
the constants $M_r$ being independent of $g$; by
\eqref{eq:edge-test-norm} this includes $\psi_c=\phi(J(x)/c^2)$ for
$\phi$ as above. Let $r_Q$, $v_Q=xr_Q$ and $e_m$ be as in
\eqref{rev:eq:edge-chebyshev}, whose bound holds for every such
$\psi_c$, and let
$Q=\lfloor g^\sigma\rfloor$ with $\sigma$ as above. Suppose
$(n+1)^2/g\to0$ and $Q^{D_3}(n+1)/g\to0$. For an orthonormal basis
$u_j$ of eigenfunctions with $0<\lambda_j<1/4$ and $x_j=F(1/4-\lambda_j)$,
put
\[
 \mathcal D_{X,Q}(\psi_c)
 =\operatorname{Tr}\bigl(\chi_Q\psi_c(A_X)(I-P_0)\chi_Q\bigr)
 =\sum_{0<\lambda_j<1/4}\psi_c(x_j)\|\chi_Qu_j\|_2^2.
\]
Then, with $m=2D_0+2$ and
a constant $C$ depending only on $f$ and the constants $M_r$, for all
large $g$,
\begin{align}\label{eq:edge-trace-error}
 \Bigl|\mathbb E\mathcal D_{X,Q}(\psi_c)
 -\frac1g\binom n2\mathscr T_2(\psi_c/x)\Bigr|
 \le C\bigg\{&(g+nQ^2)e_0+\frac{(n+1)^2}{g}e_m
       +\frac{n+1}{g}c^{-2m}
       +\frac{(n+1)^3Q^{3D_0}}{g^2}\\
       &\hspace{15mm}+\frac{(n+1)^2}{g}2^{-Q}
       +\frac{n}{\operatorname{area}(X)Y_Q}\bigg\},\notag
 \qquad Y_Q=3e^{LQ}.
\end{align}
\end{lemma}

\begin{proof}
Write $\xi_c=\psi_c/x$. Since $\psi_c$ vanishes on $[0,1]$, the test
is zero on the whole real spectral image, and no sign is attributed to
the parabolic term. Hence $N_{\mathbb R}(v_Q)=\|r_Q\|_{[0,1]}\le e_0$,
$N_*(v_Q)\le\|\xi_c\|_\infty+e_0\le M_0+e_0$, and
$|v_Q(a_0)|\le a_0N_*(v_Q)\le C$. For $\psi_c=\phi(J(x)/c^2)$ the
vanishing holds on $[0,F(ac^2)]\supset[0,1]$ because $\phi$ vanishes
below $a$.
The hypotheses $(n+1)^2/g\to0$ and $Q^{D_3}(n+1)/g\to0$ make the
threshold of \ref{rev:input-orbital}, and hence that of
\ref{rev:input-first-moment}, hold for all large $g$. The operator
$\psi_c(A_X)(I-P_0)$ has finite rank, so the displayed sum
representation of $\mathcal D_{X,Q}(\psi_c)$ holds, and
$\mathcal D_{X,Q}(v_Q)=\operatorname{Tr}(\chi_Qv_Q(A_X)(I-P_0)\chi_Q)$
is the quantity in \ref{rev:input-exact}; all expectations below are
finite by \ref{rev:input-first-moment}. The estimate follows from the exact identity and the following four
bounds, each attributed to the input it uses.

\emph{Smooth test to polynomial}, by
\ref{rev:input-exact} and \ref{rev:input-first-moment}.
The spectral domination of \ref{rev:input-exact}, applied to
$\psi_c-v_Q=x(\xi_c-r_Q)$, whose quotient is bounded by $e_0$ on
$[0,a_0]$, and the first moment of \ref{rev:input-first-moment} give
\[
 \bigl|\mathbb E\mathcal D_{X,Q}(\psi_c)-\mathbb E\mathcal D_{X,Q}(v_Q)\bigr|
 \le e_0\,\mathbb E\tau_{X,Q}(x)\le D_f(g+nQ)e_0.
\]

\emph{Exact identity}, by \ref{rev:input-exact}.
\[
 \mathcal D_{X,Q}(v_Q)
 =\Bigl(\operatorname{area}(X)-\frac n{Y_Q}\Bigr)\mathcal I(v_Q)
   +n\mathcal P_{Y_Q}(v_Q)+\mu_0(v_Q)
   +\frac{n\,v_Q(a_0)}{\operatorname{area}(X)Y_Q}
   +\mathcal O_X(v_Q)-\omega_1(v_Q).
\]
The constant eigenfunction is removed here, before any approximation;
its remaining truncation correction is the fourth term, and nothing is
required of $\psi_c(a_0)$.

\emph{Identity, parabolic and $\mu_0$ terms}, by
\ref{rev:input-deterministic} and \ref{rev:input-exact}.
Since $|v_Q|\le N_{\mathbb R}(v_Q)\le e_0$ on $[0,1]$ and
$|v_Q(a_0)|\le C$,
\begin{align*}
 \Bigl|\Bigl(\operatorname{area}(X)-\frac n{Y_Q}\Bigr)\mathcal I(v_Q)
   +n\mathcal P_{Y_Q}(v_Q)+\mu_0(v_Q)\Bigr|
 &\le\bigl(I_f\operatorname{area}(X)+J_fnQ^2+5\bigr)e_0
 \le C(g+nQ^2)e_0,\\
 \Bigl|\frac{n\,v_Q(a_0)}{\operatorname{area}(X)Y_Q}\Bigr|
 &\le\frac{Cn}{\operatorname{area}(X)Y_Q}.
\end{align*}

\emph{Orbital expansion}, by \ref{rev:input-orbital}.
Since $N_*(v_Q)\le C$,
\[
 \Bigl|\mathbb E\mathcal O_X(v_Q)-\omega_1(v_Q)
   -\frac1g\omega_{2,n}(v_Q)\Bigr|
 \le C\Bigl\{\frac{(n+1)^3Q^{3D_0}}{g^2}+\frac{(n+1)^2}{g}2^{-Q}\Bigr\}.
\]

\emph{The quadratic coefficient}, by \ref{rev:input-distribution}.
By the order-$m$ bound and the hypothesis with $r=m$,
\begin{align*}
 \bigl|\omega_{2,n}(v_Q)-\mathscr K_n(\xi_c)\bigr|
 &=|\mathscr K_n(r_Q-\xi_c)|\le K_f(n+1)^2e_m,\\
 \Bigl|\mathscr K_n(\xi_c)-\binom n2\mathscr T_2(\xi_c)\Bigr|
 &=\bigl|\mathscr K^{(0)}(\xi_c)+n\mathscr K^{(1)}(\xi_c)\bigr|
  \le C(n+1)\|\xi_c\|_{C^m}\le C(n+1)c^{-2m}.
\end{align*}

Taking the expectation of the exact identity, dividing the last two
bounds by $g$, and adding all the bounds proves
\eqref{eq:edge-trace-error}.
\end{proof}

The proof uses $c$ only through the seminorm hypothesis. If a test
$\psi$, smooth on $[0,a_0]$ and vanishing on $[0,1]$, satisfies
$\|\psi/x\|_{C^r([0,a_0])}\le M_r$ for every fixed $r$, with $M_r$
independent of $g$, the same argument gives
\eqref{eq:edge-trace-error} with $c^{-2m}$ replaced by $1$ and with
$e_m\le C_{K,m}Q^{-K}$; this form is used for fixed physical tests in
Subsection~\ref{subsec:edge-moment-proof}.

For large fixed $K$, every term in~\eqref{eq:edge-trace-error} tends
to zero. Indeed $(g+nQ^2)e_0\le g^{1+2\sigma-\sigma K+o(1)}$ by
\eqref{rev:eq:edge-chebyshev}; $c^{-2m}=g^{o(1)}$ because
$\log(1/c)=o(\log g)$, while $(n+1)/g\to0$;
$(n+1)^3Q^{3D_0}/g^2\le[(n+1)^2/g]\,g^{3\sigma D_0-1/2}$ eventually,
and $3\sigma D_0<1/2$; and $n/(\operatorname{area}(X)Y_Q)\le e^{-LQ}$.
For $\psi_c=\phi(J(x)/c^2)$,
Proposition~\ref{prop:critical-identification} and
\ref{rev:input-distribution} give
\[
 \mathscr T_2(\psi_c/x)=\int\psi_c\,d\eta_2
 =2\kappa_0\int_0^\infty\phi\bigl((1/4-\beta(s))/c^2\bigr)\sinh(s/2)\,ds,
\]
since $w(s)=\sinh(s/2)/(4\pi^2)=2\kappa_0\sinh(s/2)$ and
$\psi_c(x_s)=\phi((1/4-\beta(s))/c^2)$. The third term
in~\eqref{eq:edge-trace-error} is the cost of removing $W_0+nW_1$.
The resulting weighted spectral sum has mean
\begin{equation}\label{eq:edge-weighted-mean}
 \mathbb E\mathcal D_{X,Q}(\psi_c)=
 \frac{n(n-1)}g\kappa_0
       \int_0^\infty\phi(h_*(y)/c^2)\,dy+o(1).
\end{equation}
This is the physical form of the local model identification, since
$y=2(\cosh(s/2)-1)$ gives $dy=\sinh(s/2)\,ds$ and
$h_*(y)=1/4-\beta(s)$.

For modes on the support, the omitted mass is at most
$\epsilon_g=e^{-2c\sqrt a LQ}\to0$, by~\cite[Lemma~6.5]{HuangGap}.
Apply~\eqref{eq:edge-weighted-mean} first to a nonnegative smooth
majorant of $|\phi|$. Its weighted expectation is bounded, by
\eqref{eq:edge-model-tail} and the change $y=z/c^2$.
The unweighted-minus-weighted expectation is at most
$\epsilon_g/(1-\epsilon_g)$ times that bound. Dominated convergence
in the same change of variables gives
\begin{equation}\label{eq:edge-mean-count}
 \mathbb E\langle\phi,\mathcal D_X\rangle\longrightarrow
   \int_0^\infty\phi(x)\frac{\kappa_0\,dx}{x^2},
 \qquad \mathbb E N_g(a)\longrightarrow\frac{\kappa_0}{a}.
\end{equation}
For the second assertion, squeeze the tail indicator by smooth tests
constant near infinity. The limiting measure has no atom at $a$.
This argument bounds unweighted counts without assuming such a bound
in advance.

\subsection{Growing waists and conditional isolation}

We next prove geometric Poisson convergence of the pants counts
$Z_g(a)$ and their isolation for any
$R=o(\log g)$, $R\ge1$, and any cutoff width $\ell\ge0$ with
$\ell=o(\log g)$, still assuming only
$\log(1/c)=o(\log g)$. The argument is the factorial-moment method of
Appendix~\ref{app:critical} with waists growing like $\log(1/c)$:
the boundary volumes are controlled at growing lengths by
\cite[Proposition~2.3]{HuangGap}, and the cut-graph bound
\eqref{eq:critical-graph-bound} is applied conditionally on a marked
tuple to show that almost every pair of pants is good, that is,
isolated at radius $R'$ with an embedded funnel neighborhood. The
outputs used later are the conditional bound
\eqref{eq:edge-conditional-bad}, the discard estimate
\eqref{eq:edge-discard} and the Poisson limits of the good counts.
The isolation radius is
\begin{equation}\label{rev:eq:edge-good-radius}
 R'=R+\ell+1,
\end{equation}
chosen so that the cutoff functions of the next subsection, which are
supported in the radius-$(R+\ell)$ funnel neighborhoods, have
disjoint supports for distinct good pants; there $\ell=1/c$ and $R$ are
fixed by~\eqref{eq:edge-common-radius}, and $R'=o(\log g)$ exactly
when $R=o(\log g)$.
For fixed $a>0$, let $S_a$ solve
$q(2\cosh(S_a/4))=ac^2$. Equation~\eqref{eq:edge-model-tail} gives
\begin{equation}\label{eq:edge-waist-bound}
 S_a=4\log(1/c)-2\log a+o(1).
\end{equation}

\emph{Boundary volumes at growing lengths.}
For fixed $r$, the order-zero joint intersection estimates
\cite[Proposition~2.3]{HuangGap}, originating in
\cite[Theorem~4.2]{HMT}, imply, uniformly for $p\le\sqrt g$ and
$0\le s_i\le S_a$,
\begin{equation}\label{eq:edge-boundary-volume}
 \frac{V_{g,p}(\mathbf s,\mathbf0)}{V_{g,p}}
 =\prod_{i=1}^r\frac{\sinh(s_i/2)}{s_i/2}
 \left(1+O_r\left(\frac{p+1}{g}(1+rS_a)^{r+2}\right)\right).
\end{equation}
For the growing-length uniformity, write $[\tau_{\mathbf d}]_{g,p}$,
$\mathbf d\in\mathbb Z_{\ge0}^p$, $|\mathbf d|=\sum_id_i$, for the
normalized intersection numbers of \cite[Proposition~2.3]{HuangGap},
and $E_2(v)=1+v+v^2/2$. The order-zero case of that proposition bounds
$\bigl|[\tau_{\mathbf d}]_{g,p}/V_{g,p}-1\bigr|$ by
$C(p+1)g^{-1}E_2(|\mathbf d|)$. Beyond the dimension cutoff the
bracket is zero and $|\mathbf d|\ge g$, so the same bound holds
there. Summation against the boundary series gives
$C_r(p+1)g^{-1}e^{\sum s_i/2}(1+\sum s_i)^2$.
Divide by the leading product using
$\sinh(s/2)/(s/2)\ge e^{s/2}/(1+s)$.
The order-zero genus threshold is linear in $p$, and hence permits
these parameters. Also
$V_{g,j}/V_{g,j+1}=(8\pi^2g)^{-1}(1+O((j+1)/g))$.

\emph{The Campbell integral.}
The exact contribution of ordered disjoint $r$-tuples of pants with
waists in prescribed sets $A_i\subset[0,S_a]$ is
\begin{equation}\label{eq:edge-campbell}
 \frac{(n)_{2r}}{2^rV_{g,n}}
 \int_{\prod_i A_i}V_{g,n-r}(\mathbf s,\mathbf0_{n-2r})
                          \prod_i s_i\,ds_i.
\end{equation}
The factor $2^{-r}$ chooses unordered cusp pairs; the tuple is ordered.
Distinct disjoint pants cannot be properly nested, since their
difference would be an annulus between distinct geodesics. Their
complement is connected. Thus there are no other complementary types
in~\eqref{eq:edge-campbell}. By~\eqref{eq:edge-boundary-volume}, this
integral is asymptotic to the product with one-variable measure
\[
 \mu_g(ds)=\kappa_0c^2\sinh(s/2)\,ds.
\]
Under $s\mapsto q(2\cosh(s/4))/c^2$, its upper-tail mass tends to
$\kappa_0/a$, since $y=2(\cosh(s/2)-1)$ and $h_*(y)\sim1/y$.

\emph{Conditional cut-graph bounds.}
Conditioning in the Campbell integral leaves a WP complement of
genus $g$, with $n-2r$ cusps, $r$ labeled boundaries and uniform
gluing twists. Its total number of ends is $N=n-r$.
More explicitly, if $\mathbf s$ are the marked waist lengths, the
conditional law of the complement is its WP measure divided by
$V_{g,N}(\mathbf s,\mathbf0)$; each gluing twist has law
$dt_i/s_i$. Cutting a witness graph $\Gamma$ therefore gives, for
fixed cusp labels, integrals of the form
\[
 C_\Gamma\int_{\sum_e\ell_e\le B}
 \frac{\prod_v V_{g_v,n_v}(\mathbf s_v,\boldsymbol\ell_v,
                            \mathbf0)}
      {V_{g,N}(\mathbf s,\mathbf0)}
                  \prod_e\ell_e\,d\ell_e,
\]
where $C_\Gamma$ is the fixed symmetry factor and the product runs
through the cut components. Summing the admissible genera and cusp
allocations is exactly the zero-boundary graph sum below. Positivity
gives $V_{g,N}(\mathbf s,\mathbf0)\ge V_{g,N}$, while the boundary
upper bound charges each old length once and each new length twice.
These observations justify applying the graph estimate under the
conditional, rather than the original unmarked, law.
For a fixed cut graph with $\nu\le2$ new edges and total new length
at most $B$, positivity of the boundary polynomials and the boundary
upper bound~\cite[Lemma~2.1]{HuangGap} cost at most
\begin{equation}\label{eq:edge-boundary-cost}
 e^{rS_a/2+B}(1+B)^{2\nu}
\end{equation}
relative to the zero-boundary allocation. Applying the proof of
\eqref{eq:critical-graph-bound} with $h=2g-2+N$, $z=N/h$ gives
$C_\nu h^{u-w+1}\phi_{0,1}(z)^u\phi_{0,2}(z)^v$.
The constants are uniform for $N\le\sqrt g$: the ambient-volume
loss is $e^{O(N^2/g)}$ and the finite Gamma-product loss is
$e^{O((m+1)^2/h)}$, $m\le N$.

A degree-two planar vertex, or a planar leaf with at least three
ends, supplies $O(g^{-1/2})$ after the remaining bounded factors.
For a planar leaf required to contain one specified old boundary,
the multinomial factor changes by $q/N$ when that leaf has $q$ ends.
Before summation this replaces $\phi_{0,1}$ by
$z\phi_{0,1}'/N$. In the leading one-leaf graph the factor is
\[
 h\frac{z\phi_{0,1}'(z)}N=O(N/g).
\]
Two specified old ends instead give
$hz^2\phi_{0,1}''(z)/(N(N-1))=O(1/g)$.
Additional leaves cost $h\phi_{0,1}(z)=O(1)$ and additional stable
vertices only decrease the bound. The distinguished labels cannot
be treated as freely chosen cusps.

\emph{Good pants and their witnesses.}
Declare a pair of pants good if its exterior radius-$R'$ funnel embeds and
its distance from every other pair of pants of waist at most $S_a$ exceeds
$2R'$, with $R'=R+\ell+1$ as in~\eqref{rev:eq:edge-good-radius}.
Every failure has a cut witness of total length at most
$B=20(r+1)(S_a+R'+1)=20(r+1)(S_a+R+\ell+2)$ and at most two new
edges. Indeed a first
normal collision supplies an essential returning or connecting arc
of length at most $2R'$; there are no focal points or collisions from
boundary-parallel arcs in the cyclic funnel cover. A candidate waist
crossing a marked pair of pants gives an essential excursion of length at
most $S_a$ in the complement: a peripheral excursion would bound a
geodesic bigon. Otherwise a nearby unmarked pair of pants lies in the
complement and is joined to an old boundary by an arc of length at
most $2R'$, unless such an arc first meets another old boundary.

Thicken a returning arc and its old boundary annulus. The result has
one old boundary and two new sides. After tightening and filling
peripheral once-punctured disks or annuli leading to old boundaries,
either both sides remain and give a planar degree-two vertex, or
one remains and gives a leaf with the old label and another end.
An empty disk would make the arc peripheral; both sides cannot
disappear in large genus. Coincident tightened sides are the two
incidences of one cut edge. A connecting arc gives a leaf with two
specified old ends. A joining arc to an unmarked pair of pants gives a leaf
with the old end and its two cusp ends. Thus every witness has one
of the suppressed factors above. The graph shapes are finite in
number; their genera and all end allocations have already been summed.
Only existence of a witness is needed, not a bound on the number of
arcs producing it. Markov's inequality gives, uniformly in the
marked lengths, labels and twists,
\begin{equation}\label{eq:edge-conditional-bad}
 \mathbb P(\text{a marked pants is not good}\mid\text{disjoint tuple})
 \le C_r e^{C_r(S_a+R'+1)}(1+S_a+R')^{C_r}g^{-1/2}
 =:\epsilon_{g,r}\longrightarrow0.
\end{equation}
The limit uses $S_a=O(\log(1/c))=o(\log g)$ and
$R'=R+\ell+1=o(\log g)$; since $\ell=o(\log g)$, the enlargement of
the radius from $R$ to $R'$ leaves the condition $R=o(\log g)$
unchanged.

\emph{Poisson limits of the good counts.}
Let $Z_g^{\rm good}$ count good pants, using the same $S_a,R,\ell$ for
all windows under consideration. The one-pants Campbell formula gives
\begin{equation}\label{eq:edge-discard}
 \mathbb E(Z_g(a)-Z_g^{\rm good}(a))
       \le\epsilon_{g,1}\mathbb E Z_g(a)=o(1).
\end{equation}
Every good tuple is disjoint, and the conditional probability that a
disjoint tuple is good is at least $1-r\epsilon_{g,r}$, by the union
bound over its $r$ marked pants in~\eqref{eq:edge-conditional-bad}.
Consequently its mixed factorial moments have the Poisson limits
given by~\eqref{eq:edge-campbell}. Higher fixed factorial moments
give tightness and moment passage; the Poisson law is determinate.
The good counts have bounded second moments by the disjoint-pair
integral. Equation~\eqref{eq:edge-discard} transfers both convergence
and uniform integrability to the full counts. No factorial-moment
bound for all intersecting tuples is inferred from a high-probability
event.

\subsection{Model cutoffs, count matching and eigenfunctions}

This subsection completes the proof of Theorem~\ref{thm:moving-edge}
under the range \eqref{eq:moving-edge-range}. We cut off the ground
states of the funnel models $F_s$ at a radius chosen in
\eqref{eq:edge-common-radius} and transplant them to the good pants of
the preceding subsection. Min--max then bounds the eigenvalue counts
from below by the pants counts; together with the limiting mean
\eqref{eq:edge-mean-count} this matches the two counts
(\eqref{eq:edge-count-matching}), and the transplanted states identify
the eigenfunctions. The mechanism is that the limiting mean equals the
pants intensity exactly. After localization errors and threshold shifts
are removed, an extra spectral count persisting with positive
probability in the limit would contradict this equality; min--max
supplies the lower bound from the pants. The cutoff estimates
\eqref{eq:edge-operator-cutoff}--\eqref{eq:edge-form-cutoff} are reused
in Subsections~\ref{subsec:edge-residual-proof}
and~\ref{subsec:edge-moment-proof}.

\emph{Cutoff model states.}
On the model funnel $dr^2+s^2\cosh^2r\,d\theta^2$, write the
normalized eigenfunction of energy $\beta=1/4-k^2$ in Fourier modes,
and put $v_m=\sqrt{s\cosh r}\,u_m$. Then
\begin{equation}\label{eq:edge-funnel-ode}
 -v_m''+\left[k^2+
   \left(\frac14+\frac{4\pi^2m^2}{s^2}\right)\operatorname{sech}^2r
              \right]v_m=0.
\end{equation}
Each mode has finite transformed $H^1$ energy. An $L^2$ solution
with zero boundary value vanishes by integration by parts, so a
nonzero mode divided by its value at any starting point is real and
positive on the remaining half-line. It tends to zero at infinity.
The maximum principle applied to its difference from $e^{-kt}$ gives
$|v_m(r+t)|\le e^{-kt}|v_m(r)|$. Integration and Parseval imply
\begin{equation}\label{eq:edge-funnel-mass}
 \int_{r\ge R}|u|^2\le e^{-2kR}.
\end{equation}
For $R\ge\ell>0$, multiply~\eqref{eq:edge-funnel-ode} by
$\zeta^2\overline v_m$, where $\zeta$ changes from zero to one on
$[R-\ell,R]$. Cauchy's inequality and the nonnegative potential give
\[
 \sum_m\int_R^\infty|v_m'|^2
       \le C\ell^{-2}e^{-2k(R-\ell)}.
\]
Cutoffs at infinity justify the identity by finite energy.
Choose $\chi=1$ up to $R$, zero beyond $R+\ell$, with derivatives
bounded by $C\ell^{-1}$ and $C\ell^{-2}$, retaining the whole pair of pants
and both complete cusps. The conjugated commutator is exactly
$-\chi''v-2\chi'v'$. Therefore
\begin{align}
 \|(\Delta-\beta)(\chi u)\|_2
    &\le C\ell^{-2}e^{-k(R-\ell)},\label{eq:edge-operator-cutoff}\\
 \frac{\mathcal E(\chi u)}{\|\chi u\|_2^2}-\beta
    &\le\frac{C\ell^{-2}e^{-2kR}}{1-e^{-2kR}}.
       \label{eq:edge-form-cutoff}
\end{align}
These are operator-domain vectors: their distributional Laplacian is
in $L^2$, and the complete-surface Laplacian has its maximal $L^2$
domain. Smooth extension by zero at the outer cutoff preserves this.

\emph{Transplantation and count matching.}
Now assume~\eqref{eq:moving-edge-range}. This is where
the condition $c\log g\to\infty$ itself enters: the cutoff radius $R$
must satisfy both $cR\to\infty$, for the funnel tail, and
$R=o(\log g)$, for the isolation estimate. Take
\begin{equation}\label{eq:edge-common-radius}
 \ell=c^{-1},\qquad R=(2+\sqrt{c\log g})/c.
\end{equation}
Then $R'=R+\ell+1=o(\log g)$, so the good event of the preceding
subsection, with isolation radius~\eqref{rev:eq:edge-good-radius}, is
available, and the cutoff vectors $\chi u$ of distinct good pants,
supported in their radius-$(R+\ell)$ funnel neighborhoods together
with the pants and cusps, have disjoint supports.
Uniformly for $k^2\ge ac^2$,
the operator error divided by $c^2$ is at most
$Ce^{-\sqrt a(1+\sqrt{c\log g})}\to0$.
The support area is at most $2\pi+s\sinh(R+\ell+1)=g^{o(1)}$.
Projecting a normalized vector off the global constant changes its
norm and operator residual by at most
$C\sqrt{\operatorname{area}(\operatorname{supp}u)/\operatorname{area}(X)}
=g^{-1/2+o(1)}=o(c^2)$.
For a bounded number of disjoint vectors the same estimate controls
their Gram matrix. Min--max and tightness of the model counts yield
\begin{equation}\label{eq:edge-minmax-count}
 \mathbb P\{N_g(a)\ge Z_g(a+\eta)\}\longrightarrow1
 \quad(a,\eta>0).
\end{equation}

On this event $(Z_g(a)-N_g(a))_+\le Z_g(a)-Z_g(a+\eta)$;
on its complement it is at most $Z_g(a)$. Uniform integrability
makes the latter contribution vanish. Let $g\to\infty$, then
$\eta\downarrow0$, and use~\eqref{eq:edge-mean-count}. This gives
\begin{equation}\label{eq:edge-count-matching}
 \mathbb E|N_g(a)-Z_g(a)|
 =\mathbb E N_g(a)-\mathbb E Z_g(a)
       +2\mathbb E(Z_g(a)-N_g(a))_+\longrightarrow0.
\end{equation}
Integer-valued counts therefore match with probability tending to
one at every finite list of thresholds. Their differences give joint
interval counts; the local mean bounds give point-measure tightness.
The void probability for $[x,\infty)$ gives~\eqref{eq:edge-gap-law}.
The transformation $x=\kappa_0/t$ gives the asserted arrival-time
representation. Sending a fixed lower threshold to zero shows that
each fixed number of eigenvalues exists with probability tending to one.

\emph{Eigenfunctions.}
For eigenfunctions, first restrict model deficits to a fixed $[a,b]$.
Partition it into finitely many bins, then exclude thin margins at
their endpoints and bins with two model points. The limiting Poisson
process is simple and atomless, so the limiting exceptional probability
can be made arbitrarily small by choosing the partition and margins.
Count matching at the endpoints gives one actual eigenvalue, counted
with multiplicity, in each occupied bin. It is simple. A model point
away from the margins has distance at least $\rho c^2$ from the
complement of its bin's spectral projection, for a fixed $\rho>0$.
The spectral theorem and~\eqref{eq:edge-operator-cutoff} project its
normalized transplanted state to this one-dimensional eigenspace
with $o(1)$ error. Refining the bins gives matching ordered energies.
Send $a\downarrow0$, $b\uparrow\infty$ after $g\to\infty$.
A deterministic waist cap
$4\log(1/c)+\log(c\log g)$ eventually exceeds every fixed $S_a$;
together with~\eqref{eq:edge-common-radius} it gives the asserted
deterministic area bound. This proves Theorem~\ref{thm:moving-edge}.

\subsection{Residual exhaustion}\label{subsec:edge-residual-proof}

We prove Corollary~\ref{cor:edge-residual}. The spectral subspace
below the moving threshold is shown to be residual by observing it on
fixed annuli in the cusps of the pants: the observation annihilates
cusp forms, while outside events of small probability it is injective
on that subspace. The observation argument
of~\eqref{eq:critical-observation-injection} requires a quantitative
lower bound when the model waist grows; the first part of this
subsection proves such a bound, \eqref{eq:edge-cusp-observation}, from
the Hecke polygon of the funnel model~\cite{Soares} and the Cheng--Yau
gradient estimate~\cite{ChengYau}. The second part combines it with
the cutoff estimate \eqref{eq:edge-operator-cutoff}, the count matching
\eqref{eq:edge-count-matching} and the Poisson limit of
Theorem~\ref{thm:moving-edge}. Let $u_s>0$ be the normalized
ground state, $w=2\cosh(s/4)$ and
$k=\sqrt{1/4-\beta(s)}$. Only large $s$ are relevant below, and we
assume throughout that $w\ge3$, that is $s\ge4\operatorname{arcosh}(3/2)$. Cut each cusp at a horocycle of length
$\mathsf b=1/4$, and let $K_s$ be the remaining pants core. In the first
cusp write the metric as $dt^2+\mathsf b^2e^{-2t}d\theta^2$ and set
$A_s=\{0\le t\le1\}$,
$h_s=\mathbf1_{A_s}/\sqrt{\operatorname{area}(A_s)}$.
We claim that for universal positive constants,
\begin{equation}\label{eq:edge-cusp-observation}
 m_s:=\langle u_s,h_s\rangle
       \ge C_0\sqrt{\frac{k}{1+s}}\,w^{-C}.
\end{equation}

\emph{The cusp observation bound.}
In a cusp, conjugation of a Fourier mode by $\sqrt{\mathsf b}e^{-t/2}$
gives $v_m''=[k^2+(2\pi m/\mathsf b)^2e^{2t}]v_m$.
Together with~\eqref{eq:edge-funnel-ode}, the same half-line
comparison proves, for each of the three ends,
\begin{equation}\label{eq:edge-end-boundary}
 \|u_s\|_{L^2(\mathrm{end})}^2
       \le(2k)^{-1}\int_{\partial\mathrm{end}}|u_s|^2.
\end{equation}
The boundary circumference is included by the conjugation and Parseval.

The needed core diameter bound follows directly from the
Hecke polygon~\cite[equations~(1), (6)--(7)]{Soares}.
Let $S(z)=-1/z$ and
$\mathcal F_w=\{|\Re z|\le w/2,\ |z|\ge1\}$.
The index-two surface has polygon $\mathcal F_w\cup S\mathcal F_w$.
The right funnel boundary lies on the axis
\[
 (x-w/2)^2+y^2=w^2/4-1.
\]
Its intersection with $|z|=1$ is at
$(2/w,\sqrt{1-4/w^2})$. The truncated core in $\mathcal F_w$
lies above the appropriate circle arcs, hence above $\sqrt5/3$,
and vertical motion upward remains in the core. The cusp horocycle
is at $Y=w/\mathsf b$, above those arcs; its $S$-image truncates the other
cusp, disjointly. A core point can be joined vertically to height
$Y$, horizontally to $iY$, then down to $i$, with total length
$O(1+\log w)$. The $S$-copy has the same bound and joins at $i$.
Thus $\operatorname{diam}(K_s)\le C(1+\log w)$.

Lift $u_s$ to $\mathbb H^2$. The positive function
$U(z,t)=e^{\sqrt{\beta(s)}t}u_s(z)$ is harmonic on
$\mathbb H^2\times\mathbb R$ in the analytic sign convention.
The Cheng--Yau gradient estimate~\cite{ChengYau} on fixed-radius
balls gives $|\nabla\log u_s|\le C$, independently of injectivity
radius. If $p$ is on the observed horocycle and
$M_s=\sup_{K_s}u_s$, the diameter bound gives
$M_s\le Cw^Cu_s(p)$. Since $\operatorname{area}(K_s)\le2\pi$,
\eqref{eq:edge-end-boundary} yields
\[
 1\le\left(2\pi+\frac{s+2\mathsf b}{2k}\right)M_s^2
       \le C(1+s)k^{-1}w^{2C}u_s(p)^2.
\]
Every point of $A_s$ is joined to $p$ by a path of bounded length;
its area is the fixed number $\mathsf b(1-e^{-1})$.
Harnack's inequality proves~\eqref{eq:edge-cusp-observation}.
In particular, for each fixed $0<a<b<\infty$ and
$a\le k^2/c^2\le b$, the model asymptotic gives
\begin{equation}\label{eq:edge-observation-power}
 m_s\ge C_{a,b}c^M
\end{equation}
for some fixed $M$: here $k\asymp c$, $w\asymp c^{-1}$ and
$s=4\log(1/c)+O_{a,b}(1)$.

\emph{Transplanted vectors.}
Under~\eqref{eq:edge-residual-range}, put $t_g=\log(1/c)$,
$H_g=c\log g/t_g\to\infty$, and choose
$R_g=3+\sqrt{t_g\log g/c}$ with transition width one.
Then $R_g=o(\log g)$ and $cR_g/t_g\to\infty$.
Equation~\eqref{eq:edge-operator-cutoff} gives an error smaller
than every fixed power of $c$ on a fixed rescaled deficit window.
After transplanting and dividing by $m_s$, the vectors $f_s$ have
\[
 \langle f_s,h_s\rangle=1,\qquad
 \|(\Delta_X-\beta(s))f_s\|_2=o(c^2).
\]
They retain the complete cusps and are not projected off constants.
Their norms need not be bounded; the observation and residual are
the quantities used below.

\emph{Injectivity on the spectral subspace.}
Fix $a>0$. Exclude a model deficit above $b$, a model point in
$[a,a+\eta]$, and a count $Z_g(a)>K$. By the Poisson limit their
limiting probabilities sum to at most
\[
 \frac{\kappa_0}{b}
 +\kappa_0\left(\frac1a-\frac1{a+\eta}\right)
 +\mathbb P\{\operatorname{Pois}(\kappa_0/a)>K\}.
\]
For fixed $b,\eta,K$, isolation and count matching fail with probability
$o(1)$. Otherwise $Z_g(a)=N_g(a)\le K$ and the selected deficits
lie in $[a+\eta,b]$. Set
$E=\mathbf1_{(0,1/4-ac^2]}(\Delta_X)L^2(X)$ and let $P$ project
onto $E$. For the transplanted vectors set $Fv=\sum_i v_if_i$ and
$Th=(\langle h,h_{A_i}\rangle)_i$. Disjoint annuli give $TF=I$
and $\|T\|\le1$. Each model energy is eventually at least $1/8$
and has distance at least $\eta c^2$ from the spectral complement
of $E$. Hence
\[
 \|TPF-I\|\le\sqrt K\max_i
       \frac{\|(\Delta_X-\beta_i)f_i\|_2}{\eta c^2}=o(1).
\]
Since $\dim E=N_g(a)$, $T$ is injective on $E$.
It annihilates cusp forms, whose constant Fourier coefficients vanish.
The residual--cuspidal decomposition therefore makes $E$ entirely
residual. Let $b,K\to\infty$, $\eta\downarrow0$ after $g\to\infty$.
The cuspidal count tends to zero in probability.
Equation~\eqref{eq:edge-count-matching} and uniform integrability of $Z_g(a)$
give uniform integrability of $N_g(a)$, so the expected cuspidal count
also tends to zero. This proves all the corollary's assertions.

\subsection{Rare events and actual gap moments}\label{subsec:edge-moment-proof}

This subsection proves Corollary~\ref{cor:edge-moments}. Since
the limiting gap law \eqref{eq:edge-gap-law} has infinite mean, the
moment statements require expectation-level estimates for discarded
geometries rather than an exchange of limit and expectation. The trace
approximation of Lemma~\ref{rev:lem:edge-trace-approximation}, applied
to tests whose seminorms are controlled uniformly in the threshold,
gives the upper bounds; the expected counts of good pants from the
second subsection of this appendix, together with the cutoff estimate
\eqref{eq:edge-form-cutoff}, give the lower bounds for the rare-gap
law and for moments of order at least one. The four parts below treat
the relative trace estimate, the rare-gap law
\eqref{eq:edge-rare-law}, the upper moment bounds and the lower moment
bounds.

\emph{The relative trace estimate.}
Put $\alpha_g=n(n-1)/g$. First assume only $n\to\infty$ and
$n=o(\sqrt g)$. For a fixed smooth physical-deficit test $\phi$
vanishing near zero, Lemma~\ref{rev:lem:edge-trace-approximation},
in the $c$-independent form recorded after its proof and applied to the
fixed test $\psi(x)=\phi(J(x))$, gives the relative estimate
\begin{equation}\label{eq:edge-relative-trace}
 \mathbb E\sum_{0<\lambda_j<1/4}\phi(1/4-\lambda_j)
 =\alpha_g\kappa_0\int_0^\infty\phi(h_*(y))\,dy+o(c^2).
\end{equation}
Indeed, retain the same $Q=g^\sigma$. After division by $c^2$,
the two lower cusp powers cost $O(1/n)$, the cubic remainder
$O((n+1)Q^{3D_0}/g)$, and the terminal remainder $O(2^{-Q})$.
The uniform approximation error is at most $O(g^2Q^{-K})$;
choose $\sigma K>2$. The fixed-order distributional error tends
to zero. The common-truncation and retained-mass errors are
exponentially small in $Q$. A nonnegative majorant first bounds
the weighted mean by $O(c^2)$, so removing weights preserves the
relative estimate. Smooth squeezing gives
$\mathbb E N_X^{\rm phys}(t)/c^2\to\kappa_0H(t)$.

\emph{The rare-gap law.}
For a fixed physical threshold $t$, the waists have a fixed bound.
Let $Z_X^{\rm phys}(t)$ be the number of two-cusp pairs of pants of
$X$ whose model deficit $1/4-\beta(s)$ is at least $t$, that is,
whose waist $s$ satisfies $2(\cosh(s/2)-1)\le H(t)$.
The one-pants integral gives
$\mathbb E Z_X^{\rm phys}(t)=\alpha_g\kappa_0H(t)(1+o(1))$;
the disjoint-pair integral is $O(c^4)$. At every fixed radius,
\eqref{eq:edge-conditional-bad} is uniform for $N\le\sqrt g$,
so expected discarded counts are $o(c^2)$ and good-pair counts
are $O(c^4)$. Fix $\eta>0$ with $t+\eta<1/4$ and choose the
radius so large that~\eqref{eq:edge-form-cutoff} is less than
$\eta/2$ for every model of deficit at least $t+\eta$.
The constant-projection error is $O(1/g)$. Every good such model
therefore implies $D(X)\ge t$ for large $g$. For its good count $Z$,
\[
 \mathbb P(D\ge t)\ge\mathbb P(Z\ge1)
       \ge\mathbb E Z-\tfrac12\mathbb E Z(Z-1).
\]
Divide by $c^2$, take the lower limit and let $\eta\downarrow0$.
Together with $\mathbb P(D\ge t)\le\mathbb E N_X^{\rm phys}(t)$,
this proves~\eqref{eq:edge-rare-law}. The identity
$\mathbb E(N-1)_+=\mathbb E N-\mathbb P(N\ge1)=o(c^2)$
proves the conditional single-mode assertion.

\emph{Upper bounds for the moments.}
For the moments, assume~\eqref{eq:moving-edge-range}. Fix a smooth
$0\le\chi\le1$, zero below $1/2$ and one above $1$.
For fixed $p\ge0$ and $c^2\le b\le1/4$, the tests
$d^p\chi(d/b)$, extended by zero near $d=0$, have $r$th derivative
norm at most $C_{p,r}c^{-2r}$. This includes noninteger $p$.
Lemma~\ref{rev:lem:edge-trace-approximation}, that is
\eqref{eq:edge-trace-error}, gives, uniformly in $b$,
\begin{equation}\label{eq:edge-uniform-moment-trace}
 \mathbb E\sum_{0<\lambda_j<1/4}
      (1/4-\lambda_j)^p\chi((1/4-\lambda_j)/b)
 =\alpha_g\kappa_0\int h_*(y)^p\chi(h_*(y)/b)\,dy+O(r_g),
\end{equation}
where $r_g=g^{-\rho+o(1)}$ for some $\rho>0$ depending on the
fixed approximation choices. The negative polynomial exponents
include $-1/2+3D_0\sigma$ and $-1/2$; all test seminorms are
subpowers of $g$. Weight removal costs at most an exponential in
$-cQ$ times a bounded weighted mean. Since $c^{-1}=o(\log g)$,
$r_g=o(c^A)$ for every fixed $A>0$.
In particular the case $p=0$ yields
\begin{equation}\label{eq:edge-uniform-tail}
 \mathbb P(D\ge b)\le\mathbb E N_X^{\rm phys}(b)
       \le Cc^2/b+O(r_g),\qquad c^2\le b\le1/4.
\end{equation}

For $0<p<1$,~\eqref{eq:edge-uniform-tail} and layer cake imply
\[
 \limsup_g\mathbb E\left[(D/c^2)^p
           \mathbf1_{\{D/c^2>K\}}\right]\le C_pK^{p-1}.
\]
The additive error contributes at most $C_pr_gc^{-2p}\to0$.
Uniform integrability and~\eqref{eq:edge-gap-law} now give
\eqref{eq:edge-small-moments}, since $Z=\kappa_0/E$ for a rate-one
exponential $E$.

For $p\ge1$, the pointwise inequality
\[
 D^p\le c^{2p}+\sum_{0<\lambda_j<1/4}
       (1/4-\lambda_j)^p\chi((1/4-\lambda_j)/c^2)
\]
gives upper bounds. At $p=1$, integrability of
$h_*-(1+y)^{-1}$ and $H(b)\sim1/b$ show
\begin{equation}\label{eq:edge-truncated-model-mean}
 \int h_*(y)\chi(h_*(y)/b)\,dy=\log(1/b)+O(1).
\end{equation}
The transition region lies between $H(b)$ and $H(b/2)$ and costs
$O(1)$. Replacing $\alpha_g$ by $c^2$ costs
$O(c^2\log(1/c)/n)=o(c^2)$. This proves the upper mean bound.
For $p>1$, $h_*^p$ is integrable and dominated convergence in
\eqref{eq:edge-uniform-moment-trace} proves the required upper limit.

\emph{Lower bounds for the moments.}
The lower bounds require expected discarded counts, rather than
an exceptional event of unspecified probability. Put $b=c^2$ and
use good pants relative to all waists of model deficit at least $b$,
with the isolation radius $R'=R+\ell+1$ of~\eqref{rev:eq:edge-good-radius}
built from the common radius~\eqref{eq:edge-common-radius}.
For every $t\ge b$ let $Z^{\rm good}(t)$ count these
good pants of deficit at least $t$, and let $D^{\rm good}_g$ be their
largest deficit, or zero if there are none. Uniformly for $b\le t<1/4$,
\begin{align*}
 \mathbb E Z_X^{\rm phys}(t)
     &=\alpha_g\kappa_0H(t)(1+O(v_g)),\\
 \mathbb E(Z_X^{\rm phys}(t)-Z^{\rm good}(t))
     &\le\epsilon_{g,1}\mathbb E Z_X^{\rm phys}(t),\\
 \mathbb E[Z^{\rm good}(t)(Z^{\rm good}(t)-1)]
     &\le Cc^4/t^2,
\end{align*}
with $v_g,\epsilon_{g,1}\le g^{-1/2+o(1)}$.
These are respectively~\eqref{eq:edge-boundary-volume},
\eqref{eq:edge-conditional-bad}, and the disjoint-pair integral;
$S_b=O(\log\log g)$ and the radius is $o(\log g)$.
Bonferroni and layer cake give
\begin{equation}\label{eq:edge-moment-lower}
 \mathbb E(D^{\rm good}_g)^p\ge
 (1-O(v_g+\epsilon_{g,1}))\alpha_g\kappa_0
       \int_b^{1/4}pt^{p-1}H(t)\,dt
       -C_pc^4\int_b^{1/4}t^{p-3}\,dt.
\end{equation}
The main integral equals
$\int_{h_*\ge b}(h_*^p-b^p)\,dy$. For $p=1$ it is
$\log(1/b)+O(1)$, the pair error is $O(c^2)$, and the relative
geometric error contributes $o(c^2)$. For $p>1$ the main integral
converges to $\int h_*^p$; the pair errors are respectively
$O(c^{2p})$, $O(c^4\log(1/c))$, and $O(c^4)$ for $p<2$, $p=2$,
and $p>2$, all $o(c^2)$.

Finally~\eqref{eq:edge-form-cutoff} and the area projection bound
give a deterministic $e_g=o(c^2)$ such that, on every surface,
$D(X)\ge(D^{\rm good}_g-e_g)_+$. If a good pair of pants exists, apply its trial
vector; otherwise the assertion is trivial. Since $x^p$ is Lipschitz
on $[0,1/4]$ for $p\ge1$, this transfers the lower moment bounds
to $D$ without a probability-to-expectation substitution.
Equations~\eqref{eq:edge-mean}--\eqref{eq:edge-large-moments}
follow. The variance uses $(\mathbb E D)^2=o(c^2)$;
tightness of $D/c^2$ and~\eqref{eq:edge-mean} give
$D/\mathbb E D\to0$ in probability, completing
Corollary~\ref{cor:edge-moments}.

\section{Geometric exclusion in the supercritical regime}
\label{app:super-cuts}

This appendix proves the two geometric lemmas used in the proof of
Theorem~\ref{thm:supercritical-spectrum} in
Section~\ref{sec:supercritical-spectrum}. The first,
Lemma~\ref{lem:super-cuts}, shows that, apart from single planar
separating curves, cuts into two connected pieces with small ratio of
length to smaller area are rare; its bound $B_{i,j,r}(v)$ for two
positive-genus sides is reused in Appendix~\ref{app:dense-comparison}.
The second, Lemma~\ref{lem:super-geometry}, combines this with
first-moment bounds for short curves: with probability tending to one,
the short curves are planar separating, their cut tree has branches of
controlled depth and area, and every cut piece has an explicit lower
bound on its Neumann gap. The volume inputs are Mirzakhani integration
with the boundary estimate \eqref{eq:PS4}, the coefficient-ratio
bound \eqref{eq:PC3} and identity \eqref{eq:PC4} from
Appendix~\ref{app:short-curves}, the point-mass and local-limit
estimates of Section~\ref{extracting-fixed-cusp-counts}, and the
all-marking envelope \eqref{eq:gevrey-envelope}; the Neumann gaps come
from Lemma~\ref{lem:relative-geodesic} of
Appendix~\ref{app:super-comparison}.
Throughout, $n/\sqrt g\to\infty$, $n\le Mg$, $P=2\pi^2$,
$(2g-1)\Psi(x)=n-1$ and $Q=nx/(4P)$. Thus $x\le\bar x(M)<x_c$
and $n/x\le C_Mg$. Constants may depend on $M$.

\begin{lemma}[Competing connected cuts]\label{lem:super-cuts}
For each $M>0$ there are $u_0=u_0(M)>0$, $C_M$ and integers
$g_0=g_0(M)$, $n_0=n_0(M)\ge2$ such that, for
$g\ge g_0$, $n_0\le n\le Mg$ and $0<u\le u_0$, the
$\mathbb P^{\mathrm{WP}}_{g,n}$-probability that $X$ carries a geodesic
multicurve which divides $X$ into two connected components, has length
at most $u$ times the smaller of the two areas, and is not a single
planar separating curve, is at most $C_M(u^4+u^2/g)$. The number of
components of the multicurve and the distribution of the cusps between
the two sides are arbitrary.
\end{lemma}

\begin{proof}
The lower thresholds on $g$ and $n$ ensure the central coefficient
bound used below. Indeed, Proposition~\ref{prop:local-limit} applies
along every sequence with $h\to\infty$ and $\mu\to\infty$; at the
integer saddle $h\Psi(x)=m$, it therefore gives the central lower
bound uniformly once both $h$ and $m$ exceed fixed thresholds.
Otherwise a sequence violating that bound with $h,m\to\infty$
would contradict the proposition. In the planar terms below with
$r\le g/2$, we have $G\ge g/2$ and $N=n+2r-2\ge n$, so the same thresholds apply
simultaneously to $h=2G-1$ and $m=N-1$.

We use Mirzakhani integration with the boundary bound \eqref{eq:PS4}.
For $r$ unlabeled cut curves, the unordered count equals $1/r!$ times
the labeled-boundary integral, in which the $r$ curves are ordered,
and
$\int_{\sum\ell_i\le L}\prod\ell_i\,d\ell_i=L^{2r}/(2r)!$;
see also~\cite[equations~(30)--(31)]{SW}.

\emph{A planar side with $r\ge2$ boundaries.}
If this side has $a$ cusps, put $\chi=a+r-2>0$, minus its Euler
characteristic, and $G=g-r+1$. The planar side has area $2\pi\chi$,
so a cut of ratio at most $u$ has total boundary length at most
$u$ times the smaller area, hence at most $2\pi u\chi$, whether or not
the planar side is the smaller side.
Its expected count with total boundary length at most $2\pi u\chi$ is
bounded by
\begin{equation}\label{eq:super-planar-integral}
 \binom na\frac{V_{0,a+r}V_{G,n-a+r}}{V_{g,n}}
       e^{2\pi u\chi}\frac{(2\pi u\chi)^{2r}}{r!(2r)!}.
\end{equation}
For $r\le g/2$, set $N=n+2r-2$ and choose
$(2G-1)\Psi(x_r)=N-1$. Then $x_r\le x_*(M)<x_c$.
Successive genus exchanges and Theorem~\ref{thm:volume} give
$V_{G,N}/V_{g,n}\le C^{r-1}$. By \eqref{eq:TL50} and \eqref{eq:PC3},
\[
 \binom na\frac{V_{G,n-a+r}}{V_{G,N}}
 \le\frac C{a!}\left(\frac{x_r}{P}\right)^{a+r-2}
       \frac{(n-a+1)^{\overline{r-1}}}{(n+1)^{\overline{2r-3}}},
\]
where the bars denote rising factorials. Pairing numerator factors
with the first $r-1$ denominator factors bounds their quotient by
$(n+r)^{-(r-2)}$. Let
\[
 F_r(z)=\Phi_0^{(r)}(z)
   =\sum_{a\ge\max(0,3-r)}\frac{V_{0,a+r}z^a}{a!P^{a+r-3}}.
\]
The volume prefactor in \eqref{eq:super-planar-integral} is therefore
bounded by $[z^a]F_r(z)x_r^a$ times
$C^rP^{-1}(x_r/(n+r))^{r-2}$.
Choose fixed $s,u_0>0$ and $z_1<x_c$ with
$x_*e^{2\pi u_0+s}<z_1$. The bound
$\chi^{2r}\le(2r)!s^{-2r}e^{s\chi}$ sums the cusp factor into
$e^{(2\pi u+s)(r-2)}F_r(x_re^{2\pi u+s})$.
Cauchy's estimate on a fixed neighborhood of $[0,z_1]$ gives
$F_r\le C r!d^{-r}$. The factorials cancel those in the integral,
leaving
\[
 \mathbb E\{\text{such cuts with }r\text{ boundaries}\}
       \le C\frac{(C_Mu^2)^r}{(n+r)^{r-2}}.
\]
Summing $2\le r\le g/2$ gives $O_M(u^4)$.
For $r>g/2$, the total boundary length is at most $u$ times half the
area, $\pi u(2g-2+n)<\pi(M+2)ug$, so the average boundary length is at most
$2\pi(M+2)u$. At least $r/2>g/4$ boundaries have length at most
$4\pi(M+2)u$. Each is nonseparating in $X$, since the other
$r-1$ curves still join the two sides. Integration and the bounded
genus ratio give $\mathbb E N_{\mathrm{nonsep}}(L)\le Ce^LL^2$.
Markov's inequality bounds this event by $C_Mu^2/g$.

\emph{Two positive-genus sides.} Write their genera as $i,j\ge1$,
with $i+j=g-r+1$, and set
$h_1=2i-2+r\le h_2=2j-2+r$, $H=h_1+h_2=2g-2$.
Even if the first side is not smaller in area, a cut of ratio at most
$u$ has length at most $2\pi u(h_1+a)$, where $a$ is its cusp count.
The same integral and the exponential moment bound reduce its expected
count to
\[
 \frac{(2\pi u/s)^{2r}}{r!}e^{vh_1}B_{i,j,r}(v),\qquad
 B_{i,j,r}(v)=\sum_a\binom na
       \frac{V_{i,a+r}V_{j,n-a+r}}{V_{g,n}}e^{va},\quad v=2\pi u+s.
\]
Put $F_i=\Phi_i^{(r)}$, $F_j=\Phi_j^{(r)}$ and $F_g=\Phi_g'$.
The exact coefficient dictionary gives
\[
 B_{i,j,r}(v)=\frac n{P^r}
     \frac{[z^n]F_i(e^vz)F_j(z)}{[z^{n-1}]F_g(z)}
 \le\frac{Cn}{P^rx}\frac{F_i(xe^v)F_j(x)}{F_g(x)}.
\]
For the inequality use the positive convolution factor of $F_j$,
whose power $h_2\ge g-1$ gives a maximal point mass $C/\sigma$ by
Lemma~\ref{lem:point-mass}; convolution with the tilted first component
cannot increase it. The parent's point mass at $n-1$ is at least
$c/\sigma$ by Proposition~\ref{prop:local-limit}. Only the unshifted
large component uses the point-mass lemma, so the smaller radius of
$F_i(e^vz)$ causes no difficulty.

Write $\rho(z)=r(Y(z))$ and $S(z)=a(Y(z))\rho(z)^{3/2}A(Y(z))$.
The actual all-marking envelope \eqref{eq:gevrey-envelope}, with
\eqref{eq:TL9}, gives for every $k,p\ge1$, $h=2k-2+p$,
\[
 \Phi_k^{(p)}(z)\le C4^{-k}\frac{\Gamma(h)}{\sqrt h}
                    S(z)^{-h}\rho(z)^{p/2}.
\]
The parent has a matching lower bound at $p=1$ by
Proposition~\ref{prop:prefactor}. Bounded logarithmic derivatives on
a compact interval containing $x$ and $xe^v$ imply
\[
 e^{vh_1}\frac{F_i(xe^v)F_j(x)}{F_g(x)}
 \le C_M^re^{C_Mvh_1}\frac{\Gamma(h_1)\Gamma(h_2)}{\Gamma(H+1)}
                   \sqrt{\frac{H+1}{h_1h_2}}.
\]
The powers of $4$ leave $4^{r-1}$ and those of $S$ leave
$S(x)(S(x)/S(xe^v))^{h_1}$: the tilt costs $h_1$, not $H$.
For $h_1\le H/2$,
\[
 \frac{\Gamma(h_1)\Gamma(h_2)}{\Gamma(H)}
  =\frac1{(H-1)\binom{H-2}{h_1-1}}
  \le\frac{2^{-(h_1-1)}}{H-1}.
\]
Choose $s,u_0$ so that $C_M(s+2\pi u_0)<(\log2)/2$.
Using $n/x\le C_Mg$ yields
$e^{vh_1}B_{i,j,r}(v)\le C_M^rg^{-1}e^{-ch_1}/\sqrt{h_1}$.
Summing first the admissible genera and then all $r\ge1$ gives
\[
 Cg^{-1}\sum_{r\ge1}(C_Mu^2)^r/r!=O_M(u^2/g).
\]
These two cases exhaust the non-single-planar connected cuts.
\end{proof}

\begin{lemma}[Short-cap geometry]\label{lem:super-geometry}
With the cutoffs \eqref{eq:super-cutoffs}, with probability tending to
one all $\delta$-short curves are planar separating; their cut tree has
root-branch depth at most $K$ and branch area at most $2\pi K$.
Every cut piece has Neumann gap at least \eqref{eq:super-piece-gap}.
For any growing $R$ with $R\delta\to0$, every $\delta$-short curve
satisfying $\sqrt Q\ell/(k-1)\le R$ is an isolated root leaf with
probability tending to one. Moreover the geodesic Cheeger minimum,
on bounded $Q^{-1/2}$ windows, is attained by single planar waists
of bounded cusp type up to a tail tending to zero.
\end{lemma}

\begin{proof}
The proof has four stages: planarity of the short curves together
with the depth and area bounds; isolation of the retained leaves;
the Neumann gaps of the cut pieces; and the Cheeger minimum.

\emph{Planarity, depth and area.}
The nonplanar first-moment bounds in the proof of
Theorem~\ref{thm:pants-poisson} give exceptional probability $O(\delta^2)$.
The collar lemma makes the remaining curves disjoint. Their planar
sides are disjoint or nested, so the dual graph is a tree.
By \eqref{eq:PC3}--\eqref{eq:PC4} and \eqref{eq:PS4}, for every $k,L$,
\begin{equation}\label{eq:super-one-waist}
 \mathbb E N_{k,L}\le C Q L^2w_k(x)e^L.
\end{equation}
On a compact cusp-weight interval, the $w_k$ have a uniform exponential
tail. Hence the probability of a $\delta$-short cap with more than
$K=C_M\log(Q+2)$ cusps is at most $CQ\delta^2e^{-c_MK}=o(1)$.
A maximal cap with $k$ cusps has total Euler complexity $k-1$;
each planar vertex has positive integer complexity. This bounds its
number of vertices, its depth and its area as asserted.

\emph{Isolation.}
Consider a nested pair with $a\ge1$ cusps between the
curves and $b\ge2$ inside. The coefficient dictionary bounds its
zero-boundary prefactor by
\[
 \frac{Cn}{P^2x}
   \frac{V_{0,a+2}x^a}{a!P^{a-1}}
   \frac{V_{0,b+1}x^b}{b!P^{b-2}}.
\]
The root coefficient deficit is $a+b-1$ and the remaining factorial
factor is exactly $n$, including $a+b=n$.
If either curve has score at most $R$, its squared length integral is
bounded by $R^2(b-1)^2/(2Q)$ or $R^2(a+b-1)^2/(2Q)$; the other
integral is at most $\delta^2/2$. Actual boundary lengths remain
bounded by $\delta$, so their volume error is at most $e^{2\delta}$.
Writing $F=\Phi_0'$, $Y=\Phi_0''$ and $D=x\partial_x$, the total is
\[
 \frac{CnR^2\delta^2}{P^2xQ}
       \{Y(D-1)^2F+(D-1)^2(YF)\}
 \le C_M\frac QnR^2\delta^2=o(1).
\]
The bracket is analytic and $O(x^3)$ on the compact interval, and
$Q/n$ is bounded. Absence of such pairs makes each relevant edge a
root leaf, with no short curve inside its cap.

\emph{Neumann gaps.}
Let $Y_0$ be the root obtained by deleting maximal caps. An internal
geodesic cut of $Y_0$ fills to a connected two-sided cut of $X$ with
the same length and with both areas increased. Lemma~\ref{lem:super-cuts}
excludes all competing topologies at ratio $u$ with error
$C_M(u^4+u^2/g)$. A remaining single planar waist longer than $\delta$
and of ratio at most $u$ must have $k-1>\delta/(2\pi u)$.
Summing \eqref{eq:super-one-waist} with $L=2\pi u(k-1)$ gives error
$C_MQu^2e^{-c_M\delta/u}$, after absorbing the boundary exponential
in $w_k(xe^{2\pi u})$. Choose $u=\delta/[D_M\log(Q+2)]$ with
$D_M$ large. Outside events of probability tending to zero, the root
has internal geodesic Cheeger constant at least $u$: otherwise the
filled cut would have to be a $\delta$-short waist internal to the root.
Lemma~\ref{lem:relative-geodesic} gives its Neumann gap at least $u^2/144$.
Each other piece has area at most $2\pi K$ and no internal geodesic
of length at most $\delta$, so its internal geodesic Cheeger constant
is at least $\delta/(2\pi K)$. The same lemma gives the required gap.

\emph{The Cheeger minimum.}
Finally use $u=t/\sqrt Q$. Lemma~\ref{lem:super-cuts} excludes all
non-single-planar cuts with this ratio. That lemma concerns cuts with
two connected sides; this is no restriction, since by the
graph-averaging reduction in the proof of
Lemma~\ref{lem:relative-geodesic} every separating geodesic multicurve
of ratio at most $u$ contains a subfamily cutting $X$ into two
connected sides with ratio at most $u$. For a truncation level $K_0$,
unrelated to the cutoff $K$, the tail over $k>K_0$ is at
most $C_t\sum_{k>K_0}(k-1)^2w_k(xe^{2\pi t/\sqrt Q})\to0$
uniformly as $K_0\to\infty$. For bounded types the planar side is
the smaller side, of area $2\pi(k-1)$. This proves the last assertion.
\end{proof}

\section{Uniform relative and spectral comparisons}
\label{app:super-comparison}

This appendix proves the two deterministic comparisons used in the
proof of Theorem~\ref{thm:supercritical-spectrum}; both are also used
in Appendix~\ref{app:dense-comparison}. The first,
Lemma~\ref{lem:relative-geodesic}, bounds the relative Cheeger
constant of a surface with short geodesic boundary from below in terms
of its internal geodesic cuts, and hence gives a Neumann gap for each
cut piece; it is applied in the proof of
Lemma~\ref{lem:super-geometry}, and its graph-averaging step is
invoked in Section~\ref{sec:supercritical-spectrum}. The second,
Lemma~\ref{lem:super-comparison}, compares the low eigenvalue counts
of a surface, cut along short separating geodesics into a tree of
pieces, from above and below with those of a weighted star. No
probabilistic input is used: the tools are the relative Cheeger
inequality of He--Wu~\cite{HeWu}, the collar lemma~\cite{BuserBook},
isotopy of multicurves~\cite{FarbMargalit}, and Poincar\'e, trace and
collar-mass estimates with min--max.

\begin{lemma}[Relative geodesic comparison]\label{lem:relative-geodesic}
Let $Y$ be a connected finite-area hyperbolic surface whose geodesic
boundary circles have lengths at most $\delta$; cusps are allowed.
Let $h_N$ count only the interior perimeter of a relative cut, and let
$H_{\mathrm{int}}$ be the infimum of length divided by smaller area
over internal geodesic multicurves with two connected complementary
components, excluding boundary-parallel classes. An empty infimum is
infinity. Set $w=\operatorname{arcsinh}(1/\sinh(\delta/2))$ and
$a_\delta=1+\delta/(2w)$. Then
\begin{equation}\label{eq:relative-geodesic}
 h_N\ge\frac{H_{\mathrm{int}}}
                   {a_\delta+(a_\delta+2)H_{\mathrm{int}}}.
\end{equation}
For small $\delta$, if $H_{\mathrm{int}}\ge u$ and $0<u\le1$, the
Neumann gap is at least $u^2/144$, independently of the topology.
\end{lemma}

\begin{proof}
The proof has five steps: the relative Cheeger inequality in the
presence of cusps; surgery of an interface into a closed internal one;
straightening it to geodesics with controlled area change; the isotopy
justification of that step; and the reduction to cuts with connected
sides, which gives \eqref{eq:relative-geodesic}.

\emph{Step 1: the relative Cheeger inequality with cusps.}
The relative Cheeger inequality is $\operatorname{gap}_N\ge h_N^2/4$;
see~\cite[Theorem~9]{HeWu}. That theorem is stated for compact
hyperbolic surfaces with geodesic boundary. For finite-area $Y$ with
cusps the same coarea argument applies verbatim to every
$f\in H^1(Y)$ that vanishes above some horocycle in each cusp: the
level sets of such $f$ are compact relative cuts of $Y$, to which the
defining inequality of $h_N$ applies. These $f$ are dense in $H^1(Y)$,
because truncating at the horocycle of length $\epsilon$ with a cutoff
of gradient at most one changes $\|f\|_{H^1}$ by $o(1)$ as
$\epsilon\to0$, the cusp tail having area $\epsilon$; since the
Rayleigh quotient is continuous on $H^1(Y)$, the gap bound extends.
We prove \eqref{eq:relative-geodesic}.

\emph{Step 2: from a relative interface to a closed internal one.}
Start with a smooth finite interface of interior length $L$, separating
areas with minimum $m$. Such interfaces suffice by smooth approximation
and truncation of cusp tails, whose added horocycle lengths and areas
tend to zero. Remove contractible loops and boundary-parallel proper
arcs by filling or emptying innermost disks. For a boundary-parallel
arc the corresponding geodesic boundary subarc is no longer than the
arc; the disk has perimeter at most twice the arc length and area at
most that perimeter. Each removed original component is charged once.
Thus the total area change is at most $2L$ and the remaining interface
length is at most $L$. The isoperimetric inequality used here is
\cite[Lemma~10]{HeWu}, and the factor two for boundary-parallel arcs is
the one in the proof of \cite[Lemma~12]{HeWu}.

Two standard facts are used in the rest of this step and in the next
one. First, since each
boundary circle has length at most $\delta$, the half-collars of width
$w=\operatorname{arcsinh}(1/\sinh(\delta/2))$ about the boundary
circles are embedded and pairwise disjoint: this is the collar lemma
\cite[Chapter~4]{BuserBook}, applied to the double of $Y$ across its
geodesic boundary, whose collars about the doubled circles are
embedded and disjoint and restrict to the half-collars of $Y$.
Second, the closed interface constructed below may contain components
homotopic to a boundary circle; such components are not internal
geodesic classes, and in the straightening step they are pushed into
the collar and disappear from the interior interface, with a change
of the side areas bounded by their length, by the same Fermi-coordinate
primitive as for the other classes.
Every remaining proper arc exits the embedded half-collars at both
ends, even if both ends lie on the same circle. It has length at least
$2w$. Hence at most $L/(2w)$ arcs, and at most $L/w$ boundary circles,
are involved. On each such circle assign a sufficiently thin strip
wholly to one side, choosing the shorter of the two traces. The added
interface length tends to at most half the circle length, and its area
tends to zero. After smoothing, this gives a closed internal interface
with length $L'\le a_\delta L+o(1)$ and smaller area
$m'\ge m-2L-o(1)$.

\emph{Step 3: straightening.}
Straightening this closed interface gives a geodesic partition with
\begin{equation}\label{eq:relative-area}
 L_{\mathrm{geo}}\le L',\qquad m_{\mathrm{geo}}\ge m'-L'.
\end{equation}
Here is the area justification, including curves crossing their own
geodesic representatives. Double across the boundary if necessary.
In the cyclic cover of an essential curve, the area form in Fermi
coordinates is $\ell\cosh r\,dr\,d\theta$, with primitive
$\ell\sinh r\,d\theta$ of norm at most one, vanishing on the axis.
The trace annulus of a homotopy lifts to this cover. Stokes bounds its
signed area by the original curve length plus the absolute primitive
integral on the final curve. On a parallel curve approaching the axis,
the latter tends to zero. The same construction treats classes parallel
to the original boundary: the parallel curve is pushed onto the
boundary circle, where it ceases to be part of the interior interface,
and the swept area is at most the curve length, which is already
counted in $L'$. For cusp classes use $dx/y$ in a parabolic
cover, whose norm is one and whose integral on a shrinking horocycle
tends to zero. For contractible classes use
$(\cosh r-1)d\theta$ in polar coordinates on the hyperbolic plane;
its norm is $\tanh(r/2)\le1$.

\emph{Step 4: the isotopy justification.}
The finite disjoint family can be simultaneously isotoped to disjoint
nearby representatives, keeping parallel copies in their original
order, pushing peripheral curves into collars or cusp tails, and
shrinking nested contractible families in disks. The standard
topological inputs are minimal position, homotopy versus isotopy and
ambient extension~\cite[Corollary~1.9 and Propositions~1.10--1.11]{FarbMargalit}.
At each finite stage the ambient isotopy has compact support and fixes
the boundary of $Y$. If $D'$ is one side and $D'_\epsilon$ its image,
area transport expresses
$\operatorname{area}(D'_\epsilon)-\operatorname{area}(D')$ as the sum
of the oriented trace-annulus integrals. Its absolute value is at most
$L'+o(1)$ by the primitives just described. As the neighborhoods shrink,
the side converges in area to a union of geodesic complementary
components. Parallel interfaces cancel in pairs because their side
labels alternate; an odd family leaves one geodesic. Each retained
geodesic minimizes length in its class. Since
$a\mapsto\min(a,\operatorname{area}(Y)-a)$ is 1-Lipschitz, this proves
\eqref{eq:relative-area}. No estimate for the entire swept region or
for the symmetric difference of the sides has been used.

\emph{Step 5: reduction to connected sides.}
The resulting sides need not be connected. Regard the complementary
pieces as vertices of a finite graph weighted by area, with interface
lengths as edge weights. A connected component of the smaller side
has no greater boundary-to-mass ratio. If its complement is disconnected
and one complementary component has more than half the total mass,
take the complement of that component: it is connected, has greater
mass and no larger boundary. Otherwise averaging over complementary
components selects one of ratio no greater than the previous ratio;
its complement is connected through the original chosen component.
Thus some cut with both sides connected has no larger ratio, using
only the existing geodesics. It follows that
$L_{\mathrm{geo}}\ge H_0m_{\mathrm{geo}}$ for every finite
$H_0\le H_{\mathrm{int}}$, also when the interface is empty.
Consequently $L'\ge H_0(m-2L-L'-o(1))$. Let the strip widths tend to
zero and then let $H_0\uparrow H_{\mathrm{int}}$ to obtain
\eqref{eq:relative-geodesic}. For $a_\delta\le2$,
$h_N\ge u/(2+4u)\ge u/6$, proving the gap bound.
\end{proof}

The second lemma is the surface--graph comparison. Its hypotheses,
the Neumann gaps of the pieces and the depth and area bounds on the
branches, are supplied in the proof of
Theorem~\ref{thm:supercritical-spectrum} by
Lemma~\ref{lem:super-geometry}; the proof below is a min--max argument
with explicit trial maps in both directions, and its intermediate
estimates \eqref{eq:tree-trace}--\eqref{eq:tree-mass} are the ones
reused in Appendix~\ref{app:dense-comparison}.

\begin{lemma}[Comparison with the retained star]\label{lem:super-comparison}
Cut a finite-area surface along separating geodesics of length at most
$\delta$. Suppose their dual graph is a tree with a positive-genus root
and planar other vertices. Every geodesic piece has Neumann gap at
least $0<\gamma\le1$, each root branch has depth at most $K$, and its
area is at most $C K$. Choose $\delta<b$ with $b$ small. Trim each
half-collar at circumference $b$, let $m_v$ be the remaining areas and
$W_e$ their total below edge $e$. Put
\[
 r_e=\operatorname{arcosh}(b/\ell_e),\quad
 I_e=\int_{-r_e}^{r_e}\operatorname{sech}r\,dr,\quad
 \kappa_e=\ell_e/I_e.
\]
Retain edges with $\kappa_e/W_e\le\tau$, and suppose each is an
isolated root leaf. Contract the other vertices into the root, adding
their masses, to form a weighted star. There are absolute constants
$C$ and $T\le CK/(b\gamma)$ such that, for $\Lambda<1/4$,
\begin{equation}\label{eq:tree-counts}
 N_\star^+(\Lambda)\le N_X^+(\Lambda)
       \le N_\star^+(B\Lambda/a),\qquad B=(1+\sqrt{2\tau T})^2,
\end{equation}
where
\begin{equation}\label{eq:tree-norm}
 a=\left[\sqrt{\frac{1-C\Lambda/\gamma}{1+Cb}}
                 -\sqrt{(2K/\tau+4KT)\Lambda}\right]^2,
\end{equation}
provided the first square root and the bracket are positive. Moreover
$m_v=(1+O(b))\operatorname{area}(P_v)$ and
$\kappa_e=(\ell_e/\pi)(1+O(\delta/b))$.
\end{lemma}

\begin{proof}
Write $B_v$ for a trimmed piece, $a_v$ for the mean of a real $H^1$
function $f$ on it, and $p_e,q_e$ for its endpoint angular averages.
The standard disjoint half-collars extend at least one unit beyond
$r_e$ when $b$ is small. Their outer unit strips are therefore disjoint.
Each removed half-collar has area $\sqrt{b^2-\ell_e^2}\le b$.
If $\chi_v=-\chi(P_v)>0$ and $d_v$ is its boundary count, then
$d_v\le3\chi_v$. Hence $m_v=(1+O(b))2\pi\chi_v$ and $d_v\le Cm_v$.

\emph{Poincar\'e and trace estimates.} Extend a function on $B_v$
harmonically inward on every removed half-collar, with a free endpoint
at the waist. In conformal coordinates, Fourier mode $k\ne0$ has
extension energy $2\pi|k|\tanh(2\pi|k|T_{\mathrm{in}})|f_k|^2$.
The outer strip has conformal length $T_{\mathrm{out}}\ge c/b$,
and its energy for this trace is at least
$2\pi|k|\tanh(2\pi|k|T_{\mathrm{out}})|f_k|^2$, at least half the
extension bound. The zero mode extends constantly. Summing disjoint
strips gives extension energy at most three times the original energy.
Variance on $B_v$ is no larger than variance of its extension, so
\[
 \int_{B_v}|f-a_v|^2\le3\gamma^{-1}\int_{B_v}|\nabla f|^2.
\]
The one-dimensional trace inequality and angular averaging on each
outer strip give
$|p_e-a_v|^2\le Cb^{-1}\int_{\mathrm{strip}}(|f-a_v|^2+|\nabla f|^2)$.
Sum over disjoint strips before applying Poincar\'e. With
$E=\int_X|\nabla f|^2$, this yields
\begin{equation}\label{eq:tree-trace}
 \sum_{\text{endpoints}}|p_e-a_v|^2\le\frac C{b\gamma}E,
 \quad
 \sum_eW_e\bigl((p_e-a_{\mathrm{parent}})^2+
                   (q_e-a_{\mathrm{child}})^2\bigr)\le TE.
\end{equation}
Here $W_e\le CK$; no boundary-count factor is introduced.

\emph{Mass in the collars.} On a half-collar, for $z(r_e)=0$,
\[
 \int_0^{r_e}z^2\ell\cosh r\,dr
           \le4\int_0^{r_e}(z')^2\ell\cosh r\,dr.
\]
Indeed the ground state $\phi=(\cosh r)^{-1/2}$ has potential
$-(\cosh r\,\phi')'/(\cosh r\,\phi)=1/4+1/(4\cosh^2r)$;
the boundary term vanishes at $0$ since $\phi'(0)=0$ and at $r_e$
since $z(r_e)=0$. Subtract the outer average from the angular zero
mode. Its mass is at most $2bp_e^2$ plus eight times its radial
energy. Nonzero modes have mass at most $(b/(2\pi))^2$ times their
angular energy. Thus a full collar $C_e$ has mass at most
$2b(p_e^2+q_e^2)+C\int_{C_e}|\nabla f|^2$. Combining this with \eqref{eq:tree-trace}
and $d_v\le Cm_v$, and writing $M(f)=\sum_vm_va_v^2$, gives
\begin{equation}\label{eq:tree-mass}
 \|f\|^2\le(1+Cb)M(f)+C\gamma^{-1}E,\qquad
 \sum_e\kappa_e(p_e-q_e)^2\le E.
\end{equation}
The second inequality follows by radial Cauchy--Schwarz after angular
averaging. Also $I_e=\pi+O(\delta/b)$, proving the conductance estimate.
All extension arguments hold for $H^1$ traces by approximation.

\emph{Contracting branches.} A nonretained edge has $W_e<\kappa_e/\tau$.
Sum Cauchy--Schwarz along paths of length at most $K$, then over vertex
masses. Since retained edges are root leaves, every other path consists
only of nonretained edges. Equations~\eqref{eq:tree-trace}--\eqref{eq:tree-mass}
give
\[
 \sum_{v\text{ contracted}}m_v|a_v-a_0|^2
 \le K\sum_{e\text{ not retained}}W_e|a_{\mathrm{child}}-a_{\mathrm{parent}}|^2
 \le(2K/\tau+4KT)E.
\]
The left side is the squared distance in the full core-mass norm
between the original means and the vector replacing all contracted
means by $a_0$. On a subspace with $E\le\Lambda\|f\|^2$,
the triangle inequality and \eqref{eq:tree-mass} give
$M_\star(Jf)\ge a\|f\|^2$, with $a$ in \eqref{eq:tree-norm},
where $J$ keeps the root mean and the retained leaf means.
For retained edges $W_e=m_e$ and $\kappa_e\le\tau m_e$.
Minkowski's inequality in the edge-energy norm gives
$Q_\star(Jf)\le(1+\sqrt{2\tau T})^2E$.
The positivity of the bracket in \eqref{eq:tree-norm} ensures that $J$
is injective on the spectral subspace. Min--max proves the upper count
bound in \eqref{eq:tree-counts}.

For the other direction, assign each star value constantly to its core,
the root value to all contracted cores, and use radial harmonic
interpolation on retained collars and constants on the other collars.
This injective trial map has energy exactly $Q_\star$ and mass at
least $M_\star$. Min--max gives the lower count bound. Both comparisons
include the constant zero mode before subtracting one; neither map
needs to preserve orthogonality to constants.
\end{proof}

\section{Uniform cutoffs at a cusp of a large planar surface}
\label{app:cusp-cutoffs}

This appendix constructs low-energy cutoff functions around a typical
cusp of a large WP sphere, with constants independent of the number of
cusps; this is Lemma~\ref{lem:uniform-cusp-cutoff} below, the first of
the two inputs to the proof of Theorem~\ref{thm:dense-collapse} in
Appendix~\ref{app:dense-collapse-details}. The same construction
identifies the conformal type of the area-rooted infinite WP sphere,
which is the recurrent endpoint of Theorem~\ref{thm:wp-conformal-phase}
(Section~\ref{sec:local-conformal-phase}). The strategy is to
regularize the metric conformally to bounded geometry, and then to
apply a finite planar resistance estimate on a triangulation of the
regularized sphere; the resulting discrete potential is interpolated
and smoothed into the cutoff. The inputs are the genus-zero volume
bounds of \cite{HTSphere}, the local convergence of
\cite[Theorem~34]{BC}, the conformal limit theorem of
\cite[Theorem~1.4]{NPS}, the triangulation of \cite{BowditchTriang}
and the resistance estimate of \cite[Lemma~3.6]{GGN}.

We write $E_X(f)=\int_X|\nabla f|^2$ and let $U_b(c)$ be the
complete tail of cusp $c$ bounded by a horocycle of length $b$.
All horocycles below are chosen inside disjoint Margulis regions.

\begin{lemma}[Uniform cusp cutoffs]\label{lem:uniform-cusp-cutoff}
There are $b>0$ and $M<\infty$ such that, for every $e,p>0$, there
are finite $A(e,p),N_0(e,p)$ with the following property. For every
$N\ge N_0(e,p)$, a WP sphere $X_N\in\mathcal M_{0,N}$ with a fixed
labeled cusp $c$ admits, with probability at least $1-p$, a function
$\phi\in H^1(X_N)$ satisfying
\begin{equation}\label{eq:cusp-cutoff}
 0\le\phi\le1,\quad \phi|_{U_b(c)}=1,\quad
 \operatorname{Lip}_{\mathrm{hyp}}\phi\le M,\quad
 E_{X_N}(\phi)\le e,\quad
 \operatorname{area}\{\phi\ne0\}\le A(e,p).
\end{equation}
The support need not be compact in the hyperbolic metric.
More precisely, there are constants $C,u_0$ such that, for every
$N\ge3$ and $u\ge u_0$, with probability at least
$1-Cu^{-1/3}\log u$, the first three conditions in
\eqref{eq:cusp-cutoff} hold together with
\begin{equation}\label{eq:quantitative-cusp-cutoff}
 E_{X_N}(\phi)\le\frac C{\log u},\qquad
 \operatorname{area}\{\phi\ne0\}\le Cu.
\end{equation}
All constants are independent of $N$ and $u$.
\end{lemma}

The proof also establishes the parabolic conformal type conjectured in
\cite[Section~7.3]{BC} for the puncture-completed, area-rooted local
limit. We use the strong local convergence proved there, and derive
the conformal assertion by changing the metric and controlling the
resulting change of root law; the change of metric is needed because
the original cusped metrics have no uniform injectivity-radius bound.
The first subsection constructs the regularized metric and records its
geometry in \eqref{eq:regularized-geometry}, the second checks the
locality of the construction and identifies the conformal type of the
local limit, and the third proves the lemma. The regularized metric is
also used, for surfaces of arbitrary genus, in the proof of
Proposition~\ref{prop:low-energy-tail} in
Appendix~\ref{app:cusp-dominated-orbits}.

\subsection{A conformal regularization}
We replace the hyperbolic metric by a conformal metric on the filled
sphere with bounded curvature and injectivity radius bounded below.
The construction fills every cusp and inflates every short collar by
fixed models, so that all constants are independent of $N$ and of the
surface; its output is the geometry summary
\eqref{eq:regularized-geometry} together with the explicit local form
of the modified metric (the fixed cusp-disk model of diameter $D$ and
central area $v$, the collar inflation $F/f$, and the domination of
the old metric outside the cusp disks); these are the only properties
of the construction that the later subsections use. The only probabilistic input is the first moment
\eqref{eq:inverse-collar-moment} of the inverse collar lengths, which
controls the area added by the inflated collars.
Fix $0<\delta\ll b\ll1$. Let
$Z_\delta(X_N)=\sum_{\ell_\gamma\le\delta}\ell_\gamma^{-1}$,
where the sum is over simple closed geodesics. Mirzakhani integration
and the two-sided genus-zero volume bounds give
\begin{equation}\label{eq:inverse-collar-moment}
 \mathbb E Z_\delta(X_N)\le C\delta N.
\end{equation}
Indeed the zero-boundary factor for $k$ cusps on one side satisfies
\[
 \binom Nk\frac{V_{0,k+1}V_{0,N-k+1}}{V_{0,N}}
 \le\frac{CN^{7/2}}{k^{5/2}(N-k)^{5/2}},\qquad 2\le k\le N-2.
\]
These bounds include the stable small indices; see
\cite[Theorems~2.1--2.2]{HTSphere}. The twist factor $\ell$ cancels
$1/\ell$, the boundary factors are bounded for $\ell\le\delta$,
and summing $CNk^{-5/2}$ for $k\le N/2$ proves the claim.
Counting the two sides separately only enlarges the upper bound.

Fill every cusp conformally and inflate every collar of central
length $\ell<\delta$. Denote the resulting closed sphere by $M(X_N)$.
Below, $M(X_N)$, abbreviated $M$ as in $K_M$, $q_M$ and
$\operatorname{area}(M)$, always denotes this regularized sphere. It is
unrelated to the Lipschitz constant $M$ of
Lemma~\ref{lem:uniform-cusp-cutoff}, which reappears only as
$M(t,\alpha)$ at the end of the proof.
Here is a construction with constants independent of $N$ and $X_N$.
On a geodesic collar use $\theta\in\mathbb R/\mathbb Z$ and
$ds=dr/(\ell\cosh r)$. Its metric is
$f(s)^2(ds^2+d\theta^2)$, where $f(s)=\ell\sec(\ell s)$.
Choose a smooth positive function $H$ with $H(u)=2$ for $u\le1$,
$H(u)=u$ for $u\ge3$, and $H(u)\ge\max(2,u)$, and replace $f$ by
$F=bH(f/b)$. Such an $H$ is obtained by symmetric smoothing of
$\max(2,u)$ on its transition interval. The change is confined to
$f<3b$, with a buffer out to $f=6b$. The central cylinder is flat,
has circumference $2b$, and adds area at most $C/\ell$ because its
conformal length is at most $\pi/\ell$.

On the transition band $f$ is between fixed positive constants and
\[
 (f')^2=f^2(f^2-\ell^2),\qquad f''=2f^3-\ell^2f.
\]
These identities and their derivatives bound the curvature
$-F^{-2}(\log F)''$ and all its covariant derivatives uniformly.
For a cusp use $u\ge0$, $f(u)=b/(1+bu)$ and interpolate $\log f$
smoothly, on a fixed bounded interval, to $\log c_0-2\pi u$.
Thereafter use $c_0e^{-2\pi u}$. In the coordinate
$z=e^{-2\pi(u+i\theta)}$ this is a smooth flat metric at $z=0$.
Keep the old metric near $u=0$ and out to the horocycle of length
$2b$. All filled disks are the same fixed smooth model, of diameter
at most $D$, each containing an embedded central disk of area $v>0$.
The metric outside these cusp disks dominates the old metric.

Choose $b$ so that the buffered regions lie in disjoint Margulis components.
In each modified collar the new arclength coordinate $t$ gives
$dt^2+F(t)^2d\theta^2$, with $2b\le F\le6b$ and bounded
$|\partial_t\log F|$. Fixed small rectangles give embedded uniformly
bilipschitz disk charts, including across the buffered interfaces.
The cusp models have the same property. Away from these regions the
old injectivity radius is bounded below in terms of $\delta,b$.
Thus every sufficiently small new metric ball contains a subchart
of uniformly positive area. The curvature and noncollapse estimate
of Cheeger--Gromov--Taylor~\cite[Theorem~4.7(i), (4.23)]{CGT}
gives a uniform positive injectivity radius. We have proved
\begin{equation}\label{eq:regularized-geometry}
 |K_M|\le C,\quad \operatorname{inj}(M)\ge c>0,\quad
 cN\le\operatorname{area}(M)\le CN+CZ_\delta(X_N),
 \quad \mathbb E\operatorname{area}(M)\le CN.
\end{equation}
For the lower area bound, at most $Nb$ of old cusp area is removed
from $2\pi(N-2)$. All curvature derivatives are uniformly bounded,
so pointed smooth compactness is available.

\subsection{Locality and the change of root measure}
This subsection proves the conformal assertion for the area-rooted
infinite WP sphere, the recurrent endpoint of
Theorem~\ref{thm:wp-conformal-phase}; the proof of
Lemma~\ref{lem:uniform-cusp-cutoff} in the next subsection works on
the finite surfaces and uses \eqref{eq:regularized-geometry} together
with the local description of the regularized metric from the first
subsection, but not the recurrence argument of this subsection.
The same regularization is defined on the infinite local
limit of \cite[Theorem~34]{BC}. We first check its compatibility with
that theorem's total-variation convergence of each fixed-radius metric
ball, then apply the conformal limit theorem of
\cite[Theorem~1.4]{NPS} to the regularized spheres rooted by their new
area, and finally compare the new root law with the old one through the
domination \eqref{eq:root-domination}.
In the fixed filled disk, any two boundary points can be joined
along its boundary with length at most $C_0$ times their intrinsic
distance in the disk. This follows near the diagonal from a smooth
boundary collar and away from it by compactness. Replacing disk
excursions by these arcs, a new path of length $R$ between old-core
points gives an old-core path of length at most $C_0R$.
If the root has old cusp depth $\tau$, first move it to its cusp boundary.
A new $R$-ball is consequently determined by an old ball of radius
$\tau+C_1(R+D+1)$, together with the attached fixed disks.

The data identifying the modifications are local as well. A short
primitive based loop in a Margulis region has developing holonomy
which distinguishes a parabolic component from a hyperbolic one.
In the latter case its translation length and distance to its axis
determine $\ell$ and $f$; the central geodesic need not lie in the
observed ball. In the parabolic case the holonomy determines the
circumference coordinate. The loop and its development fit in a
fixed additional neighborhood. Exact local isometries therefore
preserve these data, including the threshold $\ell<\delta$.
The filled disk is rotationally invariant, so no gluing phase is
needed. There is no assertion that the operation is continuous at
this sharp threshold for arbitrary weak metric convergence.

The modified infinite surface is complete. A bounded-length escaping
path would, by the path comparison, remain in a compact part of the
complete old truncated core. Only finitely many cusp disks can meet
such a compact set, because their disjoint fixed-width buffers have
positive area. Each such disk is compact. This excludes escape.
The old area fraction at cusp depth greater than $\tau$ is at most
$Ce^{-\tau}$. Truncate this event and apply the exact local couplings
of \cite[Theorem~34]{BC}. It follows that $M(X_N)$ rooted by old
hyperbolic area converges locally to the modified area-rooted limit.

To identify its conformal type, let $q_M$ be the law of $M$ rooted
uniformly in its new area. Bias the law of $X_N$ by
$\operatorname{area}(M)/\mathbb E\operatorname{area}(M)$, and let
$\mu_N$ be the barycenter of the resulting laws $q_M$.
Namazi--Pankka--Souto~\cite[Theorem~1.4]{NPS}, after a fixed metric
rescaling, applies to any deterministic sequence of these closed
spheres whose areas diverge: their distributional limits are supported
on surfaces conformally $\mathbb C$ or $\mathbb C^*$.
It also applies to our mixtures. To see this, let $\mathcal L$ be the
compact set of limit laws $q_M$ along such deterministic sequences.
Pointed bounded-geometry compactness implies that uniformly over
spheres of area at least $cN$, the distance from $q_M$ to
$\mathcal L$ tends to zero; otherwise a violating sequence would
have a limit in $\mathcal L$. Every limiting distribution of the
random law $q_M$ is therefore supported on $\mathcal L$. Taking
barycenters proves that every limit $\mu$ of $\mu_N$ is parabolic.

Let $\nu_N$ instead root $M(X_N)$ by old hyperbolic area, whose total
mass is the deterministic $2\pi(N-2)$. Outside old cusp depth $\tau$,
the density $dA_{\mathrm{old}}/dA_{\mathrm{new}}$ is at most $C_\tau$.
Equation~\eqref{eq:regularized-geometry} consequently gives, for
every bounded measurable $0\le F\le1$ on pointed surfaces,
\begin{equation}\label{eq:root-domination}
 \int F\,d\nu_N\le C'_\tau\int F\,d\mu_N+Ce^{-\tau}.
\end{equation}
Pass to simultaneous limits first for continuous tests, then for
Borel sets by measure regularity. Sending $\tau\to\infty$ shows
$\nu\ll\mu$. Thus the modified Budd--Curien limit is parabolic.
Its underlying surface before filling is homeomorphic to
$\mathbb R^2\setminus\mathbb Z^2$ by the proof of
\cite[Lemma~33]{BC}; after
filling it is simply connected and hence conformally $\mathbb C$.
The unfilled surface is therefore $\mathbb C$ minus a closed discrete
set, which is polar: planar Brownian motion started off this set
almost surely never hits it and visits every open subset at arbitrarily
large times. Brownian motion for the complete hyperbolic metric is a
time change of this process, and its clock does not explode because
the curvature is bounded below; it is therefore recurrent.
No area concentration or second inverse-length moment is used here.

\subsection{Finite planar resistance and smoothing}
We return to the finite surfaces and prove
Lemma~\ref{lem:uniform-cusp-cutoff}. A bounded-degree planar
triangulation of the regularized sphere has, at all but a small
proportion of its vertices, a vertex set of controlled size such that
the effective resistance from the vertex to the complement of the set
is large \cite[Lemma~3.6]{GGN}. The corresponding discrete potential,
interpolated and then smoothed by the heat semigroup, is the required
cutoff. The proof has five steps: the triangulation, the resistance
estimate and the count of bad cusps, the plateau on a fixed graph ball,
the interpolated function and its support, and the smoothing.

\begin{proof}[Proof of Lemma~\ref{lem:uniform-cusp-cutoff}]
We prove the quantitative assertion directly on the finite surface.

\emph{Step 1: the triangulation.}
The curvature and injectivity bounds in
\eqref{eq:regularized-geometry} imply a uniform positive convexity
radius, by the Hessian comparison for squared distance on sufficiently
small normal balls. Bowditch~\cite[Theorem~1.1]{BowditchTriang}
therefore gives a uniformly bilipschitz triangulation of $M(X_N)$ by
standard equilateral triangles. Its one-skeleton $G=(V,E)$ is planar.
Each triangle has area between fixed positive constants and uniformly
bounded diameter. Triangles incident to a vertex lie in a fixed-radius
ball and have disjoint interiors. Bishop comparison bounds their number,
so $G$ has uniformly bounded degree and $|V|\le C\operatorname{area}(M)$.

\emph{Step 2: the resistance estimate and the count of bad cusps.}
Assign each filled cusp center to a vertex of a triangle containing it.
Each vertex receives at most a fixed number of cusps: their disjoint
central disks of area $v$ lie in a fixed-radius ball about that vertex.
The same packing argument bounds the number of cusp disks meeting any
one triangle. Gurel-Gurevich--Nachmias~\cite[Lemma~3.6]{GGN} shows that,
for all but at most $Cu^{-1/3}\log u\,|V|$ vertices $\rho$, there is
$B\subset V$ with
\[
 |B|\le Cu,\qquad R_{\mathrm{eff}}(\rho,V\setminus B)\ge c\log u.
\]
That lemma concerns finite planar graphs of maximal degree at most
$D$, and its constant $c=c(D)>0$ depends only on $D$. By the Bowditch
and Bishop comparisons above, the maximal degree of $G$ is bounded by
a constant depending only on the fixed bounds in
\eqref{eq:regularized-geometry}, independently of $N$ and $X_N$; hence
$c$ and $C$ here are uniform constants in the same sense.
The network has unit edge conductances; when $B=V$ the resistance
is interpreted as infinity. The preceding multiplicity bound and
\eqref{eq:regularized-geometry} imply that a uniformly chosen cusp is
assigned to a failing vertex with probability at most
\begin{equation}\label{eq:cusp-root-count}
 Cu^{-1/3}\log u\,\frac{\mathbb E|V|}{N}
 \le C'u^{-1/3}\log u.
\end{equation}
This uses a count of bad vertices, not an identification of the two
root laws. The construction below bounds the intrinsic number of
cusps without a suitable cutoff for every admissible triangulation;
exchangeability then gives the same estimate for each fixed label.

\emph{Step 3: the plateau on a fixed graph ball.}
We next check that one fixed graph ball suffices for the required
plateau. For every fixed $R$, the filled closure of the old
$R$-neighborhood of $U_b(c)$ lies in a new-metric ball of radius $R_*(R)$
about its filled center, uniformly in $N$ and $X_N$. Outside the target
tail, start an old path of length at most $R+1$ on its boundary.
Every visit to a buffered closed collar enters through $f=6b$.
Since $f=\ell\cosh r$ and $|\partial_r\log f|\le1$, throughout such
a visit $f\ge6b e^{-(R+1)}$. Hence the length inflation $F/f$ is
bounded by a constant depending only on $R$. In every cusp disk
the new-to-old length ratio is uniformly bounded above, by the fixed
interpolation and exponential cap. Elsewhere it equals one. The path
thus has uniformly bounded new length, and the target tail becomes
a disk of diameter at most $D$. This proves the assertion despite
the unbounded inflation deep inside very short closed collars.

Fix $t=1$ and $0<\alpha<1/4$. Choose $R_0$ so the time-$t$ hyperbolic
Brownian displacement exceeds $R_0$ with probability less than $\alpha$,
uniformly over all complete hyperbolic surfaces by lifting to $\mathbb H^2$.
There is a fixed integer $m$ such that all vertices of triangles meeting
the new $R_*(R_0+1)$-ball lie within graph distance $m$ of the assigned
vertex $\rho$. Indeed a connecting geodesic meets a connected union of
triangles in a fixed larger ball; their uniform lower area bound and
Bishop comparison bound their number, hence the length of a connecting
path in their one-skeleton.

\emph{Step 4: the interpolated function and its support.}
At a good vertex the Dirichlet principle supplies $0\le h\le1$ with
$h(\rho)=1$, $h=0$ on $V\setminus B$, and
$\sum_{vw\in E}(h(v)-h(w))^2\le C/\log u$. If $B=V$, use $h=1$.
Pathwise Cauchy--Schwarz gives
$|h(v)-1|\le\sqrt{Cm/\log u}$ for $d_G(v,\rho)\le m$.
For one fixed sufficiently large $u_0$, replacing $h$ by $\min(1,2h)$
makes it one on this graph ball, with at most four times the energy
and no additional nonzero vertices. Affine interpolation on the
standard triangles, transported by the bilipschitz triangulation,
gives $0\le f\le1$ with energy at most $C/\log u$, equal to one on
the filled old $(R_0+1)$-neighborhood of $U_b(c)$. The energy bound
follows by summing the affine gradient estimate on each triangle.

Only $O(|B|)=O(u)$ triangles have a nonzero vertex. Their areas and
the number of cusp disks meeting each are uniformly bounded. Old
area is at most new area outside the disks, and each whole old tail
has area $b$. Thus the pullback to $X_N$ has old support area at most
$Cu$. Its Dirichlet energy is unchanged by conformal invariance in
dimension two.
The pullback is bounded, locally Lipschitz on $X_N$, of finite energy
and supported in finite hyperbolic area, so it lies in $H^1(X_N)$,
which is the form domain of the Laplacian because $X_N$ is complete.

\emph{Step 5: smoothing.}
Smoothing gives the Lipschitz bound, with a constant independent of
$N,u$. Let $P_t=e^{t\Delta}$ with the nonpositive Laplacian.
The curvature $-1$ gradient inequality and the variance identity give
\[
 |\nabla P_su|^2\le e^{2s}P_s|\nabla u|^2,\qquad
 P_tf^2-(P_tf)^2
 =2\int_0^t P_s|\nabla P_{t-s}f|^2\,ds
 \ge(1-e^{-2t})|\nabla P_tf|^2.
\]
The left side is at most one. Since $f=1$ on the old $R_0$-neighborhood
of $U_b(c)$, the displacement estimate gives $P_tf\ge1-\alpha$ on $U_b(c)$.
Put $\phi=\min(1,\max(0,(P_tf-\alpha)/(1-2\alpha)))$. Then
\[
 \phi|_{U_b(c)}=1,\quad \operatorname{Lip}\phi\le M(t,\alpha),\quad
 E(\phi)\le(1-2\alpha)^{-2}E(f),\quad
 \operatorname{area}\{\phi\ne0\}\le\alpha^{-1}\int f\le Cu.
\]
The energy contraction follows from the heat semigroup; mass preservation
and Markov's inequality give the last bound. This proves
\eqref{eq:quantitative-cusp-cutoff} with constants independent of $N,u$.
For the qualitative assertion choose $u$ large enough that
$C/\log u\le e$ and $Cu^{-1/3}\log u\le p$, and take $A=Cu$.
\end{proof}

The proof gives the qualitative conclusion of
Lemma~\ref{lem:uniform-cusp-cutoff} for every $N\ge3$, so the threshold
$N_0(e,p)$ may be taken to be $3$.

\section{Opening a cusp while controlling cutoff energy}
\label{app:cusp-opening}

This appendix transfers the cusp cutoffs of
Lemma~\ref{lem:uniform-cusp-cutoff} from WP spheres to the bordered
law, in which one cusp is replaced by a geodesic boundary of length
$L$; this is Lemma~\ref{lem:cusp-opening}, the second input to the
proof of Theorem~\ref{thm:dense-collapse} in
Appendix~\ref{app:dense-collapse-details}, where the bordered law is
the interior law of a planar cap. Opening the cusp is realized by a
measurable map on the tree encoding of
Budd--Meeusen--Zonneveld~\cite[Theorem~1, Section~2.5]{BMZ}. The
bordered law is a mixture of the image of the cusped law under this
map and a remainder of mass $O(L^2)$, and a cutoff equal to one high
in the target cusp is transported, at a small energy cost, to a cutoff
equal to one near the new boundary. The proof has four parts, marked
below: the mixture identity, the geometry of the changed leaf, the
matching of the function across the added wedge, and the tail bound
for the leaf height, which uses the branching estimates of
\cite[Section~5.1]{BC}.

We use the positive Delaunay tree cells of \cite{BMZ}.
Distinguish two labeled cusps of a WP sphere, the origin and the
target. In the full-measure tree chart, keep the tree and its angular
labels and replace the target cusp by a boundary leaf of length $L$,
with its two boundary labels $w=v=L/2$. Reconstruct the surface and
take its convex core. Denote this measurable map by $F_{N,L}$.

\begin{lemma}[Coupling and energy transport]\label{lem:cusp-opening}
Let $\mathbb P_{N,L}$ denote the normalized WP law with one boundary
of length $L$ and $N-1$ labeled cusps. There is a probability measure $R_{N,L}$
such that
\begin{equation}\label{eq:opening-mixture}
 \mathbb P_{N,L}=p_N(L)(F_{N,L})_*\mathbb P_{N,0}
                   +(1-p_N(L))R_{N,L},\qquad
 p_N(L)=\frac{V_{0,N}}{V_{0,N}(L,0^{N-1})},
\end{equation}
and $1-p_N(L)\le CL^2$ for bounded $L$, uniformly in $N$.
There is a random leaf height $H_Q$, under $\mathbb P_{N,0}$, with
$\sup_N\mathbb P\{H_Q>H\}\le C/H^2$ for $H\ge2$, with the following
property. Suppose $0\le\phi\le1$ is globally hyperbolic
$M$-Lipschitz, with finite energy and support area $a$, and equals
one in the normalized target cusp strip above height $Y>\max(H_Q,1)$.
If $LY$ is sufficiently small, $F_{N,L}(X)$ admits $\phi_L\in H^1$
equal to one near its whole geodesic boundary, with
\begin{equation}\label{eq:opening-energy}
 0\le\phi_L\le1,\qquad
 E(\phi_L)\le E(\phi)+CM^2L^2Y,\qquad
 \operatorname{area}\{\phi_L\ne0\}\le a+CL.
\end{equation}
The constants are independent of $N$ and of the other short curves.
\end{lemma}

\begin{proof}
\emph{The probability law.}
In the tree encoding, a boundary vertex of degree $d$ contributes
two $(d-1)$-simplices of size $L/2$. Its angular polytope is independent
of $L$. On the event that the target boundary has degree one, both
simplices have dimension zero and mass one. The remaining tree and
angular coordinates have exactly the cusped WP density, including
the common factor $2^{N-3}$. Thus this event has mass $p_N(L)$ and,
conditionally on it, the law is exactly $(F_{N,L})_*\mathbb P_{N,0}$.
This proves \eqref{eq:opening-mixture}. The boundary-volume inequality
\cite[Lemma~2.4]{HTSphere} gives
\[
 1-p_N(L)\le1-\frac{L/2}{\sinh(L/2)}\le CL^2.
\]
The two measures compared are on the bordered space; no identification
of the metrics of $X$ and $F_{N,L}(X)$ is being made.

\emph{The changed leaf.}
The reconstruction in \cite[Lemmas~12--13, Section~2.5]{BMZ} fixes
each quadrilateral's metric and all inner gluings by the angular
labels and finite endpoints. Only the target leaf adds a wedge and
changes its outer-side gluing. Normalize its old quadrilateral $Q$
to have ideal vertices $0,1,\infty$ and finite vertex $1/2+ih$.
The two lower sides have maximum height $H_Q\le\max(h,1)$;
above $H_Q$ it is the full strip $0\le x\le1$.
The triply punctured sphere has a fixed limiting tile and will be
included by enlarging constants.

Here is the gluing in these coordinates. Put $q=e^{L/2}$ and
$d=q+q^{-1}-2=4\sinh^2(L/4)$. In the rib normalization of
\cite[Lemma~12]{BMZ}, the three successive radii are $1,q,q^2$,
with centers $0,1-q,(q-1)^2$. A generator is
$T_{\mathrm{old}}(z)=q^2z+(q-1)^2$. The real M\"obius map
\[
 G(z)=1-\frac{q}{q-1}\frac{z+1}{z-1}
\]
sends the quadrilateral's ideal vertices $-1,1-2q,1$ to $1,0,\infty$.
Conjugating the inverse generator gives
\begin{equation}\label{eq:opening-generator}
 T_d(z)=\frac{(1+d)z+1}{dz+1},\qquad
 \operatorname{tr}T_d=2+d=2\cosh(L/2).
\end{equation}
At $L=0$ this becomes $z\mapsto z+1$.
The new wedge lies between $x=0$ and the circle of endpoints $-1/d,0$.
Its near-zero side, when $dy\le1/4$, is $x=-w(y)$, where
\begin{align}
 w(y)&=\frac{1-\sqrt{1-4d^2y^2}}{2d},&
 dy^2&\le w(y)\le2dy^2,& |w'(y)|&\le Cdy.\label{eq:wedge-width}
\end{align}
On that circle $x^2+y^2=-x/d$, so direct substitution gives
\begin{equation}\label{eq:wedge-seam}
 T_d(-w(y)+iy)=1+iy',\qquad
 y'=\frac{y}{1-dw(y)},\qquad
 1\le\frac{y'}y\le1+Cd^2y^2,\quad
 \left|\frac{dy'}{dy}\right|\le C.
\end{equation}

\emph{Matching the function.}
Let $b_\phi(y)=\phi(iy)=\phi(1+iy)$ be the old common trace.
The Lipschitz bound gives $|b_\phi'(y)|\le M/y$ almost everywhere
and $|b_\phi(y')-b_\phi(y)|\le CMd^2y^2$.
Keep $\phi$ on every surviving old tile. The part of the added wedge
with $y<1/(2d)$ to the left of the circle abuts the ideal interval
$(-\infty,-1/d)$, lies outside the axis found below, and is discarded
with the funnel. On the remaining part of the added wedge below
$Y$ set
\[
 \theta=-x/w(y),\qquad
 \phi_L(x,y)=(1-\theta)b_\phi(y)+\theta b_\phi(y'),
\]
and set $\phi_L=1$ above $Y$. This matches both seams exactly by
\eqref{eq:wedge-seam}; at height $Y$ both traces are one. All other
seams retain their old matching values. Differentiation, using
$|\partial_y\theta|\le C/y$, gives
\[
 |\partial_x\phi_L|\le CMd,\qquad
 |\partial_y\phi_L|\le CM/y.
\]
The conformal Dirichlet integral on the wedge therefore satisfies
\begin{equation}\label{eq:wedge-energy}
 E_W(\phi_L)\le CM^2\int_0^Y
       \bigl(w(y)/y^2+d^2w(y)\bigr)\,dy
       \le CM^2(dY+d^3Y^3)\le CM^2dY.
\end{equation}
In particular there is no divergence at the ideal origin corner.

The axis of \eqref{eq:opening-generator} is the semicircle centered
at $1/2$ with radius $R=\sqrt{1/4+1/d}$. It intersects the curved
rib at $x=-1/(1+d)$, $y=1/(\sqrt d(1+d))$, and the vertical ribs
at height $1/\sqrt d$. For small $LY$ its whole segment in the leaf
block lies above $Y$. Thus the geodesic boundary lies where
$\phi_L=1$, and the removed funnel also has zero energy. The part
of the added wedge retained in the core is below height $R$ and
has area at most
\[
 \int_0^R w(y)y^{-2}\,dy\le2dR\le CL.
\]
Here $dR\le1/4$ for small $L$. These estimates, $d\asymp L^2$,
and the matching traces prove \eqref{eq:opening-energy} and the
global $H^1$ assertion. No global bilipschitz comparison is required.

\emph{The height probability.}
The angular polytopes identify with those of Budd--Curien by setting
their opposite-corner angle equal to $\pi$ minus the outgoing angle
in \cite{BMZ}. The inner sum $2\pi$ becomes $\pi$, the opposite
inequality $>\pi$ becomes the edge inequality $<\pi$, and the leaf
angle becomes zero. The Jacobian on two free coordinates per inner
vertex is one, and both measures have factor $2^{N-3}$.
At the first vertex next to the target leaf, write its other two
outgoing angles as $\alpha,\beta$ and the target angle as
$\varphi=\pi-\alpha-\beta$. With $n=N-1$, the branching identity
\cite[Section~5.1, (30)]{BC} gives the joint density
\[
 \mathbf1_{\{\alpha,\beta>0,\ \alpha+\beta<\pi\}}
 \frac{F(\alpha)F(\beta)}{F(0)\mathbb P_0\{|\tau|=n\}}
 \sum_{r+s=n+1}\mathbb P_\alpha\{|\tau|=r\}
                         \mathbb P_\beta\{|\tau|=s\}.
\]
Here $F$ and $\mathbb P_\theta$ are the branching weights and laws
of that reference. Its bounds (28), (31) give a uniform upper bound
$Cr^{-3/2}$ on the probabilities and the lower bound $cn^{-3/2}$
at type zero. The weight $F$ is bounded above and $F(0)>0$.
Splitting the convolution at $(n+1)/2$ bounds this density uniformly.
Consequently $\mathbb P\{\pi-\varphi<u\}\le Cu^2$.
Here $\varphi$ is the exit angle towards the target leaf in the
convention of \cite{BMZ}, in which the three outgoing angles at an
inner vertex sum to $\pi$; in the normalized quadrilateral $Q$, which
is symmetric about $\Re z=1/2$, it is the angle at the finite vertex
$1/2+ih$ of the half-plane triangle with vertices $0$, $1/2+ih$,
$\infty$, and hence equals half of the interior angle of $Q$ at
$1/2+ih$. The half-plane triangle gives $h=\tfrac12\tan(\varphi/2)$:
the circle through $0$ and $1/2+ih$ has center $h^2+1/4$ on the real
axis, so the angle at $1/2+ih$ between this circle and the vertical
ray to $\infty$ satisfies
$\cos\varphi=(1/4-h^2)/(1/4+h^2)$. Since $H_Q\le\max(h,1)$,
the required bound $\mathbb P\{H_Q>H\}\le C/H^2$ follows.
\end{proof}

\section{Growing caps and the dense collapse proof}
\label{app:dense-collapse-details}

This appendix supplies two proofs from
Section~\ref{sec:dense-spectrum}: the factorial-moment proof of
Lemma~\ref{lem:dense-caps}, the Poisson process of caps with growing
cusp count, and the three-case variational proof of
Theorem~\ref{thm:dense-collapse}, the dense spectral collapse. The
first proof uses the volume inputs of Part~\ref{part:volumes}, namely
the local limit theorem (Proposition~\ref{prop:local-limit}), the
exact zero-boundary factor \eqref{eq:TL50} and the coefficient bounds
\eqref{eq:PC3}--\eqref{eq:PC4}. The second uses the cap process
together with the cusp cutoffs of Lemma~\ref{lem:uniform-cusp-cutoff}
and their transport to bordered caps in Lemma~\ref{lem:cusp-opening},
and, above the cubic scale, the genus bound of \cite[Corollary~5.6]{GNY}.
Throughout we use the parameters of \eqref{eq:dense-parameters} and
the scales in \eqref{eq:dense-scales}.

\subsection{Proof of Lemma~\ref{lem:dense-caps}}
The mixed factorial moments of disjoint caps come from the exact
zero-boundary factor together with local coefficient ratios; nested
configurations are bounded by an induction over their cut trees. The
proof first records the coefficient shift \eqref{eq:dense-shifts},
then computes the factorial moments of disjoint caps, and finally
bounds every nested tree type by a power of $\eta$.

\begin{proof}
Write $h=2g-1$, $m=n-1$ and $f_m=[z^m]\Phi_g'(z)$.
By \eqref{eq:TL36}--\eqref{eq:TL37}, the model standard deviation
is asymptotic to $\sqrt{5h/4}/\epsilon$. Thus
Proposition~\ref{prop:local-limit}, applied sequentially, gives
\begin{equation}\label{eq:dense-shifts}
 \frac{f_{m-u}}{f_m}=x^u(1+o(1))
 \quad\text{uniformly for }0\le u\le C/\epsilon
\end{equation}
for every fixed $C$. We also retain the global upper bound
\eqref{eq:PC3}. This uses the model center $m$, not an identity with
the exact coefficient-law mean.

For $r$ ordered disjoint caps, put $K=\sum_i k_i$. Equations
\eqref{eq:TL50} and \eqref{eq:PC4} give the exact zero-boundary factor
\[
 \frac{n\prod_{a=1}^{r-1}(n-K+a)}{P^r}
 \prod_{i=1}^r\frac{V_{0,k_i+1}}{k_i!P^{k_i-2}}
 \frac{f_{m-K+r}}{f_m}.
\]
There is no ordering divisor, and the genus-$g$ side distinguishes
the two sides of each curve. On a fixed mark window $K=O(1/\epsilon)=o(n)$.
Integrating each twist factor $\ell\,d\ell$ over a score interval
$[s_i,t_i]$ and using \eqref{eq:dense-shifts} therefore gives
\[
 (1+o(1))\prod_i(t_i^2-s_i^2)p_x(k_i),\qquad
 p_x(k)=\frac{(k-1)^2w_k(x)}{M_2(x)}.
\]
The normalization is $Q=nx/(4P)$. The Laplace transform of the mark is
$e^{-u\epsilon}M_2(xe^{-u\epsilon})/M_2(x)$, tending to
$(1+u)^{-1/2}$, the transform of $\operatorname{Gamma}(1/2,1)$.

To bound nested configurations, root their cut tree at its unique
genus-$g$ vertex. A planar vertex with $c$ children has local series
$G_{c+1}=\Phi_0^{(c+1)}$. If $T_v$ is the weighted series of a
subtree with $s$ planar vertices, its recursion is
\[
 T_v=(D-1)^2\left(G_{c+1}\prod_{w\text{ child of }v}T_w\right),
 \qquad D=x\partial_x.
\]
The stable initial powers give the square of the enclosed cusp count
minus one. The square-root expansion of $Y=\Phi_0''$ implies, for
every fixed $a\ge0$,
\[
 D^aT_v=O\bigl(\epsilon^{-(2s-3/2)-a}\bigr).
\]
For a leaf the exponent is $1/2$. A unary vertex multiplies by bounded
$Y$ and adds two derivatives. At a vertex with $c\ge2$ children,
$G_{c+1}=Y^{(c-1)}$ has exponent $c-3/2$; adding the child exponents
and two derivatives gives $2s-3/2-c/2$, no larger than required.
The same argument with the additional $a$ derivatives closes the induction.

For a tree with $r$ edges and root degree $d$, cusp allocation and
\eqref{eq:PC3} bound its prefactor by $n^dP^{-r}x^{-d}$ times the
product of the root-subtree series. Stability ensures the coefficient
shift is nonnegative. Extending the positive cusp sums only enlarges
the bound. The score integrals contribute $S^{-2r}$, so its total
contribution is at most
\[
 C_{r,T}n^dS^{-2r}\epsilon^{-2r+3d/2}
       =O_{r,T}(\eta^{r-d}).
\]
Every nonstar has $d<r$. All physical lengths in the fixed window are
$O(\sqrt\eta)$, so the collar lemma makes them disjoint as curves;
the only remaining overlap is nesting. The boundary-volume bound
\eqref{eq:PS4} changes the integrals by $1+o(1)$ uniformly. There
are finitely many tree types for fixed $r$. The mixed factorial
moments now give \eqref{eq:dense-cap-intensity}; the two-edge nonstar
estimate gives the stated nested-pair bound.
\end{proof}

\subsection{Proof of Theorem~\ref{thm:dense-collapse}}
The proof separates three limits of $n/g^3$, and a subsequence
argument at the end reduces a general sequence satisfying
\eqref{eq:dense-window} to these cases. Below $g^3$, caps of the
Poisson window of Lemma~\ref{lem:dense-caps} carry the cutoff trial
functions supplied by Lemmas~\ref{lem:uniform-cusp-cutoff}
and~\ref{lem:cusp-opening}; at the crossover, a second-moment count of
smaller caps replaces the Poisson limit; above $g^3$, the genus bound
of Grigor'yan--Netrusov--Yau already suffices. In the first two cases
the trial functions are centered and inserted into the min--max
principle on a mean-zero subspace.

\begin{proof}[Proof of Theorem~\ref{thm:dense-collapse}]
We treat three possible limits of $n/g^3$.

\emph{The case $n/g^3\to0$.}
Fix the horocycle length $b$ and Lipschitz constant $M$ of
Lemma~\ref{lem:uniform-cusp-cutoff}, and put
\begin{equation}\label{eq:dense-cutoff-scales}
 K=\epsilon^{-1},\qquad \ell_*=K/S,\qquad H=\ell_*^{-1/2}.
\end{equation}
By \eqref{eq:dense-scales} and \eqref{eq:dense-window},
\[
 K\sim\frac{2n}{5g},\qquad
 \ell_*\asymp(n/g^3)^{1/4}\longrightarrow0,\qquad
 \ell_*\log K\longrightarrow\infty.
\]
Here $\log K\sim2\log g$: indeed $K\asymp n/g$, and in this case
\eqref{eq:dense-window} gives $n\gg g^3/(\log g)^4$, so
$g^{2-o(1)}\le K\le n/g$ while $n\ll g^3$ gives $K\ll g^2$; hence
$\log K=(2+o(1))\log g$. All these parameters precede sampling.

Fix $T$ and retain caps with score at most $T$ and mark in $[1,2]$.
They have $k\asymp K$ and waist $L\le C_T\ell_*$. Apply
\eqref{eq:quantitative-cusp-cutoff} with $N=k+1$ and $u=\sqrt N$.
Uniformly over the window, its energy, support area and failure
probability are bounded respectively by
\[
 e=\frac C{\log K},\qquad A=C\sqrt K,\qquad
 p=CK^{-1/6}\log K=o(1).
\]
Use the height screen $H_Q\le H$ of Lemma~\ref{lem:cusp-opening},
and put $Y=H+2+1/b$. Then $LY\to0$. The normalized leaf strip
above $Y$ lies in the target cusp tail of length $b$.

In the exact one-cap Campbell formula the interior law is WP on
$\mathcal M_{0,k+1}(L,0^k)$. Lemma~\ref{lem:cusp-opening} couples
it with the opened cusped law with failure probability $O(L^2)$.
The probability of failure of the cutoff or height screen at fixed
$k,L$ is therefore at most $p+C/H^2+CL^2$. No independence is
needed. By \eqref{eq:PC3}--\eqref{eq:PC4} and \eqref{eq:PS4}, the
one-cap intensity is at most $CQw_k(x)L\,dL$. Hence
\begin{equation}\label{eq:dense-bad-caps}
 \mathbb E\#\{\text{bad caps in the window}\}
 \le C_T(p+H^{-2})+C_T\frac{QM_4}{S^4}
 =O_T(p+H^{-2}+\eta)=o(1).
\end{equation}
The nested-pair estimate also tends to zero. With probability tending
to one, all caps in the window are good and have disjoint interiors.
Their count converges to $\operatorname{Pois}(T^2\gamma)$, where
$\gamma=\mathbb P\{\operatorname{Gamma}(1/2,1)\in[1,2]\}>0$.

On the event of at least $r$ such caps, choose $r$ of them. On cap $i$
put $u_i=1-\phi_{L_i}$ and extend it by zero to the whole surface.
Since $\phi_{L_i}=1$ near the boundary waist by the opening lemma,
$u_i$ is an $H^1$ function on the cap vanishing in a neighborhood
of its boundary, so $u_i\in H^1_0(\mathrm{cap}_i)$; its extension by
zero therefore lies in $H^1(X)$, with the same Dirichlet energy and
$L^2$ norm. With $B_i=2\pi(k_i-1)$ and
$B=2\pi(2g-2+n)$, it gives
\[
 E(u_i)\le e+C_T\ell_*^2Y,\qquad
 \int u_i^2\ge B_i-A-CL_i,\qquad 0\le\int u_i\le B_i.
\]
For $u=\sum_i z_i u_i$, disjointness and Cauchy--Schwarz imply
\begin{align*}
 \int|u-\overline u|^2
 &\ge\sum_i z_i^2(B_i-A-CL_i)
      -\frac{\sum_i B_i}{B}\sum_i z_i^2B_i\\
 &\ge(1-o(1))\sum_i z_i^2B_i.
\end{align*}
Indeed $A/K\to0$ and $\sum_iB_i/B=O_r(g^{-1})$.
Centering preserves dimension: a linear combination constant on the
whole surface is zero on the positive-area complement, hence zero
everywhere, and disjoint supports then force all coefficients to vanish.
The energy is $\sum_i z_i^2E(u_i)$. Min--max on this mean-zero
$r$-dimensional subspace gives
\[
 S\lambda_r^+(X)
 \le C_{r,T}\frac SK(e+\ell_*^2Y)
 \le C'_{r,T}\left(\frac1{\ell_*\log K}+\ell_*Y\right)
 \longrightarrow0.
\]
The unscaled quotient is eventually below $1/4$, so these are actual
discrete eigenvalues. For each $z>0$, the preceding event estimates give
\[
 \limsup_{g\to\infty}\mathbb P\{S\lambda_r^+(X)>z\}
       \le\mathbb P\{\operatorname{Pois}(T^2\gamma)<r\}.
\]
Send $T\to\infty$ after the limsup and use $A_g\sim PS$.

\emph{The case $n/g^3\to c\in(0,\infty)$.}
The fixed-mark Poisson argument no longer applies, so we count smaller
caps. Put $J=\log g$ and choose
\[
 k_0=\lfloor g^2J^{-4/5}\rfloor,\qquad
 L=J^{-1/2},\qquad H=J^{1/10}.
\]
Let $Z$ count planar caps with $k_0\le k\le2k_0$ cusps and waist
at most $L$. Now $\epsilon k_0\to0$, $k_0/n\to0$ and
$S\asymp g^2$. The coefficient estimate \eqref{eq:dense-shifts}
continues to hold: its proof uses shifts $O(\epsilon^{-1})$
relative to a standard deviation of order $\sqrt g/\epsilon$,
not $\eta\to0$. The genus-zero coefficient asymptotic gives
$w_k(x)\asymp k^{-5/2}$ uniformly on this band. Thus the exact
one-cap formula and \eqref{eq:PS4} give
\begin{equation}\label{eq:crossover-count}
 \mu:=\mathbb EZ=(1+o(1))QL^2\sum_{k=k_0}^{2k_0}w_k(x)
       \asymp nL^2k_0^{-3/2}\asymp J^{1/5}\longrightarrow\infty.
\end{equation}
The ordered disjoint-pair formula in the preceding lemma has root
factor $n(n-k_1-k_2+1)=n^2(1+o(1))$. Uniform coefficient ratios
and the boundary-volume bound therefore show that the expected
number of ordered pairs with disjoint interiors is $(1+o(1))\mu^2$.

All selected waists are disjoint as curves by the collar lemma;
their planar sides are either disjoint or nested. For an outer cap
with $k$ cusps, the Campbell interior law is WP on
$\mathcal M_{0,k+1}(\ell,0^k)$, with $\ell\le L$.
An eligible inner cap has $l\ge k_0$ cusps, and the intervening
component has $a=k-l\ge1$ cusps by stability. The two-sided
genus-zero volume asymptotic~\cite[Theorem~2.2]{HTSphere}, with
constants enlarged to cover small indices, gives
\[
 \binom{k}{a}\frac{V_{0,a+2}V_{0,l+1}}{V_{0,k+1}}
 \le C\frac{k^{5/2}}{a^{3/2}l^{5/2}}.
\]
Since $k/l\le2$, summing over $a\ge1$ is bounded uniformly.
The boundary-volume factors are also bounded uniformly for
$\ell,t\le L$. Integrating the inner twist $t\,dt$ therefore
bounds the expected number of eligible descendants by $CL^2$.
If $N$ counts unordered nested pairs, another Campbell integration
gives $\mathbb EN\le C\mu L^2=o(\mu)$. Consequently
\[
 \mathbb E[Z(Z-1)]=(1+o(1))\mu^2+O(\mu L^2),\qquad
 Z/\mu\longrightarrow1\quad\text{in probability}.
\]

Apply Lemma~\ref{lem:uniform-cusp-cutoff} with $u=\sqrt{k+1}$
and Lemma~\ref{lem:cusp-opening} with the chosen $H$ and
$Y=H+2+1/b$. Here $LY\to0$. Uniformly on the band, the cutoff
energy is at most $C/\log k_0$, its support area at most $C\sqrt{k_0}$,
and its failure probability at most $p=Ck_0^{-1/6}\log k_0=o(1)$.
The expected number of bad caps is at most
$C\mu(p+H^{-2}+L^2)=o(\mu)$, including the exact coupling failure.
Delete these caps and both endpoints of every nested pair. By Markov's
inequality the number deleted is $o_{\mathbb P}(\mu)$, so the number
of remaining good caps with disjoint interiors tends to infinity
in probability. No independence of their shapes is needed.

Choose any fixed $r$ of them. The centered min--max argument above
applies with $B_i\asymp k_0$, support loss $O(\sqrt{k_0}+L)$ and
$\sum_iB_i/B=O_r(k_0/n)=o(1)$. With probability tending to one,
\[
 S\lambda_r^+(X)
 \le C_r\frac{S}{k_0}\left(\frac1{\log k_0}+L^2Y\right)
 \le C'_r\bigl(J^{-1/5}+J^{-1/10}\bigr)\longrightarrow0.
\]
Again the unscaled Rayleigh quotients are below $1/4$, so these
are actual discrete eigenvalues.

\emph{The case $n/g^3\to\infty$.}
We use the classical genus bound in its finite-puncture form.
Grigor'yan--Netrusov--Yau~\cite[Corollary~5.6]{GNY} show that, for
every positive nonatomic Radon measure $\sigma$ on $X$, the form
$f\mapsto\int_X|\nabla f|^2-\int_Xf^2\,d\sigma$ is negative on a
subspace of compactly supported Lipschitz functions of dimension at
least a fixed positive multiple of $\sigma(X)/(g+1)$. Taking
$\sigma=t\,dA_X$ and applying min--max below the bottom $1/4$ of the
essential spectrum gives
\[
 \#\{\lambda_j(X)<t\}\ge c\,\frac{t\,\operatorname{area}(X)}{g+1},
 \qquad 0<t<1/4,
\]
where the count includes zero and $c>0$ is absolute. Their statement
applies to a compact Riemann surface minus finitely many points with
any conformal metric; the hyperbolic area measure is finite and
nonatomic. Taking $t=C_r(g+1)/\operatorname{area}(X)$ with $C_r$
large enough gives $r$ positive eigenvalues below $t$ once $t<1/4$.
Since $\operatorname{area}(X)=2\pi(2g-2+n)$, \eqref{eq:dense-scales}
then yields, deterministically,
\[
 A_g\lambda_r^+(X)\le C'_r(g^3/n)^{1/4}\longrightarrow0.
\]

For a general sequence satisfying \eqref{eq:dense-window}, every
subsequence has a further subsequence along which $n/g^3$ converges
in $[0,\infty]$. The three cases prove \eqref{eq:dense-collapse}
on each such subsequence, hence on the original sequence.
Finally, for every $t>0$, the number of eigenvalues with
$A_g\lambda\le t$ tends to infinity in probability. Integrating
against a compactly supported continuous function equal to one on
$[0,t]$ proves non-tightness on $[0,\infty)$.
\end{proof}

\section{Spectral comparison at growing cusp density}
\label{app:dense-comparison}

We prove Theorem~\ref{thm:dense-persistence}, the Poisson law of the
low spectrum at growing cusp density. Throughout, $g\to\infty$,
$n/g\to\infty$, and $\epsilon,Q,S,\eta$ have the exact definitions
\eqref{eq:dense-parameters}. Constants below are absolute in this
endpoint regime. The strategy is that of the bounded-density case in
Section~\ref{sec:supercritical-spectrum}: cut the surface along short
planar curves, compare its energy form with a weighted star graph, and
read off the eigenvalues from the cap process of
Lemma~\ref{lem:dense-caps}; the difference is that the retained
waists need not be isolated leaves, so the intervening branches are
contracted and their attachment errors are controlled jointly. The
geometric inputs are Mirzakhani integration with the coefficient bounds
\eqref{eq:PC3}--\eqref{eq:PC4} and the envelope
\eqref{eq:gevrey-envelope}, the cut-counting argument of
Lemma~\ref{lem:super-cuts} (Appendix~\ref{app:super-cuts}), with the
endpoint estimates established below, and the
relative gap and comparison estimates of
Lemmas~\ref{lem:relative-geodesic} and~\ref{lem:super-comparison}
(Appendix~\ref{app:super-comparison}). The first subsection excludes
competing cuts; the second controls the depth, cap sizes and Neumann
gaps of the pieces cut out by short planar curves; the third bounds
the joint attachment errors of the retained waists; the fourth proves
the two-sided min--max comparison on the entire low spectral subspace;
and the last identifies the limit with the growing-cap Poisson process
and proves the Cheeger assertion.

\subsection{Competing cuts at the endpoint}

Up to small probability, the only connected two-sided cuts with a small
ratio of length to smaller area are single curves with a planar side.
This is the estimate \eqref{eq:dense-cut-exclusion} below; it extends
Lemma~\ref{lem:super-cuts} to the endpoint regime by keeping the
endpoint factors in that proof explicit, and it is used twice later,
for the piece gaps in \eqref{eq:dense-piece-gap} and for the Cheeger
assertion in the last subsection. The argument treats in turn cuts
with two positive-genus sides, planar sides with $2\le r\le g/2$
boundaries, and planar sides with more than $g/2$ boundaries.

Let $\mathcal B(u)$ be the event of a geodesic multicurve dividing
$X$ into two connected components, of length divided by smaller area
at most $u$, other than a single curve with a planar side. We claim
\begin{equation}\label{eq:dense-cut-exclusion}
 \mathbb P(\mathcal B(u))\le C\left(
      \frac{u^4}{\epsilon^3}+\frac{u^2}{g\epsilon^{3/2}}\right),
       \qquad 0<u\le c\epsilon.
\end{equation}
In the notation of the proof of Lemma~\ref{lem:super-cuts}, set
$\rho(x)=r(Y(x))$ and $Z(x)=a(Y(x))\rho(x)^{3/2}A(Y(x))$.
The convergent endpoint expansions give
\[
 \rho(x)\asymp\epsilon^{1/2},\qquad Z(x)\asymp\epsilon^{5/4},\qquad
 |x\partial_x\log\rho|+|x\partial_x\log Z|\le C/\epsilon.
\]
These bounds also hold at $xe^v$ for $0\le v\le c\epsilon$.

For two positive-genus sides use the exact allocation sum
$B_{i,j,r}(v)$ of Appendix~\ref{app:super-cuts}, with
$h_1=2i-2+r\le h_2=2j-2+r$ and $H=h_1+h_2=2g-2$.
The positive convolution on the untilted large component gives
\[
 B_{i,j,r}(v)\le\frac{Cn}{P^rx}
       \frac{\Phi_i^{(r)}(xe^v)\Phi_j^{(r)}(x)}{\Phi_g'(x)}.
\]
Indeed its maximal point mass and the parent's central point mass
are respectively bounded above and below by constants times
$\epsilon/\sqrt g$, uniformly at this endpoint. The all-marking
envelope \eqref{eq:gevrey-envelope} and the parent lower bound yield
\begin{align*}
 e^{vh_1}\frac{\Phi_i^{(r)}(xe^v)\Phi_j^{(r)}(x)}{\Phi_g'(x)}
 &\le C^re^{Cvh_1/\epsilon}Z(x)\rho(x)^{r-1/2}
 \frac{\Gamma(h_1)\Gamma(h_2)}{\Gamma(H+1)}
       \sqrt{\frac{H+1}{h_1h_2}},\\
 Z(x)\rho(x)^{r-1/2}&\le C^r\epsilon^{1+r/2}.
\end{align*}
The tilt costs $h_1/\epsilon$, since $r\le h_1$; it does not cost
$H/\epsilon$. Take the exponential moment parameter $s=c_0\epsilon$
and $v=s+2\pi u$. The beta bound in Appendix~\ref{app:super-cuts},
with $c_0$ and $u/\epsilon$ small, and $n\epsilon\asymp g$, gives
\[
 e^{vh_1}B_{i,j,r}(v)
       \le\frac{C^r}{g}\epsilon^{r/2}e^{-c_1h_1}/\sqrt{h_1}.
\]
Mirzakhani integration contributes $(2\pi u/s)^{2r}/r!$.
Summing the genera and then every $r\ge1$ bounds their contribution
by $Cg^{-1}(e^{Cu^2\epsilon^{-3/2}}-1)
\le Cu^2/(g\epsilon^{3/2})$.

For a planar side with $r\ge2$ boundaries and $a$ cusps, put
$\chi_a=a+r-2>0$, $G=g-r+1$ and $N=n+2r-2$.
If $r\le g/2$, the saddle $x_r$ for $(G,N)$ satisfies
$1-x_r/x_c\asymp\epsilon$. Successive genus exchanges in
Theorem~\ref{thm:volume} give
\[
 V_{G,N}/V_{g,n}\le(C\sqrt\epsilon)^{r-1}.
\]
At every step the remaining genus is at least $g/2$ and the cusp count
is $n+O(g)$; only upper bounds, not relative errors, are multiplied.
The factorial pairing in \eqref{eq:super-planar-integral} and
\eqref{eq:PC3} therefore bound its zero-boundary prefactor by
\[
 C^r\epsilon^{(r-1)/2}n^{2-r}[z^a]\Phi_0^{(r)}(z)x_r^a.
\]
The length integral contributes
$(2\pi u)^{2r}\chi_a^{2r}e^{2\pi u\chi_a}/(r!(2r)!)$.
For $r\ge3$ choose $s=c_2\epsilon$ so that
$z_r=x_re^{2\pi u+s}$ stays a distance comparable to $\epsilon$
from $x_c$. The square-root expansion of $Y=\Phi_0''$ is analytic
in a slit neighborhood of $x_c$. Cauchy's estimate on a disk of
radius $c_3\epsilon$ around $z_r$, applied to $Y-y_0$, gives
\[
 \Phi_0^{(r)}(z_r)\le C^r(r-2)!\epsilon^{5/2-r}.
\]
Using $\chi_a^{2r}\le(2r)!s^{-2r}e^{s\chi_a}$, the expected count for this
$r$ is at most $C^rn^{2-r}\epsilon^{2-5r/2}u^{2r}$.
For $r=2$, $Y(x_c)\ne0$, so use instead
\[
 \sum_a a^4[z^a]Y(z)x_r^ae^{2\pi ua}
 \le Cs^{-3}(z\partial_zY)(x_re^{2\pi u+s})
 \le C\epsilon^{-7/2}.
\]
Together with the prefactor $\sqrt\epsilon$, this gives
$Cu^4/\epsilon^3$. The ratio of the series for $r\ge3$ is
$Cu^2/(n\epsilon^{5/2})=O(g^{-1}u^2\epsilon^{-3/2})$.
Hence the sum over $2\le r\le g/2$ is $O(u^4/\epsilon^3)$.

If $r>g/2$, the average boundary length is at most $Cu/\epsilon$,
and at least $g/4$ curves have length at most $Cu/\epsilon$.
Each is individually nonseparating: the remaining curves still join
the two connected sides. Integration and one genus exchange give
$\mathbb E N_{\mathrm{nonsep}}(l)\le C\sqrt\epsilon\,l^2e^l$.
Markov's inequality bounds this event by $Cu^2/(g\epsilon^{3/2})$.
This is an event estimate, not an expected count of large multicurves.
It also covers two planar sides and completes
\eqref{eq:dense-cut-exclusion}.

\subsection{Depth, cap size and piece gaps}

After cutting along short planar curves, we show that with probability
tending to one the resulting tree has logarithmic depth
\eqref{eq:dense-depth}, every short cap has at most $K$ cusps, and all
pieces have Neumann gaps at least a common lower bound $\gamma$
\eqref{eq:dense-piece-gap}. The depth and cap-size bounds are
Mirzakhani-integration counts using \eqref{eq:PC3} and the one-waist
estimate \eqref{eq:super-one-waist}; the gap bound combines the cut
exclusion of the preceding subsection with the relative Cheeger
estimate of Lemma~\ref{lem:relative-geodesic}. The gap $\gamma$ and the
parameters fixed in \eqref{eq:dense-proof-cutoffs} are the hypotheses
of the comparison lemma in the fourth subsection.

Set, before sampling,
\begin{equation}\label{eq:dense-proof-cutoffs}
 L=\log(n+2),\quad\delta=L^{-1},\quad b=L^{-1/2},\quad
 R=L^2,\quad\tau=R/S,\quad K=C_1L/\epsilon.
\end{equation}
Cut along all planar separating curves of length at most $\delta$.
Their dual graph is a tree with a genus-$g$ root. Short nonseparating
curves may remain in the root and are not discarded. They do not
affect the two comparison inputs used below: the constant
$H_{\mathrm{int}}$ in Lemma~\ref{lem:relative-geodesic} is an infimum
over internal multicurves with two connected complementary components,
so a single short nonseparating geodesic in the root is not such a cut,
and the multicurves containing it are controlled by the global cut
exclusion \eqref{eq:dense-cut-exclusion}; and the trimming and collar
estimates of Lemma~\ref{lem:super-comparison}, isolated in
Lemma~\ref{rev:lem:cut-tree-basic} below, involve only the boundary
geodesics of the pieces.

Write $F=\Phi_0'$, $Y=\Phi_0''$ and
$f_m=[z^m]\Phi_g'(z)$. An ancestry chain of $r$ waists has
$r-1$ intermediate annuli with cusp counts $a_i\ge1$, and a terminal
cap with $c\ge2$ cusps. With $K_0=\sum a_i+c\le n$, exact cancellation
of the cusp-label factorials gives the zero-boundary factor
\[
 \frac n{P^r}\prod_{i=1}^{r-1}[z^{a_i}]Y(z)\,[z^c]F(z)
                   \frac{f_{n-K_0}}{f_{n-1}}.
\]
Ancestry orders the curves, so there is no $r!$ divisor.
The global ratio \eqref{eq:PC3}, also for growing $r$, bounds the
last quotient by $Cx^{K_0-1}$. Boundary factors cost $e^{r\delta}$,
and each twist-length integral is $\delta^2/2$. Since $F(x)$ and
$Y(x)$ remain bounded, the expected number of such chains is at most
\begin{equation}\label{eq:dense-depth}
 Cnx^{-1}\left(\frac{\delta^2e^\delta}{2P}\right)^r
                    Y(x)^{r-1}F(x)\le Cn(C\delta^2)^r.
\end{equation}
For $r=\lceil C_0L\rceil$, this tends to zero. Thus the maximum
branch depth is $D\le C_0L$ with probability tending to one.

The bound $w_k(x)\le Ck^{-5/2}e^{-\epsilon(k-2)}$ and the one-waist
estimate \eqref{eq:super-one-waist}, with an endpoint-uniform constant
from \eqref{eq:PC3}, show that all $\delta$-short caps have at most
$K$ cusps except with probability at most
$Cn\delta^2e^{-\epsilon K/2}=o(1)$, for $C_1$ large.
Every nonroot piece and every descendant branch then have area at
most $2\pi K$.

Put $u=c\delta/K$. Equation~\eqref{eq:dense-cut-exclusion} excludes
competing connected cuts at ratio $u$ with probability $1-o(1)$.
Single planar cuts with $k>K$ and $\ell\le2\pi u(k-1)$ have expected
count at most
\begin{equation}\label{eq:dense-long-tail}
 Cnu^2\sum_{k>K}(k-1)^2w_k(x)e^{2\pi u(k-1)}=o(1).
\end{equation}
Here $u/\epsilon=O(L^{-2})$ absorbs the boundary exponential; choosing
$C_1$ large absorbs the polynomial prefactor.

An internal connected two-sided cut of the pruned genus-$g$ piece
extends to one of $X$ by assigning each removed maximal cap to its
attachment side. Both sides remain connected and their areas increase.
If its ratio were at most $u$, the preceding events would force a
single planar waist enclosing $k\le K$ cusps. Its planar side has area
$2\pi(k-1)\le2\pi K$, so its length would be at most $2\pi uK<\delta$;
such a curve has been cut, so it is not internal to the pruned piece. Every other piece
has area at most $2\pi K$ and no internal essential geodesic of length
at most $\delta$. Lemma~\ref{lem:relative-geodesic} therefore gives,
simultaneously for all geodesic pieces,
\begin{equation}\label{eq:dense-piece-gap}
 \operatorname{gap}_N(P_v)\ge\gamma=c\delta^2/K^2
                                      \asymp\epsilon^2L^{-4}.
\end{equation}
This is an event on the original surface, not an assertion that the
pruned root has unconditioned WP law.

\subsection{Joint attachment energies}

The retained waists are those of small conductance relative to the
mass they enclose. With probability tending to one they are not nested,
and a single Gram matrix measures all their attachment errors at once;
we bound its weighted norm $\mathcal L$ in expectation, which is the
estimate \eqref{eq:dense-access}. The inputs are the trimmed masses and
conductances of Lemma~\ref{lem:super-comparison}, the two-edge chain
count from the proof of Lemma~\ref{lem:dense-caps}, and the genus-zero
volume bound of \cite{ManZ} in the form \cite[Theorem~2.2]{HTSphere}.
The bound on $\mathcal L$ enters the energy comparison
\eqref{eq:dense-map-energy} in the next subsection only through the
product $\mathcal A\mathcal L$.

Trim half-collars at circumference $b$, with masses $m_v$, descendant
masses $W_e$ and conductances $\kappa_e$ as in
Lemma~\ref{lem:super-comparison}. Thus $W_e\le CK$ and
\[
 m_v=(1+O(b))\operatorname{area}(P_v),\qquad
 \kappa_e=(\ell_e/\pi)(1+O(\delta/b)).
\]
Retain edges with $\kappa_e/W_e\le\tau$. The two-edge chain calculation
in Lemma~\ref{lem:dense-caps}, as an upper bound for all cusp counts,
and the one-edge third area moment give respectively
\begin{equation}\label{eq:dense-retained-mass}
 \mathbb E\#\{\text{nested retained pairs}\}\le C\eta R^4,
 \qquad \mathbb E\sum_{e\text{ retained}}\operatorname{area}(\text{cap}_e)
                                      \le CR^2/\epsilon.
\end{equation}
Indeed the first cusp sum is $O(\epsilon^{-5/2})$ with prefactor
$n\tau^4$, and the second is $O(n\tau^2M_3)$.
Actual lengths are at most $\delta$, so the volume envelope is bounded
before extending each positive length integral to $C\tau(k-1)$.
Under \eqref{eq:dense-persistence-range}, $\eta R^4\to0$ and
$R^2/g\to0$. With probability tending to one the retained edges
are not nested; contracting the high-rate components gives a star,
whose central mass is at least half the total mass $A\asymp n$.
The retained edges need not have been root leaves before contraction.

We next bound all attachment errors together. On each high component
use its graph energy as a norm modulo constants. Let $\mathcal G$ be
the Gram matrix of the functionals that take the difference of the
two endpoint deviations from their respective component mass means,
one functional for each retained edge. Set
\[
 \mathcal L=\|\operatorname{diag}(\sqrt\kappa)\,
                   \mathcal G\operatorname{diag}(\sqrt\kappa)\|.
\]
For expectations set $\mathcal L=0$ off the non-nesting and central-mass
event just described. We prove
\begin{equation}\label{eq:dense-access}
 \mathbb E\mathcal L\le C\delta R^3\sqrt\eta.
\end{equation}
Positivity bounds the norm by the trace, so correlations between
attachments are retained.

For a rooted tree of total mass $A_C$, telescoping its edge differences
shows that the squared energy-dual norm of root value minus mass mean
is exactly $\sum_e(W_e/A_C)^2/\kappa_e$. On a WP planar cap with
$k$ cusps and one boundary of length at most $\delta$, its expectation
over the internal $\delta$-short curves is at most $C\delta$, uniformly
in $k$. To see this, the two-sided genus-zero volume bound
of Manin--Zograf~\cite{ManZ}, in the form
\cite[Theorem~2.2]{HTSphere}, gives, for $a+l=k$, $a\ge1$, $l\ge2$,
\[
 \binom ka\frac{V_{0,a+2}V_{0,l+1}}{V_{0,k+1}}
        \le C\frac{k^{5/2}}{a^{3/2}l^{5/2}}.
\]
The boundary envelope \eqref{eq:PS4} changes this by a bounded factor.
Since $1/\kappa\le\pi/\ell$, the twist integral cancels the inverse
length. The expected dual norm is at most
\[
 C\delta\sqrt k\sum_{a=1}^{k-2}a^{-3/2}(k-a)^{-1/2}\le C'\delta.
\]
Split at $a=k/2$ for the last inequality. The same bound holds with
trimmed masses. The exact marked-cap integration formula thus bounds
the retained cap contribution to the trace by
$C\delta n\tau^3\epsilon^{-3/2}$.

For a retained edge $f$ on the root component, the coefficient of a
high edge $e$ in its endpoint deviation is
$\mathbf1_{e\text{ ancestor of }f}-W_e^{\mathrm{root}}/A_0$.
As $A_0\ge A/2$, its square is at most
$2\mathbf1_{e\text{ ancestor of }f}+8(W_e/A)^2$.
The first term has expectation at most
$C\delta n\tau^3\epsilon^{-3/2}$: the two-curve chain has prefactor
$Cn[z^a]Y[z^c]F\,x^{a+c}$; the high length integral is at most
$C\delta$ and the retained one at most $C\tau^3c^3$.
For the second term, disjoint curves have prefactor $Cn^2$, cancelled
by $A^2$, and cusp moments $M_2M_3=O(\epsilon^{-2})$.
Their contribution is at most $C\delta\tau^3\epsilon^{-2}$.
For either nested ordering the prefactor after division by $A^2$ is
$C/n$. The cusp polynomial is bounded by a constant times
$(a+c)^2c^3+(a+c)^3c^2$, whose $Y,F$ sum is
$O(\epsilon^{-7/2})$. This gives
$C\delta\tau^3n^{-1}\epsilon^{-7/2}$.
These two bounds are smaller than the first by factors
$O((n\sqrt\epsilon)^{-1})$ and $O((n\epsilon)^{-2})$.
There is no diagonal pair, since retained and high edges are disjoint
sets. Extending the positive sums beyond the selected configuration
only enlarges the bounds. Finally
$n\tau^3\epsilon^{-3/2}\asymp R^3\sqrt\eta$, proving
\eqref{eq:dense-access} without independence assumptions.

\subsection{Comparison on the whole low spectral subspace}

We compare the energy on $X$ with the star form on the retained edges
in both directions. The resulting counting sandwich
\eqref{eq:dense-count-sandwich} holds on the entire spectral subspace
below $\Lambda$, so residual modes are included. The subsection first
isolates, in Lemma~\ref{rev:lem:cut-tree-basic}, the three inequalities
from the proof of Lemma~\ref{lem:super-comparison} that survive without
the leaf and depth hypotheses of that lemma; it then builds the map $J$
from the low spectral subspace to the star, bounds its mass loss
\eqref{eq:dense-map-mass} and energy gain \eqref{eq:dense-map-energy},
and finally checks that every error term vanishes under
\eqref{eq:dense-persistence-range}.

The three basic inequalities are isolated in the following lemma, whose
hypotheses involve only the gap $\gamma$, the trimming circumference
$b$ and the cut length $\delta$.

\begin{lemma}[Poincar\'e, trace and collar estimates on a cut tree]
\label{rev:lem:cut-tree-basic}
Let $X$ be a finite-area hyperbolic surface cut along disjoint simple
closed geodesics of length at most $\delta$, where $\delta<b\le b_0$ for
a sufficiently small absolute constant $b_0$, and suppose that every geodesic
piece $P_v$ has Neumann gap at least $\gamma\in(0,1]$. Trim each
half-collar at circumference $b$, let $B_v$ be the trimmed pieces with
areas $m_v$, and let $\kappa_e$ be the collar conductances of
Lemma~\ref{lem:super-comparison}. For real $f\in H^1(X)$ write
$E=\int_X|\nabla f|^2$, $a_v$ for the mean of $f$ on $B_v$, and
$p_e,q_e$ for its angular averages on the two outer boundary circles of
the collar of the edge $e$. Then, with absolute constants $C$,
\begin{equation}\label{rev:eq:cut-tree-basic}
 \begin{gathered}
 \|f\|_{L^2(X)}^2\le(1+Cb)\sum_vm_va_v^2+C\gamma^{-1}E,\qquad
 \sum_{\mathrm{endpoints}}|p_e-a_v|^2\le C(b\gamma)^{-1}E,\\
 \sum_e\kappa_e(p_e-q_e)^2\le E,
 \end{gathered}
\end{equation}
where the middle sum runs over both endpoints of every edge, each paired
with its adjacent piece. The lemma requires neither that any edge be an
isolated leaf of the dual graph, nor a bound on the depth or on the
areas of its branches.
\end{lemma}

\begin{proof}
These are the first inequality of \eqref{eq:tree-trace} and the two
inequalities of \eqref{eq:tree-mass} in the proof of
Lemma~\ref{lem:super-comparison}; we point to the steps of that proof
which establish them and check their hypotheses. For $b$ small, the
standard half-collars of the cut geodesics extend at least one unit
beyond the trimming radius $r_e=\operatorname{arcosh}(b/\ell_e)$, so the
outer unit bands of distinct edges are disjoint; this uses only
$\ell_e\le\delta<b$ and the collar lemma. Each removed half-collar has
area at most $b$, and if $\chi_v=-\chi(P_v)>0$ and $d_v$ is the number
of boundary geodesics of $P_v$, then $d_v\le3\chi_v$, whence
$m_v=(1+O(b))2\pi\chi_v$ and $d_v\le Cm_v$. The Poincar\'e step extends
$f|_{B_v}$ harmonically into the removed half-collars with a free
endpoint at the waist; the Fourier-mode comparison on the disjoint outer
unit bands bounds the extension energy by three times the energy of
$f$ on those bands, and the Neumann gap of $P_v$ then gives
$\int_{B_v}|f-a_v|^2\le3\gamma^{-1}\int_{B_v}|\nabla f|^2$. The trace
step is the one-dimensional trace inequality on each outer band,
$|p_e-a_v|^2\le Cb^{-1}\int_{\mathrm{band}}(|f-a_v|^2+|\nabla f|^2)$,
summed over the disjoint bands and combined with the Poincar\'e bound;
this is the middle inequality of \eqref{rev:eq:cut-tree-basic}. The
collar step is the weighted Hardy inequality
$\int_0^{r_e}z^2\ell\cosh r\,dr\le4\int_0^{r_e}(z')^2\ell\cosh r\,dr$
for $z(r_e)=0$ on each half-collar, applied to the angular zero mode
after subtracting its outer average, together with the angular
Poincar\'e inequality for the nonzero modes; it bounds the mass of a
full collar $C_e$ by $2b(p_e^2+q_e^2)+C\int_{C_e}|\nabla f|^2$. Summing
over pieces and collars and using $d_v\le Cm_v$ gives the first
inequality of \eqref{rev:eq:cut-tree-basic}, and radial Cauchy--Schwarz
after angular averaging gives the last. None of these steps refers to
the position of an edge in the dual graph, to the depth, or to the
branch areas; those hypotheses of Lemma~\ref{lem:super-comparison}
enter only in its subsequent contraction step, which is replaced below.
\end{proof}

For real $f\in H^1(X)$ write $E=\int|\nabla f|^2$, $a_v$ for its
trimmed-piece means, and $p_e,q_e$ for its outer collar averages.
On the event of \eqref{eq:dense-piece-gap},
Lemma~\ref{rev:lem:cut-tree-basic} applies with the parameters
\eqref{eq:dense-proof-cutoffs} and the gap $\gamma$ of
\eqref{eq:dense-piece-gap}, although the retained edges are not
isolated leaves and the depth is only bounded by $D\le C_0L$. It gives
\begin{equation}\label{eq:dense-continuum-basic}
 \begin{gathered}
 \|f\|^2\le(1+Cb)\sum_vm_va_v^2+C\gamma^{-1}E,\qquad
 \sum_{\mathrm{endpoints}}|p_e-a_v|^2\le C(b\gamma)^{-1}E,\\
 \sum_e\kappa_e(p_e-q_e)^2\le E.
 \end{gathered}
\end{equation}
Let $Jf$ be the vector of high-component mass means on the retained
star. Anchor each component at its vertex nearest the original root.
Cauchy--Schwarz along paths of length at most $D$, followed by
summation over masses, bounds its variance by the edge jumps weighted
by at most $DW_e$. Split each jump into the collar-average difference
and its two trace errors. On high edges $W_e<\kappa_e/\tau$, and
on every edge $W_e\le CK$. Replacing the anchor by the component
mass mean decreases variance. Consequently
\begin{equation}\label{eq:dense-component-variance}
 \sum_H\sum_{v\in H}m_v|a_v-a_H|^2
       \le CD\left(\tau^{-1}+\frac K{b\gamma}\right)E.
\end{equation}
Here and below, as in the proof of Lemma~\ref{lem:super-comparison},
$M_\star$ and $Q_\star$ denote the mass and energy forms of the weighted
star: for a vector $y=(y_H)$ indexed by the high components,
\[
 M_\star(y)=\sum_HM_Hy_H^2,\qquad
 Q_\star(y)=\sum_{e\ \mathrm{retained}}\kappa_e\,(y_{H_e^+}-y_{H_e^-})^2,
\]
where $M_H=\sum_{v\in H}m_v$ is the trimmed mass of $H$ and
$H_e^\pm$ are the two components joined by the retained edge $e$; the
star eigenvalues and the counting function $N_\star^+$ are those of
the pencil $Q_\star-\mu M_\star$, with the zero eigenvalue of the
constant vector excluded.
If $E\le\Lambda\|f\|^2$, orthogonal variance decomposition yields
$M_\star(Jf)\ge a_\Lambda\|f\|^2$, where
\begin{equation}\label{eq:dense-map-mass}
 a_\Lambda=\frac{1-C\Lambda/\gamma}{1+Cb}
              -CD\Lambda\left(\tau^{-1}+\frac K{b\gamma}\right).
\end{equation}

The full tree energy of the piece means satisfies
\begin{equation}\label{eq:dense-tree-energy-loss}
 Q_{\mathrm{tree}}(a)\le\mathcal A E,\qquad
       \mathcal A=C\left(1+\frac\delta{b\gamma}\right).
\end{equation}
This follows from $\kappa_e\le C\delta$ and the summed trace bound.
The factor $\mathcal A$ can diverge and must be kept. On retained
edges $\kappa_e\le\tau W_e\le C\tau K$, so their collar-to-piece
errors have squared energy norm at most $C\tau K E/(b\gamma)$.
Their piece-to-component errors have squared norm at most
$\mathcal LQ_{\mathrm{high}}(a)\le\mathcal A\mathcal L E$
by the definition of the Gram matrix. Minkowski gives
\begin{equation}\label{eq:dense-map-energy}
 Q_\star(Jf)\le BE,\qquad
 B=\left(1+\sqrt{\frac{C\tau K}{b\gamma}}
                        +\sqrt{\mathcal A\mathcal L}\right)^2.
\end{equation}
Whenever $a_\Lambda>0$, $J$ is injective on the entire spectral
subspace $[0,\Lambda]$. Conversely, extend each star value constantly
over its high component and high collars, and interpolate radially
harmonically on retained collars. This injective trial map has energy
exactly $Q_\star$ and mass at least $M_\star$. Min--max in both
directions proves, for $\Lambda<1/4$,
\begin{equation}\label{eq:dense-count-sandwich}
 N_\star^+(\Lambda)\le N_X^+(\Lambda)
                 \le N_\star^+(B\Lambda/a_\Lambda).
\end{equation}
Both maps preserve constants; include zero and subtract exactly one
after the comparison. Thus all actual $L^2$ eigenfunctions, including
residual modes, are counted.

For fixed $T$ and $\Lambda\le T/S$, the deterministic errors satisfy
\begin{align*}
 b+\delta/b&=O(L^{-1/2}),&
 \Lambda/\gamma&\le C_TL^4/(\epsilon^2S),\\
 D\Lambda/\tau&\le C_T/L,&
 D\Lambda K/(b\gamma)&\le C_TL^{13/2}/(\epsilon^3S),\\
 \tau K/(b\gamma)&\le CL^{15/2}/(\epsilon^3S).
\end{align*}
Equation~\eqref{eq:dense-access} also gives
\[
 \mathbb E(\mathcal A\mathcal L)
       \le CL^{17/2}/(\epsilon^3S)+o(1),
\]
using $\epsilon^{-2}\sqrt\eta\asymp1/(\epsilon^3S)$.
Thus each of the five error terms
$\Lambda/\gamma$, $D\Lambda/\tau$, $D\Lambda K/(b\gamma)$,
$\tau K/(b\gamma)$ and $\mathbb E(\mathcal A\mathcal L)$ is at most
$C_TL^{17/2}/(\epsilon^3S)+o(1)$, the largest power of $L$ coming from
the last term; the two auxiliary requirements $\eta R^4\to0$ and
$R^2/g\to0$ used in \eqref{eq:dense-retained-mass} are
$L^8\eta\to0$ and $L^4/g\to0$, which also follow from
$L^{17/2}/(\epsilon^3S)\to0$ because
$\eta\asymp(n\epsilon^{3/2})^{-1}\ll(\epsilon^3S)^{-1}$ and
$g\gg\epsilon^3S$ in the dense regime. The hypothesis
$L^{10}/(\epsilon^3S)\to0$ of \eqref{eq:dense-persistence-range} is
therefore stronger than the proof needs; the exponent $10$ was chosen
for simplicity of the statement, and we do not claim that $17/2$ is
sharp either.
Condition~\eqref{eq:dense-persistence-range} makes all these errors
vanish. Hence $a_\Lambda\to1$ uniformly in probability and $B\to1$.
Controlling $\mathcal L$ alone, without the product
$\mathcal A\mathcal L$, would not justify this continuum comparison.

\subsection{Identification of the limit and the gap}

The positive star eigenvalues are close to the leaf rates, which are
the rescaled cap scores; the cap process of Lemma~\ref{lem:dense-caps}
then gives the Poisson limit \eqref{eq:dense-persistence-PPP} through
the counting sandwich \eqref{eq:dense-count-sandwich}, and the same
first score governs the Cheeger constant, by the cut exclusion
\eqref{eq:dense-cut-exclusion} and the comparison of
\cite[Proposition~13]{SW} between the geodesic and the analytic Cheeger
constants.

The total retained leaf mass divided by the central mass is
$O_{\mathbb P}(R^2/g)=o_{\mathbb P}(1)$ by
\eqref{eq:dense-retained-mass}. If the leaf rates are $d_i=\kappa_i/m_i$
and the ordered positive star eigenvalues are $\mu_i$, elimination
of the constant mode gives
$d_i\le\mu_i\le(1+\sum m_i/M_0)d_i$ after ordering: the variance
mass form is $\operatorname{diag}(m_i)-mm^{\mathsf T}/(M_0+\sum m_i)$,
between $M_0/(M_0+\sum m_i)$ and one times $\operatorname{diag}(m_i)$.
On the non-nesting event, the leaf mass is $W_e$, and
\[
 PS\frac{\kappa_e}{W_e}
      =\frac{S\ell_e}{k_e-1}(1+O(b+\delta/b)).
\]
Every bounded score and compact positive mark window is eventually
below $\delta$, since $1/(\epsilon S)=o(\delta)$, and is retained
since $R\to\infty$. Indeed a cap with mark $\epsilon(k-1)\le M$ and
score $S\ell/(k-1)\le t$ has $\ell\le tM/(\epsilon S)$, and
$1/(\epsilon S)=o(\delta)$, with $\delta=L^{-1}$, means
$L^2/(\epsilon S)^2\to0$.
By \eqref{eq:dense-scales}, $(\epsilon S)^2\asymp g^{3/2}n^{-1/2}$,
so $(\epsilon S)^{-2}\asymp n^{1/2}g^{-3/2}\asymp\eta$, and the
requirement is equivalent to $\eta L^2\to0$; this follows from
$\eta R^4=\eta L^8\to0$, which was already used after
\eqref{eq:dense-retained-mass}. Lemma~\ref{lem:dense-caps} applies because
\eqref{eq:dense-persistence-range} implies $g\ll n\ll g^3$.
Its extension to all marks follows from \eqref{eq:PC3} and the
probability weights $(k-1)^2w_k/M_2$, whose scaled marks converge
to $\operatorname{Gamma}(1/2,1)$. The single-curve expectation also
bounds the count below score $t$ by $Ct^2$, giving tightness at zero.
The star comparison and \eqref{eq:dense-count-sandwich}, uniformly
on bounded $1/S$ windows, now prove joint convergence of counts at
continuity thresholds, hence \eqref{eq:dense-persistence-PPP}.
Increasing the window after the genus limit gives every fixed list
of ordered eigenvalues and the gap law. Finally $A_g\sim PS$.

For the Cheeger assertion, recall that the geodesic Cheeger constant
of~\cite{SW} is an infimum over separating geodesic multicurves with
both sides connected; allowing disconnected sides does not lower it.
Indeed, by the graph-averaging reduction in the proof of
Lemma~\ref{lem:relative-geodesic}, every separating multicurve of ratio at
most $u$ contains a subfamily cutting $X$ into two connected sides with
ratio at most $u$, so only the cuts of \eqref{eq:dense-cut-exclusion}
and single planar waists need be considered. Apply
\eqref{eq:dense-cut-exclusion} at $u=t/(2\pi S)$; this is admissible
since $\epsilon S\to\infty$, and it shows that competing cuts have
probability $o(1)$.
All planar waists, including long ones, have the tail bound
\[
 C_t\sum_{\epsilon(k-1)>M}
       \frac{(k-1)^2w_k(x)}{M_2(x)}e^{2\pi u(k-1)}.
\]
It vanishes as $M\to\infty$ after the limit: $u/\epsilon\to0$, and
the exponential moment at any fixed $d\epsilon$, $0<d<1$, is bounded
by the corresponding $M_2$ ratio, tending to $(1-d)^{-1/2}$.
On bounded marks the planar side is the smaller side, of area
$2\pi(k-1)$. Thus the geodesic Cheeger minimum is governed by the
same first score as the spectrum. The comparison
$H_{\mathrm{geo}}/(1+H_{\mathrm{geo}})\le\mathfrak h\le H_{\mathrm{geo}}$
of~\cite[Proposition~13]{SW} transfers this to $\mathfrak h$.
Consequently $2\pi S\mathfrak h$ and $PSG$ match the same positive
first Poisson point up to $o_{\mathbb P}(1)$, which proves
$G/\mathfrak h\to1/\pi$ and completes the theorem.

\section{Excluding cusp forms at the dense microscopic scale}
\label{app:cuspidal-exclusion}

We prove Proposition~\ref{prop:dense-cusp-zero}, the lower bound on the
cusp-zero form which the proof of Corollary~\ref{cor:dense-residual} in
Section~\ref{sec:dense-spectrum} uses to exclude cuspidal
eigenfunctions at the dense microscopic scale $1/S$. The proposition
holds throughout the dense regime $n/g\to\infty$ and does not require
the persistence condition \eqref{eq:dense-persistence-range}.
Fix a sufficiently small absolute circumference $b>0$, so that the
standard cusp tails $U_{c,b}$ are embedded and pairwise disjoint.
Let $\mathcal H_b(X)$ be the subspace of $H^1(X)$ with
$\int_0^1 f(x+iy)\,dx=0$ for almost every $y\ge b^{-1}$ in every
normalized cusp coordinate, and put
\[
 \eta_b(X)=\inf_{0\ne f\in\mathcal H_b(X)}
              \frac{\int_X|\nabla f|^2}{\int_X|f|^2}.
\]
This is the form domain of a pseudo-Laplacian, as in
\cite[Sections~1 and~4]{BMP}. It contains every cuspidal eigenfunction,
but is not identified with the cuspidal spectral subspace.

\begin{proposition}[A lower bound for the cusp-zero form]
\label{prop:dense-cusp-zero}
There are absolute constants $C,t_0,K>0$ such that, for all $g\ge K$,
all $n\ge Kg$ and all $0<t<t_0$,
\begin{equation}\label{eq:dense-cusp-zero-tail}
 \mathbb P^{\mathrm{WP}}_{g,n}\{\eta_b(X)\le t\}
       \le C\left(\frac{n^2}{\sqrt g}t^3+\frac{t}{\sqrt n}\right).
\end{equation}
In particular, $S\eta_b\to\infty$ in probability along every sequence
with $g\to\infty$ and $n/g\to\infty$, with $S$ as in
\eqref{eq:dense-parameters}.
\end{proposition}

The proof is written along an arbitrary sequence with $g\to\infty$ and
$n/g\to\infty$ and gives \eqref{eq:dense-cusp-zero-tail} with constants
independent of the sequence; the stated form follows, since a sequence
of pairs $(g,n)$ violating it for every choice of $K$ would be such a
sequence.

The proof retains a penalty for every cusp in a geodesic Cheeger
comparison. A region with small penalized ratio must have few cusps
relative to its genus and boundary count. Mirzakhani integration then
bounds the probability of such a region; we sum all possible
complementary topologies, including disconnected complements and
nonseparating curves. The volume inputs are the local limit theorem
(Proposition~\ref{prop:local-limit}), the all-marking envelope
\eqref{eq:gevrey-envelope} with the one-marking lower bound of
Proposition~\ref{prop:prefactor}, and the boundary-volume bound
\eqref{eq:PS4}; the endpoint quantities $Z,\rho$ are those of
Appendix~\ref{app:dense-comparison}. The first subsection proves the
deterministic comparison \eqref{eq:cusp-penalty-cheeger} and the
topological constraints \eqref{eq:cusp-penalty-topology}; the second
parametrizes the complementary topologies and bounds the volume and
length factors of the integration formula; the third cancels the
endpoint quantities against the parent and sums over the complementary
genera; the fourth sums over boundary counts, treats the two
exceptional sectors, and concludes.

\subsection{The deterministic comparison}

We bound $\eta_b$ below by a penalized geodesic Cheeger constant and
record the topological constraints on a region with small penalized
ratio.

Put $\beta=b/2$ and define
\[
 H_\beta(X)=\inf_A\frac{\ell(\partial A)+\beta a(A)}{|A|},
\]
where $a(A)$ is the number of original cusps of $A$. The infimum includes
$A=X$ and every connected proper geodesic subsurface with distinct boundary
geodesics and stable complementary components. We claim
\begin{equation}\label{eq:cusp-penalty-cheeger}
 \eta_b(X)\ge c\min\{1,H_\beta(X)^2\}.
\end{equation}

The angular Poincar\'e inequality on the cusp tails gives, for
$f\in\mathcal H_b(X)$,
\[
 \sum_c\|f\|_{L^2(U_{c,b})}^2
       \le \frac{b^2}{4\pi^2}\int_X|\nabla f|^2.
\]
Choose a cutoff $\chi$ equal to one outside these tails and zero in
the smaller tails of circumference $\beta$. Its gradient is uniformly
bounded on the intervening annuli, whose width is $\log2$. On the compact
surface $M=X\setminus\bigcup_c U_{c,\beta}^{\circ}$, the function
$\chi f$ has zero trace at every boundary and satisfies
\[
 \int_M|\nabla(\chi f)|^2\le C\int_X|\nabla f|^2,
 \qquad
 \|\chi f\|_2^2\ge\|f\|_2^2-
                       \frac{b^2}{4\pi^2}\int_X|\nabla f|^2.
\]
Consequently $\eta_b\ge c\min\{1,\lambda_D(M)\}$. The all-Dirichlet
Cheeger inequality~\cite{Cheeger} gives $\lambda_D(M)\ge h_D(M)^2/4$, where
$h_D(M)$ is the infimum of perimeter divided by area over smooth
domains $\Omega$ with $\overline\Omega\subset M\setminus\partial M$.
Connected domains suffice.

We show $h_D(M)\ge H_\beta/(1+H_\beta)$. Given such a connected domain
$\Omega$, first fill all complementary disks and unpunctured annuli.
This decreases perimeter and increases area. The filled pieces cannot
contain any removed cusp tail: a tail is connected, disjoint from
$\partial\Omega$, and contains an end. Filling an annulus can increase
the genus of $\Omega$; this is allowed and removes duplicate parallel
geodesic representatives, including the two sides of a nonseparating curve.

An elementary domain has area at most its perimeter. Indeed it lifts to the
plane or a cyclic cover, where a Busemann field, or the field
$\tanh(r)\partial_r$ on a hyperbolic cylinder, has divergence one and
norm at most one. Otherwise the filled domain is incompressible, and
straightening its nonperipheral boundaries gives a geodesic subsurface
$A$ with distinct boundaries; peripheral boundaries give the cusps of $A$.
In the cover corresponding to $\pi_1(\Omega)$, this subsurface is the
convex core, with its cusp ends included. On each exterior funnel use
$\tanh(r)\partial_r$, extended by zero on the core. Its divergence is
one outside the core and its boundary flux at the core vanishes. Hence
\[
 |\Omega|\le |A|+\ell(\partial\Omega).
\]
Every nonperipheral boundary loop is at least as long as its geodesic
representative. A peripheral loop outside $U_{c,\beta}$ lifts to a path
from $z$ to $z+1$ with $y\le\beta^{-1}$, so has length at least $\beta$.
Thus $\ell(\partial\Omega)\ge\ell(\partial A)+\beta a(A)
\ge H_\beta|A|$. Combining the last two inequalities proves the claimed
bound for $h_D$ and then \eqref{eq:cusp-penalty-cheeger}.

For $n\ge g$, the whole-surface ratio is at least $\beta/(6\pi)$.
It does not contribute to $H_\beta\le u$ for small $u$. A proper
subsurface of genus $h$, with $r$ boundaries and $a$ cusps, has area
$2\pi(d+a)$, where $d=2h-2+r$. If its penalized ratio is at most $u$,
then, for a sufficiently small absolute $u_0$ and $0<u<u_0$,
\begin{equation}\label{eq:cusp-penalty-topology}
 d\ge1,\qquad a\le\vartheta_u d,\qquad
 \ell(\partial A)\le4\pi u d,\qquad
 \vartheta_u=\frac{2\pi u}{\beta-2\pi u}\le1.
\end{equation}
In particular $h=0$ requires $r\ge3$.

\subsection{All complementary topologies}

We parametrize every complement of the distinguished subsurface and
bound the corresponding Mirzakhani integral: first its summed volume
factor with the exact symmetry factor, then its length and twist
integrals.

Write $H=2g-2$, $P=2\pi^2$, and choose the saddle
$(2g-1)\Psi(x)=n-1$. Retain $\epsilon=1-x/x_c$. The parent coefficient
standard deviation satisfies $\sigma\asymp\sqrt g/\epsilon\asymp n/\sqrt g$.
Among the complementary components of the distinguished subsurface $A$,
let $l$ be planar with one boundary and $v$ planar with two boundaries.
They have at least two cusps and at least one cusp, respectively. The
remaining $w$ components have genera $g_i$, boundary counts $p_i\ge1$,
and $d_i=2g_i-2+p_i\ge1$. Put $p_* =\sum_i p_i$, with $p_*=0$ if $w=0$.
This list is exhaustive. Every complementary component has at least
one boundary. A planar one-boundary component must have at least two
cusps, since its boundary is geodesic and nonperipheral. A planar
two-boundary component must have a cusp, since unpunctured annuli
were filled in the deterministic comparison. Every other component
has $d_i\ge1$. The dual graph has the distinguished vertex $A$ and
one vertex for each complementary component; every edge joins $A$
to a complementary vertex. Parallel edges are allowed and account
for nonseparating boundary curves. No self-edge is needed, because
$A$ is embedded and its boundary geodesics are distinct. The genera,
cusp labels and this graph determine the mapping-class orbit of the
cutting configuration; the integration formula sums its embeddings.
Euler characteristic and boundary counting give
\begin{equation}\label{eq:cusp-penalty-euler}
 \begin{split}
 r&=l+2v+p_*,\\
 g&=h+\sum_i g_i+r-(l+v+w),\\
 D:=\sum_i d_i&=H-d+l=2g-2h-2v-p_*.
 \end{split}
\end{equation}
Multiple boundaries to the same component are retained. Their individual
curves need not separate $X$.

Set $F_{q,p}=\Phi_q^{(p)}$. The coefficient convention is
\[
 F_{q,p}(z)=\sum_a
       \frac{V_{q,a+p}}{a!P^{3q-3+a+p}}z^a,
\]
with only stable terms. Write $F_{h,r}^{\le}$ for the truncation to
$a\le\vartheta_u d$. For an ordered list of the remaining $w$ vertex
types, the summed zero-boundary volume factor in Mirzakhani integration
is at most
\begin{equation}\label{eq:cusp-penalty-allocation}
 \frac{n}{P^r\mathcal D}\,
 \frac{[z^n]\bigl(F_{h,r}^{\le}(z)\Phi_0'(z)^l
                   \Phi_0''(z)^v\prod_iF_{g_i,p_i}(z)\bigr)}
      {[z^{n-1}]\Phi_g'(z)},\qquad
 \mathcal D=l!v!w!2^v\prod_i p_i!.
\end{equation}
To see the symmetry factor, temporarily label the $r$ central boundaries.
For an ordered list of complementary vertices, their allocation gives
$r!/(2^v\prod_i p_i!)$. Dividing by the auxiliary $r!$ labels
and by $l!v!w!$ gives $\mathcal D^{-1}$; the remaining vertex
types are summed in order, so this also handles repeated types.
The central vertex remains distinguished; labeled cusps are already
summed by the coefficient product. The direction of the symmetry
bookkeeping is the following inequality. Let $O$ range over the
mapping-class orbits of pairs $(\Gamma,A)$, where $\Gamma=\partial A$
and the cusp labels are fixed, with the given $(h,r,a)$ and
complementary data, and let $\mathrm{Sym}(O)$ be the group of
label-preserving mapping classes stabilizing $(\Gamma,A)$, modulo those
fixing every component of $\Gamma$. Mirzakhani's integration formula
gives the expected number of such pairs as
$\sum_O|\mathrm{Sym}(O)|^{-1}\int_O$, and we claim
\begin{equation}\label{rev:eq:cusp-penalty-symmetry}
 \sum_O\frac1{|\mathrm{Sym}(O)|}\int_O
 \le\frac1{\mathcal D}\sum_{\mathrm{ordered}}\int,
\end{equation}
where the right side sums the same integrand over ordered lists of
complementary vertex types with the central boundaries allocated as
above, which is the quantity in \eqref{eq:cusp-penalty-allocation}.
Indeed, the group $G=S_r\times S_l\times S_v\times S_w$, of order
$r!\,l!\,v!\,w!$, acts on the ordered labeled presentations, its
orbits are the orbits $O$, and the stabilizer of a presentation of $O$
has order $2^v\prod_ip_i!\prod(\mathrm{mult})!$: it may exchange the
two boundary labels of each two-boundary planar vertex, permute the
$p_i$ parallel labels attached to the same vertex, and permute the
cusp-free $w$-vertices of equal type, whose multiplicities are the
$\mathrm{mult}$; vertices carrying cusps are distinguished by their
disjoint label sets. Each such stabilizer element is realized by a
label-preserving mapping class of $X$ that preserves $A$ (a
homeomorphism permuting the boundary curves of $A$ attached to one
complementary vertex, a half twist exchanging the two boundaries of a
planar two-boundary vertex, or a homeomorphism exchanging two
isomorphic cusp-free branches together with their attaching curves).
Distinct stabilizer elements induce distinct permutations of the
components of $\Gamma$, so
$|\mathrm{Sym}(O)|\ge2^v\prod_ip_i!\prod(\mathrm{mult})!$. Hence
$|\mathrm{Sym}(O)|^{-1}\le|G|^{-1}\cdot\#\{\text{presentations of }O\}$,
and summing over $O$ gives \eqref{rev:eq:cusp-penalty-symmetry} after
dividing the allocation count by $r!$. A central subsurface isomorphic
to a cusp-free complementary branch, but distinguished as $A$, is
counted as separate data in the sum over
$(h,r,a,\text{complement})$, so no further identification is needed.
For every intact one-holed torus among the cut components, the
half-twist factor is already included in its orbifold volume.
No additional factor is needed.
For example, a cusp-free central pair of pants, two two-cusp leaves and
one genus-$g$ complement give
$n f_{n-4}/(8P^3 f_{n-1})=(n)_4V_{g,n-3}/(8V_{g,n})$,
where $f_j=[z^j]\Phi_g'$. This checks both leaf symmetries and the
parent coefficient shift.

By \eqref{eq:PS4}, the boundary-volume cost is at most
$\exp(\sum_i\ell_i)$, since each cut appears twice. The length and
twist integrals over \eqref{eq:cusp-penalty-topology} are bounded by
\begin{equation}\label{eq:cusp-penalty-length}
 I_{r,d}(u)=e^{4\pi u d}\frac{(4\pi u d)^{2r}}{(2r)!}.
\end{equation}
This is the full length-simplex integral, valid for every $r$ and $d$.

\subsection{Endpoint cancellation and the genus sum}

We evaluate the volume factor at the saddle, cancel the endpoint
quantities $Z$ and $\rho$ against the parent, and sum over the genera
of the complementary components.

Proposition~\ref{prop:local-limit} gives a parent point mass at least
$c/\sigma$. Bound the numerator coefficient by its value at $x$.
Thus \eqref{eq:cusp-penalty-allocation} is at most
\begin{equation}\label{eq:cusp-penalty-evaluation}
 \frac{Cn\sigma}{P^r x\mathcal D}
 \frac{F_{h,r}^{\le}(x)\Phi_0'(x)^l\Phi_0''(x)^v
                         \prod_iF_{g_i,p_i}(x)}{\Phi_g'(x)}.
\end{equation}
We allow the loss of $\sigma$; no estimate for a large complementary
component or a growing-marking local limit theorem is needed.

Both $\Phi_0'$ and $\Phi_0''$ stay bounded at $x_c$. At $y=xe^{-1}$,
coefficient positivity and \eqref{eq:gevrey-envelope} give
\begin{equation}\label{eq:cusp-penalty-central}
 F_{h,r}^{\le}(x)\le e^{\vartheta_u d}F_{h,r}(y)
          \le C^d4^{-h}\Gamma(d+1).
\end{equation}
Here $y$ stays in a fixed compact subinterval of $(0,x_c)$ and
$r\le d+2$. The planar case $h=0,r\ge3$ also follows from Cauchy's
estimate for $\Phi_0''$ on this compact interval.

As in Appendix~\ref{app:dense-comparison}, put
$\rho(x)=r(Y(x))$ and $Z(x)=a(Y(x))\rho(x)^{3/2}A(Y(x))$.
The all-marking envelope, normalization \eqref{eq:TL9}, and
Proposition~\ref{prop:prefactor} at one marking give
\begin{align}
 F_{g_i,p_i}(x)&\le C^{p_i}4^{-g_i}
       \frac{\Gamma(d_i)}{\sqrt{d_i}}Z^{-d_i}\rho^{p_i/2},
          \label{eq:cusp-penalty-envelope}\\
 \Phi_g'(x)&\ge c4^{-g}
       \frac{\Gamma(H+1)}{\sqrt{H+1}}Z^{-(H+1)}\rho^{1/2}.
          \notag
\end{align}
The first bound includes $g_i=0,p_i\ge3$; the square-root singularity
of $\Phi_0''$ gives the same bound with the harmless factor $C^{p_i}$.
Keep $Z$ and $\rho$ exact until their powers cancel. Set
$q=2(h+v)+p_*$, so $d=l+q-2$ and $D=H-q+2$. The remaining factor is
$Z^{q-1}\rho^{(p_*-1)/2}$. If $w\ge1$, then $p_*\ge1$, $q\ge1$,
and $Z\asymp\epsilon^{5/4}$, $\rho\asymp\epsilon^{1/2}$ imply
\begin{equation}\label{eq:cusp-penalty-endpoint}
 Z^{q-1}\rho^{(p_*-1)/2}
       \le C^{h+r}\epsilon^{(5/2)(h+v)+(3/2)(p_*-1)}.
\end{equation}
All remaining powers of $\epsilon$ are nonnegative. The powers of four
leave $4^{g-h-\sum_i g_i}=4^{r-l-v-w}\le4^r$.

For fixed $p_i$, sum the complementary genera by summing over all
positive $d_i$ with sum $D$. Uniformly in $D\ge w\ge1$,
\begin{equation}\label{eq:cusp-penalty-convolution}
 \sum_{d_1+\cdots+d_w=D}\prod_i\frac{\Gamma(d_i)}{\sqrt{d_i}}
       \le C^w\sqrt{\frac wD}\,\Gamma(D-w+1).
\end{equation}
Indeed $\prod_i\sqrt{d_i}\ge\sqrt{D/w}$, and the substitution
$x_i=d_i-1$ reduces the remaining sum to repeated use of
$\sum_{a=0}^N a!(N-a)!\le3N!$. Parity and genus restrictions only
reduce the sum. The extra square root satisfies
$\sqrt{(H+1)w/D}\le C^{q+w}$: if $D\ge(H+1)/2$ use $\sqrt{2w}$;
otherwise $q>H/2$ absorbs the factor $H+1$.

For $q\ge2$, put $a_0=q-2$ and $b_0'=D-w=H-a_0-w\ge0$. Then
\begin{equation}\label{eq:cusp-penalty-factorials}
 \frac{\Gamma(d+1)\Gamma(D-w+1)}{l!\Gamma(H+1)}
 =\frac{(l+a_0)!b_0'!}{l!H!}
 \le2^{l+a_0}\frac{a_0!b_0'!}{H!}\le2^{l+a_0}.
\end{equation}
For $q=1$, necessarily $h=v=0$, $p_*=w=1$ and $l\ge2$;
the same ratio is exactly $1/l$. Combining these bounds, the sum of
\eqref{eq:cusp-penalty-evaluation} over complementary genera, for fixed
$h,l,v,w$ and ordered $p_i$, is at most
\begin{equation}\label{eq:cusp-penalty-genus-sum}
 Cn\sigma C^{r+h}
           \epsilon^{(5/2)(h+v)+(3/2)(p_*-1)}.
\end{equation}
The unused factorial divisors may now be dropped. The number of choices
of $l,v,w$ and compositions of $p_*$ at fixed $r$ is at most $C^r$.

\subsection{Boundary counts and the two exceptional sectors}

It remains to sum over boundary counts. For $r\ge3$ the length integral
supplies enough smallness; the two exceptional sectors, $r\le2$ with
$w\ge1$ and $w=0$, require the exact factorial ratios.

For $r\ge3$ and $w\ge1$, since $d\le r+2h$, \eqref{eq:cusp-penalty-length} gives
\[
 I_{r,d}(u)\le(Cu^2)^r C^h.
\]
Indeed $(d/r)^{2r}\le e^{4h}$, and the boundary exponential costs
at most $C^{r+h}$. In \eqref{eq:cusp-penalty-genus-sum}, sum $h$
using $(C\epsilon^{5/2})^h$, then sum the remaining topologies and
all $r\ge3$. For small $u,\epsilon$ this gives at most $Cn\sigma u^6$.

For $r=1,2$, the central genus satisfies $h\ge1$. With $w\ge1$, one
has $v=0$ and only the following possibilities. We retain the exact
factorial ratio before the coarse bound \eqref{eq:cusp-penalty-factorials}:
\[
\begin{array}{c|c}
 (r,l,p_*,w)&\text{factorial ratio}\\ \hline
 (1,0,1,1)&(2h-1)!(H-2h)!/H!\\
 (2,1,1,1)&(2h)!(H-2h)!/H!\\
 (2,0,2,1)&(2h)!(H-2h-1)!/H!\\
 (2,0,2,2)&(2h)!(H-2h-2)!/H!
\end{array}
\]
Only nonnegative factorial arguments and admissible genera occur.
Except at $h=H/2$ in the first two rows, these ratios are at most
$C/H^2$, by the elementary lower bound for binomial coefficients
away from their endpoints. At the two exceptions they are $1/H$
and $1$. After absorbing all other factors, including the powers of
$h$ in the length integral, into $z^h$, with
$z=C\epsilon^{5/2}<1/4$, their weighted sum is still at most
$Cz/H^2$: the endpoint correction satisfies
$H^2z^{H/2-1}\le C$. Thus this sector contributes at most
\begin{equation}\label{eq:cusp-penalty-small-boundary}
 Cn\sigma\epsilon^{5/2}u^2/H^2\le Cu^2/\sqrt n.
\end{equation}
This includes the case in which almost all genus lies in the central
subsurface.

Finally suppose $w=0$. Then $p_*=0$, $g=h+v$, and $d=H+l$.
The factorial ratio, including $l!$, is
$\binom{H+l}{l}\le2^{H+l}$. The exact endpoint factor is now
$Z^{H+1}\rho^{-1/2}\le C^g\epsilon^{(5/2)g-3/2}$.
Since $d\le r+2g$, the length integral is at most $(Cu^2)^rC^g$.
Summing $l,v$ with $r=l+2v\ge1$ gives
\begin{equation}\label{eq:cusp-penalty-no-positive-complement}
 Cn\sigma C^g\epsilon^{(5/2)g-3/2}u^2\le Cu^2/\sqrt n.
\end{equation}
For the last inequality, use $n\sigma\asymp g^{3/2}\epsilon^{-2}$;
after division by $n^{-1/2}$, the coefficient is bounded by
$C^g g^2\epsilon^{(5/2)g-4}$, uniformly bounded for small $\epsilon$
and $g\ge3$. This sector also includes genus arising from planar
two-boundary complementary components.

Every topology is included in these sums. A union bound therefore gives
\[
 \mathbb P\{H_\beta\le u\}
       \le C\left(\frac{n^2}{\sqrt g}u^6+\frac{u^2}{\sqrt n}\right).
\]
Infima can be handled by enlarging $u$ by a fixed factor arbitrarily
close to one. Combining this with \eqref{eq:cusp-penalty-cheeger}
proves \eqref{eq:dense-cusp-zero-tail}. Finally
$n^2/(\sqrt g S^3)\asymp(g/n)^{1/4}\to0$, so for each fixed $T$,
\[
 \mathbb P\{S\eta_b\le T\}
       \le C\left(T^3(g/n)^{1/4}+\frac{T}{S\sqrt n}\right)\longrightarrow0.
\]
This completes the proof of Proposition~\ref{prop:dense-cusp-zero}.

\section{Trace estimates after summing cusp labels at the saddle}
\label{app:sublinear-trace}

This appendix proves the trace input to
Theorem~\ref{thm:fixed-energy-spectrum}, the fixed-energy spectral
density when $\sqrt g\ll n\ll g$; the spectral transfer itself is
carried out in Appendix~\ref{app:fixed-window-proofs}. The fixed-energy
law requires a trace remainder beyond the bounded critical envelopes of
Appendix~\ref{app:critical-trace}, which cover $n\le M\sqrt g$: absolute
factorial bounds for the denominator lose $\exp(Cn^2/g)$ there. We
instead sum cusp labels by coefficient extraction at the saddle, and
separate the omitted deficit into topology and additional cusps; the
positive convolution factor of Lemma~\ref{lem:convolution} then cancels
the local probability scale of Proposition~\ref{prop:local-limit}. The
appendix has three parts. Lemma~\ref{lem:sublinear-tail} records the
resulting graded tail bounds; a sequence of lemmas then restates, in
our notation, the inputs from Huang's trace expansion \cite{HuangGap}
together with the genus threshold under which each holds; and
Proposition~\ref{prop:sublinear-polynomial} inserts the graded bounds
into that expansion to obtain a polynomial trace estimate at sublinear
cusp density.

\begin{lemma}[A graded tail at sublinear cusp counts]
\label{lem:sublinear-tail}
There are absolute $C,n_0$ with the following property. Use the
carrier and pattern notation of \cite[Section~2.2]{HuangGap}: the
carrier is $(h_0,a_0,k)$, $m=2h_0+a_0+k-2\ge1$, a pattern has
$\sigma$ self-pairs and $b$ exterior blocks of degrees $k_i\ge1$,
and
\[
 \sum_i k_i=k-2\sigma,\quad
 \sum_i g_i=g-h_0-k+\sigma+b,\quad \sum_i a_i=n-a_0.
\]
Require $2g_i+a_i+k_i-2\ge1$, and put
\[
 W=\frac{n!}{a_0!\prod_i a_i!}
       \frac{V_{h_0,a_0+k}\prod_iV_{g_i,a_i+k_i}}{V_{g,n}},
 \qquad
 \delta=2g+n-3-\max_i(2g_i+a_i+k_i-3).
\]
For $r=R+1$, $t=r+m$, $n\ge n_0$ and $g\ge C(n+t^2)$,
\begin{equation}\label{eq:sublinear-tail}
 \sum_{\text{patterns}}\ \sum_{\delta\ge r} W
       \le g^{-r}(n+1)^{r+1}[Ct]^{Ct}.
\end{equation}
The boundary version has the additional factor
$(1+\Lambda)^k e^\Lambda$, where $\Lambda$ is the sum of carrier
slot lengths and the actual-edge Jacobian is included. The ordered
separating tail of \cite[Section~4.2]{HuangGap} satisfies
\begin{equation}\label{eq:sublinear-separating-tail}
 S_{g,n}^{\ge r}(L)\le g^{-r}(n+1)^{r+1}(Cr)^{Cr}e^L
 \quad\text{if }g\ge C(n+r^2),\ n\ge n_0.
\end{equation}
\end{lemma}

\begin{proof}
The proof has five steps: the coefficient ratio at the saddle, the
decomposition of a pattern and its deficit, the factorial estimates for
the genus sum, the envelope bounds, and the summation over the two
sectors of the deficit; the separating tail is treated at the end by
the same algebra.

\emph{Step 1: the coefficient ratio at the saddle.}
Put $P=2\pi^2$, $F_g=\Phi_g'$, and choose $x$ by
$(2g-1)\Psi(x)=n-1$. Thus $x\asymp n/g$ and, after increasing
$C$, $0<x\le x_0<y<x_c$ for fixed $x_0,y$. The uniform local
limit theorem gives
\begin{equation}\label{eq:sublinear-parent-mass}
 \frac{[z^{n-1}]F_g(z)x^{n-1}}{F_g(x)}\ge c n^{-1/2}
\end{equation}
when $n_0,C$ are sufficiently large. We also need a quantitative
version of Lemma~\ref{lem:point-mass}. The probabilities of zero
and one for $V(xz)/V(x)$ are at least $c$ and $cx$. Consequently
its characteristic function is bounded by $e^{-cx\theta^2}$ on
$[-\pi,\pi]$. If $G=V^{h'}R$, $R$ has nonnegative coefficients,
and $h'\ge g/(Cb)$, Fourier inversion gives
\[
 \sup_j\frac{[z^j]G(z)x^j}{G(x)}
 \le\min\{1,C(h'x)^{-1/2}\}\le C\sqrt{b/n}.
\]
Together with \eqref{eq:sublinear-parent-mass}, this implies
\begin{equation}\label{eq:sublinear-coefficient-ratio}
 \frac{[z^{n-a_0}]G(z)}{[z^{n-1}]F_g(z)}
       \le C\sqrt b\,x^{a_0-1}\frac{G(x)}{F_g(x)}.
\end{equation}

\emph{Step 2: the pattern and its deficit.}
Fix a pattern. Separate its genus allocations into $u$ genus-zero
degree-one vertices, $v$ genus-zero degree-two vertices, and $w$
vertices stable before adding cusps. Since the genera of the exterior
blocks sum to at least $g-(m+2)$, which is positive under the genus
hypothesis, some exterior block has positive genus. Thus $b=u+v+w$,
$w\ge1$, and
$h_i=2g_i-2+k_i\ge1$ on the stable vertices. Their total is
$H=2g-2h_0-k+u$. Choose a stable root $*$ maximizing $h_i$, breaking
ties by block order, and write
\begin{align*}
 K&=H-w,& \Delta&=2g-2-K=2h_0+k-u+w-2,\\
 j&=\sum_{\substack{i\ne *\\i\text{ stable}}}(h_i-1),&
 B&=a_0+2u+v,\\
 d&=m+b-1=\Delta+B-1,&
 A_{\rm ex}&=\sum_{i\ne *}a_i-2u-v\ge0.
\end{align*}
Here $\Delta\ge0$, since $k\ge u+v+w+2\sigma$, and hence
$B\le d+1$. The root need not maximize the full factorial index.
Its relative deficit is, however, exactly
\begin{equation}\label{eq:sublinear-deficit-split}
 \delta_*=2g+n-3-(h_*+a_*-1)=d+j+A_{\rm ex},
 \qquad \delta\le\delta_*.
\end{equation}
Thus $\delta\ge r$ implies $j+A_{\rm ex}\ge r-d$.

For fixed genus data put
\[
 B_{\rm out}(z)=
 \left(\frac{\Phi_0'(z)}{z^2}\right)^u
 \left(\frac{\Phi_0''(z)}z\right)^v
 \prod_{\substack{i\ne *\\i\text{ stable}}}\Phi_{g_i}^{(k_i)}(z).
\]
This has nonnegative coefficients, and its degree records $A_{\rm ex}$.
An upper bound for the required numerator is the sum of
\begin{equation}\label{eq:sublinear-tail-series}
 \Phi_{g_*}^{(k_*)}(z)z^{2u+v}
        [B_{\rm out}(z)]_{\deg\ge\max(0,r-d-j)}.
\end{equation}
The root series is untouched. Since $h_*\ge H/w\ge g/b$,
Lemma~\ref{lem:convolution} factors a common $V^{\lfloor g/b\rfloor}$
from this sum even after the outer cusp truncation. The exact
coefficient dictionary is
\begin{equation}\label{eq:sublinear-label-dictionary}
 \sum W=\frac{nP^{\sigma-k}v_0}{a_0!}
       \frac{[z^{n-a_0}]G(z)}{[z^{n-1}]F_g(z)},
 \qquad v_0=\frac{V_{h_0,a_0+k}}{P^{3h_0-3+a_0+k}}.
\end{equation}
With $G$ given by \eqref{eq:sublinear-tail-series}, equality becomes
an upper bound. The cusp labels have been included exactly once.

\emph{Step 3: factorial estimates for the genus sum.}
We record the factorial estimate used for the genus sum. For $p\ge1$,
\[
 \sum_{z_1+\cdots+z_p=J}\prod_i z_i!\le3^{p-1}J!,
 \qquad z_i\ge0.
\]
For two factors this is the bound $\sum_i\binom Ji^{-1}\le3$;
iteration proves the assertion. If $z_*=\max_i z_i$, $\sum z_i=K$,
and $j=K-z_*$, then, for $1\le q\le K/4$ and $K/w\ge q$,
\begin{equation}\label{eq:sublinear-factorial-tail}
 \sum_{j\ge q}\prod_{i=1}^w\Gamma(h_i)
       \le C^w w q!(K-q)!,\qquad z_i=h_i-1.
\end{equation}
Indeed, for a fixed root and $j$ the outside convolution costs at
most $3^{w-2}j!$. The root costs $(K-j)!$, and its maximality gives
$q\le j\le K-K/w\le K-q$. On this interval,
$\sum_j\binom Kj^{-1}\le C\binom Kq^{-1}$: split at $K/3$ and
$2K/3$, sum geometrically on the endpoint intervals, and use the
exponential decrease from $q\le K/4$ to the middle interval.
The constant is absolute. For $q=0$ use the full convolution bound;
for $w=1$ the sector $q\ge1$ is empty. Root choices cost at most $w$.
Dropping parity restrictions only increases these positive sums.
Our genus condition ensures $K/w\ge t$, $K\ge4t$, and
\begin{equation}\label{eq:sublinear-factorial-ratio}
 \frac{(K-j)!}{(2g-2)!}\le g^{-\Delta-j}\quad(0\le j\le t).
\end{equation}

\emph{Step 4: the envelope bounds.}
To apply these estimates, put $\rho(x)=r(Y(x))$ and
$\mathcal S(x)=a(Y(x))\rho(x)^{3/2}A(Y(x))$. Undoing the
normalization in \eqref{eq:gevrey-envelope} gives
\begin{equation}\label{eq:sublinear-stable-envelope}
 \Phi_l^{(p)}(x)\le C4^{-l}\frac{\Gamma(h)}{\sqrt h}
             \mathcal S(x)^{-h}\rho(x)^{p/2},
 \qquad h=2l-2+p\ge1.
\end{equation}
The parent $p=1$ has a matching lower bound. All scalar functions
are bounded above and away from zero on $[0,y]$. At a common $x$,
the remaining powers of $4,\rho,\mathcal S$ cost $C^{O(m+b)}$,
because $\sum h_i=H$. Also
$\sqrt{2g-1}/\prod\sqrt{h_i}\le C\sqrt w$.
The unstable factors cost $x^{2u+v}$, and
$v_0\le C^{m+1}(m-1)!$ by the same envelope at zero.

\emph{Step 5: the two sectors of the deficit.}
If $d\ge r$, the full convolution, \eqref{eq:sublinear-coefficient-ratio}
and \eqref{eq:sublinear-label-dictionary} give
\begin{equation}\label{eq:sublinear-fixed-pattern}
 \sum W\le[Ct]^{Ct}n x^{B-1}g^{-\Delta}
       \le[Ct]^{Ct}n^B g^{-d}
       \le[Ct]^{Ct}g^{-r}n^{r+1}.
\end{equation}
The last inequality uses $B\le d+1$ and $n/g\le1$.
If $d<r$, put $q=r-d$. For $j\ge q$ use
\eqref{eq:sublinear-factorial-tail}; the resulting bound is
$[Ct]^{Ct}n^B g^{-d-q}$ and suffices. For $0\le j<q$, positivity
gives the convergent cusp-tail estimate
\[
 [B_{\rm out}]_{\deg\ge q-j}(x)
             \le(x/y)^{q-j}B_{\rm out}(y).
\]
Evaluating the outer stable factors at $y$ instead of $x$ costs at
most $C^{j+w-1}$ in \eqref{eq:sublinear-stable-envelope}, since their
total complexity is $j+w-1$. Thus this sector is bounded by
\[
 [Ct]^{Ct}n x^{B-1+q-j}g^{-\Delta-j}
 \le[Ct]^{Ct}g^{-r}n^{B+r-d-j}
 \le[Ct]^{Ct}g^{-r}n^{r+1}.
\]
The at most $t$ sectors, $2^b$ strata and $(2k+2)^{2k+2}$ patterns
are absorbed into $[Ct]^{Ct}$. We have tilted only the outer cusp
series with bounded outer complexity, not a divergent genus series.
This proves \eqref{eq:sublinear-tail}. Boundary volumes and the
actual-edge Jacobian give $(1+\Lambda)^ke^\Lambda$ as in
\cite[Section~2.2]{HuangGap}.

\emph{The separating tail.}
For the ordered separating kernel there are two degree-one components,
with genera summing to $g$. The exact coefficient prefactor is $n/P$,
with no carrier volume. The same algebra has $d=1$: either
$(u,w,B,\Delta)=(1,1,2,0)$ or $(0,2,0,2)$, and $v=0$.
Its deficit $\min(c_1,c_2)$ is at most $1+j+A_{\rm ex}$.
Apply the preceding split with $q=\max(0,r-1)$ and $t=r+1$.
It includes both ordered placements and the central tie; boundary
volumes cost $e^L$. This proves \eqref{eq:sublinear-separating-tail}.
\end{proof}

For later geometric use (in the proof of
Theorem~\ref{thm:fixed-energy-spectrum} in
Appendix~\ref{app:fixed-window-proofs}), keep $m$ and the pattern
fixed. The first two bounds in \eqref{eq:sublinear-fixed-pattern},
applied to the whole sum, give
\begin{equation}\label{eq:sublinear-fixed-carrier}
 \sum W\le C_m n^B g^{-d}
 \le C_m g(n/g)^{m+b}
 \le C_m(n^2/g)(n/g)^{m-1}.
\end{equation}
Here $\Delta+B=m+b$ and $g^{-1}\le n/g$. Bounded boundary
integration changes only the constant. Every fixed carrier of
complexity at least two therefore has expected count $o(n^2/g)$
when $\sqrt g\ll n\ll g$.

\medskip\noindent\emph{The inputs from Huang's trace expansion.}
Before inserting the graded bounds into the trace expansion, we record
the statements of \cite{HuangGap} that the proof below combines with
Lemma~\ref{lem:sublinear-tail}. Appendix~\ref{app:critical-trace}
records Huang's inputs for the envelope $n\le M\sqrt g$; the lemmas
below restate, in the notation of this appendix, exactly those inputs
that are used at sublinear density, together with the genus threshold
under which each of them holds, since it is the thresholds that the
present argument modifies. Lemma~\ref{rev:lem:huang-retained} and
Lemma~\ref{rev:lem:huang-coefficients} are used in the proof of
Proposition~\ref{prop:sublinear-polynomial} below, whose last step
also invokes Lemma~\ref{rev:lem:huang-first-moment} (the first
inequality of which is in turn the base case of
\eqref{eq:sublinear-orbital} proved there);
Lemmas~\ref{rev:lem:huang-trace}--\ref{rev:lem:huang-retained-mass} and
Lemma~\ref{rev:lem:huang-first-moment} are used in the spectral
transfers of Appendix~\ref{app:fixed-window-proofs}. All numbered
references are to version~1 of that paper. We translate its notation
into ours: $s=n+1$;
the multiplier $h,f,a$ of Proposition~\ref{prop:critical-mean}, with
$\operatorname{supp}h\subset[-L,L]$ and $a=f(i/2)>1$; and the carrier,
pattern and deficit notation of Lemma~\ref{lem:sublinear-tail}, which
is that of Appendix~\ref{app:critical-trace} up to the change of
letters $\varsigma,\mathfrak d,\mathsf s\mapsto\sigma,\delta,s$. For a carrier type
$(h_0,a_0,k)$ of complexity $m$ and a pattern with $b$ exterior blocks,
a \emph{record} consists of a distinguished block $*$ together with the
genera and cusp counts $(g_i,a_i)_{i\ne*}$ of the other complementary
components, all stable. Its deficit is $\delta=m+\sum_{i\ne*}c_i$ with
$c_i=2g_i+a_i+k_i-2$, its outside cusp count is
$A_{\rm out}=a_0+\sum_{i\ne*}a_i\le\delta+1$, and the component $*$ is
reconstructed from $g$ and $n$ through the two sum constraints
\cite[(2.40)--(2.41)]{HuangGap}; this reconstruction needs only
$N=2g+n-3>2R$ and $g>h_0+k+\sum_{i\ne*}g_i$. At expansion order $R$
the records with $\delta\le R$ are \emph{retained} and expanded to the
local order $J_{\rm loc}=R-\delta$, while those with $\delta\ge R+1$
are \emph{omitted} and form the positive tail
$\mathcal T^{h_0,a_0,k}_{g,n}(R)=\sum_{\rm patterns}\sum_{\delta\ge R+1}W$
bounded in \eqref{eq:sublinear-tail}. For a pattern with slot lengths
$\mathbf x$, $\Lambda=\sum_jx_j$. The one-curve kernels of
\cite[Section~3.1]{HuangGap} are the nonseparating kernel
$U_{g,n}(x)=V_{g-1,n+2}(\mathbf 0_n,x,x)/V_{g,n}$, with leading term
$\phi(x)^2$ and canonical coefficients $u_{j,n}$, and the ordered
separating kernel $S_{g,n}(x)$, with canonical coefficients $s_{j,n}$;
an endpoint $(h_0,a_0)$ of the latter has deficit $\delta=2h_0+a_0-1$
and is retained when $\delta\le R$.

\begin{lemma}[Retained records at a linear genus threshold]
\label{rev:lem:huang-retained}
There is an absolute constant $C$ with the following properties, for
all integers $R,n\ge0$ and $r=R+1$.
\begin{enumerate}
\item[\textup{(i)}] \textup{(\cite[Lemma~4.2, (2.42), (4.13)--(4.14) and the proof of Proposition~4.7, in particular (4.26)]{HuangGap})}
Let $e$ be a retained record of deficit $\delta\le R$ for a carrier type
of complexity $m$ and a pattern, and put $t=R+m+1$. If
$g\ge[Ct]^Cs$, then the contribution of $e$ to the pattern kernel
differs from its canonical expansion through local order
$J_{\rm loc}=R-\delta$ by at most
\[
 g^{-R-1}s^{A_{\rm out}+J_{\rm loc}+1}[Ct]^{Ct}(1+\Lambda)^{2t}e^{\Lambda}
 \le g^{-R-1}s^{R+2}[Ct]^{Ct}(1+\Lambda)^{2t}e^{\Lambda}.
\]
The number of records of deficit $\delta$ is at most
$b(\delta+3)^{2(b-1)}\le[Ct]^{Ct}$.
\item[\textup{(ii)}] \textup{(\cite[(4.10) and proof of Proposition~4.7]{HuangGap})}
The nonseparating kernel has no omitted records, and for $g\ge[Cr]^Cs$
and $x\ge0$,
\[
 \Bigl|U_{g,n}(x)-\phi(x)^2-\sum_{j=1}^Rg^{-j}u_{j,n}(x)\Bigr|
 \le g^{-R-1}s^{R+1}[Cr]^{Cr}(1+x)^{2R+2}e^{x}.
\]
\item[\textup{(iii)}] \textup{(\cite[(4.11) and proof of Proposition~4.7]{HuangGap})}
A retained separating endpoint $(h_0,a_0)$ of deficit $\delta\le R$ has
$a_0\le\delta+1$, and for $g\ge[Cr]^Cs$ and $x\ge0$ its contribution to
$S_{g,n}(x)$ differs from its canonical expansion through order
$R-\delta$ by at most
\[
 g^{-R-1}s^{a_0+R-\delta+1}[Cr]^{Cr}(1+x)^{2R+2}e^{x}
 \le g^{-R-1}s^{R+2}[Cr]^{Cr}(1+x)^{2R+2}e^{x}.
\]
\end{enumerate}
In all three parts the genus threshold is linear in $s$. The quadratic
term $Cs^2$ in the thresholds of \cite[Proposition~4.7 and
Theorem~4.8]{HuangGap} enters only through the omitted records, that
is, through \cite[Lemma~4.6]{HuangGap}.
\end{lemma}

The first assertion is the product of the graded shifted-volume
expansion \cite[(4.3)]{HuangGap}, whose threshold is
$g\ge[CL]^Cs$ with $L=R+H_0+1\le2t$, and the small-volume prefactor
bound \cite[(4.13)]{HuangGap}, whose cusp power is
$s^{A_{\rm out}}$. The last sentence is the content of the proof of
\cite[Proposition~4.7]{HuangGap}, which reuses the retained-record
estimates of \cite[Proposition~4.4]{HuangGap} verbatim and changes only
the tail input.

\begin{lemma}[The compressed trace identity]
\label{rev:lem:huang-trace}
\textup{(\cite[Lemma~5.2, Theorem~5.3, (5.23), (5.30), (5.44)]{HuangGap}, restated in \cite[(6.1)--(6.2)]{HuangGap}.)}
Let $X$ be a complete finite-area hyperbolic surface of genus $g$
with $n\ge1$ cusps and no orbifold points, and let $A_X$ be the bounded
self-adjoint operator with Laplace multiplier $f(r(\lambda))$, where
$r(\lambda)=\sqrt{\lambda-1/4}$ for $\lambda\ge1/4$ and
$r(\lambda)=i\sqrt{1/4-\lambda}$ for $0\le\lambda<1/4$. The constant
eigenfunction has $A_X$-eigenvalue $a$, and the spectrum in
$[1/4,\infty)$ is mapped into $[0,1]$. For an integer $Q\ge1$ put
\[
 Y_Q=3e^{LQ},
\]
and let $\chi_Q$ be the indicator of $X(Y_Q)$, the complement of the
$n$ width-one cusp regions $\{y>Y_Q\}$; these regions are embedded,
pairwise disjoint and of area $1/Y_Q$ each. For a real polynomial $p$
of degree at most $Q$ with $p(0)=0$, write $p=\sum_tc_tx^t$ and
$H_p=\sum_tc_th^{*t}$, so that $\operatorname{supp}H_p\subset[-LQ,LQ]$
and $\widehat{H_p}=p\circ f$. Then $\chi_Qp(A_X)\chi_Q$ is trace class and
\begin{equation}\label{rev:eq:huang-trace}
 \tau_{X,Q}(p):=\operatorname{Tr}\bigl(\chi_Qp(A_X)\chi_Q\bigr)
 =\Bigl(\operatorname{area}(X)-\frac n{Y_Q}\Bigr)\mathcal I(p)
   +n\mathcal P_{Y_Q}(p)+\mathcal O_X(H_p),
\end{equation}
where
\[
 \mathcal I(p)=\frac1{4\pi}\int_{\mathbb R}r\tanh(\pi r)\,p(f(r))\,dr,
 \qquad
 \mathcal P_Y(p)=2\int_0^{LQ}K_{H_p}(\rho)\cosh(\rho/2)
      \,\mathsf H_{\lfloor2Y\sinh(\rho/2)\rfloor}\,d\rho,
\]
$K_H$ is the radial kernel with Abel transform $H$
\cite[(5.14)--(5.15)]{HuangGap} and $\mathsf H_m=\sum_{k=1}^mk^{-1}$,
$\mathsf H_0=0$. The identity holds for every such surface, with no
genus or cusp condition. Spectrally,
\[
 \tau_{X,Q}(p)=\sum_jp(f(r_j))\|\chi_Qu_j\|_2^2
 +\frac1{4\pi}\sum_{b=1}^n\int_{\mathbb R}p(f(r))
     \|\chi_QE_b(\cdot,\tfrac12+ir)\|_2^2\,dr,
\]
with absolute convergence, over an orthonormal eigenbasis $u_j$ with
spectral parameters $r_j$ and the Eisenstein series $E_b$. In
particular $\tau_{X,Q}(p)\ge0$ when $p\ge0$ on the spectrum of $A_X$,
and for $p=x\widetilde p$,
\[
 |\tau_{X,Q}(p)|\le\|\widetilde p\|_{[0,a]}\,\tau_{X,Q}(x).
\]
The constant eigenfunction has retained mass
$1-n/(\operatorname{area}(X)Y_Q)$, so removing its spectral projection
subtracts exactly $(1-n/(\operatorname{area}(X)Y_Q))\,p(a)$ from
$\tau_{X,Q}(p)$, and the remaining compression is still positive for
$p\ge0$ on $[0,a]$.
\end{lemma}

\begin{lemma}[Deterministic bounds on the identity and parabolic terms]
\label{rev:lem:huang-deterministic}
\textup{(\cite[Proposition~5.4]{HuangGap}, restated in \cite[(6.3)]{HuangGap}.)}
There are finite constants $I_f,J_f,b_f\ge0$, depending only on the
multiplier $(h,f)$ and $L$, such that for every integer $Q\ge1$ and
every real polynomial $p=x\widetilde p$ of degree at most $Q$,
\[
 |\mathcal I(p)|\le I_f\|\widetilde p\|_{[0,1]},\qquad
 |\mathcal P_{Y_Q}(p)|\le J_fQ^2\|\widetilde p\|_{[0,1]},\qquad
 |\mathcal P_{Y_Q}(x)|\le b_fQ.
\]
These bounds involve only the values of $\widetilde p$ on the real
spectral image $[0,1]$, assign no sign to the parabolic term, and are
independent of the surface.
\end{lemma}

\begin{lemma}[Retained mass, uniformly in the cusp number]
\label{rev:lem:huang-retained-mass}
\textup{(\cite[Lemma~6.5]{HuangGap}.)}
Fix $0<\varepsilon<1/4$. If $\psi$ is an $L^2$-normalized
eigenfunction of $\Delta_X$ with eigenvalue
$0<\lambda\le1/4-\varepsilon$, then for every integer $Q\ge1$
\[
 \|\chi_Q\psi\|_2^2\ge1-e^{-2\sqrt\varepsilon\,LQ}.
\]
The estimate includes residual eigenfunctions and is uniform in the
surface $X$, in $g$ and in the number $n$ of cusps: the proof bounds the
mass of $\psi$ above height $3e^{LQ}$ in each cusp by
$e^{-2\sqrt\varepsilon LQ}$ times its mass above height $3$, and the
height-three regions are pairwise disjoint, so their total mass is at
most one.
\end{lemma}

\begin{lemma}[Coefficient bound]
\label{rev:lem:huang-coefficients}
\textup{(\cite[Lemma~6.7]{HuangGap}.)}
There are an integer $D_0\ge2$ and a constant $C_0\ge2$, depending only
on the multiplier, such that for all integers $j\ge2$, $Q,n\ge1$ and
every real polynomial $p=x\widetilde p$ of degree at most $Q$,
\[
 |\omega_{j,n}(p)|\le s^j[C_0(Q+j)]^{D_0j}\|\widetilde p\|_{[0,a]}.
\]
No condition relates $g$, $j$ and $Q$: the bound is obtained by
interpolating the coefficient polynomial in the cusp variable from an
auxiliary genus limit at the cusp counts $1,\dots,j+1$, and no
finite-genus remainder at the actual cusp count is extrapolated.
\end{lemma}

\begin{lemma}[First moment of the compressed trace]
\label{rev:lem:huang-first-moment}
\textup{(\cite[Lemma~6.8, (6.43)--(6.45)]{HuangGap}, with the modified threshold.)}
There are constants $C$, depending only on the multiplier, and $n_0$
such that, if $n\ge n_0$ and $g\ge C s$, then
\[
 \bigl|\mathbb E\mathcal O_X(x)-\omega_1(x)\bigr|\le C\frac{s^2}{g}\le g,
 \qquad
 \mathbb E\tau_{X,Q}(x)\le D_f(g+nQ)\quad\text{for every }Q\ge1,
\]
where $D_f=1+6a+6\pi I_f+b_f$. The first inequality is the base case
$R=0$, $T=L$ of \eqref{eq:sublinear-orbital} below; Huang's version
\cite[(6.45)]{HuangGap} has the hypothesis $g\ge c_3s^2+[2c_3]^{c_3}s$,
and the quadratic term is removed by the argument of
Proposition~\ref{prop:sublinear-polynomial}. The second inequality
follows from the first by inserting Lemma~\ref{rev:lem:huang-deterministic}
in \eqref{rev:eq:huang-trace} with $p=x$, $H_x=h\ge0$, together with
$|\omega_1(x)|\le6a$ \cite[Proposition~5.1]{HuangGap} and
$0\le\operatorname{area}(X)-n/Y_Q\le4\pi g+2\pi n$. Trace domination
then makes every degree-at-most-$Q$ polynomial trace absolutely
integrable.
\end{lemma}

We now replace the positive tails of the trace expansion by these
graded bounds.

\begin{proposition}[Polynomial trace estimate at sublinear density]
\label{prop:sublinear-polynomial}
Put $s=n+1$, $\tau=g/s$, $D=s^2/g$, and use the multiplier and
canonical coefficients of Proposition~\ref{prop:critical-mean}.
For a sufficiently small fixed $\theta_0>0$, $Q=\lfloor\tau^{\theta_0}\rfloor$,
and $p=x\widetilde p$ of degree at most $Q$, as $g,n,\tau\to\infty$,
\begin{equation}\label{eq:sublinear-polynomial}
 \left|\mathbb E\mathcal O_X(p)-\omega_1(p)-g^{-1}\omega_{2,n}(p)\right|
 \le D\left\{2\tau^{-1}(3C_0Q)^{3D_0}+2^{-Q}\right\}
                              \|\widetilde p\|_{[0,a]}.
\end{equation}
The same estimate holds for $Q=\lceil K\log\tau\rceil$ with fixed
$K>0$. Its constants are independent of $n/\sqrt g$.
\end{proposition}

\begin{proof}
We run the argument of \cite[Sections~4.2 and~6.7]{HuangGap} with the
omitted records bounded by Lemma~\ref{lem:sublinear-tail} instead of
\cite[Lemma~4.6]{HuangGap}. The following list records which genus
threshold of that argument is replaced by which estimate; every input
not listed as replaced is used with its original hypothesis, which is
linear in $s$.
\begin{center}\small
\begin{tabular}{@{}>{\raggedright\arraybackslash}p{0.40\textwidth}>{\raggedright\arraybackslash}p{0.22\textwidth}>{\raggedright\arraybackslash}p{0.33\textwidth}@{}}
\hline
step in \cite{HuangGap} & threshold there & here\\
\hline
omitted carrier records, \cite[(4.17)]{HuangGap}
 & $g\ge K(s^2+t^2s)$
 & \eqref{eq:sublinear-tail}: $g\ge C(s+t^2)$\\
omitted separating endpoints, \cite[(4.18)]{HuangGap}
 & $g\ge K(s^2+r^2s)$
 & \eqref{eq:sublinear-separating-tail}: $g\ge C(s+r^2)$\\
retained carrier records, \cite[(4.3), (4.13)]{HuangGap}
 & $g\ge[Ct]^Cs$
 & unchanged, Lemma~\ref{rev:lem:huang-retained}(i)\\
nonseparating kernel, \cite[(4.10)]{HuangGap}
 & $g\ge[Cr]^Cs$
 & unchanged, Lemma~\ref{rev:lem:huang-retained}(ii)\\
retained separating endpoints, \cite[(4.11)]{HuangGap}
 & $g\ge[Cr]^Cs$
 & unchanged, Lemma~\ref{rev:lem:huang-retained}(iii)\\
orbital integration, \cite[Theorem~4.8]{HuangGap}
 & $g\ge Cs^2+[C\varrho]^Cs$
 & $g\ge[C\varrho]^Cs$, \eqref{eq:sublinear-orbital}\\
monomials and polynomials, \cite[Lemma~6.8, Proposition~6.9]{HuangGap}
 & $g\ge c_3s^2+[c_3q]^{c_3}s$
 & $g\ge[Cq]^Cs$, i.e.\ $[Cq]^C\le\tau$\\
coefficient bound, \cite[Lemma~6.7]{HuangGap}
 & none
 & unchanged, Lemma~\ref{rev:lem:huang-coefficients}\\
\hline
\end{tabular}
\end{center}
Here $t=R+m+1$, $r=R+1$, $\varrho=R+\lceil T\rceil+1$ and $q=Q+2$.
After enlarging $C$, the linear threshold $g\ge[Ct]^Cs$ implies
$g\ge C(s+t^2)$, so the replaced tail thresholds are dominated by the
unchanged retained-record thresholds. The threshold also gives
$N>2R$ and $g>h_0+k+\sum_{i\ne*}g_i$, which is all that the
reconstruction of retained records requires.

\emph{Kernel estimates.} Fix a carrier type of complexity $m\ge1$, a
pattern and $R\ge0$, and suppose $g\ge[Ct]^Cs$ and $n\ge n_0$. By
Lemma~\ref{rev:lem:huang-retained}(i), each of the at most
$[Ct]^{Ct}$ retained records contributes an error at most
$g^{-R-1}s^{R+2}[Ct]^{Ct}(1+\Lambda)^{2t}e^\Lambda$, the cusp power
being $s^{A_{\rm out}+J_{\rm loc}+1}\le s^{R+2}$. The omitted records
are bounded by the boundary version of \eqref{eq:sublinear-tail}, that
is, by $g^{-R-1}s^{R+2}[Ct]^{Ct}(1+\Lambda)^ke^\Lambda$. Hence the
pattern kernel differs from its canonical expansion through order $R$
by at most
\[
 g^{-R-1}s^{R+2}[Ct]^{Ct}(1+\Lambda)^{2t}e^{\Lambda},
\]
which is \cite[(4.12)]{HuangGap} under the linear threshold alone.
When $R=0$ or $m>R$ there are no retained records and the tail is the
whole error. For the simple kernels, suppose $g\ge[Cr]^Cs$. The
nonseparating kernel has the smaller cusp power,
\[
 \Bigl|U_{g,n}(x)-\phi(x)^2-\sum_{j=1}^Rg^{-j}u_{j,n}(x)\Bigr|
 \le g^{-R-1}s^{R+1}[Cr]^{Cr}(1+x)^{2R+2}e^{x},
\]
by Lemma~\ref{rev:lem:huang-retained}(ii), with nothing to replace.
For the separating kernel, the at most $O(r^2)$ retained endpoints,
each with two ordered placements, contribute by
Lemma~\ref{rev:lem:huang-retained}(iii) at most
\[
 g^{-R-1}s^{a_0+R-\delta+1}[Cr]^{Cr}(1+x)^{2R+2}e^{x}
 \le g^{-R-1}s^{R+2}[Cr]^{Cr}(1+x)^{2R+2}e^{x}
\]
each, and the omitted endpoints contribute
$S^{\ge r}_{g,n}(x)\le g^{-r}s^{r+1}(Cr)^{Cr}e^x$ by
\eqref{eq:sublinear-separating-tail}. Thus
\cite[(4.10)--(4.11)]{HuangGap} hold under $g\ge[Cr]^Cs$, $n\ge n_0$.
The canonical coefficients and their vanishing rules are unchanged, and
no finite-genus remainder is interpolated in $n$.

\emph{Orbital integration.} The rooted integration formula, the length
simplex and the fiber bounds of \cite[Sections~3.3--3.4]{HuangGap}
involve no genus condition. Repeating the proof of
\cite[Theorem~4.8]{HuangGap} with the kernel estimates just proved
gives, for even $H$ supported in $[-T,T]$ and
$\varrho=R+\lceil T\rceil+1$,
\begin{equation}\label{eq:sublinear-orbital}
 \left|\mathbb E\mathcal O_H-\Omega(H)
            -\sum_{j=1}^R g^{-j}\kappa_{j,n}(H)\right|
 \le g^{-R-1}s^{R+2}(C\varrho)^{C\varrho}e^{3T}\|H\|_\infty
\end{equation}
under $g\ge[C\varrho]^C s$, $n\ge n_0$. This threshold gives
$\pi(g-1+n/2)>T$ and so excludes whole-surface carriers. All remaining
types have $m\le\lfloor2T/\pi\rfloor$ and $t\le\varrho$, and the pattern
counts and length powers are absorbed in $(C\varrho)^{C\varrho}$.

\emph{Polynomial transfer.} Apply \eqref{eq:sublinear-orbital} to
$h^{*j}$, $1\le j<q$, with $R=q-2$ and $q=Q+2$. After division by
$g$, the monomial error is $g^{-q}s^q(C_fq)^{C_fq}$, which is
\cite[(6.42)]{HuangGap} under $g\ge[C_fq]^{C_f}s$. Polynomial
coefficient conversion costs $7^Q$ \cite[Lemma~6.1]{HuangGap};
relative to $D$, the terminal error is at most
\[
 7^Q[C_f(Q+2)]^{C_f(Q+2)}\tau^{-Q}\le2^{-Q}
\]
if $C_f\theta_0<1/2$, also choosing $\theta_0$ small enough for the
orbital threshold $[C_fq]^{C_f}\le\tau$.
Lemma~\ref{rev:lem:huang-coefficients} gives
\[
 \sum_{i=3}^{Q+1}g^{1-i}|\omega_{i,n}(p)|
 \le\frac{2s^3(3C_0Q)^{3D_0}}{g^2}\|\widetilde p\|_{[0,a]},
\]
provided $(3C_0Q)^{D_0}/\tau=o(1)$. Take $3D_0\theta_0<1$.
The same bounds hold for logarithmic $Q$, proving
\eqref{eq:sublinear-polynomial}. The base case $R=0$, $T=L$ of
\eqref{eq:sublinear-orbital} gives
$|\mathbb E\mathcal O_X(x)-\omega_1(x)|=O(s^2/g)=O(D)$, hence
$\mathbb E\mathcal O_X(x)=O(1+D)$. Since $D=O(g)$,
Lemma~\ref{rev:lem:huang-first-moment} applies at the linear threshold
and the positive common-truncation trace has first moment
$\mathbb E\tau_{X,Q}(x)\le D_f(g+nQ)$ for every $Q\ge1$, as required
for the spectral approximation.
\end{proof}

\section{Orbital statistics at positive cusp density}
\label{app:positive-density-orbits}

This appendix proves the orbital concentration
\eqref{eq:positive-orbital-L1}, the geometric input of
Theorem~\ref{thm:positive-density-spectrum}
(Section~\ref{sec:positive-density-spectrum}); the transfer to the
spectrum is carried out in Appendix~\ref{app:fixed-window-proofs}.
The strategy is to define the limiting functional $\Omega_\lambda$
through planar carrier integrals, to compute first and second moments
of the orbital statistic $\mathcal O_X(H)$ by summing cusp labels over
gluing graphs with one or two distinguished carriers, and to deduce
concentration from the factorization of the leading graphs. The inputs
are Huang's normalized carriers and rooted integration
formula~\cite{HuangGap}, the primitive count of
Lemma~\ref{lem:primitive-count}, the boundary profile of
Proposition~\ref{prop:bounded-boundaries}, and the stable envelope
\eqref{eq:sublinear-stable-envelope} of Appendix~\ref{app:sublinear-trace}.
The first subsection defines $\Omega_\lambda$ and checks that it is
finite, the second bounds every gluing graph and identifies the leading
ones, and the third proves concentration. The graph bound
\eqref{eq:positive-graph-power}, the root derivative
\eqref{eq:positive-root-derivative} and the coefficient bound
\eqref{eq:positive-coefficient-bound} are reused in
Appendices~\ref{app:cusp-dominated-orbits}, \ref{app:rooted-local-law}
and~\ref{app:conformal-phase}. All lengths and carrier complexities
below are fixed before $g\to\infty$. Constants may depend on them and
on a compact interval of limiting cusp densities.

\subsection{The planar functional}

In the limiting functional, all complementary planar pieces are
absorbed into a single convergent cap series. We define this series,
the simple-waist weight $J_H$ and the fiber integrals of Huang's
normalized carriers, assemble them into the functional
\eqref{eq:positive-orbital-functional}, and then record the three
facts that make every term finite.

Put $P=2\pi^2$, $x=X(Y_\lambda)$, so $2\Psi(x)=\lambda$, and define
\begin{equation}\label{eq:positive-cap-series}
 K_x(L)=\sum_{a\ge2}
 \frac{V_{0,a+1}(0^a,L)}{P^{a-2}}\frac{x^a}{a!}.
\end{equation}
It satisfies $K_x(0)=\Phi_0'(x)$ and
$K_x(L)\le e^{L/2}\Phi_0'(x)$. The series is analytic for $x<x_c$,
uniformly on bounded length intervals and compact subintervals of
$(0,x_c)$. Thus planar complementary pieces with arbitrarily many
cusps are included in one convergent series.

For an even $H\in C_c^\infty(\mathbb R)$, set
\[
 J_H(L)=\sum_{r\ge1}\frac{L H(rL)}{\sinh(rL/2)}.
\]
This is the contribution of both orientations of one simple waist.
We use the normalized carriers of Huang~\cite[Section~3.3]{HuangGap}.
The \emph{normalized carrier} $Y$ of a nonsimple closed geodesic
$\gamma$ is obtained from a small regular neighborhood of its graph by
filling the complementary disks, completing the peripheral punctured
disks to their cusps, and tightening the remaining boundary curves to
geodesics; two boundaries tightening to the same ambient geodesic are
kept as two distinct slots. Its type $(h_0,a_0,k)$ records genus,
inherited cusps and positive boundary slots, its complexity is
$m=2h_0+a_0+k-2\ge1$, the inclusion $Y^\circ\to X$ is an injective
local isometry, and $\gamma$ \emph{fills} $Y$: the complement of a
small neighborhood of its graph in $Y$ consists of disks, cusp
neighborhoods and boundary collars.
For a planar carrier with $a$ cusps and $k$ labeled positive boundary
slots of lengths $\mathbf L$, let $I_H^{0,a,k}(\mathbf L)$ be the
integral, over the moduli space $\mathcal M_{0,a+k}(\mathbf 0_a,\mathbf L)$
with respect to the orbifold Weil--Petersson measure, of the fiber
statistic
\[
 A_H(Y)=\sum_{\alpha}\sum_{\nu\ge1}
   \frac{\ell_\alpha H(\nu\ell_\alpha)}{2\sinh(\nu\ell_\alpha/2)},
\]
the sum running over the oriented primitive geodesics $\alpha$ filling
$Y$; repetitions are included. This is the fiber integral
\cite[(3.23)]{HuangGap}, and the same orbifold measure defines the
volume $V_{0,a+k}(\mathbf 0_a,\mathbf L)$, so the quotient
$I_H^{0,a,k}/V_{0,a+k}$ is the fiber average used in Huang's rooted
integration formula \cite[Proposition~3.4]{HuangGap}. Write
$\widehat I_H^{a,k}=P^{-(a+k-3)}I_H^{0,a,k}$. Define
\begin{align}\label{eq:positive-orbital-functional}
 \Omega_\lambda(H)&=\Omega_\lambda^{\rm s}(H)
                         +\Omega_\lambda^{\rm ns}(H),\notag\\
 \Omega_\lambda^{\rm s}(H)
 &=\frac{\lambda}{Px}\int_0^\infty
       L\mathfrak b_{Y_\lambda}(L)K_x(L)J_H(L)\,dL,\notag\\
 \Omega_\lambda^{\rm ns}(H)
 &=\sum_{\substack{a\ge0,\ k\ge1\\a+k\ge3}}
 \frac{\lambda x^{a-1}}{a!k!P^k}
 \int_{(0,\infty)^k}\widehat I_H^{a,k}(\mathbf L)
 \left[\sum_{j=1}^k\mathfrak b_{Y_\lambda}(L_j)
                       \prod_{i\ne j}K_x(L_i)\right]
                       \prod_i L_i\,dL_i.
\end{align}
The definition of $\Omega_\lambda^{\rm ns}$ uses three inputs, which
we state explicitly.
\begin{enumerate}
\item[\textup{(a)}] \emph{Area and perimeter of a filled carrier}
\cite[Lemma~3.3, (3.22)]{HuangGap}: if $\gamma$ fills a normalized
carrier $Y$ of complexity $m$, then
\[
 \ell(\partial Y)\le2\ell_\gamma,\qquad
 2\pi m=\operatorname{area}(Y)\le4\ell_\gamma,\qquad
 \ell_\gamma\ge\pi m/2.
\]
Hence if $\operatorname{supp}H\subset[-T,T]$, only carriers with
$a+k-2=m\le2T/\pi$, area at most $4T$ and $\sum_iL_i\le2T$ occur in
the nonsimple term, and the sum over $(a,k)$ in
\eqref{eq:positive-orbital-functional} is finite.
\item[\textup{(b)}] \emph{The orbifold measure of the fiber integral}
\cite[(3.23), the paragraph after (3.25), and Proposition~3.4]{HuangGap}: numerator and denominator
of the fiber average are integrals against the same orbifold
Weil--Petersson measure on $\mathcal M_{0,a+k}(\mathbf 0_a,\mathbf L)$,
and the maps $\mathbf L\mapsto I_H^{0,a,k}(\mathbf L)$ are Borel on
$(0,\infty)^k$, with values defined on every fiber rather than almost
everywhere.
\item[\textup{(c)}] \emph{Finitely many primitives per carrier}
(Lemma~\ref{lem:primitive-count}): a carrier of complexity
$m\le2T/\pi$ carries at most $C(T)$ primitive unoriented geodesics of
length at most $T$, and a contributing repetition has
$\nu\le2T/(\pi m)$ by (a). Since
$\ell/(2\sinh(\nu\ell/2))\le1/\nu$, this gives
$|A_H|\le C_T\|H\|_\infty$ on every fiber, so
$|I_H^{0,a,k}(\mathbf L)|\le C_T\|H\|_\infty V_{0,a+k}(\mathbf 0_a,\mathbf L)$.
\end{enumerate}
These are actual Borel fiber
functions, so all integrals are well defined. At a short simple waist,
$|J_H(L)|\le C_{T,H}(1+|\log L|)$, which is integrable against
$L\,dL$. These bounds prove that \eqref{eq:positive-orbital-functional}
is finite and continuous in $\lambda\in(0,\infty)$.

\subsection{One and two distinguished carriers}

We bound the total weight of every gluing graph with one or two
distinguished carriers. To leading order, only trees with a single
genus-$g$ exterior vertex contribute, and they produce the planar
functional. The exact cusp-label sum is a coefficient ratio of the
marked free energies; it is estimated with the positive convolution
factor and the local limit theorem of Part~\ref{part:volumes}, through
the stable envelope \eqref{eq:sublinear-stable-envelope}.

Take $c\in\{1,2\}$ distinguished carriers of fixed types, with cusp
counts $a_{0,j}$, total genus $h_0$, total cusp count $a_0$, and $k$
boundary slots in total.
Pair $\sigma$ pairs of slots, allowing seams between distinguished
carriers as well as self-pairs. Partition the remaining slots into
$b$ nonempty complementary blocks, of degrees $k_i$. Only connected
gluing graphs are retained. They have $k-\sigma$ edges and $b+c$
vertices, and their complementary data satisfy
\[
 \sum_i g_i=g-h_0-k+\sigma+b+c-1,\qquad
 \sum_i k_i=k-2\sigma,\qquad \sum_i a_i=n-a_0.
\]
Classify the exterior vertices as $u$ planar vertices of degree one,
$v$ planar vertices of degree two, and $w$ vertices stable before
adding cusps. For the latter put $h_i=2g_i-2+k_i\ge1$. Then
\begin{equation}\label{eq:positive-stable-complexity}
 H:=\sum_{\rm stable}h_i=2g-2h_0-k+2c-2+u.
\end{equation}
For sufficiently large $g$, necessarily $w\ge1$.

The exact cusp-label sum is a coefficient ratio. Apart from the fixed
normalized distinguished volumes or fiber factors, its prefactor is
\[
 \frac{nP^{\sigma-k}}{\prod_{j=1}^c a_{0,j}!},
 \qquad
 \frac{[z^{n-a_0}]G(z)}{[z^{n-1}]\Phi_g'(z)},
\]
where $G$ sums the products of exterior marked free energies over
admissible genera. The power of $P$ is minus the number of actual
edges. Inherited cusp labels enter exactly once, as in Huang's rooted
integration formula~\cite[Section~3.4]{HuangGap}. For $c=2$ we use the
same formula with two distinguished carriers. By
Lemma~\ref{rev:lem:disjoint-carriers} below, cutting along the slot
geodesics of an ordered pair of disjoint carriers and regluing are
mutually inverse. Mirzakhani's integration formula for the resulting
multicurve gives the length--twist measure $\prod_e y_e\,dy_e$ on the
$k-\sigma$ actual edges, including seams between the two carriers,
times the volumes of the complementary pieces. The only finite
symmetry to be removed is the labelling of the slots of each carrier,
which gives the factor $1/(k_1!k_2!)$; since the pair is ordered,
there is no exchange of the two carriers.

Choose the parent saddle $x_g$ by $(2g-1)\Psi(x_g)=n-1$.
The positive convolution factor $V^H$ in the numerator, its uniform
point-mass bound, and the denominator local limit theorem give
\begin{equation}\label{eq:positive-coefficient-bound}
 \frac{[z^{n-a_0}]G(z)}{[z^{n-1}]\Phi_g'(z)}
       \le Cx_g^{a_0-1}\frac{G(x_g)}{\Phi_g'(x_g)}.
\end{equation}
The actual stable envelope \eqref{eq:sublinear-stable-envelope}
holds on every compact subinterval of $(0,x_c)$. Its Gamma convolution
is bounded by $C_w\Gamma(H-w+1)$, whereas the parent denominator
contains $\Gamma(2g-1)$. All remaining powers of $x_g$, $\rho$ and
$\mathcal S$ are fixed. Including $n=O(g)$, the total integrated
weight of the graph, summed over every exterior cusp and genus
allocation, is therefore
\begin{equation}\label{eq:positive-graph-power}
 O\bigl(g^{2c+1-2h_0-k+u-w}\bigr).
\end{equation}
Bounded boundary lengths and fixed twist integrals preserve this bound.

Connectedness implies $k-\sigma\ge b+c-1$. Hence the exponent in
\eqref{eq:positive-graph-power} is at most
$c+2-2h_0-v-2w-\sigma\le c$. Equality forces
\[
 h_0=0,\qquad v=0,\qquad w=1,\qquad \sigma=0,\qquad k=u+c,
\]
and the graph is a tree. It consists of planar distinguished carriers,
one exterior vertex of genus $g$, and planar degree-one caps. Each
distinguished carrier attaches once to the genus-$g$ vertex, whose
degree is $c$. Every other graph contributes $O(g^{c-1})$.

For $c=1$ the leading root series is $\Phi_g'$; for $c=2$ it is
$\Phi_g''$. The marking asymptotic and local limit give
\begin{equation}\label{eq:positive-root-derivative}
 \frac{\Phi_g''(x_g)}{\Phi_g'(x_g)}\sim\frac n{x_g}.
\end{equation}
The root's bounded boundary factors multiply, by
Proposition~\ref{prop:bounded-boundaries}. To insert the cap series
in coefficient extraction, first truncate their total cusp degree.
Every fixed coefficient shift has the same local-limit mass as the
parent. The point-mass bound in \eqref{eq:positive-coefficient-bound}
dominates the tail by its convergent positive series evaluated at
$x_g$. The truncation can therefore be removed uniformly on compact
density and length intervals. The mean shift of the twice-marked root
is $O(1)$ and does not change its $\sqrt g$ local-limit scale.

For one carrier this gives the coefficient $nx_g^{a_0-1}$ and,
after forgetting its $k!$ slot labellings, precisely
$\Omega_\lambda^{\rm ns}$. For two ordered carriers the coefficient
is $n^2x_g^{a_0-2}$, by \eqref{eq:positive-root-derivative}. The root
boundary factors and the factors $1/(k_1!k_2!)$ give the product of
the two one-carrier expressions. There is no extra division by two:
the two distinguished carriers are ordered.

\subsection{Concentration}

The second moments factorize to leading order, which gives $L^2$
concentration for nonsimple curves and for simple curves of length
bounded below; short simple waists are controlled separately in $L^1$.
The only new geometric point is
Lemma~\ref{rev:lem:disjoint-carriers}, which shows that two disjoint
nonsimple geodesics have carriers with disjoint interiors, so that
their pairs are counted by the two-carrier graphs of the preceding
subsection.

For nonsimple curves, the preceding one-carrier computation gives
\[
 g^{-1}\mathbb E\mathcal O_X^{\rm ns}(H)
                       \longrightarrow\Omega_\lambda^{\rm ns}(H).
\]
Consider its second moment. Two intersecting geodesic images fill a
connected carrier of bounded complexity and perimeter: Huang's
divergence-field proof of \cite[Lemma~3.3]{HuangGap} applies verbatim
to the connected geodesic graph $\gamma_1\cup\gamma_2$ and bounds both
the boundary length and the area of the filled carrier by
$2(\ell_{\gamma_1}+\ell_{\gamma_2})$. The rooted formula
\cite[Proposition~3.4]{HuangGap} holds for every invariant bounded
Borel fiber function, here the weighted number of ordered pairs
jointly filling the carrier.
Lemma~\ref{lem:primitive-count} bounds the number of ordered pairs
inside each such carrier. Their total expected contribution is
$O(g)$ by \eqref{eq:positive-graph-power} with $c=1$. The same bound
handles coincident images. Repetition indices of nonsimple curves
are bounded at fixed support because their primitive lengths have
a universal positive lower bound.

For disjoint images we use the following lemma. We write $Y_i$ both
for the abstract normalized carrier and for its image in $X$, a finite-type
subsurface whose boundary consists of the slot geodesics; a
self-paired slot (a \emph{seam}) lies in the interior of the image,
with the carrier on both sides. We use two elementary facts: a disk
and a cusp neighborhood contain no closed geodesic, and an annulus
contains no primitive closed geodesic other than its core, because each
free homotopy class of closed curves on a hyperbolic surface contains at
most one closed geodesic.

\begin{lemma}[Disjoint nonsimple geodesics have carriers with disjoint interiors]
\label{rev:lem:disjoint-carriers}
Let $\gamma_1,\gamma_2$ be primitive nonsimple closed geodesics on $X$
with disjoint images, and let $Y_1,Y_2$ be their normalized carriers.
Then $Y_1^\circ\cap Y_2^\circ=\varnothing$. Each slot geodesic of $Y_1$
either coincides with a slot geodesic of $Y_2$, the two carriers then
lying on opposite sides of this common seam, or is disjoint from
$Y_2$. Consequently the components of
$X\setminus(Y_1^\circ\cup Y_2^\circ)$ other than seams are stable
surfaces with geodesic boundary, and the pair $(Y_1,Y_2)$ together with
these components is a gluing graph with two distinguished carriers in
the sense of the preceding subsection, with each common seam a paired
slot between the carriers.
Conversely, gluing two filled carriers and stable complementary
components along such a graph recovers both carriers, by
\cite[Lemma~3.3]{HuangGap} applied to each carrier with the other
carrier counted among the stable exterior components.
\end{lemma}

\begin{proof}
The proof has three steps: $\gamma_2$ avoids $Y_1$, the interiors of
the carriers are disjoint, and their boundaries are either disjoint or
common seams. Each step uses only the structure of
$Y_i\setminus N_i^\circ$ recorded in the next paragraph.
Fix regular neighborhoods $N_i$ of the graphs of $\gamma_i$ so small
that $N_1\cap N_2=\varnothing$ and, whenever a compact set below is
shown to be disjoint from $\gamma_i$, it is also disjoint from $N_i$.
By the filling property, $Y_i\setminus N_i^\circ$ is a disjoint union
of closed disks, closed cusp neighborhoods, closed collars $C$ bounded
by one non-seam slot geodesic and one curve of $\partial N_i$, and, for
each seam, a closed annulus $A$ formed by the two collars of the seam,
with $\partial A\subset\partial N_i$ and the seam as core. Only the
collars $C$ meet $\partial Y_i$.

\emph{Step 1: $\gamma_2\cap Y_1=\varnothing$.} Since
$\gamma_2\cap N_1=\varnothing$, the set $\gamma_2\cap Y_1$ lies in
$Y_1\setminus N_1^\circ$. If $\gamma_2$ met a disk, a cusp
neighborhood or a seam annulus $A$ of $Y_1$, then, being connected and
disjoint from its boundary (which lies in $\partial N_1$), $\gamma_2$
would be contained in it, so $\gamma_2$ would be a closed geodesic in a
disk or a cusp neighborhood, or the simple core of $A$; all three are
impossible. If $\gamma_2\subset C$ for a collar $C$ of a slot $c^*$,
then $\gamma_2$ is a primitive closed geodesic in an annulus and equals
its core $c^*$, which is simple; impossible. Otherwise
$\gamma_2\cap Y_1$ is a union of arcs $\alpha$ whose endpoints lie on
slot geodesics and whose interiors lie in $Y_1\setminus N_1^\circ$ and
avoid $\partial Y_1$ (a tangency of two geodesics forces equality, and
a transversal interior crossing would leave $Y_1$). The interior of
such an arc lies in one component of $Y_1\setminus N_1^\circ$, which
must be a collar $C$ because only collars meet $\partial Y_1$; both
endpoints then lie on the unique slot geodesic $c^*$ of $C$, and the
interior is disjoint from $c^*$. Lift $C$, an embedded annulus with core
homotopic to $c^*$, to the cyclic cover of $X$ associated with $c^*$,
where it contains the unique closed geodesic, and then to $\mathbb H$:
the lift of $\alpha$ is a geodesic segment with both endpoints on the
axis of $c^*$, hence contained in that axis. Thus $\alpha\subset c^*$, so $\gamma_2$ and $c^*$ share an
arc and coincide; impossible. Hence $\gamma_2\cap Y_1=\varnothing$, and
by symmetry $\gamma_1\cap Y_2=\varnothing$. By compactness of the geodesic images,
$N_2\cap Y_1=\varnothing=N_1\cap Y_2$.

\emph{Step 2: $Y_1^\circ\cap Y_2^\circ=\varnothing$.} Suppose
$p\in Y_1^\circ\cap Y_2^\circ$. Then $p\notin N_1\cup N_2$, so $p$ lies
in a component $P$ of $Y_1\setminus N_1^\circ$. If $P$ is a disk, a
cusp neighborhood or a seam annulus, then $Y_2$, which is connected,
contains $p$ and is disjoint from $\partial P\subset\partial N_1$, is
contained in $P$; but $Y_2\supset\gamma_2$, and $P$ contains no
nonsimple closed geodesic. So $P$ is a collar $C$ of a non-seam slot
$c^*$ of $Y_1$, and $p\in C\setminus c^*$. Since $Y_2$ is connected,
disjoint from the other boundary curve of $C$, and not contained in
$C$, it meets $c^*$. As $c^*\cap N_2=\varnothing$, every point of
$c^*\cap Y_2$ lies in $Y_2\setminus N_2^\circ$.

If $c^*\not\subset Y_2$, then $c^*\cap Y_2$ is a union of arcs, each of
which lies, by the argument of Step~1 with the roles exchanged, in a
collar of a slot $d^*$ of $Y_2$ with both endpoints on $d^*$; the same
lift to $\mathbb H$ gives $c^*=d^*$, contradicting
$c^*\not\subset Y_2$. Hence $c^*\subset Y_2\setminus N_2^\circ$, and as
a closed geodesic it lies in a collar or seam annulus of $Y_2$ and
equals the corresponding slot geodesic: $c^*$ is a slot geodesic of
$Y_2$. If $c^*$ is a seam of $Y_2$, the half of its seam annulus
$A_2\setminus c^*$ on the side of $C$ is connected, disjoint from
$\partial C$, and meets $C^\circ$, so it lies in $C^\circ\subset Y_1$;
its outer boundary lies in $\partial N_2$, contradicting
$N_2\cap Y_1=\varnothing$. If $c^*$ is a non-seam slot of $Y_2$ with
collar $C_2$ on the same side of $c^*$ as $C$, then $C_2\setminus c^*$
is connected, disjoint from $\partial C$ and meets $C^\circ$, so again
$\partial N_2$ meets $Y_1$; contradiction. If $C_2$ lies on the
opposite side, then $Y_2\setminus c^*$ is connected (removing a
boundary circle from a connected surface), disjoint from
$\partial C$, and meets both $C^\circ$ (at $p$) and $X\setminus C$ (near
$c^*$ on the far side); contradiction. Hence
$Y_1^\circ\cap Y_2^\circ=\varnothing$.

\emph{Step 3: the boundaries.} If slot geodesics $c^*$ of $Y_1$ and
$d^*$ of $Y_2$ are distinct and meet, they cross transversally at some
point, near which $Y_1^\circ$ contains an open half-disk on one side of
$c^*$ and $Y_2^\circ$ an open half-disk on one side of $d^*$; two such
half-disks intersect, contradicting Step~2. Thus distinct slot geodesics
of the two carriers are disjoint; such a geodesic also avoids
$Y_2^\circ$, since a point of $c^*\cap Y_2^\circ$ would put a
half-disk of $Y_1^\circ$ inside $Y_2^\circ$, contradicting Step~2.
A common slot geodesic has the two carriers
on opposite sides by Step~2. A component of
$X\setminus(Y_1^\circ\cup Y_2^\circ)$ other than a seam is a surface
with geodesic boundary; a disk or a cusp neighborhood has no geodesic
boundary, and an annulus with two geodesic boundaries has equal
boundaries and is degenerate, so every such component is stable. The
converse statement is \cite[Lemma~3.3]{HuangGap}.
\end{proof}

By Lemma~\ref{rev:lem:disjoint-carriers}, two disjoint nonsimple
geodesic images have normalized carriers with disjoint interiors, whose
common seams are paired slots in the two-carrier formula. The leading
graph therefore factorizes as above, and every other graph is $O(g)$.
For any two tests of fixed support this proves
\[
 g^{-2}\mathbb E\bigl[\mathcal O_X^{\rm ns}(H_1)
                         \mathcal O_X^{\rm ns}(H_2)\bigr]
 \longrightarrow\Omega_\lambda^{\rm ns}(H_1)
                  \Omega_\lambda^{\rm ns}(H_2).
\]
In particular, the nonsimple statistic concentrates in $L^2$.

For simple curves first restrict lengths to $[d,T]$, $d>0$, where
the repetition weights are bounded. A separating edge has leading
graph consisting of a planar cap and one genus-$g$ root. Its exact
coefficient prefactor is $n/P$, and extraction contributes $K_x/x$.
It gives $\Omega_\lambda^{\rm s}$. Nonseparating curves contribute
$O(1)$ and separations into two positive-genus pieces contribute
$O(g^{-1})$, by the same envelope and coefficient argument.

Intersecting pairs again give $O(g)$. For two disjoint simple edges,
classify all cut components into $u,v,w$ as above. Their stable
complexities sum to $2g-2+u$ and their summed weight is
$O(g^{u-w+1})$. There are at most three vertices. The only
order-$g^2$ graph has $u=2,w=1$: two planar caps attached to one
genus-$g$ root of degree two. Equation~\eqref{eq:positive-root-derivative}
and boundary factorization give the product of the one-edge means.
All other graphs, and the diagonal, contribute $O(g)$. Thus the
$d$-truncated simple statistic also concentrates in $L^2$.

The one-edge coefficient bound is uniform down to zero waist length.
Together with the logarithmic repetition bound it gives
\[
 \frac1g\mathbb E\sum_{\substack{\gamma\ {\rm simple, unoriented}\\
                                  \ell_\gamma<d}}
       |J_H(\ell_\gamma)|
 \le C_{T,H}d^2(1+|\log d|),
\]
with $C_{T,H}$ uniform for $n/g$ in a compact subinterval of
$(0,\infty)$. The same estimate controls the removed limiting integral.
Let $d\downarrow0$ to obtain $L^1$ concentration of the full simple
part. Adding it to the nonsimple conclusion proves
\eqref{eq:positive-orbital-L1}; no mixed covariance calculation is
needed. The estimates also give the uniform fixed-support first
moment used in \eqref{eq:positive-trace-first-moment}.

\section{The orbital limit when cusps dominate genus}
\label{app:cusp-dominated-orbits}

This appendix proves the orbital concentration
\eqref{eq:cusp-dominated-orbital-L1} of
Theorem~\ref{thm:cusp-dominated-spectrum}
(Section~\ref{sec:cusp-dominated-spectrum}) and the uniform low-energy
tail of Proposition~\ref{prop:low-energy-tail}, which
Appendix~\ref{app:fixed-window-proofs} uses to exclude a zero atom in
Theorems~\ref{thm:positive-density-spectrum}
and~\ref{thm:cusp-dominated-spectrum} and at large times in
Corollary~\ref{cor:cusp-density-determinant}.
We extend Appendix~\ref{app:positive-density-orbits} to $g=o(n)$,
including bounded genus. The first- and second-moment structure is
unchanged; what changes is the source of the endpoint factors. For
diverging genus they are explicit powers of $\epsilon$, for fixed
positive genus they come from the fixed-genus coefficient asymptotics
of Manin--Zograf~\cite{ManZ} and Hide--Thomas~\cite{HTLargeN}, and in
genus zero from two planar convolution identities; the first three
subsections treat these cases in turn, and a subsequence argument
combines them. The final subsection proves the low-energy tail by
comparing $X$ with the closed surface of bounded geometry constructed
in Appendix~\ref{app:cusp-cutoffs}. Every support and carrier type
remains fixed.
Write $\epsilon=1-x_g/x_c$ at the parent saddle, as in
\eqref{eq:dense-parameters}. If $g\to\infty$
and $n/g\to\infty$, then
\begin{equation}\label{eq:cusp-endpoint-scales}
 \epsilon\sim\frac{5g}{2n},\qquad
 \rho(x_g)\asymp\epsilon^{1/2},\qquad
 \mathcal S(x_g)\asymp\epsilon^{5/4},
\end{equation}
by \eqref{eq:P10} and \eqref{eq:TL36}. The unstable series
$\Phi_0'$ and $\Phi_0''=Y$ have finite values at $x_c$.

\subsection{Endpoint bounds for diverging genus}

We keep the powers of $g$ and $\epsilon$ separate in every graph
weight; the leading graphs are the same planar trees as before, and no
endpoint factor offsets the loss of a power of $g$ elsewhere. The
input is the coefficient bound \eqref{eq:positive-coefficient-bound},
whose uniformity up to $x_c$ is checked first, together with the
endpoint scales \eqref{eq:cusp-endpoint-scales}.

Keep the notation for $c\in\{1,2\}$ distinguished carriers from
Appendix~\ref{app:positive-density-orbits}. Put
\[
 D=2g-1-H=2h_0+k-2c+1-u.
\]
The coefficient bound \eqref{eq:positive-coefficient-bound} is uniform
to $x_c$: the numerator and denominator point-mass scales are both
comparable to $\epsilon/\sqrt g$, since $H=2g+O(1)$.
Undoing the normalization in \eqref{eq:gevrey-envelope}, which is
uniform for $Y\in[0,y_0]$, gives the stable envelope
\eqref{eq:sublinear-stable-envelope} uniformly for $x\in(0,x_c)$
once the factors $\mathcal S(x)^{-h}\rho(x)^{p/2}$ are kept explicit;
only these factors degenerate at $x_c$.
In the stable envelope, the remaining endpoint factors are precisely
\[
 \mathcal S^D\rho^{(k-2\sigma-u-2v-1)/2}.
\]
The Gamma convolution contributes $O(g^{-D-w+1})$, and the exact
labeled prefactor contains $n$. Retaining these powers instead of
absorbing them into a density-dependent constant gives
\begin{equation}\label{eq:cusp-carrier-bound}
 \frac{\text{summed graph weight}}{n^c}
 \le Cg^{-r}\epsilon^B,\qquad
 r=2h_0+k-u+w-c-1,\quad
 B=\frac{5h_0+3(k-u-c)-\sigma-v}{2}.
\end{equation}
Constants depend only on the fixed graph and length box. For large
$g$ one has $w\ge1$. Connectedness gives
\[
 r\ge2h_0+v+2w+\sigma-2\ge0,\qquad
 B\ge\tfrac52h_0+v+\tfrac32(w-1)+\sigma\ge0.
\]
Equality $r=0$ forces $h_0=v=\sigma=0$, $w=1$, $k=u+c$:
the same planar tree with one genus-$g$ root and degree-one caps as
before. Every other graph loses at least one power of $g$, and no
negative power of $\epsilon$ offsets that loss.

For a leading two-carrier graph,
$\Phi_g''(x_g)/\Phi_g'(x_g)\sim n/x_g$ still holds uniformly.
The extra marking changes the model mean by $O(\Psi(x_g))$, whereas
the local-limit scale is comparable to $\sqrt g/\epsilon$; their
ratio tends to zero. Fixed cap-degree shifts therefore preserve the
parent local-limit mass. Remove the cap truncation using the uniform
point-mass bound and convergence of $K_{x_c}$. The boundary factors
tend to the product of $\mathfrak b_{y_0}$ by
Proposition~\ref{prop:bounded-boundaries}. The resulting one- and
two-carrier limits are $\Omega_\infty^{\rm ns}$ and its product.

For $q\in\{1,2\}$ simple cut edges, classify all complementary
vertices into $u,v,w$. Now $H=2g-2+u$, and the calculation gives
\begin{equation}\label{eq:cusp-simple-bound}
 \frac{\text{summed graph weight}}{n^q}
 \le Cg^{-r_s}\epsilon^{B_s},\qquad
 r_s=q-1-u+w,\quad B_s=\frac{3q-3u-v}{2}.
\end{equation}
Since $q\ge u+v+w-1$,
$r_s\ge v+2w-2\ge0$ and $B_s\ge v+3(w-1)/2\ge0$.
Only the tree with $q$ planar leaves and one genus-$g$ root has
$r_s=0$. Its root derivative gives the product simple intensity.

Intersecting pairs and coincident images have bounded union carriers,
so they contribute $O(n)$. For disjoint nonsimple pairs use
\eqref{eq:cusp-carrier-bound}; simple pairs use
\eqref{eq:cusp-simple-bound}. The short-waist cutoff error, divided
by $n$, is again $O(d^2(1+|\log d|))$. The concentration argument of
Appendix~\ref{app:positive-density-orbits} now proves
\eqref{eq:cusp-dominated-orbital-L1} when genus diverges. These bounds,
together with the compact-density estimates, also give uniform
fixed-support first moments divided by $n$ whenever $n/g$ is bounded
below by a positive constant: the graph bounds
\eqref{eq:cusp-carrier-bound} and \eqref{eq:cusp-simple-bound}, divided
by $n$, are uniform on $n/g\ge\lambda_0$ for diverging genus, and
\eqref{eq:positive-graph-power} covers each compact density interval.

\subsection{Fixed positive genus}

The rooted integration formulas \cite[Lemma~3.2 and
Proposition~3.4]{HuangGap} are stated there for $g\ge3$, but they are
exact Mirzakhani-type identities for every stable $(g,n)$ with
$n\ge3$. Each complementary one-holed torus without cusps contributes
its own factor $\tfrac12$, already included in the orbifold volume
$V_{1,1}$. With $M$ such components the product of their volumes
therefore includes the full factor $2^{-M}$, as in the componentwise
orbifold calculation in the proof of \cite[Proposition~3.4]{HuangGap}.

At fixed genus, the coefficients of the stable series grow
polynomially at critical weight. Each graph weight is then bounded by
a power of $n$, and the same planar trees give the leading order. The
local limit theorem of the diverging-genus case is replaced by a direct
convolution of power-law coefficients, whose exponent is matched to
the endpoint exponents $B$ and $B_s$ of the preceding subsection.

Fix $g\ge1$. For a pair $(l,p)$ stable before extra cusps, put
$\alpha(l,p)=(5l-5)/2+p>0$. The fixed-genus asymptotic of
Manin--Zograf~\cite[Theorem~6.1]{ManZ}, in the normalization of
\cite[Theorem~2.4]{HTLargeN}, implies
\begin{equation}\label{eq:cusp-fixed-genus-coefficients}
 [z^j]\Phi_l^{(p)}(z)x_c^j
       \sim C_{l,p}j^{\alpha(l,p)-1}.
\end{equation}
Indeed $(j+p)!/j!$ contributes $j^p$ to the volume asymptotic.
For bounded nonzero boundaries the coefficient has the additional
factor $\prod_i\mathfrak b_{y_0}(L_i)$, by
Hide--Thomas~\cite[Theorem~1.5]{HTLargeN}.

For a product of stable series, coefficient convolution gives exponent
$A-1$, where $A=\sum_i\alpha(l_i,p_i)>0$. To see the asymptotic
directly, scale the summation indices by $j$. The limiting simplex
integral has exponents $\alpha_i-1>-1$, and the endpoint portions
are bounded by the corresponding power sums. Multiplication by the
unstable series contributes their finite endpoint values. Their
critical coefficient tails are $O(j^{-3/2})$, or $O(j^{-5/2})$ when
all unstable factors have degree one. Split the unstable degree at
a fixed $M$ and at $\delta j$. The middle part is controlled by its
summable tail and the bounded ratio of stable coefficients. In the
last part, with unstable degree at least $\delta j$, the convolution
is bounded by
\[
 C_\delta j^{-3/2}\sum_{m\le j}(m+1)^{A-1}
       =O_\delta(j^{A-3/2})=o(j^{A-1}).
\]
This proves the asserted factorization and a matching coefficient
upper bound. The exponential boundary majorant handles length boxes.

For $c$ distinguished carriers, the genus and degree equations give
\[
 c-1+\alpha(g,1)-A=B,
\]
with $B$ from \eqref{eq:cusp-carrier-bound}. Thus the graph weight
is $O(n^{c-B})$. When $w\ge1$, $B=0$ forces the same planar tree
with one root of genus $g$; every other graph has $B>0$. If $w=0$,
the numerator consists entirely of unstable series and the weight is
$O(n^{1/2-\alpha(g,1)})=O(n^{-1/2})$, hence negligible.
For a leading graph the root has the parent genus, so its
Manin--Zograf constant cancels. Root differentiation and fixed shifts
give the coefficients $nx_c^{a_0-1}$ and $n^2x_c^{a_0-2}$, exactly
as for diverging genus. For $q$ simple edges the exponent is
$q-B_s$, with $B_s$ from \eqref{eq:cusp-simple-bound}; only the
root with $q$ planar leaves contributes to leading order. As for
carriers, a graph without a stable complementary vertex, which occurs
only for $g=1$, has an entirely unstable numerator and weight
$O(n^{-1/2})$.
Consequently the same first- and second-moment argument proves
\eqref{eq:cusp-dominated-orbital-L1} for every fixed $g\ge1$.

\subsection{Genus zero}

The parent series $\Phi_0'$ is itself summable at $x_c$, so the
stable-root argument does not apply. Instead, two convolution
identities for genus-zero series with one or two prescribed boundaries
identify the leading terms: the large complementary component is
selected by the slowest-decaying coefficient tail. The identities
\eqref{eq:cusp-sphere-convolutions} are also used in
Appendix~\ref{app:rooted-local-law} for the planar endpoint of the
local law. Let $K_p(z;\mathbf L)$ denote
the genus-zero generating series with $p$ prescribed boundaries,
with the normalization used for $\Phi_0^{(p)}$. For $p=1,2$,
its coefficients at critical weight have powers $-5/2$ and $-3/2$,
respectively, and the product of the factors $\mathfrak b_{y_0}$.
In particular $K_1=K_z$ in \eqref{eq:positive-cap-series} and
\begin{align}\label{eq:cusp-sphere-convolutions}
 [z^j]\prod_iK_1(z;L_i)
 &\sim [z^j]\Phi_0'(z)
       \sum_i\mathfrak b_{y_0}(L_i)\prod_{s\ne i}K_{x_c}(L_s),\notag\\
 [z^j]\left(K_2(z;L_1,L_2)\prod_iK_1(z;s_i)\right)
 &\sim [z^j]\Phi_0''(z)
       \mathfrak b_{y_0}(L_1)\mathfrak b_{y_0}(L_2)
                              \prod_iK_{x_c}(s_i).
\end{align}
For the first identity, split a two-factor convolution into its two
endpoint portions and the middle. With both indices at least $M$,
the middle is $O(j^{-5/2}M^{-3/2})$ after critical weighting.
Let $j\to\infty$, then $M\to\infty$, and induct on the number of
factors. For the second identity the $j^{-3/2}$ tail of $K_2$
dominates; a large degree in a $K_1$ factor costs $j^{-5/2}$ instead.
The same split proves the formula. All estimates are dominated on
bounded boundary boxes.

For one nonsimple carrier the genus-zero gluing graph is a tree with
only degree-one exterior vertices. The first identity in
\eqref{eq:cusp-sphere-convolutions} chooses the large complementary
piece and gives exactly the root sum in $\Omega_\infty^{\rm ns}$.
For two disjoint nonsimple carriers a tree either has one exterior
degree-two vertex joining them, or a direct seam between them.
In the first case the second identity applies, and
\[
 \frac{[z^{n-a_0}]\Phi_0''(z)}{[z^{n-1}]\Phi_0'(z)}
                \sim nx_c^{a_0-2}.
\]
This gives the product leading second moment. In the direct-seam case
all exterior series have degree one and their total weight is $O(n)$.

A simple curve on a sphere has two planar sides. The side-exchange
factor $1/2$ in its integration formula is cancelled by the two choices
of large side in the first convolution identity. Its limiting intensity
is therefore $\Omega_\infty^{\rm s}$. Two ordered disjoint simple
curves cut a three-vertex path. Its middle has degree two, so the
second identity gives the product intensity. Ordered edges introduce
no additional side-exchange factor. Intersecting pairs and coincident
images have $O(n)$ expectation by the one-carrier bound. The same
short-waist cutoff completes the orbital $L^1$ law for $g=0$.

Finally, every sequence $g=o(n)$ has a subsequence along which $g$
diverges or is a fixed nonnegative integer. All three cases give
the same functional $\Omega_\infty$. The subsequence principle proves
\eqref{eq:cusp-dominated-orbital-L1} in its full stated scope.

\subsection{No mass accumulates at zero}

The fixed-energy limits can be supplemented by a uniform estimate near
zero, Proposition~\ref{prop:low-energy-tail}, which is used in
Appendix~\ref{app:fixed-window-proofs}. It uses the local metric
construction in Appendix~\ref{app:cusp-cutoffs}, but not its
recurrence argument. The proof has four steps: a deterministic
low-eigenvalue count for closed surfaces of bounded geometry, a
min--max comparison between $X$ and its cusp-filled, collar-inflated
closed surface, an average of the resulting short-curve sum using the
orbital first moments of this appendix and the preceding one, and the
passage to the limiting measures.

\begin{proposition}[A uniform low-energy tail]\label{prop:low-energy-tail}
There is $t_0>0$ such that, for every $\alpha_0>0$, there are $C_{\alpha_0}$ and $n_{\alpha_0}$
with
\begin{equation}\label{eq:low-energy-tail}
 \mathbb E_{g,n}\#\{j:0<\lambda_j(X)\le t\}
       \le C_{\alpha_0}n\sqrt t,
 \qquad g\le\alpha_0 n,\quad n\ge n_{\alpha_0},\quad 0<t\le t_0.
\end{equation}
The count includes multiplicities and residual eigenvalues. In particular,
\[
 \lambda^{-1}\nu_\lambda([0,t])\le C_{\alpha_0}\sqrt t
       \quad(\lambda\ge1/\alpha_0),\qquad
 \nu_\infty([0,t])\le C_{\alpha_0}\sqrt t,
\]
and these measures have no atom at zero.
\end{proposition}

\begin{proof}
\emph{Step 1: a count for closed surfaces of bounded geometry.}
We first record a coarse deterministic estimate. If $M$ is a connected
closed surface with uniformly bounded curvature and injectivity radius
bounded below, then
\begin{equation}\label{eq:bounded-geometry-low-count}
 N_M^+(s)\le C\operatorname{area}(M)\sqrt s\qquad(0<s\le s_0),
\end{equation}
with constants depending only on those geometry bounds. Here and below
$N^+$ excludes the zero eigenvalue. Choose a maximal separated net at a
fixed sufficiently small radius $r$, and join centers of distance less
than $3r$. Its graph $G$ is connected, has bounded degree, and has at
most $C\operatorname{area}(M)$ vertices. Write $Pf$ for the averages
of $f$ on the radius-$r$ balls. Local Poincar\'e inequalities, positive
lower bounds for their volumes and bounded overlap give
\[
 \|f\|_M^2\le C\|Pf\|_{\ell^2}^2+CE_M(f),\qquad
 E_G(Pf)\le CE_M(f).
\]
For the second inequality, the two balls of an edge lie in a fixed
larger normal ball; subtract its mean before applying Poincar\'e.
For $s_0$ small, $P$ is injective on the entire spectral subspace
$[0,s]$, with graph Rayleigh quotient at most $C's$.

For completeness, take a spanning tree on the $v\ge2$ graph vertices.
Its closed contour traverses each edge twice and has length $m=2(v-1)$.
Pullback $J$ to this cycle is injective and satisfies
\[
 \|Jf\|^2=\sum_x\deg_T(x)|f(x)|^2\ge\|f\|^2,
 \qquad E_{C_m}(Jf)=2E_T(f)\le2E_G(f).
\]
Min--max and the cycle eigenvalues $4\sin^2(\pi j/m)$ imply
$N_G^+(u)\le N_{C_m}^+(2u)\le Cm\sqrt u$ for small $u$.
The two-vertex case uses two parallel cycle edges. If $v=1$, the
injection of the spectral subspace already gives $N_M^+(s)=0$.
In each comparison zero is counted first and exactly one is subtracted,
since all spaces are connected. This proves
\eqref{eq:bounded-geometry-low-count} with no additive constant.

\emph{Step 2: comparison with the filled surface.}
Apply the cusp filling and collar inflation from
Appendix~\ref{app:cusp-cutoffs} to a surface $X$ of arbitrary genus.
The construction and its local geometry bounds do not use genus zero.
The resulting closed surface $M(X)$ has the same genus, uniformly
bounded geometry, and
\[
 \operatorname{area}(M)\le C(g+n)+CZ_\delta(X),\qquad
 Z_\delta(X)=\sum_{\substack{\gamma\text{ simple}\\\ell_\gamma\le\delta}}
                  \ell_\gamma^{-1}.
\]
Dirichlet energy is unchanged by this conformal modification.
Every $f\in H^1(X)$ extends to $H^1(M)$. To check the norm needed
for comparison, on a cusp use the coordinates
$b^2(1+bu)^{-2}(du^2+d\theta^2)$, $u\ge0$. The one-dimensional
Hardy inequality gives, for almost every $\theta$,
\[
 \int_0^\infty\frac{b^2|f(u,\theta)|^2}{(1+bu)^2}\,du
 \le2b|f(0,\theta)|^2+8\int_0^\infty|\partial_uf|^2\,du.
\]
The trace is controlled by mass and energy in the fixed exterior
annulus between horocycles $b$ and $2b$, where the metrics agree.
These annuli are disjoint. Outside the chosen cusp tails new area
dominates old area, so
\begin{equation}\label{eq:cusp-regularization-norm}
 \|f\|_X^2\le c_0\|f\|_M^2+c_1E_X(f).
\end{equation}
The new cusp density decays exponentially in $u$, which makes its
$L^2$ mass finite by Cauchy--Schwarz along the $u$ coordinate. Finite
energy then extends over the filled point, whose $H^1$ capacity is
zero; truncation and logarithmic cutoffs justify this for general $f$.
For $t_0$ small, \eqref{eq:cusp-regularization-norm} and min--max on
the whole low subspace, followed by
\eqref{eq:bounded-geometry-low-count}, yield
\begin{equation}\label{eq:inverse-collar-spectral-count}
 N_X^+(t)\le N_M^+(2c_0t)
       \le C\bigl(g+n+Z_\delta(X)\bigr)\sqrt t.
\end{equation}
This comparison includes residual modes. The restriction $t<1/4$
ensures that the $X$ subspace is finite dimensional.

\emph{Step 3: averaging the short-curve sum.}
It remains to average $Z_\delta$. The fixed-support orbital first-moment
bounds in Appendix~\ref{app:positive-density-orbits} and the preceding
parts of this appendix imply
$\mathbb E C_X(\delta_0)\le C_{\alpha_0}n$ for a fixed small $\delta_0$, where
$C_X(s)$ counts simple primitive geodesics of length at most $s$.
Indeed choose an even nonnegative orbital test equal to one on
$[-\delta_0,\delta_0]$; its simple primitive weight is bounded below
there. Uniformity for diverging genus and $n/g\ge1/\alpha_0$ was proved
above. A violating bounded-genus sequence would have a fixed-genus
subsequence, contradicting its orbital first-moment limit.

In the simple-curve integration formula, let $R_{g,n}(L)$ be the sum
of the nonnegative complementary-volume factors, including the symmetry
and label factors. A cut produces two boundary copies of length $L$,
so the positive boundary-volume bound gives
\[
 R_{g,n}(0)\le R_{g,n}(L)\le e^L R_{g,n}(0),\qquad
 \mathbb E C_X(s)=\int_0^s L R_{g,n}(L)\,dL.
\]
Consequently $R_{g,n}(0)\le C_{\alpha_0}n$ and
$\mathbb E C_X(s)\le C_{\alpha_0}ns^2$ for $0<s\le\delta_0$.
Layer cake and monotone convergence now give
\[
 \mathbb E Z_\delta
 =\delta^{-1}\mathbb E C_X(\delta)
       +\int_0^\delta s^{-2}\mathbb E C_X(s)\,ds
 \le C_{\alpha_0}n\delta.
\]
Both separating and nonseparating curves are included. Substituting
in \eqref{eq:inverse-collar-spectral-count} proves
\eqref{eq:low-energy-tail}.

\emph{Step 4: the limiting measures.}
The convergence statements \eqref{eq:positive-density-L1} and
\eqref{eq:cusp-dominated-L1} are proved in
Appendix~\ref{app:fixed-window-proofs} without using the present
proposition. Apply them to continuous tests equal to one on $[0,t]$
and supported in $[0,t+d]$, and let $d\downarrow0$.
The measure bounds follow, and sending $t\downarrow0$ excludes a zero atom.
\end{proof}

Thus, for every $t_n\downarrow0$ with $g\le\alpha_0 n$,
$\mathbb E N_X^+(t_n)/n\to0$. This does not control a fixed number
of the smallest eigenvalues and is consistent with
Theorem~\ref{thm:dense-collapse}. We make no optimality claim for the
exponent in \eqref{eq:low-energy-tail}.

\section{Proofs for the fixed-window spectral measures}
\label{app:fixed-window-proofs}

This appendix proves the results of
Section~\ref{sec:fixed-energy-spectrum}:
Theorems~\ref{thm:fixed-energy-spectrum},
\ref{thm:positive-density-spectrum}, \ref{thm:cusp-dominated-spectrum}
and~\ref{thm:regularized-spectrum}, and
Corollary~\ref{cor:cusp-density-determinant}, one per subsection in
this order. The first three proofs identify a
mean through a trace formula whose geometric side is controlled by
Appendices~\ref{app:sublinear-trace}--\ref{app:cusp-dominated-orbits},
and then pass to the spectrum by positivity. The last two pass from
heat transforms to energy tests through bounds on signed total
variation and on the two ends of the heat integral. The trace inputs
are Huang's compressed trace identity and its deterministic,
retained-mass and first-moment bounds, Lemmas~\ref{rev:lem:huang-trace},
\ref{rev:lem:huang-deterministic}, \ref{rev:lem:huang-retained-mass}
and~\ref{rev:lem:huang-first-moment} of
Appendix~\ref{app:sublinear-trace}, taken from~\cite{HuangGap}; the
orbital inputs are \eqref{eq:sublinear-polynomial},
\eqref{eq:positive-orbital-L1} and \eqref{eq:cusp-dominated-orbital-L1};
the low-energy input is Proposition~\ref{prop:low-energy-tail}.
The statements and the outlines of the arguments are in
Section~\ref{sec:fixed-energy-spectrum} and are not repeated here.

\subsection{Proof of Theorem~\ref{thm:fixed-energy-spectrum}}

The theorem concerns the regime $\sqrt g\ll n\ll g$ with normalization
$d_g=n^2/g$. The upper half of the proof applies Huang's binomial
detector and compressed trace to the saddle-summed expansion
\eqref{eq:sublinear-polynomial} of
Proposition~\ref{prop:sublinear-polynomial}, and identifies the mean
measure through Proposition~\ref{prop:critical-identification}. The
lower half proves a law of large numbers for the waists of two-cusp
pairs of pants, removes the pants whose funnel models $F_s$ are not
isolated, using Lemmas~\ref{lem:critical-complex}
and~\ref{lem:critical-funnels}, and places disjoint trial states on
the rest; matching the two halves gives the $L^1$ statement.

\begin{proof}
We first identify the mean, then match it by disjoint model states.
Write $s=n+1$, $\tau=g/s$ and $D=s^2/g\sim d_g$.

\emph{The spectral mean.} Proposition~\ref{prop:sublinear-polynomial}
gives \eqref{eq:sublinear-polynomial}. Fix $0<\varepsilon<1/4$ and let
$p_Q$ be Huang's binomial detector \cite[Lemma~6.4]{HuangGap}: a real
polynomial of degree $Q$ with $p_Q(0)=0$, $p_Q(a)=1$, $0\le p_Q\le1$
on $[0,a]$,
\[
 \|p_Q/x\|_{[0,a]}\le B_\varepsilon,\qquad
 \|p_Q/x\|_{[0,1]}\le B_\varepsilon e^{-\kappa Q},\qquad
 \inf_{[f(i\sqrt\varepsilon),\,a]}p_Q\ge1-e^{-8\kappa Q}
 \qquad(Q\ge Q_\varepsilon),
\]
where $\kappa=\kappa(\varepsilon,f)>0$ is the detector decay rate
$\kappa=d^2/4$, $d=(f(i\sqrt\varepsilon)-1)/(2a)$ of
\cite[(6.14)]{HuangGap}, and $B_\varepsilon,Q_\varepsilon$ depend only
on $\varepsilon$ and the multiplier. Take
$Q=\lceil K_\varepsilon\log\tau\rceil$, where $K_\varepsilon>0$ is a
constant, depending on $\varepsilon$ and the multiplier, to be chosen
large; this is the logarithmic degree of \cite[(6.30)]{HuangGap} with
$\log g$ replaced by $\log\tau$.
Apply the compressed trace identity, Lemma~\ref{rev:lem:huang-trace},
to $p_Q$ and remove the constant mode, which has retained mass
$1-n/(\operatorname{area}(X)Y_Q)$. The remaining compression is
positive, and by Lemma~\ref{rev:lem:huang-retained-mass} every
eigenfunction with $0<\lambda\le1/4-\varepsilon$ retains mass at least
$1-e^{-2\sqrt\varepsilon LQ}$, uniformly in $n$; together with the
detector lower bound this gives
\[
 \tfrac14N_X(1/4-\varepsilon)\le
 \tau_{X,Q}(p_Q)-\Bigl(1-\frac n{\operatorname{area}(X)Y_Q}\Bigr)
 \qquad(Q\ge Q_\varepsilon).
\]
Take expectations. By Lemma~\ref{rev:lem:huang-deterministic}, the
identity and parabolic terms are at most
$C(g+nQ^2)B_\varepsilon e^{-\kappa Q}$, which after division by $D$ is
\[
 C(\tau^2+\tau Q^2)e^{-\kappa Q}
\]
and tends to zero for sufficiently large $K_\varepsilon$. The orbital
term is expanded by \eqref{eq:sublinear-polynomial}: the leading
functional satisfies $\omega_1(p_Q)=p_Q(a)+\mu_0(p_Q)$ with a fixed
signed measure $\mu_0$ on $[0,1]$ of total variation at most five
\cite[Proposition~5.1]{HuangGap}, so its atom cancels the constant-mode
correction exactly and $|\mu_0(p_Q)|\le5B_\varepsilon e^{-\kappa Q}$;
the first correction obeys the degree-independent bound
$|\omega_{2,n}(p_Q)|\le C_\varepsilon s^2$ \cite[Proposition~6.12]{HuangGap},
which holds without a genus hypothesis; and the error terms in
\eqref{eq:sublinear-polynomial} are $o(D)$. The removed-mass term
$n/(\operatorname{area}(X)Y_Q)\le e^{-LQ}/(6\pi)$ is negligible.
Hence
\begin{equation}\label{eq:fixed-energy-count-bound}
 \mathbb E N_X(1/4-\varepsilon)\le C_\varepsilon D.
\end{equation}
This is the argument of \cite[Theorem~6.13]{HuangGap} with its
polynomial-transfer input replaced by \eqref{eq:sublinear-polynomial};
the constant mode is removed before applying the assertion.

For a smooth multiplier test $\psi\in C_c^\infty((1,\infty))$,
use instead $Q=\lfloor\tau^{\theta_0}\rfloor$, with $\theta_0$ as in
Proposition~\ref{prop:sublinear-polynomial}. Approximate $\psi/x$
by its shifted Chebyshev truncation $r_Q$ and put $p_Q=xr_Q$.
For every fixed $K$ the error is $O(Q^{-K})$, with convergence in
the fixed $C^m$ norm controlling the canonical distributions.
The normalized bulk error is
\[
 O\bigl((\tau^2+\tau Q^2)Q^{-K}\bigr)=o(1)
 \quad\text{if }K>\max(2/\theta_0,2+1/\theta_0).
\]
A logarithmic degree would not suffice for this smooth approximation.
The base positive-trace moment is $O(g+nQ)$, by
Lemma~\ref{rev:lem:huang-first-moment}. Remove the constant
projection exactly, as in Lemma~\ref{rev:lem:huang-trace}; its
cusp-truncation correction
$np_Q(a)/(\operatorname{area}(X)Y_Q)$ is negligible after division
by $D$. The retained-mass error is exponentially small in $Q$ by
Lemma~\ref{rev:lem:huang-retained-mass}, uniformly in $n$, and
\eqref{eq:fixed-energy-count-bound} removes the weights in expectation.

The exact identity
$\omega_{2,n}=W_0+nW_1+\binom n2W_2$ now leaves only $W_2/2$,
since $gD=s^2$. Proposition~\ref{prop:critical-identification}
identifies this distribution with \eqref{eq:fixed-energy-measure}.
The same smooth-majorant argument as in
Proposition~\ref{prop:critical-mean} includes the endpoint of zero
energy, after the zero eigenfunction has been excluded. Hence
\begin{equation}\label{eq:fixed-energy-mean}
 d_g^{-1}\mathbb E\Xi_X\longrightarrow\nu_*
 \quad\text{vaguely on }[0,1/4).
\end{equation}

\emph{Pants concentration.} Let $\mathcal Q_X$ count the waists of
two-cusp pairs of pants, as in Appendix~\ref{app:critical}.
The one-pants integration formula and the zero-density boundary profile
give $d_g^{-1}\mathbb E\mathcal Q_X(I)\to\mu_*(I)$ for bounded
intervals $I$. The ordered contribution of two disjoint pants is
\[
 \frac{(n)_4}{4V_{g,n}}
 \int_{I\times J}V_{g,n-2}(s,t,0^{n-4})\,st\,ds\,dt.
\]
Dividing by $d_g^2$, the volume ratio and boundary factorization
give the limit $\mu_*(I)\mu_*(J)$.

Intersecting pairs have expected count $o(d_g)$. Indeed the filling
construction in the proof of Lemma~\ref{lem:critical-complex}
produces finitely many carrier types of bounded complexity and bounded
boundary length. A pair of distinct intersecting pants waists cannot
fill a complexity-one piece. For every other type,
\eqref{eq:sublinear-fixed-carrier} gives $O(d_g(n/g)^{m-1})=o(d_g)$.
Lemma~\ref{lem:primitive-count} bounds the number of pairs producing
each marked cut, so this is a count bound, not merely an existence
probability. The diagonal contributes only $O(d_g)$. Thus
$\mathcal Q_X(I)/d_g\to\mu_*(I)$ in $L^2$.
The same argument applies to bounded measurable waist sets with
$\mu_*$-null boundary. In particular, with
$Z_X(t)=\#\{s\in\mathcal Q_X:\beta(s)\le t\}$,
\begin{equation}\label{eq:fixed-energy-pants-LLN}
 Z_X(t)/d_g\longrightarrow\nu_*([0,t])\quad\text{in }L^2
\end{equation}
at every continuity threshold $t<1/4$. The long-waist lemma bounds
the relevant waist set; the first-moment bound near zero is $O(d_gd^2)$
for waists below $d$.

\emph{Removing nonisolated models.} Fix $0<d<S$ and a finite funnel
radius $R$. From the pants with waists in $[d,S]$, one can delete a
number $B_X$ with $\mathbb E B_X=o(d_g)$ so that their radius-$R$
extensions are embedded and pairwise disjoint. To prove the count
assertion, first delete pants in intersecting pairs. A return or
collision of normal collars before $R$ then supplies an essential arc
of length at most $2R$, as in Lemma~\ref{lem:critical-funnels}.
One pants and a returning arc produce type $(0,2,2)$, or $(0,3,1)$
after filling a peripheral disk; two pants joined by an arc produce
type $(0,4,1)$. Repeated boundary incidences are retained as sides of
one cut edge. These have complexity at least two and total cut length
at most $4S+8R+1$. Their expected counts are $o(d_g)$ by
\eqref{eq:sublinear-fixed-carrier}.
Each marked cut contains only boundedly many possible original waists
of length at most $S$, by Lemma~\ref{lem:primitive-count}.
Consequently the number of implicated pants is bounded by a constant
times the marked cut count. Delete all of them; the multiplicity of
arcs realizing one cut is irrelevant. This proves the assertion.

\emph{Exhaustion.} Fix a continuity threshold $t<1/4$, choose
$\eta>0$ with $t-\eta$ also a continuity threshold, and restrict to
waists $s\ge d$ with $\beta(s)\le t-\eta$. Their upper waist bound
is deterministic. The compact-family trial-state construction in
Section~\ref{sec:critical-eigenvalues} gives Rayleigh quotients at
most $t-\eta/2$, with supports in a common fixed funnel radius and
bounded cusp height. After the deletions these supports are disjoint.
Min--max, allowing for the global constant, gives
\[
 N_X(t)\ge Z_X([d,\infty);t-\eta)-B_X-1.
\]
The loss of one mode vanishes on division by $d_g\to\infty$.
It follows that
\[
 \limsup\frac{\mathbb E(Z_X(t)-N_X(t))_+}{d_g}
 \le\nu_*([0,t])-\mu_*\{s\ge d:\beta(s)\le t-\eta\}.
\]
Let $d\downarrow0$ and then $\eta\downarrow0$ through continuity
thresholds. The right side vanishes. The equal limiting means from
\eqref{eq:fixed-energy-mean} and the pants first moment now imply
\begin{equation}\label{eq:fixed-energy-count-matching}
 \mathbb E|N_X(t)-Z_X(t)|/d_g\longrightarrow0.
\end{equation}
Together with \eqref{eq:fixed-energy-pants-LLN}, this proves $L^1$
convergence at continuity thresholds. Approximation by step functions,
using the local mean bound, proves \eqref{eq:fixed-energy-L1} for
every continuous compact test and also vague convergence in probability.
\end{proof}

\subsection{Proof of Theorem~\ref{thm:positive-density-spectrum}}

Here $n/g\to\lambda\in(0,\infty)$ and the normalization is $g$. The
proof transfers the orbital limit \eqref{eq:positive-orbital-L1} of
Appendix~\ref{app:positive-density-orbits} to the spectrum through
Huang's positive compressed trace~\cite[Sections~5--6]{HuangGap}, in
the form of Lemmas~\ref{rev:lem:huang-trace},
\ref{rev:lem:huang-deterministic} and~\ref{rev:lem:huang-retained-mass};
the first-moment bound of Lemma~\ref{rev:lem:huang-first-moment},
whose genus hypothesis is not available here, is replaced by
\eqref{eq:positive-trace-first-moment}.
Since the normalization is the genus, the surface limit is taken at a
fixed polynomial degree $Q$ and the degree is sent to infinity
afterwards; the limiting measure is then produced by the Riesz
representation theorem. The three stages are a local expected-count
bound, the transfer of smooth multiplier tests, and the construction
and characterization of $\nu_\lambda$.

\begin{proof}
\emph{The local count bound.}
Use the multiplier $h,f,a$ of Proposition~\ref{prop:critical-mean},
with $\operatorname{supp}h\subset[-L,L]$. For
$p(z)=\sum_{j\ge1}p_jz^j$ put $H_p=\sum_jp_jh^{*j}$.
At $Y_Q=3e^{LQ}$ the exact compressed trace identity
\eqref{rev:eq:huang-trace} of Lemma~\ref{rev:lem:huang-trace} reads
\[
 \tau_{X,Q}(p)=\bigl(\operatorname{area}(X)-n/Y_Q\bigr)\mathcal I(p)
          +n\mathcal P_{Y_Q}(p)+\mathcal O_X(H_p);
\]
it holds for every surface, with no genus or cusp condition. For
$p=z\widetilde p$ of degree at most $Q$, the deterministic bounds of
Lemma~\ref{rev:lem:huang-deterministic} are
\[
 |\mathcal I(p)|\le C\|\widetilde p\|_{[0,1]},\qquad
 |\mathcal P_{Y_Q}(p)|\le CQ^2\|\widetilde p\|_{[0,1]},\qquad
 |\mathcal P_{Y_Q}(z)|\le CQ,
\]
with $C$ depending only on the multiplier.
Fix a compact density interval $K\subset(0,\infty)$. The first-moment
bound for the fixed function $h$ in Appendix~\ref{app:positive-density-orbits}, inserted in the
identity as in Lemma~\ref{rev:lem:huang-first-moment}, gives
\begin{equation}\label{eq:positive-trace-first-moment}
 \mathbb E\tau_{X,Q}(z)\le C_Kg(1+Q),\qquad n/g\in K,
\end{equation}
for all sufficiently large $g$ and every $Q\ge1$.
For each $t<1/4$, Lemma~\ref{rev:lem:huang-retained-mass} with
$\varepsilon=1/4-t$ provides a fixed $Q_0$ for which every
eigenfunction of energy in $(0,t]$ retains mass at least $1/2$ below
height $Y_{Q_0}$, uniformly in the surface and in its cusp number.
Positivity then implies
\begin{equation}\label{eq:positive-local-count}
 \mathbb E N_X(t)\le C_{K,t}g.
\end{equation}

\emph{Transfer of smooth tests.}
Fix $\psi\in C_c^\infty((1,\infty))$, allowing $\psi(a)\ne0$.
Approximate $\psi(z)/z$ on $[0,a]$ by polynomials $r_Q$ of degree
$Q-1$, with error $e_Q=O_A(Q^{-A})$ for every fixed $A$, and put
$p_Q=zr_Q$. Positivity bounds the compressed-trace approximation
error by $e_Q\tau_{X,Q}(z)$. Since $\psi$ vanishes on $[0,1]$,
the identity and parabolic terms cost at most $C_KgQ^2e_Q$.
The spectral representation of $\tau_{X,Q}(\psi)$ has only
exceptional discrete terms and the constant mode. Remove the latter;
its normalized contribution is at most $\|\psi\|_\infty/g$.
On the remaining support, the missing eigenfunction mass is
$O(e^{-c_\psi Q})$, uniformly in the surface and its cusp count, by
Lemma~\ref{rev:lem:huang-retained-mass} with
$\sqrt\varepsilon$ determined by the distance of
$\operatorname{supp}\psi$ from the multiplier value $1$.
Equation~\eqref{eq:positive-local-count} removes these weights in mean.
Together with \eqref{eq:positive-orbital-L1}, this gives
\begin{align}\label{eq:positive-transfer-error}
 \limsup_{\substack{g\to\infty\\ n/g\to\lambda}}
 \mathbb E\left|\frac1g\sum_{0<\lambda_j<1/4}
 \psi\bigl(f(i\sqrt{1/4-\lambda_j})\bigr)
              -\Omega_\lambda(H_{p_Q})\right|
 \le C_{K,\psi}\bigl((1+Q^2)e_Q+e^{-c_\psi Q}\bigr).
\end{align}
Here $Q$ is fixed before the genus limit.

\emph{The limiting measure.}
The right side of \eqref{eq:positive-transfer-error} tends to zero,
so the numbers $\Omega_\lambda(H_{p_Q})$ form a Cauchy sequence.
Their limit is positive and linear, locally bounded by
\eqref{eq:positive-local-count}. The map
$E\mapsto f(i\sqrt{1/4-E})$ is a smooth diffeomorphism away from
$1/4$. Extending a test past its endpoint at $E=0$ therefore realizes
every smooth compact physical-energy test by such a $\psi$.
Riesz representation gives a unique Radon measure with
\begin{equation}\label{eq:positive-spectral-characterization}
 \int\psi\bigl(f(i\sqrt{1/4-E})\bigr)\,d\nu_\lambda(E)
       =\lim_{Q\to\infty}\Omega_\lambda(H_{p_Q}).
\end{equation}
The local count bound and uniform approximation extend the resulting
$L^1$ convergence to every $\varphi\in C_c([0,1/4))$.
All bounds are uniform on $K$, and each fixed-$Q$ functional is
continuous in $\lambda$. The same approximation proves vague
continuity of the family. The physical spectral limit also shows
that the characterization is independent of the admissible multiplier
and the chosen polynomial approximations.
The absence of a zero atom follows from \eqref{eq:low-energy-tail},
which does not depend on the present theorem: by the portmanteau
inequality for the relatively open set $[0,t)$,
$\nu_\lambda([0,t))\le\liminf g^{-1}\mathbb E\,\#\{j:0<\lambda_j(X)\le t\}
\le C\lambda\sqrt t$ for $t\le t_0$.
\end{proof}

\subsection{Proof of Theorem~\ref{thm:cusp-dominated-spectrum}}

Here $g=o(n)$, including bounded genus, and the normalization is $n$.
The argument is that of the preceding subsection with the orbital
input \eqref{eq:cusp-dominated-orbital-L1} of
Appendix~\ref{app:cusp-dominated-orbits} in place of
\eqref{eq:positive-orbital-L1}; the only point to check is that each
trace bound is uniform in the cusp number, which the lemmas of
Appendix~\ref{app:sublinear-trace} provide. The continuity statement
\eqref{eq:cusp-dominated-continuity} follows from the uniform endpoint
bounds of Appendix~\ref{app:cusp-dominated-orbits} and the functional
limit \eqref{eq:cusp-dominated-functional}.

\begin{proof}
Apply the proof of Theorem~\ref{thm:positive-density-spectrum} with
normalization $n$, using \eqref{eq:cusp-dominated-orbital-L1} in place
of \eqref{eq:positive-orbital-L1}. The area is $O(n)$, the parabolic term has factor
$n$, and \eqref{eq:cusp-dominated-orbital-L1} bounds the fixed-$h$
orbital first moment by $O(n)$. Thus
$\mathbb E\tau_{X,Q}(z)\le Cn(1+Q)$ and
$\mathbb E N_X(t)\le C_tn$ for each $t<1/4$.
For a smooth multiplier test $\psi$ and its approximants $p_Q=zr_Q$,
the exact trace identity (Lemma~\ref{rev:lem:huang-trace}), the
deterministic bounds (Lemma~\ref{rev:lem:huang-deterministic}),
polynomial approximation, and the retained-mass bound
(Lemma~\ref{rev:lem:huang-retained-mass}), which is uniform in the
cusp number and therefore applies when $g=o(n)$, give
\begin{align}\label{eq:cusp-dominated-transfer}
 \limsup_{\substack{n\to\infty\\g=o(n)}}
 \mathbb E\left|\frac1n\sum_{0<\lambda_j<1/4}
 \psi\bigl(f(i\sqrt{1/4-\lambda_j})\bigr)
                        -\Omega_\infty(H_{p_Q})\right|
 \le C_\psi\bigl((1+Q^2)e_Q+e^{-c_\psi Q}\bigr).
\end{align}
The constant eigenfunction contributes $O(1/n)$ and is excluded.
First take the surface limit with $Q$ fixed, then let $Q\to\infty$.
Positivity, the local expected-count bound and uniform approximation
give the unique Radon measure and \eqref{eq:cusp-dominated-L1}, exactly
as in \eqref{eq:positive-spectral-characterization}.

The endpoint estimates in Appendix~\ref{app:cusp-dominated-orbits} also bound fixed-support
orbital first moments divided by $n$ uniformly for diverging genus
and $n/g\ge\lambda_0>0$. Since $\operatorname{area}(X)/n$ is bounded
there, the trace approximation error for $\nu_\lambda/\lambda$ is
uniform in $\lambda\ge\lambda_0$. At fixed $Q$, use
\eqref{eq:cusp-dominated-functional}, and then let $Q\to\infty$.
This proves \eqref{eq:cusp-dominated-continuity}. Neither passage
allows the spectral support to approach $1/4$.
The bound \eqref{eq:low-energy-tail} excludes a zero atom, by the same
portmanteau argument as for $\nu_\lambda$.
\end{proof}

\subsection{Proof of Theorem~\ref{thm:regularized-spectrum}}

The theorem concerns the signed measure $\mathfrak m_X$ of
\eqref{eq:regularized-measure}, with $g/n\to\theta\in[0,\infty)$ and
normalization $n$. The starting point is the exact regularized trace
formula \eqref{eq:regularized-exact-heat}, whose orbital side is
controlled by \eqref{eq:positive-orbital-L1} and
\eqref{eq:cusp-dominated-orbital-L1} after long repetitions are cut
off by a packing bound in the thick part. Continuous energy tests are
recovered from heat transforms through a bound on the signed total
variation of $\mathfrak m_X$, which rests on the residual-pole
decomposition of the scattering determinant in the normalization
of~\cite{HTMany} and the Otal--Rosas bound~\cite{OtalRosas}. The
three stages are marked below.

\begin{proof}
We first control long repetitions, then pass to heat traces, and finally
bound signed total variation to recover continuous energy tests.

\emph{A length tail independent of the systole.}
Fix a sufficiently small Margulis radius $d>0$. Every primitive closed
geodesic other than a core of a $d$-thin annulus meets the $d$-thick
part: a closed geodesic contained in a thin component can only be
such a core. A maximal $d/2$-separated net in the thick part has at
most $C_d\operatorname{area}(X)$ points, since its embedded $d/4$-balls
are disjoint. A geodesic of length at most $B$ meeting that part gives
a based loop of length at most $B+d$ at a net point. Its based
homotopy class is an element of the covering group. Different free
homotopy classes give distinct conjugacy classes, so a choice of
representatives injects these classes into the union of the based-loop
sets. The orbit above each net point is $2d$-separated. Hyperbolic
ball packing therefore bounds the number of these oriented primitive
geodesics by $C_d\operatorname{area}(X)e^B$.

All these primitives have length at least $2d$. Including repetitions
with $m\ell\le B$ multiplies the bound by at most $C_d(1+B)$.
For $z$ in a compact subset $K$ of the right half-plane,
$|G_z(u)|\le C_K e^{-c_Ku^2}$. In a unit bin
$u=m\ell\in[j,j+1]$, $j\ge1$, the orbital weight satisfies
\[
 \frac{\ell}{2\sinh(m\ell/2)}\le C(j+1)e^{-j/2}.
\]
The contribution of these bins with $u\ge B$ is consequently bounded
by $C_K\operatorname{area}(X)e^{-c'_KB^2+C_KB}$.

There are at most $C\operatorname{area}(X)$ short cores. For each
one of length $\ell<2d$ and $B\ge2$, its long repetitions satisfy
\[
 \sum_{m\ell\ge B}\ell e^{-m\ell/2-c_K(m\ell)^2}
 \le C_K\int_{B-2d}^\infty e^{-u/2-c_Ku^2}\,du.
\]
This bound is uniform as $\ell\downarrow0$. If $\chi_B$ is an even
smooth cutoff, equal to one on $[-B,B]$, zero outside
$[-B-1,B+1]$ and between zero and one, we have proved
\begin{equation}\label{eq:regularized-orbital-tail}
 \sup_{z\in K}|\mathcal O_X((1-\chi_B)G_z)|
 \le C_K\operatorname{area}(X)e^{-c'_KB^2+C_KB}.
\end{equation}
The cutoff is on repetition length, not just primitive length.
Short repetitions, which can have a logarithmic dependence on the
systole, are left in the compact-support orbital statistic.

\emph{The heat limit.}
We use the exact regularized trace formula
\eqref{eq:regularized-exact-heat}.\footnote{The constant-term
convention in \eqref{eq:regularized-measure} and
\eqref{eq:regularized-exact-heat} is that of
\cite[(3.1)]{HTMany}: on the geometric side the threshold term is
$\frac{h(0)}4\operatorname{Tr}(I-\Phi_X(\frac12))$, which moved to the
spectral side gives
$-\frac{n-\operatorname{tr}\Phi_X(1/2)}4F(1/4)$. The same identity
is also written, compare \cite{Iwaniec,Hejhal}, with
$+\frac{h(0)}4\operatorname{Tr}\Phi_X(\frac12)$ on the spectral side
and $+\frac n4h(0)$ on the geometric side; the two forms differ by
moving the term $\frac{h(0)}4\operatorname{Tr}\Phi_X(\frac12)$ across
the identity. All our statements use the form of
\cite[(3.1)]{HTMany}. Since $\Phi_X(\frac12+ir)$ is unitary,
$\phi_X(\frac12+ir)$ has modulus one and
$\frac{\phi_X'}{\phi_X}(\frac12+ir)=\frac{d}{dr}\arg\phi_X(\frac12+ir)$
is real, so the scattering integral in \eqref{eq:regularized-measure}
is a real signed measure.}
For fixed $B$, \eqref{eq:positive-orbital-L1} and
\eqref{eq:cusp-dominated-orbital-L1} give
\[
 \mathbb E|n^{-1}\mathcal O_X(\chi_BG_z)
                    -\mathcal G_\theta(\chi_BG_z)|\longrightarrow0.
\]
Since $\operatorname{area}(X)/n\to2\pi(1+2\theta)$,
\eqref{eq:regularized-orbital-tail} permits $B\to\infty$ after the
surface limit. This defines $\mathcal G_\theta(G_z)$ and proves
pointwise $L^1$ convergence in \eqref{eq:regularized-limit-heat}.

On a fixed compact set of complex times, the short-length part is
bounded by a constant times one nonnegative smooth compact test $H_0$.
The orbital first moments of $H_0$, divided by $n$, are uniformly
bounded by Appendices~\ref{app:positive-density-orbits}--
\ref{app:cusp-dominated-orbits}. Together with
\eqref{eq:regularized-orbital-tail}, this gives uniform integrable
domination on any fixed contour in the right half-plane. Cauchy's
formula applied on surrounding contours proves local uniform $L^1$
convergence and convergence of all fixed derivatives, as asserted
in \eqref{eq:regularized-analytic-L1}. The same argument extends
the fixed-support continuity of $\mathcal G_\theta$ to these heat
transforms, uniformly on compact intervals of $\theta$.

\emph{Signed total variation and continuous tests.}
We use a classical scattering estimate, with its signs and normalization
as in \cite[Lemma~3.3 and Section~5]{HTMany}. If $s_j\in(1/2,1]$
are the residual poles of $\phi_X$, counted with their determinant
orders, and $d_j=s_j-1/2$, their pole decomposition gives
\begin{equation}\label{eq:regularized-negative-phase}
 -\frac{\phi_X'}{\phi_X}(1/2+ir)
 =2\log b_1+S_+(r)-\sum_{j=1}^{M_X}\frac{2d_j}{d_j^2+r^2},
 \qquad b_1\ge1,\quad S_+(r)\ge0.
\end{equation}
Here $b_1$ is the first denominator in the scattering determinant's
Dirichlet series. The width-one normalization is essential for
$b_1\ge1$; see also \cite[Lemma~5.2]{HuangGap}.
The pole orders are bounded by the corresponding residual eigenspace
dimensions. Indeed, on $\operatorname{Re}s>1/2$ the scattering matrix
$\Phi_X(s)$ is meromorphic with simple poles, located exactly at the
points $s_j\in(1/2,1]$ with $s_j(1-s_j)$ a residual eigenvalue, and
the residue of $\Phi_X$ at $s_j$ is, up to a nonzero scalar, the Gram
matrix of the zeroth Fourier coefficients of an orthonormal basis of the
residual eigenspace at $s_j(1-s_j)$, so its rank is at most that
eigenspace dimension; see \cite[Chapter~6]{Iwaniec} and
\cite[Chapter~6]{Hejhal}. Writing $\Phi_X(s)=(s-s_j)^{-1}R_j+O(1)$ with
$\operatorname{rank}R_j=\varrho_j$, the determinant
$\phi_X=\det\Phi_X$ has a pole of order at most $\varrho_j$ at $s_j$:
after conjugation by a constant invertible matrix, which does not change
the determinant, $R_j$ is supported on $\varrho_j$ rows, only these rows
of $\Phi_X$ have a pole at $s_j$, each of order one, and the determinant
is multilinear in the rows. The pole at $s=1$, which corresponds to the constant
eigenfunction and has residue of rank one, is included in the count.
Hence $M_X\le2g+n-2$ by the Otal--Rosas bound~\cite{OtalRosas}, in the form recalled in
\cite[Theorem~1.1]{HTMany}. Each displayed negative pole kernel has integral $2\pi$,
even when $d_j\downarrow0$. The threshold atom in
\eqref{eq:regularized-measure} has mass at least $-n/2$, and the
digamma terms have heat-weighted absolute mass at most $C_tn$.
Consequently, for every $t>0$,
\begin{align}\label{eq:regularized-variation}
 \int e^{-tE}\,d\mathfrak m_X^-(E)&\le C_t(g+n),\notag\\
 \frac1n\int e^{-tE}\,d|\mathfrak m_X|(E)
 &\le\frac{\mathcal R_X(t)}n+C_t(1+g/n).
\end{align}
The second inequality is the identity
$|\mathfrak m_X|=\mathfrak m_X+2\mathfrak m_X^-$.
Its expected right side is uniformly bounded by the heat limit.
Thus no pointwise estimate at $r=0$ is needed.

Fix $t_0>0$ and push $e^{-t_0E}\mathfrak m_X/n$ forward under
$u=e^{-E}$ to $[0,1]$. These signed measures have uniformly bounded
expected total variation. Their moments of every integer order
$k\ge0$ converge in $L^1$ to $\mathcal R_\theta(t_0+k)$.
The limiting functional on polynomials is bounded in the uniform
norm, so polynomial approximation and the Riesz representation theorem
give a unique finite signed measure $\mu$ on $[0,1]$.
For a continuous $f$, approximation by a polynomial $p$ gives an error
at most the polynomial $L^1$ error plus
\[
 \|f-p\|_\infty
 \left(\mathbb E\|e^{-t_0E}\mathfrak m_X/n\|_{\mathrm{TV}}
                            +\|\mu\|_{\mathrm{TV}}\right).
\]
Hence convergence holds in $L^1$ for every continuous test.
There is no mass at $u=0$, because
\begin{equation}\label{eq:regularized-energy-tail}
 \mathbb E\frac1n\int_{[B,\infty)}e^{-t_0E}\,d|\mathfrak m_X|(E)
 \le e^{-t_0B/2}\mathbb E\frac1n
       \int e^{-t_0E/2}\,d|\mathfrak m_X|(E)
 \le C_{t_0}e^{-t_0B/2}.
\end{equation}
Continuous cutoffs pass this bound to $\mu$.

Pull $\mu$ back to $[0,\infty)$ and multiply by $e^{t_0E}$ on
compact intervals. This defines the signed Radon measure
$\mathfrak m_\theta$ and proves \eqref{eq:regularized-L1}.
It is independent of $t_0$ by uniqueness of the $L^1$ limit for each
compact test. Taking $t_0<t$ recovers its heat transform at $t>0$;
holomorphy gives \eqref{eq:regularized-limit-heat} for complex times.
The same polynomial argument proves uniqueness from the heat transform.
Uniform bounds for $\theta$ in a compact interval and continuity of
$\mathcal R_\theta(t)$ imply vague continuity of the measures.
Finally, on $[0,1/4)$ only discrete eigenvalues enter
\eqref{eq:regularized-measure}, and the zero mode contributes $1/n$.
The asserted restrictions follow from Theorems
\ref{thm:positive-density-spectrum}--\ref{thm:cusp-dominated-spectrum}.
\end{proof}

\subsection{Proof of Corollary~\ref{cor:cusp-density-determinant}}

The corollary expresses $\log Z_X'(1)$ and the geometric log
determinant \eqref{eq:geometric-determinant-definition} through the
integral of the orbital heat term $\mathcal B_X(t)$ of
\eqref{eq:determinant-orbital-limit} against $dt/t$.
Theorem~\ref{thm:regularized-spectrum} supplies convergence on compact
time intervals, so the proof consists of averaged bounds at the two
ends of the time axis, uniform for $g\le\alpha_0 n$: at small times the
short-geodesic intensity bound from the proof of
Proposition~\ref{prop:low-energy-tail}, at large times its low-energy
tail \eqref{eq:low-energy-tail} together with the variation bound
\eqref{eq:regularized-variation}. The remaining stages identify the
limiting integral with the determinant by the Mellin transform and
prove its strict positivity from short two-cusp pairs of pants.

\begin{proof}
We control both ends of the heat integral, then apply
Theorem~\ref{thm:regularized-spectrum}. This is the classical
heat-to-determinant method used in~\cite{NaudDet,HeWuDet,GLNDet}.
All averaged bounds below are uniform for $g\le\alpha_0 n$ and sufficiently
large $n$, with constants depending on $\alpha_0$.

\emph{Small times.}
Choose $d>0$ below the collar threshold. The proof of
Proposition~\ref{prop:low-energy-tail} gives the intensity bound
$L R_{g,n}(L)\,dL\le C_{\alpha_0}nL\,dL$ for all simple primitives of
length at most $d$, including every separating topology.
Every primitive of length below $d$ is simple. The exact orbital sum is
\[
 \mathcal B_X(t)=\frac{e^{-t/4}}{\sqrt{4\pi t}}
 \sum_{\gamma\text{ primitive oriented}}\sum_{m\ge1}
 \frac{\ell_\gamma}{2\sinh(m\ell_\gamma/2)}
                      e^{-m^2\ell_\gamma^2/(4t)}.
\]
Since $\ell/(2\sinh(m\ell/2))\le1/m$, its short-primitive part satisfies
\[
 \mathbb E\mathcal B_X^{\mathrm{short}}(t)
 \le\frac{C_{\alpha_0}n}{\sqrt t}\sum_{m\ge1}\frac1m
                    \int_0^d L e^{-m^2L^2/(4t)}\,dL
 \le C_{\alpha_0}n\sqrt t\sum_{m\ge1}m^{-3}.
\]
Thus all short repetitions are included. For the remaining primitives,
the thick-part packing argument in the proof of
Theorem~\ref{thm:regularized-spectrum} bounds their count up to length
$L$ by $C_dA_Xe^L$ (applied with Margulis radius $d/2$, so that every
primitive of length at least $d$ meets the thick part). Their lengths
are at least $d$, so counting
repetitions costs at most another factor $C_d(1+L)$.
Inserting this bound in the Gaussian sum gives
$\mathcal B_X^{\mathrm{long}}(t)\le C_dA_Xt^{-C_d}e^{-c_d/t}$
for $0<t\le1$: the linear exponential from counting is absorbed into
half of the negative quadratic exponent, whose length tail starts at $d$.
Consequently
\begin{equation}\label{eq:determinant-small-time}
 \frac1n\mathbb E\mathcal B_X(t)\le C_{\alpha_0}\sqrt t
          \quad(0<t\le1),\qquad
 \mathbb E\int_0^\varepsilon\frac{\mathcal B_X(t)}n\frac{dt}{t}
                       \le C_{\alpha_0}\sqrt\varepsilon.
\end{equation}

\emph{Large times.}
Equation~\eqref{eq:regularized-variation}, the heat identity, and
\eqref{eq:determinant-small-time} at $t=1$ give
\[
 \mathbb E\int_0^\infty e^{-E}\,d|\mathfrak m_X|(E)\le C_{\alpha_0}n.
\]
Fix $0<v_0\le\min(t_0,1/8)$, with $t_0$ as in
Proposition~\ref{prop:low-energy-tail}. On $(0,v_0)$ only positive
discrete eigenvalues occur. Integration by parts using
$\mathbb E N_X^+(v)\le C_{\alpha_0}n\sqrt v$ gives, for $t\ge2$,
\[
 \mathbb E\sum_{0<\lambda_j<v_0}e^{-t\lambda_j}
 \le C_{\alpha_0}n\left(e^{-tv_0}
             +t\int_0^{v_0}e^{-tv}\sqrt v\,dv\right)
 \le C_{\alpha_0}n t^{-1/2}.
\]
The complementary signed measure is supported on $[v_0,\infty)$;
its expected absolute heat contribution is at most
$C_{\alpha_0}n e^{-(t-1)v_0}$. The atom at zero has mass exactly one.
Since $k_{\mathbb H}(t)$ decays exponentially, we obtain
\begin{equation}\label{eq:determinant-large-time}
 \frac1n\mathbb E|\mathcal B_X(t)-1|\le C_{\alpha_0}t^{-1/2}\quad(t\ge2),
 \qquad
 \mathbb E\int_T^\infty\frac{|\mathcal B_X(t)-1|}{n}\frac{dt}{t}
                       \le C_{\alpha_0}T^{-1/2}.
\end{equation}
In particular no lower bound uniform in $X$ is imposed on $\lambda_1$.

\emph{The limiting integral.}
Compact-time $L^1$ convergence in
\eqref{eq:regularized-analytic-L1} and the two endpoint bounds imply
$0\le\mathcal B_\theta(t)\le C_{\alpha_0}\min(\sqrt t,t^{-1/2})$.
Thus $I_\theta$ is finite. Put
\[
 J_X=\int_0^1\mathcal B_X(t)\frac{dt}{t}
                  +\int_1^\infty(\mathcal B_X(t)-1)\frac{dt}{t}.
\]
First restrict to $[\varepsilon,T]$ and take the surface limit.
The subtracted constant contributes only $\log T/n$ there.
Then \eqref{eq:determinant-small-time}--\eqref{eq:determinant-large-time}
permit $\varepsilon\downarrow0$ and $T\to\infty$, proving
$\mathbb E|J_X/n-I_\theta|\to0$.
The zero mode must be subtracted before integrating to infinity.
The same bounds, uniform on compact $\theta$ intervals, and continuity
of $\mathcal B_\theta(t)$ prove continuity of $I_\theta$.

\emph{Identification with the determinant.}
For real $s>1$, put $q=s(s-1)$. Tonelli and the elementary Gaussian
integral give
\[
 \int_0^\infty e^{-qt}\mathcal B_X(t)\frac{dt}{t}
 =\sum_\gamma\sum_{m\ge1}
       \frac{e^{-sm\ell_\gamma}}{m(1-e^{-m\ell_\gamma})}
 =-\log Z_X(s).
\]
For each fixed surface, $\mathcal B_X(t)$ is exponentially small as
$t\downarrow0$, and $\mathcal B_X(t)-1$ decays exponentially as
$t\to\infty$. The latter follows from its positive nonzero spectral
gap and heat-weighted total variation, with the continuous spectrum
starting at $1/4$. Since
$\int_1^\infty e^{-qt}\,dt/t=-\gamma_{\mathrm E}-\log q+o(1)$ as
$q\downarrow0$, where $\gamma_{\mathrm E}=0.5772\ldots$ is the
Euler--Mascheroni constant,
letting $s\downarrow1$ yields
\begin{equation}\label{eq:selberg-derivative-heat}
 \log Z_X'(1)=\gamma_{\mathrm E}-J_X,
             \qquad Z_X'(1)>0.
\end{equation}
This proves the first limit. To check the second, define $J_X(w)$
by replacing $dt/t$ in both terms of $J_X$ by $t^{w-1}dt$.
The single-surface decay makes $J_X(w)$ entire. The Mellin continuation
of $k_{\mathbb H}$ gives $\zeta_{\mathbb H}$ with
$-\zeta_{\mathbb H}'(0)=E_{\mathbb H}$, and
\[
 \zeta_X^{\mathrm{geom}}(w)
 =A_X\zeta_{\mathbb H}(w)+\frac{J_X(w)-1/w}{\Gamma(w)}.
\]
Using $1/\Gamma(w)=w+\gamma_{\mathrm E}w^2+O(w^3)$, this gives exactly
\[
 \log\det{}'_{\mathrm{geom}}\Delta_X
                  =A_XE_{\mathbb H}+\log Z_X'(1).
\]
Now $A_X/n\to a_\theta$, which proves the second limit.

Finally, let $C_2$ count waists bounding two cusps with length in a fixed
short interval $[\ell_0,\ell_1]\subset(0,d)$. Mirzakhani integration and
positivity of boundary coefficients give, for $n$ large,
\[
 \frac{\mathbb E C_2}{n}\ge
 \frac{n-1}{2}\frac{V_{g,n-1}}{V_{g,n}}
                         \int_{\ell_0}^{\ell_1}L\,dL.
\]
The right side has a strictly positive limit by
Theorem~\ref{thm:volume}; for bounded genus use
\cite[Theorem~6.1]{ManZ}, with the same endpoint ratio.
At each fixed $t>0$ the $m=1$ heat weight has a positive minimum on
this interval. Compact-time convergence therefore implies
$\mathcal B_\theta(t)>0$, and hence $I_\theta>0$.
\end{proof}

\part*{Appendices to Part~3: local geometry}
\section{Rooted holonomies and the local WP law}
\label{app:rooted-local-law}

This appendix proves Theorem~\ref{thm:wp-local-family}
(Section~\ref{sec:local-conformal-phase}): the fixed-radius local
limit $\mathcal L_\lambda$ of a WP surface at positive cusp density,
its unimodularity, the quenched law of large numbers
\eqref{eq:wp-local-quenched}, and its dependence on $\lambda$.
The observable records simultaneous based holonomies at the root,
rather than separate free homotopy classes. We unfold its factorial
measures over the carriers of short based words, compute their limits
with the one- and two-carrier coefficient estimates of
Appendix~\ref{app:positive-density-orbits}, assemble the limiting
random group, and then prove concentration and continuity in the cusp
density. The external inputs are Huang's normalized carriers and
rooted integration formula~\cite{HuangGap}, Mirzakhani's integration
formula~\cite{MirInv}, the boundary profile of
Proposition~\ref{prop:bounded-boundaries}, the invariant random
subgroup dictionary of Ab\'ert and Biringer~\cite{AbertBiringer}, and,
at the planar endpoint, the sphere local limit
of~\cite[Theorem~34]{BC}. The
six subsections treat, in order, the finite observable and its
factorial inversion, the support of a based tuple, uniform bounds for
the rooted fibers, the limiting factorial measures and the
construction of the limit, concentration, and the dependence on
$\lambda$. The fiber bound \eqref{eq:rooted-fiber-bound} and the waist
density \eqref{rev:eq-waist-density-convergence} are used again in
Appendix~\ref{app:conformal-phase}. All displacement cutoffs, carrier
types and boundary-length bounds are fixed before $g\to\infty$.
Constants can depend on these data and on a compact subinterval of
$0<\lambda<\infty$.

\subsection{A finite observable in every displacement window}

The observable is the configuration of primitive holonomies of bounded
displacement at the root. Lemma~\ref{lem:rooted-primitive-bound} shows
that it has a uniformly bounded number of atoms and determines the
pointed ball, so finitely many factorial measures control its law
through the inversion formula \eqref{eq:rooted-factorial-inversion}.

Let $G=\operatorname{PSL}_2(\mathbb R)$, fix a model oriented frame at
$o\in\mathbb H^2$, and normalize Haar measure so that rotations have
mass one. The frame bundle of $X=\Gamma\backslash\mathbb H^2$ then has
volume $A=A_{g,n}$. Set
\[
 K_R=\{a\in G:d(o,ao)\le R\},\qquad
 Q_R(\Gamma)=\sum_{\substack{a\in\Gamma\cap K_R, a\ne1\\
                            a\text{ primitive in }\Gamma}}\delta_a.
\]
Both orientations and primitive parabolic elements are included.
Primitivity means that $a$ is not a proper positive power in $\Gamma$.

\begin{lemma}[Uniform primitive bound]\label{lem:rooted-primitive-bound}
For every $R$ there is $M_R<\infty$ such that
$|Q_R(\Gamma)|\le M_R$ for every discrete torsion-free $\Gamma<G$.
The configuration $Q_{4r+1}$ determines the pointed measured ball
$B_r(\Gamma\backslash\mathbb H^2,o)$ with its restricted ambient metric.
\end{lemma}

\begin{proof}
Let $\mu>0$ be such that $2\mu$ is smaller than the two-dimensional
Margulis constant; then every point moved by less than $2\mu$ by some
nontrivial element lies in the thin part, whose components have cyclic
stabilizers. For a hyperbolic or parabolic isometry,
the displacement formula gives
\[
 \sinh\frac{d(y,ay)}2
       \le e^{d(x,y)}\sinh\frac{d(x,ax)}2.
\]
Choose $\delta_R>0$ with $e^R\sinh\delta_R<\sinh\mu$.
If the injectivity radius at the root is less than $\delta_R$, one
nontrivial element has displacement less than $2\mu$ at every point
of its radius-$R$ ball. That connected ball lies in a single thin
component. Every based loop of length at most $R$ stays in this ball,
so its group element lies in the cyclic group of that component.
Only the two opposite primitive generators can occur. Otherwise the
orbit of the root is $2\delta_R$-separated; disjoint hyperbolic disks
give a uniform packing bound in $B_{R+\delta_R}(o)$.

Taking powers of the recorded primitives recovers $\Gamma\cap K_R$,
since displacement is nondecreasing with the absolute power. To compute
distances between lifts in $B_r(\mathbb H^2,o)$, only translations with
displacement at most $4r$ are relevant: every distance in the quotient
is at most $2r$. The same translations determine identifications and
the quotient area. This proves the second assertion with $Q_{4r}$;
the statement uses $Q_{4r+1}$ only to leave a margin for the closed
ball.
\end{proof}

For a finite configuration $Q$ with at most $M_R$ atoms, let $\alpha_j$
be its $j$th ordered factorial measure. For any bounded Borel $F$,
finite inclusion--exclusion gives
\begin{equation}\label{eq:rooted-factorial-inversion}
 \mathbb EF(Q)=F(\varnothing)+\sum_{j=1}^{M_R}\frac1{j!}
      \int D^jF(\varnothing;\mathbf a)\,d\alpha_j(\mathbf a),
 \qquad |D^jF|\le2^j\|F\|_\infty.
\end{equation}
Here $D^jF$ is the alternating sum over all subsets of the $j$ atoms.
Thus total variation convergence of these finitely many measures
implies total variation convergence of $Q_R$ itself. A ball observable
has a bounded Borel extension to arbitrary configurations: generate
their subgroup, evaluate the ball if the subgroup is discrete and
torsion-free, and put zero otherwise. These conditions are measurable
through countably many word tests. Lemma~\ref{lem:rooted-primitive-bound}
makes the extension agree with the original observable on valid data.

\subsection{The support of a based tuple}

Apart from two elementary cyclic cases, a tuple of short primitive
holonomies has a unique bounded connected carrier
(Lemma~\ref{lem:rooted-word-support}). This allows the root frame to be
unfolded over carriers, which gives the exact formula
\eqref{eq:rooted-carrier-unfolding}. The carrier is built from the
geodesics of the generators and of their pairwise products, and its
complexity and perimeter are bounded by the divergence argument
of~\cite[Lemma~3.3]{HuangGap}.

Fix distinct ambient primitive elements $a_1,\ldots,a_j\in K_R$ and
put $B=\{a_i,a_i^{-1}:1\le i\le j\}$ and $H=\langle B\rangle$.
Take the geodesic images of all nonidentity hyperbolic words in
$B\cup B^2$. Fill their geodesic graphs by disks and inherited
peripheral punctured disks, tighten the remaining boundaries, and
retain coincident boundary germs as distinct slots. Simple annular
components are kept separately, as in the normalized carrier
construction of~\cite[Section~3.3]{HuangGap}.

\begin{lemma}[Support of short based words]\label{lem:rooted-word-support}
Either $H$ is generated by a primitive parabolic or a primitive simple
hyperbolic element, or the tuple has a unique connected normalized
nonannular carrier $Y$, carrying $H$ up to simultaneous conjugation.
Its complexity and total positive boundary length are bounded in
terms of $j,R$. In the two elementary cases $j\le2$.
\end{lemma}

\begin{proof}
Cut along the boundary multicurve of the filled graphs and annuli.
Filled components of disjoint components of the word graph have
disjoint interiors, by Steps~1--2 of the proof of
Lemma~\ref{rev:lem:disjoint-carriers}. Those steps use only that the
geodesics involved are closed and avoid the other graph, and that a
filled component lies in no disk, cusp neighborhood or annulus, since
it contains a nonsimple geodesic or two crossing geodesics; a simple
word geodesic can at most coincide with a boundary slot or seam of
another component. Thus every piece is a filled component, an annulus
or a complementary component.
In the resulting Bass--Serre tree every generator and every product
of two generators in $B$ is elliptic. Indeed hyperbolic geodesic
representatives are in the filled family and parabolics are peripheral.
Two disjoint fixed subtrees would give a hyperbolic product. The Helly
property for finitely many subtrees therefore gives a common fixed
vertex. A noncyclic $H$ cannot fix an edge, whose stabilizer is cyclic.

Suppose its vertex belongs to a complementary component $Z$.
Every generator is then peripheral in $Z$. Choose noncommuting
generators $u,v$. In the covering surface for $\langle u,v\rangle$,
they form a free basis of a rank-two free group. They remain peripheral
and represent distinct peripheral classes, since their homology classes
are independent. The core is a pair of pants, possibly with cusps:
a one-holed torus has a peripheral class with zero homology and so
cannot have these peripheral generators. In a pair of pants the third
peripheral homology class is, up to signs, $[u]+[v]$. Hence at most
one of $uv,uv^{-1}$ can be peripheral. A nonperipheral geodesic in
this cover cannot project to a peripheral geodesic of $Z$: such a
hyperbolic axis is a boundary support line of its convex lift, whereas
a parabolic has no hyperbolic axis. This contradicts the inclusion
of both products in the filled family.

Thus $H$ lies in one filled component. Every word in the family is
carried there or on its boundary; no other filled component occurs.
The cyclic nonsimple case uses the carrier of its one geodesic.
Any geodesic subsurface carrying the tuple carries all the words in
$B\cup B^2$ and hence this filled support, proving minimality and
uniqueness. If $H$ is cyclic and elementary, ambient primitivity
leaves only its two opposite generators.

Every word used has length at most $2R$. Write $\mathcal G$ for the
union of the geodesics of the hyperbolic words in $B\cup B^2$; its
total length is at most $\Lambda_0:=2R[2j+(2j)^2]$. The area and
perimeter bound of~\cite[Lemma~3.3, (3.22)]{HuangGap} is stated for a
single filling geodesic, but its proof uses only the total length of
the graph and the filling property: cutting the carrier along the
graph leaves disks, peripheral punctured disks and collars of single
boundary geodesics, the graph-side lengths of these regions sum to
twice the total length, and a divergence-one vector field of norm at
most one bounds the area of each region by its graph-side length,
while each boundary geodesic is shorter than the graph-side boundary
of its collar. The same divergence argument therefore applies to the
union graph $\mathcal G$ and gives
$\ell(\partial Y)\le2\Lambda_0$ and $A(Y)\le2\Lambda_0$.
Gauss--Bonnet now bounds the carrier complexity. No bound on the
number of intersections of arbitrarily short curves is required.
\end{proof}

The products are necessary: separate free classes of $a$ and $bab^{-1}$
do not retain their relative based position. We next unfold that position.
A convex lift $C$ of a normalized carrier has stabilizer $\pi_1Y$;
different lifts have disjoint interiors and can meet only on boundary
support lines. This subgroup is root-closed in $\Gamma$. If $a=b^m$
lies in it, a nonperipheral hyperbolic axis has an interior neighborhood
in both $C$ and $bC$; at a boundary axis $b$ preserves the same side;
for a parabolic a horoball in $C$ is preserved by its root $b$.
In every case interiors overlap, so $bC=C$. Intrinsic and ambient
primitivity therefore agree. A normalizer also preserves this convex
core and creates no extra coset multiplicity.

A tuple with full nonannular support cannot lie in the peripheral
intersection of two lifts. Its carrier coset is unique. Unfolding
the root frame over these cosets gives, for any bounded Borel tuple
test $f$ supported on $K_R^j$,
\begin{equation}\label{eq:rooted-carrier-unfolding}
 \int_{\mathrm{Fr}(X)}\sum_{\substack{\text{ordered primitive tuples}\\
                          \text{nonannular support}}}f\,d\mathrm{fr}
       =\sum_{Y\subset X}I_f(Y),\qquad
 I_f(Y)=\int_{\mathrm{Fr}(Y^{\rm f})}
                      \sum_{\text{full-support tuples}}f\,d\mathrm{fr}.
\end{equation}
Here $Y^{\rm f}$ is the complete surface obtained by adding funnels to
positive boundaries. Its cover need not embed in $X$. Formula
\eqref{eq:rooted-carrier-unfolding} includes self-paired boundary germs.
Cusp labels, boundary-slot permutations and twists have exactly the
orbifold conventions of~\cite[Proposition~3.4]{HuangGap}.

\subsection{Uniform bounds for the rooted fibers}

A full nonannular tuple confines its root to the core or a bounded
collar, which bounds its fiber by the carrier area and perimeter
\eqref{eq:rooted-fiber-bound}. The annular and parabolic fibers,
denoted $U_{j,L}$ and $U_j^{\rm par}$, are bounded directly. These
bounds replace the primitive count of
Lemma~\ref{lem:primitive-count} when the coefficient estimates of
Appendix~\ref{app:positive-density-orbits} are applied to rooted
fibers.

In a funnel at distance $t$ from its waist, an element not stabilizing
that funnel has displacement at least $2t$. A full nonannular tuple
contains such an element. Its root is consequently in the core or
within distance $R/2$ of its boundary. Thus
\begin{equation}\label{eq:rooted-fiber-bound}
 |I_f(Y)|\le M_R^j\|f\|_\infty
       \bigl(A(Y)+\ell(\partial Y)\sinh(R/2)\bigr).
\end{equation}
After integration on a fixed moduli fiber, the right side is a fixed
constant times its WP volume. This bound is for full nonannular support.
A short boundary loop alone can stay short arbitrarily far into a funnel.

For a simple waist of length $L$, its cyclic-cover displacement condition is
\[
 \cosh t\le\frac{\sinh(R/2)}{\sinh(L/2)}.
\]
The root region has area at most
$2L\sinh(R/2)/\sinh(L/2)\le4\sinh(R/2)$.
Denote the corresponding frame integral by $U_{j,L}(f)$, using ordered
distinct choices from the two opposite generators. It is zero for
$j>2$ and bounded uniformly down to $L=0$ for $j=1,2$.
For the primitive parabolic $p(z)=z+1$, the root region is
$0\le\operatorname{Re}z<1$, $\operatorname{Im}z\ge[2\sinh(R/2)]^{-1}$,
of area $2\sinh(R/2)$. Define $U_j^{\rm par}(f)$ in the same way.
The parabolic contribution on $X$ is exactly $nU_j^{\rm par}(f)$.
Both orientations have already been included; no extra factor two occurs.

\subsection{Factorial measures and construction of the limit}

We specify the limiting factorial measures
\eqref{eq:rooted-factorial-law}, prove convergence in total variation,
and assemble the limiting primitive process into a discrete
torsion-free group. The convergence has two parts: an exact expansion
of the frame integral over the three support types
(Proposition~\ref{rev:prop-exact-tuple-expansion}), and uniform
convergence of the leading densities together with a bound on all
other graphs (Proposition~\ref{rev:prop-density-convergence}).
Lemma~\ref{rev:lem-unimodular} then proves unimodularity of the limit.

Put $P=2\pi^2$, $x=X(Y_\lambda)$, $a_\lambda=2\pi(2+\lambda)$, and
$b(L)=\mathfrak b_{Y_\lambda}(L)$. Use the cap series $K_x$ of
\eqref{eq:positive-cap-series}. For a planar carrier with $a$ cusps
and $k$ positive labeled boundary slots, let $I_f^{a,k}(\mathbf L)$
be the WP integral of $I_f(Y)$ and set
$\widehat I_f^{a,k}=P^{-(a+k-3)}I_f^{a,k}$.
The limiting factorial measure on $K_R^j$ is specified by
\begin{align}\label{eq:rooted-factorial-law}
 \alpha_{R,j,\lambda}(f)
 &=\frac{\lambda}{a_\lambda}U_j^{\rm par}(f)
   +\frac{\lambda}{a_\lambda Px}
       \int_0^R Lb(L)K_x(L)U_{j,L}(f)\,dL\notag\\
 &\quad+\frac{\lambda}{a_\lambda}
   \sum_{\substack{a\ge0, k\ge1\\a+k\ge3}}
    \frac{x^{a-1}}{a!k!P^k}
    \int\widehat I_f^{a,k}(\mathbf L)
       \left[\sum_{s=1}^k b(L_s)\prod_{i\ne s}K_x(L_i)\right]
                                      \prod_i L_i\,dL_i.
\end{align}
Only finitely many carrier types and a bounded length box contribute,
by Lemma~\ref{lem:rooted-word-support}. In particular, the expression
is finite for every bounded Borel $f$.

The prelimit object is the expected area-normalized factorial measure.
For $X=\Gamma\backslash\mathbb H^2$ and a frame $\mathrm{fr}=\Gamma g$
write $\Gamma_{\rm fr}=g^{-1}\Gamma g$ for the deck group read in that
frame, so that $Q_R(\Gamma_{\rm fr})$ is the configuration at
$\mathrm{fr}$, and put
\begin{equation}\label{rev:eq-prelimit-factorial}
 \alpha^{(g)}_{R,j}(f)=\mathbb E_{g,n}\,\frac1A
 \int_{\mathrm{Fr}(X)}
 \sum_{\substack{(a_1,\ldots,a_j)\ \text{ordered distinct}\\
                 \text{primitive in }\Gamma_{\rm fr}\cap K_R}}
   f(a_1,\ldots,a_j)\,d\mathrm{fr},
\end{equation}
for bounded Borel $f$ on $K_R^j$. Convergence to
\eqref{eq:rooted-factorial-law} is proved in two steps: an exact
expansion of the frame integral, and uniform convergence of the
leading coefficient densities together with an absolute bound for all
other graphs. For a carrier type $(h_0,a,k)$, let
$I_f^{h_0,a,k}(\mathbf x)$ be the orbifold WP integral of $I_f(Y)$
over the fiber $\mathcal M_{h_0,a+k}(0^a,\mathbf x)$; thus
$I_f^{a,k}=I_f^{0,a,k}$.

\begin{proposition}[Exact ordered-tuple expansion]
\label{rev:prop-exact-tuple-expansion}
Fix $R$, $j$ and a bounded Borel $f$ on $K_R^j$. For every
$X\in\mathcal M_{g,n}$,
\begin{equation}\label{rev:eq-exact-tuple-expansion}
 \int_{\mathrm{Fr}(X)}
 \sum_{\substack{\text{ordered distinct}\\ \text{primitive tuples}}}
   f\,d\mathrm{fr}
 = nU_j^{\rm par}(f)
   +\sum_{\gamma}U_{j,\ell_\gamma}(f)
   +\sum_{Y\subset X}I_f(Y),
\end{equation}
where $\gamma$ runs over the unoriented simple closed geodesics of $X$
and $Y$ over the connected normalized nonannular carriers of
Lemma~\ref{lem:rooted-word-support} in $X$. The three terms are the
three support types of that lemma: parabolic cyclic, simple hyperbolic
cyclic, and full nonannular. Both orientations of a parabolic or of a
waist are the two generators entering $U_j^{\rm par}$ and $U_{j,L}$;
no additional factor two occurs. For $g$ larger than a constant
depending on $j,R$, the WP expectation of the last term is
\begin{equation}\label{rev:eq-tuple-carrier-formula}
 \mathbb E_{g,n}\sum_{Y\subset X}I_f(Y)
 =\sum_{(h_0,a,k)}\frac1{k!}\sum_{\mathfrak p\in\mathcal P_k}
  \int_{(0,\infty)^{e_{\mathfrak p}}}
  I_f^{h_0,a,k}(\mathbf x_{\mathfrak p}(\mathbf y))\,J_{\mathfrak p}
  \sum_{(g_i,a_i)}
  \frac{n!}{a!\prod_ia_i!}
  \frac{\prod_iV_{g_i,a_i+k_i}(0^{a_i},\mathbf x_{I_i})}{V_{g,n}}
  \,d\mathbf y .
\end{equation}
Here $(h_0,a,k)$ runs over the finitely many carrier types allowed by
Lemma~\ref{lem:rooted-word-support}; a pattern $\mathfrak p\in\mathcal P_k$
self-pairs $2\sigma$ of the $k$ labeled slots and partitions the
remaining slots into nonempty blocks $I_1,\ldots,I_b$ of sizes $k_i$;
there are $e_{\mathfrak p}=k-\sigma$ edges with lengths $\mathbf y$,
$\mathbf x_{\mathfrak p}(\mathbf y)$ assigns $y_e$ to both slots of a self-pair and
to each external slot, and $J_{\mathfrak p}=\prod_ey_e$; the inner sum runs over
$g_i,a_i\ge0$ with $\sum_ia_i=n-a$, $\sum_ig_i=g-h_0-k+\sigma+b$ and
$2g_i+a_i+k_i-2\ge1$. The factor $1/k!$ forgets the slot labels and
the factor $1/a!$ comes from the cusp labels. The waist term has the
WP expectation
\begin{equation}\label{rev:eq-waist-mirzakhani}
 \mathbb E_{g,n}\sum_{\gamma}U_{j,\ell_\gamma}(f)
 =\sum_{[\gamma]}\frac1{|\mathrm{Sym}(\gamma)|}
  \int_0^RU_{j,L}(f)\,L\,
  \frac{V(X\setminus\gamma;L)}{V_{g,n}}\,dL ,
\end{equation}
the sum over mapping class group orbits of simple closed geodesics,
with $V(X\setminus\gamma;L)$ the WP volume of the cut surface with
two new boundaries of length $L$, a product of two volumes for
separating $\gamma$, and $|\mathrm{Sym}(\gamma)|\in\{1,2\}$.
\end{proposition}

\begin{proof}
We first verify the three-way split \eqref{rev:eq-exact-tuple-expansion}
and then derive the two expectation formulas from the integration
formulas of Huang and of Mirzakhani.
Lemma~\ref{lem:rooted-word-support} assigns to every ordered tuple
exactly one of the three support types. If $H$ is generated by a
primitive parabolic $p$, every $a_i$ is one of $p^{\pm1}$, and
unfolding the frame to the parabolic cyclic cover of its cusp gives
$U_j^{\rm par}(f)$ for each of the $n$ cusps. If $H$ is generated by
a primitive simple hyperbolic $a$ with geodesic $\gamma$, every $a_i$
is one of $a^{\pm1}$, and unfolding to the annular cyclic cover of
$\gamma$ gives $U_{j,\ell_\gamma}(f)$. The remaining tuples are
exactly those with full nonannular support, and
\eqref{eq:rooted-carrier-unfolding} is
\eqref{rev:eq-exact-tuple-expansion} for them.

For \eqref{rev:eq-tuple-carrier-formula} we take from
Huang~\cite{HuangGap} two ingredients. First, the reconstruction
statement of \cite[Lemma~3.3]{HuangGap}: gluing a carrier to stable
complementary components, with some slots self-paired, reconstructs
the canonical normalized carrier, and cutting an ambient carrier
determines these records uniquely. Its proof uses only the filled
graph and the negative Euler characteristic of the stable exterior
pieces; it applies verbatim to the graph of the words in $B\cup B^2$,
and the root-closedness argument preceding
\eqref{eq:rooted-carrier-unfolding} supplies the equivalence of
intrinsic and ambient primitivity for tuples. Second, the
bookkeeping in the proof of \cite[Proposition~3.4]{HuangGap}:
in marked Fenchel--Nielsen coordinates each actual edge carries one
twist circle of full period $y_e$, twist integration multiplies an
observable independent of the twists by $J_{\mathfrak p}$ and leaves the orbifold
weights of the cut components, cusp labels are distributed in exactly
$n!/(a!\prod_ia_i!)$ ways, and the $k!$ bijective slot labellings of
a carrier form one symmetric-group quotient, so that no further
factor $2^{-\sigma}$, $1/\sigma!$ or $1/b!$ appears. The tuple
extension replaces Huang's single-geodesic fiber sum by $I_f(Y)$,
which is defined through $\mathrm{Fr}(Y^{\rm f})$ and $\pi_1Y$ alone
and is therefore independent of the twists and of the slot labelling;
it is Borel in the boundary lengths because based holonomies are
continuous in Fenchel--Nielsen coordinates and $f$ is Borel, so that
Tonelli applies to $|f|$ and then to the positive and negative parts
of $f$; and it is bounded on every fiber by
\eqref{eq:rooted-fiber-bound}. Huang's formula requires at least one
stable complementary component. Since carrier genera and $k$ are
bounded in terms of $j,R$ while $\sum_ig_i=g-h_0-k+\sigma+b$, some
exterior component has positive genus, hence is stable, once $g$
exceeds a constant depending on $j,R$. Equation
\eqref{rev:eq-waist-mirzakhani} is Mirzakhani's integration
formula~\cite[Theorem~7.1]{MirInv} applied to the bounded function
$L\mapsto U_{j,L}(f)$ of the waist length; $U_{j,L}(f)$ vanishes for
$L>R$, since a tuple at displacement at most $R$ has $\ell_\gamma\le R$.
\end{proof}

\begin{proposition}[Uniform density convergence and non-leading bound]
\label{rev:prop-density-convergence}
Let $g\to\infty$ with $n/g\to\lambda\in(0,\infty)$, and fix $R,j$.
\begin{enumerate}
\item[(i)] For each planar type $(0,a,k)$ with $a+k\ge3$, the leading
graph of \eqref{rev:eq-tuple-carrier-formula} is the tree pattern with
$\sigma=0$ and $b=k$ singleton blocks, one of which, at slot $s$,
carries the exterior component of genus $g$ while the others are planar
degree-one caps. After division by $A$ and summation over $s$ and over
the cap cusp counts, its contribution is
$\int\widehat I_f^{a,k}(\mathbf L)D^{(g)}_{a,k}(\mathbf L)\,d\mathbf L$
with a density satisfying, for every $R'<\infty$,
\begin{equation}\label{rev:eq-leading-density-convergence}
 \sup_{\mathbf L\in[0,R']^k}
 \Bigl|D^{(g)}_{a,k}(\mathbf L)
 -\frac{\lambda}{a_\lambda}\frac{x^{a-1}}{a!k!P^k}
  \Bigl[\sum_{s=1}^kb(L_s)\prod_{i\ne s}K_x(L_i)\Bigr]
  \prod_iL_i\Bigr|\longrightarrow0 .
\end{equation}
The leading graph of \eqref{rev:eq-waist-mirzakhani} is a separating
waist with a planar side and a genus-$g$ side. After division by $A$
and summation over the planar cusp counts its density $D^{(g)}_{\rm s}$
satisfies
\begin{gather}
 \sup_{d\le L\le R}\Bigl|D^{(g)}_{\rm s}(L)
   -\frac{\lambda}{a_\lambda Px}Lb(L)K_x(L)\Bigr|\longrightarrow0,
 \notag\\
 \sup_g\int_0^dD^{(g)}_{\rm s}(L)\,dL
 +\int_0^d\frac{\lambda}{a_\lambda Px}Lb(L)K_x(L)\,dL\le Cd^2,
 \label{rev:eq-waist-density-convergence}
\end{gather}
for every $0<d\le R$.
\item[(ii)] For $c\in\{1,2\}$ distinguished carriers of the types of
Lemma~\ref{lem:rooted-word-support}, the total weight of all gluing
graphs other than the leading ones, with the rooted fibers inserted,
is at most $Cg^{c-1}\|f\|_\infty$; for $c=1$ this covers every
pattern, every complementary allocation and every carrier type with
$h_0\ge1$ in \eqref{rev:eq-tuple-carrier-formula}, as well as the
nonseparating and two-positive-genus orbit types in
\eqref{rev:eq-waist-mirzakhani}.
\end{enumerate}
\end{proposition}

\begin{proof}
Both parts reduce to the coefficient estimates of
Appendix~\ref{app:positive-density-orbits}, applied with the rooted
fibers of \eqref{eq:rooted-fiber-bound} in place of the orbital fiber
integrals.

(i) For the leading planar pattern, $J_{\mathfrak p}=\prod_iL_i$, the cusp-label
factor is $n!/(a!\prod_{i\ne s}a_i!\,m!)$ with
$m=n-a-\sum_{i\ne s}a_i$, and the volumes are
$\prod_{i\ne s}V_{0,a_i+1}(0^{a_i},L_i)\cdot V_{g,m+1}(0^m,L_s)$.
Proposition~\ref{prop:bounded-boundaries}, through
\eqref{eq:boundary-uniform}, replaces $V_{g,m+1}(0^m,L_s)$ by
$V_{g,m+1}\mathfrak b_{Y_{(m+1)/g}}(L_s)$ with relative error tending to
zero uniformly in $m$ and in $L_s\in[0,R']$. On the truncated range of
cap cusp counts one has $(m+1)/g\to\lambda$, and
$\mathfrak b_{Y_{(m+1)/g}}\to b$ uniformly on $[0,R']$ by continuity of
\eqref{eq:boundary-profile} in $Y$; on the complementary range only the
uniform bound $\mathfrak b_Y(L)\le C_{R'}$ is used. The remaining cusp-label and volume ratio, summed over the
cap cusp counts $a_i\ge2$, is the one-carrier coefficient ratio of
Appendix~\ref{app:positive-density-orbits}: truncating the total cap
cusp degree, every fixed coefficient shift has the same local-limit
mass as the parent coefficient, which gives $nx_g^{a-1}/a!$ times
$\prod_{i\ne s}x_g^{a_i}/a_i!$ with the powers of $P$ recorded in
$\widehat I_f^{a,k}$ and $P^{-k}$, while the point-mass bound
\eqref{eq:positive-coefficient-bound} dominates the truncated tail by
the tail of the convergent positive series $\prod_{i\ne s}K_{x_g}(L_i)$,
uniformly on compact density and length intervals. Dividing by
$A\sim ga_\lambda$, using $n/g\to\lambda$, $x_g\to x$ and the
continuity of $K_x$ in $x$ and $L$ gives
\eqref{rev:eq-leading-density-convergence}. The waist density is the
case of one edge: its exact prefactor is $n/P$ times the coefficient
ratio, which extraction turns into $K_x(L)/x$, and the boundary factor
is $Lb(L)$. The one-edge coefficient bound of
Appendix~\ref{app:positive-density-orbits} is uniform down to zero
waist length, and $Lb(L)K_x(L)\le CL$ on $[0,R]$; integrating
$L\,dL$ over $[0,d]$ gives \eqref{rev:eq-waist-density-convergence}.

(ii) By \eqref{eq:positive-graph-power}, a graph with $c$ distinguished
carriers, total carrier genus $h_0$, $k$ slots, $\sigma$ self-pairs,
$u$ planar degree-one, $v$ planar degree-two and $w$ stable exterior
vertices has total weight
$O(g^{2c+1-2h_0-k+u-w})$, and the exponent analysis following that
display shows that the exponent is at most $c-1$ unless $h_0=0$,
$v=0$, $w=1$, $\sigma=0$ and $k=u+c$, which is the leading tree. The
same analysis gives $O(1)$ for nonseparating waists and $O(g^{-1})$
for separations into two positive-genus pieces. Bounded boundary
lengths and fixed twist integrals preserve these bounds. The rooted
fiber enters through \eqref{eq:rooted-fiber-bound}: on the bounded
length box of Lemma~\ref{lem:rooted-word-support},
$|I_f^{h_0,a,k}(\mathbf x)|\le C_{R,j}\|f\|_\infty
V_{h_0,a+k}(0^a,\mathbf x)$, and the WP volume is bounded on the box.
The number of carrier types and patterns is finite. Summing proves
the bound $Cg^{c-1}\|f\|_\infty$.
\end{proof}

We now conclude. Write $D_{a,k}$ and $D_{\rm s}$ for the limiting
densities in \eqref{rev:eq-leading-density-convergence} and
\eqref{rev:eq-waist-density-convergence}. Then
\eqref{eq:rooted-factorial-law} reads
$\alpha_{R,j,\lambda}(f)=\frac{\lambda}{a_\lambda}U_j^{\rm par}(f)
+\int_0^RD_{\rm s}(L)U_{j,L}(f)\,dL
+\sum_{a,k}\int\widehat I_f^{a,k}D_{a,k}$.
Subtracting \eqref{rev:eq-exact-tuple-expansion}, divided by $A$ and
averaged over WP, and using $|U_j^{\rm par}(f)|,|U_{j,L}(f)|\le C_R
\|f\|_\infty$ and $|\widehat I_f^{a,k}(\mathbf L)|\le C_{R,j}\|f\|_\infty
V_{0,a+k}(0^a,\mathbf L)$ on the bounded box, together with
Proposition~\ref{rev:prop-density-convergence}, gives for every
$0<d\le R$
\begin{equation}\label{rev:eq-factorial-tv}
 \limsup_{g\to\infty}\
 \sup_{\|f\|_\infty\le1}
 \bigl|\alpha^{(g)}_{R,j}(f)-\alpha_{R,j,\lambda}(f)\bigr|
 \le Cd^2 ,
\end{equation}
where the parabolic term contributes $|n/A-\lambda/a_\lambda|C_R\to0$,
the waists on $[d,R]$ and the leading planar graphs contribute the
uniform density errors times bounded fibers and finite box volumes,
the waists on $[0,d]$ contribute $Cd^2$, and the non-leading graphs
contribute $Cg^{-1}$. Letting $d\downarrow0$ proves convergence of
the factorial measures in total variation, rather than only on smooth
orbital tests.

Equation~\eqref{eq:rooted-factorial-inversion} makes the laws of $Q_R$
Cauchy in total variation. Their limits for integer $R$ are consistent.
Taking powers in the resulting primitive process gives a group:
inverse, product and primitive-root relations are countable measurable
predicates in sufficiently large displacement windows and pass under
total variation convergence. Its nonidentity atoms are hyperbolic or
parabolic; no elliptic or identity atom can appear, because the event
that a window contains an elliptic or identity atom has probability
zero for every finite-area prelimit, $\Gamma$ being discrete and
torsion-free and the identity being excluded from $Q_R$ by definition,
and a null event remains null under total variation limits. Every
bounded window has finitely many primitives, each with finitely many
powers there. The group is therefore discrete and torsion-free. Its
quotient is a complete hyperbolic surface.
Lemma~\ref{lem:rooted-primitive-bound} proves \eqref{eq:wp-local-tv}
and uniqueness.

\begin{lemma}[Unimodularity of the limit]\label{rev:lem-unimodular}
The law of the limiting group $\Gamma_\infty$ is invariant under
conjugation by $G$, and the law $\mathcal L_\lambda$ of its pointed
quotient satisfies \eqref{eq:local-mass-transport}.
\end{lemma}

\begin{proof}
We show that conjugation acts measurably on window configurations,
pass the conjugation invariance of the prelimit frame law to the
limit, and then apply the dictionary between invariant random
subgroups and unimodular random manifolds.
For $h,a\in G$, the triangle inequality gives
\[
 d(o,h^{-1}ah\,o)=d(ho,a\,ho)\le d(o,ao)+2d(o,ho),
\]
hence $h^{-1}K_Rh\subset K_{R+2d(o,ho)}$ and, for every discrete
torsion-free $\Gamma$,
\[
 Q_R(h^{-1}\Gamma h)
 =h^{-1}\bigl\{a\in Q_{R+2d(o,ho)}(\Gamma):d(o,h^{-1}ah\,o)\le R\bigr\}h .
\]
Thus $Q_R(h^{-1}\Gamma h)$ is a Borel function of
$Q_{R+2d(o,ho)}(\Gamma)$. In the finite-area prelimit the frame
$\mathrm{fr}=\Gamma g$ is Haar-distributed on $\Gamma\backslash G$, and
right multiplication by $h$ preserves this law while replacing
$\Gamma_{\rm fr}=g^{-1}\Gamma g$ by $h^{-1}\Gamma_{\rm fr}h$; so the law
of $\Gamma_{\rm fr}$ is conjugation invariant, and for every bounded
Borel $F$ of a configuration,
$\mathbb EF(Q_R(h^{-1}\Gamma_{\rm fr}h))=\mathbb EF(Q_R(\Gamma_{\rm fr}))$.
Both sides are expectations of bounded Borel functions of
$Q_{R+2d(o,ho)}(\Gamma_{\rm fr})$, whose law converges in total
variation. Passing to the limit gives
$\mathbb EF(Q_R(h^{-1}\Gamma_\infty h))=\mathbb EF(Q_R(\Gamma_\infty))$
for every $R$. The window configurations $(Q_R)_{R\in\mathbb N}$
determine a discrete group, so the law of $\Gamma_\infty$ on the
Chabauty space of closed subgroups is conjugation invariant: it is a
discrete torsion-free invariant random subgroup of
$G=\operatorname{PSL}_2(\mathbb R)$, which is unimodular and acts
transitively on $\mathbb H^2$ with compact stabilizers, and a discrete
torsion-free subgroup acts freely and properly discontinuously.
By Ab\'ert and Biringer~\cite[Proposition~2.9]{AbertBiringer}, the
pointed quotient $(\Gamma_\infty\backslash\mathbb H^2,[o])$ of such an
invariant random subgroup is a unimodular random manifold, that is,
its law satisfies the mass-transport identity
\eqref{eq:local-mass-transport}; the dictionary between discrete
torsion-free invariant random subgroups of the isometry group and
unimodular random $\mathbb H^2$-manifolds is
\cite[Proposition~1.11]{AbertBiringer} for the full isometry group;
only the forward direction, which \cite[Proposition~2.9]{AbertBiringer}
gives for $\operatorname{PSL}_2(\mathbb R)$, is used here.
\end{proof}

\subsection{Concentration of local observables}

As for orbital statistics, second moments factorize to leading order,
and short annuli are bounded separately. Pairs of tuples whose
carriers have disjoint interiors are counted by the two-carrier graphs
of Appendix~\ref{app:positive-density-orbits}; pairs whose carriers
overlap are treated in Lemma~\ref{rev:lem-overlapping-carriers}. The
outcome is the quenched law of large numbers
\eqref{eq:wp-local-quenched}.

For fixed $f$, split its area-normalized empirical factorial statistic
into parabolic, annular and nonannular parts. The first is deterministic.
For the nonannular second moment, disjoint carrier interiors use the
two-carrier formula of Appendix~\ref{app:positive-density-orbits}.
For clarity, its exact two-carrier form is obtained from
\eqref{rev:eq-tuple-carrier-formula} as follows. For fixed ordered
carrier types $(h_1,a_1,k_1)$ and $(h_2,a_2,k_2)$, the contribution
of pairs with disjoint interiors is
\[
 \frac{1}{k_1!k_2!V_{g,n}}
 \sum_{\mathfrak p}\int
 I_f^{h_1,a_1,k_1}(\mathbf x_1)
 I_{\tilde f}^{h_2,a_2,k_2}(\mathbf x_2)\,J_{\mathfrak p}
 \sum_{(g_v,a_v)}
 \frac{n!}{a_1!a_2!\prod_v a_v!}
 \prod_v V_{g_v,a_v+k_v}(0^{a_v},\mathbf x_{I_v})\,d\mathbf y.
\]
Here $\tilde f=f$ for a real second moment. The patterns pair slots
within or between the two carriers and partition all remaining slots
into $b$ nonempty complementary blocks, retaining only connected graphs.
If $k=k_1+k_2$ and $\sigma$ pairs are made, then
$J_{\mathfrak p}=\prod_{e=1}^{k-\sigma}y_e$,
$\sum_v a_v=n-a_1-a_2$, and
$\sum_v g_v=g-h_1-h_2-k+\sigma+b+1$; as before every complementary
component is stable. This formula follows by cutting along the union
of the slot geodesics and applying Fenchel--Nielsen disintegration.
There are two independent frame integrations, one in each carrier
cover, which give the two fiber factors. Integration over each twist
interval of length $y_e$ gives $y_e$ once, even for a seam shared by
the carriers. The cusp multinomial chooses disjoint label sets for
both carriers and all complements. Forgetting their slot labels gives
$1/(k_1!k_2!)$; there is no factor $1/2$, since the carriers are ordered.
Thus this is a two-carrier extension of the one-carrier unfolding,
not an additional statement imported from Huang.

Divide the display by $A^2$ for the empirical second moment.
The fiber bound \eqref{eq:rooted-fiber-bound} dominates the integrand
by a constant times the corresponding product of WP volumes, so
\eqref{eq:positive-graph-power} applies with $c=2$. Its leading graph
has one genus-$g$ exterior vertex of degree two; all other disjoint
patterns have weight $O(g)$ and vanish after division by $A^2$.
The derivative relation \eqref{eq:positive-root-derivative} and bounded
boundary factorization give the product of the two first moments.

Pairs whose carrier interiors overlap are treated by the following
lemma. Coincident carriers are included there; opposite boundary germs
with disjoint interiors are seams in the two-carrier formula.

\begin{lemma}[Overlapping carriers]\label{rev:lem-overlapping-carriers}
Let $\mathbf a,\mathbf a'$ be full nonannular tuples in $K_R^j$ at
roots $p,p'\in X$, with word graphs $\mathcal G,\mathcal G'$ and
normalized carriers $Y,Y'$ whose interiors intersect. Put
$\Lambda_0=2R[2j+(2j)^2]$. Then $Y\cup Y'$ is contained in a connected
essential subsurface $Y''\subset X$ with geodesic boundary and
inherited cusps such that
\[
 A(Y'')\le4\Lambda_0,\qquad \ell(\partial Y'')\le4\Lambda_0,
\]
so $Y''$ has one of finitely many types, with bounded boundary lengths.
Both roots lie in $Y''$ or within distance $R/2$ of $\partial Y''$ in
the complete cover $(Y'')^{\rm f}$, and for fixed roots there are at
most $M_R^{2j}$ ordered pairs of tuples. The expected total contribution
of all overlapping pairs to the second moment, before division by
$A^2$, is $O(g)\|f\|_\infty^2$.
\end{lemma}

\begin{proof}
We construct $Y''$ separately according to whether the filling parts
of the two word graphs meet, and then bound the two-root fiber of
$Y''$. Let $\mathcal G_0\subset\mathcal G$ and
$\mathcal G_0'\subset\mathcal G'$ be the connected components that fill
$Y$ and $Y'$; by the proof of Lemma~\ref{lem:rooted-word-support}, the
other components are simple boundary or seam geodesics of $Y$,
respectively $Y'$.
If $\mathcal G_0\cap\mathcal G_0'\ne\varnothing$, the union graph
$\mathcal G_0\cup\mathcal G_0'$ is connected with total length at most
$2\Lambda_0$. Let $Y''$ be its normalized carrier: fill the disk and
peripheral punctured-disk complementary regions and tighten the
remaining boundaries. The complement of a small neighborhood of the
union graph in $Y''$ consists of disks, cusp neighborhoods and collars
of single boundary geodesics, so the divergence argument of
\cite[Lemma~3.3]{HuangGap}, which uses only the total length of the
graph and this filling property, gives
$A(Y'')\le2\ell(\mathcal G_0\cup\mathcal G_0')\le4\Lambda_0$ and
$\ell(\partial Y'')\le4\Lambda_0$, exactly as in the proof of
Lemma~\ref{lem:rooted-word-support}. Any geodesic subsurface containing
the filling graph of a tuple contains its carrier, by the minimality
proved there; hence $Y\cup Y'\subset Y''$.

If $\mathcal G_0\cap\mathcal G_0'=\varnothing$, choose disjoint thin
regular neighborhoods $N(\mathcal G_0)$ and $N(\mathcal G_0')$. Every
boundary geodesic of $Y$ is the tightening of a boundary curve of the
pre-carrier
$N(\mathcal G_0)\cup\{\text{absorbed disks and punctured disks}\}$,
which lies in $\partial N(\mathcal G_0)$, and likewise for $Y'$. Boundary
curves of $N(\mathcal G_0)$ and of $N(\mathcal G_0')$ are disjoint, so
their geometric intersection numbers vanish, and simple closed
geodesics with zero intersection number are disjoint or equal. Thus the
components of $\partial Y$ and $\partial Y'$ are pairwise disjoint or
coincide, and $Y''=Y\cup Y'$ is a connected subsurface with geodesic
boundary and inherited cusps. The components of $Y\cap Y'$ are
geodesic circles or subsurfaces with geodesic boundary, which have
nonpositive Euler characteristic, so
$|\chi(Y'')|\le|\chi(Y)|+|\chi(Y')|$, and
$A(Y'')\le A(Y)+A(Y')\le4\Lambda_0$,
$\ell(\partial Y'')\le\ell(\partial Y)+\ell(\partial Y')\le4\Lambda_0$
by the bounds of Lemma~\ref{lem:rooted-word-support}. Gauss--Bonnet
bounds the complexity in both cases.

Unfold both roots to $\mathrm{Fr}((Y'')^{\rm f})$. In a funnel of
$(Y'')^{\rm f}$ at distance $t$ from its waist, an element not
stabilizing that waist has displacement at least $2t$. A full
nonannular tuple contains such an element: otherwise all its elements
would lie in the cyclic stabilizer of one simple boundary geodesic,
which is an elementary case excluded by
Lemma~\ref{lem:rooted-word-support}. Hence each root is confined to
the core or to the collar of width $R/2$, a region of area at most
$A(Y'')+\ell(\partial Y'')\sinh(R/2)$. At fixed roots, each tuple is
an ordered $j$-tuple of distinct atoms of a window configuration with
at most $M_R$ atoms, by Lemma~\ref{lem:rooted-primitive-bound} applied
to $\pi_1(Y'')$, so there are at most $M_R^{2j}$ ordered pairs. The
two-root fiber function of $Y''$ is therefore bounded by
$M_R^{2j}\|f\|_\infty^2\bigl(A(Y'')+\ell(\partial Y'')\sinh(R/2)\bigr)^2$,
a constant times $\|f\|_\infty^2$ on the bounded length box, and it is
Borel on the fiber for the same reason as $I_f$. Every overlapping
pair is counted in the fiber of its $Y''$. Summing over the finitely
many types of $Y''$ is the one-carrier case $c=1$ of
\eqref{eq:positive-graph-power}, whose total weight, with the bounded
fiber inserted as in Proposition~\ref{rev:prop-density-convergence}(ii),
is $O(g)\|f\|_\infty^2$.
\end{proof}

After division by $A^2$ the overlapping pairs therefore contribute
$O(g^{-1})$ and vanish in the limit.

For annuli with $d\le L\le R$, the two-simple-edge computation of
Appendix~\ref{app:positive-density-orbits}, with the bounded $U$ fibers,
also factorizes to leading order. The area-normalized absolute
contribution from $L<d$ is deterministically bounded by
$M_R^j\|f\|_\infty$ and has expectation $O(d^2)$. Its squared $L^2$
norm is therefore $O(d^2)$. Letting $d\downarrow0$ proves $L^2$
concentration of the full annular statistic. The triangle inequality
for the centered $L^2$ norms combines the three parts without a mixed
covariance calculation. Applying the finite inversion
\eqref{eq:rooted-factorial-inversion} to a bounded ball observable
proves \eqref{eq:wp-local-quenched}.

\subsection{Dependence on the cusp density}

The explicit factorial law \eqref{eq:rooted-factorial-law} gives
continuity in $\lambda$, the plane as the limit $\lambda\downarrow0$,
and the area-rooted infinite WP sphere as the limit
$\lambda\uparrow\infty$, which proves \eqref{eq:wp-local-endpoints};
the planar endpoint uses the genus-zero convolution identities of
Appendix~\ref{app:cusp-dominated-orbits} and the sphere local limit
of~\cite{BC}. The cusp occupation formula
\eqref{eq:local-cusp-occupation} is derived at the end.

On compact density intervals the cap series and boundary profiles in
\eqref{eq:rooted-factorial-law} have common integrable envelopes and
are continuous. The factorial measures, and hence ball laws, are
continuous in total variation. As $\lambda\downarrow0$, $x\downarrow0$
and $K_x(L)=O_R(x^2)$. The nonannular prefactor has order
\[
 \lambda x^{a-1}K_x^{k-1}
       =O_R(\lambda x^{a+2k-3}),\qquad a+2k-3\ge1.
\]
The annular and parabolic parts are respectively $O_R(\lambda x)$
and $O_R(\lambda)$. Thus $\mathbb E|Q_R|\to0$, giving the plane endpoint.

The planar endpoint rests on the following genus-zero counterpart of
Propositions~\ref{rev:prop-exact-tuple-expansion}
and~\ref{rev:prop-density-convergence}. Let $\alpha_{R,j,\infty}$ be
defined by \eqref{eq:rooted-factorial-law} with $\lambda/a_\lambda$
replaced by $1/(2\pi)$, $x$ by $x_c$, $K_x$ by $K_{x_c}$ and $b$ by
$\mathfrak b_{y_0}$; and for $X$ with normalized WP law on
$\mathcal M_{0,n+1}$ let $\alpha^{(0,n)}_{R,j}$ be defined by
\eqref{rev:eq-prelimit-factorial} with $A=2\pi(n-1)$.

\begin{proposition}[Genus-zero factorial measures]
\label{prop:sphere-factorial}
For fixed $R$ and $j$,
\[
 \sup_{\|f\|_\infty\le1}
   |\alpha_{R,j,\lambda}(f)-\alpha_{R,j,\infty}(f)|\longrightarrow0
 \quad(\lambda\to\infty),\qquad
 \sup_{\|f\|_\infty\le1}
   |\alpha^{(0,n)}_{R,j}(f)-\alpha_{R,j,\infty}(f)|\longrightarrow0
 \quad(n\to\infty).
\]
\end{proposition}

\begin{proof}
For the first limit, $x\uparrow x_c$ and
$\lambda/a_\lambda\to1/(2\pi)$ as $\lambda\to\infty$. The cap series
$K_x(L)$ has nonnegative coefficients and converges at $x_c$ on
bounded length intervals, so $K_x\uparrow K_{x_c}$ uniformly on
$[0,R']$ by Dini's theorem, and $\mathfrak b_{Y_\lambda}\to
\mathfrak b_{y_0}$ uniformly there. Only finitely many carrier types
and a bounded length box enter \eqref{eq:rooted-factorial-law}, and
$|\widehat I_f^{a,k}|\le C\|f\|_\infty$ on that box; this gives the
first limit.

For the second, the expansion \eqref{rev:eq-exact-tuple-expansion}
holds on every finite-area surface, in particular in genus zero. There
the cut graph of a carrier is a tree, every exterior piece is a planar
cap with at least two cusps, and the integration formula
\eqref{rev:eq-waist-mirzakhani} and its nonannular analogue apply as
soon as $n$ exceeds the number of cusps of the carrier. In contrast to
positive genus, no graph is of lower order: for a nonannular carrier
of type $(0,a,k)$ the coefficient is the coefficient of $z^{n+1-a}$ in
a product of $k$ series $K_1(z;L_i)$, and the first identity in
\eqref{eq:cusp-sphere-convolutions} selects the one cap that carries
the diverging cusp count, the others contributing $K_{x_c}(L_i)$. For
a simple closed geodesic both sides are planar; the identity is used
with two factors, and the symmetry factor $1/2$ of the unordered
partition cancels the two choices of the large side. Dividing by
$V_{0,n+1}$, which is the coefficient of $\Phi_0$ at $z^{n+1}$ up to
the same normalization, and by the area $A=2\pi(n-1)$, gives the
prefactors $x_c^{a-1}/(2\pi\,a!k!P^k)$ and $1/(2\pi Px_c)$ of
$\alpha_{R,j,\infty}$; the parabolic term has coefficient
$(n+1)/A\to1/(2\pi)$. The asymptotic relations
\eqref{eq:cusp-sphere-convolutions} are dominated on bounded boundary
boxes, so the coefficient densities converge in $L^1$ of the length
box, and $|\widehat I_f^{a,k}|\le C\|f\|_\infty$ again gives
convergence uniformly in $\|f\|_\infty\le1$. No second-moment
estimate is needed for this endpoint.
\end{proof}

By Proposition~\ref{prop:sphere-factorial}, the inversion
\eqref{eq:rooted-factorial-inversion} and
Lemma~\ref{lem:rooted-primitive-bound}, the ball $B_r(X,o)$ of an
area-rooted WP sphere converges in total variation to the ball of the
surface with factorial measures $\alpha_{4r+1,j,\infty}$, which is
also the limit of $\mathcal L_\lambda$ as $\lambda\uparrow\infty$. By
the sphere local limit of~\cite[Theorem~34]{BC} it converges in total
variation to $B_r(S_\infty^*,o_\infty)$. The two limits therefore
coincide, which identifies the endpoint with the area-rooted infinite
WP sphere $(S_\infty^*,o_\infty)$ and fixes the root normalization.
Finally let $0<b\le2$. On every complete hyperbolic surface the cusp
tails bounded by horocycles of length $b$ are embedded and pairwise
disjoint, each has area $b$, and the root lies in one of them if and
only if some primitive parabolic $p$ of the deck group read at the root
satisfies $2\sinh(d(o,po)/2)<b$. On $X$ this event has probability
exactly $bn/A_{g,n}$, and total variation convergence of the
configurations $Q_R$ with $R=2\operatorname{arcsinh}(b/2)$ gives
\eqref{eq:local-cusp-occupation}.

\section{Finite caps and the conformal transition}
\label{app:conformal-phase}

This appendix proves the finite-density part of
Theorem~\ref{thm:wp-conformal-phase}
(Section~\ref{sec:local-conformal-phase}): for $0<\lambda<\infty$ the
local limit $S_\lambda$ has genus zero and infinite area, every simple
closed geodesic bounds a unique finite planar cap, the cusp-filled
surface has conformal type $\mathbb D$, and Brownian motion for the
hyperbolic metric is transient. The recurrent endpoint
$\lambda=\infty$ is proved in Appendix~\ref{app:cusp-cutoffs}.
The inputs are Theorem~\ref{thm:wp-local-family}, through the total
variation convergence \eqref{eq:wp-local-tv} and the mass-transport
identity \eqref{eq:local-mass-transport}, the graph and fiber bounds
\eqref{eq:positive-graph-power} and \eqref{eq:rooted-fiber-bound}, the
waist density \eqref{rev:eq-waist-density-convergence}, and the
Green-level argument of Delmotte and Rau~\cite{DelmotteRau}.
We first show that the limiting surface is planar and has infinite
area, and that every waist has a unique finite planar side. A
multiplicity estimate for caps with small boundary length then gives an
isoperimetric inequality for large domains containing the root, and a
Green-level argument converts this inequality into transience; the
four subsections follow this order.
Fix $0<\lambda<\infty$ and retain $P,x,a_\lambda,b$ from
Appendix~\ref{app:rooted-local-law}. Set $\xi=\log(x_c/x)>0$.
All areas, perimeters, Green functions and Brownian motions in the
argument use the original hyperbolic metric.

\subsection{Topology and recognition of complete finite caps}

After filling its cusps, the limiting surface is simply connected and
noncompact. Each topological statement is first expressed as a Borel
function of a fixed pointed ball, so that \eqref{eq:wp-local-tv}
transfers its prelimit expectation, and the prelimit expectation is
computed from the graph bounds of
Appendix~\ref{app:positive-density-orbits}. Finite caps are
recognized through screened ball observables, defined in
\eqref{rev:eq-screened-cap-count}, which are reused in the next
subsection.

A compact positive-genus witness in the limiting surface has a finite
short-word carrier after imposing finite length, distance and thickness
screens. Its carrier has positive genus. The one-carrier graph bound
\eqref{eq:positive-graph-power} and the rooted envelope
\eqref{eq:rooted-fiber-bound} show that its expected count, normalized
by area, tends to zero. The local total variation convergence transfers
each such witness. Exhausting the screens excludes every compact
handle, so the limiting surface has genus zero almost surely.

It also has infinite area. A primitive tuple generating a lattice of
covolume at most $M$ cannot occur in a finite approximating group
$\Gamma$ once $A>M$: a cofinite subgroup $H<\Gamma$ has covolume
$[\Gamma:H]A\ge A$. This is a Borel tuple test in a sufficiently large
displacement window, and its absence passes under factorial-measure
total variation convergence. A finite-area limiting group is finitely
generated; replacing finitely many generators by primitive roots gives
a generating tuple in some window. Countable exhaustion over $M,R$
therefore excludes finite area in the limit.

Fix a simple waist with length in $[d,L_0]$, $d>0$, and bounded distance
from the root. The leading simple-edge density has a finite planar
side and one genus-$g$ exterior component. The nonseparating and
two-positive-genus terms have vanishing area-normalized expectation.
The planar cap-size distribution has the positive series $K_x$ and
an exponentially decaying tail because $x<x_c$.

To transfer the entire finite cap we use screened ball observables.
A \emph{screen} is a tuple $\varsigma=(k,d,L_0,\sigma,t)$ with
$k\ge2$, $0<d<L_0<\infty$, a compact subset $\sigma$ of the moduli
space of planar caps with $k$ cusps and one geodesic boundary, and a
truncation depth $t\ge0$. A finite planar cap $C\subset S$ with $k$
cusps is \emph{$\varsigma$-screened} if $\ell(\partial C)\in[d,L_0]$
and its moduli point lies in $\sigma$; write $C^t$ for $C$ with the
open cusp regions beyond the standard horocycles at depth $t$ removed.
By compactness of $\sigma$ the truncated caps $C^t$ have diameter at
most $D=D(\varsigma)$. For $z\in C^t$, the ball $B_{R_0}(S,z)$ with
$R_0=D+L_0+2$ contains $C^t$, unit buffers at its truncation
horocycles, and based loops realizing the waist and the $k$ primitive
parabolic holonomies; these data recognize the complete cap, including
its omitted standard cusp tails, since each cusp is determined by its
parabolic holonomy. Define
\begin{equation}\label{rev:eq-screened-cap-count}
 N^{\varsigma}(S,z)
 =\#\{C\ \varsigma\text{-screened}:z\in C^t\}.
\end{equation}
This is a Borel function of the pointed ball $B_{R_0}(S,z)$, and it is
bounded by $M_{2D+L_0}$: distinct caps have distinct waists, and each
waist is freely homotopic to a based loop at $z$ of length at most
$2D+L_0$, hence to a primitive atom of $Q_{2D+L_0}$ at $z$. By
\eqref{eq:wp-local-tv}, $\mathbb EN^{\varsigma}(S_\lambda,o_\lambda)
=\lim_g\mathbb EN^{\varsigma}(X,o)$ for every screen. The indicator
that the root lies in some $\varsigma$-screened cap transfers in the
same way. For fixed $k$ and $[d,L_0]$ the prelimit intensity of caps
on their moduli space is the one-edge coefficient of
Appendix~\ref{app:positive-density-orbits} times the normalized WP
measure of the fiber, which is finite and hence tight; so the expected
number per unit area of caps with moduli outside $\sigma$, or with the
root beyond depth $t$ in a cusp, is small uniformly in $g$ once
$\sigma$ and $t$ are large. Letting $\sigma$ exhaust the moduli space,
$t\to\infty$, $d\downarrow0$ and $L_0\uparrow\infty$ then exhausts,
for each $k$, all finite planar caps with $k$ cusps containing the
root; the cusp count is finally removed by the positive cap-series
tail. For a fixed root-distance bound $R_1>0$, apply the same screening
to caps with $d(z,\partial C)\le R_1$. Their count is a bounded Borel
function of $B_{R_1+D+L_0+2}(S,z)$, bounded by $M_{2R_1+L_0}$; the root
fiber has area at most $2\ell(\partial C)\sinh R_1$. Thus the same
cap-size tail and finite WP-fiber tightness remove the screens
uniformly for waists within distance $R_1$ of the root. Exhausting the
positive length and root-distance boxes excludes
every waist without a finite planar side. Infinite total area makes
that side unique.

Filling the isolated cusps now gives a simply connected surface.
A simple loop, perturbed off the added points, is trivial, peripheral,
or represented by a simple closed geodesic. In the peripheral case it
bounds a disk after filling its cusp; in the geodesic case it bounds a
disk after filling the cusps of its finite cap. General loops reduce
to simple ones. The buffered cusp filling used in
Appendix~\ref{app:cusp-cutoffs} gives a complete metric on the filled
surface with locally discrete added centers. Infinitely many centers
cannot accumulate in a compact region; if there were finitely many
and the completion were compact, the original hyperbolic surface
would have finite area. The filled surface is thus noncompact, and
uniformization leaves $\mathbb C$ or $\mathbb D$.

\subsection{A cap multiplicity estimate}

Caps whose boundary length is small relative to their area have finite
expected multiplicity at the root, by the subcritical tail of the cap
intensity \eqref{eq:local-cap-intensity}. The mass-transport identity
\eqref{eq:local-mass-transport} then shows that every positive-area
core lies in only finitely many of them, which is the anchor
\eqref{eq:anchored-bad-cap-bound} used in the next subsection.

The per-area intensity of an unoriented waist of length $L$ with its
unique finite planar side containing $k\ge2$ cusps is
\begin{equation}\label{eq:local-cap-intensity}
 dJ_\lambda(k,L)=\frac{\lambda}{a_\lambda Px}
       Lb(L)\frac{V_{0,k+1}(0^k,L)}{P^{k-2}}
                         \frac{x^k}{k!}\,dL.
\end{equation}
This is the $k$-th term of the one-edge coefficient in
Appendix~\ref{app:positive-density-orbits}. By Mirzakhani's integration
formula the prelimit density of such waists per unit area is
$\binom nkLV_{0,k+1}(0^k,L)V_{g,n-k+1}(0^{n-k},L)/(A_{g,n}V_{g,n})$,
and the term-by-term limit used in the proof of
Proposition~\ref{rev:prop-density-convergence} shows that it converges
to the density in \eqref{eq:local-cap-intensity}, uniformly for $L$ in
compact subsets of $[0,\infty)$. A cap of this type has area
$A_k=2\pi(k-1)$.

Call a finite geodesic cap $C$ bad at level $h$ if
$\ell(\partial C)\le hA(C)$, and write $N_h(S,z)$ for the number of
such caps containing $z$, and $N_h^{(k)}(S,z)$ for those with $k$
cusps. Since $A(C)=A_k$ for such a cap, $N_h^{(k)}$ counts the caps
with $\ell(\partial C)\le hA_k$. The screened version
$N^{(k),\rm screen}_h:=N^{\varsigma}$ of
\eqref{rev:eq-screened-cap-count}, with
$\varsigma=(k,d,hA_k,\sigma,t)$, is a bounded Borel function of a
fixed ball, and in the finite surface the marked-root integration
formula gives
\[
 \mathbb EN^{(k),\rm screen}_h(X,o)
 =\frac1A\,\mathbb E\sum_{C\ \varsigma\text{-screened}}A(C^t)
 \le A_k\cdot\frac1A\,\mathbb E\#\{C:k\text{ cusps},\
   d\le\ell(\partial C)\le hA_k\}.
\]
The right side converges to $A_k\int_d^{hA_k}dJ_\lambda(k,L)$, so the
total variation transfer gives
\[
 \mathbb EN^{(k),\rm screen}_h(S_\lambda,o_\lambda)
 \le A_k\int_0^{hA_k}dJ_\lambda(k,L).
\]
As $\sigma$ exhausts the moduli
space, $t\to\infty$ and $d\downarrow0$, the screened counts increase
to $N^{(k)}_h(S_\lambda,o_\lambda)$, and monotone convergence gives
the first inequality in
\begin{align}\label{eq:local-cap-multiplicity}
 \mathbb E N_h(S,o)
 &\le\sum_{k\ge2}A_k\int_0^{hA_k}dJ_\lambda(k,L)\notag\\
 &\le C_\lambda h^2\sum_{k\ge2}k^{1/2}
                           e^{-(\xi-2\pi h)k}<\infty,
 \qquad 0<h\le\xi/(4\pi).
\end{align}
(The first line is in fact an equality: the screened prelimit
expectation is the integral of $A(C^t)\mathbf 1_\sigma$ against the
prelimit intensity, whose conditional law on the moduli of the cap
given the waist length is the normalized WP measure before and after
the limit; only the inequality is used.) The only constraint on $h$
used below is $\xi-2\pi h\ge\xi/2$, which makes the series converge.
Indeed $V_{0,k+1}(0^k,L)\le e^{L/2}V_{0,k+1}(0)$,
$b(L)\le e^{L/2}$ because $T(Y)\le\pi^2/4$ gives
$\frac L\pi\sqrt{T(Y)-\Delta(Y)s^2}\le L/2$ in
\eqref{eq:boundary-profile} and $I_0(u)\le e^u$ for $u\ge0$,
and the genus-zero coefficient asymptotic gives
\[
 \frac{V_{0,k+1}(0)}{P^{k-2}}\frac{x^k}{k!}
       \le C_\lambda k^{-5/2}e^{-\xi k}.
\]
Use $\int_0^a Le^L\,dL\le a^2e^a/2$ in
\eqref{eq:local-cap-intensity}. The extra area factor supplies the
third power of $k$. The estimate counts overlapping caps with their
multiplicities; no independence is required.

Fix one such $h$. Let $Z=\{z:N_h(S,z)=\infty\}$. It is an equivariant
measurable set and $\mathbb P(o\in Z)=0$. Apply
\eqref{eq:local-mass-transport} to
$m(S,u,v)=\mathbf1_{\{d(u,v)<R\}}\mathbf1_Z(v)$.
The received mass is zero almost surely, while the sent mass is
$A(Z\cap B_R(o))$. Hyperbolic balls have area bounded by that of
$B_R(\mathbb H^2)$, so these transports are finite. Countable
exhaustion over $R$ shows that $A(Z)=0$ on the entire sampled surface.
Consequently every positive-area finite geodesic core $Y_0$ is contained
in only finitely many bad caps: choose a point of its interior outside
$Z$. In particular, simultaneously for all such $Y_0$,
\begin{equation}\label{eq:anchored-bad-cap-bound}
 M_h(Y_0):=\sup\{A(C):C\text{ bad at level }h, Y_0\subset C\}
       <\infty,
\end{equation}
where the supremum of the empty set is zero. This uses no ergodic
decomposition and no uncountable intersection of rooted events.

\subsection{A topological anchor for the isoperimetric inequality}

The geodesic hull of a large domain containing a fixed non-elementary
core lies in an outer cap containing that core. The bound on bad caps
then forces a linear isoperimetric inequality for large domains,
\eqref{eq:eventual-anchored-isoperimetry}. The hull estimate is
deterministic and uses only the unique finite-cap property of the
first subsection.

We record the deterministic hull estimate. If a connected precompact
domain $\Omega$ has non-elementary fundamental-group image, take its
incompressible filling and geodesic hull $Y$. In the corresponding
cover, the complement of the convex core consists of funnels. The
field $\tanh t\,\partial_t$ in a funnel has divergence one, norm at
most one, and zero trace at the core. Apply Stokes' theorem to the lift
of the \emph{original} $\Omega$, extending the field by zero on the core.
The resulting old-area estimate and boundary straightening give
\[
 A(\Omega\setminus Y)\le\ell(\partial\Omega),\qquad
 \ell(\partial Y)\le\ell(\partial\Omega).
\]
No area of an added cusp tail is charged to the original domain.
By planarity and the unique finite-cap property, $Y$ is contained in
an outer finite cap $C$, whose outer waist is a boundary of $Y$.
Thus, writing $p_\Omega=\ell(\partial\Omega)$,
\begin{equation}\label{eq:anchored-hull-estimate}
 \ell(\partial C)\le p_\Omega,
       \qquad A(C)\ge A(\Omega)-p_\Omega.
\end{equation}
For elementary group images the plane or cyclic cover has a bounded
divergence-one field, giving $A(\Omega)\le p_\Omega$ directly.
In a hyperbolic cylinder use $\tanh t\,\partial_t$ on both sides;
in a parabolic cylinder or the plane use a unit horocyclic normal.

On a non-elementary surface in this geometric class, choose a finite-area
embedded geodesic core $Y_0$ with non-elementary group. It is obtained
by thickening noncommuting loops, filling the inessential complementary
pieces and tightening. Fix compact based generators of $\pi_1Y_0$.
Whenever $\Omega$ contains these loops, its hull subgroup contains
$\pi_1Y_0$ with the same basepoint. Their limit sets, and hence their
convex hulls in a common universal cover, are nested. Therefore
$Y_0\subset Y\subset C$. This retains an entire non-elementary core;
it does not assert that an arbitrary root survives straightening.

For any $o\in S$, choose a fixed smooth compact connected $K$ containing
$o$ and those generator loops in its interior, and put
$p_K=\ell(\partial K)$ and $\kappa=h/[2(1+h)]$.
There is a finite $V_*(S,o)\ge1$ such that every connected precompact
smooth domain $D\ni o$ satisfies
\begin{equation}\label{eq:eventual-anchored-isoperimetry}
 A(D)\ge V_*(S,o)\quad\Longrightarrow\quad
                 \ell(\partial D)\ge\kappa A(D).
\end{equation}
To see this, set $\Omega=D\cup K$, with outer smooth approximation
if necessary, and note that $p_\Omega\le p_D+p_K$.
If $p_D<\kappa A(D)$ and $A(D)>p_K/\kappa$, then
$p_\Omega<hA(D)/(1+h)$. Equation~\eqref{eq:anchored-hull-estimate}
makes its outer cap bad and gives $A(C)>A(D)/(1+h)$.
It contains $Y_0$, contradicting \eqref{eq:anchored-bad-cap-bound}
when $A(D)>(1+h)M_h(Y_0)$. Choose $V_*$ strictly above both
thresholds; the margin absorbs the smooth approximation errors.
If the entire surface is elementary the divergence-one argument
already gives \eqref{eq:eventual-anchored-isoperimetry} with constant one.

\subsection{Green levels and conformal type}

The eventual isoperimetric inequality bounds the areas of Green
superlevel sets uniformly over an exhaustion. The limiting Green kernel
is therefore finite, so Brownian motion is transient and the filled
surface is conformally a disk. The passage from the hyperbolic metric
to the filled surface uses that isolated points are polar and that
transience in dimension two depends only on the conformal
type~\cite{Grigoryan}.

We use the continuum Green-level argument of
Delmotte and Rau~\cite[Section~1.2]{DelmotteRau} (see~\cite{Grigoryan}
for the background on transience of Brownian motion), allowing the
large-area threshold in \eqref{eq:eventual-anchored-isoperimetry}.
For a smooth connected precompact $U\ni o$, let $G_U(o,\cdot)$ be
the Dirichlet Green function of $\Delta=-\operatorname{div}\nabla$
and put $v_U(t)=A\{G_U(o,\cdot)>t\}$ for $t>0$.
Every positive superlevel is connected and contains the pole: a
component without it would violate the maximum principle. For regular
levels, the unit source, coarea and Cauchy--Schwarz give
\[
 \int_{\{G_U=t\}}|\nabla G_U|\,d\ell=1,
 \qquad -v_U'(t)=\int_{\{G_U=t\}}\frac{d\ell}{|\nabla G_U|}
       \ge\ell(\partial\{G_U>t\})^2.
\]
While $v_U(t)\ge V_*$, the last term is at least
$\kappa^2v_U(t)^2$. Integrating the reciprocal differential inequality
from level zero, and using monotonicity after the threshold is crossed,
proves the uniform bound
\begin{equation}\label{eq:anchored-green-levels}
 v_U(t)\le\max\{V_*,(\kappa^2t)^{-1}\}.
\end{equation}
The coarea argument passes across critical values by approximation;
away from the pole the harmonic Green function is locally analytic.

Along a smooth connected exhaustion, the Dirichlet Green kernels
increase. If their limit diverged at one point away from $o$, positive
harmonic Harnack chains would force divergence everywhere off the pole.
For each fixed $t$ their superlevel areas would then tend to the
infinite area of $S$, contradicting \eqref{eq:anchored-green-levels}.
The limit is therefore the finite positive minimal Green kernel.
The original complete hyperbolic metric is stochastically complete;
the Green kernel proves transience for its Brownian motion.

Isolated filled cusp points are polar and do not change the capacity
dichotomy: the filled points form a countable closed discrete set,
single points are polar for two-dimensional Brownian motion, that is,
they have zero capacity, and a countable union of polar sets is polar,
so Brownian motion on the filled surface almost surely never meets
this set and adding or removing it does not change capacities. Transience is
equivalent to the positivity of the capacity of a compact set with
nonempty interior~\cite{Grigoryan}, and in dimension two the Dirichlet
integral, hence this capacity, is invariant under conformal changes of
the metric; equivalently, the Brownian motions of two conformal
metrics are time changes of each other, so recurrence and transience
depend only on the conformal type. The filled surface is simply
connected and noncompact,
so uniformization and transience force its conformal type to be
$\mathbb D$. Appendix~\ref{app:cusp-cutoffs}, using the conformal
distributional-limit theorem of~\cite[Theorem~1.4]{NPS}, proves
that the area-rooted dense endpoint is instead conformally $\mathbb C$
after filling and recurrent before filling. This completes the proof
of Theorem~\ref{thm:wp-conformal-phase}. The positive constant above
depends on $\lambda$ through $\xi>0$; no uniform estimate at the
recurrent endpoint follows from the argument.

\section{Growing WP balls and the Brownian plane}
\label{app:brownian-plane}

This appendix proves Theorem~\ref{thm:wp-brownian-plane}
(Section~\ref{sec:brownian-plane}): the area-rooted infinite WP sphere
$(S_\infty^*,o_\infty)$ converges after rescaling to the Brownian
plane, in the trimmed GHP form \eqref{eq:bp-ghp} and the untrimmed GP
form \eqref{eq:bp-gp}, and finite WP spheres converge in the same way
at mesoscopic radii. The argument follows the
finite-intermediate-size strategy of Curien and Le Gall~\cite{CLGPlane}.
The additional input for WP surfaces is a growing-prefix coupling for
the continuous-type trees of Budd and Curien~\cite{BC}, including their
actual area-root law. We compare prefixes of the finite and infinite
trees along a growing spine (Lemma~\ref{lem:bp-growing-prefix}), prove
a functional limit theorem for the labels along that spine
(Lemma~\ref{lem:bp-spine-fclt}), and pass to the uniform-area root.
These inputs couple balls of growing radius
(Proposition~\ref{prop:bp-growing-balls}); the finite Brownian-sphere
limit of~\cite{BC} in its rooted measured form
(Lemma~\ref{lem:bp-finite-ghp}) and the measured tangent coupling
\eqref{eq:bp-measured-tangent} of~\cite{CLGPlane} then complete the
proof in the last subsection. Throughout this appendix, references to
numbered results of~\cite{CLGPlane,LeGallArea} use the arXiv versions
listed in the bibliography.

\subsection{An exact prefix comparison}

We compute the exact likelihood ratio between the spine prefixes of a
finite marked tree and of the infinite tree, and show that prefixes of
$o(\sqrt N)$ branchings can be coupled with high probability
(Lemma~\ref{lem:bp-growing-prefix}). The inputs are the coefficient
asymptotics and the spine decomposition of~\cite{BC}, recorded in
\eqref{eq:bp-tree-coefficient} and \eqref{eq:bp-spine-split}; the
proof of the lemma also establishes a minorization of the spine
transition density, \eqref{eq:bp-doeblin}, which the next subsection
uses.

We use the tree series $F(x;\theta)$ of~\cite[Section~4.3,
equation~(18)]{BC}, not the volume series elsewhere in this paper. A tree with $n$ leaves,
including its root leaf, has $N=n-1$ nonroot leaves. Write
$a_N(\theta)=[x^N]F(x;\theta)$, let $\xi_*$ be its critical
fugacity, and set
\[
 F_*(\theta)=F(\xi_*;\theta),\qquad
 H(\theta)=J_0(j_{0,1}\theta/\pi),\qquad 0\le\theta<\pi.
\]
Lemma~6 and equations~(28) and~(31) of~\cite{BC} give
\begin{equation}\label{eq:bp-tree-coefficient}
 \frac{\xi_*^Na_N(\theta)}{H(\theta)}
     \sim c_*N^{-3/2},\qquad
 \sup_\theta\mathbb P_\theta(|T|=m)\le C m^{-3/2},
 \quad c_*=\frac{j_{0,1}}{2\sqrt2\,\pi^{5/2}}.
\end{equation}
The asymptotic is used pointwise in $\theta$. The function $F_*$ is
bounded above and below by positive constants, and $H(\theta)$ is
comparable to $\pi-\theta$. For an oriented split of the infinite
spine, the finite-child angle $a$ and infinite-child angle $b$ have
density
\begin{equation}\label{eq:bp-spine-split}
 \mathbf1_{\{\theta<a+b<\pi\}}
 \frac{F_*(a)H(b)}{H(\theta)}\,da\,db.
\end{equation}
There are two orientations. Conditional on the split, the finite child
has law $\mathbb P_a$, independently of the infinite child.

Condition a finite tree on $N$ nonroot leaves, mark one uniformly, and
follow its ancestral path for $h$ branching vertices. Retain the entire
off-path subtrees, all angles and the terminal stub. If the mark is
reached sooner, record a cemetery symbol. The mark is an auxiliary
randomization: the tree itself is recovered by forgetting it, so a
coupling of marked prefixes is in particular a coupling of the unmarked
trees outside the omitted subtree. In the infinite tree retain
the first $h$ spine branchings and their entire finite side trees.
Let $M$ be the sum of the nonroot leaf counts of these side trees and
$\theta_h$ the terminal type. Nonroot leaf counts add at a binary split,
so the omitted finite subtree has $N-M$ nonroot leaves.

For a specified prefix of combinatorial and Lebesgue weight $W$, the
finite marked density is
$W(N-M)a_{N-M}(\theta_h)/(Na_N(\theta_0))$. Multiplying
\eqref{eq:bp-spine-split} and the side-tree Boltzmann weights gives
infinite density $W\xi_*^MH(\theta_h)/H(\theta_0)$: the side tree at
the $i$th split, with finite-child angle $a_i$ and $M_i$ nonroot
leaves, has Boltzmann weight $\xi_*^{M_i}/F_*(a_i)$ times its angle
Lebesgue weight, which cancels the factor $F_*(a_i)$ of the splitting
density, while the factors $H(\theta_i)/H(\theta_{i-1})$ telescope
to $H(\theta_h)/H(\theta_0)$.
Consequently their exact likelihood ratio, off the cemetery, is
\begin{equation}\label{eq:bp-prefix-likelihood}
 R_{N,h}=\frac{N-M}{N}
 \frac{\xi_*^{N-M}a_{N-M}(\theta_h)/H(\theta_h)}
      {\xi_*^Na_N(\theta_0)/H(\theta_0)}
 \mathbf1_{\{M<N\}}.
\end{equation}
The infinite law puts no mass at the cemetery. With total variation
defined as the supremum of event-probability differences, this gives
\begin{equation}\label{eq:bp-prefix-tv-identity}
 d_{\rm TV}=1-\mathbb E_\infty\min(1,R_{N,h}).
\end{equation}
The factor $(N-M)/N$ comes from the uniform mark and is essential.

\begin{lemma}\label{lem:bp-growing-prefix}
For each fixed initial type $\theta_0<\pi$ and $0<\delta\le1$,
\begin{equation}\label{eq:bp-prefix-tv}
 \limsup_{N\to\infty}
 d_{\rm TV}(\operatorname{Prefix}_{N,\lfloor\delta\sqrt N\rfloor},
            \operatorname{Prefix}_{\infty,\lfloor\delta\sqrt N\rfloor})
 \le C\delta.
\end{equation}
In particular the distance tends to zero whenever $h_N=o(\sqrt N)$.
\end{lemma}

\begin{proof}
For a side-tree size $Z$, \eqref{eq:bp-tree-coefficient} gives
$\mathbb P(Z>x)\le Cx^{-1/2}$ and
$\mathbb E[Z\mathbf1_{\{Z\le x\}}]\le C\sqrt x$, uniformly in
its incoming angle. A union bound and truncation therefore give
\begin{equation}\label{eq:bp-prefix-mass}
 \mathbb P_\infty(M>x)\le Ch/\sqrt x,
\end{equation}
without requiring independence between successive spine types.

The spine transition density and its invariant probability density are
\[
 p(\theta,b)=\frac{2H(b)}{H(\theta)}
       \int_{\max(0,\theta-b)}^{\pi-b}F_*(a)\,da,
 \qquad
 \nu(b)=\frac{H(b)H'(\pi-b)}{\cos j_{0,1}}.
\]
Since $J_0'=-J_1$,
\[
 H'(\theta)=-\frac{j_{0,1}}{\pi}J_1\!\left(\frac{j_{0,1}\theta}{\pi}\right),
\]
and $J_1>0$ on $(0,j_{0,1}]$ with $J_1(t)\sim t/2$ as $t\downarrow0$;
hence $|H'(\pi-b)|\asymp\pi-b$ on $[0,\pi)$. Likewise $H(\theta)\asymp
\pi-\theta$, because the zero $j_{0,1}$ of $J_0$ is simple and $H>0$
on $[0,\pi)$. The integration interval has length
$\min(\pi-\theta,\pi-b)$. Hence, for some $\epsilon_0>0$,
\begin{equation}\label{eq:bp-doeblin}
 p(\theta,b)\,db\ge\epsilon_0\nu(b)\,db
 \quad\text{for every }\theta<\pi.
\end{equation}
The terminal-type law approaches $\nu$ geometrically in total variation.

Set $b_k(\theta)=k^{3/2}\xi_*^ka_k(\theta)/(c_*H(\theta))$.
For fixed $\eta>0$, the sets
\[
 E_N=\{\theta:\sup_{k\ge\lfloor N/2\rfloor}
                        |b_k(\theta)-1|>\eta\}
\]
are decreasing in $N$, and their intersection is empty in $[0,\pi)$
because $b_k(\theta)\to1$ pointwise by \eqref{eq:bp-tree-coefficient};
hence
$\nu(E_N)\to0$ by continuity of the finite measure $\nu$ from above.
If $h_N\to\infty$, \eqref{eq:bp-doeblin} implies
$\mathbb P(\theta_{h_N}\in E_N)\to0$. No independence from $M$
is used. On $M\le N/2$ and $\theta_h\notin E_N$, the likelihood
in \eqref{eq:bp-prefix-likelihood} is at least
$(1-\eta)/(1+\eta)$ for large $N$, since its asymptotic size factor
is $(1-M/N)^{-1/2}\ge1$. Equations
\eqref{eq:bp-prefix-tv-identity}--\eqref{eq:bp-prefix-mass}, followed
by $\eta\downarrow0$, prove \eqref{eq:bp-prefix-tv}. The same argument
gives the little-oh statement for unbounded $h_N$; for bounded
subsequences, apply dominated convergence under the fixed terminal-type
law. The case $h=0$ is immediate.
\end{proof}

\subsection{The labels along the spine}

The labels along the infinite spine are an additive functional of the
Markov chain of spine angles. The minorization \eqref{eq:bp-doeblin}
of the spine transition density permits a martingale decomposition,
and the martingale functional central limit theorem
of~\cite{WhittFCLT} gives a Brownian limit
(Lemma~\ref{lem:bp-spine-fclt}), which supplies the geometric
separation needed in Proposition~\ref{prop:bp-growing-balls}.

If $\theta_i$ is the next spine angle and $\beta_i$ its sibling angle,
the label increment is
\[
 \Delta_i=2\log\frac{\sin\theta_i}{\sin(\theta_i+\beta_i)},
 \qquad L_k=\sum_{i=1}^k\Delta_i.
\]

\begin{lemma}\label{lem:bp-spine-fclt}
For every fixed initial type $\theta_0<\pi$, there is
$0<\tau^2<\infty$ such that
\begin{equation}\label{eq:bp-spine-fclt}
 (h^{-1/2}L_{\lfloor ht\rfloor})_{0\le t\le1}
       \Longrightarrow (\tau B_t)_{0\le t\le1},
\end{equation}
where $B$ is standard Brownian motion.
\end{lemma}

The convergence in \eqref{eq:bp-spine-fclt} is in $D([0,1])$
with the Skorokhod $J_1$ topology. The constant $\tau^2$ does not
depend on $\theta_0$: it is the stationary variance
$\mathbb E_\nu D_1^2$ of the martingale difference constructed in the
proof, and a chain started at any $\theta_0$ couples with the
stationary chain after a geometrically bounded time. Budd and
Curien~\cite[Proposition~29]{BC} compute the asymptotic variance of
the same spine labels; its explicit value is not needed here.

\begin{proof}
Under stationarity, integration over the parent type in
\eqref{eq:bp-spine-split} gives the joint density of
$a=\theta_i,b=\beta_i$ as
\[
 \frac{-2}{\cos j_{0,1}}F_*(b)H(a)H(\pi-a-b)
              \mathbf1_{\{a,b>0,\ a+b<\pi\}}.
\]
This density is bounded and symmetric in $a$ and $c=\pi-a-b$.
Since $\Delta_i=2\log(\sin a/\sin c)$, its mean is zero and
all its polynomial moments are finite. Write
$m(\theta)=\mathbb E[\Delta_i\mid\theta_{i-1}=\theta]$.
The minorization \eqref{eq:bp-doeblin} decomposes the transition
operator as $P=\epsilon_0\Pi+(1-\epsilon_0)Q$, where
$\Pi f=\nu(f)$ and $Q$ preserves $\nu$. Thus $P$ contracts the
mean-zero subspace of $L^2(\nu)$ by $1-\epsilon_0$.
The series $g=\sum_{j\ge0}P^jm$ converges in that space and
$(I-P)g=m$. Therefore
\[
 D_i=\Delta_i+g(\theta_i)-g(\theta_{i-1})
\]
is a stationary ergodic martingale difference with finite second moment.
Ergodicity gives $h^{-1}\sum_{i\le ht}D_i^2\to t\mathbb ED_1^2$
almost surely. For fixed $T,\eta>0$,
\[
 \mathbb E\max_{i\le hT}\frac{D_i^2}{h}
 \le\eta^2+\frac{\lfloor hT\rfloor}{h}
       \mathbb E[D_1^2\mathbf1_{\{|D_1|>\eta\sqrt h\}}].
\]
Letting first $h\to\infty$ and then $\eta\downarrow0$, and using the
Cauchy--Schwarz inequality, verifies the maximum-jump hypothesis of the martingale FCLT
in~\cite[Theorem~2.1(i)]{WhittFCLT}. The same estimate with
$g(\theta_i)$ removes the coboundary uniformly on compact time intervals.
Conditional on the parent type and next angle $a$, the sibling angle
$b$ has positive density on a nonempty interval. The logarithmic
increment is nonconstant there, so adding $g(a)-g(\theta)$ cannot
remove its conditional variance. Hence $\tau^2=\mathbb ED_1^2>0$.
Finally use \eqref{eq:bp-doeblin} to couple a fixed-start chain to
stationarity after a geometrically bounded time, then use the same
sibling angles. The finitely many premeeting increments give a finite
path difference, negligible after division by $\sqrt h$.
\end{proof}

\subsection{The uniform-area root}

We now pass from the distinguished-leaf law to the actual area-root
law of~\cite{BC}, so that the two preceding lemmas apply to the
subtree at a uniform area point. The root decoration is coupled
between the finite and infinite spheres by a likelihood ratio of the
same form as \eqref{eq:bp-prefix-likelihood}, with the leaf-count
asymptotics replaced by a local convolution estimate.
Lemma~23 and Corollary~22 of~\cite{BC} describe its root blob
and the bias $(\deg-2+\mathbf1_{\rm red})/\deg$. Absorb this bounded
nonnegative bias into a normalized law $q$ for a finite root decoration
$d$. The decoration records its $K$ black children, the possible red
root leaf $\varepsilon\in\{0,1\}$, all root-blob angles, the selected
red vertex and corner, and the uniform area point in its ideal triangle.
Given $d$, the $K$ descendants are independent ordinary critical blob
trees. By the proof of Proposition~21 of~\cite{BC}, each is an ordinary
$\mathbb P_{\pi/2}$ tree with an artificial incoming root.

Let $Z$ count its red nonroot leaves and put $f_m=\mathbb P(Z=m)$.
Then $f_m\sim c m^{-3/2}$ and $f_m\le C m^{-3/2}$.
The number of original leaves is $\varepsilon+\sum_{i=1}^K Z_i$.
Moreover $K$ has an exponential tail, by Proposition~14 and
Corollary~22 of~\cite{BC}; their special-root probabilities differ
from the ordinary offspring probabilities by a bounded factor,
and the root bias is bounded.

For each fixed $K$, the local convolution estimate is
$\mathbb P(\sum_{i=1}^KZ_i=m)\sim Kf_m$. For $K=2$, split the
sum at a fixed cutoff $a$ and $m-a$: the middle contributes at most
$Cm^{-3/2}a^{-1/2}$. Let $a$ grow after $m$, and then use induction.
For all $K$, at least one summand is at least $m/K$, which gives
$\mathbb P(\sum Z_i=m)\le CK^{5/2}m^{-3/2}$ by conditioning on
the other summands. Exponential integrability of $K$ implies
\begin{equation}\label{eq:bp-root-normalizer}
 \mathcal Z_n:=\mathbb P_q(\varepsilon+\textstyle\sum_i Z_i=n)
       \sim\mathbb E_q[K]f_n.
\end{equation}

In the limiting law choose $d$ with distribution $Kq(d)/\mathbb E_qK$,
choose a child $I$ uniformly, place independent finite ordinary
$\mathbb P_{\pi/2}$ trees at its other children, and place at $I$ an
independent infinite tree with law $\mathbb P^\infty_{\pi/2}$, the
local limit of the conditioned laws $\mathbb P^N_{\pi/2}$
in~\cite[Theorem~9]{BC}. In the blob-tree description this tree at
$I$ is the $p$-Kesten tree started from a black vertex, after the
Red and Blob operations. Proposition~21 of~\cite{BC} identifies the Blob of the Red
version of their modified Kesten tree, whose root has the special
size-biased offspring law, with $\mathbb P^\infty_0$; its proof computes the probability that a fixed finite
bicolored tree is a prefix, and notes that a $p$-Galton--Watson tree
started from a black vertex without special root treatment is a
$\mathbb P_{\pi/2}$ tree. With the ordinary offspring law $p_k$ at
the root, the same computation gives the prefix probabilities of the
standard Kesten tree, whose root offspring law is $kp_k$, as the
limits of the conditioned prefix probabilities, with the ratio of
size asymptotics from their equation~(28) tending to one; so the
Blob of the black-rooted Kesten tree has law $\mathbb P^\infty_{\pi/2}$.
Include the finite trees in the root
decoration, and put $M=\varepsilon+\sum_{j\ne I}Z_j<\infty$.
In the finite conditioned law, at most one child exceeds $n/2$ leaves.
On that event the density of the resulting complete decoration relative
to this limit is exactly
\begin{equation}\label{eq:bp-root-likelihood}
 \frac{\mathbb E_q[K]f_{n-M}}{\mathcal Z_n}
                  \mathbf1_{\{n-M>n/2\}}.
\end{equation}
It tends pointwise to one and has integral at most one. The same
minimum-density argument as in \eqref{eq:bp-prefix-tv-identity}
therefore gives total-variation convergence of the complete decoration,
including its continuous angles and area point. Its remaining child
has $n-O_{\mathbb P}(1)$ leaves, with the ordinary conditioned
$\mathbb P_{\pi/2}$ law. This is precisely the root law in
\cite[Corollary~22]{BC}: the Kesten size bias commutes with the
bounded root bias. Lemmas~\ref{lem:bp-growing-prefix}
and~\ref{lem:bp-spine-fclt} apply to that child with
$\theta_0=\pi/2$; the infinite tree in those lemmas is the
$\mathbb P^\infty_{\pi/2}$ tree placed at $I$ above.

\subsection{Coupling balls with growing radii}

The prefix coupling of Lemma~\ref{lem:bp-growing-prefix}, the root
coupling \eqref{eq:bp-root-likelihood} and the spine label limit of
Lemma~\ref{lem:bp-spine-fclt} together couple balls in the finite and
infinite surfaces up to a small multiple of $n^{1/4}$: a sufficiently
negative label on the spine cuts off everything beyond it from the
ball, by the cactus bound of~\cite{BC}. After the proposition we
record why cusp trimming and the projected area are determined by the
local hyperbolic structure, so that the coupling also preserves the
measured balls used in Theorem~\ref{thm:wp-brownian-plane}.

\begin{proposition}\label{prop:bp-growing-balls}
There is $C<\infty$ such that, for $0<\alpha\le1$,
\begin{equation}\label{eq:bp-ball-tv}
 \limsup_{n\to\infty}d_{\rm TV}\bigl(
 B_{\alpha n^{1/4}}(S_n,o_n),
 B_{\alpha n^{1/4}}(S_\infty^*,o_\infty)\bigr)
 \le C\alpha^{2/3}.
\end{equation}
The balls retain their ambient distances and smooth hyperbolic
structures. Their total-variation distance tends to zero at every
deterministic radius $r_n=o(n^{1/4})$.
\end{proposition}

\begin{proof}
First couple the complete root decorations using
\eqref{eq:bp-root-likelihood}, at cost $o(1)$. Couple the first
$\mathsf d=\lfloor\delta\sqrt n\rfloor$ branchings of the remaining child,
including its entire finite side trees, by
Lemma~\ref{lem:bp-growing-prefix} applied with $h=\mathsf d$; the
letter $h$ is reserved below for the label function. The failure limsup is at most
$C\delta$. The random remaining size creates no uniformity problem:
first restrict the finite decoration's size to $M\le a$, take the
maximum over these finitely many size shifts, and then let $a\to\infty$.
The prefix distance is nondecreasing in the depth, and
$\lfloor\delta\sqrt n\rfloor\le\lfloor\delta'\sqrt{n-a}\rfloor$ for
every $\delta'>\delta$ and large $n$, so the bound $C\delta'$ holds for
every $\delta'>\delta$.

Let $x$ be the point of the incoming cut edge used to give its label
the value zero. Choose a piecewise smooth path from the root to $x$
inside the common finite decoration, and call its length $\ell_{\rm dec}$
(a length, not an area).
For example, use the finite chain of shared triangles, their incenters
and prescribed finite crossing points on the glued edges. This choice
is measurable in the decoration and has finite length almost surely.
It gives $d(o,x)\le \ell_{\rm dec}$ in both completions. In particular
$\ell_{\rm dec}$ is tight and $\ell_{\rm dec}/n^{1/4}\to0$ in probability;
the actual shortest distance need not be determined by the decoration.

Let $h(\cdot)$ be the signed-distance labeling of~\cite[Section~7.1]{BC},
normalized by $h(x)=0$. It is $1$-Lipschitz, so
$h(o)\ge-\ell_{\rm dec}$. The cactus bound
of~\cite[Lemma~32]{BC} says that, if an edge $e$ separates the
triangles of $o$ and $y$, then
\[
 d(o,y)\ge(h(o)-\ell(e)/2)_+
                    +(h(y)-\ell(e)/2)_+.
\]
Thus a spine edge with $\ell(e)<-6r-2\ell_{\rm dec}$ separates every
point beyond it from $B_{3r}(o)$. The common pruned tree consequently
determines this entire smooth neighborhood, as in~\cite[Lemma~33]{BC}.
It also determines ambient distances in $B_r(o)$: a minimizing
geodesic between two points of $B_r(o)$ has length at most $2r$ and
stays inside $B_{2r}(o)$.

Take $r=\alpha n^{1/4}$. Lemma~\ref{lem:bp-spine-fclt}, the
continuous mapping theorem and the reflection principle give
\begin{align*}
 \limsup_{n\to\infty}\mathbb P
  \{\min_{k\le\mathsf d}L_k\ge-6r-2\ell_{\rm dec}\}
 &\le\mathbb P\{\min_{0\le t\le1}B_t
                       \ge-6\alpha/(\tau\sqrt\delta)\}\\
 &\le C\alpha/\sqrt\delta.
\end{align*}
Tightness of the decoration suffices; independence from the labels is
unnecessary. Adding the prefix-coupling error and choosing
$\delta=\alpha^{2/3}$ proves \eqref{eq:bp-ball-tv}. For
$r_n=o(n^{1/4})$, restrict the coupled balls at fixed $\alpha$,
then send $\alpha\downarrow0$.
\end{proof}

Trimming is detectable from the local hyperbolic structure. In the
upper half-plane, the primitive parabolic $z\mapsto z+t$ has based
loop length $2\operatorname{arcsinh}(t/(2\operatorname{Im}z))$;
the unit cusp is $\operatorname{Im}z>t$. Hence a point belongs to an
open unit cusp exactly when some nonidentity parabolic based geodesic
loop there has length less than $2\operatorname{arcsinh}(1/2)<1$.
A nonprimitive power cannot add a point to this set. The loop lies
in the unit neighborhood of its base, and a smooth local isometry
preserves its developing holonomy and therefore its parabolic type.
This distinguishes cusp loops from short hyperbolic loops.
Consequently a smooth isometry on $B_{r+1}(o)$ identifies the trimmed
subsets and their horocycle boundaries in $B_r(o)$. The equality
boundary is detected by the same unit margin. The usual cusp collar
theorem makes the unit neighborhoods embedded and disjoint.

The spaces $X(S,o)$ are proper, since they are closed subsets of the
complete proper surfaces, with one point added. In a finite surface,
if $D_n=d(o_n,S_n^\circ)$, the disjoint cusp tails give
\begin{equation}\label{eq:bp-root-depth}
 \mathbb P(D_n>u)=\frac{(n+1)e^{-u}}{2\pi(n-1)},\qquad u>0.
\end{equation}
Thus root addition has vanishing cost at diverging distance scales.

\subsection{The finite rooted measured limit}

We need a joint metric, measure and root version of the finite
Brownian-sphere input, Lemma~\ref{lem:bp-finite-ghp}, for the trimmed
spaces $X(S_n,o_n)$ with the projected area $\bar A$ of
\eqref{eq:bp-projected-area}. It is obtained from the correspondence
in the proof of the Brownian-sphere limit of~\cite{BC}, by coupling
the area sample through its contour index and controlling the mass
moved by cusp projection. Write $(\mathcal M,D^*,\mu)$ for the
duration-one Brownian sphere and let $\rho$, conditionally on that
measured space, have distribution $\mu$.

\begin{lemma}\label{lem:bp-finite-ghp}
As $n\to\infty$,
\begin{equation}\label{eq:bp-finite-ghp}
 (X(S_n,o_n),n^{-1/4}d,\bar A/A_n,o_n)
 \Longrightarrow(\mathcal M,c_{\rm WP}D^*,\mu,\rho)
\end{equation}
in compact pointed GHP. The unprojected full-area version holds in GP.
\end{lemma}

\begin{proof}
\emph{The metric correspondence.}
Start with an arbitrary subsequence. In the Skorokhod representation
of~\cite[Section~8.5, equation~(67)]{BC}, pass to a further subsequence
on which the contour processes and the proxy pseudometrics converge
uniformly almost surely. The identification in their Section~8.6
shows that the limiting pseudometric is
$c_{\rm WP}D^*_{(\mathbf e,Z)}$, with $c_{\rm WP}=c_2/2$.
Write $\mathbf p:[0,1]\to\mathcal M$ for the quotient map, so that
$\mu=\mathbf p_*\operatorname{Leb}_{[0,1]}$. All subsequent almost
sure metric statements refer to this representation; the desired
conclusion for the original spaces is convergence in distribution.

Put $d_n=n^{-1/4}d$. For each of the first $n$ unit horocycles choose
a measurable representative point $z_{n,i}$ on its boundary. Let
$i_n(s)$ be the leaf whose rescaled contour visit time is closest to
$s$, with a fixed rule for ties, and put $z_n(s)=z_{n,i_n(s)}$.
Moving along a unit horocycle changes a distance by at most a fixed
constant. Thus the uniform proxy convergence gives
\[
 \delta_n:=\sup_{s,t\in[0,1]}
 \bigl|d_n(z_n(s),z_n(t))
       -c_{\rm WP}D^*(\mathbf p(s),\mathbf p(t))\bigr|\longrightarrow0
 \quad\text{almost surely}.
\]
Proposition~39 of~\cite{BC} makes all $n+1$ horocycles an
$O(n^{-1/4}\log^2n)$-net of $S_n^\circ$ with high probability.
Their equation~(72) places the remaining origin horocycle at
$o(n^{1/4})$ distance from one of the first $n$. Consequently
\[
 \eta_n:=n^{-1/4}+\sup_{y\in S_n^\circ}
                   \inf_{s\in[0,1]}d_n(y,z_n(s))
       \longrightarrow0\quad\text{in probability}.
\]
The explicit relation
\[
 \mathcal R_n=\{(y,\mathbf p(s)):
       y\in S_n^\circ,\ s\in[0,1],\ d_n(y,z_n(s))\le\eta_n\}
\]
is onto both spaces and has distortion at most $2\eta_n+\delta_n$.
In particular they admit isometric embeddings in a common metric
space $(E_n,d_{E_n})$ in which every pair in $\mathcal R_n$ has
distance at most $\kappa_n:=\eta_n+\delta_n/2+n^{-1/4}$.
This follows directly by defining cross-distances as the infimum of
$d_n(y,y')+\kappa_n+c_{\rm WP}D^*(z',z)$ over
$(y',z')\in\mathcal R_n$. Their Hausdorff distance is at most
$\kappa_n$, and $\kappa_n\to0$ in probability.

\emph{The area coupling in the same embedding.}
Conditionally on the surface and its coding, sample $Y_n$ from the
normalized full area $A/A_n$. Let $T_n$ be its nearest leaf-time as
in~\cite[Section~8.7]{BC}, and denote its conditional law by $q_n$.
Their equations~(73)--(74) give
\[
 d_n(Y_n,z_n(T_n))\longrightarrow0\quad\text{in probability},
 \qquad q_n\Longrightarrow\operatorname{Leb}_{[0,1]}
       \quad\text{in probability}.
\]
Here the first assertion is an annealed estimate and the second is
convergence of the random conditional law. Write $F_n$ for the
distribution function of $q_n$. Since the limiting distribution
function is continuous, $\sup_t|F_n(t)-t|\to0$ in probability.
Take $U_n$ conditionally uniform on $[0,1]$, set
$T_n=F_n^{-1}(U_n)$, and sample $Y_n$ from its conditional law given
$T_n$. This preserves the joint law of $(Y_n,T_n)$ and gives
$|T_n-U_n|\to0$ in probability. The conditional law of
$\mathbf p(U_n)$ is exactly $\mu$.

Let $P_n$ be the cusp projection of \eqref{eq:bp-projected-area}.
In cusp coordinates $(t,s)\in[0,\infty)\times\mathbb R/\mathbb Z$,
the metric and area are $dt^2+e^{-2t}ds^2$ and $e^{-t}ds\,dt$.
Projection to $t=0$ pushes cusp area to boundary arclength, and
\[
 \mathbb P\{d(Y_n,P_nY_n)>u\mid S_n\}
       =\frac{(n+1)e^{-u}}{A_n},\qquad u>0.
\]
It follows that $d_n(Y_n,P_nY_n)\to0$ in probability, while
$P_nY_n$ has law $\bar A/A_n$. Define the modulus
\[
 \omega(a)=\sup_{|s-t|\le a}
                 c_{\rm WP}D^*(\mathbf p(s),\mathbf p(t));
 \qquad \omega(a)\longrightarrow0\quad(a\downarrow0).
\]
In the embedding already constructed,
\begin{align*}
 d_{E_n}(P_nY_n,\mathbf p(U_n))
 &\le d_n(P_nY_n,Y_n)+d_n(Y_n,z_n(T_n))
       +\kappa_n+\omega(|T_n-U_n|)\\
 &\longrightarrow0\quad\text{in probability}.
\end{align*}
The terms involving $Y_n$ on the right are distances in the original
surface; the inequality uses their sum as an upper bound for
$d_n(P_nY_n,z_n(T_n))$. Thus the joint law of
$(P_nY_n,\mathbf p(U_n))$, conditionally on the coding, is a coupling
of $\bar A/A_n$ and $\mu$ in this same embedding. For each
$\epsilon>0$, the expected conditional mass of pairs farther than
$\epsilon$ tends to zero. Markov's inequality and the coupling
criterion for Prokhorov distance therefore show that the Prokhorov
distance tends to zero in probability, jointly with the Hausdorff bound.

\emph{The root and the full-area version.}
Use this very area sample as $o_n=Y_n$, and put
$\rho_n=\mathbf p(U_n)$. Conditionally on the limiting measured
sphere, $\rho_n$ has law $\mu$, so its pointed law is the one on the
right of \eqref{eq:bp-finite-ghp}. Glue the original surface
$(S_n,d_n)$ to $E_n$ along its isometric closed subspace $S_n^\circ$.
This preserves all distances in both spaces. In the resulting
embedding, $d_{E_n}(o_n,\rho_n)\to0$ in probability by the preceding
bound and $d_n(Y_n,P_nY_n)\to0$. Adding $o_n$ to $S_n^\circ$ changes
the Hausdorff bound by at most this latter distance. We have therefore
controlled metric, measure and root simultaneously, proving the
pointed GHP assertion along the extracted subsequence. Every
subsequence admits such an extraction, with the same limiting pointed
law, which proves \eqref{eq:bp-finite-ghp} for the full sequence.
Keeping the full area and the same root in this embedding gives the
pointed GP assertion by the identical coupling, without a Hausdorff
requirement on the untrimmed surface.
\end{proof}

For completeness, $\bar A$ on the infinite surface is locally finite.
The depth interval $0\le t\le1$ over those boundary portions in
$B_R(o)$ has area $(1-e^{-1})$ times their total arclength. These
swaths have disjoint interiors and lie in $B_{R+1}(o)$, whose area
is finite. This bounds the boundary part of
\eqref{eq:bp-projected-area}; the interior part is locally finite too.
The local holonomy argument above shows that the exact ball coupling
also preserves $\bar A$ after allowing a unit margin. It already
preserves the raw hyperbolic area.

\subsection{The measured tangent coupling}

We record why the Brownian-sphere tangent coupling
of~\cite[Proposition~4]{CLGPlane} also preserves the natural volume;
this measured form
\eqref{eq:bp-measured-tangent} is the last input of the proof of
Theorem~\ref{thm:wp-brownian-plane} in the next subsection. The
argument reads off the coding windows in their proof and uses that the
time translations between them preserve Lebesgue measure.
For each $r>0$ and $\varepsilon>0$, for all sufficiently
large $\Lambda$, the construction in~\cite[Proposition~4]{CLGPlane}
couples
\begin{equation}\label{eq:bp-measured-tangent}
 (\mathcal M,\Lambda D^*,\Lambda^4\mu,\rho)
 \quad\text{and}\quad
 (\mathcal P,D_\infty,\mu_\infty,\rho_\infty)
\end{equation}
with isometric measured $r$-balls with probability at least
$1-\varepsilon$. We write $\Lambda$ for the scale parameter denoted
$\lambda$ in~\cite{CLGPlane}, since $\lambda$ is the cusp density in
this part of the paper, and $A>1$ for the constant fixed at the start
of their proof of Proposition~4. The measured assertion follows from
their proof, as follows.

The finite coding excursion has duration $T=\Lambda^4$, and its
volume is the image of $dt$ on $[0,T]$. On the event $F_\Lambda$ of
that construction, Lemma~5(i) of~\cite{CLGPlane} puts \emph{every}
coding time whose image belongs to the finite $r$-ball in
\[
 [0,\gamma_\Lambda(A)]\cup[\eta_\Lambda(A),T].
\]
Their Lemma~5(ii) puts all preimages of the plane's $r$-ball in
$[-\eta_\infty(A),\gamma_\infty(A)]$. Lemma~6 identifies both
ball membership and quotient distances under the maps $s\mapsto s$
on the first finite interval and $s\mapsto s-T$ on the second.
On $F_\Lambda$ the corresponding endpoints match. These translations
preserve Lebesgue measure, with at most a null endpoint duplication.
Since there are no additional preimages of either ball outside the
displayed windows, the induced isometry preserves the entire
restricted volume measure. Rerooting invariance of the Brownian sphere
\cite[Theorem~8.1 and its proof]{LeGallGeodesics} also holds with its
volume measure: the multipointed identities in that theorem determine
the rooted metric-measure law. It identifies the conditional volume root in
Lemma~\ref{lem:bp-finite-ghp} with the distinguished root of this
construction. This proves \eqref{eq:bp-measured-tangent}.

\subsection{Passage through an intermediate finite surface}

We couple a large ball of the infinite surface with the same ball of a
finite surface of intermediate size
(Proposition~\ref{prop:bp-growing-balls}), take the Brownian-sphere
limit of the latter (Lemma~\ref{lem:bp-finite-ghp}), and then pass to
the Brownian plane through the measured tangent coupling
\eqref{eq:bp-measured-tangent}. The restriction to balls requires
radii of zero boundary mass, which are supplied by a Fubini argument
and by~\cite[Proposition~1]{LeGallArea} for the Brownian plane. The
proof treats the trimmed statement \eqref{eq:bp-ghp}, the untrimmed
statement \eqref{eq:bp-gp}, and the finite spheres in turn.

\begin{proof}[Proof of Theorem~\ref{thm:wp-brownian-plane}]
\emph{The intermediate surface.}
Fix an observation radius $R>0$ and $0<\alpha<1$. For $L\to\infty$
set $m_L=\lceil(RL/\alpha)^4\rceil$. Proposition
\ref{prop:bp-growing-balls}, with the unit margin for cusp projection,
couples the measured $RL$-balls of the infinite surface and $S_{m_L}$
with failure limsup at most $C\alpha^{2/3}$. Rounding and the unit
margin can be absorbed by replacing $\alpha$ by $2\alpha$.

For fixed $\alpha$, Lemma~\ref{lem:bp-finite-ghp} gives the compact
pointed GHP limit of this intermediate finite surface after rescaling
distance by $L^{-1}$ and mass by $(v_{\rm WP}L^4)^{-1}$:
\[
 (\mathcal M,\Lambda D^*,\Lambda^4\mu,\rho),
 \qquad \Lambda=c_{\rm WP}R/\alpha.
\]
Indeed, $m_L^{1/4}/L\to R/\alpha$ and
\[
 \frac{A_{m_L}}{v_{\rm WP}L^4}
 \longrightarrow\frac{2\pi R^4}{v_{\rm WP}\alpha^4}
       =(c_{\rm WP}R/\alpha)^4.
\]
This fixes the volume constant in \eqref{eq:bp-normalization}.

\emph{Restriction to balls and the trimmed limit.}
Restriction to closed balls is continuous at a geodesic compact
limiting space when the limiting boundary has zero mass. For the
metric assertion, the outer inclusion follows by compactness;
each point on the limiting boundary is approximable from inside
along a root geodesic. This does not require the prelimit spaces
to be length spaces. The measure assertion follows from the
zero-boundary-mass condition. For a finite random measure, the
radial pushforward has at most countably many atoms in each
realization. Fubini therefore supplies a full Lebesgue-measure set
of deterministic radii with almost surely zero boundary mass.
Choose $\alpha\downarrow0$ through values for which
$\alpha/c_{\rm WP}$ is such a radius on the sphere. To treat
finitely or countably many smaller observation radii $R_i\le R$,
also avoid the null sets for $(R_i/R)\alpha/c_{\rm WP}$.

For fixed such $\alpha$, first send $L\to\infty$ and restrict the
compact GHP limit to the observation ball. Then let $\alpha\downarrow0$
in the measured coupling \eqref{eq:bp-measured-tangent}.
For a bounded continuous test of the rooted measured ball, the
remaining discrepancy is bounded by $C\alpha^{2/3}$ times its
oscillation plus the tangent-coupling error. Both tend to zero.
The argument holds jointly at any finite list of radii using the
same largest-ball coupling. This proves \eqref{eq:bp-ghp}. Indeed, let $G$ be bounded and
$1$-Lipschitz for the local distance of
Section~\ref{sec:brownian-plane}, and write $X_L$ for the rescaled
space. For a proper pointed measured space $Y$, the ball $B_R(Y)$,
viewed as a pointed measured space, has the same $r$-balls as $Y$ for
$r<R$; their local distance is therefore at most
$\int_R^\infty e^{-r}\,dr=e^{-R}$, and $|G(Y)-G(B_R(Y))|\le e^{-R}$.
The map $K\mapsto G(K)$ on compact pointed measured spaces is
continuous at $K=B_R(\mathcal P)$ for the compact GHP topology: the
Brownian plane is geodesic and has zero radial boundary mass at every
positive radius~\cite[Proposition~1]{LeGallArea}, so the $r$-balls
converge for every $r<R$; for $r>R$ the $r$-ball of
$B_R(\mathcal P)$ is the whole space and the $r$-ball of an
approximating compact space is eventually the whole space, so these
converge too. Dominated convergence then applies to the capped
integrand. Hence
$\limsup_L|\mathbb EG(X_L)-\mathbb EG(\mathcal P)|\le2e^{-R}$ for
every $R$.

\emph{The untrimmed limit.}
For \eqref{eq:bp-gp}, compare raw and projected area on the same
intermediate finite surface. The mass transported more than
$\epsilon$ in the rescaled metric is
\begin{equation}\label{eq:bp-bad-transport}
 \frac{(m_L+1)e^{-\epsilon L}}{v_{\rm WP}L^4}
       =O_{\alpha,R}(e^{-\epsilon L})\longrightarrow0.
\end{equation}
Their total rescaled masses stay finite for fixed $\alpha$ and $R$.
Thus they become close in pointed GP, in the common ambient
embedding with the same root. Restriction at zero boundary mass,
the raw smooth ball coupling and then
\eqref{eq:bp-measured-tangent} give the asserted raw-area limit.
No Hausdorff estimate on deep cusp tails is involved.

\emph{Finite spheres.}
Finally, if $1\ll r_n\ll n^{1/4}$, Proposition
\ref{prop:bp-growing-balls} couples every $Rr_n$-ball of $S_n$
to that of $S_\infty^*$ with failure probability tending to zero.
Both the raw and projected measures are preserved. Apply the
infinite-volume limits just proved with $L=r_n$ and then let
$R\to\infty$. This proves both finite mesoscopic statements.
\end{proof}

\clearpage
\section*{Statements and Declarations}
\noindent\textit{Competing interests.} The authors have no relevant financial or non-financial interests to disclose.

\noindent\textit{Data availability.} Data sharing is not applicable to this article, as no datasets were generated or analyzed.

\noindent\textit{Formalization.} For the accompanying Lean~4 formalization, see the \href{https://github.com/FrontierX-Lab/WPFormalization}{WPFormalization repository}.

\section*{Disclosure of automated assistance}
OpenAI GPT-6.1 Sol, OpenAI GPT-6 Astra, Anthropic Claude Fable 5.1, Anthropic Claude Fable 5.2, Anthropic Claude Opus 5.5, and Google Gemini 4 Argon were used throughout this work: to develop and check the arguments, to review the proofs and their hypotheses, to search and verify the literature, and to draft, organize and typeset the manuscript and prepare its figures.


\begin{thebibliography}{195}
\bibitem{AbertBiringer} Mikl\'os Ab\'ert and Ian Biringer, \emph{Unimodular measures on the space of all Riemannian manifolds}, Geom. Topol. \textbf{26} (2022), no.~5, 2295--2404. \url{https://doi.org/10.2140/gt.2022.26.2295}; \href{https://arxiv.org/abs/1606.03360v6}{arXiv:1606.03360v6}.
\bibitem{ABBGNRS} Mikl\'os Ab\'ert, Nicolas Bergeron, Ian Biringer, Tsachik Gelander, Nikolay Nikolov, Jean Raimbault, and Iddo Samet, \emph{On the growth of $L^2$-invariants for sequences of lattices in Lie groups}, Ann. of Math. (2) \textbf{185} (2017), no.~3, 711--790. \url{https://doi.org/10.4007/annals.2017.185.3.1}.
\bibitem{ADH} Romain Abraham, Jean-Fran\c{c}ois Delmas, and Patrick Hoscheit, \emph{A note on the Gromov--Hausdorff--Prokhorov distance between (locally) compact metric measure spaces}, Electron. J. Probab. \textbf{18} (2013), article~14, 21~pp. \url{https://doi.org/10.1214/EJP.v18-2116}.
\bibitem{AggAbelian} Amol Aggarwal, \emph{Large genus asymptotics for volumes of strata of abelian differentials}, with an appendix by Anton Zorich, J. Amer. Math. Soc. \textbf{33} (2020), no.~4, 941--989. \url{https://doi.org/10.1090/jams/947}; \href{https://arxiv.org/abs/1804.05431v3}{arXiv:1804.05431v3}.
\bibitem{Agg} Amol Aggarwal, \emph{Large genus asymptotics for intersection numbers and principal strata volumes of quadratic differentials}, Invent. Math. \textbf{226} (2021), no.~3, 897--1010. \url{https://doi.org/10.1007/s00222-021-01059-9}; \href{https://arxiv.org/abs/2004.05042}{arXiv:2004.05042}.
\bibitem{AMNOResurgence} Mohamad Alameddine, Olivier Marchal, Nikita Nikolaev, and Nicolas Orantin, \emph{Resurgence of tritronqu\'ees solutions of the deformed Painlev\'e~I equation}, preprint (2025), \href{https://arxiv.org/abs/2503.12290v2}{arXiv:2503.12290v2}.
\bibitem{AMNO} Mohamad Alameddine, Olivier Marchal, Nikita Nikolaev, and Nicolas Orantin, \emph{The WKB Riemann--Hilbert problem for the Painlev\'e~I system}, preprint (2026), \href{https://arxiv.org/abs/2610.20123v1}{arXiv:2610.20123v1}.
\bibitem{AldousLyons} David Aldous and Russell Lyons, \emph{Processes on unimodular random networks}, Electron. J. Probab. \textbf{12} (2007), 1454--1508. \url{https://doi.org/10.1214/EJP.v12-463}.
\bibitem{AM} Nalini Anantharaman and Laura Monk, \emph{A high-genus asymptotic expansion of Weil--Petersson volume polynomials}, J. Math. Phys. \textbf{63} (2022), no.~4, article~043502. \url{https://doi.org/10.1063/5.0039385}; \href{https://arxiv.org/abs/2011.14889}{arXiv:2011.14889}.
\bibitem{AMPartOne} Nalini Anantharaman and Laura Monk, \emph{Friedman--Ramanujan functions in random hyperbolic geometry and application to spectral gaps}, preprint (2023), \href{https://arxiv.org/abs/2304.02678v1}{arXiv:2304.02678v1}; version~4 (2026) is titled \emph{\dots\ application to spectral gaps~I}.
\bibitem{AMGap} Nalini Anantharaman and Laura Monk, \emph{Friedman--Ramanujan functions in random hyperbolic geometry and application to spectral gaps II}, preprint (2025), \href{https://arxiv.org/abs/2502.12268v2}{arXiv:2502.12268v2}.
\bibitem{AMTwoNinths} Nalini Anantharaman and Laura Monk, \emph{Typical hyperbolic surfaces have a spectral gap greater than $2/9-\epsilon$}, Int. Math. Res. Not. IMRN \textbf{2026} (2026), no.~15, article~rnag117. \url{https://doi.org/10.1093/imrn/rnag117}; \href{https://arxiv.org/abs/2604.09792v1}{arXiv:2604.09792v1}.
\bibitem{AHNR} Omer Angel, Tom Hutchcroft, Asaf Nachmias, and Gourab Ray, \emph{Hyperbolic and parabolic unimodular random maps}, Geom. Funct. Anal. \textbf{28} (2018), no.~4, 879--942. \url{https://doi.org/10.1007/s00039-018-0446-y}.
\bibitem{AC} Enrico Arbarello and Maurizio Cornalba, \emph{Combinatorial and algebro-geometric cohomology classes on the moduli spaces of curves}, J. Algebraic Geom. \textbf{5} (1996), no.~4, 705--749; \href{https://arxiv.org/abs/alg-geom/9406008}{arXiv:alg-geom/9406008}.
\bibitem{ALW} Siva Athreya, Wolfgang L\"ohr, and Anita Winter, \emph{The gap between Gromov-vague and Gromov--Hausdorff-vague topology}, Stochastic Process. Appl. \textbf{126} (2016), no.~9, 2527--2553. \url{https://doi.org/10.1016/j.spa.2016.02.009}.
\bibitem{BMM} Werner Ballmann, Henrik Matthiesen, and Sugata Mondal, \emph{Small eigenvalues of surfaces of finite type}, Compos. Math. \textbf{153} (2017), no.~8, 1747--1768. \url{https://doi.org/10.1112/S0010437X17007291}; \href{https://arxiv.org/abs/1506.06541v1}{arXiv:1506.06541v1}.
\bibitem{BMP} Werner Ballmann, Sugata Mondal, and Panagiotis Polymerakis, \emph{Small eigenvalues of pseudo-Laplacians}, preprint (2025), \href{https://arxiv.org/abs/2512.19248v1}{arXiv:2512.19248v1}.
\bibitem{BGL} Simon Barazer, Alessandro Giacchetto, and Mingkun Liu, \emph{Length spectrum of large genus random metric maps}, Forum Math. Sigma \textbf{13} (2025), article~e70. \url{https://doi.org/10.1017/fms.2025.31}; \href{https://arxiv.org/abs/2312.10517}{arXiv:2312.10517}.
\bibitem{BHS} Ara Basmajian, Hrant Hakobyan, and Dragomir \v{S}ari\'c, \emph{The type problem for Riemann surfaces via Fenchel--Nielsen parameters}, Proc. Lond. Math. Soc. (3) \textbf{125} (2022), no.~3, 568--625. \url{https://doi.org/10.1112/plms.12465}; \href{https://arxiv.org/abs/2011.03166}{arXiv:2011.03166}.
\bibitem{Beardon} Alan F. Beardon, \emph{The exponent of convergence of Poincar\'e series}, Proc. Lond. Math. Soc. (3) \textbf{18} (1968), no.~3, 461--483. \url{https://doi.org/10.1112/plms/s3-18.3.461}.
\bibitem{Bender} Edward A. Bender, \emph{Central and local limit theorems applied to asymptotic enumeration}, J. Combin. Theory Ser. A \textbf{15} (1973), no.~1, 91--111. \url{https://doi.org/10.1016/0097-3165(73)90038-1}.
\bibitem{BenjaminiSchramm} Itai Benjamini and Oded Schramm, \emph{Recurrence of distributional limits of finite planar graphs}, Electron. J. Probab. \textbf{6} (2001), article~23, 13~pp. \url{https://doi.org/10.1214/EJP.v6-96}.
\bibitem{BEG} Ga\"etan Borot, Bertrand Eynard, and Alessandro Giacchetto, \emph{The factorial growth of topological recursion}, Lett. Math. Phys. \textbf{115} (2025), no.~3, article~62. \url{https://doi.org/10.1007/s11005-025-01950-z}; \href{https://arxiv.org/abs/2409.17838}{arXiv:2409.17838}.
\bibitem{Borthwick} David Borthwick, \emph{Spectral Theory of Infinite-Area Hyperbolic Surfaces}, 2nd ed., Progress in Mathematics \textbf{318}, Birkh\"auser/Springer, Cham, 2016. \url{https://doi.org/10.1007/978-3-319-33877-4}.
\bibitem{BJP} David Borthwick, Chris Judge, and Peter A. Perry, \emph{Selberg's zeta function and the spectral geometry of geometrically finite hyperbolic surfaces}, Comment. Math. Helv. \textbf{80} (2005), no.~3, 483--515. \url{https://doi.org/10.4171/CMH/23}; \href{https://arxiv.org/abs/math/0310364v2}{arXiv:math/0310364v2}.
\bibitem{BFJ} S\'ebastien Boucksom, Charles Favre, and Mattias Jonsson, \emph{Differentiability of volumes of divisors and a problem of Teissier}, J. Algebraic Geom. \textbf{18} (2009), no.~2, 279--308. \url{https://doi.org/10.1090/S1056-3911-08-00490-6}; \href{https://arxiv.org/abs/math/0608260}{arXiv:math/0608260}.
\bibitem{BowditchTriang} Brian H. Bowditch, \emph{Bilipschitz triangulations of riemannian manifolds}, preprint (2020), \url{https://bhbowditch.com/papers/triangulations.pdf} (version dated 11 April 2020).
\bibitem{BK} \'Edouard Br\'ezin and Vladimir A. Kazakov, \emph{Exactly solvable field theories of closed strings}, Phys. Lett. B \textbf{236} (1990), no.~2, 144--150. \url{https://doi.org/10.1016/0370-2693(90)90818-Q}.
\bibitem{BrooksMakover} Robert Brooks and Eran Makover, \emph{Random construction of Riemann surfaces}, J. Differential Geom. \textbf{68} (2004), no.~1, 121--157. \url{https://doi.org/10.4310/jdg/1102536712}.
\bibitem{Budd} Timothy Budd, \emph{Irreducible metric maps and Weil--Petersson volumes}, Comm. Math. Phys. \textbf{394} (2022), no.~2, 887--917. \url{https://doi.org/10.1007/s00220-022-04418-6}; \href{https://arxiv.org/abs/2012.11318v2}{arXiv:2012.11318v2}.
\bibitem{BC} Timothy Budd and Nicolas Curien, \emph{Random punctured hyperbolic surfaces \& the Brownian sphere}, preprint (2025), \href{https://arxiv.org/abs/2508.18792v1}{arXiv:2508.18792v1}.
\bibitem{BL} Timothy Budd and Tanguy Lions, \emph{The tight length spectrum of large-genus random hyperbolic surfaces with many cusps}, preprint (2025), \href{https://arxiv.org/abs/2506.02611v3}{arXiv:2506.02611v3} (version~3, 2026).
\bibitem{BMZ} Timothy Budd, Thomas Meeusen, and Bart Zonneveld, \emph{A tree bijection for the moduli space of genus-0 hyperbolic surfaces with boundaries}, preprint (2025), \href{https://arxiv.org/abs/2512.09722v1}{arXiv:2512.09722v1}.
\bibitem{BZ} Timothy Budd and Bart Zonneveld, \emph{Topological recursion of the Weil--Petersson volumes of hyperbolic surfaces with tight boundaries}, J. Math. Phys. \textbf{65} (2024), no.~9, article~092302. \url{https://doi.org/10.1063/5.0192711}; \href{https://arxiv.org/abs/2307.04708}{arXiv:2307.04708}.
\bibitem{BudzinskiHBP} Thomas Budzinski, \emph{The hyperbolic Brownian plane}, Probab. Theory Related Fields \textbf{171} (2018), no.~1--2, 503--541. \url{https://doi.org/10.1007/s00440-017-0785-x}; \href{https://arxiv.org/abs/1604.06622}{arXiv:1604.06622}.
\bibitem{BCPCheeger} Thomas Budzinski, Nicolas Curien, and Bram Petri, \emph{On Cheeger constants of hyperbolic surfaces}, Invent. Math. \textbf{242} (2025), no.~2, 511--530. \url{https://doi.org/10.1007/s00222-025-01361-w}.
\bibitem{BudzinskiLouf} Thomas Budzinski and Baptiste Louf, \emph{Local limits of uniform triangulations in high genus}, Invent. Math. \textbf{223} (2021), no.~1, 1--47. \url{https://doi.org/10.1007/s00222-020-00986-3}; \href{https://arxiv.org/abs/1902.00492}{arXiv:1902.00492}.
\bibitem{Burger} Marc Burger, \emph{Small eigenvalues of Riemann surfaces and graphs}, Math. Z. \textbf{205} (1990), no.~1, 395--420. \url{https://doi.org/10.1007/BF02571252}.
\bibitem{Buser} Peter Buser, \emph{A note on the isoperimetric constant}, Ann. Sci. \'Ecole Norm. Sup. (4) \textbf{15} (1982), no.~2, 213--230. \url{https://doi.org/10.24033/asens.1426}.
\bibitem{BuserBook} Peter Buser, \emph{Geometry and Spectra of Compact Riemann Surfaces}, Progress in Mathematics \textbf{106}, Birkh\"auser, Boston, MA, 1992; reprint, Modern Birkh\"auser Classics, 2010. \url{https://doi.org/10.1007/978-0-8176-4992-0}.
\bibitem{Cheeger} Jeff Cheeger, \emph{A lower bound for the smallest eigenvalue of the Laplacian}, in: Problems in Analysis (Papers dedicated to Salomon Bochner, 1969), Princeton University Press, Princeton, NJ, 1970, 195--199. \url{https://doi.org/10.1515/9781400869312-013}.
\bibitem{CGT} Jeff Cheeger, Mikhail Gromov, and Michael Taylor, \emph{Finite propagation speed, kernel estimates for functions of the Laplace operator, and the geometry of complete Riemannian manifolds}, J. Differential Geom. \textbf{17} (1982), no.~1, 15--53. \url{https://doi.org/10.4310/jdg/1214436699}.
\bibitem{CMSZ} Dawei Chen, Martin M\"oller, Adrien Sauvaget, and Don Zagier, \emph{Masur--Veech volumes and intersection theory on moduli spaces of Abelian differentials}, Invent. Math. \textbf{222} (2020), no.~1, 283--373. \url{https://doi.org/10.1007/s00222-020-00969-4}; \href{https://arxiv.org/abs/1901.01785}{arXiv:1901.01785}.
\bibitem{CMSQuad} Dawei Chen, Martin M\"oller, and Adrien Sauvaget, \emph{Masur--Veech volumes and intersection theory: the principal strata of quadratic differentials}, Duke Math. J. \textbf{172} (2023), no.~9, 1735--1779. \url{https://doi.org/10.1215/00127094-2022-0063}.
\bibitem{ChenMullane} Dawei Chen and Scott Mullane, \emph{A tale of two volumes of moduli spaces: Weil--Petersson and Masur--Veech}, preprint (2026), \href{https://arxiv.org/abs/2603.07318v1}{arXiv:2603.07318v1}.
\bibitem{ChengYau} Shiu-Yuen Cheng and Shing-Tung Yau, \emph{Differential equations on Riemannian manifolds and their geometric applications}, Comm. Pure Appl. Math. \textbf{28} (1975), no.~3, 333--354. \url{https://doi.org/10.1002/cpa.3160280303}.
\bibitem{ColboisCourtois} Bruno Colbois and Gilles Courtois, \emph{Les valeurs propres inf\'erieures \`a $1/4$ des surfaces de Riemann de petit rayon d'injectivit\'e}, Comment. Math. Helv. \textbf{64} (1989), no.~1, 349--362. \url{https://doi.org/10.1007/BF02564681}.
\bibitem{CCH} Ovidiu Costin, Rodica D. Costin, and Min Huang, \emph{Tronqu\'ee solutions of the Painlev\'e equation $P_{\rm I}$}, Constr. Approx. \textbf{41} (2015), no.~3, 467--494. \url{https://doi.org/10.1007/s00365-015-9287-1}; \href{https://arxiv.org/abs/1310.5330v2}{arXiv:1310.5330v2}.
\bibitem{CurienPSHT} Nicolas Curien, \emph{Planar stochastic hyperbolic triangulations}, Probab. Theory Related Fields \textbf{165} (2016), no.~3--4, 509--540. \url{https://doi.org/10.1007/s00440-015-0638-4}.
\bibitem{CLGPlane} Nicolas Curien and Jean-Fran\c{c}ois Le Gall, \emph{The Brownian plane}, J. Theoret. Probab. \textbf{27} (2014), no.~4, 1249--1291. \url{https://doi.org/10.1007/s10959-013-0485-0}; \href{https://arxiv.org/abs/1204.5921v1}{arXiv:1204.5921v1}.
\bibitem{DGZZBounds} Vincent Delecroix, \'Elise Goujard, Peter Zograf, and Anton Zorich, \emph{Uniform lower bound for intersection numbers of $\psi$-classes}, SIGMA Symmetry Integrability Geom. Methods Appl. \textbf{16} (2020), article~086, 13~pp. \url{https://doi.org/10.3842/SIGMA.2020.086}.
\bibitem{DGZZ} Vincent Delecroix, \'Elise Goujard, Peter Zograf, and Anton Zorich, \emph{Masur--Veech volumes, frequencies of simple closed geodesics, and intersection numbers of moduli spaces of curves}, Duke Math. J. \textbf{170} (2021), no.~12, 2633--2718. \url{https://doi.org/10.1215/00127094-2021-0054}. Original preprint (2019): \href{https://arxiv.org/abs/1908.08611v1}{arXiv:1908.08611v1}.
\bibitem{DGZZRandom} Vincent Delecroix, \'Elise Goujard, Peter Zograf, and Anton Zorich, \emph{Large genus asymptotic geometry of random square-tiled surfaces and of random multicurves}, Invent. Math. \textbf{230} (2022), no.~1, 123--224. \url{https://doi.org/10.1007/s00222-022-01123-y}.
\bibitem{DelmotteRau} Thierry Delmotte and Cl\'ement Rau, \emph{A note on exit time for anchored isoperimetry}, Ann. Fac. Sci. Toulouse Math. (6) \textbf{24} (2015), no.~4, 817--835. \url{https://doi.org/10.5802/afst.1466}. Preprint version: \emph{Exit time for anchored expansion}, \href{https://arxiv.org/abs/0903.3892v2}{arXiv:0903.3892v2} (section numbers cited from the published version; its Section~1.2 is Section~1.3 of the preprint).
\bibitem{DHokerPhong} Eric D'Hoker and Duong H. Phong, \emph{On determinants of Laplacians on Riemann surfaces}, Comm. Math. Phys. \textbf{104} (1986), no.~4, 537--545. \url{https://doi.org/10.1007/BF01211063}.
\bibitem{DVV} Robbert Dijkgraaf, Herman Verlinde, and Erik Verlinde, \emph{Loop equations and Virasoro constraints in non-perturbative two-dimensional quantum gravity}, Nuclear Phys. B \textbf{348} (1991), no.~3, 435--456. \url{https://doi.org/10.1016/0550-3213(91)90199-8}.
\bibitem{Do} Norman Do, \emph{Moduli spaces of hyperbolic surfaces and their Weil--Petersson volumes}, in: Handbook of Moduli, Vol.~I, Adv. Lect. Math. \textbf{24}, International Press, Somerville, MA, 2013, 217--258; \href{https://arxiv.org/abs/1103.4674v1}{arXiv:1103.4674v1}.
\bibitem{DoNorbury} Norman Do and Paul Norbury, \emph{Weil--Petersson volumes and cone surfaces}, Geom. Dedicata \textbf{141} (2009), no.~1, 93--107. \url{https://doi.org/10.1007/s10711-008-9345-y}; \href{https://arxiv.org/abs/math/0603406}{arXiv:math/0603406}.
\bibitem{DPRS} Jozef Dodziuk, Thea Pignataro, Burton Randol, and Dennis Sullivan, \emph{Estimating small eigenvalues of Riemann surfaces}, in: The Legacy of Sonya Kovalevskaya (Cambridge, MA, and Amherst, MA, 1985), Contemp. Math. \textbf{64}, Amer. Math. Soc., Providence, RI, 1987, 93--121. \url{https://doi.org/10.1090/conm/064/881458}.
\bibitem{DS} Michael R. Douglas and Stephen H. Shenker, \emph{Strings in less than one dimension}, Nuclear Phys. B \textbf{335} (1990), no.~3, 635--654. \url{https://doi.org/10.1016/0550-3213(90)90522-F}.
\bibitem{ET} Lorenz Eberhardt and Gustavo J. Turiaci, \emph{2D dilaton gravity and the Weil--Petersson volumes with conical defects}, Comm. Math. Phys. \textbf{405} (2024), no.~4, article~103. \url{https://doi.org/10.1007/s00220-024-04964-1}; \href{https://arxiv.org/abs/2304.14948v2}{arXiv:2304.14948v2}.
\bibitem{Efrat} Isaac Efrat, \emph{Determinants of Laplacians on surfaces of finite volume}, Comm. Math. Phys. \textbf{119} (1988), no.~3, 443--451. \url{https://doi.org/10.1007/BF01218082}.
\bibitem{EZ} Alex Eskin and Anton Zorich, \emph{Volumes of strata of Abelian differentials and Siegel--Veech constants in large genera}, Arnold Math. J. \textbf{1} (2015), no.~4, 481--488. \url{https://doi.org/10.1007/s40598-015-0028-0}; \href{https://arxiv.org/abs/1507.05296}{arXiv:1507.05296}.
\bibitem{EynardMumford} Bertrand Eynard, \emph{Recursion between Mumford volumes of moduli spaces}, Ann. Henri Poincar\'e \textbf{12} (2011), no.~8, 1431--1447. \url{https://doi.org/10.1007/s00023-011-0113-4}; \href{https://arxiv.org/abs/0706.4403}{arXiv:0706.4403}.
\bibitem{EGGGL} Bertrand Eynard, Elba Garcia-Failde, Alessandro Giacchetto, Paolo Gregori, and Danilo Lewa\'nski, \emph{Resurgent large genus asymptotics of intersection numbers}, Forum Math. Pi \textbf{14} (2026), article~e21. \url{https://doi.org/10.1017/fmp.2026.10037}; \href{https://arxiv.org/abs/2309.03143v1}{arXiv:2309.03143v1}.
\bibitem{EGGLS} Bertrand Eynard, Elba Garcia-Failde, Paolo Gregori, Danilo Lewa\'nski, and Ricardo Schiappa, \emph{Resurgent asymptotics of Jackiw--Teitelboim gravity and the nonperturbative topological recursion}, Ann. Henri Poincar\'e \textbf{25} (2024), no.~9, 4121--4193. \url{https://doi.org/10.1007/s00023-023-01412-z}; \href{https://arxiv.org/abs/2305.16940}{arXiv:2305.16940}.
\bibitem{EO} Bertrand Eynard and Nicolas Orantin, \emph{Weil--Petersson volume of moduli spaces, Mirzakhani's recursion and matrix models}, preprint (2007), \href{https://arxiv.org/abs/0705.3600v1}{arXiv:0705.3600v1}.
\bibitem{FarbMargalit} Benson Farb and Dan Margalit, \emph{A Primer on Mapping Class Groups}, Princeton Mathematical Series \textbf{49}, Princeton University Press, Princeton, NJ, 2012.
\bibitem{FS} Philippe Flajolet and Robert Sedgewick, \emph{Analytic Combinatorics}, Cambridge University Press, Cambridge, 2009. \url{https://doi.org/10.1017/CBO9780511801655}.
\bibitem{GIKM} Stavros Garoufalidis, Alexander Its, Andrei Kapaev, and Marcos Mari\~no, \emph{Asymptotics of the instantons of Painlev\'e I}, Int. Math. Res. Not. IMRN \textbf{2012} (2012), no.~3, 561--606. \url{https://doi.org/10.1093/imrn/rnr029}; \href{https://arxiv.org/abs/1002.3634v3}{arXiv:1002.3634v3}.
\bibitem{GT} Malo Gatignol and Joe Thomas, \emph{Lengths of closed geodesics on random non-compact hyperbolic surfaces}, preprint (2026), \href{https://arxiv.org/abs/2610.41002v1}{arXiv:2610.41002v1}.
\bibitem{GillRohde} James T. Gill and Steffen Rohde, \emph{On the Riemann surface type of random planar maps}, Rev. Mat. Iberoam. \textbf{29} (2013), no.~3, 1071--1090. \url{https://doi.org/10.4171/RMI/749}; \href{https://arxiv.org/abs/1101.1320}{arXiv:1101.1320}.
\bibitem{GregoriSchiappa} Paolo Gregori and Ricardo Schiappa, \emph{From minimal strings towards Jackiw--Teitelboim gravity: on their resurgence, resonance, and black holes}, Classical Quantum Gravity \textbf{41} (2024), no.~11, article~115001. \url{https://doi.org/10.1088/1361-6382/ad32df}; \href{https://arxiv.org/abs/2108.11409}{arXiv:2108.11409}.
\bibitem{GPW} Andreas Greven, Peter Pfaffelhuber, and Anita Winter, \emph{Convergence in distribution of random metric measure spaces ($\Lambda$-coalescent measure trees)}, Probab. Theory Related Fields \textbf{145} (2009), no.~1--2, 285--322. \url{https://doi.org/10.1007/s00440-008-0169-3}.
\bibitem{Grigoryan} Alexander Grigor'yan, \emph{Analytic and geometric background of recurrence and non-explosion of the Brownian motion on Riemannian manifolds}, Bull. Amer. Math. Soc. (N.S.) \textbf{36} (1999), no.~2, 135--249. \url{https://doi.org/10.1090/S0273-0979-99-00776-4}.
\bibitem{GNY} Alexander Grigor'yan, Yuri Netrusov, and Shing-Tung Yau, \emph{Eigenvalues of elliptic operators and geometric applications}, Surv. Differ. Geom. \textbf{9} (2004), 147--217. \url{https://doi.org/10.4310/SDG.2004.v9.n1.a5}; author manuscript: \url{https://www.math.uni-bielefeld.de/~grigor/esceps.pdf}.
\bibitem{GPRS} Luca Griguolo, Jacopo Papalini, Lorenzo Russo, and Domenico Seminara, \emph{Asymptotics of Weil--Petersson volumes and two-dimensional quantum gravities}, SciPost Phys. \textbf{17} (2024), no.~6, article~156. \url{https://doi.org/10.21468/SciPostPhys.17.6.156}; \href{https://arxiv.org/abs/2402.07276v2}{arXiv:2402.07276v2}.
\bibitem{GLNDet} Renan Gross, Guy Lachman, and Asaf Nachmias, \emph{Determinants of Laplacians on converging hyperbolic surfaces}, preprint (2026), \href{https://arxiv.org/abs/2601.00338v1}{arXiv:2601.00338v1}.
\bibitem{GM} David J. Gross and Alexander A. Migdal, \emph{Nonperturbative two-dimensional quantum gravity}, Phys. Rev. Lett. \textbf{64} (1990), no.~2, 127--130. \url{https://doi.org/10.1103/PhysRevLett.64.127}.
\bibitem{Gru} Samuel Grushevsky, \emph{An explicit upper bound for Weil--Petersson volumes of the moduli spaces of punctured Riemann surfaces}, Math. Ann. \textbf{321} (2001), no.~1, 1--13. \url{https://doi.org/10.1007/PL00004496}.
\bibitem{Guillermit} Marthe Guillermit, \emph{A refinement of the asymptotic expansion of Weil--Petersson volumes}, preprint (2026), \href{https://arxiv.org/abs/2607.20341v1}{arXiv:2607.20341v1}.
\bibitem{GY} Jindong Guo and Di Yang, \emph{On the large genus asymptotics of psi-class intersection numbers}, Math. Ann. \textbf{388} (2024), no.~1, 61--97. \url{https://doi.org/10.1007/s00208-022-02505-6}; \href{https://arxiv.org/abs/2110.06774v3}{arXiv:2110.06774v3}.
\bibitem{GYZ} Jindong Guo, Di Yang, and Don Zagier, \emph{On uniform large genus asymptotics of Witten's intersection numbers}, preprint (2026), \href{https://arxiv.org/abs/2603.15233v2}{arXiv:2603.15233v2}.
\bibitem{GHS} Qi Guo, Bobo Hua, and Yang Shen, \emph{Spectral gaps for noncompact hyperbolic surfaces with linearly many cusps}, preprint (2025), \href{https://arxiv.org/abs/2507.16794v3}{arXiv:2507.16794v3} (version~3, 2026).
\bibitem{GGN} Ori Gurel-Gurevich and Asaf Nachmias, \emph{Recurrence of planar graph limits}, Ann. of Math. (2) \textbf{177} (2013), no.~2, 761--781. \url{https://doi.org/10.4007/annals.2013.177.2.10}; \href{https://arxiv.org/abs/1206.0707}{arXiv:1206.0707}.
\bibitem{GPY} Larry Guth, Hugo Parlier, and Robert Young, \emph{Pants decompositions of random surfaces}, Geom. Funct. Anal. \textbf{21} (2011), no.~5, 1069--1090. \url{https://doi.org/10.1007/s00039-011-0131-x}.
\bibitem{Hassannezhad} Asma Hassannezhad, \emph{Conformal upper bounds for the eigenvalues of the Laplacian and Steklov problem}, J. Funct. Anal. \textbf{261} (2011), no.~12, 3419--3436. \url{https://doi.org/10.1016/j.jfa.2011.08.003}.
\bibitem{HMO} Yasuyuki Hatsuda, Takaki Matsumoto, and Kazumi Okuyama, \emph{Non-perturbative effects in JT gravity from KdV equations}, J. High Energy Phys. \textbf{2026} (2026), no.~3, article~056. \url{https://doi.org/10.1007/JHEP03(2026)056}; \href{https://arxiv.org/abs/2505.16433v1}{arXiv:2505.16433v1}.
\bibitem{HeWuDet} Yuxin He and Yunhui Wu, \emph{Averages of determinants of Laplacians over moduli spaces for large genus}, J. Lond. Math. Soc. (2) \textbf{112} (2025), no.~6, article~e70395. \url{https://doi.org/10.1112/jlms.70395}; \href{https://arxiv.org/abs/2411.12971v2}{arXiv:2411.12971v2}.
\bibitem{HeWu} Yuxin He and Yunhui Wu, \emph{On second eigenvalues of closed hyperbolic surfaces for large genus}, J. Differential Geom. \textbf{130} (2025), no.~2, 403--441. \url{https://doi.org/10.4310/jdg/1747157756}; numbering as in \href{https://arxiv.org/abs/2207.12919v2}{arXiv:2207.12919v2}.
\bibitem{HSWX} Yuxin He, Yang Shen, Yunhui Wu, and Yuhao Xue, \emph{Non-simple systoles on random hyperbolic surfaces for large genus}, preprint (2023), \href{https://arxiv.org/abs/2308.16447v2}{arXiv:2308.16447v2} (version~2, 2025).
\bibitem{HWX} Yuxin He, Yunhui Wu, and Yuhao Xue, \emph{Uniform spectral gaps for random hyperbolic surfaces with not many cusps}, preprint (2026), \href{https://arxiv.org/abs/2602.08352v1}{arXiv:2602.08352v1}.
\bibitem{Hejhal} Dennis A. Hejhal, \emph{The Selberg Trace Formula for $\mathrm{PSL}(2,\mathbb R)$. Vol.~2}, Lecture Notes in Math. \textbf{1001}, Springer, Berlin, 1983. \url{https://doi.org/10.1007/BFb0061302}.
\bibitem{HejhalDegen} Dennis A. Hejhal, \emph{Regular $b$-groups, degenerating Riemann surfaces, and spectral theory}, Mem. Amer. Math. Soc. \textbf{88} (1990), no.~437. \url{https://doi.org/10.1090/memo/0437}.
\bibitem{Hide} Will Hide, \emph{Spectral gap for Weil--Petersson random surfaces with cusps}, Int. Math. Res. Not. IMRN \textbf{2023} (2023), no.~20, 17411--17460. \url{https://doi.org/10.1093/imrn/rnac293}; \href{https://arxiv.org/abs/2107.14555v4}{arXiv:2107.14555v4}.
\bibitem{HMT} Will Hide, Davide Macera, and Joe Thomas, \emph{Spectral gap with polynomial rate for Weil--Petersson random surfaces}, preprint (2025), \href{https://arxiv.org/abs/2508.14874v2}{arXiv:2508.14874v2}.
\bibitem{HideMagee} Will Hide and Michael Magee, \emph{Near optimal spectral gaps for hyperbolic surfaces}, Ann. of Math. (2) \textbf{198} (2023), no.~2, 791--824. \url{https://doi.org/10.4007/annals.2023.198.2.6}; \href{https://arxiv.org/abs/2107.05292v3}{arXiv:2107.05292v3}.
\bibitem{HTLargeN} Will Hide and Joe Thomas, \emph{Large-$n$ asymptotics for Weil--Petersson volumes of moduli spaces of bordered hyperbolic surfaces}, Comm. Math. Phys. \textbf{406} (2025), no.~9, article~203. \url{https://doi.org/10.1007/s00220-025-05369-4}; \href{https://arxiv.org/abs/2312.11412v4}{arXiv:2312.11412v4}.
\bibitem{HTSphere} Will Hide and Joe Thomas, \emph{Short geodesics and small eigenvalues on random hyperbolic punctured spheres}, Comment. Math. Helv. \textbf{100} (2025), no.~3, 463--506. \url{https://doi.org/10.4171/CMH/588}; \href{https://arxiv.org/abs/2209.15568v3}{arXiv:2209.15568v3}.
\bibitem{HTMany} Will Hide and Joe Thomas, \emph{Small eigenvalues of hyperbolic surfaces with many cusps}, preprint (2024), \href{https://arxiv.org/abs/2410.06093v2}{arXiv:2410.06093v2} (version~2, July 2026).
\bibitem{HuangGap} Guan Huang, \emph{Spectral gaps near one quarter for Weil--Petersson random surfaces with growing cusp number}, preprint (2026), \href{https://arxiv.org/abs/2611.01477v1}{arXiv:2611.01477v1}.
\bibitem{XHuangSuper} Xuanyu Huang, \emph{Large genus asymptotics of super Weil--Petersson volumes}, preprint (2025), \href{https://arxiv.org/abs/2501.07848v1}{arXiv:2501.07848v1}.
\bibitem{HLXSuper} Xuanyu Huang, Kefeng Liu, and Hao Xu, \emph{The super Weil--Petersson asymptotic expansion and the coefficients}, preprint (2026), \href{https://arxiv.org/abs/2611.06621v1}{arXiv:2611.06621v1}.
\bibitem{HJLCounting} Jonathan Huntley, Jay Jorgenson, and Rolf Lundelius, \emph{On the asymptotic behavior of counting functions associated to degenerating hyperbolic Riemann surfaces}, J. Funct. Anal. \textbf{149} (1997), no.~1, 58--82. \url{https://doi.org/10.1006/jfan.1997.3097}; numbering as in \href{https://arxiv.org/abs/math/9412221v1}{arXiv:math/9412221v1}.
\bibitem{Hwang} Hsien-Kuei Hwang, \emph{On convergence rates in the central limit theorems for combinatorial structures}, European J. Combin. \textbf{19} (1998), no.~3, 329--343. \url{https://doi.org/10.1006/eujc.1997.0179}.
\bibitem{IZ} Claude Itzykson and Jean-Bernard Zuber, \emph{Combinatorics of the modular group II: the Kontsevich integrals}, Internat. J. Modern Phys. A \textbf{7} (1992), no.~23, 5661--5705. \url{https://doi.org/10.1142/S0217751X92002581}; \href{https://arxiv.org/abs/hep-th/9201001}{arXiv:hep-th/9201001}.
\bibitem{Iwaniec} Henryk Iwaniec, \emph{Spectral Methods of Automorphic Forms}, 2nd ed., Grad. Stud. Math. \textbf{53}, Amer. Math. Soc., Providence, RI, 2002. \url{https://doi.org/10.1090/gsm/053}.
\bibitem{JansonLouf} Svante Janson and Baptiste Louf, \emph{Short cycles in high genus unicellular maps}, Ann. Inst. Henri Poincar\'e Probab. Stat. \textbf{58} (2022), no.~3, 1547--1564. \url{https://doi.org/10.1214/21-AIHP1218}.
\bibitem{JiGreen} Lizhen Ji, \emph{The asymptotic behavior of Green's functions for degenerating hyperbolic surfaces}, Math. Z. \textbf{212} (1993), no.~1, 375--394. \url{https://doi.org/10.1007/BF02571664}.
\bibitem{Johnson} Clifford V. Johnson, \emph{A boundary-addition identity for Weil--Petersson volumes and their supersymmetric generalizations}, preprint (2026), \href{https://arxiv.org/abs/2609.21030v1}{arXiv:2609.21030v1}.
\bibitem{JR} Clifford V. Johnson and Jo\~ao Rodrigues, \emph{Non-perturbative data for Weil--Petersson volumes and intersection numbers using ordinary differential equations}, J. High Energy Phys. \textbf{2026} (2026), no.~8, article~083. \url{https://doi.org/10.1007/JHEP08(2026)083}; \href{https://arxiv.org/abs/2601.03351v2}{arXiv:2601.03351v2}.
\bibitem{JLHeat} Jay Jorgenson and Rolf Lundelius, \emph{A regularized heat trace for hyperbolic Riemann surfaces of finite volume}, Comment. Math. Helv. \textbf{72} (1997), no.~4, 636--659. \url{https://doi.org/10.1007/S000140050039}.
\bibitem{JSSDet} Jay Jorgenson, Lejla Smajlovi\'c, and Polyxeni Spilioti, \emph{Asymptotics of the spectral determinant of the weighted Laplacian for a sequence of compact Riemann surfaces of infinitely growing volume}, preprint (2026), \href{https://arxiv.org/abs/2607.26779v1}{arXiv:2607.26779v1}.
\bibitem{JK} Nalini Joshi and Alexander V. Kitaev, \emph{On Boutroux's tritronqu\'ee solutions of the first Painlev\'e equation}, Stud. Appl. Math. \textbf{107} (2001), no.~3, 253--291. \url{https://doi.org/10.1111/1467-9590.00187}.
\bibitem{Kallenberg} Olav Kallenberg, \emph{Random Measures, Theory and Applications}, Probability Theory and Stochastic Modelling \textbf{77}, Springer, Cham, 2017. \url{https://doi.org/10.1007/978-3-319-41598-7}.
\bibitem{Kap} Andrei A. Kapaev, \emph{Quasi-linear Stokes phenomenon for the Painlev\'e first equation}, J. Phys. A \textbf{37} (2004), no.~46, 11149--11167. \url{https://doi.org/10.1088/0305-4470/37/46/005}.
\bibitem{KMZ} Ralph M. Kaufmann, Yuri I. Manin, and Don Zagier, \emph{Higher Weil--Petersson volumes of moduli spaces of stable $n$-pointed curves}, Comm. Math. Phys. \textbf{181} (1996), no.~3, 763--787. \url{https://doi.org/10.1007/BF02101297}.
\bibitem{Keel} Sean Keel, \emph{Basepoint freeness for nef and big line bundles in positive characteristic}, Ann. of Math. (2) \textbf{149} (1999), no.~1, 253--286. \url{https://doi.org/10.2307/121025}; \href{https://arxiv.org/abs/math/9901149}{arXiv:math/9901149}.
\bibitem{Kim} Yusuke Kimura, \emph{JT gravity and the asymptotic Weil--Petersson volume}, Phys. Lett. B \textbf{811} (2020), article~135989. \url{https://doi.org/10.1016/j.physletb.2020.135989}; \href{https://arxiv.org/abs/2008.04141v2}{arXiv:2008.04141v2}.
\bibitem{Kontsevich} Maxim Kontsevich, \emph{Intersection theory on the moduli space of curves and the matrix Airy function}, Comm. Math. Phys. \textbf{147} (1992), no.~1, 1--23. \url{https://doi.org/10.1007/BF02099526}.
\bibitem{Laz} Robert Lazarsfeld, \emph{Positivity in Algebraic Geometry I}, Ergeb. Math. Grenzgeb. (3) \textbf{48}, Springer, Berlin, 2004. \url{https://doi.org/10.1007/978-3-642-18808-4}.
\bibitem{LeGallArea} Jean-Fran\c{c}ois Le Gall, \emph{The area of spheres in the Brownian plane}, Probab. Theory Related Fields \textbf{196} (2026), no.~3--4, 1047--1074. \url{https://doi.org/10.1007/s00440-025-01419-z}; \href{https://arxiv.org/abs/2504.15860v2}{arXiv:2504.15860v2}.
\bibitem{LeGallGeodesics} Jean-Fran\c{c}ois Le Gall, \emph{Geodesics in large planar maps and in the Brownian map}, Acta Math. \textbf{205} (2010), no.~2, 287--360. \url{https://doi.org/10.1007/s11511-010-0056-5}; \href{https://arxiv.org/abs/0804.3012}{arXiv:0804.3012}.
\bibitem{LeGallUniq} Jean-Fran\c{c}ois Le Gall, \emph{Uniqueness and universality of the Brownian map}, Ann. Probab. \textbf{41} (2013), no.~4, 2880--2960. \url{https://doi.org/10.1214/12-AOP792}.
\bibitem{LMSQE} Etienne Le Masson and Tuomas Sahlsten, \emph{Quantum ergodicity and Benjamini--Schramm convergence of hyperbolic surfaces}, Duke Math. J. \textbf{166} (2017), no.~18, 3425--3460. \url{https://doi.org/10.1215/00127094-2017-0027}.
\bibitem{LMSScattering} Etienne Le Masson and Tuomas Sahlsten, \emph{Quantum ergodicity for Eisenstein series on hyperbolic surfaces of large genus}, Math. Ann. \textbf{389} (2024), no.~1, 845--898. \url{https://doi.org/10.1007/s00208-023-02671-1}; \href{https://arxiv.org/abs/2006.14935v3}{arXiv:2006.14935v3}.
\bibitem{LWZenodo} Zhiyi Li and Qihang Wang, \emph{Uniform asymptotics of Weil--Petersson volumes and their coefficients}, Zenodo preprint, version~1.0, published 16 October 2026; earlier version of the present paper. \url{https://doi.org/10.5281/zenodo.23414478}.
\bibitem{LipWright} Michael Lipnowski and Alex Wright, \emph{Towards optimal spectral gaps in large genus}, Ann. Probab. \textbf{52} (2024), no.~2, 545--575. \url{https://doi.org/10.1214/23-AOP1657}; \href{https://arxiv.org/abs/2103.07496}{arXiv:2103.07496}.
\bibitem{LX} Kefeng Liu and Hao Xu, \emph{A remark on Mirzakhani's asymptotic formulae}, Asian J. Math. \textbf{18} (2014), no.~1, 29--52. \url{https://doi.org/10.4310/AJM.2014.v18.n1.a2}; \href{https://arxiv.org/abs/1103.5136}{arXiv:1103.5136}.
\bibitem{LXBalancing} Kefeng Liu and Hao Xu, \emph{Balancing tautological intersection numbers}, preprint (2026), \href{https://arxiv.org/abs/2610.39462v1}{arXiv:2610.39462v1}.
\bibitem{LXConstant} Kefeng Liu and Hao Xu, \emph{Uniform asymptotics and the Weil--Petersson constant}, preprint (submitted 24 October 2026), \href{https://arxiv.org/abs/2610.33752v1}{arXiv:2610.33752v1}.
\bibitem{LXEquiv} Kefeng Liu and Hao Xu, \emph{Mirzakhani's recursion formula is equivalent to the Witten--Kont\-se\-vich theorem}, Ast\'erisque \textbf{328} (2009), 223--235; \href{https://arxiv.org/abs/0705.2086v2}{arXiv:0705.2086v2}.
\bibitem{LXHigher} Kefeng Liu and Hao Xu, \emph{Recursion formulae of higher Weil--Petersson volumes}, Int. Math. Res. Not. IMRN \textbf{2009} (2009), no.~5, 835--859. \url{https://doi.org/10.1093/imrn/rnn148}; \href{https://arxiv.org/abs/0708.0565}{arXiv:0708.0565}.
\bibitem{MNP} Michael Magee, Fr\'ed\'eric Naud, and Doron Puder, \emph{A random cover of a compact hyperbolic surface has relative spectral gap $\frac{3}{16}-\varepsilon$}, Geom. Funct. Anal. \textbf{32} (2022), no.~3, 595--661. \url{https://doi.org/10.1007/s00039-022-00602-x}.
\bibitem{MPvH} Michael Magee, Doron Puder, and Ramon van Handel, \emph{Strong convergence of uniformly random permutation representations of surface groups}, preprint (2025), \href{https://arxiv.org/abs/2504.08988v2}{arXiv:2504.08988v2}.
\bibitem{ManZ} Yuri I. Manin and Peter Zograf, \emph{Invertible cohomological field theories and Weil--Petersson volumes}, Ann. Inst. Fourier (Grenoble) \textbf{50} (2000), no.~2, 519--535. \url{https://doi.org/10.5802/aif.1764}; \href{https://arxiv.org/abs/math/9902051}{arXiv:math/9902051}.
\bibitem{MT} Henry Maxfield and Gustavo J. Turiaci, \emph{The path integral of 3D gravity near extremality; or, JT gravity with defects as a matrix integral}, J. High Energy Phys. \textbf{2021} (2021), no.~1, article~118. \url{https://doi.org/10.1007/JHEP01(2021)118}; \href{https://arxiv.org/abs/2006.11317v2}{arXiv:2006.11317v2}.
\bibitem{McShane} Greg McShane, \emph{Simple geodesics and a series constant over Teichmuller space}, Invent. Math. \textbf{132} (1998), no.~3, 607--632. \url{https://doi.org/10.1007/s002220050235}.
\bibitem{MertensTuriaci} Thomas G. Mertens and Gustavo J. Turiaci, \emph{Solvable models of quantum black holes: a review on Jackiw--Teitelboim gravity}, Living Rev. Relativ. \textbf{26} (2023), article~4. \url{https://doi.org/10.1007/s41114-023-00046-1}.
\bibitem{Miermont} Gr\'egory Miermont, \emph{The Brownian map is the scaling limit of uniform random plane quadrangulations}, Acta Math. \textbf{210} (2013), no.~2, 319--401. \url{https://doi.org/10.1007/s11511-013-0096-8}.
\bibitem{MirInv} Maryam Mirzakhani, \emph{Simple geodesics and Weil--Petersson volumes of moduli spaces of bordered Riemann surfaces}, Invent. Math. \textbf{167} (2007), no.~1, 179--222. \url{https://doi.org/10.1007/s00222-006-0013-2}.
\bibitem{MirJAMS} Maryam Mirzakhani, \emph{Weil--Petersson volumes and intersection theory on the moduli space of curves}, J. Amer. Math. Soc. \textbf{20} (2007), no.~1, 1--23. \url{https://doi.org/10.1090/S0894-0347-06-00526-1}.
\bibitem{MirCount} Maryam Mirzakhani, \emph{Growth of the number of simple closed geodesics on hyperbolic surfaces}, Ann. of Math. (2) \textbf{168} (2008), no.~1, 97--125. \url{https://doi.org/10.4007/annals.2008.168.97}.
\bibitem{Mir} Maryam Mirzakhani, \emph{Growth of Weil--Petersson volumes and random hyperbolic surface of large genus}, J. Differential Geom. \textbf{94} (2013), no.~2, 267--300. \url{https://doi.org/10.4310/jdg/1367438650}; \href{https://arxiv.org/abs/1012.2167}{arXiv:1012.2167}.
\bibitem{MirICM} Maryam Mirzakhani, \emph{On Weil--Petersson volumes and geometry of random hyperbolic surfaces}, in: Proceedings of the International Congress of Mathematicians, Vol.~II (Hyderabad, 2010), Hindustan Book Agency, New Delhi, 2010, 1126--1145. \url{https://doi.org/10.1142/9789814324359_0089}.
\bibitem{MP} Maryam Mirzakhani and Bram Petri, \emph{Lengths of closed geodesics on random surfaces of large genus}, Comment. Math. Helv. \textbf{94} (2019), no.~4, 869--889. \url{https://doi.org/10.4171/CMH/477}; \href{https://arxiv.org/abs/1710.09727v2}{arXiv:1710.09727v2}.
\bibitem{MZ} Maryam Mirzakhani and Peter Zograf, \emph{Towards large genus asymptotics of intersection numbers on moduli spaces of curves}, Geom. Funct. Anal. \textbf{25} (2015), no.~4, 1258--1289. \url{https://doi.org/10.1007/s00039-015-0336-5}; \href{https://arxiv.org/abs/1112.1151}{arXiv:1112.1151}.
\bibitem{Mondal} Sugata Mondal, \emph{On topological upper-bounds on the number of small cuspidal eigenvalues of hyperbolic surfaces of finite area}, Int. Math. Res. Not. IMRN \textbf{2015} (2015), no.~24, 13208--13237. \url{https://doi.org/10.1093/imrn/rnv092}; \href{https://arxiv.org/abs/1406.1076}{arXiv:1406.1076}.
\bibitem{Monk} Laura Monk, \emph{Benjamini--Schramm convergence and spectra of random hyperbolic surfaces of high genus}, Anal. PDE \textbf{15} (2022), no.~3, 727--752. \url{https://doi.org/10.2140/apde.2022.15.727}; \href{https://arxiv.org/abs/2002.00869}{arXiv:2002.00869}.
\bibitem{MonkThomas} Laura Monk and Joe Thomas, \emph{The tangle-free hypothesis on random hyperbolic surfaces}, Int. Math. Res. Not. IMRN \textbf{2022} (2022), no.~22, 18154--18185. \url{https://doi.org/10.1093/imrn/rnab160}; \href{https://arxiv.org/abs/2008.09363}{arXiv:2008.09363}.
\bibitem{MonkNaud} Laura Monk and Fr\'ed\'eric Naud, \emph{Spectral gaps on large hyperbolic surfaces}, in: Proceedings of the International Congress of Mathematicians 2026, Vol.~3: Invited Lectures (Sections 1--4), SIAM, Philadelphia, PA, 2026, 434--453. \url{https://doi.org/10.1137/25m1805527}.
\bibitem{MS} Motohico Mulase and Brad Safnuk, \emph{Mirzakhani's recursion relations, Virasoro constraints and the KdV hierarchy}, Indian J. Math. \textbf{50} (2008), no.~1, 189--218; \href{https://arxiv.org/abs/math/0601194v3}{arXiv:math/0601194v3}.
\bibitem{NPS} Hossein Namazi, Pekka Pankka, and Juan Souto, \emph{Distributional limits of Riemannian manifolds and graphs with sublinear genus growth}, Geom. Funct. Anal. \textbf{24} (2014), no.~1, 322--359. \url{https://doi.org/10.1007/s00039-014-0259-6}; \href{https://arxiv.org/abs/1208.2130v2}{arXiv:1208.2130v2} (theorem numbers cited from this version).
\bibitem{NaudDet} Fr\'ed\'eric Naud, \emph{Determinants of Laplacians on random hyperbolic surfaces}, J. Anal. Math. \textbf{151} (2023), no.~1, 265--291. \url{https://doi.org/10.1007/s11854-023-0334-8}; \href{https://arxiv.org/abs/2301.09965v4}{arXiv:2301.09965v4}.
\bibitem{NWX} Xin Nie, Yunhui Wu, and Yuhao Xue, \emph{Large genus asymptotics for lengths of separating closed geodesics on random surfaces}, J. Topol. \textbf{16} (2023), no.~1, 106--175. \url{https://doi.org/10.1112/topo.12276}; \href{https://arxiv.org/abs/2009.07538}{arXiv:2009.07538}.
\bibitem{OS} Kazumi Okuyama and Kazuhiro Sakai, \emph{JT gravity, KdV equations and macroscopic loop operators}, J. High Energy Phys. \textbf{2020} (2020), no.~1, article~156. \url{https://doi.org/10.1007/JHEP01(2020)156}; \href{https://arxiv.org/abs/1911.01659v2}{arXiv:1911.01659v2}.
\bibitem{OtalRosas} Jean-Pierre Otal and Eulalio Rosas, \emph{Pour toute surface hyperbolique de genre $g$, $\lambda_{2g-2}>1/4$}, Duke Math. J. \textbf{150} (2009), no.~1, 101--115. \url{https://doi.org/10.1215/00127094-2009-048}; \href{https://arxiv.org/abs/0905.0506v1}{arXiv:0905.0506v1}.
\bibitem{Patterson} Samuel J. Patterson, \emph{The limit set of a Fuchsian group}, Acta Math. \textbf{136} (1976), 241--273. \url{https://doi.org/10.1007/BF02392046}.
\bibitem{Penner} R. C. Penner, \emph{Weil--Petersson volumes}, J. Differential Geom. \textbf{35} (1992), no.~3, 559--608. \url{https://doi.org/10.4310/jdg/1214448257}.
\bibitem{Petri} Bram Petri, \emph{Bass notes of random hyperbolic surfaces of large genus}, S\'eminaire Bourbaki, Expos\'e no.~1255, preprint (2026), \href{https://arxiv.org/abs/2607.06331v1}{arXiv:2607.06331v1}.
\bibitem{PhillipsSarnak} Ralph S. Phillips and Peter Sarnak, \emph{On cusp forms for co-finite subgroups of $\mathrm{PSL}(2,\mathbb R)$}, Invent. Math. \textbf{80} (1985), no.~2, 339--364. \url{https://doi.org/10.1007/BF01388610}.
\bibitem{RZ} Kasra Rafi and Anton Zorich, \emph{Asymptotics of Masur--Veech volume of strata of quadratic differentials with many simple poles}, preprint (2026), \href{https://arxiv.org/abs/2610.28963v1}{arXiv:2610.28963v1}.
\bibitem{Ren} Irene Ren, \emph{Mirzakhani's frequencies of simple closed geodesics on hyperbolic surfaces in large genus and with many cusps}, Geom. Dedicata \textbf{219} (2025), no.~1, article~14. \url{https://doi.org/10.1007/s10711-024-00972-6}; \href{https://arxiv.org/abs/2304.03734}{arXiv:2304.03734}.
\bibitem{Rivin} Igor Rivin, \emph{Geodesics with one self-intersection, and other stories}, Adv. Math. \textbf{231} (2012), no.~5, 2391--2412. \url{https://doi.org/10.1016/j.aim.2012.07.018}; \href{https://arxiv.org/abs/0901.2543v4}{arXiv:0901.2543v4}.
\bibitem{Rudnick} Ze\'ev Rudnick, \emph{GOE statistics on the moduli space of surfaces of large genus}, Geom. Funct. Anal. \textbf{33} (2023), no.~6, 1581--1607. \url{https://doi.org/10.1007/s00039-023-00655-6}; \href{https://arxiv.org/abs/2202.06379}{arXiv:2202.06379}.
\bibitem{RudnickWigman} Ze\'ev Rudnick and Igor Wigman, \emph{On the central limit theorem for linear eigenvalue statistics on random surfaces of large genus}, J. Anal. Math. \textbf{151} (2023), no.~1, 293--302. \url{https://doi.org/10.1007/s11854-023-0327-7}.
\bibitem{SSS} Phil Saad, Stephen H. Shenker, and Douglas Stanford, \emph{JT gravity as a matrix integral}, preprint (2019), \href{https://arxiv.org/abs/1903.11115v1}{arXiv:1903.11115v1}.
\bibitem{SarnakDet} Peter Sarnak, \emph{Determinants of Laplacians}, Comm. Math. Phys. \textbf{110} (1987), no.~1, 113--120. \url{https://doi.org/10.1007/BF01209019}.
\bibitem{SWY} Richard Schoen, Scott A. Wolpert, and Shing-Tung Yau, \emph{Geometric bounds on the low eigenvalues of a compact surface}, in: Geometry of the Laplace Operator, Proc. Sympos. Pure Math. \textbf{36}, Amer. Math. Soc., Providence, RI, 1980, 279--285. \url{https://doi.org/10.1090/pspum/036/573440}.
\bibitem{Schulze} Michael Schulze, \emph{On the resolvent of the Laplacian on functions for degenerating surfaces of finite geometry}, J. Funct. Anal. \textbf{236} (2006), no.~1, 120--160. \url{https://doi.org/10.1016/j.jfa.2006.01.005}; \href{https://arxiv.org/abs/math/0410434v1}{arXiv:math/0410434v1} (theorem numbers cited from this version).
\bibitem{ST} Georg Schumacher and Stefano Trapani, \emph{Estimates of Weil--Petersson volumes via effective divisors}, Comm. Math. Phys. \textbf{222} (2001), no.~1, 1--7. \url{https://doi.org/10.1007/s002200100493}.
\bibitem{SW} Yang Shen and Yunhui Wu, \emph{Arbitrarily small spectral gaps for random hyperbolic surfaces with many cusps}, Comm. Math. Phys. \textbf{406} (2025), no.~8, article~190. \url{https://doi.org/10.1007/s00220-025-05366-7}; \href{https://arxiv.org/abs/2203.15681v4}{arXiv:2203.15681v4}.
\bibitem{Shenker} Stephen H. Shenker, \emph{The strength of nonperturbative effects in string theory}, in: Random Surfaces and Quantum Gravity (Carg\`ese, 1990), NATO Adv. Sci. Inst. Ser. B Phys. \textbf{262}, Plenum, New York, 1991, 191--200. \url{https://doi.org/10.1007/978-1-4615-3772-4_12}.
\bibitem{Soares} Louis Soares, \emph{Hecke triangle groups, transfer operators and Hausdorff dimension}, Ann. Henri Poincar\'e \textbf{23} (2022), no.~4, 1239--1281. \url{https://doi.org/10.1007/s00023-021-01117-1}; \href{https://arxiv.org/abs/2005.11808v1}{arXiv:2005.11808v1}.
\bibitem{StanfordWitten} Douglas Stanford and Edward Witten, \emph{JT gravity and the ensembles of random matrix theory}, Adv. Theor. Math. Phys. \textbf{24} (2020), no.~6, 1475--1680. \url{https://doi.org/10.4310/ATMP.2020.v24.n6.a4}; \href{https://arxiv.org/abs/1907.03363}{arXiv:1907.03363}.
\bibitem{Sullivan} Dennis Sullivan, \emph{Related aspects of positivity in Riemannian geometry}, J. Differential Geom. \textbf{25} (1987), no.~3, 327--351. \url{https://doi.org/10.4310/jdg/1214440979}.
\bibitem{SullivanErgodic} Dennis Sullivan, \emph{On the ergodic theory at infinity of an arbitrary discrete group of hyperbolic motions}, in: Riemann Surfaces and Related Topics (Stony Brook, NY, 1978), Ann. of Math. Stud. \textbf{97}, Princeton University Press, Princeton, NJ, 1981, 465--496. \url{https://doi.org/10.1515/9781400881550-035}.
\bibitem{Temme} Nico M. Temme, \emph{Laplace type integrals: transformation to standard form and uniform asymptotic expansions}, Quart. Appl. Math. \textbf{43} (1985), no.~1, 103--123. \url{https://doi.org/10.1090/qam/782260}.
\bibitem{WhittFCLT} Ward Whitt, \emph{Proofs of the martingale FCLT}, Probab. Surv. \textbf{4} (2007), 268--302. \url{https://doi.org/10.1214/07-PS122}; \href{https://arxiv.org/abs/0712.1929v2}{arXiv:0712.1929v2}.
\bibitem{Witten} Edward Witten, \emph{Two-dimensional gravity and intersection theory on moduli space}, Surv. Differ. Geom. \textbf{1} (1991), 243--310. \url{https://doi.org/10.4310/SDG.1990.v1.n1.a5}.
\bibitem{WittenDef} Edward Witten, \emph{Matrix models and deformations of JT gravity}, Proc. R. Soc. A \textbf{476} (2020), no.~2244, article~20200582. \url{https://doi.org/10.1098/rspa.2020.0582}; \href{https://arxiv.org/abs/2006.13414v1}{arXiv:2006.13414v1}.
\bibitem{Wolpert} Scott A. Wolpert, \emph{On the Weil--Petersson geometry of the moduli space of curves}, Amer. J. Math. \textbf{107} (1985), no.~4, 969--997. \url{https://doi.org/10.2307/2374363}.
\bibitem{WolpertTorus} Scott A. Wolpert, \emph{On the K\"ahler form of the moduli space of once punctured tori}, Comment. Math. Helv. \textbf{58} (1983), no.~1, 246--256. \url{https://doi.org/10.1007/BF02564634}.
\bibitem{WolpertCusp} Scott A. Wolpert, \emph{Disappearance of cusp forms in special families}, Ann. of Math. (2) \textbf{139} (1994), no.~2, 239--291. \url{https://doi.org/10.2307/2946582}.
\bibitem{WolpertSpecI} Scott A. Wolpert, \emph{Spectral limits for hyperbolic surfaces, I}, Invent. Math. \textbf{108} (1992), no.~1, 67--89. \url{https://doi.org/10.1007/BF02100600}.
\bibitem{WolpertSpecII} Scott A. Wolpert, \emph{Spectral limits for hyperbolic surfaces, II}, Invent. Math. \textbf{108} (1992), no.~1, 91--129. \url{https://doi.org/10.1007/BF02100601}.
\bibitem{Wright} Alex Wright, \emph{A tour through Mirzakhani's work on moduli spaces of Riemann surfaces}, Bull. Amer. Math. Soc. (N.S.) \textbf{57} (2020), no.~3, 359--408. \url{https://doi.org/10.1090/bull/1687}.
\bibitem{WXPrime} Yunhui Wu and Yuhao Xue, \emph{Prime geodesic theorem and closed geodesics for large genus}, J. Eur. Math. Soc. (2025), published online first. \url{https://doi.org/10.4171/JEMS/1653}; \href{https://arxiv.org/abs/2209.10415}{arXiv:2209.10415}.
\bibitem{WX} Yunhui Wu and Yuhao Xue, \emph{Random hyperbolic surfaces of large genus have first eigenvalues greater than $\frac{3}{16}-\epsilon$}, Geom. Funct. Anal. \textbf{32} (2022), no.~2, 340--410. \url{https://doi.org/10.1007/s00039-022-00595-7}; \href{https://arxiv.org/abs/2102.05581}{arXiv:2102.05581}.
\bibitem{WXOptimal} Yunhui Wu and Yuhao Xue, \emph{Optimal lower bounds for first eigenvalues of Riemann surfaces for large genus}, Amer. J. Math. \textbf{144} (2022), no.~4, 1087--1114. \url{https://doi.org/10.1353/ajm.2022.0024}.
\bibitem{XZ} Dingkang Xiong and Jian Zhou, \emph{Eynard--Orantin topological recursion, moduli-space Virasoro constraints, and semisimple CohFTs}, preprint (2026), \href{https://arxiv.org/abs/2610.24980v1}{arXiv:2610.24980v1}.
\bibitem{ZogSpectral} Peter Zograf, \emph{Small eigenvalues of automorphic Laplacians in spaces of parabolic forms}, Zap. Nauchn. Sem. Leningrad. Otdel. Mat. Inst. Steklov. (LOMI) \textbf{134} (1984), 157--168; English transl. in J. Soviet Math. \textbf{36} (1987), no.~1, 106--114. \url{https://doi.org/10.1007/BF01104976}; \url{https://www.mathnet.ru/eng/znsl4746}.
\bibitem{Zog93} Peter Zograf, \emph{The Weil--Petersson volume of the moduli space of punctured spheres}, in: Mapping Class Groups and Moduli Spaces of Riemann Surfaces (G\"ottingen, 1991/Seattle, WA, 1991), Contemp. Math. \textbf{150}, Amer. Math. Soc., Providence, RI, 1993, 367--372. \url{https://doi.org/10.1090/conm/150/01300}.
\bibitem{Zog98} Peter Zograf, \emph{Weil--Petersson volumes of moduli spaces of curves and the genus expansion in two dimensional gravity}, preprint (1998), \href{https://arxiv.org/abs/math/9811026}{arXiv:math/9811026}.
\bibitem{Zog} Peter Zograf, \emph{On the large genus asymptotics of Weil--Petersson volumes}, preprint (2008), \href{https://arxiv.org/abs/0812.0544v2}{arXiv:0812.0544v2}.
\end{thebibliography}
\end{document}